\documentclass[12pt,twosides]{amsart}
\usepackage{amssymb,amsmath,amsthm, amscd, enumerate, mathrsfs}
\usepackage{graphicx, hhline}
\usepackage[all]{xy}
\usepackage{braket}
\usepackage[usenames]{color}
\usepackage{hyperref}
\usepackage{fancyhdr}
\usepackage{bm}
\usepackage{url}
\usepackage{enumitem}
\usepackage{comment}
\usepackage{tikz}
\usepackage{caption}
\usepackage[top=30truemm,bottom=30truemm,left=30truemm,right=30truemm]{geometry}

\usepackage{here}

\hypersetup{colorlinks=true}

\title{Normal stable degenerations of Noether-Horikawa surfaces}
\date{\today}

\author{Hiroto Akaike} 
\address{Institute of Mathematics, Tohoku University,6-3 Aramaki Aza, Aoba, Aoba Ward, Sendai, Miyagi 980-8578, Japan}
\email{hiroto.akaike.c6@tohoku.ac.jp}

\author{Makoto Enokizono} 
\address{Graduate School of Mathematical Sciences, University of Tokyo,
3-8-1 Komaba, Meguro-ku, Tokyo 153-8914, Japan}
\email{enokizono@g.ecc.u-tokyo.ac.jp}

\author{Masafumi Hattori} 
\address{School of 
Mathematical Sciences, University of Nottingham, Nottingham, England}
\email{masafumi.hattori@nottingham.ac.uk}

\author{Yuki Koto} 
\address{Institute of Mathematics, Academia Sinica, Astronomy-Mathematics Building, No.\ 1, Sec.\ 4, Roosevelt Road, Taipei 10617, Taiwan.}
\email{ykoto@gate.sinica.edu.tw}

\newcommand{\Supp}[0]{{\operatorname{Supp}}}

\newtheorem{thm}{Theorem}[section]

\newtheorem{lem}[thm]{Lemma}
\newtheorem{cor}[thm]{Corollary}
\newtheorem{prop}[thm]{Proposition}

\newtheorem{ques}[thm]{Question}

\theoremstyle{definition}
\newtheorem{defn}[thm]{Definition}
\newtheorem{lem-defn}[thm]{Lemma-Definition} 
\newtheorem{rem}[thm]{Remark}
\newtheorem{note}[thm]{Notation}
\newtheorem{const}[thm]{Construction}
\newtheorem{exam}[thm]{Example}

\newtheorem{ass}[thm]{Assumption}
\newtheorem*{ack}{Acknowledgments} 
 
\newtheorem*{b-divisor}{b-divisors} 
 
\newtheorem*{g-pair}{Generalized pairs} 
\newtheorem*{adj-g-pair}{Divisorial adjunction for generalized pairs} 
\newtheorem*{mmp-g-pair}{MMP for generalized pairs}

\newtheorem*{claim*}{Claim}

\def\C{\mathbb C}

\def\Q{\mathbb Q}
\def\Z{\mathbb Z}

\def\PP{\mathbb P}

\def\cA{\mathcal{A}}

\def\O{\mathcal{O}}

\def\t{t}

\def\x{x}

\newcommand{\cAT}{\cA_{\operatorname{T}}}
\newcommand{\cAR}{\cA_{\operatorname{R}}}
\newcommand{\cASR}{\cA_{\operatorname{SR}}}
\newcommand{\cAE}{\cA_{\operatorname{E}}}
\newcommand{\hor}{{\operatorname{hor}}}
\newcommand{\Ell}{{\operatorname{Ell}}}
\newcommand{\SR}{{\operatorname{SR}}}

\usepackage[colorinlistoftodos]{todonotes}
\DeclareFontFamily{U}{mathx}{}
\DeclareFontShape{U}{mathx}{m}{n}{<-> mathx10}{}
\DeclareSymbolFont{mathx}{U}{mathx}{m}{n}
\DeclareMathAccent{\widecheck}{0}{mathx}{"71}

\usepackage{pgfplots}
\pgfplotsset{compat=1.15}
\usepackage{mathrsfs}
\usetikzlibrary{arrows}

\usepackage{subcaption}
\begin{document}

\begin{abstract}
We classify all normal stable Horikawa surfaces with only $\mathbb{Q}$-Gorenstein smoothable log canonical singularities.
Furthermore, we provide a criterion for their global $\Q$-Gorenstein smoothability and describe the boundary strata of the moduli space of $\Q$-Gorenstein smoothable normal stable Horikawa surfaces.
\end{abstract}

\maketitle

\tableofcontents

\section{Introduction}

Noether-Horikawa surfaces, or simply Horikawa surfaces, are minimal surfaces of general type satisfying the numerical relation
$$
K_{X}^{2}=2p_g(X)-4.
$$
The aim of this paper is to classify all normal stable Horikawa surfaces with only $\Q$-Gorenstein smoothable log canonical singularities, and to describe the boundary of the moduli space of normal stable degenerations of smooth Horikawa surfaces.

\subsection{Background}\label{sec--background}

In algebraic geometry, the classification of projective varieties and their degenerations is a central but difficult problem. 
A key perspective from birational geometry is to understand varieties as built from three fundamental types: Fano, Calabi-Yau, and general type varieties. 
Among these, the classification of surfaces of general type is particularly challenging due to the infinite range of their numerical invariants. 
One classical tool is the Noether inequality, which asserts that
$$
K_{X}^{2}\ge 2p_g(X)-4
$$
for smooth minimal surfaces of general type.

Enriques \cite{Enriques} discovered examples of surfaces for which the above inequality becomes an equality.
In the 1970s, Horikawa conducted a detailed study of surfaces near the Noether line (\cite{horikawa,horikawa-ii,horikawa-iii,horikawa-iv,horikawa-v}), particularly those satisfying the equality in the above inequality. 
He classified such surfaces with only Du Val singularities, showing that they are often genus-two fibrations and can be constructed as double covers of surfaces of minimal degree with smooth or mildly singular branch divisors. His work also described their moduli spaces.

In parallel, advances in moduli theory---especially the pioneering work of Koll\'ar-Shepherd-Barron \cite{KSB} and Alexeev \cite{Al}---have led to the construction of the {\em KSBA moduli space}: 
a proper moduli space parametrizing $\Q$-Gorenstein deformations of {\em stable varieties}, that is, projective schemes with ample canonical class and only semi-log canonical (slc) singularities (see \cite{kollar-moduli} for an overview).
In particular, this moduli space compactifies the Gieseker moduli space \cite{gieseker} parametrizing canonical models of smooth surfaces of general type, which generalizes the Deligne-Mumford compactification of the moduli space of curves \cite{DM}.
The theory of KSBA moduli has become a central tool to study moduli spaces of higher dimensional varieties: not only of general type but also of broader classes, toric and Abelian varieties \cite{Alab}, $K3$ surfaces \cite{AET,Al-En,AE2,AEDGS,ADH}, Enriques surfaces \cite{AEDGS}, and elliptic surfaces \cite{ABI,AB,ABE,modulielliptic}, considering boundary divisors.
On the other hand, explicit descriptions of KSBA moduli spaces without boundary divisor remain rare, even in dimension two. 
Indeed, only a few examples are known: Campedelli and Burniat surfaces \cite{Al-Pa2,AH} and certain quotients or finite covers of products of curves \cite{vO,Liu2,Rol}.

These challenges naturally lead to the problem of studying the KSBA compactification of the moduli of Horikawa surfaces. 
This problem has been approached from various perspectives by many researchers, and there have been many significant developments since Horikawa's celebrated work on the moduli of Horikawa surfaces with only Du Val singularities.
For example, Lee-Park \cite{LP} constructed Kawamata log terminal (klt) Horikawa surfaces with non-Du Val singularities, and showed that they are $\Q$-Gorenstein smoothable. 
Anthes \cite{Anthes} and Rana-Rollenske \cite{RR} studied components of the KSBA moduli space parametrizing certain (not necessarily normal) Gorenstein stable Horikawa surfaces.
Recently, Evans-Simonetti-Urz\'ua \cite{ESU} investigated $\Q$-Gorenstein smoothable normal stable Horikawa surfaces with $p_g=3$, and Monreal-Negrete-Urz\'ua \cite{MNU} studied non-Gorenstein stable Horikawa surfaces with only T-singularities, that is, $\Q$-Gorenstein smoothable klt singularities.

Despite substantial progress made in recent years, the classification of stable Horikawa surfaces with log canonical (lc) singularities and the structure of their KSBA moduli space remain incomplete. 
One major difficulty lies in the complexity of the deformation theory of stable surfaces. 
For example, not all such surfaces are $\Q$-Gorenstein smoothable (cf.~\cite{RR}), and it remains challenging to control the singularities that appear in degenerations, despite general boundedness results \cite{Al,Lee,hmx-boundgentype}.

This paper aims to address this gap by providing a comprehensive description of the locus of normal surfaces in the KSBA compactification of the moduli space of Horikawa surfaces.

\subsection{Main results}

We begin by formulating the central problems addressed in this paper.

\begin{ques}\label{question}
    \phantom{A}
    \begin{itemize}
        \item[$(1)$] 
        Can one give a list of all normal stable Horikawa surfaces with only $\Q$-Gorenstein smoothable singularities?
        \item[$(2)$]
        Can one give a list of all $\Q$-Gorenstein smoothable normal stable Horikawa surfaces?
        \item[$(3)$]
        Can one describe a stratification of the KSBA boundary of the moduli space of normal stable Horikawa surfaces?
    \end{itemize}
\end{ques}

Note that the second question is different from the first question because stable surfaces with only singularities that are locally $\Q$-Gorenstein smoothable may not be globally $\Q$-Gorenstein smoothable.
The first question is answered by Theorems \ref{thm:standard_Horikawa} and \ref{thm:classification_non-standard_Horikawa_noninv}. 
We explain these results in Section \ref{intro-subsubsec:result_Q1}.
The answer to the second question (Theorem \ref{intro-thm:Q2}) is explained in Section \ref{intro-subsubsec:result_classification}, and the answer to the third question (Theorem \ref{thm--stratification--intro}) is given in Section \ref{intro-subsubsec:result_moduli}.

\subsubsection{Classification of normal stable Horikawa surfaces}\label{intro-subsubsec:result_Q1}

To explain our classification of normal stable Horikawa surfaces with only $\Q$-Gorenstein smoothable singularities, we divide the class of such surfaces into two types: those whose canonical maps give double covers (\emph{standard Horikawa surfaces}\footnote{Standard Horikawa surfaces are slightly different from \emph{standard stable Horikawa surfaces} introduced by Rana-Rollenske \cite{RR}.
The latter refers to (not necessarily normal) stable Horikawa surfaces that are obtained as double covers of Hirzebruch surfaces.}), and those whose canonical maps are composed with a pencil (\emph{non-standard Horikawa surfaces}).

Standard Horikawa surfaces are natural generalizations of the surfaces studied by Horikawa \cite{horikawa}, and were later investigated by Chen \cite{Che}.
In particular, the canonical map of a standard Horikawa surface gives rise to a double cover over a minimal degree surface (Proposition \ref{prop:standard_Horikawa}).
Following Horikawa's notation, we say that it is of type $(\infty)$ if the canonical image is $\PP^2$, of type $(d)$ if the canonical image is the $d$-th Hirzebruch surface $\Sigma_d$, and of type $(d)'$ if the canonical image is $\overline{\Sigma}_d$, a projective cone over a smooth rational curve of degree $d$.
Note that a geometric genus of a standard Horikawa surface of type $(d)'$ is $d+2$.
Furthermore, we classify all possible singularities of a branch divisor on a minimal degree surface whose associated double cover gives a normal stable Horikawa surface.
This completes the classification of standard Horikawa surfaces ({\cite[Theorem~4.1]{Che}}, Theorem~\ref{thm:standard_Horikawa}).

We now turn our attention to the classification of non-standard Horikawa surfaces.  
In Theorem \ref{thm:classification_non-standard_Horikawa_noninv}, we provide a list of all possible combinations of non-Du Val singularities that can appear on non-standard Horikawa surfaces with only $\Q$-Gorenstein smoothable singularities.
Two key tools of the proof are the log Noether inequality (Theorem \ref{normalstable}) and a combinatorial analysis of \emph{extended T-chains} (Appendix \ref{app:drill}).

\begin{rem}
    We make a few remarks on our classification of non-standard Horikawa surfaces (Theorem \ref{thm:classification_non-standard_Horikawa_noninv}).
    \begin{itemize}
        \item[$(1)$]
        Our classification extends the classification of non-Gorenstein stable Horikawa surfaces with only T-singularities given by Monreal-Negrete-Urz\'ua \cite[Theorem 4.3]{MNU}, which was established independently.
        The main difficulty in the lc case lies in the log Noether inequality.
        See Section \ref{intro-subsec:classification} for details.
        \item[$(2)$]
        Except for a few cases, we do not verify the existence of a non-standard Horikawa surface whose set of non-Du Val singularities coincides with a given candidate in the list provided in Theorem \ref{thm:classification_non-standard_Horikawa_noninv}. 
        Nevertheless, we provide a recipe for constructing each candidate from an elliptic fibration with specified Euler characteristic and prescribed singular fibers.
        In particular, we observe that there exists a non-standard Horikawa surface with exactly one strictly lc singularity of type $(2,4,4)[3]$; see Lemma-Definition \ref{lem:smoothable_rational_strict_lc} (iii) for notation, and Remark \ref{rem:triangle_Horikawa} for the proof.
        This seems to be the first example of a stable Horikawa surface with a strictly lc singularity of type $(2,4,4)[3]$.
        \item[$(3)$] 
        Theorem \ref{thm:classification_non-standard_Horikawa_noninv} provides valuable insights into components of the KSBA moduli space of stable Horikawa surfaces beyond the Gieseker component.
        We expect that, once the classification of lc singularities that can appear on surfaces parametrized by the connected component containing the Gieseker component is established, our classification techniques can be applied to study other components of the KSBA moduli space.
    \end{itemize}
\end{rem}

\subsubsection{Classification of $\Q$-Gorenstein smoothable normal stable Horikawa surfaces}\label{intro-subsubsec:result_classification}

The following theorem provides a necessary condition for a normal stable Horikawa surface to be $\Q$-Gorenstein smoothable.
The strictly lc singularity of type $(2,2,2,2)$ appearing in the statement is introduced in Lemma-Definition~\ref{lem:smoothable_rational_strict_lc}~(i).

\begin{thm}\label{intro-qGsmHor}
Let $X$ be a $\Q$-Gorenstein smoothable normal stable Horikawa surface.
Then exactly one of the following holds:

\begin{itemize}
    
\item[$(1)$]
$X$ is Gorenstein and the canonical map gives a double covering $\varphi_{K_X}\colon X\to W$ over a minimal degree surface $W \subset \mathbb{P}^{p_g-1}$. 

\item[$(2)$]
$(p_g=4)$
$X$ is Gorenstein and the canonical linear system is composed with a pencil.
In this case, the bicanonical map gives a double covering $\varphi_{2K_X}\colon X\to W$ over an elliptic cone $W \subset \mathbb{P}^{8}$ of degree $8$.

\item[$(3)$]
$X$ has two T-singularities of type $\frac{1}{(p_g-1)^2}(1,p_g-2)$ and no other non-Gorenstein singularities.

\item[$(4)$]
$(p_g=3)$
$X$ has a T-singularity of type $\frac{1}{50}(1,29)$ and no other non-Gorenstein singularities.

\item[$(5)$]
$(p_g=3)$
$X$ has a strictly lc singularity of type $(2,2,2,2)$ and no other non-Gorenstein singularities.

\item[$(6)$]
$(p_g=3)$
$X$ has a T-singularity of type $\frac{1}{50}(1,29)$, two T-singularities of type $\frac{1}{4}(1,1)$ and no other non-Gorenstein singularities.

\item[$(7)$]
$(p_g=3)$
$X$ has a T-singularity of type $\frac{1}{50}(1,29)$, a strictly lc singularity of type $(2,2,2,2)$ and no other non-Gorenstein singularities.
\end{itemize}
\end{thm}

This theorem follows immediately from Theorem \ref{thm:standard_Horikawa} and Theorem \ref{thm:classification_non-standard_Horikawa}.
The proof mainly relies on the classification of non-standard Horikawa surfaces (Theorem \ref{thm:classification_non-standard_Horikawa_noninv}) and on further combinatorial analysis with respect to certain involutions. 
We now give a brief description of each surface in the above classification list. 
\begin{itemize}
    \item[$(1)$]
    These are exactly standard Horikawa surfaces explained in Section \ref{intro-subsubsec:result_Q1}.
    In contrast, surfaces in (2)--(7) are non-standard.
    \item[$(2)$]
    These have at least one simple elliptic singularity of degree $4$.
    These surfaces appear only when $p_g=4$.
    In cases (3)--(7), $X$ is not Gorenstein.
    \item[$(3)$]
    These are generalizations of the examples found by Lee-Park \cite{LP}, which had been constructed by Fintushel-Stern \cite{FS} as rational blowdowns of smooth fourfolds in differential topology.
    We refer to these as \emph{Horikawa surfaces of Lee-Park type}.
    These surfaces form a boundary divisor in the moduli space.
    \item[$(4)$]
    These surfaces form a boundary divisor in the moduli space.
    \item[$(5)$]
    These arise as degenerations of those in case (3), and form a codimension $2$ locus in the moduli space.
    \item[$(6)$]
    These are common degenerations of surfaces in cases (3) and (4).
    \item[$(7)$]
    These are common degenerations of surfaces in cases (5) and (6), and form a codimension $3$ locus in the moduli space.
\end{itemize}
Each of these surfaces can be constructed by birationally modifying a double cover of a Hirzebruch surface or $\mathbb{P}^2$ with a specially chosen branch divisor (see Theorem~\ref{thm:standard_Horikawa} for case~$(1)$ and Construction~\ref{construction} for the other cases).
Note that surfaces in case (2) are first discovered and studied in this paper.
We also remark that some of these had already been found in \cite{Anthes,RR,ESU,MNU}.
See Remark \ref{intro-rem:other_works_Q2} for the relationship between recent research and Theorem \ref{intro-qGsmHor}.

In fact, Theorem \ref{intro-qGsmHor} is a classification of normal stable Horikawa surfaces that admit a \emph{good involution} (Definition~\ref{defn--good-involution}) and have only $\Q$-Gorenstein smoothable singularities.
The notion of a good involution is defined specifically for our purposes and has the following elegant properties (Proposition~\ref{prop:good_involution_uniqueness}~$(1)$, Proposition~\ref{prop--involution}): 
\begin{itemize}
    \item
    For any normal stable Horikawa surface $X$ with only $\Q$-Gorenstein smoothable singularities, if a good involution on $X$ exists, then it is unique. 
    \item 
    Every $\Q$-Gorenstein smoothable normal stable Horikawa surface admits a unique good involution.
\end{itemize}
By construction, each surface in the list is naturally equipped with an involution, which is in fact a good involution.
However, we remark that the existence of the good involution on $X$ does not imply that $X$ is $\Q$-Gorenstein smoothable.
Therefore, to answer Question~\ref{question}~(2), we need to study the $\Q$-Gorenstein smoothability of Horikawa surfaces listed in Theorem \ref{intro-qGsmHor}.

To simplify the discussion, we focus here on a standard Horikawa surface $X$. 
By assumption, $X$ is equipped with the good involution $\sigma$ induced by the canonical map, and the quotient $W:=X/\sigma$ is either a projective plane, a Hirzebruch surface, or a cone over a rational normal curve.
We divide the singularities on $X$ into two types according to whether their images in $W$ are smooth: those whose images are smooth are called \emph{mild singularities} (Definition~\ref{defn--mild--sing}), and those whose images are cone singularities are called \emph{double cone singularities} (Definition~\ref{defn:(d;k_1,k_2)_sing}).
Note that $X$ has at most one double cone singularity because the quotient $W$ has at most one cone singularity.
Since $X$ is Gorenstein, a double cone singularity on $X$ is either a Du Val singularity or an elliptic singularity; see Proposition \ref{prop:classification_cone_sing} for a more detailed classification.
Note that among the cases listed in Theorem \ref{intro-qGsmHor}, only cases $(1)$ and $(3)$ may admit an elliptic double cone singularity.
The smoothing problem of $X$ with no elliptic double cone singularities was already studied; see e.g. \cite{horikawa,ESU,MNU}.
In this paper, we investigate the case where $X$ has an elliptic double cone singularity.

The following theorem reduces the \emph{global} $\Q$-Gorenstein smoothing problem for normal stable Horikawa surfaces to the study of the \emph{local} equivariant smoothability of elliptic double cone singularities. 

\begin{thm}\label{thm--main--ii}
    Let $X$ be a normal stable Horikawa surface with only $\Q$-Gorenstein smoothable singularities that admits the good involution.
    Then, $X$ is $\Q$-Gorenstein smoothable if and only if all elliptic double cone singularities on $X$ are equivariantly smoothable with respect to the good involution.
    In particular, if $X$ has no elliptic double cone singularity, then $X$ is $\Q$-Gorenstein smoothable.
\end{thm}
We prove Theorem \ref{thm--main--ii} by considering the deformation property of the log pair $(W,\frac{1}{2}B)$, where $W$ is the quotient of $X$ by the good involution and $B$ is the associated branch divisor, instead of the Horikawa surface $X$ itself.
The main tools for proving Theorem \ref{thm--main--ii} are the vanishing theorems of Fujino \cite{fujino--slc--vanishing} and the second author \cite{Enokizono}, deformation theories of divisorial sheaves under $\mathbb{Q}$-Gorenstein deformations of T-singularities, and an \emph{anti-P-resolution}, a new notion of a birational transformation introduced in this paper.
We refer to Section \ref{intro-subsec:smoothing_and_anti_P} for the definition of anti-P-resolutions and a detailed discussion.

We note that all elliptic double cone singularities over the cone singularity of type $\frac{1}{p_g-2}(1,1)$ with $p_g\le 6$ are equivariantly smoothable,
and that for each $p_g\ge 7$, there exists an equivariantly smoothable elliptic double cone singularity over the cone singularity of type $\frac{1}{p_g-2}(1,1)$ (cf. Section~\ref{subsec:conclusion--equiv--smoothability}). 

Combining Theorem~\ref{intro-qGsmHor} with Theorem~\ref{thm--main--ii}, we obtain a classification of $\Q$-Gorenstein smoothable normal stable Horikawa surfaces, that is, the answer to Question \ref{question} (2).

\begin{thm}\label{intro-thm:Q2}
Let $X$ be a normal stable Horikawa surface with only $\Q$-Gorenstein smoothable singularities that admits a good involution.
Then exactly one of the following holds, where the labels $(1)$\text{--}$(7)$ correspond to those in Theorem~\ref{intro-qGsmHor}.
\begin{itemize}
    
\item[$(1)$]
$X$ has some mild singularities and at most one double cone singularity. 
Furthermore, $X$ is smoothable if and only if either $X$ has no elliptic double cone singularity, or the elliptic double cone singularity on $X$ is locally equivariantly smoothable with respect to the good involution.

\item[$(2)$]
$X$ has a simple elliptic singularity of degree $4$ and some mild singularities.
Furthermore, $X$ is smoothable.

\item[$(3)$]
$X$ has two T-singularities of type $\frac{1}{(p_g-1)^2}(1,p_g-2)$, some mild singularities and at most one double cone singularity.
Furthermore, $X$ is $\Q$-Gorenstein smoothable if and only if either $X$ has no elliptic double cone singularity, or the elliptic double cone singularity on $X$ is locally and equivariantly smoothable with respect to the good involution.

\item[$(4)$]
$X$ has a T-singularity of type $\frac{1}{50}(1,29)$ and some mild singularities.
Furthermore, $X$ is $\Q$-Gorenstein smoothable.

\item[$(5)$]
$X$ has a strictly lc singularity of type $(2,2,2,2)$ and some mild singularities.
Furthermore, $X$ is $\Q$-Gorenstein smoothable.

\item[$(6)$]
$X$ has a T-singularity of type $\frac{1}{50}(1,29)$, two T-singularities of type $\frac{1}{4}(1,1)$ and some mild singularities.
Furthermore, $X$ is $\Q$-Gorenstein smoothable.

\item[$(7)$]
$X$ has a T-singularity of type $\frac{1}{50}(1,29)$, a strictly lc singularity of type $(2,2,2,2)$ and some mild singularities.
Furthermore, $X$ is $\Q$-Gorenstein smoothable.
\end{itemize}
\end{thm}

\begin{rem}[Recent research]\label{intro-rem:other_works_Q2}
    While we were writing this paper, two notable developments were made by Evans-Simonetti-Urz\'ua \cite{ESU} and Monreal-Negrete-Urz\'ua \cite{MNU}.
    Including these works, we describe how recent research relates to Theorems \ref{intro-qGsmHor}, \ref{intro-thm:Q2}.
    \begin{itemize}
        \item[(1)]
        Anthes \cite{Anthes} studied (not necessarily normal) Gorenstein stable Horikawa surfaces with $p_g=3$, which can be realized as double covers of $\PP^2$ branched along octics via the canonical linear systems on them \cite[Corollary 2.7]{Anthes}. 
        Such surfaces are standard Horikawa surfaces if they are normal.
        \item[(2)]
        Rana-Rollenske \cite{RR} studied (not necessarily normal) stable Horikawa surfaces that arise as double covers of Hirzebruch surfaces branched along certain divisors.
        Such surfaces are standard if they are normal.
        Note that the quotient of a standard Horikawa surface may not be a Hirzebruch surface.
        \item[(3)] 
        Evans-Simonetti-Urz\'ua \cite{ESU} have developed a tropical-geometric approach and investigated $\Q$-Gorenstein smoothable normal stable Horikawa surfaces with $p_g=3$.
        In particular, they achieved Theorem \ref{intro-thm:Q2} for the case $p_g=3$.
        Note that the condition on equivariant smoothability appearing in Theorem \ref{intro-thm:Q2} becomes vacuous when $p_g=3$. 
        \item[(4)]
        Motivated by the \emph{Horikawa problem} (see Section~\ref{subsec--other--relate} for details), Monreal-Negrete-Urz\'ua \cite{MNU} classified stable non-Gorenstein Horikawa surfaces with only T-singularities and studied their $\Q$-Gorenstein smoothability.
        They showed that no klt stable non-Gorenstein Horikawa surfaces with $p_g\geq10$ except \emph{Lee-Park examples}\footnote{
            Monreal-Negrete-Urz\'ua \cite{MNU} introduced what they call Lee-Park examples: klt surfaces with two Wahl singularities admitting very similar configurations to Lee-Park type Horikawa surfaces.
            Note that Lee-Park examples do not require the existence of a good involution.} 
        are $\Q$-Gorenstein smoothable.
        Theorem \ref{intro-qGsmHor} refines their result by showing that among normal stable non-Gorenstein Horikawa surfaces with $p_g\geq4$, only those of Lee-Park type can be $\Q$-Gorenstein smoothable.
        Moreover, Theorem \ref{intro-thm:Q2} provides an equivalent condition for Lee-Park type Horikawa surfaces to be $\Q$-Gorenstein smoothable in terms of local conditions for any $p_g\geq3$.
        In particular, we can deduce that all klt Horikawa surfaces of Lee-Park type are $\Q$-Gorenstein smoothable.
    \end{itemize} 
\end{rem}

In a different direction, we also highlight the following interesting observation. 

\begin{thm}\label{intro-thm:involution}
    Every Gorenstein normal stable Horikawa surface admits the good involution.
\end{thm}

This follows from the log Noether inequality---more precisely, from its consequence Proposition~\ref{prop:standard_Horikawa}---and from Proposition~\ref{prop--canonical--pencil--invol}.

\subsubsection{Boundary of the KSBA moduli}\label{intro-subsubsec:result_moduli}
To describe a stratification of the KSBA boundary, a finer analysis of Theorem \ref{intro-thm:Q2} is required.
We introduce subclasses of the standard Horikawa surfaces (Theorem \ref{intro-thm:Q2} (1)) and the Horikawa surfaces of Lee-Park type (Theorem \ref{intro-thm:Q2} (3)) into several types, respectively.

We first discuss Horikawa surfaces of Lee-Park type. 
They can be subdivided into three classes according to the quotient $W=X/\sigma$ by the good involution: 
those without double cone singularities and with $\rho(W)=2$ (\emph{general Lee-Park type}), those with Gorenstein double cone singularities (\emph{special Lee-Park type}),
and those without double cone singularities and with $\rho(W)=1$ (\emph{infinite Lee-Park type}).
This subdivision corresponds to the subdivision of standard Horikawa surfaces into of type $(d)$, of type $(p_g-2)'$, and of type $(\infty)$.
If $p_g\geq4$, these are further subdivided according to the branch divisors in the quotients, or, equivalently, according to the configurations of the exceptional curves of non-Du Val singularities.

\begin{exam}[general Lee-Park type]\label{exam--lp--general}
    Let $X$ be a Horikawa surface of general Lee-Park type with $p_g=p_g(X)\geq 4$.
    We provide two types of constructions of $X$.
    
    By Theorem \ref{intro-thm:Q2}, $X$ has exactly two T-singularities, whose exceptional curves form the chain $[p_g+1,2^{p_g-3}]$. 
    Here, $2^{p_g-3}$ in the string represents that the number $2$ appears consecutively $p_g-3$ times at that position.
    Futhermore, the quotient of $X$ by the good involution has a unique singularity, which is a T-singularity (away from the branch divisor) corresponding to the chain $[p_g+1,2^{p_g-3}]$.
    We follow the Kodaira notation for singular fibers of relatively minimal elliptic fibrations.
    \begin{itemize}
        \item[(1)]
        Let $\widetilde{X}\to\PP^1$ be a relatively minimal elliptic fibration with a specific singular fiber $F$ of type $\mathrm{I}_n$ ($n\geq0$) and a fiberwise involution $\sigma$.
        Suppose there exist two distinct sections $S_1$ and $S_2$, each with self-intersection number $-(p_g+1)$, such that $\sigma$ permutes $S_1$ and $S_2$.
        As shown in Figure \ref{fig:LP_gen} (1), we can take two distinct chains $[p_g+1,2^{p_g-3}]$ in $F+S_1+S_2$.
        We impose the condition that $\sigma$ interchanges the two chains.
        Finally, we define $X$ to be a Horikawa surface obtained from $\widetilde{X}$ by contracting the two chains and the components of fibers disjoint from these chains.
        Such $X$ is called a \emph{Horikawa surface of general Lee-Park type $\mathrm{I}$}. 
        \item[(2)]
        Let $p_g\geq5$.
        We carry out the same construction as (1) by replacing $F$ with a singular fiber of type $\mathrm{I}^*_n$.
        The two distinct chains of the form $[p_g+1,2^{p_g-3}]$ in $F+S_1+S_2$ are arranged as shown in Figure \ref{fig:LP_gen} (2).
        Such $X$ is called a \emph{Horikawa surface of general Lee-Park type $\mathrm{I}^*$}.  
    \end{itemize}
    
We note that the examples found by Lee-Park \cite{LP} are included in the class of Horikawa surfaces of general Lee-Park type I.
\end{exam}

\begin{exam}[Horikawa surfaces of infinite Lee-Park type]
    We set $p_g=6$, and carry out the same construction as Example \ref{exam--lp--general} (2) by replacing the two chains with those shown in Figure \ref{fig:LP_gen} (3).
    Then, the resulting surface is of infinite Lee-Park type.
    Note that such a surface was first found in \cite[Example 4.4]{MNU}.
\end{exam}

\begin{exam}[Horikawa surfaces of special Lee-Park type]
    Let $X$ be a Horikawa surface of special Lee-Park type with $p_g=p_g(X)\geq4$.
    By Theorem \ref{intro-thm:Q2}, $X$ has two T-singularities corresponding to $[p_g+1,2^{p_g-3}]$ and one Gorenstein double cone singularity.
    The quotient of $X$ by the good involution has exactly two singularities: a T-singularity corresponding to $[p_g+1,2^{p_g-3}]$ and a cone singularity.
    Although $X$ may have a Du Val double cone singularity when $p_g\le6$, we focus here on the case where $X$ has an elliptic double cone singularity.
    
    Let $Y\to\PP^1$ be a relatively minimal elliptic fibration with a specific singular fiber $F$ of type $\mathrm{I}_n$ ($n\geq0$) and a fiberwise involution $\sigma$.
    Suppose there exist two distinct sections $S_1$ and $S_2$, each with self-intersection number $-p_g$, such that $\sigma$ permutes $S_1$ and $S_2$.
    For each section, we perform $(p_g-2)$ successive blow-ups at the intersection points of the proper transforms of $F$ and the total transforms of the section, which yields a surface $\widetilde{X}$.
    Finally, we define $X$ to be a Horikawa surface obtained from $\widetilde{X}$ by contracting all solid lines in Figure~\ref{fig:LP_special}, which are the union of the two distinct chains of the form $[p_g+1,2^{p_g-3}]$ and the proper transform of $F$, and the components of fibers disjoint from these chains.
    Such $X$ is called a \emph{Horikawa surface of special Lee-Park type $\mathrm{I}$}.


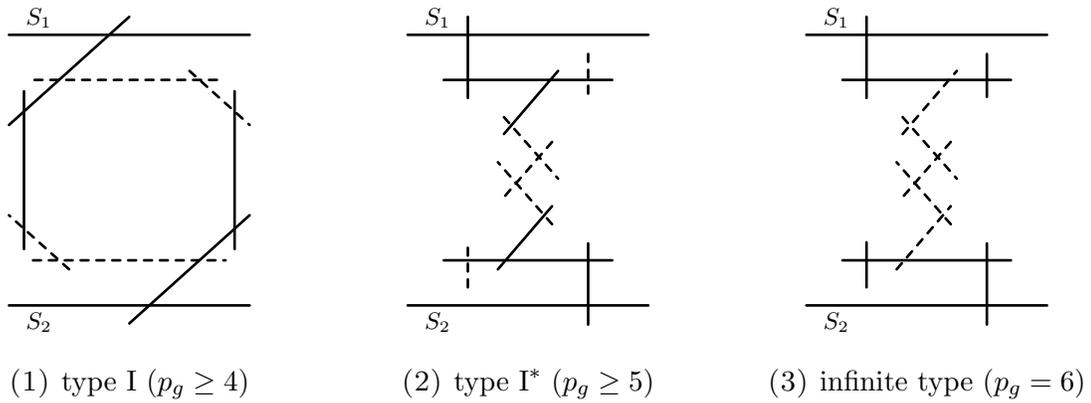
\begin{figure}[H]
\centering
\begin{subfigure}[t]{0.32\linewidth}
\centering
\begin{tikzpicture}[line cap=round,line join=round,>=triangle 45,x=0.4cm,y=0.3cm]
\clip(-5,-8) rectangle (5,8);
\draw [line width=1pt] (-4,6)-- (4,6);
\draw [line width=1pt] (-4,-6)-- (4,-6);
\draw [line width=1pt,] (0,6.8)-- (-4,2);
\draw [line width=1pt] (4,-2)-- (0,-6.8);
\draw [line width=1pt,dashed] (-3.15,4)-- (3.15,4);
\draw [line width=1pt,dashed] (2,4.4)-- (4,2);
\draw [line width=1pt,dashed] (-2,-4.4)-- (-4,-2);
\draw [line width=1pt,dashed] (-3.2,-4)-- (3.2,-4);
\draw [line width=1pt] (3.5,-3.5)-- (3.5,3.5);
\draw [line width=1pt] (-3.5,-3.5)-- (-3.5,3.5);
\begin{scriptsize}
\draw [color=black] (-3,6.75) node {$S_1$};
\draw [color=black] (-3,-6.75) node {$S_2$};
\end{scriptsize}
\end{tikzpicture}
\caption{type I ($p_g\geq4$)}
\label{fig:type_I}
\end{subfigure}
\hfill
\begin{subfigure}[t]{0.32\linewidth}
\centering
\begin{tikzpicture}[line cap=round,line join=round,>=triangle 45,x=0.4cm,y=0.3cm]
\clip(-5,-8) rectangle (5,8);
\draw [line width=1pt] (-4,6)-- (4,6);
\draw [line width=1pt] (-4,-6)-- (4,-6);
\draw [line width=1pt] (-2,6.8)-- (-2,3.2);
\draw [line width=1pt] (-2.8,4)-- (2.8,4);
\draw [line width=1pt,dashed] (2,5.15)-- (2,3.25);
\draw [line width=1pt] (2,-6.85)-- (2,-3.25);
\draw [line width=1pt] (-2.8,-4)-- (2.8,-4);
\draw [line width=1pt,dashed] (-2,-5.2)-- (-2,-3.2);
\draw [line width=1pt] (-0.8,1.6)-- (1,4.4);
\draw [line width=1pt] (-1,-4.4)-- (0.8,-1.6);
\draw [line width=1pt,dashed] (-0.8,2.35)-- (1,-0.38);
\draw [line width=1pt,dashed] (0.8,-2.4)-- (-1,0.35);
\draw [line width=1pt,dashed] (0.8,1.25)-- (-0.8,-1.2);
\begin{scriptsize}
\draw [color=black] (-3,6.75) node {$S_1$};
\draw [color=black] (-3,-6.75) node {$S_2$};
\end{scriptsize}
\end{tikzpicture}
\caption{type $\textrm{I}^*$ ($p_g\geq5$)}
\label{fig:type_I^*}
\end{subfigure}
\hfill
\begin{subfigure}[t]{0.32\linewidth}
\centering
\begin{tikzpicture}[line cap=round,line join=round,>=triangle 45,x=0.4cm,y=0.3cm]
\clip(-5,-8) rectangle (5,8);
\draw [line width=1pt] (-4,6)-- (4,6);
\draw [line width=1pt] (-4,-6)-- (4,-6);
\draw [line width=1pt] (-2,6.8)-- (-2,3.2);
\draw [line width=1pt] (-2.8,4)-- (2.8,4);
\draw [line width=1pt] (2,5.15)-- (2,3.25);
\draw [line width=1pt] (2,-6.85)-- (2,-3.25);
\draw [line width=1pt] (-2.8,-4)-- (2.8,-4);
\draw [line width=1pt] (-2,-5.2)-- (-2,-3.2);
\draw [line width=1pt,dashed] (-0.8,1.6)-- (1,4.4);
\draw [line width=1pt,dashed] (-1,-4.4)-- (0.8,-1.6);
\draw [line width=1pt,dashed] (-0.8,2.35)-- (1,-0.38);
\draw [line width=1pt,dashed] (0.8,-2.4)-- (-1,0.35);
\draw [line width=1pt,dashed] (0.8,1.25)-- (-0.8,-1.2);
\begin{scriptsize}
\draw [color=black] (-3,6.75) node {$S_1$};
\draw [color=black] (-3,-6.75) node {$S_2$};
\end{scriptsize}
\end{tikzpicture}
\caption{infinite type ($p_g=6$)}
\label{fig:type_infty}
\end{subfigure}
\caption{Horikawa surfaces of general/infinite Lee-Park type}
\label{fig:LP_gen}
\end{figure}

\begin{figure}[H]
\centering

\begin{subfigure}[t]{0.4\linewidth}
\centering
\begin{tikzpicture}[line cap=round,line join=round,>=triangle 45,x=0.4cm,y=0.3cm]
\clip(-5,-8) rectangle (5,8);
\draw [line width=1pt] (-4,6)-- (4,6);
\draw [line width=1pt] (-4,-6)-- (4,-6);
\draw [line width=1pt] (-2,7)-- (-2,-4);
\draw [line width=1pt] (2,-7)-- (2,4);
\draw [line width=1pt] (-4,2)-- (4,2);
\draw [line width=1pt] (-4,-2)-- (4,-2);
\begin{scriptsize}
\draw [color=black] (-3,6.75) node {$5$};
\draw [color=black] (-3,-6.75) node {$5$};
\draw [color=black] (0,2.5) node {$2$};
\draw [color=black] (0,-2.5) node {$2$};
\draw [color=black] (2.5,0) node {$2$};
\draw [color=black] (-2.5,0) node {$2$};
\end{scriptsize}
\end{tikzpicture}
\end{subfigure}
%
\hspace{0.02\linewidth}
\raisebox{5.5em}{$\xleftarrow{\text{blow-ups}}$}
\hspace{0.02\linewidth}
%
\begin{subfigure}[t]{0.4\linewidth}
\centering
\begin{tikzpicture}[line cap=round,line join=round,>=triangle 45,x=0.4cm,y=0.3cm]
\clip(-5,-8) rectangle (5,8);
\draw [line width=1pt] (-4,6)-- (4,6);
\draw [line width=1pt] (-1,7)-- (-4,5);
\draw [line width=1pt] (-4,5.5)-- (-1,3);
\draw [line width=1pt,dashed] (-1,4)-- (-4,2);
\draw [line width=1pt] (-2.5,3.5)-- (-2.5,-2.5);
\draw [line width=1pt] (-3,2)-- (3,2);
\clip(-5,-8) rectangle (5,8);
\draw [line width=1pt] (4,-6)-- (-4,-6);
\draw [line width=1pt] (1,-7)-- (4,-5);
\draw [line width=1pt] (4,-5.5)-- (1,-3);
\draw [line width=1pt,dashed] (1,-4)-- (4,-2);
\draw [line width=1pt] (2.5,-3.5)-- (2.5,2.5);
\draw [line width=1pt] (3,-2)-- (-3,-2);
\begin{scriptsize}
\draw [color=black] (2,6.5) node {$6$};
\draw [color=black] (-2,7) node {$2$};
\draw [color=black] (-2,4.5) node {$2$};
\draw [color=black] (-3,0) node {$5$};
\draw [color=black] (0,2.5) node {$2$};
\draw [color=black] (-4,3) node {$1$};
\draw [color=black] (-2,-6.5) node {$6$};
\draw [color=black] (2,-7) node {$2$};
\draw [color=black] (2,-4.5) node {$2$};
\draw [color=black] (3,0) node {$5$};
\draw [color=black] (0,-2.5) node {$2$};
\draw [color=black] (4,-3) node {$1$};
\end{scriptsize}
\end{tikzpicture}
\end{subfigure}
\caption{Horikawa surfaces of special Lee-Park type $\mathrm{I}$ with $p_g=5$}
\label{fig:LP_special}
\end{figure}
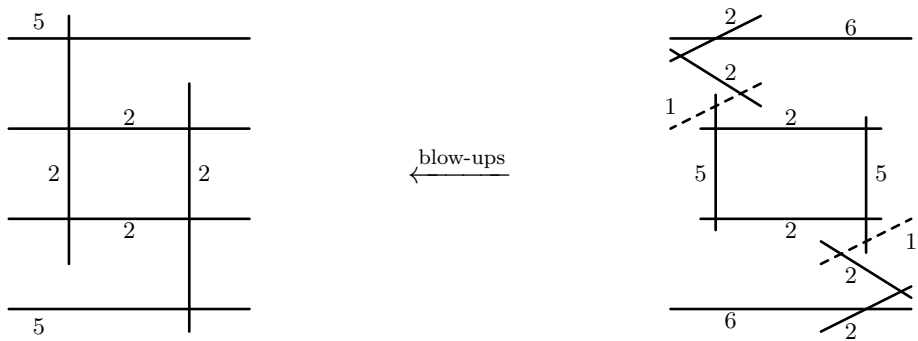
\end{exam}

In Section \ref{sec:deformation}, it turns out that the KSBA boundary is stratified as follows: 
\begin{itemize}
    \item  All Horikawa surfaces of special Lee-Park type with $p_g=5$ and an elliptic double cone singularity are partially $\Q$-Gorenstein smoothable to Horikawa surfaces of general Lee-Park type $\mathrm{I}^*$, that is, there is a $\Q$-Gorenstein deformation of  such a Horikawa surface of special Lee-Park type whose general fibers are Horikawa surfaces of general Lee-Park type $\mathrm{I}^*$.
    \item  All Horikawa surfaces of general Lee-Park type $\mathrm{I}^*$ and all Horikawa surfaces of special Lee-Park type are partially $\mathbb{Q}$-Gorenstein smoothable to Horikawa surfaces of general Lee-Park type I.
    \item If $p_g\geq6$, neither a Horikawa surface of general Lee-Park type $\mathrm{I}^*$ nor a Horikawa surface of special Lee-Park type degenerates to the other. 
\end{itemize}
These phenomena arise from the difference of the positivity of an anticanonical divisor of the quotient $W$ of $X$.
In fact, $-K_W$ is ample if $X$ is of type I, while $-K_W$ is not nef if $X$ is of type $\mathrm{I}^*$.
Since the ampleness of an anticanonical divisor is an open condition, this observation leads to the subtlety in the description of the KSBA boundary as explained above.

Next, we turn to the deformation theory of standard Horikawa surfaces.
The most technically challenging part of the stratification problem is 
to understand the deformations of standard Horikawa surfaces of type $(p_g-2)'$ for $p_g\ge 7$, which have elliptic double cone singularities.
We further divide the class of standard Horikawa surfaces of type $(p_g-2)'$ with $p_g\ge 7$ into three subclasses, according to the geometry of the branch divisor:
{\em of double Fano type}, {\em of double non-Fano type}, and {\em supersingular}.
This subdivision is essential for describing the stratification of the moduli space. 
See Definition \ref{defn--of--fano} for the precise definitions. 

To describe examples of the three classes, we introduce some notations.
Let $\eta\colon W^{+}=\Sigma_{p_g-2}\to W=\overline{\Sigma}_{p_g-2}$ be the minimal resolution of the quotient $W=X/\sigma$, where $\sigma$ denotes a good involution.
Let $B^{+}$ be the divisor on $W^{+}$ defined by the log crepant condition: 
$$
K_{W^{+}}+\frac{1}{2}B^{+}=\eta^{*}\left(K_W+\frac{1}{2}B\right),
$$
where $B$ is the branch divisor of the double cover $X\to W$.
According to the classification of elliptic double cone singularities (Proposition~\ref{prop:classification_cone_sing}),
the divisor $B^{+}$ decomposes as $B^{+}=\widehat{B}+2\Delta_0$, where $\Delta_{0}$ is the exceptional curve of $\eta$, and $\widehat{B}$ is the proper transform of $B$ such that any local intersection number of $\widehat{B}$ and $\Delta_0$ is at most $2$.

\begin{exam}[Standard Horikawa surfaces of double Fano type]
Consider the case where the divisor $\widehat{B}$ intersects the ruling fibers transversally along $\Delta_0$, as illustrated in Figure~\ref{fig:fano_non-fano_super-sing}~(1).
Such a surface $X$ is called a {\em standard Horikawa surface of double Fano type}.
In the process of resolving the singularities of $\widehat{B}+\Gamma$,
the proper transform of the fiber $\Gamma$ is blown up exactly once.
This property characterizes these surfaces by the existence of a certain log crepant birational model (i.e., an anti-P-resolution) that satisfies the weak Fano condition.
In fact, every smoothable Horikawa surface of special Lee-Park type with $p_g\ge 7$ deforms to a standard Horikawa surface of double Fano type.
\end{exam}

\begin{exam}[Standard Horikawa surfaces of double non-Fano type]
Consider the case where $\widehat{B}$ becomes tangent to some ruling fiber $\Gamma$ along $\Delta_0$, as illustrated in Figure~\ref{fig:fano_non-fano_super-sing}~(2).
We call such an $X$ a {\em standard Horikawa surface of double non-Fano type}.
In this case, to resolve the singularities of $\widehat{B}+\Gamma$, the proper transform of the fiber $\Gamma$ must be blown up twice.
This behavior implies the nonexistence of any weak Fano anti-P-resolution that is an isomorphism outside the fiber $\Gamma$.
Consequently, standard Horikawa surfaces of double non-Fano type do not degenerate to any Horikawa surface of special Lee-Park type for $p_g\ge7$.
Nevertheless, we will show that such surfaces can deform to standard Horikawa surfaces of double Fano type.
\end{exam}

\begin{exam}[Supersingular standard Horikawa surfaces]
Consider the case where the divisor $\widehat{B}$ contains a fiber $\Gamma$, and $\widehat{B}-\Gamma$ meets $\Gamma$ with multiplicity at least $3$, as illustrated in Figure~\ref{fig:fano_non-fano_super-sing}~(3).
To resolve the singularities of $\widehat{B}$, the proper transform of the fiber $\Gamma$ must be blown up at least three times.

For $p_g\neq 10$, supersingular standard Horikawa surfaces behave similarly to those of double non-Fano type.
In fact, they do not degenerate to any Horikawa surface of special Lee-Park type, but deform to standard Horikawa surfaces of double Fano type.

However, the case $p_g=10$ is exceptional.
According to Horikawa \cite{horikawa}, the moduli space of Horikawa surfaces with only Du Val singularities with $p_g=10$ has two connected components:
one generically parameterizing surfaces of type $(0)$, and the other of type $(6)$.
We will show that supersingular standard Horikawa surfaces deform to those of type $(6)$, while smoothable non-supersingular standard Horikawa surfaces deform only to type $(0)$, and not to type $(6)$.
This phenomenon reflects the subtlety of the deformation theory.
For instance, when a log crepant anti-P-resolution $(W^-,\frac{1}{2}B^-)\to (W,\frac{1}{2}B)$ is non-isomorphic along the fiber $\Gamma$, the cohomology group $H^1(\mathcal{O}_{W^-}(B^-))$ may not vanish, which obstructs log $\Q$-Gorenstein smoothings.
As discussed in Section~\ref{intro-subsec:smoothing_and_anti_P}, to address this difficulty, we study log $\Q$-Gorenstein smoothings for such pairs $(W^-,\frac{1}{2}B^-)$.
Although there exists a supersingular Horikawa surface that is smoothable to those of type $(0)$ (Example~\ref{exam:connectingGieseker}, Corollary~\ref{cor--p_g=10--connected}), it remains unknown whether all supersingular standard Horikawa surfaces with $p_g=10$ are smoothable to type $(0)$. 
This example plays a crucial role in the proof of Theorem \ref{thm--connectedness} as we will see later (Example~\ref{exam:connectingGieseker}). 
\end{exam}

\begin{figure}[H]
\centering
\begin{subfigure}[t]{0.32\linewidth}
\centering
\begin{tikzpicture}[line cap=round,line join=round,>=triangle 45,x=0.35cm,y=0.25cm]
\clip(-8,-10) rectangle (8,10);
\draw [line width=1pt] (-7,8)-- (-7,-8);
\draw [line width=1pt] (-7,8)-- (7,8);
\draw [line width=1pt] (7,8)-- (7,-8);
\draw [line width=1pt] (-7,-8)-- (7,-8);
\draw [line width=1pt] (-6,4)-- (6,4);
\draw [line width=1pt] (4,6)-- (4,-6);
\draw [line width=1pt,red] (3,-2)-- (5,-2);
\draw [line width=1pt,red] (3,-4)-- (5,-4);
\draw[line width=1pt,color=red,smooth,samples=100,domain=3:5,rotate around={45:(4,4)}] plot(\x,{(\x - 4)^4+4});
\draw[line width=1pt,color=red,smooth,samples=100,domain=3:5,rotate around={45:(4,4)}] plot(\x,{0-(\x - 4)^4+4});
\draw[line width=1pt,color=red,smooth,samples=100,domain=-5:-3,rotate around={45:(-4,4)}] plot(\x,{(\x + 4)^4+4});
\draw[line width=1pt,color=red,smooth,samples=100,domain=-5:-3,rotate around={45:(-4,4)}] plot(\x,{0-(\x + 4)^4+4});
\begin{scriptsize}
\draw [fill] (4.5,0) node {$\Gamma$};
\draw [fill] (0,5) node {$\Delta_0$};
\end{scriptsize}
\end{tikzpicture}
\caption{of double Fano type}
\label{fig:intro_Fano_type}
\end{subfigure}
\hfill
\begin{subfigure}[t]{0.32\linewidth}
\centering
\begin{tikzpicture}[line cap=round,line join=round,>=triangle 45,x=0.35cm,y=0.25cm]
\clip(-8,-10) rectangle (8,10);
\draw [line width=1pt] (-7,8)-- (-7,-8);
\draw [line width=1pt] (-7,8)-- (7,8);
\draw [line width=1pt] (7,8)-- (7,-8);
\draw [line width=1pt] (-7,-8)-- (7,-8);
\draw [line width=1pt] (-6,4)-- (6,4);
\draw [line width=1pt] (4,6)-- (4,-6);
\draw[line width=1pt,color=red,smooth,samples=100,domain=3:5,rotate around={90:(4,4)}] plot(\x,{(\x - 4)^2+4});
\draw[line width=1pt,color=red,smooth,samples=100,domain=3:5,rotate around={90:(4,4)}] plot(\x,{0-(\x - 4)^2+4});
\draw[line width=1pt,color=red,smooth,samples=100,domain=-5:-3,rotate around={45:(-4,4)}] plot(\x,{(\x + 4)^4+4});
\draw[line width=1pt,color=red,smooth,samples=100,domain=-5:-3,rotate around={45:(-4,4)}] plot(\x,{0-(\x + 4)^4+4});
\begin{scriptsize}
\draw [fill] (4.5,0) node {$\Gamma$};
\draw [fill] (0,5) node {$\Delta_0$};
\end{scriptsize}
\end{tikzpicture}
\caption{of double non-Fano type}
\label{fig:intro_non-Fano_type}
\end{subfigure}
\hfill
\begin{subfigure}[t]{0.32\linewidth}
\centering
\begin{tikzpicture}[line cap=round,line join=round,>=triangle 45,x=0.35cm,y=0.25cm]
\clip(-8,-10) rectangle (8,10);
\draw [line width=1pt] (-7,8)-- (-7,-8);
\draw [line width=1pt] (-7,8)-- (7,8);
\draw [line width=1pt] (7,8)-- (7,-8);
\draw [line width=1pt] (-7,-8)-- (7,-8);
\draw [line width=1pt] (-6,4)-- (6,4);
\draw [line width=1pt,red] (4,6)-- (4,-6);
\draw[line width=1pt,color=red,smooth,samples=100,domain=-5:-3,rotate around={45:(-4,4)}] plot(\x,{(\x + 4)^4+4});
\draw[line width=1pt,color=red,smooth,samples=100,domain=-5:-3,rotate around={45:(-4,4)}] plot(\x,{0-(\x + 4)^4+4});
\draw[line width=1pt,color=red,smooth,samples=100,domain=3:5,rotate around={90:(4,4)}] plot(\x,{0-(\x - 4)^4+4});
\begin{scriptsize}
\draw [fill] (4.5,0) node {$\Gamma$};
\draw [fill] (0,5) node {$\Delta_0$};
\end{scriptsize}
\end{tikzpicture}
\caption{supersingular}
\label{fig:intro_super_singular}
\end{subfigure}
\caption{Standard Horikawa surfaces of type $(p_g-2)'$}
\label{fig:fano_non-fano_super-sing}
\end{figure}

Taking into account the subdivisions introduced above, we now describe the stratifications of the boundary of the moduli space and an answer to Question~\ref{question}~(3).
For each integer $p_g\ge 3$, let $\overline{M}^{\mathrm{Gie}}_{2p_g-4,p_g}$ be the coarse moduli space of $\Q$-Gorenstein smoothable normal stable Horikawa surfaces of geometric genus $p_g$.
Let $(\overline{M}^{\mathrm{Gie}}_{2p_g-4,p_g})^\mathrm{sn}$ denote its seminormalization.
Recall that taking the seminormalization does not alter the topological type (see Definition \ref{defn--seminormal}).
The following theorem answers Question~\ref{question}~(3).

\begin{thm}[{cf.~Theorem \ref{thm--stratification}}]\label{thm--stratification--intro}
For any $p_g\neq 10$, 
    The stratification of the boundary of $(\overline{M}^{\mathrm{Gie}}_{2p_g-4,p_g})^\mathrm{sn}$ is described in Figures~\ref{fig:hasse-p_g=3}, \ref{fig:hasse-p_g=4}, \ref{fig:hasse-p_g=5}, \ref{fig:hasse-p_g=6}, \ref{fig:hasse-p_gge7}, and \ref{fig:hasse-p_gge10mod4} (see Notation \ref{note--stratification}). 
\end{thm}

In the case $p_g=10$, we also give a complete description of the moduli stratification, excluding the supersingular standard Horikawa surfaces, which require separate treatment; see Theorem~\ref{thm--stratification}.

\begin{rem}[Recent research]
We make a few remarks on previous work concerning stratifications of moduli spaces.
    \begin{itemize}
        \item[$(1)$]
        Anthes \cite{Anthes} described a stratification of the moduli space parametrizing (possibly non-normal) Gorenstein stable Horikawa surfaces with $p_g=3$ in terms of elliptic singularities and conductors.
        \item[$(2)$]
        Rana-Rollenske \cite{RR} described a stratification of the moduli space parameterizing (possibly non-normal) Gorenstein stable Horikawa surfaces that arise as double covers over Hirzebruch surfaces.
        Note that the intersection of this moduli space and the moduli space $(\overline{M}^{\mathrm{Gie}}_{2p_g-4,p_g})^\mathrm{sn}$ parametrizes standard Horikawa surfaces of type $(d)$ for some $d$.
        \item[$(3)$]
        After completing our proof, we noticed that Evans-Simonetti-Urz\'ua revised their paper \cite{ESU} and dealt with the stratification problem in the case where $p_g=3$.
    \end{itemize}
\end{rem}

As explained earlier, Theorem \ref{thm--main--ii} was proven by studying the log $\Q$-Gorenstein deformation of the quotient pair $(W,\frac{1}{2}B)$. 
The same strategy applies to Theorem \ref{thm--stratification--intro}.
Let $(M^{\mathrm{nq}}_{p_g})^{\mathrm{sn}}$ denote the seminomalization of the coarse moduli space of such quotient pairs, obtained from a Horikawa surface $X$ with a good involution.
The next result shows that understanding the stratification of the boundary of $(\overline{M}^{\mathrm{Gie}}_{2p_g-4,p_g})^\mathrm{sn}$
reduces to studying the corresponding stratification for $(M^{\mathrm{nq}}_{p_g})^{\mathrm{sn}}$.

\begin{thm}[{$=$Theorem \ref{prop--moduli--seminormaliazation}}]\label{thm--moduli--seminormaliazation--intro}
    There exists a finite morphism $\overline{\alpha}'\colon(\overline{M}^{\mathrm{Gie}}_{2p_g-4,p_g})^{\mathrm{sn}}\to (M^{\mathrm{nq}}_{p_g})^{\mathrm{sn}}$ which restricts to an isomorphism over the open locus parameterizing all Horikawa surfaces except those described in Theorem~\ref{intro-qGsmHor}~(2). 
\end{thm}

When the Gieseker moduli space has two connected components, equivalently in the case where the canonical volume is divisible by eight, it is natural to ask whether these components are joined in the KSBA compactification.
Rana-Rollenske \cite{RR} proved that the complete KSBA moduli space of $\Q$-Gorenstein smoothable Horikawa surfaces is connected by allowing degenerations to non-normal stable surfaces.
On the other hand, Monreal-Negrete-Urz\'ua \cite{MNU} showed that they are not connected via degenerations to $\Q$-Gorenstein smoothable klt stable Horikawa surfaces.
Inspired by their results, we investigate a similar phenomenon under the restriction that only $\Q$-Gorenstein smoothable lc stable Horikawa surfaces are allowed, and obtain the following result.

\begin{thm}\label{thm--connectedness}
    $\overline{M}^{\mathrm{Gie}}_{2p_g-4,p_g}$ is connected if and only if $p_g=6,10$ or $p_g-2\not\in 4\mathbb{Z}$. 
\end{thm}

To prove Theorem \ref{thm--connectedness} in the case $p_g=10$, we construct supersingular standard Horikawa surfaces that smooth to both type $(0)$ and type $(6)$ surfaces (Example~\ref{exam:connectingGieseker}).
Note that any normal stable Horikawa surface other than supersingular standard Horikawa surfaces does not lie in the closure of the component parameterizing type $(6)$ surfaces.

\subsection{Classification of normal stable Horikawa surfaces}\label{intro-subsec:classification}
We give an overview of the classification method of normal stable Horikawa surfaces with only $\Q$-Gorenstein smoothable singularities and those admitting the good involution.
These results provide (partial) answers to Question \ref{question} $(1)$ and $(2)$.
It is worth noting that the list of normal stable Horikawa surfaces with only $\Q$-Gorenstein smoothable singularities (without imposing the existence of good involutions) is of independent interest.
Indeed, this list is expected to serve as a key guide in the investigation of components other than the Gieseker component within the KSBA moduli space.

The key tools developed in this paper are the \emph{log Noether inequality} (Theorem \ref{normalstable}) and an analysis of \emph{extended T-chains} (Appendix \ref{app:drill}).
Our classification strategy is as follows: First, we establish the log Noether inequality and structure theorems for standard and non-standard Horikawa surfaces.
Then, by analyzing extended T-chains combinatorially, we obtain our classification results (Theorems \ref{thm:standard_Horikawa}, \ref{thm:classification_non-standard_Horikawa_noninv}, \ref{thm:classification_non-standard_Horikawa}).
In what follows, we explain the key tools in turn.

\subsubsection{Log Noether inequality}\label{intro-subsubsec:logNoether}
Let $X$ be a normal stable Horikawa surface.
We consider the following diagram:
\[
\xymatrix{
 & \ar[ld]_{\rho} \widetilde{X} \ar[rd]^{\pi}&  \\
Y  &    & X \\
}
\]
where $\pi\colon\widetilde{X}\to X$ denotes the minimal resolution, and $Y$ denotes a minimal model of $\widetilde{X}$.
If $X$ has only rational singularities, then the geometric genus of $\widetilde{X}$ equals that of $X$, but the canonical volume satisfies $K_Y^{2}\le K_{X}^{2}$.
If $X$ has non-Du Val singularities, the inequality is strict and $Y$ does not satisfy the Noether inequality, implying that $Y$ has Kodaira dimension $1$.
In particular, $Y$ admits an elliptic fibration, which provides information of the structure of $X$.

However, if $X$ has elliptic singularities, the above approach does not apply directly.
Indeed, the geometric genus of $\widetilde{X}$ may be less than that of $X$. 

To overcome this, we establish a {\em logarithmic version of Noether's inequality} (Theorem~\ref{normalstable}, see also \cite{Che}), which is of independent interest.
Instead of the diagram above, we consider the following diagram:
\[
\xymatrix{
 & \ar[ld]_{\rho} (\widetilde{X},\widetilde{\Delta}) \ar[rd]^{\pi}&  \\
(Y,\Delta)  &    & X \\
}
\]
where $\widetilde{\Delta}$ is the reduced exceptional divisor over all elliptic singularities of $X$, and $(\widetilde{X}, \widetilde{\Delta})\to (Y, \Delta)$ denotes a log minimal model.
Theorem \ref{normalstable} (2) (applied to standard Horikawa surfaces) provides a structure theorem for standard Horikawa surfaces (Proposition \ref{prop:standard_Horikawa}), whereas Theorem \ref{normalstable} (1) together with \eqref{eqn:ellfundeq}, the formula relating invariants of $X$ and singularities on $X$, provides a structure theorem for non-standard Horikawa surfaces (Proposition~\ref{prop:non-std_Horikawa}~$(2)$).

Using the structure theorem for standard Horikawa surfaces, we obtain their classification (Theorem \ref{thm:standard_Horikawa}).
These correspond to the Horikawa surfaces appearing in Theorem \ref{intro-qGsmHor} (1).

On the other hand, the structure theorem for non-standard Horikawa surfaces indicates that the surface $Y$ must be an elliptic surface over either a smooth rational curve or an elliptic curve.
In order to determine the singularities on the non-standard Horikawa surface $X$, we examine the divisor $C\subset Y$, which is the pushforward along $\widetilde{X}\to Y$ of the reduced inverse image $\widetilde{C}\subset\widetilde{X}$ of the non-Du Val singularities on $X$.
Since $C$ may involve \emph{extended T-chains}, which are introduced below, we study extended T-chains to examine $C$.

\begin{rem}
    Variants of log Noether inequalities have been established in \cite{TZ,LR,Liu,Che}, all of which are not equivalent to Theorem~\ref{normalstable}.
    While the result in \cite{Che} is very close to Theorem~\ref{normalstable},
    the formulation of Theorem~\ref{normalstable} is more applicable to our situation. 
\end{rem}

\subsubsection{Extended T-chains}\label{intro-subsubsec:extT}
Throughout this paper, the terms “string” and “chain” are used interchangeably.

An \emph{extended T-chain} is a generalization of a \emph{T-chain}.
We begin by recalling the notion of a T-chain.
A T-chain is the string associated with a non-Du Val T-singularity via the Hirzebruch-Jung continued fraction.
According to \cite[Proposition 3.11]{KSB}, all T-chains are obtained from $[4]$ or $[3,2,\dots,2,3]$ with $2$ repeated $n$ times ($n\geq0$) by repeatedly applying the following two operations to strings:
\[
[b_1,\dots,b_r]\mapsto[2,b_1,\dots,b_r+1], \quad [b_1,\dots,b_r]\mapsto[b_1+1,\dots,b_r,2].
\]

We are now ready to introduce an extended T-chain.
We call a string $B=[b_1,\dots,b_r]$ an \emph{extended T-chain} if $B$ corresponds to the chain obtained from a chain $\widetilde{B}$ consisting of T-chains and $(-1)$-curves by contracting all $(-1)$-curves; see Definition \ref{def:extT} for the precise definition.
The chain $\widetilde{B}$ is called a \emph{T-train} associated with $B$.
For example, the string $[3,2,2,3,4,3,2]$ is an extended T-chain since it can be realized as the blow-down of the T-train $[3,3,1,3,5,2,1,3,5,3,2]$.
We write the T-train as 
\begin{equation}\label{intro-example_of_T-train}
    [3,3]-1-[3,5,2]-1-[3,5,3,2].
\end{equation}
This example appears in Theorem \ref{thm:one-section-T}, a classification of configurations on certain normal stable surfaces.

We provide a geometric interpretation of extended T-chains and associated T-trains.
Let $B$ be an extended T-chain, let $\widetilde{B}$ be a T-train associated with $B$, and let $p\in U$ be a cyclic quotient singularity whose exceptional set corresponds to $B$.
Then, there exists the resolution $U'\to (p\in U)$ whose exceptional set corresponds to $\widetilde{B}$, and the partial resolution $\widetilde{U}\to U$ such that $\widetilde{U}$ is obtained from $U'$ by contracting all T-chains in the T-train $\widetilde{B}$.
We say that a T-train $\widetilde{B}$ is \emph{ample} if $K_{\widetilde{U}}$ is ample over $U$, and that an extended T-chain is \emph{P-admissible} if there exists an associated ample T-train.
For example, the T-train \eqref{intro-example_of_T-train} is ample, and $[3,2,2,3,4,3,2]$ is P-admissible; this can be checked by a combinatorial computation (Lemma \ref{lem_ampleness}).
Note that, if $\widetilde{B}$ is ample, the partial resolution $\widetilde{U}\to U$ is nothing but a P-resolution \cite{KSB}.
See Proposition \ref{charP-resol} for futher details.

During the classification process, it becomes necessary to verify, in many instances, whether a given string is a P-admissible chain.
The necessary results are collected in Appendix \ref{app:drill}.

\begin{rem}
    The notion of extended T-chains and ample T-trains have already been studied in several works.
    We briefly review these here.
    \begin{itemize}
        \item[$(1)$]
        Figueroa-Rana-Urz\'ua provided a necessary combinatorial condition for a T-train to be nef \cite[Corollary 2.4]{FRU}.
        \item[$(2)$]
        Monreal-Negrete-Urz\'ua \cite{MNU} gave a classification of certain extended T-chains that is very similar to ours.
        In fact, if we restrict our classification to extended T-chains that can arise from Horikawa surfaces with only T-singularities, these results coincide.
    \end{itemize}
\end{rem}

\subsubsection{Completion of the classifications}
As we already explained the classification of standard Horikawa surfaces in Section \ref{intro-subsubsec:logNoether}, we focus on the classification of non-standard Horikawa surfaces here.

According to the structure theorem (Proposition \ref{prop:non-std_Horikawa}), a non-standard Horikawa surface $X$ with only $\Q$-Gorenstein smoothable singularities falls into the following three classes:
\begin{itemize}
    \item[(A)] 
    $Y$ is a relatively minimal elliptic fibration over $\PP^1$, and the divisor $C\subset Y$ contains exactly one section of $Y\to\PP^1$.
    \item[(B)] 
    $Y$ is a relatively minimal elliptic fibration over $\PP^1$, and the divisor $C\subset Y$ contains exactly two sections of $Y\to\PP^1$.
    \item[(C)] 
    $Y$ is a relatively minimal elliptic fibration over an elliptic curve, and the divisor $C\subset Y$ contains exactly one section.
\end{itemize}
We point out that the case~(C) never occurs when $X$ is klt.

If $X$ has a T-singularity, an extended T-chain naturally appears in the image of $\widetilde{C}\subset\widetilde{X}$ during the process of contracting $(-1)$-curves in $\widetilde{X}$.
In each case described above, we analyze the configuration of $C$ and $\widetilde{C}$ using a combinatorial analysis of extended T-chains and the list of $\Q$-Gorenstein smoothable singularities (Lemma \ref{lem:smoothable_elliptic_singularity}, \ref{lem:smoothable_rational_strict_lc}, and Appendix \ref{app:cusp}).
As a result, we identify all possible $C\subset Y$ and $\widetilde{C}\subset\widetilde{X}$, and obtain the classification (Theorem \ref{thm:classification_non-standard_Horikawa_noninv}).
This provides a list of all possible non-standard Horikawa surfaces with only $\Q$-Gorenstein smoothable singularities; that is, it gives an answer to Question \ref{question}~(1).

If we further assume that $X$ admits a good involution, we obtain the list of surfaces among the classification of non-standard Horikawa surfaces admitting a good involution (Theorem \ref{thm:classification_non-standard_Horikawa}).
The main tools are the classification of non-standard Horikawa surfaces (Theorem \ref{thm:classification_non-standard_Horikawa_noninv}) and the list of singular fibers of elliptic surfaces with an involution (Appendix \ref{app:sing_fiber_inv}).
This proves Theorem \ref{intro-qGsmHor}.

We conclude this section by explaining which case in the above three classes each surface in Theorem \ref{intro-qGsmHor} corresponds to.
The surfaces in Theorem \ref{intro-qGsmHor} $(4)$, $(5)$, $(6)$, and $(7)$ correspond to the first case; those in Theorem \ref{intro-qGsmHor} $(3)$, to the second case; and those in Theorem \ref{intro-qGsmHor} $(2)$, to the third case.

\begin{rem}\label{rem:broader}
    Although our purpose is to classify normal stable Horikawa surfaces, we also obtain a classification of certain normal stable surfaces by the same procedure.
    See Section \ref{subsec:overview_non-std} for a detailed explanation. 
    If we restrict our classification to surfaces with only klt singularities, we recover the classification of \emph{small surfaces} \cite{MNU}.  
    Note that the class of small surfaces extends the class of klt non-standard Horikawa surfaces in case (A); see Remark \ref{rem:small_surfaces} for details.
\end{rem}

\begin{rem}
Here, we briefly discuss how our arguments in Section~\ref{sec:non-std_Horikawa} differ from those by Monreal-Negrete-Ur\'zua \cite{MNU}.
\begin{itemize}
    \item[(1)]
    In \cite[Corollary 2.12]{MNU}, they provide a structure theorem for non-standard Horikawa surfaces with only T-singularities by using \cite[Proposition 3.1]{RU2}, a classification of certain klt projective surfaces.
    In particular, they showed that non-standard Horikawa surfaces with only T-singularities fall into the cases (A) and (B).
    This is the specialization of Proposition \ref{prop:non-std_Horikawa} to the klt setting.
    \item[(2)] 
    The case (A) in the klt setting is studied in \cite[Sections 3, 4]{MNU}.
    They classified small surfaces \cite[Theorem 3.12]{MNU}, and showed that their constructions fall into $14$ building blocks.
    In contrast, we classify certain surfaces similar to the surfaces in the case (A), as explained in Remark~\ref{rem:broader}, and obtain $27$ building blocks (Theorems \ref{thm:one-section-T}, \ref{thm:one-section-halfcusp}, and \ref{thm:one-section-triangle}).
    \item[(3)] 
    In the klt setting, the classification of the case (B) is straightforward, as shown in \cite[Corollary~2.12]{MNU}.
    In fact, it is shown that $\widetilde{X} = Y$ if $X$ is klt.
    However, if $X$ has at least one elliptic singularity, the morphism $\widetilde{X} \to Y$ is no longer an isomorphism and a further analysis is needed; see Section \ref{subsec:0,2,0}.
    \item[(4)] 
    In \cite[Theorem~5.12]{MNU}, they provide candidates for $\Q$-Gorenstein smoothable klt non-standard Horikawa surfaces under the assumption that $p_g \ge 10$. 
    Their approach relies on the observation that, if the surface is $\Q$-Gorenstein smoothable, then it must admit an involution whose quotient is a degeneration of a Hirzebruch surface.  
    To prove \cite[Theorem~5.12]{MNU}, they make use of \cite[Theorem~5.6]{MNU}, which states that degenerations of Hirzebruch surfaces with only T-singularities have at most four singular points.
    
    In the present paper, we focus more on involutions.
    We classify the possible actions of involutions on singularities of non-standard Horikawa surfaces, and determine singularities that can occur on quotient surfaces.
    To achieve this, we utilize the classification of singular fibers of elliptic surfaces with involution (see Appendix~\ref{app:sing_fiber_inv}).  
    As a result, we are able to narrow down the possible candidates for $\Q$-Gorenstein smoothable non-standard Horikawa surfaces without assuming any condition on the geometric genus (Theorem~\ref{thm:classification_non-standard_Horikawa}).
\end{itemize}
\end{rem}

\subsection{$\Q$-Gorenstein smoothing and anti-P-resolution}
\label{intro-subsec:smoothing_and_anti_P}

To address Question~\ref{question}~(2), we investigate the $\Q$-Gorenstein smoothing problem for the Horikawa surfaces appearing in Theorem~\ref{intro-qGsmHor}.
Note that any $\Q$-Gorenstein smoothable normal stable Horikawa surface $X$ must carry a (unique) good involution as the limit of the good involutions of smooth Horikawa surfaces.
Thus, the existence of a good involution is a necessary condition for $\Q$-Gorenstein smoothability.
Given a Horikawa surface $X$ with a good involution $\sigma$, we can express $X$ as a double cover of the quotient $W:=X/\sigma$, branched along a divisor $B\in |2L|$.
Then, $X$ is $\Q$-Gorenstein smoothable if and only if the pair $(W, \frac{1}{2}B)$ is log $\Q$-Gorenstein smoothable.
This reduces the smoothing problem for $X$ to the smoothing problem for the pair $(W, \frac{1}{2}B)$. 

In most cases, $-K_W$ is big, so the local-to-global obstruction for deforming $W$ vanishes (cf.~\cite{HP}).
Consequently, $W$ can be smoothed. 
Since $W$ has only smoothable singularities by the classification, we can construct a one-parameter smoothing $\mathscr{W}\to C$ of $W$.
Then, the problem of the $\Q$-Gorenstein smoothability for $X$ reduces to the problem of the liftability of $B\in |2L|$ on $\mathscr{W}$. 
When $\mathscr{W}$ is a $\Q$-Gorenstein smoothing of T-singularities, we provide a criterion for deforming the divisorial sheaf $L$ (Proposition~\ref{prop--smoothability--criterion}).
In this case, the extendability of the branch divisor $B$ to the total space of the deformation can be checked by the vanishing of the cohomology $H^{1}(\O_{W}(B))$.

However, we must also address the cases in which $\mathscr{W}$ is not $\Q$-Gorenstein.
This case occurs when $X$ has a double cone singularity, and so $\mathscr{W}$ is a smoothing of the cone singularity of $W=\mathscr{W}_0$.
Log $\Q$-Gorenstein deformation theory is more complicated and the vanishing $H^{1}(\O_{W}(B))=0$ is not sufficient to deduce the extendability of $B$ (cf.~Remark \ref{rem--necessity-of-P-or-anti-p}).
Fortunately, since $W$ is klt in all cases of Theorem \ref{intro-qGsmHor}, we obtain a small birational morphism $\mathscr{W}^{+}\to \mathscr{W}$ such that the central fiber $\mathscr{W}^{+}_{0}$ has only T-singularities (cf.~\cite{fujino--some--remarks}).
The induced map on central fibers $\mathscr{W}^{+}_{0}\to \mathscr{W}_{0}$ is a P-resolution in the sense of \cite{KSB}.
Thus, one can attempt to first take a P-resolution $W^{+}\to W$, and then study the log $\Q$-Gorenstein smoothability of the pair $(W^{+}, \frac{1}{2}B^{+})$, which is log crepant with respect to $(W, \frac{1}{2}B)$.
However, handling the deformation of the branch divisor $B^{+}$ can be subtle.
Indeed, $H^{1}(\O_{W^{+}}(B^{+}))$ may fail to vanish, making the deformation theory more difficult. 

To deal with this issue, we introduce the notion of an {\em anti-P-resolution}.
For a klt surface $W$, a projective birational morphism $\mu\colon W^{-}\to W$ is called an anti-P-resolution if $W^{-}$ has only $\Q$-Gorenstein smoothable slc singularities and $-K_{W^{-}}$ is ample over $W$.
In this paper, we consider anti-P-resolutions $\mu\colon (W^{-}, \frac{1}{2}B^{-})\to (W, \frac{1}{2}B)$ satisfying the log crepant condition:
$$
K_{W-}+\frac{1}{2}B^{-}=\mu^{*}\left(K_W+\frac{1}{2}B\right),
$$
where both $B$ and $B^{-}$ are effective divisors.
By this condition, $B^-$ is automatically $\mu$-ample and $H^1(\mathcal{O}_{W^-}(B^-))$ is more likely to vanish. This is an advantage of an anti-P-resolution.
Anti-P-resolutions naturally arise in the following context.
Let $\mathscr{X}\to C$ be a $\Q$-Gorenstein smoothing of Horikawa surfaces, and consider $\mathscr{W}=\mathscr{X}/\sigma$ and $\mathscr{W}^{+}\to \mathscr{W}$ as before, where $\sigma$ is the fiberwise good involution.
Since $\mathscr{W}$ is potentially klt, we can perform a $-K_{\mathscr{W}^{+}}$-flip $\mathscr{W}^{+}\dasharrow \mathscr{W}^{-}$. 
Since the morphism $\mathscr{W}^{-}\to \mathscr{W}$ is small, the branch divisor $\mathscr{B}$ on $\mathscr{W}$ lifts uniquely to an effective divisor $\mathscr{B}^{-}$ on $\mathscr{W}^{-}$.
The induced map on central fibers $(\mathscr{W}^{-}_{0},\frac{1}{2}\mathscr{B}^{-}_{0})\to (\mathscr{W}_{0},\frac{1}{2}\mathscr{B}_0)$ is an anti-P-resolution which is log crepant.

If we aim to prove the non-existence of a $\Q$-Gorenstein smoothing, we instead assume contrary that such a smoothing exists.
We then obtain an anti-P-resolution $(W^{-}, \frac{1}{2}B^{-})\to (W,\frac{1}{2}B)$ as described above.
Such an anti-P-resolution has some restrictions on its exceptional locus and the self-intersection number of the relative canonical class.
We refer to them as {\em admissible anti-P-resolutions} (Definition~\ref{def--adm--anti--P}).
If we can classify all admissible anti-P-resolutions, then in some cases, we may reach a contradiction.
In this paper, we classify all admissible anti-P-resolutions (Theorem~\ref{thm--anti-P-k_1-k_2}).
Note that admissible anti-P-resolutions are classified into three types:
those of normal type (Example~\ref{exam--anti-P-type-C}), those of non-normal type $\mathrm{I}$ (Example~\ref{exam--anti-P-type-B}), and those of non-normal type $\mathrm{II}$ (Example~\ref{exam--anti-P-type-A}).

On the other hand, we propose the following strategy for proving the existence of a $\Q$-Gorenstein smoothing of a given Horikawa surface $X$.
Start by taking an admissible anti-P-resolution $\mu\colon (W^{-},\frac{1}{2}B^{-})\to (W, \frac{1}{2}B)$ for the quotient $W=X/\sigma$ by the good involution and the associated branch divisor $B$.
Then, consider a $\Q$-Gorenstein smoothing $\mathscr{W}^{-} \to C$ and a lifting $\mathscr{B}^{-}$ of the branch divisor $B^{-}$.
Finally, take the ample model $(\mathscr{W}, \frac{1}{2}\mathscr{B})$ of $K_{\mathscr{W}^{-}/C}+\frac{1}{2}\mathscr{B}^{-}$ and the double cover $\mathscr{X}\to \mathscr{W}$ branched along $\mathscr{B}$.
The resulting family $\mathscr{X}\to C$ is in fact a $\Q$-Gorenstein smoothing of $X$.

Based on the above methods, we prove Theorem \ref{thm--main--ii} using the following strategy. 
We consider the case where $X$ has a double cone singularity that is locally equivariantly smoothable with respect to the good involution. 
This assumption ensures the existence of an admissible anti-P-resolution $\mu\colon (W^-,\frac{1}{2}B^-)\to (W,\frac{1}{2}B)$ such that $W^-$ is $\Q$-Gorenstein smoothable in a neighborhood of the exceptional curve.
Since the local-to-global obstraction vanishes (by Theorem~\ref{thm--hacking--prokhorov}), the surface $W^-$ admits a $\Q$-Gorenstein smoothing family 
$\mathscr{W}^-\to C$.
Thus, to verify the $\Q$-Gorenstein smoothability of $X$, it suffices to prove the vanishing $H^1(\mathcal{O}_{W^-}(B^-))=0$, because $\mathscr{W}^-$ is $\Q$-Gorenstein.

To establish this vanishing, we use vanishing theorems from \cite{fujino--slc--vanishing,Enokizono}.
We rely on Fujino's vanishing theorem for slc schemes \cite{fujino--slc--vanishing} since admissible anti-P-resolutions are often non-normal.
To apply this theorem, we make use of our explicit classification of admissible anti-P-resolutions and verify the positivity of $B^-$ case by case.
However, in certain situations---such as for supersingular standard Horikawa surfaces with $p_g\neq 10$---the divisor $B^--K_{W^-}$ fails to be big and nef, making Fujino's vanishing theorem inapplicable.
In these cases, we instead use the vanishing theorem of the second author \cite{Enokizono}, which requires only a weaker condition on the positivity of $B^--K_{W^-}$, known as $\mathbb{Z}$-positivity (see Definition \ref{defn-z-positive}). 

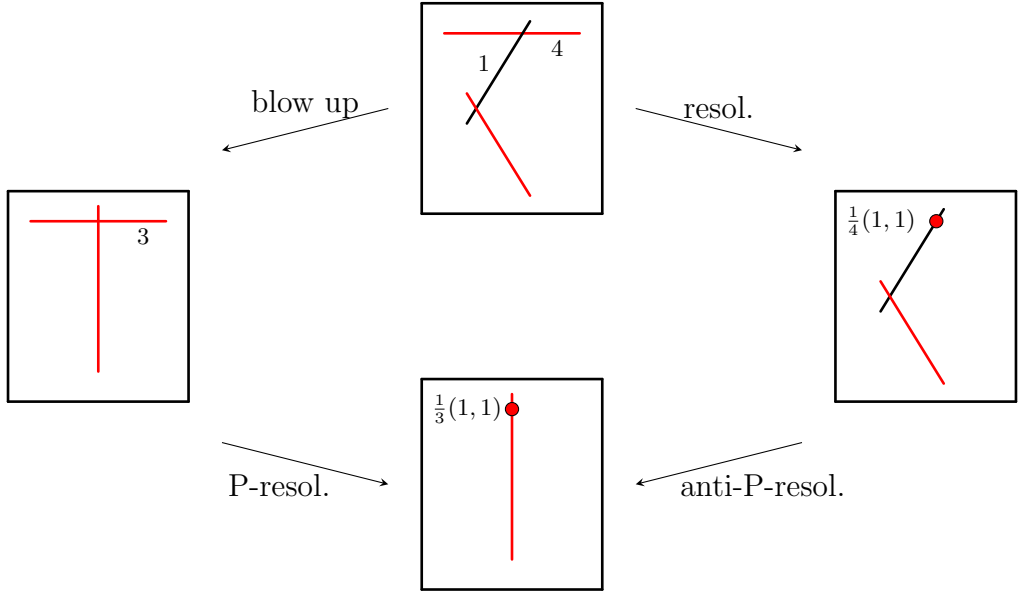
\begin{figure}[H]
\centering
\begin{tikzpicture}[>=stealth]


\node (A) at (0,2.5) {
  \begin{tikzpicture}[line cap=round,line join=round,>=triangle 45,x=0.3cm,y=0.2cm]
\clip(-5,-8) rectangle (5,10);
\draw [line width=1pt] (-4,8)-- (-4,-6);
\draw [line width=1pt] (-4,-6)-- (4,-6);
\draw [line width=1pt] (-4,8)-- (4,8);
\draw [line width=1pt] (4,8)-- (4,-6);
\draw [line width=1pt,red] (-3,6)-- (3,6);
\draw [line width=1pt] (0.8,6.8)-- (-2,0);
\draw [line width=1pt,red] (-2,2)-- (0.8,-4.8);
\begin{scriptsize}
\draw [fill] (2,5) node {$4$};
\draw [fill] (-1.25,4) node {$1$};
\end{scriptsize}
\end{tikzpicture}
};

\node (B) at (-5.5,0) {
  \begin{tikzpicture}[line cap=round,line join=round,>=triangle 45,x=0.3cm,y=0.2cm]
    \clip(-5,-8) rectangle (5,10);
\draw [line width=1pt] (-4,8)-- (-4,-6);
\draw [line width=1pt] (-4,-6)-- (4,-6);
\draw [line width=1pt] (-4,8)-- (4,8);
\draw [line width=1pt] (4,8)-- (4,-6);
    \draw [line width=1pt,red] (-3,6)-- (3,6);
    \draw [line width=1pt,red] (0,7)-- (0,-4);
    \begin{scriptsize}
      \draw [fill] (2,5) node {$3$};
    \end{scriptsize}
  \end{tikzpicture}
};

\node (C) at (5.5,0) {
 \begin{tikzpicture}[line cap=round,line join=round,>=triangle 45,x=0.3cm,y=0.2cm]
\clip(-5,-8) rectangle (5,10);
\draw [line width=1pt] (-4,8)-- (-4,-6);
\draw [line width=1pt] (-4,-6)-- (4,-6);
\draw [line width=1pt] (-4,8)-- (4,8);
\draw [line width=1pt] (4,8)-- (4,-6);
\draw [line width=1pt] (0.8,6.8)-- (-2,0);
\draw [line width=1pt,red] (-2,2)-- (0.8,-4.8);
\begin{scriptsize}
\draw [fill=red] (0.47,6) circle (2.5pt);
\draw [fill] (-2,6) node {$\frac{1}{4}(1,1)$};
\end{scriptsize}
\end{tikzpicture}
};

\node (D) at (0,-2.5) {
  \begin{tikzpicture}[line cap=round,line join=round,>=triangle 45,x=0.3cm,y=0.2cm]
    \clip(-5,-8) rectangle (5,10);
\draw [line width=1pt] (-4,8)-- (-4,-6);
\draw [line width=1pt] (-4,-6)-- (4,-6);
\draw [line width=1pt] (-4,8)-- (4,8);
\draw [line width=1pt] (4,8)-- (4,-6);
    \draw [line width=1pt,red] (0,7)-- (0,-4);
    \begin{scriptsize}
    \draw [fill=red] (0,6) circle (2.5pt);
      \draw [fill] (-2,6) node {$\frac{1}{3}(1,1)$};
    \end{scriptsize}
  \end{tikzpicture}
};

\draw[->]  (A.west) -- (B.north east) node[midway, above] {blow up};
\draw[->] (A.east) -- (C.north west) node[midway, above] {resol.};
\draw[->] (B.south east) -- (D.west) node[midway, below] {P-resol.\;\;\;\;\;\;};
\draw[->] (C.south west) -- (D.east) node[midway, below] {\;\;\;\;\;\;\;\;\; anti-P-resol.} ;

\end{tikzpicture}
\caption{Anti-P-resolution over $\overline{\Sigma}_3$}
\label{fig:anti-P_example}
\end{figure}

\begin{rem}\label{rem--anti--p}
The notion of anti-P-resolutions is implicitly inspired by Horikawa's argument in \cite{horikawa}.
In that paper, a log $\Q$-Gorenstein smoothing of the pair $(\overline{\Sigma}_3,16\overline{\Gamma})$ to $(\Sigma_1,6\Delta_0+10\Gamma)$ is constructed.
Here, $\Sigma_3$ serves as the P-resolution of $\overline{\Sigma}_3$, and we have $h^1(\mathcal{O}_{\Sigma_3}(6\Delta_0+16\Gamma))=1$.
To proceed, Horikawa considered the blow-up of $\Sigma_3$ at the point where the two components of the branch divisor intersect, and effectively reduced the problem to the $\mathbb{Q}$-Gorenstein smoothability of $\overline{\Sigma}_4$.
This procedure can be reinterpreted in terms of anti-P-resolutions:
that is, Horikawa's method can essentially be viewed as constructing an anti-P-resolution of $\overline{\Sigma}_3$ by contracting the $(-4)$-curve arising from the blow-up of $\Sigma_3$ at the intersection point.
See Figure \ref{fig:anti-P_example}.
\end{rem}

\phantom{A}

\phantom{A}

\phantom{A}

\phantom{A}

\phantom{A}

As an application of anti-P-resolutions, we give the following explicit examples of normal Horikawa surfaces of geometric genus $10$ which are degenerations of smooth Horikawa surfaces of both types $(0)$ and $(6)$:

\begin{exam} \label{exam:connectingGieseker}
Let us fix two distinct fibers $\Gamma$ and $\Gamma'$ of the ruling $\Sigma_8 \to \mathbb{P}^1$.  
We consider a divisor $\widehat{B} = B_0 + \Gamma \in |4\Delta_0 + 36\Gamma|$, constructed as shown in the following figure.

\begin{figure}[H]
\centering
\begin{tikzpicture}[line cap=round,line join=round,>=triangle 45,x=0.4cm,y=0.4cm]
\clip(-2,-3) rectangle (6,2);
\draw [line width=1pt] (-2,0)-- (8,0);
\draw [line width=1pt] (0,2)-- (0,-2);
\draw [line width=1pt,red] (4,-2)-- (4,2);
\draw[line width=1pt,color=red,smooth,samples=100,domain=2.75:5.25,rotate around={90:(4,0)}] plot(\x,{0-(\x - 4)^2});
\draw[line width=1pt,color=red,smooth,samples=100,domain=-0.75:0.75,rotate around={45:(0,0)}] plot(\x,{(\x )^4});
\draw[line width=1pt,color=red,smooth,samples=100,domain=-0.75:0.75,rotate around={45:(0,0)}] plot(\x,{0-(\x )^4});
\begin{scriptsize}
\draw [fill] (2,0.5) node {$\Delta_0$};
\draw [fill] (4,-2.5) node {$\Gamma$};
\draw [fill] (0,-2.5) node {$\Gamma'$};
\end{scriptsize}
\end{tikzpicture}
\caption{}
\end{figure}

The existence of such a divisor is guaranteed by Lemma~\ref{lem--supersingular--basepointfree}~(2).
Let $\eta \colon \Sigma_8 \to \overline{\Sigma}_8$ denote the contraction of $\Delta_0$, and define $B := \eta_*\widehat{B} \in |36\overline{\Gamma}|$.  
Then the double cover $X \to W := \overline{\Sigma}_8$ branched along $B$ defines a standard Horikawa surface of type $(8)'$ with $p_g = 10$ and an elliptic double cone singularity.  
We show that this surface deforms to smooth Horikawa surfaces of both types $(0)$ and $(6)$ by using the results in Section~\ref{sec:deformation}:

\smallskip

\noindent
\textbf{Smoothing to Type $(0)$.}
We first construct a smoothing of $X$ to surfaces of type $(0)$.  
Consider the sequence of six blow-ups
\[
W_6 \to \cdots \to W_1 \to \Sigma_8
\]
constructed as follows (see Figure~\ref{fig:anti_P_p_g=10_1}).  
The first four blow-ups give an embedded resolution of the 4-fold node $B_0 \cap \Gamma'$ on $\widehat{B}$.  
The proper transform of $B_0$ then intersects the exceptional curve of $\psi_4$ at two smooth points.  
Choose one of these points, say $x$, and let $\psi_5$ be the blow-up at $x$.  
Let $\psi_6$ be the blow-up at the intersection of the proper transform of $B_0$ and the exceptional curve of $\psi_5$.
The proper transforms of $\Delta_0$ and the exceptional curves of $\psi_i$ for $i=1,\ldots,5$ form a T-chain of type $[9,2^5]$.  
Let $c \colon W_6 \to W^{-}$ denote the contraction of this chain.  
It is easy to check that $W^{-}$ is a log del Pezzo surface with Picard number $2$ and a T-singularity of type $\frac{1}{49}(1,6)$.
Let $B^{-}$ be the effective divisor on $W^{-}$ such that $K_{W^{-}} + \frac{1}{2}B^{-}$ is the pullback of $K_{W} + \frac{1}{2}B$. 
The birational morphism $\mu \colon (W^{-}, \frac{1}{2}B^{-}) \to (W, \frac{1}{2}B)$ is then an admissible anti-P-resolution.  
By Proposition~\ref{prop--smoothing--c-type}, we obtain a $\mathbb{Q}$-Gorenstein smoothing $\mathscr{W}^{-} \to C$ of $W^{-}$ together with an extension $\mathscr{B}^{-}$ of $B^{-}$, such that $\mathscr{W}^{-}$ is a $\mathbb{Q}$-factorial threefold, and for general $c \in C$, the fiber $\mathscr{W}^{-}_c \cong \Sigma_0$ and $\mathscr{B}^{-}_c$ is smooth.
Taking a log canonical model $(\mathscr{W}, \frac{1}{2}\mathscr{B})$ of $(\mathscr{W}^{-}, \frac{1}{2}\mathscr{B}^{-})$, and considering the double cover $\mathscr{X} \to \mathscr{W}$ branched along $\mathscr{B}$, we obtain a smoothing $\mathscr{X} \to C$ of $X$ to smooth Horikawa surfaces of type $(0)$.

\smallskip

\noindent
\textbf{Smoothing to Type $(6)$.}
Next, we construct a smoothing of $X$ to surfaces of type $(6)$.  
We again consider six blow-ups
\[
W_6 \to \cdots \to W_1 \to \Sigma_8,
\]
this time beginning with an embedded resolution of the 4-fold node $B_0 \cap \Gamma$ on $\widehat{B}$ by the first four blow-ups.  
Let $\psi_5$ be the blow-up at the intersection of the proper transforms of $B_0$ and $\Gamma$, and $\psi_6$ the blow-up at the intersection of the proper transform of $B_0$ and the exceptional curve of $\psi_5$.
The sequence of blow-ups is as illustrated in Figure~\ref{fig:anti_P_p_g=10_2}.
As before, the proper transforms of $\Delta_0$ and the exceptional curves $\psi_i$ for $i=1,\ldots,5$ form a T-chain $[9,2^5]$.  
Let $c \colon W_6 \to W^{-}$ be the contraction of this chain, and let $B^{-}$ be defined similarly.  
Then $W^{-}$ is again a log del Pezzo surface with Picard number $2$ and a T-singularity of type $\frac{1}{49}(1,6)$.
In this case, the proper transform of $\Gamma$ on $W^{-}$ becomes a $(-6)$-curve contained in $B^{-}$, disjoint from all other components.  
By the proof of Proposition~\ref{prop--super--singular--typeC}~(2), we can construct a $\mathbb{Q}$-Gorenstein smoothing $\mathscr{W}^{-} \to C$ and an extension $\mathscr{B}^{-}$ of $B^{-}$, such that the singularity $\frac{1}{49}(1,6)$ is smoothed, while the $(-6)$-curve remains.  
Moreover, $\mathscr{W}^{-}$ is a $\mathbb{Q}$-factorial threefold, with general fiber $\mathscr{W}^{-}_c \cong \Sigma_6$, and $\mathscr{B}^{-}_c$ consists of two disjoint smooth curves: $\Delta_0$ and $\mathscr{B}^{-}_c - \Delta_0$.
Taking the log canonical model $(\mathscr{W}, \frac{1}{2}\mathscr{B})$ of $(\mathscr{W}^{-}, \frac{1}{2}\mathscr{B}^{-})$ and the double cover $\mathscr{X} \to \mathscr{W}$ branched along $\mathscr{B}$, the family $\mathscr{X} \to C$ gives a smoothing of $X$ to smooth Horikawa surfaces of type $(6)$.

\newpage

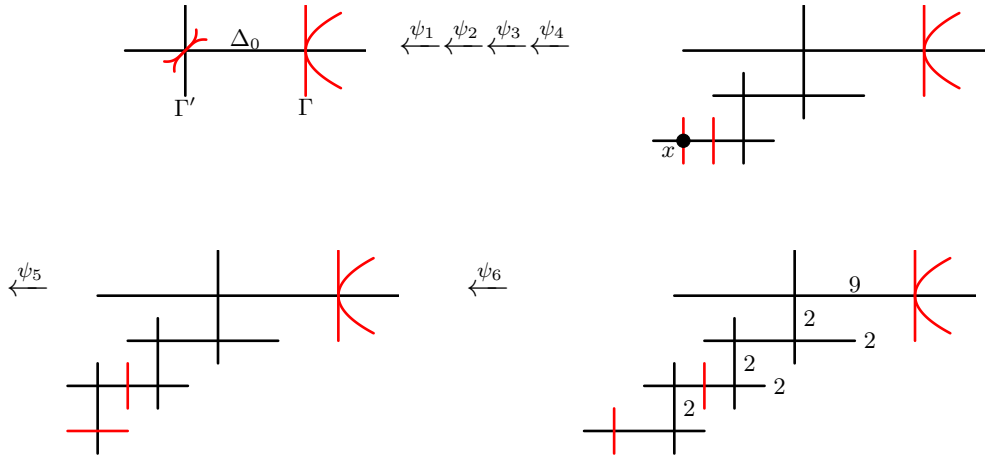
\begin{figure}[H]
\centering
\begin{subfigure}[t]{0.32\linewidth}
\centering
\begin{tikzpicture}[line cap=round,line join=round,>=triangle 45,x=0.4cm,y=0.3cm]
\clip(-8,-6) rectangle (6,2);
\draw [line width=1pt] (-2,0)-- (8,0);
\draw [line width=1pt] (0,2)-- (0,-2);
\draw [line width=1pt,red] (4,-2)-- (4,2);
\draw[line width=1pt,color=red,smooth,samples=100,domain=2.75:5.25,rotate around={90:(4,0)}] plot(\x,{0-(\x - 4)^2});
\draw[line width=1pt,color=red,smooth,samples=100,domain=-0.75:0.75,rotate around={45:(0,0)}] plot(\x,{(\x )^4});
\draw[line width=1pt,color=red,smooth,samples=100,domain=-0.75:0.75,rotate around={45:(0,0)}] plot(\x,{0-(\x )^4});

\begin{scriptsize}
\draw [fill] (2,0.5) node {$\Delta_0$};
\draw [fill] (4,-2.5) node {$\Gamma$};
\draw [fill] (0,-2.5) node {$\Gamma'$};
\end{scriptsize}
\end{tikzpicture}
\end{subfigure}
\hspace{0.05\linewidth}
\raisebox{4.25em}{$\xleftarrow{\psi_1}\xleftarrow{\psi_2}\xleftarrow{\psi_3}\xleftarrow{\psi_4}$}
\hspace{0.02\linewidth}
\begin{subfigure}[t]{0.32\linewidth}
\centering
\begin{tikzpicture}[line cap=round,line join=round,>=triangle 45,x=0.4cm,y=0.3cm]
\clip(-6,-6) rectangle (6,2);
\draw [line width=1pt] (-4,0)-- (8,0);
\draw [line width=1pt] (0,2)-- (0,-3);
\draw [line width=1pt] (-3,-2)-- (2,-2);
\draw [line width=1pt] (-2,-1)-- (-2,-5);
\draw [line width=1pt] (-5,-4)-- (-1,-4);
\draw [line width=1pt,red] (-4,-3)-- (-4,-5);
\draw [line width=1pt,red] (-3,-3)-- (-3,-5);
\draw [line width=1pt,red] (4,-2)-- (4,2);
\draw[line width=1pt,color=red,smooth,samples=100,domain=2.75:5.25,rotate around={90:(4,0)}] plot(\x,{0-(\x - 4)^2});

\begin{scriptsize}
\draw [fill=black] (-4,-4) circle (2.5pt);
\draw [fill] (-4.5,-4.5) node {$x$};
\end{scriptsize}
\end{tikzpicture}
\end{subfigure}

\vspace{2em} 

\raisebox{5.8em}{$\xleftarrow{\psi_5}$}
\hspace{-0.1\linewidth}
\begin{subfigure}[t]{0.4\linewidth}
\centering
\begin{tikzpicture}[line cap=round,line join=round,>=triangle 45,x=0.4cm,y=0.3cm]
\clip(-8,-8) rectangle (6,2);
\draw [line width=1pt] (-4,0)-- (8,0);
\draw [line width=1pt] (0,2)-- (0,-3);
\draw [line width=1pt] (-3,-2)-- (2,-2);
\draw [line width=1pt] (-2,-1)-- (-2,-5);
\draw [line width=1pt] (-5,-4)-- (-1,-4);
\draw [line width=1pt] (-4,-3)-- (-4,-7);
\draw [line width=1pt,red] (-3,-6)-- (-5,-6);
\draw [line width=1pt,red] (-3,-3)-- (-3,-5);
\draw [line width=1pt,red] (4,-2)-- (4,2);
\draw[line width=1pt,color=red,smooth,samples=100,domain=2.75:5.25,rotate around={90:(4,0)}] plot(\x,{0-(\x - 4)^2});

\begin{scriptsize}

\end{scriptsize}
\end{tikzpicture}
\end{subfigure}
\hspace{0.02\linewidth}
\raisebox{5.8em}{$\xleftarrow{\psi_6}$}
\hspace{0.02\linewidth}
\begin{subfigure}[t]{0.32\linewidth}
\centering
\begin{tikzpicture}[line cap=round,line join=round,>=triangle 45,x=0.4cm,y=0.3cm]
\clip(-8,-8) rectangle (6,2);
\draw [line width=1pt] (-4,0)-- (8,0);
\draw [line width=1pt] (0,2)-- (0,-3);
\draw [line width=1pt] (-3,-2)-- (2,-2);
\draw [line width=1pt] (-2,-1)-- (-2,-5);
\draw [line width=1pt] (-5,-4)-- (-1,-4);
\draw [line width=1pt] (-4,-3)-- (-4,-7);
\draw [line width=1pt] (-3,-6)-- (-7,-6);
\draw [line width=1pt,red] (-3,-3)-- (-3,-5);
\draw [line width=1pt,red] (-6,-5)-- (-6,-7);
\draw [line width=1pt,red] (4,-2)-- (4,2);
\draw[line width=1pt,color=red,smooth,samples=100,domain=2.75:5.25,rotate around={90:(4,0)}] plot(\x,{0-(\x - 4)^2});

\begin{scriptsize}
\draw [fill] (2,0.5) node {$9$};
\draw [fill] (-0.5,-4) node {$2$};
\draw [fill] (2.5,-2) node {$2$};
\draw [fill] (-3.5,-5) node {$2$};
\draw [fill] (-1.5,-3) node {$2$};
\draw [fill] (0.5,-1) node {$2$};
\end{scriptsize}
\end{tikzpicture}
\end{subfigure}
\caption{Smoothing to Type $(0)$.}
\label{fig:anti_P_p_g=10_1}
\end{figure}

\begin{figure}[H]
\centering

\begin{subfigure}[t]{0.32\linewidth}
\centering
\begin{tikzpicture}[line cap=round,line join=round,>=triangle 45,x=0.4cm,y=0.3cm]
\clip(-2,-6) rectangle (6,2);
\draw [line width=1pt] (-4,0)-- (8,0);
\draw [line width=1pt] (0,2)-- (0,-2);
\draw [line width=1pt,red] (4,-2)-- (4,2);
\draw[line width=1pt,color=red,smooth,samples=100,domain=2.75:5.25,rotate around={90:(4,0)}] plot(\x,{0-(\x - 4)^2});
\draw[line width=1pt,color=red,smooth,samples=100,domain=-0.75:0.75,rotate around={45:(0,0)}] plot(\x,{(\x )^4});
\draw[line width=1pt,color=red,smooth,samples=100,domain=-0.75:0.75,rotate around={45:(0,0)}] plot(\x,{0-(\x )^4});

\begin{scriptsize}
\draw [fill] (2,0.5) node {$\Delta_0$};
\draw [fill] (4,-2.5) node {$\Gamma$};
\draw [fill] (0,-2.5) node {$\Gamma'$};
\end{scriptsize}
\end{tikzpicture}
\end{subfigure}
\hspace{0.02\linewidth}
\raisebox{4.25em}{$\xleftarrow{\psi_1}\xleftarrow{\psi_2}\xleftarrow{\psi_3}\xleftarrow{\psi_4}$}
\hspace{0.02\linewidth}
\begin{subfigure}[t]{0.32\linewidth}
\centering
\begin{tikzpicture}[line cap=round,line join=round,>=triangle 45,x=0.4cm,y=0.3cm]
\clip(-2,-6) rectangle (12,2);
\draw [line width=1pt] (-2,0)-- (8,0);
\draw [line width=1pt] (0,2)-- (0,-2);
\draw [line width=1pt] (4,2)-- (4,-3);
\draw [line width=1pt] (7,-2)-- (2,-2);
\draw [line width=1pt] (6,-1)-- (6,-5);
\draw [line width=1pt] (9,-4)-- (5,-4);
\draw [line width=1pt,red] (8,-3)-- (8,-7);
\draw [line width=1pt,red] (7,-3)-- (7,-5);
\draw[line width=1pt,color=red,smooth,samples=100,domain=-0.75:0.75,rotate around={45:(0,0)}] plot(\x,{(\x )^4});
\draw[line width=1pt,color=red,smooth,samples=100,domain=-0.75:0.75,rotate around={45:(0,0)}] plot(\x,{0-(\x )^4});

\begin{scriptsize}

\end{scriptsize}
\end{tikzpicture}
\end{subfigure}

\vspace{2em} 

\raisebox{5.75em}{$\xleftarrow{\psi_5}$}
\hspace{0.02\linewidth}
\begin{subfigure}[t]{0.35\linewidth}
\centering
\begin{tikzpicture}[line cap=round,line join=round,>=triangle 45,x=0.4cm,y=0.3cm]
\clip(-2,-8) rectangle (12,2);
\draw [line width=1pt] (-2,0)-- (8,0);
\draw [line width=1pt] (0,2)-- (0,-2);
\draw [line width=1pt] (4,2)-- (4,-3);
\draw [line width=1pt] (7,-2)-- (2,-2);
\draw [line width=1pt] (6,-1)-- (6,-5);
\draw [line width=1pt] (9,-4)-- (5,-4);
\draw [line width=1pt] (8,-3)-- (8,-7);
\draw [line width=1pt,red] (7,-6)-- (11,-6);
\draw [line width=1pt,red] (7,-3)-- (7,-5);
\draw[line width=1pt,color=red,smooth,samples=100,domain=-0.75:0.75,rotate around={45:(0,0)}] plot(\x,{(\x )^4});
\draw[line width=1pt,color=red,smooth,samples=100,domain=-0.75:0.75,rotate around={45:(0,0)}] plot(\x,{0-(\x )^4});

\begin{scriptsize}

\end{scriptsize}
\end{tikzpicture}
\end{subfigure}
\hspace{0.05\linewidth}
\raisebox{5.75em}{$\xleftarrow{\psi_6}$}
\hspace{0.02\linewidth}
\begin{subfigure}[t]{0.35\linewidth}
\centering
\begin{tikzpicture}[line cap=round,line join=round,>=triangle 45,x=0.4cm,y=0.3cm]
\clip(-2,-8) rectangle (12,2);
\draw [line width=1pt] (-2,0)-- (8,0);
\draw [line width=1pt] (0,2)-- (0,-2);
\draw [line width=1pt] (4,2)-- (4,-3);
\draw [line width=1pt] (7,-2)-- (2,-2);
\draw [line width=1pt] (6,-1)-- (6,-5);
\draw [line width=1pt] (9,-4)-- (5,-4);
\draw [line width=1pt] (8,-3)-- (8,-7);
\draw [line width=1pt] (7,-6)-- (11,-6);
\draw [line width=1pt,red] (10,-5)-- (10,-7);
\draw [line width=1pt,red] (7,-3)-- (7,-5);
\draw[line width=1pt,color=red,smooth,samples=100,domain=-0.75:0.75,rotate around={45:(0,0)}] plot(\x,{(\x )^4});
\draw[line width=1pt,color=red,smooth,samples=100,domain=-0.75:0.75,rotate around={45:(0,0)}] plot(\x,{0-(\x )^4});
\begin{scriptsize}
\draw [fill] (2,0.5) node {$9$};
\draw [fill] (7.5,-2) node {$2$};
\draw [fill] (4,-3.5) node {$2$};
\draw [fill] (9.25,-4) node {$2$};
\draw [fill] (6,-5.5) node {$2$};
\draw [fill] (8,-7.5) node {$2$};
\draw [fill] (10,-7.5) node {$6$};
\end{scriptsize}
\end{tikzpicture}
\end{subfigure}
\caption{Smoothing to Type $(6)$.}
\label{fig:anti_P_p_g=10_2}
\end{figure}
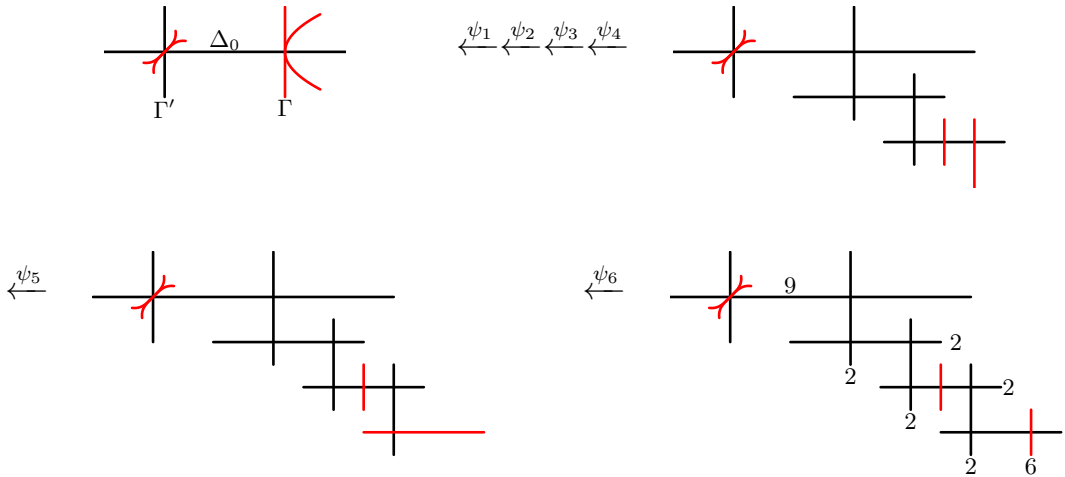

\end{exam}

\subsection{Relation to other works and further questions}\label{subsec--other--relate}

Describing the KSBA moduli of slc Horikawa surfaces is not only interesting in itself but also related to various problems. 
We will explain these relationships in turn.

\par
\vspace{1\baselineskip}

\noindent\textbf{Horikawa problem}\par
\vspace{0.5\baselineskip} 

From the point of view of differential topology, Horikawa surfaces are interesting objects.
Horikawa \cite{horikawa} showed that there are two Gieseker components if and only if $K^2$ is a multiple of eight, and that, if $K^2\equiv 8$ mod $16$, two Horikawa surfaces, each lying in a different connected component of the moduli space, are not diffeomorphic.
His observation leads to a long-standing question, the {\it Horikawa problem}, asking whether there exist two diffeomorphism types for Horikawa surfaces with the same geometric genus in the remaining cases. 
The difficulty arises from the fact that all known invariants in differential topology coincide. 
This problem has been discussed extensively; see e.g. \cite{Cat3,CM,FS,Au,RR,MNU}. 

One possible way to address the Horikawa problem is to find common degenerations of two surfaces, respectively lying in the two components of the Gieseker moduli space.
Manneti \cite[Theorem 1.5]{Manetti2} showed that all smooth surfaces in a connected component of the moduli space parametrizing stable surfaces with only T-singularities have the same diffeomorphism type.
Motivated by his result, Monreal-Negrete-Urz\'ua \cite{MNU} classified $\Q$-Gorenstein smoothable klt Horikawa surfaces (with $p_g\geq10$).
They showed that the two components of the Gieseker moduli space of Horikawa surfaces cannot be connected using only $\mathbb{Q}$-Gorenstein smoothable klt surfaces for any case where the canonical volume is divisible by eight.
In contrast, Theorem~\ref{thm--connectedness} states that, when $p_g=10$, the two components are connected using $\mathbb{Q}$-Gorenstein smoothable lc surfaces.
Although it is not yet fully understood how degenerations to lc surfaces affect the differential topology, we expect that our result will help address the Horikawa problem.
\par
\vspace{1\baselineskip}

\noindent\textbf{Surfaces near the Noether line}\par
\vspace{0.5\baselineskip} 

The main result of this paper can be viewed as a generalization of Horikawa's first paper \cite{horikawa} in the series~\cite{horikawa,horikawa-ii,horikawa-iii,horikawa-iv,horikawa-v}.
The subsequent papers in the series investigate surfaces of general type near the Noether line, especially those satisfying $K^2=2p_g-3$ and $K^2=2p_g-2$.
In this context, it is natural to consider the classification of normal stable degenerations of surfaces with $K^2=2p_g-3$ as the next step.
Partial progress in this direction has been made by various researchers \cite{FPR, FPRR,CFPRR,GPSZ,Rana,FG1,CP},
mostly focusing on klt or Gorenstein surfaces.
Using the techniques developed in this paper, we expect to achieve a complete classification of normal stable degenerations of surfaces with $K^2=2p_g-3$:
Let $X$ be a normal stable surface with $K^2=2p_g-3$, assuming that all singularities of $X$ are $\Q$-Gorenstein smoothable.
Then, by the log Noether inequality (Theorem~\ref{normalstable}), $X$ must fall into one of the following three cases:
\begin{itemize}
\item 
$X$ is Gorenstein,
\item 
The minimal resolution $\widetilde{X}$ of $X$ contracts to a Gorenstein Horikawa surface,
\item 
$\widetilde{X}$ admits an elliptic fibration.
\end{itemize}
The first case is expected to be treatable in a similar way to Horikawa's original work \cite{horikawa-ii}.
The classification of Horikawa surfaces in the present paper contributes to the analysis of the second case.
Compared to the case of Horikawa surfaces, non-Gorenstein double cone singularities play a crucial role in this case (\cite{CP}).
Most instances of the third case have already been addressed in Section~\ref{sec:non-std_Horikawa}, except for the case where $n_X=2$.
Note that $X$ may not admit a ``good'' involution even if it is $\Q$-Gorenstein smoothable, since, for example, general quintic surfaces do not carry an involution canonically.
We plan to pursue the classification of normal stable surfaces with $K^2=2p_g-3$ and only $\Q$-Gorenstein smoothable singularities in future work.
The case $K^2=2p_g-2$ also appears to be of interest, as these surfaces lie on the borderline of the Noether inequality for the case of canonical pencils (Theorem~\ref{normalstable}~(1.2)).

\newpage

\par
\vspace{1\baselineskip}

\noindent\textbf{Geography of threefolds of general type and slope inequality}\par
\vspace{0.5\baselineskip} 

Horikawa~\cite{horikawa-genus2} investigated the connection between surfaces lying near the Noether line and genus-two fibrations: 
For any relatively minimal genus-two fibration $f\colon X\to B$, the {\em slope inequality}
$K_{X/B}^{2}\ge 2\deg f_{*}\omega_{X/B}$
holds, which specializes to the Noether inequality $K_X^{2}\ge 2\chi(\O_X)-6$ when $B=\mathbb{P}^1$.
In addition, he introduced a refined version of this inequality known as the {\em slope equality}:
$$
K_{X/B}^{2}=2\deg f_{*}\omega_{X/B}+\sum_{p\in B}\mathrm{H}(f^{-1}(p)),
$$
where each $\mathrm{H}(f^{-1}(p))$ is a non-negative integer determined by the fiber germ over $p$, called the {\em Horikawa index}.
He also explicitly computed these indices for all possible fiber germs.
Because most surfaces near the Noether line admit a natural genus-two fibration over $\mathbb{P}^1$, the study of such surfaces is closely tied to the study on genus $2$ fibrations.

This theory has been generalized in various directions.
For example, Cornalba-Harris~\cite{CoHa} and Xiao~\cite{Xia} extended the slope inequality to fibrations of arbitrary genus $g\ge 2$, giving:
$$
K_{X/B}^{2}\ge \frac{4(g-1)}{g}\deg f_{*}\omega_{X/B}.
$$
Slope equalities and Horikawa indices have also been developed for specific classes such as hyperelliptic fibrations, non-hyperelliptic genus-three fibrations and so on (see the survey \cite{AsKo}).
These developments have significant applications in the geography of surfaces of general type \cite{Rei, Kon2, Pardini2}.
A key insight is that the lower bound of the slope $K_{X/B}^{2}/\deg f_{*}\omega_{X/B}$ reflects the geometry of the general fiber.

On the threefold side, the {\em Noether inequality for algebraic threefolds} has recently been completely established \cite{CCJ,CCJ+,CHJ}: 
$$
K_{X}^{3}\ge \frac{4}{3}p_g(X)-\frac{10}{3}
$$
holds for minimal threefolds of general type $X$.
According to \cite{CHPZ}, any threefold achieving the equality of the above inequality---so called {\em Noether threefolds}---admits a fibration $f\colon X\to \mathbb{P}^1$ whose general fiber is a surface with $(K^2, p_g)=(1,2)$.
Furthermore, a slope inequality of the form $K_{X/B}^{3}\ge \frac{4}{3}\deg f_{*}\omega_{X/B}$ has been studied for fibrations by $(1,2)$-surfaces in \cite{HuZh}.
Thus, as in the surface case, understanding degenerations of $(1,2)$-surfaces is essential for studying threefolds of general type near the Noether line.
We expect that our classification techniques for stable Horikawa surfaces with only $\Q$-Gorenstein smoothable singularities are useful for analyzing threefolds near the Noether line.

In another direction, Kobayashi \cite{Kob} showed that for a minimal threefold of general type $X$, the inequality 
$$
K_X^{3}\ge 2p_g(X)-6
$$
holds when the canonical map is generically finite onto its image.
This also can be seen as an analogue of the Noether inequality for surfaces.
Moreover, most threefolds that achieve the equality in this inequality admit a fibration $X\to \mathbb{P}^1$ whose general fiber is a surface with $(K^2,p_g)=(2,3)$.
According to \cite{HuZh}, the corresponding slope inequality $K_{X/B}^{3}\ge 2\deg f_{*}\omega_{X/B}$ is satisfied for any relatively minimal fibration $f\colon X\to B$ of $(2,3)$-surfaces.
Therefore, degenerations of Horikawa surfaces also play a crucial role in the study of threefolds of general type.

While slope inequalities for fibered threefolds have been studied extensively (\cite{Ohn, Bar, HuZh, Aka}), and even in higher dimensions (e.g., \cite{Bar2, HuZh2, CTV}),  the theory of slope equalities, classification of degenerate fibers, and Horikawa indices in the threefold setting remains largely undeveloped.
Since the theory of slope (in)equalities is closely connected to the enumerative geometry of moduli spaces via classifying maps of fibrations (cf. \cite{CoHa, Mor, Enokizono2}), we believe that our results on the KSBA moduli spaces of Horikawa surfaces, together with computations involving families of Horikawa surfaces, can contribute significantly to the study of the enumerative geometry of their moduli spaces.
For Campedelli and Burniat surfaces, similar attempts have been made in \cite{Alexeev}.

\par
\vspace{1\baselineskip}

\noindent\textbf{Torelli-type problems for surfaces of general type}\par
\vspace{0.5\baselineskip} 

Torelli-type problems for surfaces of general type have been studied extensively by many researchers (see e.g. \cite{CFPR,FG,Murakami,Pardini,PZ,Reider1,Reider2,Usui2,Cat1,Cat2,Ky}). 
They ask how much the variations of Hodge structures recover the deformations of complex structures.
Surfaces of general type exhibit a wide range of behaviors for each class, 
and accordingly, Torelli-type problems for them become complicated. 
In fact, while Torelli-type theorems have been established for various classes of surfaces of general type (see e.g. \cite{FG,Murakami,Pardini,PZ,Reider1,Reider2,Usui2}), counterexamples are also known (see e.g.~\cite{Cat1,Cat2,Ky}).

For algebraic curves, there exists an approach to Torelli-type theorems that makes use of the degenerations of curves (see e.g.~the chapter by Friedman in \cite{Griffiths}, \cite{Usui}).
Consider the Deligne-Mumford compactification of the moduli space of curves of genus $g$. 
One analyzes the asymptotic behavior of the period map for curves degenerating to an irreducible stable curve $C$ with a single node.
The geometric genus of $C$ is $g-1$, and the limit of the period map at the point corresponding to $C$ can be described in terms of the periods of the normalization $\widetilde{C}$ and the conductor.
This allows one to reduce the injectivity of the period map at $C$ to that of $\widetilde{C}$, a curve of smaller genus.

It is natural to ask whether such an idea can be applied to Torelli-type problems for surfaces.
In fact, Friedman \cite{Fri1} proved the global Torelli theorem for polarized K3 surfaces using this approach. 
In the case of surfaces of general type, one approach is to analyze the asymptotic behavior of the period map near the boundary of the KSBA compactification, and to attempt to reduce the Torelli-type problems to those for surfaces with smaller numerical invariants.
Such an approach has been explored by Friedman and Griffiths \cite{FG} in the case of $(1,2)$-surfaces.
They gave a detailed description of deformations of $(1,2)$-surfaces with certain elliptic singularities in \cite{FG1}, 
and studied the asymptotic behavior of the period map when such singularities acquire. 
We expect that our classification method will serve as a first step in applying this approach to Horikawa surfaces.

\newpage

\par
\vspace{1\baselineskip}

\noindent\textbf{Application of anti-P-resolutions I: wall-crossing theories for moduli spaces}\par
\vspace{0.5\baselineskip}

The wall-crossing theory for the KSBA moduli, established by \cite{ABIP,MZ}, states that if one considers the KSBA moduli of pairs $(W,cB)$ and varies the coefficient $c$ continuously, changes of the moduli space, so-called wall-crossings, occur only finitely many times.
While this technique is effectively used to describe some moduli spaces (cf.~\cite{Al-Pa}), explicit descriptions of KSBA moduli wall-crossing are rare when $K_W$ and $B$ are not proportional.
In Theorem \ref{thm--moduli--seminormaliazation--intro}, we can regard $\mathcal{M}^{\mathrm{nq}}_{p_g}$ as the moduli corresponding to the wall $c=\frac{1}{2}$.
In the case of $\mathcal{M}^{\mathrm{nq}}_{p_g}$, the wall $c=\frac{1}{2}$, anti-P-resolutions are actually applicable to describe the locus in the moduli space parameterizing normal pairs, as we explained in \S\ref{intro-subsec:smoothing_and_anti_P}.
As in this case, if $B$ is not proportional to $K_W$, we have to deal with a log $\Q$-Gorenstein deformation of the pair that induces a non-$\Q$-Gorenstein deformation of $W$ itself on some walls.
We expect that the theory of anti-P-resolutions will be helpful to describe the KSBA moduli wall-crossing phenomenon in some explicit cases.

We also expect that we can apply the theory of anti-P-resolutions to describe K-moduli wall-crossing phenomena explicitly. 
Here, K-stability is one of the main tools to construct a moduli space of klt log Fano pairs, known as K-moduli. See \cite{X} for more details.
Since Ascher-DeVleming-Liu \cite{ADL} established the K-moduli wall-crossing of plane curves---describing how the K-moduli of log Fano pairs change as we vary the coefficients---there has been growing interest in the explicit structure of K-moduli spaces and their wall-crossings.
This development is particularly important for understanding concrete moduli problems, such as those arising in the Hassett-Keel program (see \cite{Zhao,ADLW}).
For general theory in the proportional case, we refer to \cite{Z2,Zhou}.
Recently, Liu-Zhou \cite{LZ} partially established the theory of non-proportional K-moduli wall-crossing.
Although the K-moduli wall-crossing when $K_W$ and $B$ are not proportional is also useful to describe some moduli spaces explicitly (see \cite{DV+} for example), explicit descriptions of such K-moduli wall-crossing are rare.
To describe such K-moduli wall-crossing explicitly for klt log Fano surfaces, there would be the same difficulty as we explained in the non-proportional KSBA moduli wall-crossing.
To address these problems, anti-P-resolutions should also be useful.
\par
\vspace{1\baselineskip}

\noindent\textbf{Application of anti-P-resolutions II: equivariant cusp smoothing}\par
\vspace{0.5\baselineskip} 

We observe that anti-P-resolutions can be effectively applied to the problem of equivariant cusp smoothability in cases where a finite group acts non-freely outside the cusp singularity.
Although Simonetti \cite{Simonetti} and Jiang \cite{J} investigated equivariant smoothability under the assumption that the group acts freely outside the singularity with additional conditions (see Remark~\ref{rem--equiv--smooth--detail} for more details), the above cases have not been fully understood.
For example, the group $\mathbb{Z}/2\mathbb{Z}$ naturally acts non-freely outside the singularity on double cone singularities.
In this setting, we apply the theory of anti-P-resolutions to construct examples of equivariantly smoothable cusp singularities (see Proposition~\ref{prop--smoothing--c-type}, Remark~\ref{rem--smoo}).
We expect that anti-P-resolutions offer a useful approach to studying equivariant smoothability more broadly for cusp singularities with finite group actions.
\par
\vspace{1\baselineskip}

\subsection{Structure of the paper}

In Section~\ref{sec:preliminaries}, we collect fundamental concepts and results used throughout this paper.  
More specifically, we review basic definitions and relevant facts on the following topics: log pairs, $\Q$-Gorenstein deformations, KSBA moduli, log canonical singularities, and double covers.

In Section \ref{sec:extT}, we introduce the notion of extended T-chains.
Extended T-chains appear naturally in the study of normal stable Horikawa surfaces.
While their analysis is necessary for the classification of Horikawa surfaces, we avoid it here and instead introduce basic concepts such as T-trains and P-admissible chains, along with their geometric interpretation.

In Section \ref{sec:geography_normal_stable}, we investigate relations among the numerical invariants of normal stable surfaces.
In particular, we establish two fundamental results: a log Noether inequality for normal stable surfaces (Theorem~\ref{normalstable}) and a fundamental equality (Theorem~\ref{ellfundeq}) for normal stable surfaces with only $\Q$-smoothable singularities admitting an elliptic pencil.
These (in)equalities play a crucial role in the classification of normal stable Horikawa surfaces with only $\Q$-Gorenstein smoothable singularities.

From Section \ref{sec:normal_stbale_surface} onward, we focus on normal stable Horikawa surfaces (Definition \ref{def:Horikawa}).
We introduce several concepts that will play significant roles in studying classifications and deformations of these surfaces.

In Section \ref{subsec:Horikawa}, we divide the class of normal stable Horikawa surfaces into two classes based on the dimension of their canonical images: \emph{standard} (dimension two canonical image) and \emph{non-standard} (dimension one canonical image).
Standard Horikawa surfaces directly generalize the surfaces studied by Horikawa \cite{horikawa}, as shown in Proposition \ref{prop:standard_Horikawa}.
On the other hand, non-standard Horikawa surfaces with only $\Q$-Gorenstein smoothable singularities fall into three types using the log Noether inequality (Proposition \ref{prop:non-std_Horikawa}).

In Section \ref{subsec:good_involution}, we introduce the notion of \emph{good involutions} (Definition \ref{defn--good-involution}).
Although the definition of a good involution is ad hoc, we observe that all good involutions share useful properties (Proposition~\ref{prop:good_involution_uniqueness}~(1), Proposition~\ref{prop--involution}).
We provide a structure theorem for non-standard Horikawa surfaces with good involutions (Proposition~\ref{prop:good_involution_uniqueness}~$(2)$).

We classify standard Horikawa surfaces in Section \ref{sec:std_Horikawa}, and non-standard Horikawa surfaces with only $\Q$-Gorenstein smoothable singularities, as well as those admitting good involutions, in Section \ref{sec:non-std_Horikawa}.
These classifications provide partial answers to Questions \ref{question} (1) and (2).
The remaining task is to determine when a surface in these classifications is $\Q$-Gorenstein smoothable; this will be addressed in Section \ref{sec:deformation}.
We explain the classifications in order.

In Section \ref{sec:std_Horikawa}, we introduce two classes of singularities, namely \emph{mild singularities} and \emph{double cone singularities}, in order to analyze the singularities appearing on a normal stable Horikawa surface $X$ with a good involution.

In Sections~\ref{sec:mild_singularity} and~\ref{sec:double_cone}, we define mild singularities (Definition~\ref{defn--mild--sing}) and double cone singularities (Definition~\ref{defn:(d;k_1,k_2)_sing}); 
we classify mild singularities (Proposition~\ref{lem:classify_mild1}) and Gorenstein double cone singularities (Proposition~\ref{prop:classification_cone_sing}).
In Section~\ref{sec:classification_standard_Horikawa}, we prove a structural theorem for standard Horikawa surfaces (Theorem~\ref{thm:standard_Horikawa}).

In Section~\ref{sec:non-std_Horikawa}, we classify non-standard Horikawa surfaces with only $\Q$-Gorenstein smoothable singularities, as well as those admitting good involutions.

In Section~\ref{subsec:overview_non-std}, we present an overview of our classification strategy.
According to Proposition \ref{prop:standard_Horikawa}, the structure theorem for non-standard Horikawa surfaces with only $\Q$-Gorenstein smoothable singularities, these surfaces fall into three types.
In Sections~\ref{subsec:0,1,0}, \ref{subsec:0,2,0}, \ref{subsec:1,1,1}, we separately address the classification of surfaces of each type, while in Sections~\ref{subsec:0,1,0,involution}, \ref{subsec:0,2,0,involution}, \ref{subsec:1,1,1,involution}, we do so for those admitting good involutions.
We also explain that we can classify a broader class of surfaces than Horikawa surfaces by carrying out our strategy.

In Section~\ref{subsec:non-standard_Horikawa}, we apply the results from Sections~\ref{subsec:0,1,0}--\ref{subsec:1,1,1,involution} to obtain classifications of non-standard Horikawa surfaces with only $\Q$-Gorenstein smoothable singularities and those admitting good involutions (Theorems~\ref{thm:classification_non-standard_Horikawa_noninv}, \ref{thm:classification_non-standard_Horikawa}). 
Note that the former provides an answer to Question~\ref{question} (1), and the latter shows Theorem~\ref{intro-qGsmHor}.
In Construction~\ref{construction}, we also describe explicit constructions of the surfaces listed in Theorem~\ref{intro-qGsmHor}. 

In Section \ref{sec--deformation-arekore-honke}, we establish several general results on $\Q$-Gorenstein deformations of normal surfaces to deduce our results on $\Q$-Gorenstein smoothability for normal stable Horikawa surfaces collected in Sections \ref{sec:std_Horikawa} and \ref{sec:non-std_Horikawa}.
This section consists of four parts.
In Section \ref{sec--deformation--arekore}, we generalize methods of \cite{DVS} to provide a complete criterion for divisorial sheaves on a one-parameter $\mathbb{Q}$-Gorenstein smoothing of T-singularities.
In Section \ref{subsec:invol-deform-8}, we provide a criterion for the log $\Q$-Gorenstein smoothability for surface pairs with some conditions by using the above criterion.
In Section \ref{subsec:rational-surf-deform-8}, we provide several lemmas to investigate the properties of the log $\Q$-Gorenstein smoothing family of rational surfaces.
Finally, we define anti-P-resolutions and admissible anti-P-resolutions, observe their properties, and classify all admissible anti-P-resolutions in Section \ref{subsec:anti-P-8}.

In Section \ref{sec:deformation}, we discuss $\Q$-Gorenstein deformation of normal stable Horikawa surfaces described in Sections \ref{sec:std_Horikawa} and \ref{sec:non-std_Horikawa}.
In Section~\ref{subsec:eliip-cone--smoothable}--\ref{subsec:LPS-p_g>6--smoothable}, we examine the $\Q$-Gorenstein smoothability for each class of normal stable Horikawa surfaces listed in Theorem~\ref{intro-qGsmHor}.
In Section~\ref{subsec:conclusion--equiv--smoothability}, we prove Theorems~\ref{thm--main--ii} and \ref{intro-thm:Q2}, answering Quenstion~\ref{question}~(2).
We also remark the equivariant cusp smoothing problem.
In Section \ref{subsec:comparison}, we discuss seminormalizations of the two moduli spaces $\overline{M}^{\mathrm{Gie}}_{2p_g-4,p_g}$ and $M^{\mathrm{nq}}_{p_g}$ and show Theorem \ref{thm--moduli--seminormaliazation--intro}.
Finally, we show Theorems \ref{thm--stratification--intro} and \ref{thm--connectedness} in Section \ref{subsec:stratification}.

In Appendix~\ref{app:sing_fiber_inv}, we classify involutions on elliptic fiber germs $Y \to T$ that satisfy certain conditions.
We use this classification in Section~\ref{sec:non-std_Horikawa} to study the quotient of a non-standard Horikawa surface by a good involution.

In Appendix \ref{app:cusp}, we determine which cusp singularities that appear on the surfaces studied in Section \ref{sec:non-std_Horikawa} are smoothable.

In Appendix \ref{app:drill}, we collect the classifications of certain extended T-chains required in Sections~\ref{sec:non-std_Horikawa}, \ref{sec--deformation-arekore-honke}.

\begin{ack}
The authors especially thank Kenta Hashizume for fruitful discussions and a lot of insightful comments from the earliest stage of this work.
We are grateful to Giancarlo Urz\'ua for introducing us the notion of wormhole singularities and for providing helpful comments.
The authors also thank Hamid Abban, Tadashi Ashikaga, Jonathan Evans, Kazuhiro Konno, Yuji Odaka, Rita Pardini, S\"{o}nke Rollenske and Aline Zanardini for useful comments.

H.A.~is supported by JSPS KAKENHI No.24KJ0011.
M.E.~is supported by JSPS KAKENHI No.25K06926.
M.H.~is partially supported by Royal Society International Collaboration Award ICA\textbackslash1\textbackslash231019.
\end{ack}

\section{Preliminaries}
\label{sec:preliminaries}

\subsection*{Notations and convention}
We work over the field of complex numbers $\mathbb{C}$ throughout this paper.
We collect the following terminologies, which we will use in this paper.

\begin{enumerate}
    \item
    In this paper, we assume that all schemes are locally Noetherian over $\mathbb{C}$.
    By a {\em variety}, we mean a quasi-projective scheme over $\mathbb{C}$ that is irreducible and reduced.
    Unless stated otherwise, a {\em surface} refers to a reduced scheme of dimension $2$ that is quasi-projective over $\mathbb{C}$.
    \item
    Let $S$ be a scheme. For any point $s\in S$, we define  $\kappa(s)$ to be the residue field $\mathcal{O}_{S,s}/\mathfrak{m}_s$, where $\mathfrak{m}_s$ is the maximal ideal at $s$.
    \item
    Let $X$ be a normal surface and let $L$ be a $\Q$-Cartier divisor on $X$. 
    Then, the {\it Iitaka dimension} of $L$ is denoted by $\kappa(L)$.
    The {\it Kodaira dimension} of $X$ is defined as $\kappa(K_{\tilde{X}})$, where $\tilde{X}$ is the minimal resolution of $X$.
    \item
    We freely use fundamental notions of stacks from \cite{Ols}.
    If $\mathfrak{X}$ is a stack over $\mathbb{C}$, then let $\mathfrak{X}(\mathbb{C})$ denote the collection of objects of $\mathbb{C}$-valued points of $\mathfrak{X}$.
    \item
    $\Sigma_d$ denotes the Hirzebruch surface $\mathbb{P}_{\mathbb{P}^1}(\mathcal{O}\oplus\mathcal{O}(-d))$ of degree $d$.
    Let $\Delta_0$ denote a section satisfying $\Delta_0^2=-d$ and let $\Gamma$ denote a fiber over $\mathbb{P}^1$.
    If $d\ge1$, then $\Delta_0$ is uniquely determined. 
    In this case, we denote by $\overline{\Sigma}_{d}$ the normal surface obtained by contracting $\Delta_0$.
    Let $\overline{\Gamma}$ denote the image of $\Gamma$ by this contraction, which is a generator of the class group of $\overline{\Sigma}_{d}$.
    \item
    For singular fibers of relatively minimal elliptic surfaces, we follow Kodaira's notation (cf.~\cite[V.~Section 7]{BHPV}). 
    \item
    Let $f\colon X\to Y$ be a projective morphism.
    If $f_*\mathcal{O}_X\cong \mathcal{O}_Y$, then we say that $f$ is a {\it contraction}.
    In this case, the relative Picard group is defined as $\mathrm{Pic}(X/Y)=\mathrm{Pic}(X)/\sim_Y$, and if $X$ is normal, the relative class group is defined as $\mathrm{Cl}(X/Y)=\mathrm{Cl}(X)/\sim_Y$, where $\sim_Y$ denotes linear equivalence over $Y$.
    We write $D_1\sim_{Y,\mathbb{Q}}D_2$ if $mD_1\sim_{Y} mD_2$ for some positive integer $m$. 
    When $Y=\mathrm{Spec}\,\mathbb{C}$, we simply write $D_1\sim_{\mathbb{Q}}D_2$.
    We write $D_1\equiv D_2$ to mean $D_1$ and $D_2$ are numerically equivalent over $Y=\mathrm{Spec}\,\mathbb{C}$.
    \item 
    If $X$ is a projective surface with only slc singularities, we set $p_g(X):=h^0(X,\omega_X)$ (resp.~$q(X):=h^1(X,\mathcal{O}_X)$) and call this the {\em geometric genus} (resp.~the {\em irregularity}) of $X$.
    If there is no fear of confusion, we will simply write $p_g(X)$, $\chi(\mathcal{O}_X)$ and $q(X)$ as $p_g$, $\chi$ and $q$ respectively.
    \item 
    Let $X$ be a proper reduced and irreducible complex space.
    The {\it Picard number} $\rho(X)$ of $X$ is defined as the rank of the N\'eron-Severi group $\mathrm{NS}(X)$, which is the image of $\mathrm{Pic}(X)\to H^2(X,\mathbb{Z})$.
    We denote $\mathrm{Eff}(X)\subset \mathrm{NS}(X)$ as the {\it effective cone}.
    If $f\colon X\to S$ be a projective contraction, the {\em relative Picard number} $\rho(X/S)$ is the rank of the numerical classes of curves contracted by $f$ (cf.~\cite[2.16]{KoMo}). 
    \item
    Let $f\colon X\to S$ be a flat projective morphism.
    Then, we denote by $\mathrm{Aut}_S(X)$ the {\it automorphism group scheme} of $X$ over $S$.
    See \cite[Section 5]{FGA} for more details.
    \item
    Let $X$ be a projective variety and let $|L|$ be the complete linear system.
    Then we denote by $\varphi_L\colon X\dasharrow  \mathbb{P}^{h^0(L)-1}$ the rational map induced by $|L|$.
\end{enumerate}

\subsection{Fundamental notions of log pairs}

We collect fundamental notions of log pairs.

\begin{defn}[Singularities of normal pairs]
    Let $X$ be a normal variety with an effective $\mathbb{Q}$-divisor $\Delta$ on $X$ such that $K_X+\Delta$ is $\mathbb{Q}$-Cartier.
    We say that a projective birational morphism $\pi\colon Y\to X$ is a {\em log resolution} of $(X,\Delta)$ if $Y$ is smooth and $\pi_*^{-1}\Delta+\mathrm{Ex}(\pi)$ is simple normal crossing, where $\mathrm{Ex}(\pi)$ is the support of the exceptional locus of $\pi$.
    Then, we define the following conditions.
        Let $\pi\colon Y\to X$ be a log resolution of $(X,\Delta)$.
        Then, we can express $K_Y$ as 
        \[
        K_Y=\pi^*(K_X+\Delta)+\sum_{i=1}^ra_i E_i,
        \]
        where $\pi_*(\sum_{i=1}^ra_iE_i)=-\Delta$.
    We call $a_i+1$ the {\em log discrepancy} of $(X,\Delta)$ with respect to $E_i$ and write it as $A_{X,\Delta}(E_i)$. If $\Delta=0$, we write it simply as $A_{X}(E_i)$. Note that this value is independent of the choice of $\pi$.
\begin{enumerate}
    \item If $a_i\ge-1$ for all $i$, we say that $(X,\Delta)$ is {\em log canonical} ({\em lc} for short).
    If $a_i>-1$ for all $i$, we say that $(X,\Delta)$ is {\em Kawamata log terminal} ({\em klt} for short).
    \item Suppose that $(X,\Delta)$ is lc.
    If $a_i\ge0$ (resp.~$a_i>0$) for all $i$ such that $E_i$ is $\pi$-exceptional, we say that $(X,\Delta)$ is {\em canonical} (resp.~ {\em terminal}).
\end{enumerate}
We note that the above notions are independent of the choice of $\pi$. For more details, see \cite[Section 2]{KoMo}.
We say that $X$ has only lc (resp.~klt, canonical or terminal) singularities if so does $(X,0)$ as a pair.
    \end{defn}

    \begin{exam}
    As discussed in \cite[Proposition 5.20]{KoMo}, quotient singularities are typical examples of klt singularities.
    In the dimension two case, it is known by \cite{kawamata} that every klt singularity is a quotient singularity.
    Two-dimensional canonical singularities are called {\it Du Val singularities} and see \cite[Section 4.2]{KoMo} for more details.

    On the other hand, cyclic quotient singularities are important in this paper.
        Let $p\in X$ be a $2$-dimensional cyclic quotient singularity.
    It is well known that this singularity is analytically isomorphic to the singularity of $\mathbb{A}^2/(\mathbb{Z}/n\mathbb{Z})$, where $\zeta$ is a primitive $n$-th root and the action of $\mathbb{Z}/n\mathbb{Z}$ is given by
    \[
    \mathbb{Z}/n\mathbb{Z}\cong\left\langle\begin{pmatrix}
\zeta & 0 \\
0 & \zeta^a \\
\end{pmatrix}\right\rangle\subset GL(2,\mathbb{C})
    \]
    for some $n,a\in\mathbb{Z}$ such that $a$ is coprime to $n$.
    We say that $p\in X$ is {\em of type $\frac{1}{n}(1,a)$}. 
    \end{exam}

For fundamental notions of the minimal model program (MMP), see \cite{KoMo, fujino-foundation}.

\begin{defn}[Semi-log canonical and stable pairs]
    Let $Y$ be an equidimensional reduced scheme.
    If $Y$ satisfies Serre's $S_2$-condition and $Y$ has only nodes as singularities of codimension one, we say that $Y$ is {\it deminormal}.
    
    Let $\pi\colon X\to Y$ be the normalization of a deminormal scheme.
    Then, we set the conductor divisor $D_X$ on $X$ whose ideal sheaf is defined as $\mathcal{H}om_{\mathcal{O}_Y}(\pi_*\mathcal{O}_X,\mathcal{O}_Y)\subset \pi_{*}\mathcal{O}_X$.
    Let $\Delta$ be an effective $\mathbb{Q}$-divisor on $Y$ whose support is not contained in the non-normal locus of $Y$.
    We say that $(Y,\Delta)$ is {\em semi-log canonical} ({\em slc} for short) if $K_Y+\Delta$ is $\mathbb{Q}$-Cartier and $(X,D_X+\pi^{-1}_*\Delta)$ is lc.
Here, we note that 
$$
K_X+D_X+\pi^{-1}_*\Delta=\pi^*(K_Y+\Delta).
$$
Let $E$ be a prime divisor over $X$ such that $A_{X,D_X+\pi^{-1}_*\Delta}(E)=0$.
Then, we say that the image of $E$ to $Y$ is an {\it slc center} of $(Y,\Delta)$. If further $Y$ is normal, we call the slc center defined by $E$ an {\it lc center}. See \cite[Section 5]{kollar-mmp} for more details.

We say that $(Y,\Delta)$ is {\em stable}
if $(Y,\Delta)$ is slc, $Y$ is proper and $K_Y+\Delta$ is ample.
We say that $Y$ is {\em slc} (resp.~{\em stable}) if so is $(Y,0)$.

On the other hand, if there exists a projective birational morphism $f\colon Y_1\to Y_2$ of two slc pairs $(Y_1,\Delta_1)$ and $(Y_2,\Delta_2)$ such that $K_{Y_1}+\Delta_1=f^*(K_{Y_2}+\Delta_2)$, then we say that $(Y_1,\Delta_1)$ and $(Y_2,\Delta_2)$ are {\it log crepant}.
\end{defn}

The rest of the subsection devotes to introduce the useful notions of divisors that are needed to apply vanishing theorems in \cite{fujino--slc--vanishing} and \cite{Enokizono}.

\begin{defn}[Log big $\mathbb{Q}$-Cartier $\mathbb{Q}$-divisor]
    Let $(Y,\Delta)$ be a slc pair that is projective over a quasi-projective scheme $S$, and let $H$ be a $\mathbb{Q}$-Cartier $\mathbb{Q}$-divisor such that any irreducible component of $\mathrm{Supp}(H)$ is not contained in $\mathrm{Sing}(X)$.  
    If $H|_W$ is relatively big over $S$ for any slc center $W$ and for any irreducible component $W$ of $X$, we say that $H$ is {\it relatively log big over $S$}.
    If $S=\mathrm{Spec}\,\mathbb{C}$, then we simply say that $H$ is {\em log big}.
\end{defn}

We explain the following notion invented by the second author.

\begin{defn}[$\Z$-positive divisor]\label{defn-z-positive}
Let $X$ be a normal proper surface.
We say that a non-zero effective divisor $E$ is {\em negative definite} if $C^2<0$ for any Weil divisor $C$ such that $\emptyset\ne\mathrm{Supp}(C)\subset \mathrm{Supp}(E)$.
A pseudo-effective Weil divisor $D$ is {\em $\Z$-positive} if for any effective negative definite divisor $E$, there exists an irreducible component $C$ of $E$ such that $(D-E)\cdot C>0$.

For any pseudo-effective Weil divisor $D$, we say that $D=P_{\mathbb{Z}}+N_{\mathbb{Z}}$ is {\it the $\mathbb{Z}$-Zariski decomposition} if the following holds:
\begin{enumerate}
    \item $P_{\mathbb{Z}}$ is a $\mathbb{Z}$-positive divisor,
    \item $N_{\mathbb{Z}}=0$ or $N_{\mathbb{Z}}$ is a non-zero effective negative definite divisor, and
    \item $P_{\mathbb{Z}}\cdot C\le 0$ for each irreducible component $C$ of $N_{\Z}$.
\end{enumerate}

From \cite[Theorem 3.5 and Lemma 3.6]{Enokizono}, there uniquely exists a $\Z$-Zariski decomposition of any pseudo-effective divisor $D$.
We say that $N_{\mathbb{Z}}$ is the {\it $\mathbb{Z}$-negative part} of $D$ and $P_{\mathbb{Z}}$ is the {\it $\mathbb{Z}$-positive part} of $D$.
\end{defn}

The above notion is important for applying \cite[Theorems 1.3, 1.4]{Enokizono}.

\begin{defn}[Chain-connected divisor]
    Let $X$ be a normal proper surface.
    A non-zero effective divisor $D$ on $X$ is {\em chain-connected} if for any decomposition $D=A+B$ with $A$ and $B$ non-zero effective, there exists an irreducible component $C$ of $B$ such that $A\cdot C>0$.
    For any non-zero effective divisor $D$ with connected support, there uniquely exists a chain-connected subdivisor $D_{c}\le D$ such that $\Supp(D_c)=\Supp(D)$ and $D_{c}\cdot C\le 0$ for each irreducible component $C$ of $D-D_c$ (\cite[Proposition~1.5~(3)]{Kon}).
    We call $D_c$ the {\em chain-connected component} of $D$.
    If $D$ is a non-zero effective divisor which may not have connected support, then we decompose $D=\sum_{\lambda}D_{\lambda}$ into the connected components and set $D_{c}:=\sum_{\lambda}D_{\lambda,c}$.
    We refer to the decomposition 
    $$
    D=D_c+(D-D_c)
    $$
    as the {\em chain-connected component decomposition} of $D$.
\end{defn}

The following lemma is useful to detect the $\Z$-positity of divisors:

\begin{lem} \label{lem:Z-poschar}
    Let $D$ be a divisor on a normal proper surface $X$.
    Then $D$ is $\Z$-positive if and only if for any negative definite and chain-connected divisor $E>0$ on $X$, there exists an irreducible component $C$ of $E$ such that $(D-E)\cdot C>0$.
\end{lem}

\begin{proof}
    The only if part is trivial.
    Let $D$ be a divisor on $X$ satisfying the assumption of the if part.
    Let $E>0$ be a negative definite effective divisor.
    It suffces to show that there exists an irreducible component $C$ of $E$ such that $(D-E)\cdot C>0$.
    Take the chain-connected decomposition $E=E_{c}+E_1$, where $E_c$ is a disjoint union of chain-connected divisors with $\Supp(E)=\Supp(E_{c})$ and $E_1\ge 0$ satisfies that $E_c\cdot C\le 0$ for any irreducible component $C$ of $E_1$.
    By assumption, there exists a curve $C\subset E_c$ such that $(D-E_c)\cdot C>0$.
    If $(D-E)\cdot C>0$, then the claim holds.
    Assume that $(D-E)\cdot C\le 0$.
    Then $E_1>0$.
    We take the chain-connected component decomposition $E_1=E_{1,c}+E_2$.
    By assumption, there exists a curve $C_1\subset E_{1,c}$ such that $(D-E_{1,c})\cdot C_1>0$.
    Then 
    $$
    (D-E_c-E_{1,c})\cdot C_1\ge (D-E_{1,c})\cdot C_1>0.
    $$
    If $(D-E)\cdot C_1> 0$, then the claim holds.
    If $(D-E)\cdot C_1\le 0$, then $E_2>0$ holds.
    Repeating this process, we can find a curve $C\subset E$ such that $(D-E)\cdot C>0$.
\end{proof}

\begin{rem}\label{rem:Z-poschar}
    Let $S$ be a negative definite reduced curve on a normal proper surface $X$.
    By \cite[Lemma A.12]{Enokizono}, there exists an effective divisor $Z$ such that $\mathrm{Supp}(Z)=S$ and $Z\cdot E<0$ for any irreducible curve $E\subset S$.
    By the same discussion of \cite[Lemma 7.2.2]{Ishii}, we can choose the smallest $Z$ such that $\mathrm{Supp}(Z)=S$ and $Z\cdot E<0$ for any irreducible curve $E\subset S$. 
    Then we can see that any chain-connected divisor $C$ with $\Supp(C)=S$ is a subdivisor of $Z$ as the same discussion of \cite[Proposition~7.2.4]{Ishii}. 
    In particular, the number of such divisors $C$ is finite.
\end{rem}

\subsection{$\mathbb{Q}$-Gorenstein deformation}\label{subsec--q--goren}
In this subsection, we summarize the fundamental concepts of $\mathbb{Q}$-Gorenstein deformations over schemes.

First, we consider $\Q$-Gorenstein deformations over reduced schemes.
Let $S$ be a reduced scheme over $\mathbb{C}$, and let $f\colon \mathscr{X}\to S$ be a flat morphism of schemes with geometrically connected and deminormal fibers of dimension $n$.
We say that a closed subscheme $D\subset \mathscr{X}$ is an {\em effective relative Mumford divisor} if the following hold:
\begin{enumerate}
    \item There exists an open subset $U\subset \mathscr{X}$ such that $\mathrm{codim}_{\mathscr{X}_{\bar{s}}}(\mathscr{X}_{\bar{s}}\setminus U_{\bar{s}})\ge2$ for any geometric point $\bar{s}\in S$,
    \item $U_{\bar{s}}$ contains all generic points of $D_{\bar{s}}$ for any geometric point $\bar{s}\in S$,
    and
    \item $D|_U$ is a Cartier divisor of $U$ and is flat over $S$.
\end{enumerate}
A {\em relative Mumford divisor} (resp.\ {\em relative Mumford $\mathbb{Q}$-divisor}) is a $\Z$-linear (resp.\ $\Q$-linear) combination of 
reduced, irreducible effective relative Mumford divisors, written as $\Delta=\sum_{i=1}^na_iD_i$.
For any morphism $\varphi\colon T\to S$ and any relative Mumford $\mathbb{Q}$-divisor $\Delta$ on $\mathscr{X}$, we define the pullback $\Delta_T$ of $\Delta$ to $\mathscr{X}_T:=\mathscr{X}\times_ST$ as follows:
Since $\Delta|_{U}$ is a $\mathbb{Q}$-Cartier divisor, we can define $\varphi_T^*(\Delta|_U)$, where $\varphi_T\colon U\times_ST\to U$.
Then, we take $\Delta_T$ as the closure of $\varphi_T^*(\Delta|_U)$.
If $T=\mathrm{Spec}(\kappa(s))$ (resp.~$T=\mathrm{Spec}(\overline{\kappa(s)})$) for some $s\in S$, then we write $\Delta_T$ as $\Delta_s$ (resp.~$\Delta_{\bar{s}}$). 

Let $\omega_{\mathscr{X}/S}$ be the relative dualizing sheaf.
We can choose $U$ such that it satisfies the above conditions for each $D_{i}$ and $\omega_{\mathscr{X}/S}|_U$ is a line bundle.
Let $\iota\colon U\hookrightarrow \mathscr{X}$ denote the canonical open immersion.
The morphism $f\colon (\mathscr{X},\Delta)\to S$ is said to be {\em log $\mathbb{Q}$-Gorenstein} if there exists $m\in\mathbb{Z}_{>0}$ such that $\iota_*(\omega_{U/S}^{\otimes m}(m\Delta|_U))$ is a line bundle. 
In this case, we define 
$$
\mathcal{O}_\mathscr{X}(mK_{\mathscr{X}/S}+m\Delta):=\iota_*(\omega_{U/S}^{\otimes m}(m\Delta|_U)).
$$
If $\Delta=0$, $S$ is irreducible and general fibers of $f$ are smooth, we say that $\mathscr{X}_s$ or its singularities are {\em $\mathbb{Q}$-Gorenstein smoothable} for any $s\in S$.
For Gorenstein singularities, the $\mathbb{Q}$-Gorenstein smoothability is equivalent to the smoothability by \cite[Theorem 9.1.6]{Ishii}.
On the other hand, fix a $\mathbb{C}$-valued point $s\in S$ and a point $x\in \mathscr{X}_s$.
If there exists an open neighborhood $U\subset \mathscr{X}$ of $x$ such that $U$ is \'etale equivalent to $U_s\times S$ locally around $x$, then we say that $f$ is {\em locally trivial} around $x$.
If $f$ is locally trivial around any closed point of $X$, then we say that $f$ is {\em locally trivial}.

The important special case occurs when $S$ is a smooth curve.
Fix a closed point $s\in S$ and a $\mathbb{Q}$-Gorenstein family $f\colon \mathscr{X}\to S$.
Let $X:=\mathscr{X}_s$.
We note in this case that $K_\mathscr{X}$ is $\mathbb{Q}$-Cartier.
The morphism $f$ is then called a {\em partial $\mathbb{Q}$-Gorenstein smoothing} of $X$.
We note that if $\mathscr{X}_s$ is slc, then there is an open neighborhood $s\in U$ such that the pair $(f^{-1}(U),X)$ is slc by \cite[2.3]{kollar-moduli}.
Furthermore, if $\mathscr{X}$ is a $\mathbb{Q}$-Gorenstein smoothing, we see that then $f^{-1}(U)$ is klt.

\begin{defn}[T-singularity and Milnor number] \label{def--T-singularity}
Let $p\in X$ be a germ of a normal surface singularity.  

\begin{itemize}
    \item[$(1)$]
    We say that $p\in X$ is a {\em T-singularity} if it is klt and $\Q$-Gorenstein smoothable.
    \item[$(2)$]
    Let $\mathscr{X}\to C$ be a smoothing of $p\in X$.
    We call the second Betti number of the Milnor fiber of $\mathscr{X}\to C$ around $p$ the {\it Milnor number} of $p\in X$ with respect to $\mathscr{X}\to C$ (see \cite[Section 3]{Manetti} for more details).
    Since there is no local obstruction for $\mathbb{Q}$-Gorenstein deformations of T-singularities (see the proof of \cite[Theorem 8.2]{Hac}), we see that the Milnor number of a T-singularity $p\in X$ does not depend on the choice of $\mathbb{Q}$-Gorenstein smoothing $\mathscr{X}\to C$.
    Therefore, we let $\mu_p$ denote the Milnor fiber of a T-singularity $p\in X$ with respect to every $\mathbb{Q}$-Gorenstein smoothing.
    As \cite[Lemma 2.4]{HP}, for a Du Val singularity $p\in X$ of type $A_r$, $D_r$ or $E_r$, then $\mu_p=r$.
    \item[$(3)$]
    We say that $p\in X$ is a {\em Wahl singularity} if it is a T-singularity with the Milnor number zero.    
\end{itemize}

\end{defn}

\begin{defn}[Universal hull]
Let $f\colon \mathscr{X}\to S$ be a flat morphism of schemes of finite type over $\mathbb{C}$ such that each geometric fiber is a deminormal scheme of dimension $d$.
Let $n$ be an integer, and let $\mathscr{F}$ be a coherent sheaf on $\mathscr{X}$ that is flat over $S$.
We say that $\mathscr{F}$ is the {\it universal hull} of $\omega_{\mathscr{X}/ S}^{\otimes n}$, denoted by $\omega_{\mathscr{X}/ S}^{[n]}$, if the following holds:
\begin{itemize}
    \item for any $s\in S$, $\varphi_s$ is isomorphic to the natural morphism $\omega_{\mathscr{X}_s}^{\otimes n}\to \omega_{\mathscr{X}_s}^{[n]}$, where $\omega_{\mathscr{X}_s}^{[n]}$ is the $S_2$-hull of $\omega_{\mathscr{X}_s}^{\otimes n}$.
\end{itemize}
Note that if a universal hull exists, then the isomorphic class of $\mathscr{F}$ is uniquely determined.
See \cite[\S9]{kollar-moduli} for more details.
\end{defn}
 
Next, we collect fundamental terminologies of the $\mathbb{Q}$-Gorenstein deformation theory of surfaces over not necessarily reduced schemes.
Let $X$ be an slc surface with a closed point $p\in X$ and $N$ be the minimal positive integer such that $\omega_{X}^{[N]}$ is a line bundle around $p$.
We call this $N$ the {\it canonical index} of the singularity $p\in X$.
Take an open neighborhood $U$ of $p$ such that $\omega_{U}^{[N]}$ is a line bundle.
Then, we define the canonical covering of $\pi\colon Z\to U$ by
$$
Z=\mathrm{Spec}_U(\mathcal{O}_U\oplus\mathcal{O}_U(K_U)\oplus\ldots\oplus\mathcal{O}_U((N-1)K_U)).
$$
Such $Z$ is uniquely determined \'etale locally over $X$ and we can define the {\em canonical covering stack} $\mathfrak{X}\to X$ by gluing the data of the canonical coverings \'etale locally.
We note that $\mathfrak{X}$ is a Deligne-Mumford stack with the coarse moduli scheme $X$ (cf.~\cite[\S3]{Hac}).
Hacking defined the following notion in \cite{Hac} (see \cite[\S3.2]{Hac}).
\begin{defn}
    Let $f\colon\mathscr{X}\to S$ be a flat family of slc surfaces over a scheme $S$ with a point $s\in S$ such that $\mathscr{X}_s=X$.
    We say that $f$ is {\em $\mathbb{Q}$-Gorenstein} if there exists $g\colon\widetilde{\mathfrak{X}}\to S$ a flat family of Deligne-Mumford stack with the coarse moduli space $\mathscr{X}$ such that $\widetilde{\mathfrak{X}}_s=\mathfrak{X}$, the canonical covering stack of $X$.
    When $S$ is reduced, $f$ is $\Q$-Gorenstein if and only if $f\colon(\mathscr{X}, 0)\to S$ is log $\Q$-Gorenstein.
\end{defn}

Under the situation as above, suppose that $S=\mathrm{Spec}\,A$, where $A$ is a local Artinian $\mathbb{C}$-algebra and $s$ is the unique closed point of $S$. Let $L_{\widetilde{\mathfrak{X}}/A}$ denote the cotangent complex of $\widetilde{\mathfrak{X}}$ over $A$.
For the definition of cotangent complex, see \cite{I1} and \cite{I2}.
Suppose that $p\colon \widetilde{\mathfrak{X}}\to \mathscr{X}$ the canonical morphism.
For any finite $A$-module $M$, we set 
\begin{align*}
    T^i_{QG}(\mathscr{X}/A,M)&:=\mathrm{Ext}^i(L_{\widetilde{\mathfrak{X}}/A},\mathcal{O}_{\widetilde{\mathfrak{X}}}\otimes_AM),\\
    \mathscr{T}^i_{QG}(\mathscr{X}/A,M)&:=p_*\mathcal{E}xt^i(L_{\widetilde{\mathfrak{X}}/A},\mathcal{O}_{\widetilde{\mathfrak{X}}}\otimes_AM).
\end{align*}
By \cite[Theorem 3.9]{Hac}, $T^2_{QG}(\mathscr{X}/A,M)$ is the space of the obstructions.
Note that there is the following natural spectral sequence
\begin{equation}\label{eq--T^i--spectral}
E_2^{i,j}:=H^i(\mathscr{X},\mathscr{T}^j_{QG}(\mathscr{X}/A,M))\Rightarrow {T}^{i+j}_{QG}(\mathscr{X}/A,M).
\end{equation}
It follows from \cite[Lemma 3.8]{Hac} that $\mathscr{T}^0_{QG}(X/\mathbb{C},\mathbb{C})=T_X:=\mathcal{H}om_{\mathcal{O}_X}(\Omega_X,\mathcal{O}_X)$.
For any slc surface $X$, we write $\mathscr{T}^i_{QG}(X)$ for $\mathscr{T}^i_{QG}(X/\mathbb{C},\mathbb{C})$.
If there is no fear of confusion, we abbreviate it to $\mathscr{T}^i_{QG}$.

The following theorem implies that slc surfaces with big anticanonical divisors have no local-to-global obstractions:

\begin{thm}[{\cite[Theorem 8.2]{Hac}, \cite[Proposition 3.1]{HP}}]\label{thm--hacking--prokhorov}
    Let $W$ be an irreducible projective surface with only slc singularities such that $-K_W$ is big.
    Then $H^2(W,T_W)=0$.
    In particular, if $W$ is normal and has only $\mathbb{Q}$-Gorenstein smoothable singularities, then $W$ admits a global $\mathbb{Q}$-Gorenstein smoothing.
    Furthermore, if $W$ has only T-singularities, then $W$ has an unobstructed $\mathbb{Q}$-Gorenstein deformation.
\end{thm}

\begin{proof}
If $W$ is normal, all the assertions are dealt with in \cite[Theorem 8.2]{Hac} and \cite[Proposition 3.1]{HP}.
Therefore, it suffices to deal with the case when $W$ is not normal.
    Let $\nu\colon Y\to W$ be the normalization and $\Delta$ the conductor on $Y$.
    Let $D$ be the non-normal locus on $W$ with the reduced structure.
    We have the following exact sequence 
    \[
    0\to \nu_*(\mathcal{O}_{Y}(-\Delta))\to \mathcal{O}_W\to \mathcal{O}_{D}\to 0.
    \]
    By this, we have the following exact sequence
    \[
    0\to \nu_*T_Y(-\Delta)\to T_W\to \mathscr{F}\to0,
    \]
    where $\mathscr{F}$ is the cokernel sheaf supported on $D$.
    Thus, $H^2(Y,T_Y(-\Delta))\to H^2(W,T_W)$ is surjective.
    It suffices to show that $H^2(Y,T_Y(-\Delta))=0$.
    By the Serre duality, the dual space of $H^2(Y,T_Y(-\Delta))$ is isomorphic to 
    \[
    \mathrm{Hom}(\mathcal{O}_Y(-K_Y-\Delta),\Omega_Y^{\vee\vee}).
    \]
    By the Bogomolov-Sommese vanishing theorem \cite[Theorem 1.4]{GKK}, \[\mathrm{Hom}(\mathcal{O}_Y(-K_Y-\Delta),\Omega_Y^{\vee\vee})=0\] since $-K_Y-\Delta=-\nu^*K_W$ is big.
    Thus, we have the assertions. (See also \cite[Notation 3.6 and Lemma 3.8]{Hac}).
\end{proof}

\subsection{KSBA moduli}
For definitions and fundamental notions of Deligne-Mumford stacks, see \cite{Ols} and \cite{Stacks} for example.
Fix $v\in\mathbb{Q}_{>0}$, $c\in \mathbb{Q}\cap[0,1)$ and $d\in\mathbb{Z}_{>0}$.
Then, we define the following quasi-functor $\mathcal{M}^{\mathrm{KSBA}}_{d,v,c}$ to the category of groupoids such that for any scheme $S$, $\mathcal{M}^{\mathrm{KSBA}}_{d,v,c}(S)$ is the groupoid whose objects are families $f\colon (\mathscr{X},\Delta)\to S$ such that 
\begin{enumerate}
    \item $f\colon \mathscr{X}\to S$ is a projective flat morphism of relative dimension $d$,
    \item $\Delta=c\mathscr{D}$, where $\mathscr{D}$ is a family of K-flat divisors (cf.~\cite[Definition 7.1]{kollar-moduli}),
    \item if $m$ is the denominator of $c$, then $\omega_{\mathscr{X}/S}^{[m]}(m\Delta)$ is a flat family of divisorial sheaves,
    \item for any geometric point $\bar{s}\in S$, $(\mathscr{X}_s,\Delta_s)$ is slc and $\mathrm{vol}(K_{\mathscr{X}_s}+\Delta_s)=v$, and 
    \item $K_{\mathscr{X}/S}+\Delta$ is an $f$-ample $\mathbb{Q}$-line bundle,
\end{enumerate}
where an isomorphism $\alpha$ from $f\colon (\mathscr{X},\Delta_{\mathscr{X}})\to S$ to $g\colon (\mathscr{Y},\Delta_{\mathscr{Y}})\to S$ is defined as an isomorphism $\alpha\colon \mathscr{X}\to\mathscr{Y}$ such that $\Delta_{\mathscr{X}}=\alpha^*\Delta_{\mathscr{Y}}$.
\begin{thm}[{\cite{Al,hmx-boundgentype,fujino-semi-positivity,KP,kollar-moduli}}]\label{thm-ksba}
    $\mathcal{M}^{\mathrm{KSBA}}_{d,v,c}$ is a Deligne-Mumford stack proper over $\mathbb{C}$ with the projective coarse moduli scheme $\pi\colon\mathcal{M}^{\mathrm{KSBA}}_{d,v,c}\to M^{\mathrm{KSBA}}_{d,v,c}$. 
\end{thm}

$M^{\mathrm{KSBA}}_{2,v,0}$ is a projective scheme containing the open locus $M^{\mathrm{Gie}}_{v}$, which parametrizes normal surfaces with only canonical singularities (see \cite[Corollary 4.10]{kollar-mmp}).  
Gieseker \cite{gieseker} first constructed the moduli space $M^{\mathrm{Gie}}_{v}$ by applying the geometric invariant theory.
Since the geometric genus is an invariant for $\mathbb{Q}$-Gorenstein deformations (cf.~\cite[Corollary 2.69]{kollar-moduli}), we can take a closed and open subscheme $M^{\mathrm{Gie}}_{v,p}$ of $M^{\mathrm{Gie}}_{v}$ that parametrizes $X$ with $p_g(X)=p$ for any $p\in\mathbb{Z}_{\ge0}$. 
On the other hand, let $\mathcal{M}^{\mathrm{KSBA,normal}}_{2,v,0}$ denote the open substack of $\mathcal{M}^{\mathrm{KSBA}}_{2,v,0}$ parametrizing normal varieties.

\subsection{Log canonical surface singularities}\label{subsec:lc_sing}
In this section, we collect basic facts on $\Q$-Gorenstein smoothable log canonical surface singularities.

Let $p\in X$ be a germ of a normal surface
and $\pi\colon \widetilde{X}\to X$ the minimal resolution.
If $K_X$ is $\mathbb{Q}$-Cartier, then the canonical divisor is written by 
$$
K_{\widetilde{X}}=\pi^{*}K_{X}+K_{p},
$$
where $K_{p}=\sum_{i}(\alpha_i-1)E_i$ is a $\Q$-linear combination of exceptional divisors, which is called the {\em canonical cycle} of $p\in X$,
and $\alpha_i=A_X(E_i)$.
Note that $-K_p$ is effective since $\pi$ is minimal. 
As noted in \cite[Section 4.1]{KoMo}, we can define the canonical cycle $K_p$ even if $K_X$ is not $\mathbb{Q}$-Cartier.
Moreover, the condition that all $\alpha_i$'s are non-negative implies that $K_X$ is $\mathbb{Q}$-Cartier.
Assume that $p\in X$ is lc, that is, $\alpha_i\ge 0$ for any $i$.
Let $E=\sum_{i}E_{i}$ denote the reduced $\pi$-exceptional divisor.

If a surface singularity $p\in X$ is a T-singularity but not a Du Val singularity, then this is a cyclic singularity of the form $$\frac{1}{dn^2}(1,adn-1),$$
where $d,n$ and $a\in\mathbb{Z}_{>0}$ such that $1\leq a \leq n-1$ and $a$ is coprime to $n$ by \cite[Proposition 3.10]{KSB}.
 By \cite[Lemma 2.4]{HP}, $p\in X$ has Milnor number $\mu_p=d-1$ in this case.
It is known that the $\pi$-exceptional curves of a non-Du Val T-singularity form a chain $E_{1}+E_{2}+\cdots+E_{r}$ of smooth rational curves.
The self-intersection numbers $-b_i:=E_{i}^{2}$ satisfy certain arithmetic conditions ({\cite[Proposition~3.11]{KSB}}, see also the next section).
This chain is denoted by $[b_1,b_2,\ldots,b_r]$ and called a {\em T-chain}.
It is easy to see that $K_{p}^{2}+K_{p}\cdot E=1$ for a T-singularity $p\in X$.

If $p\in X$ is lc but not klt, which is called {\em strictly lc}, then the $\pi$-exceptional divisors are classified by {\cite[Chapter~3]{Koetal}}.
We first deal with non-rational singularities.
Here, we say that $p\in X$ is a {\it rational singularity} if $p\in X$ is normal and for any resolution of singularities $\pi\colon \widetilde{X}\to X$, $R^i\pi_*\mathcal{O}_{\widetilde{X}}=0$ for any $i>0$.
Note that any klt singularity is rational by \cite[Theorem 5.22]{KoMo}.
It is well known that a strictly lc surface singularity $p\in X$ 
is not rational if and only if it is Gorenstein, that is, it is a {\em simple elliptic} singularity or a {\em cusp} singularity (cf.~\cite[\S7.6]{Ishii}, \cite[\S4.4]{KoMo}).
We call these singularities {\it elliptic}.
In this case, $K_{p}=-E$ and hence $K_{p}^{2}+K_{p}\cdot E=0$.

We use the following notation for elliptic singularities:
\begin{itemize}
     \item  A simple elliptic singularity whose exceptional curve $E$ is a smooth elliptic curve with self-intersection number $-b$.
     We call such a singularity a {\em simple elliptic singularity of type} $\operatorname{Ell}_b$.

     \item A cusp singularity whose exceptional divisor $E$ consists of smooth rational curves $C_1, \ldots, C_r$  arranged in a cycle, such that each $C_i$ intersects transversally with $C_{i-1}$ and $C_{i+1}$ (with indices taken modulo $r$, so $C_r$ meets $C_1$). Let $-b_i$ be the self-intersection number of $C_i$. 
     We call such a singularity a {\em cusp singularity of type} $[b_1, \ldots, b_r]^{\circ}$ (resp. {\em cusp singularity of type} $[b_1+2]^{\circ}$) if $r\geq2$ (resp. $r=1$).
     For cusp singularities, we follow the same rule as introduced in Notation~\ref{note--T-chain} as well.
\end{itemize}
Let $X$ be a proper normal surface with an elliptic singularity $p\in X$, and let $\pi\colon \widetilde{X}\to X$ be a minimal resolution at the singularity $p$.
Then the Leray spectral sequence implies that $\chi(\widetilde{X},\mathcal{O}_{\widetilde{X}})=\chi(X,\mathcal{O}_X)-1$.

We review the $\Q$-Gorenstein smoothability for simple elliptic and strictly lc rational singularities. For cusp singularities, see Appendix~\ref{app:cusp}.

\begin{lem}[{\cite[\S5]{Pink}}]\label{lem:smoothable_elliptic_singularity}
Let $p\in X$ be a simple elliptic singularity of type $\operatorname{Ell}_b$.
Then, $p\in X$ is smoothable
if and only if $b \le 9$.
\end{lem}

\begin{lem-defn}[{\cite[Theorem~1.2]{Pro}}] \label{lem:smoothable_rational_strict_lc}
Let $p\in X$ be a strictly lc rational surface singularity which admits a $\Q$-Gorenstein smoothable.
Then the dual graph of the $\pi$-exceptional curves is one of the following:

\smallskip

\noindent
(i) type $(2,2,2,2)[b_1,\ldots,b_n]$, where $b_i\ge 2$ and $\sum_{i=1}^{n}(b_i-3)\le 3$:
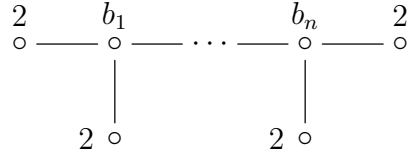
\[
\xygraph{
    \circ ([]!{-(0,-.3)} {2}) - [r]
    \circ ([]!{-(0,-.3)} {b_1}) (
        - [d] \circ ([]!{-(.3,0)} {2}),
        - [r] \cdots ([]!{-(-0,-.3)} {})
         - [r] \circ ([]!{-(0,-.3)} {b_n}) (
         - [d] \circ ([]!{-(.3,0)} {2}),
        - [r] \circ ([]!{-(0,-.3)} {2})
)}
\]
We call such a singularity $p \in X$ a \emph{strictly lc singularity of type} $(2,2,2,2)[b_1,\dots,b_r]$, or a \emph{singularity of type} $(2,2,2,2)$.
\smallskip

\noindent
(ii) type $(3,3,3)[b]$, where $b=2,3,4$:
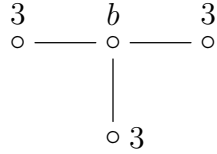
\[
\xygraph{
    \circ ([]!{-(0,-.3)} {3}) 
    - [r]  \circ ([]!{-(0,-.3)} {b})(
        - [r] \circ ([]!{-(0,-.3)} {3}),
        - [d] \circ ([]!{-(-0.25,-.0)} {3}),
}
\]
We call such a singularity $p \in X$ a \emph{strictly lc singularity of type} $(3,3,3)[b]$, or a \emph{singularity of type} $(3,3,3)$.

\smallskip

\noindent
(iii) type $(2,4,4)[b]$, where $b=2,3$:
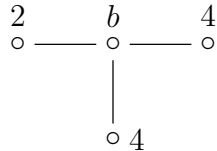
\[
\xygraph{
    \circ ([]!{-(0,-.3)} {2}) 
    - [r]  \circ ([]!{-(0,-.3)} {b})(
        - [r] \circ ([]!{-(0,-.3)} {4}),
        - [d] \circ ([]!{-(-0.25,-.0)} {4}),
}
\]
We call such a singularity $p \in X$ a \emph{strictly lc singularity of type} $(2,4,4)[b]$, or a \emph{singularity of type} $(2,4,4)$.

\smallskip

\noindent
(iv) type $(2,3,6)[2]$:
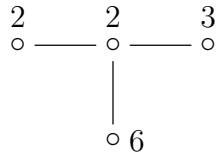
\[
\xygraph{
    \circ ([]!{-(0,-.3)} {2}) 
    - [r]  \circ ([]!{-(0,-.3)} {2})(
        - [r] \circ ([]!{-(0,-.3)} {3}),
        - [d] \circ ([]!{-(-0.25,-.0)} {6}),
}
\]
We call such a singularity $p \in X$ a \emph{strictly lc singularity of type} $(2,3,6)[2]$, or a \emph{singularity of type} $(2,3,6)$.

\smallskip

\noindent
Furthermore, $K_{p}^{2}+K_{p}\cdot E=0$ for the case of (i)
and $K_{p}^{2}+K_{p}\cdot E=1$ otherwise.
\end{lem-defn}

For non-normal surface singularity, $\Q$-Gorenstein smoothability is more complicated. For the slt case, that is, the case where the normalization is plt, $\Q$-Gorenstein smoothability is completely classified by Koll\'ar-Shepherd-Barron.

\begin{thm}[{\cite[Theorems 4.23, 5.2]{KSB}}]
Let $p\in X$ be a non-normal slt singularity that is $\Q$-Gorenstein smoothable. Then, $p\in X$ is one of the following:
\begin{itemize}
\item[$(1)$] $p\in X$ is a node,
\item[$(2)$] $p\in X$ is a pinch point, that is, \'etale equivalent to $(x^2=yz^2)\subset \mathbb{A}^3$, and
\item[$(3)$] $p\in X$ is \'etale equivalent to the quotient of $(xy=0)\subset\mathbb{A}^3$ by $\mathbf{\mu}_r$ for some $r\in\mathbb{Z}_{>1}$, where $x,y$ and $z$ are the canonical coordinates of $\mathbb{A}^3$ and $\mathbf{\mu}_r$ acts on $\mathbb{A}^3$ in the way that 
    \[
    \mu_r\cong\left\langle\begin{pmatrix}
\zeta & 0 & 0 \\
0 & \zeta^{-1}&0 \\
0 & 0&\zeta^{a}\\
\end{pmatrix}\right\rangle\subset GL(3,\mathbb{C})
    \]
    for some $a\in\mathbb{Z}$ such that $a$ is coprime to $r$, where $\zeta$ is a primitive $r$-th root.
    We write this type singularity as $(xy=0)\subset \frac{1}{r}(1,-1,a)$.
\end{itemize}
Conversely, the above singularities are $\Q$-Gorenstein smoothable.
\end{thm}

We say that a non-normal Gorenstein surface singularity is a {\it degenerate cusp} if the exceptional divisor of a semi-resolution (cf.~\cite[Definition 4.3]{KSB}) with no $(-1)$-curve contracted is a cycle of rational curves or a nodal rational curve.
We know that all degenerate cusps are slc and smoothable by \cite[Theorem 4.21]{KSB} and \cite[(3.5)]{JSte}.

\subsection{Double covering and even resolution}

Let $\varphi\colon X\to W$ be a double covering between normal varieties.
Let $B$ be the branch divisor of $\varphi$.
Then, it is known that there exists a Weil divisorial sheaf $L$ on $W$ such that $B\in |L^{[2]}|$ and $X\cong \mathbf{Spec}_{W}(\O_{W}\oplus L^{[-1]})$, where the $\O_{W}$-algebra structure on $\O_{W}\oplus L^{[-1]}$ is naturally induced by multiplying the section defining $B$. 
This relation induces a one-to-one correspondence between the set of isomorphism classes of double coverings over $W$ and the set of isomorphism classes of pairs $(B, L)$ where $L$ is a Weil divisorial sheaf on $W$ and $B$ is an effective reduced divisor contained in $|L^{[2]}|$.

Let $\psi\colon W'\to W$ be a proper birational morphism between normal varieties.
For a double covering $X'\to W'$ branched along $B'\in |L'^{[2]}|$, the Stein factorization of the composition $X'\to W'\to W$ induces a double covering $X\to W$ branched along $B\in |L^{[2]}|$, where $B=\psi_{*}B'$ and $L=\psi_{*}L'$.
For a double covering $\varphi\colon X\to W$ branched along $B\in |L^{[2]}|$, we give a double covering $X'\to W'$ from a normal variety $X'$ by taking the normalized base change of $\varphi$ by $\psi$.
If $W$ and $W'$ are smooth surfaces, the branch divisor $B'\in |L'^{[2]}|$ can be described explicitly by the even resolution:

\begin{defn}[Even resolution] \label{defn--evenresol} 
Let $\psi\colon W'\to W$ be a proper birational morphism between smooth surfaces.
We decompose this into a sequence of blow-ups 
$$
\psi=\psi_{1}\circ \cdots \circ\psi_{N}\colon W'=W_{N}\to W_{N-1}\to \cdots \to W_{0}=W,
$$
where $\psi_{i}\colon W_i\to W_{i-1}$ is the blow-up at a point $x_i\in W_{i-1}$ with the exceptional curve $E_i$.
Let $X\to W$ denote the double cover branched along $B\in |2L|$.
Normalized base changes define double covers $X_i\to W_i$.
Then, the branch locus $B_i\in |2L_i|$ on $W_i$  is described inductively as
$$
B_{i}=\psi_i^*B_{i-1}-2\left\lfloor\frac{m_i}{2}\right\rfloor E_i,\quad L_{i}=\psi_i^*L_{i-1}-\left\lfloor\frac{m_i}{2}\right\rfloor E_i,
$$
where $m_i:=\mathrm{mult}_{x_i}B_{i-1}$.
Namely, $B_i$ is the proper transform of $B_{i-1}$ if $m_i$ is even, and is the sum of the proper transform of $B_{i-1}$ and $E_i$ if $m_i$ is odd.
If we start with a singular double cover $X$, then we can take a sequence of blow-ups as above with $m_i\ge 2$ so that $B_N$ is smooth for some $N$ and hence $X_N$ is smooth.
This resolution process is called an {\em even resolution} or a {\em canonical resolution}, which was firstly introduced in \cite{horikawa-quintic}.
\end{defn}

\section{Extended T-chain}\label{sec:extT}  
In this section, we introduce the notion of an extended T-chain. 
It will play an important role in various classifications of certain surfaces (Section \ref{sec:non-std_Horikawa}).

Throughout this paper, we use the following notation.

\begin{note}\label{note--T-chain}
    \phantom{A}
    \begin{itemize}
        \item[$(1)$] The notation $2^n$ ($n\geq0$) in a string represents the number $2$ repeated $n$ times.
        For example, $[3,2^3,3]$ denotes the string $[3,2,2,2,3]$.
        \item[$(2)$] The expression $a,2^{-1},b$ in a string represents a single number $a+b-2$.
        For example, $[3,2^{-1},3]$ denotes the string $[4]$.
    \end{itemize}
\end{note}

We recall the fundamental combinatorial operartions for chains of rational curves:

\begin{defn}[Operations $L_i$ and $R_j$]
Let $\sum_{i=1}^{r}C_{i}$ be a chain of smooth rational curves on a smooth proper surface $Y$ with $b_{i}:=-C_{i}^{2}\ge 2$ and $C_{i}\cdot C_{i+1}=1$, which is denoted as $[b_{1},\ldots,b_{r}]$ for simplicity.

It is well known that the singularity $p\in X$ obtained by contracting all $C_{i}$'s is a cyclic quotient singularity of type $\frac{1}{n}(1,a)$, where $a$ and $n$ are coprime to each other and determined by the negative continued fraction expansion
$$
\frac{n}{a}=b_1-\frac{1}{b_2-\frac{1}{b_3-\cdots}}.
$$

For a chain $[b_1,\ldots,b_r]$, we define
\begin{align*}
&L_2[b_1,\ldots,b_r]:=[2,b_1,\ldots,b_r], \\
&L_1[b_1,\ldots,b_r]:=[b_1+1,b_2,\ldots,b_r], \\ 
&[b_1,\ldots,b_r]R_1:=[b_1,\ldots,b_{r-1}, b_r+1],\\
&[b_1,\ldots,b_r]R_2:=[b_1,\ldots,b_r,2].
\end{align*}
Since the operations $L_{i}$ and $R_{j}$ are commutative,
a description like 
$$
L_{i_p}\cdots L_{i_1}[b_1,\ldots,b_r]R_{j_1}\cdots R_{j_q}
$$
is well-defined.
\end{defn}

According to {\cite[Proposition~3.11]{KSB}}, all T-chains are nothing but of the form
$$
L_{i_p}\cdots L_{i_1}(L_1[2^n]R_1)R_{3-i_1}\cdots R_{3-i_p}. 
$$

\begin{defn}[Log discrepancies]\label{defn--T-sing--log--disc}
Let $[b_1,\ldots,b_r]$ be a chain corresponding to a cyclic quotient singularity $p\in X$.
We define the {\em log discrepancy} $\alpha_i$ of $b_i$ in $[b_1,\ldots,b_r]$ as the log discrepancy of the corresponding exceptional curve $E_i$ with respect to $X$.
\end{defn}

If $[b'_{0},b'_{1},\ldots,b'_{r}]=L_2[b_1,\ldots,b_r]R_1$,
the log discrepancies $\alpha'_{0},\ldots,\alpha'_{r}$ in $[b'_{0},b'_{1},\ldots,b'_{r}]$ are determined by using the log discrepancies $\alpha_1,\ldots,\alpha_{r}$ in $[b_1,\ldots,b_r]$ as 
$$
\alpha'_{0}=\frac{1}{1+\alpha_r}, \alpha'_{1}=\frac{\alpha_1}{1+\alpha_r},\ldots,\alpha'_{r}=\frac{\alpha_r}{1+\alpha_r}.
$$
The case where $[b'_{0},b'_{1},\ldots,b'_{r}]=L_1[b_1,\ldots,b_r]R_2$ is similar.
Hence the log discrepancies in a T-chain can be computed inductively.

\begin{defn}[Core]
Let $[b_1,\ldots,b_r]$ be a T-chain corresponding to a T-singularity $p\in X$ of type $\frac{1}{dn^2}(1,adn-1)$.
A component $b_i$ with the smallest log discrepancy in $[b_1,\ldots,b_r]$ is called a {\em core}.
Note that if we write 
\[
[b_1,\ldots,b_r]=L_{i_p}\cdots L_{i_1}(L_1[2,\ldots,2]R_1)R_{3-i_1}\cdots R_{3-i_p},
\]
the cores in $[b_1,\ldots,b_r]$ are the components corresponding to $[2,\ldots,2]$ and the number of cores is equal to $d$.
\end{defn}

\begin{rem}
The notion of the core of a T-chain coincides with  the center of a T-chain defined in \cite[Definition~2.2]{FRU}. 
In this paper, we also use the term \emph{core} to refer to this notion, in order to avoid confusion with (lc) centers.
\end{rem}

\begin{defn}[Extended T-chains]\label{def:extT}
Let $[b_1,\ldots,b_r]$ be a chain corresponding to a chain of rational curves $E_1,\ldots,E_{r}$ on a surface $Y$.
Then $[b_{1},\ldots,b_{r}]$ is said to be an {\em extended T-chain} if after taking a finite number of blow-ups at nodes $\rho\colon \widetilde{X}\to Y$, its reduced total transform is of the form
$$
\sum_{j=1}^{r_{1}}C_{1,j}+E_{1}+\sum_{j=1}^{r_{2}}C_{2,j}+E_{2}+\cdots+E_{l-1}+\sum_{j=1}^{r_{l}}C_{l,j},
$$
where each $\sum_{j=1}^{r_{i}}C_{i,j}$ is a T-chain, which is denoted by $[b_{i,1},\ldots,b_{i,r_{i}}]$ with $b_{i,j}:=-C_{i,j}^{2}$, and each $E_{i}$ is a $(-1)$-curve such that $C_{i,r_{i}}\cdot E_{i}=C_{i+1,1}\cdot E_{i}=1$.
This chain is simply denoted by
$$
[b_{1,1},\ldots,b_{1,r_{1}}]-1-\cdots-1-[b_{l,1},\ldots,b_{l,r_{l}}]
$$
and said to be a {\em T-train} associated to the extended T-chain $[b_{1},\ldots,b_{r}]$.
We call the chain $[b_{i,1},\ldots,b_{i,r_{i}}]$ the {\em $i$-th vehicle}.
The log discrepancies corresponding to $[b_{i,1},\ldots,b_{i,r_{i}}]$ is denoted by $\alpha_{i,1},\ldots,\alpha_{i,r_i}$.

The T-train satisfies {\em the ample condition} if the sum $\alpha_{i,r_{i}}+\alpha_{i+1,1}$ of the log discrepancies of $C_{i,r_{i}}$ and $C_{i+1,1}$  is less than $1$ for $1\leq i\leq l-1$.
This is equivalent to that $K_{X}\cdot\pi(E_{i})>0$ for each $i$, where $\pi\colon \widetilde{X}\to X$ is the contraction of all $C_{i,j}$.
In this case, any core of any T-chain $\sum_{j=1}^{r_{i}}C_{i,j}$ is not contracted by $\rho$ from {\cite[Corollary~2.4]{FRU}}.
We call a T-train that satisfies the ample condition an \emph{ample T-train}.  
An extended T-chain is said to be \emph{P-admissible} if it admits an associated ample T-train.
\end{defn}

Consider an extended T-chain $C$ and fix an associated T-train $\widetilde{C}$.
Let $r$ denote the length of $C$, and let $l$ denote the number of vehicles of $\widetilde{C}$.
Set $d=\sum_{i=1}^{l}d_{i}$, where $d_i$ denotes the length of the cores in the $i$-th vehicle.
Then similarly to T-chains, the equality
$$
r-d+2=\sum_{k=1}^{r}(b_{k}-2)
$$
holds by simple calculation.

Recall that a {\em P-resolution $\pi\colon X'\to X$} is a partial resolution of $X$ such that $X'$ has only T-singularities and $K_{X'}$ is $\pi$-ample.
Note that if $X$ has only T-singularities, the identity map $X \to X$ is regarded as a P-resolution.
The following is a relation between extended T-chains and P-resolutions.

\begin{prop} \label{charP-resol}
Let $p\in X$ be a cyclic quotient singularity and $E_1+\cdots+E_r$ the corresponding chain of exceptional curves with $E_{i}^{2}=-b_i$.

\begin{itemize}

\item[$(1)$]
A partial resolution $\pi\colon X'\to X$ is a P-resolution such that all $E_i$'s are not $\pi$-exceptional if and only if 
it corresponds to an ample T-train associated to $[b_1,\ldots,b_r]$, that is, $X'$ is obtained by contracting all T-chains in the T-train.

\item[$(2)$]
There exists a one-to-one correspondence between the P-resolutions of $p\in X$ and the tuples of T-trains with the ample condition associated to disjoint sub-chains in $[b_1,\ldots,b_r]$.
\end{itemize}

\end{prop}

\begin{proof}
(1) If part is straightforward by the definition of a P-resolution.
Let $\pi\colon X'\to X$ be a P-resolution such that none of the $E_i$'s are $\pi$-exceptional.
Let $\widetilde{X}\to X$ and $\widetilde{X}'\to X'$ denote the minimal resolutions of $X$ and $X'$, respectively.
Let $X_{m}\to X$ be the maximal resolution as described in {\cite[Lemma~3.13]{KSB}}.
By {\cite[Lemma~3.14]{KSB}}, the map factors as
$X_m\to \widetilde{X}'\to \widetilde{X}\to X$.
In particular, the exceptional set of $\widetilde{X}'\to X$ forms a chain.
This implies that the birational morphism $\widetilde{X}'\to \widetilde{X}$ is a sequence of blow-ups at nodes within the exceptional chain.
Since $K_{X'}$ is ample over $X$ and $\widetilde{X}'\to X'$ contracts all $E_i$'s, the exceptional curves on $\widetilde{X}'$ over $X$ which are not exceptional over $X'$ are disjoint unions of $(-1)$-curves.
Thus, the exceptional curves of $\widetilde{X}'\to X$ form a T-train associated to $[b_1,\ldots,b_r]$.
This establishes the claim (1).

(2) Let $\pi\colon X'\to X$ be a P-resolution.
Let $Y$ denote the surface over $X$ obtained by contracting the $E_i$'s which are not $\pi$-exceptional from the minimal resolution of $X$.
By construction, $\pi$ factors through $Y$, and $X'\to Y$ is a P-resolution of $Y$.
The singularities of $Y$ correspond to disjoint sub-chains within $[b_1,\ldots,b_r]$.
Therefore, the claim (2) follows directly from the claim (1).
\end{proof}

\section{Geography of normal stable surfaces}\label{sec:geography_normal_stable}
In this section, we introduce a useful formula for normal stable surfaces and a classification result that will be needed in later sections.

\subsection{Inequalities of normal stable surfaces}
In this subsection, we show Noether-type inequalities for normal stable surfaces. 
Let $X$ be a normal stable surface.
Let $\pi\colon \widetilde{X}\to X$ denote the minimal resolution of $X$, and let $\widetilde{\Delta}$ denote the reduced inverse image of all elliptic singularities on $X$.
We can run a $(K_{\widetilde{X}}+\widetilde{\Delta})$-minimal model program (cf. {\cite[\S~4.10]{fujino-foundation}})
$$
(\widetilde{X},\widetilde{\Delta})=(Y_{0},\Delta_{0})\xrightarrow{\rho_{1}}(Y_{1},\Delta_{1})\xrightarrow{\rho_{2}}\cdots \xrightarrow{\rho_{m}} (Y_{m},\Delta_{m})=(Y,\Delta),
$$
where $(Y,\Delta)$ is a log minimal model of $(\widetilde{X},\widetilde{\Delta})$.

\begin{lem}
Each contraction $\rho_{i}\colon Y_{i-1}\to Y_{i}$ is a blowing down of a $(-1)$-curve $E_{i}$ disjoint from $\Delta_{i-1}$.
In particular, $Y$ is smooth and $\Delta$ is normal crossings.
\end{lem}

\begin{proof}
Let $E_{0}$ be a $\rho_{1}$-exceptional curve on $Y_{0}$.
Then $(K_{Y_0}+\Delta_{0})\cdot E_{0}<0$.
Since $K_{Y_{0}}+\Delta_{0}$ is numerically trivial over $\Delta_{0}$, the exceptional curve $E_{0}$ is not contained in $\Delta_{0}$.
It follows that $K_{Y_{0}}\cdot E_{0}<0$, in other words, $\rho_{0}$ is a $K_{Y_{0}}$-negative extremal contraction.
Hence $E_{0}$ is a $(-1)$-curve and $\rho_{0}$ is the blowing down of $E_{0}$.
Moreover, $E_{0}$ is disjoint from $\Delta_{0}$ since 
$$
0>(K_{Y_{0}}+\Delta_{0})\cdot E_{0}\ge K_{Y_{0}}\cdot E_{0}=-1.
$$
Hence $Y_{1}$ is smooth, $\Delta_{1}$ is normal crossings and $K_{Y_1}+\Delta_{1}$ is numerically trivial over $\Delta_{1}$.
By induction on $m$, the claim follows.
\end{proof}

In this case, we have
$$
h^{0}(\omega_{Y}(\Delta))=h^{0}(\omega_{\widetilde{X}}(\widetilde{\Delta}))=h^{0}(\omega_{X})=p_g(X)
$$
and
$$
(K_{Y}+\Delta)^{2}=\mathrm{vol}(K_{\widetilde{X}}+\widetilde{\Delta})\le \mathrm{vol}(K_{X})=K_{X}^{2},
$$
where the last inequality is strict if and only if there exists a non-Gorenstein singularity on $X$. 
The following Noether-type inequalities will be used later.

\begin{thm}[cf.\ \cite{Che}] \label{normalstable}
Let $X$ be a normal stable surface and let $(Y,\Delta)$ be the pair described as above.
Then the following holds.

\begin{itemize}
    \item[$(1)$] 
    Assume that $|K_{Y}+\Delta|$ is composed with a pencil.
    Then one of the following holds.
    \begin{itemize}
        \item[$(1.1)$] 
        $\kappa(K_{Y}+\Delta)=1$, and $|K_{Y}+\Delta|$ gives an elliptic fibration over $\mathbb{P}^1$.
        The contraction from $Y$ to its relatively minimal model $Y_{\min}$ is a successive blowing down of vertical $(-1)$-curves intersecting with or contained in the proper transform of $\Delta$.
        The image of each connected component of $\Delta$ to $Y_{\min}$ coincides with a multiple fiber, or a fiber of type $\mathrm{I}$, that is, $\mathrm{I}_n$ for some $n$. 
        Moreover, $q(X)=0$ holds.
        \item[$(1.2)$]
        $\kappa(K_{Y}+\Delta)=2$ and $(K_{Y}+\Delta)^{2}\ge 2p_{g}(X)-2$.
        \item[$(1.3)$] 
        $\kappa(K_{Y}+\Delta)=2$, $(K_{Y}+\Delta)^{2}=1$, $p_{g}(X)=2$, and $|K_{Y}+\Delta|$ is a pencil of genus $2$ with one base point.
        \item[$(1.4)$]
        $\kappa(K_{Y}+\Delta)=2$ and  $(K_{Y}+\Delta)^{2}=p_{g}(X)$.
        In this case, $Y$ admits an elliptic fibration over an elliptic curve and the horizontal part of $\Delta$ is a section.
        The contraction from $Y$ to its minimal model $Y_{\min}$ is a successive blowing down of vertical $(-1)$-curves intersecting with or contained in the proper transform of $\Delta$ and the image of each connected component of the vertical part of $\Delta$ to $Y_{\min}$ coincides with a fiber of type $\mathrm{I}$.
        Moreover, $q(X)=0$.
    \end{itemize}

    \item[$(2)$] 
    Assume that $|K_{Y}+\Delta|$ is not composed with a pencil.
    Then we have 
    \[
    (K_{Y}+\Delta)^{2}\ge 2p_{g}(X)-4
    \]
    with the equality holding only if $|K_{Y}+\Delta|$ is basepoint-free and defines a generically finite morphism of degree $2$ onto a minimal degree surface in $\mathbb{P}^{p_g(X)-1}$.
\end{itemize}
\end{thm}

\begin{proof}
First, assume that $|K_{Y}+\Delta|$ is composed with a pencil.
It is written as
$$
|K_{Y}+\Delta|=|M|+V,
$$
where $|M|=|F_1+\cdots+F_{n}|$ is the movable part and $V$ is the fixed part.
Let $Y'\to Z$ be the Stein factorization of a resolution of indeterminacy of the rational map associated to $|M|$.
Then $|M|$ is induced by some complete linear system $|M_{Z}|$ of degree $n$ on the curve $Z$.
Applying the Riemann-Roch theorem to $M_{Z}$, it follows that
\begin{equation} \label{pgn}
p_{g}(X)=1+n-g(Z)+h^{1}(M_{Z})\le 1+n
\end{equation}
with the equality holding if and only if $g(Z)=0$.
Let $F$ denote the numerical class such that $M\equiv nF$.
Then we have
$$
(K_Y+\Delta)^{2}=(K_Y+\Delta)\cdot (nF+V)=n(K_Y+\Delta)\cdot F+(K_Y+\Delta)\cdot V
$$
and both $(K_Y+\Delta)\cdot F$ and $(K_Y+\Delta)\cdot V$ are non-negative since $K_Y+\Delta$ is nef.

Suppose that $(K_Y+\Delta)\cdot F=0$.
Then $F^2=0$ and $V\cdot F=0$ since $F$ is nef.
Hence, $Y'=Y$ and $V$ consists of vertical components.
Since $K_{Y}+\Delta$ is nef, it follows that
$$
0\le (K_Y+\Delta)^2=(nF+V)^2=V^2\le 0,
$$
whence $V^2=0$.
Since $V$ is a fixed part, we have $V\equiv aF$ for some $a\in\mathbb{Q}_{\ge0}$.
In particular, 
$\kappa(K_Y+\Delta)=1$.
Since $X$ is stable, there exists a $\pi$-exceptional rational curve on $\widetilde{X}$ horizontal over $Z$.
This shows that $Z\cong\mathbb{P}^1$ and $p_g(X)=n+1$.
Hence the map $Y\to Z$ coincides with an Iitaka fibration with respect to $K_{Y}+\Delta$.
Note that $\Delta$ consists of vertical components since $K_{Y}+\Delta$ is numerically trivial over $\Delta$.
In particular, a general fiber $F$ has genus $1$.
Conversely, $\kappa(K_{Y}+\Delta)=1$ implies $(K_Y+\Delta)\cdot F=0$ since $(K_Y+\Delta)^2=0$.
Hence the first part of the claim in (1.1) holds.

Next, we prove the remaining part of the statement (1.1).
Note that $\Delta\cdot E=-K_{Y}\cdot E=1$ for each vertical $(-1)$-curve $E$ on $Y$ since $K_{Y}+\Delta$ is nef.
Let $Y_{\min}$ be the relatively minimal model of $Y$ and $\Delta_{\min}$ the proper transform of $\Delta$.
Then, $Y_{\min}$ is a minimal elliptic surface over $Z$, $\Delta_{\min}$ is a finite union of fibers of type $\mathrm{I}$ and reduced parts of multiple fibers, and $(Y_{\min},\Delta_{\min})$ is log crepant to $(Y,\Delta)$.
Using this and the canonical bundle formula for elliptic surfaces, it is easy to see that the support of $V$ does not contain any component of $\Delta$.
Therefore, letting $l_E$ denote the number of elliptic singularities of $X$, we obtain $p_g(Y)=n+1-l_E$ and $\chi(\mathcal{O}_Y)=\chi(\mathcal{O}_X)-l_E$.
This implies that $q(Y)=q(X)$.
Moreover, the fact that $q(Y)=0$ follows from the canonical bundle formula for elliptic fibrations over $\mathbb{P}^1$ (see \cite[V, Section 12]{BHPV}).
Hence, we obtain the latter assertion of (1.1).

Next, suppose that $(K_{Y}+\Delta)\cdot F=1$.
Then either $F^2=1$, $F\cdot V=0$ and $n=1$ or $F^2=0$ and $F\cdot V=1$ hold.

If the former holds, then $p_g(X)=2$ and $Z\cong \mathbb{P}^1$ from \eqref{pgn}.
Since $F\cdot V=0$ and $F$ is nef and big, $V^2\le 0$ with the equality holding if and only if $V=0$ by the Hodge index theorem.
On the other hand, $V^2=(K_{Y}+\Delta)\cdot V\ge 0$,
whence $V=0$ and so $K_{Y}+\Delta=F$.
In particular, $(K_Y+\Delta)^{2}=1$.
By using adjunction, it follows that $|K_{Y}+\Delta|$ is a pencil of genus $2$ with one base point, whence (1.3) holds.

If the latter holds, then $|M|$ is basepoint-free.
Let $\varphi\colon Y\to Z$ denote the corresponding fibration with a general fiber $F$.
By adjunction, we have
$$
2g(F)-2=(K_Y+F)\cdot F=K_Y\cdot F=(K_Y+\Delta)\cdot F-\Delta \cdot F=1-\Delta \cdot F.
$$
Hence either $g(F)=0$ and $\Delta \cdot F=3$ or $g(F)=1$ and $\Delta \cdot F=1$.
In particular, the horizontal part $\Delta_{\mathrm{hor}}$ of $\Delta$ is non-zero. 
For each irreducible component $C\subset \Delta_{\mathrm{hor}}$,
we have
$$
V\cdot C=(K_Y+\Delta)\cdot C-nF\cdot C=-nF\cdot C<0.
$$
Since $\Delta$ is reduced, we have $\Delta_{\mathrm{hor}} \le V$.
Thus $\Delta \cdot F\le V\cdot F=1$,
which implies that $\Delta\cdot F=1$ and $g(F)=1$.
Hence $\Delta_{\mathrm{hor}}$ is a section of $\varphi$, which is denoted by $S$.
It is an elliptic curve since each connected component of $\Delta$ is a cycle of rational curves or an elliptic curve and $\Delta_{\mathrm{vert}}=\Delta-S$.
It follows from \eqref{pgn} and $g(Z)=g(S)=1$ that $p_g(X)=n$.
Let $Y_{\mathrm{min}}$ be the minimal model of $Y$.
The contraction $\psi\colon Y\to Y_{\mathrm{min}}$ decomposes into 
$$
Y\xrightarrow{\psi_{1}}Y_{1}\xrightarrow{\psi_{2}}Y_{2}=Y_{\mathrm{min}},
$$
where $\psi_{2}$ is a sequence of blow-ups at $\psi(S)\cap \psi(\Delta_{\mathrm{vert}})$.
Let $\Delta_{i}$ (resp.\ $S_{i}$, $F_{i}$) denote the pushforward of $\Delta$ (resp.\ $S$, $F$) to $Y_{i}$.
Then we have $K_{Y}+\Delta=\psi_{1}^{*}(K_{Y_1}+\Delta_{1})$ and
$$
K_{Y_1}+\Delta_{1}=\psi_{2}^{*}(K_{Y_2}+\Delta_{2})-\sum_{i=1}^{l}E_{i}=\psi_{2}^{*}(K_{Y_2}+\Delta_{2}-S_{2})+S_{1},
$$
where $E_{i}$'s are $\psi_{2}$-exceptional $(-1)$-curves and the number $l$ coincides with the number of connected components of $\Delta_{\mathrm{vert}}$.

Since $\Delta_{2,\mathrm{vert}}=\Delta_{2}-S_{2}$ consists of $l$ fibers of type $\mathrm{I}$ and the canonical bundle formula, $K_{Y_2}+\Delta_{2}-S_{2}$ is numerically equivalent to $(\chi(\O_{Y})+l)F_{2}$.
It follows from $K_{Y}+\Delta\equiv nF+V$ that $n=\chi(\O_{Y})+l$ and $V=\psi_{1}^{*}S_{1}$.
Hence we have
$$
(nF_{1}+S_{1})\cdot S_{1}=\psi_{1}^{*}(nF_{1}+S_{1})\cdot S=(K_{Y}+\Delta)\cdot S=0,
$$
where the last equality follows from the numerical triviality of $K_{Y}+\Delta$ over $\Delta$.
Thus 
$$
(K_{Y}+\Delta)^{2}=\psi_{1}^{*}(K_{Y_1}+\Delta_{1})^{2}=(nF_{1}+S_{1})^{2}=(nF_{1}+S_{1})\cdot nF_{1}=n=p_g(X),
$$
whence the case (1.4) holds.
Note that $p_{g}(X)=n=\chi(\O_{Y})+l=\chi(\O_{X})-1$ implies $q(X)=0$.

If $(K_Y+\Delta)\cdot F\ge 2$, then it follows from \eqref{pgn} that
$$
(K_Y+\Delta)^{2}\ge 2n+(K_Y+\Delta)\cdot V\ge 2p_g(X)-2+(K_Y+\Delta)\cdot V\ge 2p_g(X)-2,
$$
which implies the case (1.2).

Finally we assume that $|K_{Y}+\Delta|$ is not composed with a pencil.
If $p_g(X)\le 2$, then the claim~(2) is clear.
Thus, we may assume that $p_g(X)\ge 3$.
Let $C\in |M|$ be a general member, which can be taken as an irreducible and reduced curve by Bertini's theorem.
Let $\Delta=\sum_{p}\Delta_{p}$ denote the connected component decomposition.
Let $Y\to Y'$ denote the birational contraction of all the $\Delta_{p}$'s which are disjoint from $C$.
Let $\Delta'$ denote the pushforward of $\Delta$ to $Y'$.
Note that the contraction $Y\to Y'$ is log crepant and hence $(K_Y+\Delta)^{2}=(K_{Y'}+\Delta')^{2}$.
Since $K_{Y'}+\Delta'$ is numerically trivial over $\Delta'$ and $C$ intersects with any connected component of $\Delta'$, it follows that $\Delta'$ is contained in the fixed part $V'$ of $|K_{Y'}+\Delta'|$.
In particular, $K_{Y'}-C\sim V'-\Delta'$ is effective.
Since
$$
|K_{C}|\ge |K_{Y'}+C||_{C}\ge |2C||_{C}\ge |C||_{C}+|C||_{C}
$$
and $\dim |C||_{C}=p_g(X)-2$,
we have $p_{a}(C) \ge 2p_g(X)-3$.
By adjunction, $p_a(C)=\frac{1}{2}(K_{Y'}+C)\cdot C+1$.
Since $C$ and $K_{Y'}+\Delta'$ are nef,
we have
$$
(K_{Y'}+\Delta')^{2}-\frac{1}{2}(K_{Y'}+C)\cdot C
=\frac{1}{2}C\cdot (V'+\Delta')+(K_{Y'}+\Delta')\cdot V'\ge 0.
$$
Therefore,
\begin{align*}
(K_{Y}+\Delta)^{2}&=
(K_{Y'}+\Delta')^{2} \\
&\ge \frac{1}{2}(K_{Y'}+C)\cdot C \\
&=p_a(C)-1 \\
&\ge 2p_g(X)-4.
\end{align*}
The rest of the claim is straightforward (see the proof of \cite[Lemma 1.1]{horikawa}).
\end{proof}

\subsection{Formula for normal stable surfaces}
Let $X$ be a normal stable surface with only $\Q$-Gorenstein smoothable singularities.
In this case, any singularity of $X$ is either a T-singularity, a strictly lc rational singularity (Lemma~\ref{lem:smoothable_rational_strict_lc}), or an elliptic singularity (Lemma \ref{lem:smoothable_elliptic_singularity}, Appendix~\ref{app:cusp}).
The aim of this and next subsections is to study the minimal model of the minimal resolution of $X$ and to provide a useful formula for $X$. 

We introduce the following notations.

\begin{note}\label{note:normal_stable_surface}
Let $X$ be a normal stable surface with only $\Q$-Gorenstein smoothable singularities.
    \begin{itemize}
        \item[$(1)$] 
        Let $\pi\colon \widetilde{X}\to X$ denote the minimal resolution, and let $\widetilde{C}$ denote the (reduced) inverse image on $\widetilde{X}$ of the non-canonical singularities on $X$.
        \item[$(2)$]
        Let 
        \[
        \rho\colon \widetilde{X}=Y_{0}\xrightarrow{\rho_{1}} Y_{1}\xrightarrow{\rho_{2}} \cdots \xrightarrow{\rho_{m}} Y_{m}=Y
        \]
        be any sequence of contractions of $(-1)$-curves,
        where each $\rho_{k}\colon Y_{k-1}\to Y_{k}$ is the blowing down of a $(-1)$-curve $E_{k}$ to a point $y_{k}\in Y_{k}$.
        Let $C_{k}$ denote the pushforward of $\widetilde{C}$ to $Y_{k}$ and $C:=C_{m}$.
        Let $M_{k}$ denote the multiplicity of $C_{k}$ at $y_{k}$. 
        If $y_{k}$ is not contained in $C_{k}$, we set $M_{k}=0$.
        
    \end{itemize}
\end{note}

The following lemma, which generalizes {\cite[Lemmas~3.11 and 3.12]{FRU}}, is a key observation.

\begin{lem} \label{multi}
$M_{k}\ge 2$ for each $k$.
\end{lem}

Before proceeding to the proof, we introduce a useful notion.

\begin{defn} \label{completelyseparated}
Let $D$ be a (topologically) connected subdivisor of $C_k$ and let $E$ be an irreducible component of $C_{k}$ that meet $D$ transversally at exactly one point $p$.
Let $\widetilde{p}$ denote the point in the proper transform of $E$ on $\widetilde{X}$ that corresponds to $p$.
We say that \textit{$D$ is completely separated from $E$ after the blow-ups} if $\widetilde{C}$ is smooth at $\widetilde{p}$. 
\end{defn}

 \begin{proof}[Proof of Lemma \ref{multi}]
 We show the claim by induction on the number $m$ of blow-ups in $\rho$.
Since $\pi$ is a minimal resolution, the $(-1)$-curve $E_{1}$ is not contained in $\widetilde{C}$.
If $M_{1}=0$, then $\pi(E_{1})$ is still a $(-1)$-curve in the smooth locus of $X$, which contradicts to the ampleness of $X$.
If $M_{1}=1$, then $E_{1}$ intersects only one irreducible component of $\widetilde{C}$, say $D$, at one point transversally.
It follows that
$$
K_{X}\cdot \pi(E_{1})=\pi^{*}K_{X}\cdot E_{1}=(K_{\widetilde{X}}-aD)\cdot E_{1}=-1-a\le 0,
$$
which is also a contradiction, where the last inequality follows from the log-canonicity of $X$.

Assume that $m\ge 2$ and $M_{k}\ge 2$ for $k=1,\ldots,m-1$.
If $E_{m}$ is not contained in $C_{m-1}$, then the same argument as above shows that $M_{m}\ge 2$.
Indeed, if $M_{m}\le 1$, then there are no blow-ups over
$E_{m}$ by inductive hypothesis which leads to a contradiction similarly.
Hence we may assume that $E_{m}$ is contained in $C_{m-1}$.
If $M_{m}=0$, the same argument as above also leads to a contradiction since $\widetilde{C}$ has no $(-1)$-curves.
Hence we may assume that $M_{m}=1$.
Then $E_{m}$ and other components of $C_{m-1}$ intersects at one point transversally.
A blow-up must occur at this point, say $y_{m-1}$,
since $E_{m}$ is a $(-1)$-curve and $\pi$ is a minimal resolution.
Let $D$ denote the irreducible component of $C_{m}$ passing through $y_{m}$ and $\pi^{-1}(p_1)$ denote the connected component of $\widetilde{C}$ containing the proper transform of $D$.
By inductive hypothesis, the union of $\pi^{-1}(p_1)$ and the $\rho$-exceptional curves over $y_{m}$ are described as 
$$
\pi^{-1}(p_{1})+E^{1}+\cdots+\pi^{-1}(p_{s-1})+E^{s-1}+\pi^{-1}(p_{s})
$$
for some $s\ge 2$,
where $p_{i}$'s are non-Du Val singularities on $X$ and each $E^{i}$ is a $(-1)$-curve connecting $\pi^{-1}(p_{i})$ and $\pi^{-1}(p_{i+1})$.
Note that $\pi^{-1}(p_{i})$ is a chain and so $p_{i}$ is a T-singularity for each $i\ge 2$, but $p_{1}$ may be strictly lc.
From {\cite[Corollary~2.4]{FRU}}, we may assume that $p_{1}$ is strictly lc.
If $E_{m}$ is completely separated from the proper transform of $D$ after blow-ups,
then $\pi^{-1}(p_{2})+E^{2}+\cdots+\pi^{-1}(p_{s})$ comes from the extended T-chain $[2,\ldots,2]$, which is a contradiction due to Lemma~\ref{lem:extT_2m2}.
Hence $\pi^{-1}(p_{1})$ has a sub-chain which is exceptional over $y_{m-1}$.
In particular, $p_{1}$ is rational.
By Lemma~\ref{lem:smoothable_rational_strict_lc}, the above sub-chain consists of a single $(-r)$-curve for some $r=2,3,4$ or $6$.
Then $\pi^{-1}(p_{2})+E^{2}+\cdots+\pi^{-1}(p_{s})$ comes from the extended T-chain 
$$
L_2^{r-2}L_1[2^{\beta+1}]=[2^{r-2},3,2^{\beta}]
$$
for some $\beta\ge 0$, which is a contradiction by Lemma~\ref{lem:extT_2m2}. 
 \end{proof}

We define a number $\delta_{k}$ as $1$ if $E_{k}$ is contained in $C_{k-1}$ and $0$ otherwise.
Let $l_{T}$, $l_{R_{0}}$, $l_{R_{1}}$ and $l_{E}$ respectively denote the number of T-singularities on $X$ which are not Du Val, the number of strictly lc rational singularities of type $(2,2,2,2)$ on $X$, 
the number of strictly lc rational singularities on $X$ as in Lemma~\ref{lem:smoothable_rational_strict_lc}~(ii--iv) and the number of elliptic singularities on $X$. 
It follows from the Leray spectral sequence that
\begin{equation} \label{ellsingeuler}
\chi(\O_{Y})=\chi(\O_{\widetilde{X}})=\chi(\O_{X})-l_{E}.
\end{equation}

\begin{lem} \label{fundeq}
We have 
$$
K_{Y}^{2}+K_{Y}\cdot C+\sum_{k=1}^{m}(M_{k}-\delta_{k}-1)=K_{X}^2+l_{T}+l_{R_{1}}.
$$
\end{lem}

\begin{proof}
Note that for any point $p$ of $X$, the equality 
$$
K_{p}\cdot \widetilde{C}+K_{p}^{2}=\delta_{p}
$$
holds, where $K_{p}$ is the canonical cycle of $p$ and
$\delta_{p}$ is defined as $1$ if $p$ is either a non-Du Val T-singularity or a strictly lc rational singularity as in Lemma~\ref{lem:smoothable_rational_strict_lc}~(ii--iv), and $0$ otherwise.
Hence we have

\begin{align*}
K_{\widetilde{X}}\cdot \widetilde{C}&=(\pi^{*}K_{X}+\sum_{p}K_{p})\cdot \widetilde{C} \\
&=\sum_{p}K_{p}\cdot \widetilde{C} \\
&=l_{T}+l_{R_{1}}-\sum_{p}K_{p}^{2} \\
&=l_{T}+l_{R_1}+K_{X}^{2}-K_{\widetilde{X}}^{2} \\
&=l_{T}+l_{R_1}+K_{X}^{2}-K_{Y}^{2}+m.
\end{align*}

On the other hand, the equations
$$
C_{k-1}=\rho_{k}^{*}C_{k}-(M_{k}-\delta_{k})E_{k}
$$
 for $k=1,\ldots,m$ imply
$$
K_{\widetilde{X}}\cdot \widetilde{C}=K_{Y}\cdot C+\sum_{k=1}^{m}(M_{k}-\delta_{k}).
$$
Combining the above two equalities, the claim follows.
\end{proof}

Let us define some notations used throughout the remaining sections.
Let $C=\sum_{\alpha}C^{(\alpha)}$ denote the connected component decomposition.
Let $\widetilde{C}=\sum_{\alpha}\widetilde{C}^{(\alpha)}$
and $C_{k}=\sum_{\alpha}C_{k}^{(\alpha)}$ denote the corresponding decompositions.
For a connected component $C^{(\alpha)}$, let $V_{\alpha}$ denote the set of non-Du Val singularities on $X$ such that $\rho(\pi^{-1}(p))$ is contained in $C^{(\alpha)}$.
Let $l_{T}^{(\alpha)}$ (resp.\  $l_{R_0}^{(\alpha)}$, $l_{R_1}^{(\alpha)}$, $l_{E}^{(\alpha)}$) denote the number of non-Du Val T-singularities
(resp.\ strictly lc rational singularities of type $(2,2,2,2)$, strictly lc rational singularities as in Lemma~\ref{lem:smoothable_rational_strict_lc}~(ii--iv), elliptic singularities) in $V_{\alpha}$.
Then 
$$
\#V_{\alpha}=l_{T}^{(\alpha)}+l_{R_0}^{(\alpha)}+l_{R_1}^{(\alpha)}+l_{E}^{(\alpha)},\quad l_{T}=\sum_{\alpha}l_{T}^{(\alpha)},
\quad l_{R_i}=\sum_{\alpha}l_{R_i}^{(\alpha)},
\quad l_{E}=\sum_{\alpha}l_{E}^{(\alpha)}.
$$
For a non-Du Val singularity $p$ on $X$, let $\widetilde{C}^{(p)}$, $C_{k}^{(p)}$ and $C^{(p)}$ respectively denote the subdivisors of $\widetilde{C}$, $C_{k}$ and $C$  which are exceptional over $p$.
Clearly we have
$$
\widetilde{C}^{(\alpha)}=\sum_{p\in V_{\alpha}}\widetilde{C}^{(p)},\quad C_{k}^{(\alpha)}=\sum_{p\in V_{\alpha}}C_{k}^{(p)},\quad C^{(\alpha)}=\sum_{p\in V_{\alpha}}C^{(p)}.
$$
Let 
$$
M(\alpha):=\sum_{y_{k}\in C_{k}^{(\alpha)}}(M_{k}-\delta_{k}-1).
$$
Since the decomposition $C_{k}=\sum_{\alpha}C_{k}^{(\alpha)}$ is disjoint, we have 
\begin{equation} \label{multdecomp}
\sum_{k=1}^{m}(M_{k}-\delta_{k}-1)=\sum_{\alpha}M(\alpha).
\end{equation}

\begin{lem} \label{multbound}
For any connected component $C^{(\alpha)}$, we have
$$
M(\alpha)\ge \#V_{\alpha}-1.
$$
\end{lem}

\begin{proof}
By reordering the blow-ups if necessary, we may assume that $y_{k}$ is contained in $C_{k}^{(\alpha)}$ if $k=1,\ldots,m_{\alpha}$ and not contained in $C_{k}^{(\alpha)}$ if $k> m_{\alpha}$ for some $m_{\alpha}\le m$. 
Let $\gamma_{k}$ be the number of connected components of  $C^{(\alpha)}_{k}$ on $Y_{k}$
for $k=0,\ldots,m_{\alpha}$.
Note that 
$$
\#V_{\alpha}=\gamma_{0}\ge \gamma_{1}\ge \cdots \ge \gamma_{m_{\alpha}}=1.
$$
If $\delta_{k}=1$, then $\gamma_{k-1}=\gamma_{k}$ holds since $E_{k}$ is contained in $C^{(\alpha)}_{k-1}$ and intersects with the proper transform of any analytic branch of $C^{(\alpha)}_{k}$ passing through $y_{k}$.
If $\delta_{k}=0$, then $M_{k}-1\ge \gamma_{k-1}-\gamma_{k}$ holds.
Hence we have 
$$
M_{k}-\delta_{k}-1\ge \gamma_{k-1}-\gamma_{k}
$$
for each $k=1,\ldots,m_{\alpha}$.
Summing it up, the claim holds.
\end{proof}

\subsection{Elliptic case}\label{elliptic_case}
In this subsection, we further assume that $Y$ admits a relatively minimal elliptic fibration $f\colon Y\to Z$.
Let $F$ denote a general fiber of $f$.
By the canonical bundle formula (see \cite[V.~Theorem 12.1]{BHPV}), we have
$$
K_{Y}\equiv \left(\chi(\O_{Y})+2g(Z)-2+\sum_{\lambda}\left(1-\frac{1}{m_{\lambda}}\right)\right)F,
$$
where $m_{\lambda}$'s are the multiplicities of multiple fibers of $f$.
If there exists a multiple fiber, let $m_X$ be the maximum multiplicities of multiple fibers.
Otherwise, set $m_X=1$.
For a non-Du Val singularity $p$ of $X$, let $n^{(p)}:=F\cdot C^{(p)}$.
Then we have
\begin{equation} \label{ellKC}
K_{Y}^{2}=0,\quad K_{Y}\cdot C=\left(\chi(\O_{Y})+2g(Z)-2+\delta_{f} \right)n_{X},
\end{equation}
where $n_X:=\sum_{p}n^{(p)}$ and $\delta_{f}:=\sum_{\lambda}(1-1/m_{\lambda})$.
The following lemma implies that each connected component of $C$ has at least one horizontal component.

\begin{lem} \label{degbound}
For any connected component $C^{(\alpha)}$ of $C$, we have
$$
n^{(\alpha)}:=\sum_{p\in V_{\alpha}}n^{(p)}\ge m_X.
$$
\end{lem}

\begin{proof}
Note that clearly $n^{(\alpha)}$ is divisible by $m_X$, and that $n_X>0$ holds since otherwise $X$ would admit an elliptic fibration which contradicts the ampleness of $K_{X}$.
Assume on the contrary that $n^{(\alpha)}=0$.
Then the connected curve $C^{(\alpha)}$ is contained in some fiber $f^{-1}(x)$ of $f$.
Since $n^{(\beta)}\ge m_X$ for some $\beta\neq \alpha$, the curve $C^{(\beta)}$ intersects with the fiber $f^{-1}(x)$.
Since $C^{(\alpha)}$ and $C^{(\beta)}$ are disjoint, $C^{(\alpha)}$ does not contain all components of $f^{-1}(x)$
and so the intersection form on $C^{(\alpha)}$ is negative definite.
Let $\widetilde{\Delta}:=-\sum_{p}K_{p}$.
Then $K_{\widetilde{X}}+\widetilde{\Delta}=\pi^{*}K_{X}$ is nef.
Since $Y$ is a surface, the pushforward $K_{Y}+\rho_{*}\widetilde{\Delta}=\rho_{*}(K_{\widetilde{X}}+\widetilde{\Delta})$ is also nef.
Let $\rho_{*}\widetilde{\Delta}=\sum_{\gamma}\Delta^{(\gamma)}$ be a decomposition such that the support of $\Delta^{(\gamma)}$ equals $C^{(\gamma)}$ for each $\gamma$.
Since $K_{Y}$ is $f$-numerically trivial,
we have
$$
(K_{Y}+\rho_{*}\widetilde{\Delta})\cdot \Delta^{(\alpha)}=(\rho_{*}\widetilde{\Delta})\cdot \Delta^{(\alpha)}=(\Delta^{(\alpha)})^{2}<0,
$$
which is a contradiction to the nefness of $K_{Y}+\rho_{*}\widetilde{\Delta}$.
\end{proof}

Let $l_{E_{v}}$ denote the number of elliptic singularities $p\in X$ such that the exceptional set $\pi^{-1}(p)$ is contained in some fiber of the elliptic fibration $\widetilde{X}\to Z$,
and let $l_{E_{h}}:=l_{E}-l_{E_{v}}$.
Note that if $l_{E_{v}}>0$, then $C$ is connected since the set $\rho(\pi^{-1}(p))$, where $p\in X$ is an elliptic singularity contributing to $l_{E_{v}}$, coincides with a fiber of type $\mathrm{I}$ or a multiple fiber and hence intersects with any connected component of $C$ by Lemma~\ref{degbound}.
Here, we define $l_{E_v,m}$ as the number of elliptic singularities $p$ such that $C^{(p)}$ corresponds to a multiple fiber with multiplicity $m$ for $1\le m\le m_X$.

\begin{lem} \label{multbound2}
In the above situations, we have
$$
r^{(\alpha)}:=M(\alpha)-\#V_{\alpha}+\frac{n^{(\alpha)}}{m_X}-\sum_{1\le m\le m_X}\left(\frac{n^{(\alpha)}}{m}-1\right)l_{E_{v},m}\ge 0.
$$
\end{lem}

\begin{proof}
If $l_{E_{v}}=0$, then the claim follows from Lemmas~\ref{multbound} and \ref{degbound}.
Hence we may assume that $l_{E_{v}}>0$.
In this case, we note that $C$ is connected, and hence $C=C^{(\alpha)}$.
Let $p_{1},\ldots,p_{l_{E_{v}}}$ be the elliptic singularities contributing to $l_{E_{v}}$.
Let $m_i$ be the multiplicity of the fiber corresponding to $C^{(p_i)}$.
In the proof of Lemma~\ref{multbound}, the blow-ups over $C^{(p_i)}$ contribute at least $\frac{n^{(\alpha)}}{m_i}$ to $M(\alpha)$ for each $i$.
However, after the blow-ups over $\sum_{i=1}^{l_{E_{v}}-1}C^{(p_i)}$,
$C_{k}$ has exactly $l_{E_v}$ connected components, and each component contains a unique $C^{(p_j)}_k$ associated with a unique elliptic singularity $p_j$.
Hence, the claim follows from the similar arguments as in the proof of Lemma~\ref{multbound}. 
\end{proof}

The following fundamental equality is important to classify normal stable Horikawa surfaces with only $\Q$-Gorenstein smoothable singularities.

\begin{thm} \label{ellfundeq}
Let $X$ be a normal stable surface with only $\Q$-Gorenstein smoothable singularities,
and assume that the minimal resolution of $X$ admits an elliptic fibration over a curve $Z$.
Then we have
\begin{equation}
\label{eqn:ellfundeq}
K_{X}^{2}=\left(\chi(\O_{X})+2g(Z)-3+\delta'_{f}\right)n_{X}+r+l_{R_0}-(n_X-1)l_{E_{h}},
\end{equation}
where 
$$
r:=\sum_{\alpha}r^{(\alpha)},\quad \delta'_f:=\delta_f-\sum_{1\le m\le m_X}\left(1-\frac{1}{m}\right)l_{E_v,m}+1-\frac{1}{m_X}
$$
and these numbers are non-negative.
In particular, if $X$ belongs to (1.1) (resp.\ (1.4)) in Theorem~\ref{normalstable}, then we have
$K_{X}^{2}\ge (p_g(X)-2+\delta'_f)n_X$ (resp.\ $K_{X}^{2}\ge (p_g(X)-1)n_X+1$).
\end{thm}

\begin{proof}
The first assertion follows from \eqref{ellsingeuler}, \eqref{multdecomp}, \eqref{ellKC} and Lemma~\ref{fundeq}.
Indeed, 
\begin{align*} 
K_{X}^{2}&=\left(\chi(\O_{Y})+2g(Z)-2+\delta_{f}\right)n_X+\sum_{\alpha}M(\alpha)-l_{T}-l_{R_{1}} \\ \nonumber
&=\left(\chi(\O_{X})+2g(Z)-3+\delta_{f}\right)n_X+\sum_{\alpha}(n^{(\alpha)}+M(\alpha))-n_Xl_{E}-l_{T}-l_{R_1} \\ \nonumber
&=\left(\chi(\O_{X})+2g(Z)-3+\delta'_{f}\right)n_X+\sum_{\alpha}r^{(\alpha)}+l_{R_0}-(n_X-1)l_{E_{h}}.
\end{align*}
In the third line, we make use of the equality $n^{(\alpha)}+M(\alpha)=(1-\frac{1}{m_X})n^{(\alpha)}+\# V_{\alpha}+\sum_{m>0}(\frac{n^{(\alpha)}}{m}-1)l_{E_v,m}+r^{(\alpha)}$.
The last assertion follows directly from Theorem~\ref{normalstable} and Lemma~\ref{multbound2}.
Indeed, note that condition (1.1) (resp.\ (1.4)) in Theorem~\ref{normalstable} implies that $(g(Z), l_{E_{h}})=(0,0)$ (resp.\ $(g(Z), l_{E_{h}},m_X)=(1,1,1)$) and 
$\chi(\O_{X})+2g(Z)-3+\delta'_{f}=p_{g}(X)-2+\delta'_f$ (resp.\ $\chi(\O_{X})+2g(Z)-3+\delta'_{f}=p_{g}(X)$).
\end{proof}

As a corollary, the following Noether-type inequality can be reproved for locally $\Q$-Gorenstein smoothable normal stable surfaces,
which is a special version of {\cite[Theorem~1.1]{Liu}}:

\begin{cor}
Let $X$ be a normal stable surface with only $\Q$-Gorenstein smoothable singularities.
Then we have
$K_{X}^{2}\ge p_{g}(X)-2$.
\end{cor}

\begin{proof}
If $p_{g}(X)\le 2$, the claim follows since $K_{X}$ is ample.
Then we may assume that $p_{g}(X)>2$.
If $K_{X}^{2}\ge 2p_{g}(X)-4$, then the claim trivially holds.
Hence we may also assume that $K_{X}^{2}<2p_{g}(X)-4$.
Applying Theorem~\ref{normalstable}, the minimal resolution $\widetilde{X}\to X$ admits an elliptic fibration over a rational or elliptic curve $B$ and either (1.1) or (1.4) in Theorem~\ref{normalstable} holds.
Hence the claim follows from Theorem~\ref{ellfundeq}.
\end{proof}

Finally, we state the following lemma. 
Its proof is a direct computation.

\begin{lem}\label{lem:no_bisection}
    If $\delta'_f=0$ and $r=0$, then $\widetilde{C}$ contains exactly $n_X$ horizontal components.
\end{lem}

\section{Normal stable Horikawa surface}\label{sec:normal_stbale_surface}
From now on, we focus on Horikawa surfaces.
We define normal stable Horikawa surfaces, and divide them into two classes: standard and non-standard.
We also introduce the notion of \emph{good involutions}.

\subsection{Standard/Non-standard Horikawa surface}\label{subsec:Horikawa}


\begin{defn}\label{def:Horikawa}
A \emph{stable} (resp. \emph{normal stable}) \emph{Horikawa surface} $X$ is a stable (resp. normal stable) surface with $K_X^2=2p_g(X)-4$.
We say that a normal stable Horikawa surface $X$ is \emph{standard} if its canonical linear system is not composed with a pencil, and is \emph{non-standard} if $X$ is not standard.
\end{defn}

The following proposition is a direct consequence of Theorem \ref{normalstable}~(2), which justifies the use of the term ``standard''.

\begin{prop}[cf. {\cite[Theorem~3.3]{Che}}] \label{prop:standard_Horikawa}
Let $X$ be a standard Horikawa surface. 
Then, $X$ is Gorenstein and its canonical map gives a double covering $\varphi_{K_{X}}\colon X\to W$ over a minimal degree surface $W \subset \mathbb{P}^{p_g(X)-1}$.
\end{prop}

\begin{proof}
    From Theorem~\ref{normalstable}~(2), the linear system $|K_{Y}+\Delta|$ is basepoint-free and defines a generically finite morphism $\varphi\colon Y\to W$ to a minimal degree surface $W\subset \mathbb{P}^{p_g-1}$.
    Note that $\widetilde{X}=Y$ and $\varphi$ contracts $\Delta$ to finitely many points since $X$ is Gorenstein and stable.
    Hence $\varphi$ factors through the canonical map $\varphi_{K}\colon X\to W$, which is finite since $K_{X}$ is ample.
    Thus we complete the proof.
\end{proof}

Horikawa \cite{horikawa} made a similar observation for normal stable Horikawa surfaces with only canonical singularities. 
In this sense, standard Horikawa surfaces can be seen as a straightforward generalization of the surfaces studied in \cite{horikawa}. 

In the non-standard case, the canonical linear system is composed with a pencil, exhibiting fundamentally different behavior compared to the standard case.
In particular, non-standard Horikawa surfaces have the following structure.

\begin{prop} \label{prop:non-std_Horikawa}
Let $X$ be a non-standard Horikawa surface, let $\widetilde{X}$ be the minimal resolution of $X$, and let $n_X$ and $l_{E_h}$ be as in Section \ref{elliptic_case}.
\begin{itemize}
    \item[$(1)$] 
    $|K_{\widetilde{X}}+\widetilde{\Delta}|$ gives an elliptic fibration over a curve $Z$.
    The genus $g:=g(Z)$ is either $0$ or $1$. 
    \item[$(2)$]
    If $X$ has only $\Q$-Gorenstein smoothable singularities, then $q(X)=0$, and exactly one of the following holds:
    \begin{itemize}
        \item[$(2.1)$] 
        $(g,n_X,l_{E_h})=(0,1,0)$ and $\chi(\mathcal{O}_X)=r+l_{R_0}+3$.
        \item[$(2.2)$] 
        $(g,n_X,l_{E_h})=(0,2,0)$ and $\delta'_f=r=l_{R_0}=0$.
        The divisor $\widetilde{C}$ contains exactly two horizontal components, each of which is a section of $\widetilde{X}\to Z$.
        All elliptic singularities of $X$ are obtained from $\widetilde{X}$ by contracting components of elliptic fibers.
        \item[$(2.3)$] 
        $(g,n_X,l_{E_h})=(1,1,1)$ and $\chi(\mathcal{O}_X)=r+l_{R_0}+5$.
        $X$ has exactly one simple elliptic singularity obtained from $\widetilde{X}$ by contracting a section, and all other elliptic singularities of $X$ are obtained by contracting components of elliptic fibers.
    \end{itemize}
\end{itemize}

\end{prop}

\begin{proof}
We first show $(1)$.
Since $K_{X}^{2}=2p_g(X)-4>0$ by definition, either $(1.1)$ or $(1.4)$ in Theorem \ref{normalstable} holds.
In particular, we have an elliptic fibration $\widetilde{X}\to Z$ with $g=0$ or $g=1$.

Suppose that $n_{X}=0$.
Then, Lemma \ref{degbound} implies that $\widetilde{C}$ is empty.
Hence, $X$ becomes a Horikawa surface with only canonical singularities, which is a standard Horikawa surface \cite{horikawa}.
This contradicts the assumption that $X$ is non-standard.
This shows $n_{X}>0$.

From now on, we assume that $X$ has only $\Q$-Gorenstein smoothable singularities.
We show that exactly one of the three cases $(2.1)$, $(2.2)$, and $(2.3)$ holds. 
Note that the case (1.1) or (1.4) in Theorem~\ref{normalstable} holds according to $g=0$ or $g=1$.
In particular, $q(X)=0$.
Applying Theorem~\ref{ellfundeq} to $X$, we obtain the following equality:
\[
(2-n_{X})(\chi(\mathcal{O}_X)-3) + (n_{X}-1)l_{E_h} = (2g+\delta'_{f})n_{X} + r + l_{R_0}.
\]
When $g=0$, the case $(1.1)$ in Theorem \ref{normalstable} implies
that $l_{E_h}=0$.
Since $n_{X} \geq1$ and the right-hand side is non-negative, it holds either $n_{X}=1$ or $(n_{X}, r,l_{R_0}, \delta'_{f})=(2,0,0,0)$.
If $n_{X}=1$, we have $\chi(\mathcal{O}_X)=r+l_{R_0}+3$.
If $n_{X}=2$, we have $\delta'_f=0$ and $r=0$, so it follows from Lemma \ref{lem:no_bisection} that $\widetilde{C}$ must contains two sections.
When $g=1$, the case $(1.4)$ in Theorem \ref{normalstable} holds.
In particular, it implies that $n_{X}=1$ and $l_{E_h}=1$.
Hence we conclude that $\chi(\mathcal{O}_X)=r+l_{R_0}+5$.
\end{proof}

\subsection{Good involution}\label{subsec:good_involution}

In this subsection, we introduce the notion of good involutions on normal stable Horikawa surfaces, and study their properties.
They will play a significant role in the study of $\Q$-Gorenstein smoothable normal stable Horikawa surfaces.

\begin{defn}\label{defn--good-involution}
Let $X$ be a normal stable Horikawa surface and let $\sigma$ be an involution on $X$.
If $\sigma$ satisfies one of the following conditions, then we call $\sigma$ a \emph{good involution}.
\begin{itemize}
    \item 
    The Horikawa surface $X$ is standard.  
     The involution $\sigma$ coincides with the covering transformation of the canonical double covering $\varphi_{K_{X}}$.
    \item 
    The Horikawa surface $X$ is non-standard. 
    The involution $\sigma$ induces a fiberwise involution $\widetilde{\sigma}$ on $\widetilde{X}\to Z$, and the projection $\widetilde{X}/ \widetilde{\sigma} \to Z$ is a ruling over $Z$.
    Furthermore, if the quotient  $W=X/\sigma$ is klt, then $W$ admits a (possibly trivial) P-resolution $W^{\dagger} \to W$ with $K_{W^{\dagger}}^2\geq8$.
\end{itemize}
\end{defn}

Although the definition of a good involution on $X$ differs depending on whether $X$ is standard or non-standard, good involutions share common properties.
The following lemma asserts the uniqueness of a good involution, and describes the good involution on a non-standard Horikawa surface.

\begin{prop}\label{prop:good_involution_uniqueness}
Let $X$ be a normal stable Horikawa surface with only $\Q$-Gorenstein smoothable singularities.
\begin{itemize}
    \item[$(1)$]
    A good involution on $X$ is unique if it exists.
    \item[$(2)$]
    Suppose that $X$ is non-standard and the good involution $\sigma$ on $X$ exists.
    Then, one of the following holds:
    \begin{itemize}
        \item[$(2.1)$]
        $(g,n_X,l_{E_h})=(0,1,0)$ and the involution $\widetilde{\sigma}$ fixes the unique horizontal component of $\widetilde{C}$, which is a section of $\widetilde{X}\to Z$.
        \item[$(2.2)$]
        $(g,n_X,l_{E_h})=(0,2,0)$, $\delta'_f=r=l_{R_0}=0$, and the involution $\widetilde{\sigma}$ interchanges the two horizontal components of $\widetilde{C}$, each of which is a section of $\widetilde{X}\to Z$.
        \item[$(2.3)$]
        $(g,n_X,l_{E_h})=(1,1,1)$ and the involution $\widetilde{\sigma}$ fixes the unique horizontal component of $\widetilde{C}$, which is a section of $\widetilde{X}\to Z$.
    \end{itemize}
\end{itemize}
\end{prop}

\begin{proof}
    We first show $(2)$.
    Since the involution $\widetilde{\sigma}$ induced by $\sigma$, it preserves $\widetilde{C}\subset\widetilde{X}$.
    Due to Proposition \ref{prop:non-std_Horikawa} $(2)$, one sees that $n_X$ is either $1$ or $2$.
    If $n_X=1$, then $\widetilde{C}$ has exactly one horizontal component $S$.
    Since $\widetilde{\sigma}$ acts on $\widetilde{X}$ fiberwise, it must fix $S$.
    Suppose $n_X=2$.
    Again, from Proposition \ref{prop:non-std_Horikawa} $(2.2)$, we see that $\widetilde{C}$ has exactly two horizontal components, denoted by $S_1$ and $S_2$, each of which is a section of $\widetilde{X}\to Z$.
    It will be proved in Proposition \ref{prop:interchange} that $\widetilde{\sigma}$ interchanges $S_1$ and $S_2$.
    This shows $(2)$; see Remark \ref{rem:uniqueness_good_involution}.
    
    We prove the uniqueness of $\sigma$.
    If $X$ is a standard Horikawa surface, then it is trivial by definition.
    Hence, we may assume that $X$ is a non-standard Horikawa surface with a good involution $\sigma$.
    We claim that, in each case listed in $(2)$, there exists a $\widetilde{\sigma}$-fixed section of the elliptic surface $\widetilde{X}\to Z$.
    Note that this implies the uniqueness of $\sigma$ as a non-trivial involution on an elliptic curve with one marking point is unique.
    In cases $(2.1)$ and $(2.3)$, the claim follows immediately.
    In case $(2.2)$, a good involution $\widetilde{\sigma}$ must interchange the two sections contained in $\widetilde{C}$.
    By taking the midpoint of the two sections fiberwise, we obtain a $\widetilde{\sigma}$-fixed section.
    This completes the proof.
\end{proof}

\begin{rem}\label{rem:uniqueness_good_involution}
    In the proof above, we make use of Proposition \ref{prop:interchange}.
    It will be proved by classifying non-standard Horikawa surfaces of type $(2.2)$ in Proposition \ref{prop:non-std_Horikawa}, regardless of good involutions.
    The classification is independent of Proposition \ref{prop:good_involution_uniqueness}, and hence there is no logical issue.
\end{rem}

The following proposition asserts that certain deformation families of a normal stable Horikawa surface $X$ (including $\Q$-Gorenstein smoothing families) induce the good involution on $X$. 
In particular, we see that any $\Q$-Gorenstein smoothable normal stable Horikawa surface must have the good involution.

\begin{prop}\label{prop--involution}
    Let $X$ be a $\mathbb{Q}$-Gorenstein smoothable stable Horikawa surface, and let $f\colon\mathscr{X}\to C$ be a projective and flat morphism with a closed point $0\in C$ such that $C$ is a smooth affine curve and $\mathscr{X}_0=X$.
    Suppose that $\mathscr{X}$ is $\mathbb{Q}$-Gorenstein and $\mathscr{X}_c$ has only Du Val singularities for any $c\in C\setminus\{0\}$.
Then there exists an involution $\sigma_{\mathscr{X}}\in \mathrm{Aut}_C(\mathscr{X})$ that satisfies the following:
\begin{itemize}
    \item The restriction $\sigma_{\mathscr{X}}|_{\mathscr{X}_c}$ coincides with the covering transformation of the canonical map $\varphi_{K_{\mathscr{X}_c}}$ for $c \in C \setminus \{ 0 \}$
    \item The restriction $\sigma_{\mathscr{X}}|_{X}$ acts trivially on $H^0(X,\omega_{X})$.
    \item The restriction $\sigma_{\mathscr{X}}|_{X}$ is a good involution on $X$.
\end{itemize}

\end{prop}

\begin{proof}

Since all fibers are lc, $C \ni c \mapsto h^0(\mathscr{X}_c, \omega_{\mathscr{X}_c}) \in \Z$ is a constant function and the direct image sheaf $f_*\omega_{\mathscr{X}/C}$ is locally free (\cite[Corollary 2.69]{kollar-moduli}).
For each $c \in C \setminus \{0\}$, the fiber $\mathscr{X}_c$ is a Horikawa surface with only Du Val singularities.
Hence, $|\omega_{\mathscr{X}_c}|$ is basepoint-free for any $c \in C \setminus \{0\}$.
Therefore, we obtain the relative canonical morphism
\[
\varphi_{\mathscr{X}^{\circ}/C} \colon \mathscr{X}^\circ := \mathscr{X} \times_C (C \setminus \{0\}) \to \mathbb{P}_{C \setminus \{0\}}(f_*\omega_{\mathscr{X}/C}|_{C \setminus \{0\}})
\]
over $C \setminus \{0\}$.

Let $\mathscr{W}^\circ$ denote the image of $\mathscr{X}^\circ$ under $\varphi_{\mathscr{X}^\circ/C}$.
Then $\varphi_{\mathscr{X}^\circ/C} \colon \mathscr{X}^\circ \to \mathscr{W}^\circ$ is a double covering.
We denote by $\sigma^\circ$ the covering involution of the double covering $\mathscr{X}^\circ \to \mathscr{W}^\circ$.
By construction, the restriction $\sigma^\circ|_{\mathscr{X}_c}$ for $c \in C \setminus \{0\}$ coincides with the covering involution of the canonical morphism of $\mathscr{X}_c$.

Since $\mathcal{M}^{\mathrm{KSBA}}_{2,2p_g(X)-4,0}$ is a proper Deligne-Mumford stack (cf.~\cite[Proposition 2.50]{kollar-moduli}),
the group scheme $\mathrm{Aut}_C(\mathscr{X})$ is finite over $C$.
By the valuative criterion of properness, the involution $\sigma^\circ$ extends to an involution $\sigma_{\mathscr{X}} \in \mathrm{Aut}_C(\mathscr{X})$.
Since the restriction $\sigma_{\mathscr{X}}|_{\mathscr{X}_c}$ for $c \in C \setminus \{0\}$ coincides with the covering involution of the canonical morphism of $\mathscr{X}_c$,
it acts trivially on $H^0(\mathscr{X}_c, \omega_{\mathscr{X}_c})$.
Therefore, $\sigma_{\mathscr{X}}|_{X}$ also acts trivially on $H^0(X, \omega_X)$.

Let $\sigma := \sigma_{\mathscr{X}}|_{X}$ be the involution on $X$.
We will show that this involution $\sigma$ is a good involution.
If $X = \mathscr{X}_0$ is a standard stable Horikawa surface, then $|\omega_{\mathscr{X}_0}|$ is base point free, and hence the morphism $\varphi_{\mathscr{X}^{\circ}/C}$ over $C \setminus \{0\}$ extends to a morphism on $\mathscr{X}$.
It follows that $\sigma$ is a good involution.

Assume that $X = \mathscr{X}_0$ is a non-standard stable Horikawa surface.
First, we show that the quotient $W = X / \sigma$ admits a (possibly trivial) P-resolution $W^{+} \to W$ such that $K_{W^{+}}^2 \geq 8$ if $W$ is klt.
Let $\mathscr{W}:=\mathscr{X}/\sigma_{\mathscr{X}}$ be the quotient of $\mathscr{X}$ by the involution $\sigma_{\mathscr{X}}$ and $\mathscr{W} \to C$ be the induced family.
By the construction of the involution $\sigma$, 
the surface $W$ is a central fiber of the family $\mathscr{W} \to C$.
We also note that the general fiber of the family $\mathscr{W} \to C$ is either $\mathbb{P}^2$ or a Hirzebruch surface.
By \cite[Proposition 5.20]{KoMo} and \cite{kawakita}, it is easy to see that $(\mathscr{W},\frac{1}{2}\mathscr{B}+\mathscr{W}_0)$ is lc. Since the generic fiber of $\mathscr{X} \to C$ has only Du Val singularities, we may assume that $(\mathscr{W},\frac{1}{2}\mathscr{B})$ is a klt pair.
Let $\mathscr{W}^{+} \to \mathscr{W}$ be the log canonical modification of $\mathscr{W}$ by \cite[Theorem 5.41]{kollar-moduli}.
Since $K_{\mathscr{W}^{+}}$ is ample over $\mathscr{W}$, we see that $W^+:=\mathscr{W}_0^{+}$ is klt. 
Therefore, $W^{+} \to W$ is a P-resolution.
Note that the restricted morphism $\mathscr{W}^{+}\setminus W^+   \to \mathscr{W} \setminus W$ is an isomorphism.
Since the general fiber of the family $\mathscr{W}^{+} \to C$ is a $\mathbb{P}^2$ or a Hirzebruch surface, we get $K_{W^{+}}^2 \ge 8$. 
We may assume that $\mathscr{W}^{+}$ is $\Q$-factorial by \cite[Lemma 2]{Manetti} (see also Lemma \ref{lem--lower-semi--rho} below).

Finally, we show that the involution $\widetilde{\sigma}$ on the elliptic surface $\widetilde{X}$, induced by $\sigma$, acts fiberwise, and that the natural fibration $\widetilde{X}/\widetilde{\sigma} \to Z$ is a ruled surface.
Since the elliptic fibration $\widetilde{X}\to Z$ is given by $|K_{\widetilde{X}}+\widetilde{\Delta}|$, and $\widetilde{\sigma}$ acts trivially on $H^{0}(K_{\widetilde{X}}+\widetilde{\Delta})=H^{0}(K_X)$, 
it follows that $\widetilde{\sigma}$ acts fiberwise on $\widetilde{X}$.
Thus, it suffices to show that the branch locus of the quotient map  $\widetilde{X}\to \widetilde{W}:=\widetilde{X}/\widetilde{\sigma}$ contains a horizontal component.
If the horizontal part of $\widetilde{C}$, as in Notation~\ref{note:normal_stable_surface}, consists of a single section, then this section is a horizontal component of the branch locus.
Hence, we may assume that $X$ belongs to the case (2.2) of Proposition \ref{prop:good_involution_uniqueness}.
Suppose, for contradiction, that the branch locus is entirely vertical.
Then the fibration $\widetilde{W}\to Z=\mathbb{P}^1$ is also elliptic.
Note that the horizontal part of $\widetilde{C}$ consists of two sections that are exchanged by $\widetilde{\sigma}$, and their image in the quotient  becomes a section of $\widetilde{W}\to Z$.

We claim that $W$ has only rational singularities and that $\widetilde{W}$ is a rational surface.
Since $W$ is a degeneration of a Hirzebruch surface or $\mathbb{P}^2$, we have $p_g(W)=q(W)=0$, and hence $\chi(\O_{W})=1$.
On the other hand, since $\widetilde{W}$ has only $A_1$-singularites and the map $\widetilde{W}\to W$ is proper and birational, we have $\chi(\O_{\widetilde{W}})\le \chi(\O_{W})$ with equality holding if and only if $W$ has only rational singularities.
From the proof of Theorem~\ref{normalstable}~(1.1), we have $q(\widetilde{X})=0$, and hence $q(\widetilde{W})=0$.
Therefore, we conclude that $\chi(\O_{\widetilde{W}})=1$, $p_g(\widetilde{W})=0$ and $W$ has only rational singularities.

Let $W_{\mathrm{min}}\to Z$ be the relatively minimal model of a resolution of $\widetilde{W}$ over $Z$, and
let $Y'\to W_{\mathrm{min}}$ be the double cover induced by $\widetilde{X}\to \widetilde{W}$.
We claim that the branch locus $B$ of the map $Y'\to W_{\mathrm{min}}$ decomposes into four connected components $B=B_1+B_2+B_3+B_4$, each $B_i$ being a $(-2)$-curve contained in a fiber $F_i$ of type $\mathrm{I}_{2n}$, $\mathrm{I}^{*}_{m}$ or $\mathrm{III}$ for some $0<n\le 4$ and $m\ge 0$.
Since $B$ is a branch divisor which is vertical with respect to the elliptic fibration, there exists a divisor $L$ on $W_{\mathrm{min}}$ such that $B\in |2L|$ and $L^{2}\le 0$.
Moreover, since $Y'\to Z$ has no multiple fibers (due to the existence of sections), $B$ does not contain any entire fiber.
Let $B=\sum_{i=1}^{l}B_i$ be the decomposition into connected components.
If $B_i$ is a chain of $(-2)$-curves of length at least two, then for an end component $D$ of the chain, we have $2L\cdot D=B\cdot D=B_{i}\cdot D=-1$, which is a contradiction.
For similar reasons, each $B_i$ must be an irreducible component of a fiber of type $\mathrm{I}_{2n}$ or $\mathrm{III}$ for some $0<n\le l$ or one of the branch of a fiber of type $\mathrm{I}^{*}_{m}$ for some $m\ge 0$.
Applying the Riemann-Roch theorem and the canonical bundle formula to $W_{\mathrm{min}}$, we obtain
$$
\chi(\O_{Y'})=\chi(\O_{W_{\mathrm{min}}})+\chi(-L)=2+\frac{1}{2}L^{2}=2-\frac{l}{4}.
$$
Since $Y'$ is smooth and birational to $\widetilde{X}$, we have $\chi(\O_{Y'})=\chi(\O_{\widetilde{X}})=1+p_g(\widetilde{X})\ge 1$.
Here, the number $l$ is positive since $W_{\mathrm{min}}$ is a rational surface.
Therefore, we conclude that $l=4$ and that $\widetilde{X}$ is also a rational elliptic surface.
Note that $B$ is not contained in a single fiber of type $\mathrm{I}^{*}_m$, since otherwise the intersection number of $B$ and a section equals one, which contradicts to $B\sim 2L$.

It is straightforward to see that the natural contraction $Y'\to Y_{\mathrm{min}}$ contracts the ramification locus to four isolated fixed points of the involution, and the fiber containing each fixed point is of type $\mathrm{I}_n$, $\mathrm{I}^{*}_{2m}$ or $\mathrm{III}$, depending on whether the corresponding fiber of $W_{\mathrm{min}}$ is of type $\mathrm{I}_{2n}$, $\mathrm{I}^{*}_{m}$ or $\mathrm{III}$.
We now use Theorem~\ref{thm:two-sections}, the classification of non-standard Horikawa surfaces that belong to Proposition~\ref{prop:non-std_Horikawa}~(2.2), which will be proved in Section~\ref{sec:non-std_Horikawa} without taking involutions into account.
Let $F$ be a fiber of $Y_{\mathrm{min}}\to\PP^1$ that corresponds to a one of $Y'\to\PP^1$ containing a ramification locus.
Then Theorem~\ref{thm:two-sections} and the above observation show that one of the following holds:
\begin{itemize}
    \item $F$ is of type $\mathrm{III}$. In this case, the morphism $\widetilde{X}\to Y_{\mathrm{min}}$ is an isomorphism along $F$ and the image of the unique isolated fixed point on $F$ under $\widetilde{X}\to W$ is an $A_1$-singularity.
    \item $F$ is of type $\mathrm{I}^{*}_{2m}$. In this case, the morphism $\widetilde{X}\to Y_{\mathrm{min}}$ is an isomorphism along $F$ and the image of the two isolated fixed point on $F$ under $\widetilde{X}\to W$ is a $D_{l}$-singularity for some $2\le l\le m+4$, where a $D_2$-singularity means two $A_1$-singularities and a $D_3$-singularity means an $A_3$-singularity.
    \item $F$ is of type $\mathrm{I}_{n}$ which does not contribute to $l_{E_v}$.
    In other words, $F$ is not included in $\widetilde{C}$.
    In this case, the morphism $\widetilde{X}\to Y_{\mathrm{min}}$ is an isomorphism along $F$ and the image of the $n$ isolated fixed points on $F$ under $\widetilde{X}\to W$ is an $A_{2n-1}$-singularity.
    \item $F$ is of type $\mathrm{I}_{n}$ which contributes to $l_{E_v}$.
    In this case, the image of the corresponding cusp singularity under $X\to W$ is a cusp singularity of type $[2^{2n-1},l+2]^{\circ}$ for some $l>0$.
\end{itemize}
Note that the last case never occurs since $W$ has no elliptic singularities.
Hence $W$ must have Du Val singularities and the sum of their Milnor numbers is at least $4$.
This contradicts to the fact that $W$ is a degeneration of a Hirzebruch surface or $\mathbb{P}^{2}$.
Indeed, as shown above, there exists a P-resolution $W^{+} \to W$ with $K_{W^{+}}^{2}\ge 8$.
The surface $W^{+}$ violates the formula \cite[Proposition 6.1]{Pro} (see also Lemma~\ref{lem:restrict_involution}).
\end{proof}

\section{Classification of standard Horikawa surfaces}\label{sec:std_Horikawa}
If a normal stable Horikawa surface admits a good involution, then the quotient map $\varphi: X \to W$ gives $X$ the structure of a double cover.
Let $B$ be the branch locus of $\varphi$.
In this setting, the singularities of $X$ can be regarded as singularities of the pair $(W, \tfrac{1}{2}B)$.
Motivated by this perspective, we introduce two classes of singularities: mild singularities and double cone singularities.

The aim of this section is twofold:
\begin{itemize}
    \item To introduce and classify two classes of singularities: mild singularities and double cone singularities.
    \item To classify standard Horikawa surfaces with only $\Q$-Gorenstein smoothable singularities.
\end{itemize}

\subsection{Mild singularity}
\label{sec:mild_singularity}
We introduce mild singularities for log surfaces.

\begin{defn}[Mild singularity]\label{defn--mild--sing}
Let $W$ be a smooth surface, let $L$ be a line bundle on $W$, and let $B\in |L^{\otimes 2}|$ be an effective divisor on $W$.
The triple $(W,\frac{1}{2}B,L)$ is \emph{mild} if the double covering $X$ induced by the triple is lc.
When $L$ is clear from the context, we simply say that the log pair $(W,\frac{1}{2}B)$ is mild.
The corresponding singularity on the double cover $X$ over $W$ branched along $B\in |L^{\otimes 2}|$ is also called a mild singularity.
\end{defn}

Since whether a given triple is mild can be checked locally, it suffices to understand which germs are mild.
The following is a classification of a germ of a mild pair.

\begin{prop}[{\cite[Theorem~3.8]{Al-Pa}}, {\cite[Proposition~2.15]{Anthes}}] \label{lem:classify_mild1}
Let $(0\in W,\frac{1}{2}B)$ be a germ of a mild pair.
Let $W_1 \to W$ be the blow-up at $0\in W$
and let $E_1$ be its exceptional curve.
Let $m:=\mathrm{mult}_{0}B$ and $B_1$ the proper transform of $B$ on $W_1$.
Then the corresponding double covering $X\to W$ has an equivariant $\Q$-Gorenstein smoothing, and
the singularities of $X$ are classified into one of the following:

\begin{itemize}
    \item[$(1)$] 
    $m\le 1$. In this case, $X$ is smooth.

    \item[$(2)$] 
    $m=2$. In this case, $X$ has a Du Val singularity of type $A$.

    \item[$(3)$] 
    $m=3$ and $B_1$ has at worst double points.
    In this case, $X$ has a Du Val singularity of type $D$ or $E$.
    
    \item[$(4)$] 
    $m=3$ and $B_1$ has a triple point.
    Let $W_2\to W_1$ be the blow-up at this point and let $E_2$ be its exceptional curve.
    Then the proper transform $B_2$ of $B_1$ has at worst double points and it contacts with $E_2$ of order at most $2$.
    In this case, $X$ has either a simple elliptic singularity with degree $1$  or a cusp singularity of type $[3,2^k]^{\circ}$.

    \item[$(5)$] 
    $m=4$, $B_1$ has at worst double points and it contacts with $E_1$ of order at most $2$.
    In this case, $X$ has either a simple elliptic singularity with degree $2$  or a cusp singularity of type $[3,2^{k_1},3,2^{k_2}]^{\circ}$ for $k_1 \ge -1$, $k_2 \ge 0$.
\end{itemize}

\end{prop}

\begin{proof}
Let 
$$
\psi=\psi_{1}\circ \cdots \circ \psi_{N}\colon W'=W_{N}\to W_{N-1}\to \cdots \to W_{0}=W
$$
be an even resolution described in Definition~\ref{defn--evenresol} which is a log resolution of the pair $(W, B)$.
We freely use notations in Definition~\ref{defn--evenresol}.
Then, we have 
$$
K_{W'}+\frac{1}{2}B'=\psi^{*}\left(K_{W}+\frac{1}{2}B\right)+\sum_{i=1}^{N}\left(1-\left\lfloor \frac{m_i}{2} \right\rfloor\right)\mathbb{E}_{i},
$$
where $\mathbb{E}_i$ is the total toransform of $E_i$.
By taking the double cover $\varphi'\colon X'\to W'$ banched along $B'$, we have a log resolution $\xi\colon X'\to X$ and 
$$
K_{X'}=\xi^{*}K_{X}+\sum_{i=1}^{N}\left(1-\left\lfloor \frac{m_i}{2} \right\rfloor\right)\varphi'^{*}\mathbb{E}_{i}.
$$
Note that the pullback of the proper transform of $E_{i}$ via $\varphi'$ has geometric multiplicity $1$ (resp.\ $2$) if $m_i$ is even (resp.\ odd).
From the above observations and the log canonicity of $X$, the claim follows by direct computations.
\end{proof}

\subsection{Double cone singularity}
\label{sec:double_cone}
In this subsection, we introduce double cone singularities and provide their classification.
Roughly speaking, a double cone singularity is a lc singularity that admits the structure of a branched double cover over the germ of a cone singularity.

Gorenstein double cone singularities can be various types of elliptic singularities, depending on the branch locus.
In order to describe this precisely, we first introduce the following classes of double points of curves.

\begin{defn}
Let $B$ be a curve on a smooth surface $W$ and let $k \geq 0$ be a non-negative integer.
Then, $B$ has an {\em $A_{k}$-singularity} at $q \in W$ 
if the local equation of $B$ at $q$ is $x^2-y^{k+1}=0$ for some local coordinate $(x,y)$ around $q \in W$.
This definition makes sense when $k\ge 0$.
An $A_{2k-1}$-singularity (resp.\ $A_{2k}$-singularity) is also called a $k$-fold node (resp.\ $k$-fold cusp).
Note that a $k$-fold node (resp.\ $k$-fold cusp) becomes a $(k-1)$-fold node (resp.\ $(k-1)$-fold cusp) after one blow-up.
\end{defn}

\begin{defn}[Double cone singularity] \label{defn:(d;k_1,k_2)_sing}
    Let $0\in W$ be a germ of the cyclic quotient singularity of type $\frac{1}{d}(1,1)$, where $d\ge 2$.
    Let $L$ be a divisorial sheaf on $W$, and let $B\in |L^{[2]}|$ be an effective Weil divisor.
    Assume that the pair $(W, \frac{1}{2}B)$ is lc and consider the double cover $X\to W$ branched along $B\in |L^{[2]}|$.
    Then we say that $X$ has a {\em double cone singularity} over $0\in W$.
    Let $\eta\colon W^{+}\to W$ be the minimal resolution at $0$ and consider the normalized base change $(X\times_W W^{+})^\nu\to W^{+}$ of the double covering $X\to W$ defined by $B\in |L^{[2]}|$.
    Let $\widehat{B}$ denote its branch divisor on $W^{+}$.
    Moreover, we define an effective divisor $B^{+}$ on $W^{+}$ as the log crepant condition:
    $$
    K_{W^{+}}+\frac{1}{2}B^{+}=\mu^{*}\left(K_{W}+\frac{1}{2}B\right).
    $$
    Let $\Delta_{0}$ denote the reduced $\eta$-exceptional curve, which is a $(-d)$-curve.
    Since the pair $(W, \frac{1}{2}B)$ is lc, the multiplicity $m(X)$ of $B^{+}$ along $\Delta_{0}$ belongs to the subset $[0.2]\cap \Q$.
    Moreover, assuming that $X$ is Gorenstein, then $m(X)=0,1$ or $2$.
    \end{defn}

    From now on, we focus on Gorenstein double cone singularities.

    \begin{exam}[Double cone singularity of type $(d;k_1,k_2)$]  
    We say that $(0\in W,\frac{1}{2}B)$ has a \emph{cone singularity of type $(d;k_1,k_2)$} at $0\in W$ if the following hold: 
    \begin{itemize}
        \item 
        $m(X)=2$.  
        \item
        $\widehat{B}$ has two distinct singular points on $\Delta_0$, which are an $A_{(k_1 -1)}$-singularity and an $A_{(k_2 -1)}$-singularity.
        The definition makes sense for $k_i=0$ by replacing an $A_{(k_i-1)}$-singularity with two smooth points of $\widehat{B}$ intersecting transversally with $\Delta_0$.
        \item The both tangent cones of $A_{(k_1 -1)}$-singularity and $A_{(k_2 -1)}$-singularity meet $\Delta_0$ transversally.
    \end{itemize}
    In this case, the corresponding singularity on $X$ is called a {\em double cone singularity of type $(d;k_1,k_2)$}.
    This is a simple elliptic singularity of degree $2d$ ($k_1=k_2=0$) or a cusp singularity of type $[d+2,2^{k_1-1},d+2,2^{k_2-1}]^{\circ}$ ($k_1+k_2>0$).
    \end{exam}

The following is the classification of Gorenstein double cone singularities:

\begin{prop}\label{prop:classification_cone_sing}
    Let $X \to (0\in W, \frac{1}{2}B)$ be a Gorenstein double cone singularity.
    Then, it satisfies one of the following:

    \begin{itemize} 
    \item[$(1)$] $d=2$, $m(X)=1$ and $X\to W$ has isolated branch point at $0$.
    In this case, $X$ is smooth.

    \item[$(2)$] $d=2$ and $X\to W$ is \'{e}tale.
     In this case, $X$ has two $A_1$-singularities.

    \item[$(3)$] $d=2$, $m(X)=1$, $(B^{+}-\Delta_0)\cdot \Delta_0=2$ and $B^{+}-\Delta_0$ has an $A_{k-1}$-singularity on $\Delta_0$.
    In this case, $X$ has a $D_{k+3}$-singularity, where we regard this as an $A_3$-singularity if $k=0$.

    \item[$(4)$] $d=3$, $m(X)=1$ and $(B^{+}-\Delta_0)\cdot \Delta_{0}=1$.
    In this case, $X$ has an $A_2$-singularity.
    
    \item[$(5)$] $d=4$, $m(X)=1$ and $X\to W$ has isolated branch point at $0$.
    In this case, $X$ has an $A_1$-singularity.
    
    \item[$(6)$] $m(X)=2$, $(B^{+}-2\Delta_0)\cdot \Delta_0=4$ and $B^{+}-2\Delta_0$ has an $A_{k_1-1}$-singularity and an $A_{k_2-1}$-singularity on $\Delta_0$ for some $k_1, k_2$.
    In this case, $X$ has a double cone singularity of type $(d;k_1,k_2)$.  
    \end{itemize}
\end{prop}

\begin{proof}
we consider the even resolution of $\widehat{B}$ described in Definition~\ref{defn--evenresol} as follows:
\[
\psi:W'=W_{N} \xrightarrow{\psi_{N}} W_{N-1} \xrightarrow{\psi_{N-1}} \cdots \xrightarrow{\psi_{2}} W_{1} \xrightarrow{\psi_{1}} W_{0} := W^{+},
\]
where $\psi_{i}:W_{i} \to W_{i-1}$ is the blow-up at a point $p_{i} \in W_{i-1}$.
Let $E_{i}$ be the exceptional divisor of $\psi_{i}$ and let $m_{i}=\mathrm{mult}_{p_{i}}B_{i-1}$.
We have the following diagram:
\[
\xymatrix{
X' \ar[r]^{\tau} \ar[d]^{\varphi'}& X^{+}  \ar[r]^{\mu} \ar[d]_{\varphi^{+}} &  X \ar[d]^{\varphi} 	\\
W' \ar[r]_{\psi} & W^{+}  \ar[r]_{\eta} &  W  \\
}	
\]
Then, we have
\begin{align*}
K_{W'}=\psi^{\ast}K_{W^{+}} + \sum_{i=1}
^{N}\mathbb{E}_{i}, \quad 
B'=\psi^{\ast}\widehat{B} - \sum_{i=1}
^{N}2 \left\lfloor \frac{m_{i}}{2} \right\rfloor \mathbb{E}_{i},
\end{align*}
where $\mathbb{E}_{i}$ is the pullback of $E_{i}$ on $W'$.
Therefore, it holds that
\[
K_{W'}+\frac{1}{2}B' =\psi^{\ast}\left( K_{W^{+}}+\frac{1}{2}\widehat{B} \right)-\sum_{i=1}
^{N} \left(  \left\lfloor \frac{m_{i}}{2} \right\rfloor-1 \right) \mathbb{E}_{i}.
\]
Pulling this back this via $\varphi'$, we obtain
\[
K_{X'} =\tau^{\ast}K_{X^{+}}- \varphi'^{\ast}\left(\sum_{i=1}
^{N} \left(  \left\lfloor \frac{m_{i}}{2} \right\rfloor -1 \right) \mathbb{E}_{i}\right).
\]

On the other hand, if we write $B^{+}=\widehat{B}+a\Delta_{0}$ for some $a\in \Q$, we have
$$
K_{X^{+}}=\varphi^{+*}\left(K_{W^{+}}+\frac{1}{2}\widehat{B}\right)=\varphi^{+*}\left(\eta^{*}\left(K_W+\frac{1}{2}B\right)-\frac{a}{2}\Delta_0\right)=\mu^{*}K_X-\frac{a}{2}\varphi^{+*}\Delta_0.
$$
Substituting this to the above formula, we get
\[
K_{X'} =(\mu \circ \tau)^{\ast}K_{X}- \varphi'^{\ast}\left(
 \frac{a}{2} \psi^{\ast}\Delta_{0}
+\sum_{i=1}^{N} \left(  \left\lfloor \frac{m_{i}}{2} \right\rfloor -1 \right) \mathbb{E}_{i}\right).
\]
Since $X$ is Gorenstein and lc, we have $a\in 2\Z$ and $a\le 2$ (resp.\ $a\in \Z$ and $a\le 1$) if $\Delta_{0}$ is not contained in $\widehat{B}$ (resp.\ $\Delta_{0}$ is contained in $\widehat{B}$).

Moreover, since $(K_{W^+}+\frac{1}{2}B^{+})\cdot \Delta_0=0$ and $K_{W^{+}}\cdot \Delta_{0}=d-2$, we have $\widehat{B}\cdot \Delta_0=4-(2-a)d$.
If $\Delta_{0}$ is not contained in $\widehat{B}$ (resp.\ $\Delta_{0}$ is contained in $\widehat{B}$), we have $0\le \widehat{B}\cdot \Delta_0=4-(2-a)d$ (resp.\ $0\le (\widehat{B}-\Delta_0)\cdot \Delta_0=4-(1-a)d$). 

First, suppose that $a=2$.
Then $\Delta_{0}$ is not contained in $\widehat{B}$ and so $m(X)=2$.
Since the proper transformation of $E_{i}$ is contained in $\psi^{\ast}\Delta_{0}$,
it holds that
\[
1 \geq 
\mathrm{ord}_{E_{i}} \left(
\psi^{\ast}\Delta_{0}
+\sum_{i=1}^{N} \left(  \left\lfloor \frac{m_{i}}{2} \right\rfloor -1 \right) \mathbb{E}_{i} \right) \geq \left\lfloor \frac{m_{i}}{2} \right\rfloor,
\]
where the first inequality follows from the log-canonicity of $X$.
Therefore, $m_{i} \leq 3$ for any $i$.
If $m_i=3$, then the proper transform of $E_{i}$ is contained in the branch locus $B'$,
which contradicts the assumption that $X$ is lc.
Thus, we have $m_{i} \leq 2$ for any $i$.
Therefore, we conclude that $(W,\frac{1}{2}B)$ has a cone singularity of type $(d; k_1,k_2)$ at $0 \in W$ for some $k_{1}, k_{2} \in \Z_{\geq 0}$.

Next, we suppose that that $a=1$.
Then, $\Delta_0$ is contained in $\widehat{B}$.
Since the proper transformation of $E_{i}$ is contained in $\psi^{\ast}\Delta_{0}$ and $X$ is Gorenstein lc,
we have $m_{i} \leq 3$ and $E_i$ is contained in the branch locus, that is, $m_i=3$ for each $i$.
However, since $\psi$ is an even resolution, $m_1=3$ implies that $m_i=2$ for some $i$.
Hence, $\widehat{B}$ has no singularities, which is a contradiction since $\widehat{B}-\Delta_0$ intersects $\Delta_0$.

Finally, suppose that $a\le 0$.
Take the minimal resolution $\widetilde{X}\to X^{+}$.

Assume also that the composition $\widetilde{X}\to X^{+}\to X$ is not a minimal resolution.
Then, the proper transform of $\varphi^{+*}\Delta_0$ contains a $(-1)$-curve.
If $\Delta_0\subset \widehat{B}$, then $\varphi^{+*}\Delta_0=2E$ for an irreducible and reduced curve $E$ with $E^2=-d/2$.
Hence, $-d/2\ge -1$ holds.
This implies $d=2$ and hence $\widetilde{X}\to X^{+}$ is an isomorphism.
Since $E$ is a $(-1)$-curve, $X$ is smooth.
This case implies (1).
If $\Delta_0$ is not contained in $\widehat{B}$, then  
$\varphi^{+*}\Delta_0$ is reduced and at most two irreducible components.
Hence, $-2d=(\varphi^{+*}\Delta_0)^2\ge -2$, which contradicts to $d\ge 2$.

Assume that the composition $\widetilde{X}\to X^{+}\to X$ is minimal.
This condition implies that $a=0$ and $\mu$ is crepant.
Then we conclude that $x\in X$ is Du Val.
Indeed, otherwise, any exceptional curve on $\widetilde{X}$ would have a negative discrepancy, which is a contradiction.
If $\Delta_0$ is not contained in $\widehat{B}$, then $0\le 4-2d$ implies that $d=2$ and $X\to W$ is \'{e}tale.
This case implies (2).
If $\Delta_0$ is contained in $\widehat{B}$, then $0\le 4-d$ implies $d=2,3$ or $4$.
It is straightforward to see that the case $d=2$ (resp.\ $3$, $4$) implies (3) (resp.\ (4), (5)).
\end{proof}

\begin{rem}
Gorenstein double cone singularities that appear on standard Horikawa surfaces were already classified in {\cite[Theorem~4.1]{Che}} by using linear systems on minimal degree surfaces.
The proof of Proposition~\ref{prop:classification_cone_sing} relies only on local arguments.
Since such singularities occur not only on standard Horikawa surfaces but also those of Lee-park type,
Proposition~\ref{prop:classification_cone_sing} is also applicable to Horikawa surfaces of Lee-park type.
\end{rem}

\subsection{Result}
\label{sec:classification_standard_Horikawa}
Let $X$ be a standard Horikawa surface.
According to Proposition~\ref{prop:standard_Horikawa}, its canonical map $\varphi_{K_{X}}\colon X\to W$ is a double covering over a minimal degree surface $W \subset \mathbb{P}^{p_{g}-1}$.
Let $B$ denote its branch divisor on $W$.

The following theorem gives the classification of standard Horikawa surfaces, which is an extension of Horikawa's work \cite[Theorem 1.6]{horikawa}.

\begin{thm}[{\cite[Theorem~4.1]{Che}}]\label{thm:standard_Horikawa}
    Let $X$ be a standard Horikawa surface, and let $W$ and $B$ be as above.
    Then, $q(X)=0$ and one of the following holds.
\begin{itemize}
    \setlength{\itemindent}{24pt}
    \item[$\operatorname{Type}$ $(\infty)$]
    $W=\mathbb{P}^2$, $p_g(X)\in\{3,6\}$, $B\in |\mathcal{O}_{\mathbb{P}^2}(\frac{2}{3}p_g(X)+6)|$, and
     a log pair $(W,\frac{1}{2}B)$ is mild.
    In particular, $X$ has only mild singularities.
    \item[$\operatorname{Type}$ $(d)$\ ] $W=\Sigma_d$ and $B\in |6\Delta_0+(p_g(X)+3d+2)\Gamma|$, where $d\in\mathbb{Z}_{\ge0}$, $d\leq\min\left\{p_g(X)-4,\frac{1}{2}p_g(X)+1\right\}$, $p_g(X)-d$ is even and
     a log pair $(W,\frac{1}{2}B)$ is mild.
    In particular, $X$ has only mild singularities.
    \item[$\operatorname{Type}$ $(d)'$] $W=\overline{\Sigma}_d$, $d=p_g(X)-2\ge 2$, and $B \in |(4p_{g}(X)-4)\overline{\Gamma}|$.  
    $X$ has mild singularities and a Gorenstein double cone singularity as in classified in Proposition~\ref{prop:classification_cone_sing}~(2)--(6).
\end{itemize}
    
\end{thm}
\begin{proof}
According to Proposition~\ref{prop:standard_Horikawa}, the canonical map $\varphi_{|K_{X}|}:X \to W$ is a double covering of the minimal degree surface $W \subset \mathbb{P}^{p_g(X) -1}$.
The classification of minimal degree surfaces \cite{Nagata} implies that one of the following holds:
\begin{itemize}
    \item
    $(W,\mathcal{O}_{W}(1))=(\mathbb{P}^2,\mathcal{O}_{\mathbb{P}^2}(1))$ and $p_{g}=3$. 
    \item
$(W,\mathcal{O}_{W}(1))=(\mathbb{P}^2,\mathcal{O}_{\mathbb{P}^2}(2))$  and $p_{g}=6$.
    \item
$(W,\mathcal{O}_{W}(1))=(\Sigma_{d},\mathcal{O}_{\Sigma_{d}}(\Delta_{0}+\frac{p_g-2+d}{2} \Gamma))$, $d\leq p_g(X)-4$, and $p_g(X)-d$ is even.
    \item
$(W,\mathcal{O}_{W}(1))=(\overline{\Sigma}_{p_{g}-2},\mathcal{O}_{\overline{\Sigma}_{p_{g}-2}}((p_g-2)\overline{\Gamma}))$, and $d=p_g(X)-2\ge 2$.
\end{itemize}
Since $\mathcal{O}_{W}(1) \cong \mathcal{O}_{W}(K_{W}+ \frac{1}{2}B)$, the branch divisor $B$ belongs to the desired linear system.
Furthermore, if $W=\Sigma_d$, then the inequality $p_g(X) \geq 2d-2$ holds, since otherwise $B-2\Delta_0 \geq 0$, contradicting the normality of $X$.

Finally, we show $q(X)=0$.
Since $X$ is a double covering of $W$ with the branch divisor $B\in |L^{[2]}|$,
we have $\chi(\O_{X})=\chi(\O_W)+\chi(L^{[-1]})$.
Since $W$ is rational with only rational singularities, $\chi(\O_{W})=1$ holds.
Thus, it suffices to show that $\chi(L^{[-1]})=p_g(X)$.
If $W$ is smooth, this follows from the Riemann-Roch theorem.
Assume $W=\overline{\Sigma}_{p_g-2}$ and let
$\widehat{\psi}\colon \widehat{W}=\Sigma_{p_g-2}\to W$ be the minimal resolution.
Put 
$$
\widehat{L}:=\O_{\widehat{W}}(3\Delta_{0}+(2p_g-2)\Gamma).
$$
Then $\widehat{\psi}_{*}\widehat{L}^{-1}=L^{[-1]}$ and $R^{1}\widehat{\psi}_{*}\widehat{L}^{-1}=R^{1}\widehat{\psi}_{*}\omega_{\widehat{W}}(-\Delta_{0}+(3p_g-2)\Gamma)=0$ hold from the Kawamata-Viehweg vanishing theorem. 
This implies that $\chi(L^{[-1]})=\chi(\widehat{L}^{-1})$.
Applying the Riemann-Roch theorem again, this equals $p_g(X)$.
\end{proof}

\section{Classification of non-standard Horikawa surfaces}\label{sec:non-std_Horikawa}

\subsection{Overview}\label{subsec:overview_non-std}
The purpose of this section is to classify non-standard Horikawa surfaces with only $\Q$-Gorenstein smoothable singularities and those with a good involution (Definition \ref{defn--good-involution}).

Throughout this section, we use Notation \ref{note:normal_stable_surface}.
Recall that Proposition \ref{prop:non-std_Horikawa} provides a rough classification of non-standard Horikawa surfaces with only $\Q$-Gorenstein smoothable singularities into the following three types:
\begin{itemize}
\item $(g,n_X,l_{E_h})=(0,1,0)$ (Section \ref{subsec:0,1,0}).

\item $(g,n_X,l_{E_h})=(0,2,0)$, $\delta'_{f}=r=l_{R_0}=0$ (Section \ref{subsec:0,2,0}).

\item $(g,n_X,l_{E_h})=(1,1,1)$ (Section \ref{subsec:1,1,1}).
    
\end{itemize}
In this section, we classify the minimal resolution $\widetilde{X}$ and the divisor $\widetilde{C}$ for each case.
Furthermore, we identify those surfaces that admit a good involution (Section \ref{subsec:0,1,0,involution}, \ref{subsec:0,2,0,involution}, \ref{subsec:1,1,1,involution}), including an explicit construction (Construction \ref{construction}).
The result will be provided in Section \ref{subsec:non-standard_Horikawa}.
Since the existence of a good involution is a necessary condition for a non-standard Horikawa surface $X$ to be smoothable (Proposition \ref{prop--involution}), this classification provides candidates for $\Q$-Gorenstein smoothable $X$.
The smoothablity and deformations of non-standard Horikawa surfaces with good involutions will be studied in Section \ref{sec:deformation}.

Since the arguments in the cases that $(g,n_X,l_{E_h})=(0,1,0)$ and $(g,n_X,l_{E_h})=(1,1,1)$ remain valid even when $X$ is not a Horikawa surface,  we will give classifications for a broader class than non-standard Horikawa surfaces. 
More precisely, we treat a surface $X$ satisfying the following.
\begin{ass}\label{ass:normal_stable_elliptic}
    The surface $X$ is a normal stable surface with only $\Q$-Gorenstein smoothable singularities, and its minimal resolution admits an elliptic fibration over a smooth curve of genus $g$.
    Furthermore, $\chi(\mathcal{O}_X)\geq2$.
\end{ass}

We introduce the following notations and definitions.

\begin{note}
\label{notation_7}
    \phantom{A}
    \begin{enumerate}[label=$(\arabic*)$]
        \item 
        We take $Y$ in Notation \ref{note:normal_stable_surface} to be the relatively minimal elliptic fibration.
        \item 
        We decompose $C$ as follows:
        \begin{equation}\label{eqn:decomp}
            C=C_\hor + \sum_{j=1}^J D_j,
        \end{equation}
        where $C_\hor$ is the union of the horizontal components, and $D_{j}$'s are contained in pairwise disjoint fibers $F_{j}$.
        By abuse of notation, the proper transforms of $F_j$ are also denoted by $F_j$.
        \item 
        For $0\leq j\leq J$ and $0\leq k\leq m$, let $D_{j,k}$ denote the maximal subdivisor of $C_k$ whose support is contained in the pullback of $F_j$ by $Y_k \to X$.
        \item 
        For $0\leq j\leq J$ and $0\leq k\leq m$, let $B_{j,k}$ denote the connected components of $D_{j,k}$ that intersect with $C_\hor$.
        For $0\leq k\leq m$, let $B_k$ denote the connected components of $C_k$ that contain $C_\hor$.
        \item 
        We set $\chi:=\chi(\O_X)$.
        \item 
        When $n_{X} = 1$, the horizontal component $C_{\mathrm{hor}}=:S$ consists of an irreducible section. Let $D'_{j}$ denote the irreducible component of $D_{j}$ intersecting $S$, let $x_{j}$ denote the intersection point of $D_{j}$ and $S$.
        We write $p\in X$ for the singularity whose exceptional set $C^{(p)}$ contains $S$.
    \end{enumerate}
\end{note}

\begin{defn}
\phantom{A}
\begin{enumerate} 
\item
For a reduced normal crossing divisor $E=\sum_{i}E_i$ on a proper smooth surface,
an irreducible component $E_j$ is said to be a {\em fork} if $(E-E_j)\cdot E_j\ge 3$.
In this case, each connected component of $E-E_j$ is called a {\em branch} with respect to $E_j$.

\item 
Let $E\subset C_k$ be a subdivisor on $Y_k$.
For $0\le l<k$, let $E_l\subset C_l$ denote the maximal subdivisor  whose image on $Y_k$ is contained in $E$.
We say that $E$ {\em becomes $E_l$ after the blow-ups} $Y_l\to Y_k$.
Note that if a connected component of $C_k$ becomes a chain $E_l$ after the blow-ups, then
$E_l$ is an extended T-chain with the ample condition, since $E_0$ is a disjoint union of chains that are connected components of $C_0$, which are nothing but T-chains.
\end{enumerate}
\end{defn}

\begin{defn}[decorated graph]\label{defn:graph}
\phantom{A}
\begin{itemize}
    \item[$(1)$] A \emph{decorated graph} $G$ consists of the following triple: a graph $\Gamma$, a map $V(\Gamma)\to\Z$, and a decomposition of $V(\Gamma)$ into three subsets, where $V(\Gamma)$ is the set of vertices of $\Gamma$.
    \item[$(2)$] We define the decorated graph $G_{j,k}$ by the following data:
    \begin{itemize}
        \item the dual graph (also denoted by $G_{j,k}$) of the union of the total transform of $D_{j}$ under $Y_k \to X$ and the horizontal part $C_\hor$;
        \item the map $V(G_{j,k})\to\Z$ which assigns $-E^2$ to the vertex of $G_{j,k}$ corresponding to a vertical curve $E$, and assigns to the vertex corresponding to a component $S\subset C_\hor$ the number of blow-ups occurring on (the proper transform of) $S$ over the fiber $F_j$;
        \item the decomposition of the vertex set
        \[
        V(G_{j,k})=\{C \mid C \subset D_{j,k}\} \sqcup \{C \mid C \subset \overline{F_j\setminus D_{j,k}}\} \sqcup \{S\mid S\subset C_\hor \}.
        \]
    \end{itemize}
    When illustrating the graph $G_{j,k}$, vertices in $\{C \mid C \subset D_{j,k}\}$, $\{C \mid C \subset \overline{F_j\setminus D_{j,k}}\}$ and $\{S\mid S\subset C_\hor\}$ are represented by the symbols $\circ$, $\bullet$ and $\diamond$, respectively. 
    The assigned numbers are displayed next to the corresponding vertices.
    \end{itemize}
\end{defn}

To simplify the graph $G_{j,k}$, we will use the following abbreviated notation.
\begin{itemize}
    \item If $G_{j,k}$ contains a chain consisting of $\beta$ vertices with all assigned $2$, then we abbreviate it as 
\[\xygraph{
    \circ ([]!{-(0,-.3)} {2}) - [r]
      \circ ([]!{-(0,-.3)} {2 }) - [r]
       \cdots ([]!{-(0,-.3)} {}) - [r]
       \circ ([]!{-(0,-.3)} {2}) - [r]
        \circ ([]!{-(0,-.3)} {2})
        }
= \xygraph{
    \circ ([]!{-(0,-.3)} {2^{\beta}})}.
    \] 
\end{itemize}

\subsection{Case of $(g,n_X,l_{E_h})=(0,1,0)$}\label{subsec:0,1,0}
In this subsection, we classify $C$ and its resolution $\widetilde{C}$ assuming $(g,n_X,l_{E_h})=(0,1,0)$ under Assumption \ref{ass:normal_stable_elliptic}.
One sees that $S^2=-\chi(\O_Y)$ by the adjunction formula.
By abuse of notation, we also denote the proper transform by $S$.
We will proceed with the classification in the following steps:
\begin{enumerate}[label=\textbf{(Step~\arabic*)}]
    \setlength{\itemindent}{20pt}
    \item We classify blow-ups that may occur during the initial steps for each $D_j$, referred to as the {\it initial blow-up process}. 
    In fact, if we allow reordering of the blow-ups, the initial blow-up process can be classified into finite patterns based on the following condition: each connected component of $D_{j,0}$, except for the one intersecting $S$, consists of the exceptional curves over a $\Q$-Gorenstein smoothable lc singularity. 
    Furthermore, the component that intersects $S$, if it exists, must be a subset of the exceptional curves over $p$.
    For a precise statement, refer to Proposition \ref{amulet}.
    \item We classify the blow-up processes that give a branch with respect to $S$ that is completely separated from $S$.
    Thanks to the classification in Step $1$, it is sufficient to examine only a finite number of cases. 
    The classification is given in Proposition \ref{compsepamulet}.
    \item We classify the configurations of curves in $C$ and their resolutions when $p$ is a T-singularity. 
    See Theorem \ref{thm:one-section-T} for details.
    \item We classify the configurations of curves in $C$ and their resolutions when $p$ is a strictly lc rational singularity. For more information, refer to Theorems \ref{thm:one-section-halfcusp} and \ref{thm:one-section-triangle}.
\end{enumerate}

\subsubsection{Step~1: initial blow-up process}\label{subsubsec:0,1,0-step1}
In this step, we classify all possible initial blow-up processes.
The main tools are the classification of singular fibers of elliptic fibrations and the list of $\Q$-Gorenstein smoothable singularities (Section \ref{subsec:lc_sing}, Appendix \ref{app:cusp}).
The result of this step is the following proposition.  
The notations used in the statement are introduced in Lemmas~\ref{chain}-\ref{modifyIV*} below.

\begin{prop} \label{amulet}
    Assume that $(g,n_X,l_{E_h})=(0,1,0)$ and Assumption \ref{ass:normal_stable_elliptic}.
    If $D_{j,0}$ does not intersect with $S$, set $l=1$ and assume that $D_{j,1}$ intersects with $S$ by reordering the blow-ups; otherwise, set $l=0$.
    Then, exactly one of the following holds:
    \begin{itemize}
        \item[$(1)$] 
        $S+B_{j,l}$ is a chain of smooth rational curves and each connected component of the effective divisor $D_{j,l}-B_{j,l}$ is a T-chain. 
        By reordering the blow-ups if necessary, there exists some $k$ such that no blow-ups occur on $D_{j,k}-B_{j,k}$, and the graph $G_{j,k}$ coincides with exactly one of the following:
        \[
        S- [\mathrm{CT}]_\beta\ (\beta\geq0),\quad S-[\mathrm{I}_n \mathrm{T}]_\beta\ (n\geq1, \beta\geq0),\quad S-[\mathrm{IIT1}], 
        \]
        \[
        S-[\mathrm{IIT2}],\quad S-[\mathrm{IIT3}],\quad S-[\mathrm{IIT4}],\quad S-[\mathrm{IIIT1}],
        \]
        \[
        S-[\mathrm{IIIT2}],\quad S-[\mathrm{IIIT3}],\quad S-[\mathrm{IVT1}],\quad S-[\mathrm{IVT2}].
        \]
        
        \item[$(2)$] 
        $S+B_{j,l}$ is a chain of smooth rational curves, the effective divisor $D_{j,l}-B_{j,l}$ consists of T-chains and exceptional curves of strictly lc rational singularities, and it contains at least one of the latter.
        By reordering the blow-ups if necessary, there exists some $k$ such that no blow-ups occur on the effective divisor $D_{j,k}-B_{j,k}$, and the graph $G_{j,k}$ coincides with exactly one of the following:
        \[
        S-[\mathrm{I}_n \mathrm{R}]_\beta\ (n\geq3, \beta=0,1),\quad S-[\mathrm{IIR}],\quad S-[\mathrm{IIIR}]_\beta\ (\beta=0,1),
        \]
        \[
        S-[\mathrm{IVR}]_\beta\ (\beta=0,1,2),\quad S-[\mathrm{I}_n^* \mathrm{R}]_\beta\ (n\geq1, 0\leq\beta\leq n+3).
        \]
        
        \item[$(3)$] 
        $S+B_{j,l}$ has a fork but no cycles, and $D_{j,l}=B_{j,l}$ or the effective divisor $D_{j,l}-B_{j,l}$ consists of T-chains.
        By reordering the blow-ups if necessary, there exists some $k$ such that no blow-ups occur on the effective divisor $D_{j,k}-B_{j,k}$, and the graph $G_{j,k}$ coincides with exactly one of the following:
        \[
        S-[\mathrm{CSR}],\quad S-[\mathrm{I}_n \mathrm{SR1}]\ (\chi=3, n\geq3),\quad
        S-[\mathrm{I}_n \mathrm{SR2}]\ (\chi=4, n\geq3),
        \]
        \[
        \quad S-[\mathrm{I}_n^*\mathrm{SR1}]\ (n\geq0), \quad
        S-[\mathrm{I}_n^*\mathrm{SR2}]\ (n\geq1),\quad S-[\mathrm{II}^* \mathrm{SR}],\quad S-[\mathrm{III}^*\mathrm{SR}],\quad S-[\mathrm{IV}^*\mathrm{SR}].
        \]

        \item[$(4)$] 
        $D_{j,0}$ contains the exceptional curves of an elliptic singularity.
        By reordering the blow-ups if necessary, there exists some $k$ such that no blow-ups occur on the effective divisor $D_{j,k}-B_{j,k}$, and the graph $G_{j,k}$ coincides with exactly one of the following:
        \[
        S-[\mathrm{I}_n \mathrm{E}]_\beta\ (n\geq0, 0\leq\beta\leq n+8).
        \]
    \end{itemize}
\end{prop}

Proposition \ref{amulet} follows immediately from the subsequent nine lemmas (and their proofs).
We consider the following three cases for $D_j$ and classify possible initial blow-up processes in each case:
\begin{itemize}
    \item 
    $D_{j}$ forms a chain. (Lemma \ref{chain})
    \item 
    $D_{j}$ has a cycle.
    In this case, $D_{j}$ coincides with the support of $F_{j}$ and the fiber $F_{j}$ is of type $\mathrm{I}_{n}$ (Lemma \ref{modifyI}), $\mathrm{II}$ (Lemma \ref{modifyII}), $\mathrm{III}$ (Lemma \ref{modifyIII}) or $\mathrm{IV}$ (Lemma \ref{modifyIV}).
    \item 
    $D_{j}$ has a fork.
    In this case, the fiber $F_{j}$ is of type $\mathrm{I}^{*}_{n}$ (Lemma \ref{modifyI*}), $\mathrm{II}^{*}$ (Lemma \ref{modifyII*}), $\mathrm{III}^{*}$ (Lemma \ref{modifyIII*}) or $\mathrm{IV}^{*}$ (Lemma \ref{modifyIV*}).
\end{itemize}

We will frequently use the following lemma.

\begin{lem} \label{strlcratbranch}
Let $D$ be an irreducible component of $C_k$ for some $k$, and let $B\subset C_{k}$ be a connected component of $C_{k}-D$ intersecting $D$ transversally at one point $y$. 
Then we have the following:
\begin{itemize}
    \item[$(1)$] If $B$ consists of a single $(-b)$-curve and is completely separated from $D$, then it becomes a T-chain $[2^{b-3}, b+1]$ after the blow-ups. 
    If $B$ is a chain of $(-2)$-curves, then $B$ is not completely separated from $D$.
\end{itemize}
Moreover, if $D$ is a fork of a strictly lc rational singularity $q$, then the following hold:
\begin{itemize}
    \item[$(2)$] $B$ is not a chain of $(-2)$-curves of length greater than $1$.
    \item[$(3)$] If $B$ consists of a single $(-b)$-curve and not completely separated from $D$, then there are no blow-ups over $B$.
    \end{itemize}
\end{lem}

\begin{proof}
We prove (1).  
Assume that $B$ is a single $(-b)$-curve.  
If $B$ is completely separated from $D$, then $B$ becomes the chain $L_2^{\beta}L_1[b] = [2^{\beta}, b+1]$ for some $\beta \geq 0$.  
Since $[2^{\beta}, b+1]$ must be P-admissible, it follows that $\beta = b - 3$.  
Assume next that $B$ is a chain of $(-2)$-curves of length $l \geq 2$.  
If $B$ is completely separated from $D$, then after the blow-ups, $B$ becomes a chain of the form $[2^{\beta}, 3, 2^{l-1}]$ as a P-admissible chain, which contradicts Lemma~\ref{lem:extT_2m2}.  
Thus, $B$ is not completely separated from $D$.

From now on, we assume that $D$ is a fork of a strictly lc rational singularity $q$.
Assume that $B$ is a chain of $(-2)$-curves of length $l \geq 2$. 
Since $q$ is strictly lc rational as in Lemma~\ref{lem:smoothable_rational_strict_lc}, at least one blow-up occurs over $B$.
Furthermore, by (1), at least one blow-up occurs at $y$.
Hence, after all blow-ups over $y$, the branch $B$ becomes a P-admissible chain
\[
L_{2}^{r-2}L_1L_{2}^{\beta-1}L_1[2^{l}]=[2^{r-2},3,2^{\beta-2},3,2^{l-1}]
\]
where $-r$ is the self-intersection number of the resulting branch of $\pi^{-1}(q)$.
This implies that $r=2$ and $l=1$, which contradicts the assumption that $l\geq2$.
Thus, (2) is proved.

We assume that $B$ is a $(-b)$-curve and that $B$ is not completely separated from $D$.  
Assume to the contrary that the blow-up at $y$ occurs.  
After all the blow-ups over $y$, the branch $B$ becomes a P-admissible chain  
\[
L_{2}^{r-2}L_1L_{2}^{\beta-1}L_1[b] = [2^{r-2}, 3, 2^{\beta-2}, b+1]
\]
where $-r$ is the self-intersection number of the resulting branch $E_1$ of $\pi^{-1}(q)$.  
By computing the log discrepancies of $E_1$ and the leftmost component $E_2$ of T-trains associated with the above chain, one sees that this configuration violates the ampleness of $K_X$. 
Indeed, the $(-1)$-curve $E$ connecting $E_1$ and $E_2$ satisfies $K_{X} \cdot \pi(E) < 0$.  
Therefore, we have shown (3).
\end{proof}

We first deal with the case that $D_j$ forms a chain.

\begin{lem} \label{chain}
If $D_j$ is a chain of length $\beta(\geq1)$, then no blow-ups occur on $D_j$ and the dual graph $G_{j,0}$ is one of the following:

$(\mathrm{T})\;
S- [\mathrm{CT}]_\beta:=\xygraph{
    \diamond ([]!{-(0,-.3)} {+0}) - [r]
    \circ ([]!{-(0,-.3)} {2^{\beta}}) }
\quad
(\mathrm{SR})\; 
    S-[\mathrm{CSR}]
    :=\xygraph{
    \diamond ([]!{-(0,-.3)} {+0}) - [r]
    \circ ([]!{-(0,-.3)} {2}) (
        - [d] \circ ([]!{-(-0.25,-.0)} {2}),
        - [r] \circ ([]!{-(0,-.3)} {2}),
}
$

\noindent
In the latter case, the fiber $F_j$ is of type $\mathrm{I}_n$ with $n\geq4$.
\end{lem}

\begin{proof}
Let $q\in X$ be the singularity such that $D'_{j}\subset C^{(q)}$.
Since $D_j$ does not contain an elliptic curve and a cycle of rational curves, $q$ is a rational singularity.
If $D_j-D'_{j}$ is connected, then the dual graph $G_{j,m}$ coincides with the graph described in the case $(\mathrm{T})$.
Indeed, no blow-ups occur on $D_j$ in this case; if $q$ is a T-singularity, it follows from Corollary \ref{cor_drill_2}, and if $q$ is a strictly lc rational singularity, it follows from Lemma \ref{strlcratbranch}.

From now on, we assume that $D_j-D'_{j}$ is not connected.
In this case, this consists of two disjoint chains of $(-2)$-curves.
We show that the length of each chain must be one.

Suppose that $q$ is a T-singularity.
Then at least one of the three branches of the fork $D'_{j}$ is completely separated from $D'_{j}$.
If a branch not containing $S$ is completely separated after the blow-ups, then this branch becomes a P-admissible chain of the form $[2^{k},3,2^{l}]$ ($k\ge 0, l\ge 0$).
This contradicts Lemma~\ref{lem:extT_2m2}.
If the branch containing $S$ is completely separated after the blow-ups, the chain $D_j$ becomes a P-admissible chain of the form $[2^{k},n,2^l]$,
where $k, l\ge 1$ and $n\ge 3$.
This also contradicts Lemma~\ref{lem:extT_2m2}.

Suppose that $q$ is strictly lc rational.
Applying Lemma~\ref{strlcratbranch} to the fork $D'_{j}$,
each connected component of $D_j-D'_j$ consists of a single $(-2)$-curve, and there are no blow-ups over $D_j$.
Hence, this is the case $(\mathrm{SR})$.
\end{proof}

\begin{lem} \label{modifyI}
Assume that $D_j$ has a cycle and $F_j$ is of type $\mathrm{I}_{n}$.
By reordering the blow-ups if necessary, there exists some $k$ such that the graph $G_{j,k}$ coincides with one of the following:

$(\mathrm{T})\;S-[\mathrm{I}_n \mathrm{T}]_\beta\;
:=\xygraph{
    \diamond ([]!{-(0,-.3)} {+0}) - [r]
    \circ ([]!{-(0,-.3)} {3+\beta}) (
        - [rd] \bullet ([]!{-(.3,0)}  {1}, -[r] \circ ),
        - [r] \circ ([]!{-(0,-.3)} {2^{n-2}})
        - [r] \circ ([]!{-(0,-.3)} {3})
        - [d] \circ ([]!{-(-0.3,-.0)} {2^{\beta}})
}
$

\noindent
where $n\geq1$ and $\beta\geq0$.
Here, if $n=1$, the subchain $[3+\beta, 2^{n-2}, 3]$ is regarded as $[4+\beta]$.

$(\mathrm{R})\;
S-[\mathrm{I}_n \mathrm{R}]_{\beta}\;
:=\xygraph{
    \diamond ([]!{-(0,-.3)} {+1}) - [r]
     \circ ([]!{-(0,-.3)} {2^{\beta-1}}) - [r]
     \circ ([]!{-(0,-.3)} {3}) - [r]
      \circ ([]!{-(0,-.3)} {2}) - [r]
      \bullet ([]!{-(0,-.3)} {1}) - [r]
       \circ ([]!{-(0,-.3)} {3}) - [r]
    \circ ([]!{-(-0.5,-.0)} {3+\beta}) ( 
        - [u] \circ ([]!{-(.3,0)} {3}) - [r] \bullet ([]!{-(-0,-.3)} {1}) - [r] \circ ([]!{-(-0.3,-.0)} {3}) - [d] \circ ([]!{-(-0.5,-.0)} {2^{n-5}}),
        - [d] \circ ([]!{-(.3,0)} {3})- [r] \bullet ([]!{-(-0,-.3)} {1}) - [r] \circ ([]!{-(-0.3,-.0)} {3}) - [u] \circ ([]!{-(.3,0)} {}),       
)}
$
where $\beta=0, 1$ and $n\ge 3$.
Here, if $n=3$ or $n=4$, the subchain $-1-[3,2^{n-5},3]-1-$ is regarded as $-1-$ or $-1-[4]-1-$, respectively.
If $\beta=0$, then the subchain 
$\xygraph{
    \diamond ([]!{-(0,-.3)} {+1}) - [r]
     \circ ([]!{-(0,-.3)} {2^{\beta-1}}) - [r]
     \circ ([]!{-(0,-.3)} {3}) - [r]
      \circ ([]!{-(0,-.3)} {2}) - [r]
      \bullet ([]!{-(0,-.3)} {1})}$ is regarded as 
$\xygraph{
\diamond ([]!{-(0,-.3)} {+2}) - [r]
\circ ([]!{-(0,-.3)} {2}) - [r]
\bullet ([]!{-(0,-.3)} {1})}$ .

$(\mathrm{SR1})\;S-[\mathrm{I}_n \mathrm{SR1}]:=
\xygraph{
    \diamond ([]!{-(0,-.3)} {+0}) - [r]
    \circ ([]!{-(-0.3,-.0)} {2}) ( 
        - [u] \circ ([]!{-(.3,0)} {3}) - [r] \bullet ([]!{-(-0,-.3)} {1}) - [r] \circ ([]!{-(-0.3,-.0)} {3}) - [d] \circ ([]!{-(-0.5,-.0)} {2^{n-5}}),
        - [d] \circ ([]!{-(.3,0)} {3})- [r] \bullet ([]!{-(-0,-.3)} {1}) - [r] \circ ([]!{-(-0.3,-.0)} {3}) - [u] \circ ([]!{-(.3,0)} {}),      )}
$

where $\chi=3$ and $n\ge 3$.
Here, if $n=3$ or $n=4$, the subchain $-1-[3,2^{n-5},3]-1-$ is regarded as $-1-$ or $-1-[4]-1-$, respectively.

$(\mathrm{SR}2)\;S-[\mathrm{I}_n \mathrm{SR2}]:=
\xygraph{
    \diamond ([]!{-(0,-.3)} {+0}) - [r]
    \circ ([]!{-(-0.3,-.0)} {3}) ( 
        - [u] \circ ([]!{-(.3,0)} {4})
        - [r] \bullet ([]!{-(-0,-.3)} {1}) %
        - [r] \circ ([]!{-(-0,-.3)} {2}) %
        - [r] \circ ([]!{-(-0.3,-.0)} {3}) 
        - [d] ,
        - [d] \circ ([]!{-(.3,0)} {2})
        - [r] \bullet ([]!{-(-0,-.3)} {1})
        - [r] \circ ([]!{-(-0,-.3)} {4}) %
        - [r] \circ ([]!{-(-0.5,-.0)} {2^{n-4}}) 
        - [u]  ([]!{-(.3,0)} {}),      )}
$

where $\chi=4$ and $n\ge 3$.

$(\mathrm{E})\;S-[\mathrm{I}_n \mathrm{E}]_\beta:=
\xygraph{
    \diamond ([]!{-(0,-.3)} {+1}) 
    - [r] \circ ([]!{-(0,-.3)} {2^{\beta}})
    - [r] \bullet ([]!{-(0,-.3)} {1})
     - [r] \circ ([]!{-(-0.5,-.0)} {3+\beta})( 
        - [u] \circ ([]!{-(.3,0)} {2}) - [r] \circ ([]!{-(-0,-.3)} {2}) - [r] \circ ([]!{-(-0.3,-.0)} {2}) - [d] \circ ([]!{-(-0.5,-.0)} {2^{n-7}}),
        - [d] \circ ([]!{-(.3,0)} {2})- [r] \circ ([]!{-(-0,-.3)} {2}) - [r] \circ ([]!{-(-0.3,-.0)} {2}) - [u] \circ ([]!{-(.3,0)} {}),      )}
$

where $n\geq0$ and $0\leq\beta\leq n+8$.
Here, if $n = 0$, the cycle is regarded as a smooth elliptic curve with self-intersection number $-(\beta+1)$, and, if $n = 1$, the cycle is regarded as a rational curve with self-intersection number $-(\beta+1)$ that has exactly one node.
\end{lem}

\begin{proof}
For simplicity, we assume that $n\ge 2$. The cases of $n = 0$ and $n = 1$ are similar.

We first assume that $D_{j}'$ is not a fork of exceptional curves over any strictly lc rational singularity.
To resolve the fork, at least one of three branches connected to $D'_{j}$ is completely separated from $D'_{j}$ after the blow-ups.
If an irreducible component next to $D_j'$ in $F_j$ is completely separated from $D_j'$ after $\beta+1$ blow-ups, it gives the case $(\mathrm{T})$.

Suppose that $S$ is completely separated from $D'_{j}$ after $\beta+1$ times blow-ups. 
After the blow-ups, $D_{j}$ becomes a cycle consisting of one $(-3-\beta)$-curve which is the proper transform of $D'_{j}$ and $n-1$ $(-2)$-curves.
If no further blow-ups occur on this cycle, this gives a case $(\mathrm{E})$.
The upper bound of $\beta$ can be obtained from Corollary \ref{cor:cusp1}. 
We now assume that this cycle becomes a chain after some blow-ups (at some node).
After the blow-ups, this cycle becomes a P-admissible chain
$$
L_{i_m}\cdots L_{i_1}L_1[2^{k},3+\beta,2^{n-1-k}]R_1R_{j_1}\cdots R_{j_m}
$$
for some $k$ and $i_1,\ldots,i_m, j_1,\ldots, j_m$ with $\{i_p, j_p\}=\{1,2\}$ for each $p=1,\ldots,m$.
One sees from a simple calculation that this is not possible.

Assume that $D_{j}'$ is a fork of exceptional curves over some strictly lc rational singularity $q \in X$.
Let $z_1$ and $z_2$ be two intersection points of $D'_{j}$ and $D_j-D'_{j}$.
By Lemma~\ref{lem:smoothable_rational_strict_lc}, we may assume that the branches of $\pi^{-1}(q)$ over $z_1$ and $z_2$ are a $(-r_1)$-curve and a $(-r_2)$-curve for some $r_1, r_2 \in \{2,3,4,6\}$, respectively.
If blow-ups occur over both $z_1$ and $z_2$, then after blow-ups over $z_1$ and $z_2$, $D_j-D'_j$ must become an extended T-chain
\begin{align}
\label{align:Iext_T}
L_2^{r_1-2}L_1L_2^{\beta_1}L_1[2^{n-1}]R_1R_2^{\beta_2}R_1R_2^{r_2-2}=[2^{r_1-2},3,2^{\beta_1-1},3,2^{n-3},3,2^{\beta_2-1},3,2^{r_2-2}].
\end{align}
If the exceptional divisor of $q\in X$ is of the form described in Lemma~\ref{lem:smoothable_rational_strict_lc} (i), then $r_1 = r_2 = 2$.  
If it is as in Lemma~\ref{lem:smoothable_rational_strict_lc} (ii)--(iv),  
then the singularity $q\in X$ must be of type $(3,3,3)[4]$.  
Indeed, when the extended T-chain (\ref{align:Iext_T}) is formed in the fiber, the proper transform of $D'_j$ becomes a $-(\beta_1 + \beta_2 + 4)$-curve.  
According to the classification in Lemma~\ref{lem:smoothable_rational_strict_lc}, only the singularity of type $(3,3,3)[4]$ is allowed.
Therefore, we have $r_1 = r_2 = 3$ and $\beta_1 = \beta_2 = 0$. 

Consequently, P-admissible chain~\eqref{align:Iext_T} must be of one of the following forms:
\[
[3,2^{\beta_1-1},3,2^{n-3},3,2^{\beta_2-1},3], \quad [2,4,2^{n-3},4,2].
\]
We show that each case is ruled out by the ampleness of $K_X$.
By Lemma~\ref{lem:extT_I}~(1), the former chain is not P-admissible.
By Lemma~\ref{lem:extT_I}~(2), the chain $C=[2,4,2^{n-3},4,2]$ is a P-admissible chain and its associated T-train is $[2,5]-1-[3,2^{n-5},3]-1-[5,2]$.
In this case, the component $D'_j$ becomes the fork of the exceptional divisor of the singularity of type $(3,3,3)[4]$, and $\widetilde{X}$ has the following configuration of curves:
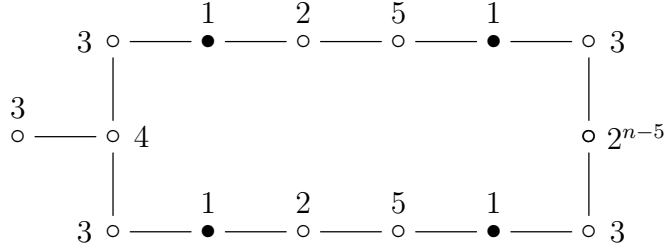
\[
\xygraph{
    \circ ([]!{-(0,-.3)} {3}) - [r]
    \circ ([]!{-(-0.3,-.0)} {4}) ( 
        - [u] \circ ([]!{-(.3,0)} {3})
        - [r] \bullet ([]!{-(-0,-.3)} {1}) %
        - [r] \circ ([]!{-(-0,-.3)} {2}) %
        - [r] \circ ([]!{-(-0,-.3)} {5}) %
        - [r] \bullet ([]!{-(-0,-.3)} {1}) 
        - [r] \circ ([]!{-(-0.3,-.0)} {3}) 
        - [d] \circ ([]!{-(-0.5,-.0)} {2^{n-5}}),
        - [d] \circ ([]!{-(.3,0)} {3})
        - [r] \bullet ([]!{-(-0,-.3)} {1})
        - [r] \circ ([]!{-(-0,-.3)} {2}) %
        - [r] \circ ([]!{-(-0,-.3)} {5}) %
        - [r] \bullet ([]!{-(-0,-.3)} {1}) %
        - [r] \circ ([]!{-(-0.3,-.0)} {3}) 
        - [u] \circ ([]!{-(.3,0)} {}),      )}
\]
Note that the leftmost vertex represents the proper transform of the section $S$.
The log discrepancy on $X$ with respect to each component of the branch of the exceptional divisor over the singularity of type $(3,3,3)[4]$ is $1/3$,  
while the log discrepancy on $X$ of the $(-2)$-curve in the T-chain $[2,5]$ is $2/3$.  
Therefore, if $E$ denotes a $(-1)$-curve connecting $(3,3,3)[4]$ and $[2,5]$, then we have $\pi(E) \cdot K_X = 0$.
This contradicts the ampleness of $K_X$.

Hence we may assume that there are no blow-ups at $z_2$.
In particular, $r_2\ge 3$.
Let $D'_{j2}$ denote the irreducible component of $D_{j}-D'_{j}$ passing through $z_{2}$, and let $z'_{2}$ denote the intersection point of $D'_{j2}$ and $D-D'_{j}-D'_{j2}$.
If blow-ups over $z_1$ occur, then after the blow-ups over $z_1$ and $z'_2$, 
we get an extended T-chain with the ample condition
\begin{align*}
L_2^{r_1-2}L_1L_2^{\beta_1}L_1[2^{n-2}]R_1R_2^{r_2-3}=[2^{r_1-2},3,2^{\beta_1-1},3,2^{n-4},3,2^{r_2-3}].
\end{align*}
By Lemma~\ref{lem:extT_I}~$(4)$, we must have either $(r_1,r_2)=(3,3)$ or $(r_1,r_2)=(2,4)$.

If $(r_1,r_2)=(3,3)$, then the associated T-train is $[2,3,2^{\beta_1 -1},4]-1-[3,2^{n-5},3]$.
In this case, $\widetilde{X}$ has the following configuration of curves:
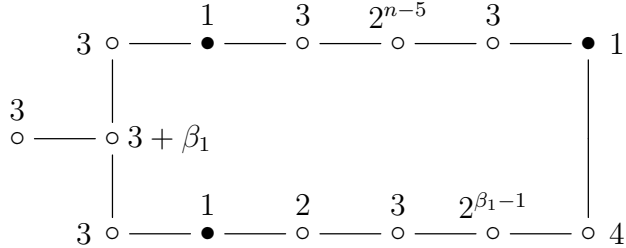
\[
\xygraph{
    \circ ([]!{-(0,-.3)} {3}) - [r]
    \circ ([]!{-(-0.6,-.0)} {3+\beta_1}) ( 
        - [u] \circ ([]!{-(.3,0)} {3})
        - [r] \bullet ([]!{-(-0,-.3)} {1}) %
        - [r] \circ ([]!{-(-0,-.3)} {3}) %
        - [r] \circ ([]!{-(-0,-.3)} {2^{n-5}}) %
        - [r] \circ ([]!{-(-0,-.3)} {3}) 
        - [r] \bullet ([]!{-(-0.3,-.0)} {1}) 
        - [d] ,
        - [d] \circ ([]!{-(.3,0)} {3})
        - [r] \bullet ([]!{-(-0,-.3)} {1})
        - [r] \circ ([]!{-(-0,-.3)} {2}) %
        - [r] \circ ([]!{-(-0,-.3)} {3}) %
        - [r] \circ ([]!{-(-0,-.3)} {2^{\beta_1-1}}) %
        - [r] \circ ([]!{-(-0.3,-.0)} {4}) 
        - [u]  ([]!{-(.3,0)} {}),      )}
\]
If we write $E$ for the $(-1)$-curve connecting $(3,3,3)[3+\beta_1]$
and the chain $[2,3,2^{\beta_1-1},4]$, then $K_{X}\cdot \pi(E)=0$. 
It contradicts the ampleness of $K_{X}$.

If $(r_1,r_2)=(2,4)$, then the associated T-train is $[3,2^{\beta_1-2},3]-1-[4,2^{n-3},3,2]$, and the self-intersection number of the proper transform of $D'_j$ on $\widetilde{X}$ is $-(3+\beta_1)$.
Lemma \ref{lem:smoothable_rational_strict_lc} shows that $\beta_1=0$, and $\widetilde{X}$ has the following configuration of curves:
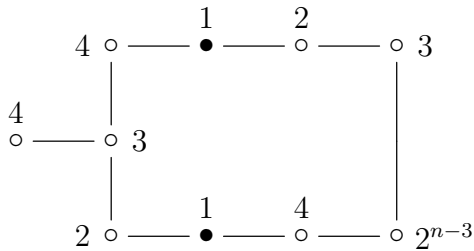
\[
\xygraph{
    \circ ([]!{-(0,-.3)} {4}) - [r]
    \circ ([]!{-(-0.3,-.0)} {3}) ( 
        - [u] \circ ([]!{-(.3,0)} {4})
        - [r] \bullet ([]!{-(-0,-.3)} {1}) %
        - [r] \circ ([]!{-(-0,-.3)} {2}) %
        - [r] \circ ([]!{-(-0.3,-.0)} {3}) 
        - [d] ,
        - [d] \circ ([]!{-(.3,0)} {2})
        - [r] \bullet ([]!{-(-0,-.3)} {1})
        - [r] \circ ([]!{-(-0,-.3)} {4}) %
        - [r] \circ ([]!{-(-0.5,-.0)} {2^{n-3}}) 
        - [u]  ([]!{-(.3,0)} {}),      )}
\]
Since there exists no blow-ups over the intersection point of $D'_j$ and $S$, it becomes a case $(\mathrm{SR}2)$.

Finally, we assume there are no blow-ups at $z_1$ and $z_2$.
In particular, $n\ge 3$ and $r_1, r_2\ge 3$.
After appropriate blow-ups, we get a connected component of $C_k$ of the form
$$
L_2^{r_1-3}L_1[2^{n-3}]R_1R_2^{r_2-3}=[2^{r_1-3},3,2^{n-5},3,2^{r_2-3}].
$$
Since this must be P-admissible, it holds that $r_1=r_2=3$ by Lemma \ref{lem:extT_I}. If there exist blow-ups over the intersection point of $D'_{j}$ and $S$, then it becomes the case $(\mathrm{R})$, while otherwise it becomes the case $(\mathrm{SR}1)$.

\end{proof}

\begin{lem} \label{modifyII}
Assume that $D_{j}$ has a cycle and $F_{j}$ is of type $\mathrm{II}$.
By reordering the blow-ups if necessary, there exists some $k$ such that the graph $G_{j,k}$ coincides with one of the following:

$(\mathrm{T1})\;
 S-[\mathrm{IIT1}] 
    :=\xygraph{!~:{@{=}}
    \diamond ([]!{-(0,-.3)} {+0}) - [r]
      \circ ([]!{-(0,-.3)} {4}) :  [r]
       \bullet ([]!{-(0,-.3)} {1}) 
}
\quad \quad \;\;\; \;
(\mathrm{T2})\;
   S-[\mathrm{IIT2}] 
    :=\xygraph{
    \diamond ([]!{-(0,-.3)} {+0}) - [r]
      \circ ([]!{-(0,-.3)} {5}) 
      (- [d]   \bullet ([]!{-(-0.25,-.0)}  {1} ),
       - [r]  \circ ([]!{-(0,-.3)} {2}) 
       -[ld]  \circ ([]!{-(0,-.3)} 
       }
$

$
(\mathrm{T3})\;
    S-[\mathrm{IIT3}]
    :=\xygraph{
    \diamond ([]!{-(0,-.3)} {+0}) - [r]
    \circ ([]!{-(0,-.3)} {6}) (
        - [d] \circ ([]!{-(.3,0)} {2}) ( 
         - [r] \bullet ([]!{-(0,-.3)} {1}) 
        - [d] \circ ([]!{-(.3,0)} {4}),
         - [d] \circ ([]!{-(.3,0)} {2})
    )}
 (\mathrm{T4})\;
    S-[\mathrm{IIT4}]:=
    \xygraph{
    \diamond ([]!{-(0,-.3)} {+0}) - [r]
    \circ ([]!{-(0,-.3)} {7}) (
        - [d] \circ ([]!{-(.3,0)} {2^3}) (
         - [r] \bullet ([]!{-(0,-.3)} {1}) 
         - [r] \circ ([]!{-(-0.25,-.0)} {5}) (
        - [d] \circ ([]!{-(-0.25,-.0)} {3}),
        - [u] \circ ([]!{-(-0.25,-.0)} {2}))
    )}
$

$
(\mathrm{R})\;
S-[\mathrm{IIR}]:=
    \xygraph{
    \diamond ([]!{-(0,-.3)} {+0}) - [r]
      \circ ([]!{-(0,-.3)} {8}) - [r]
       \circ ([]!{-(0,-.3)} {2^4}) - [r]
        \bullet ([]!{-(0,-.3)} {1}) - [r]
    \circ ([]!{-(0,-.3)} {6}) - [r]
    \circ ([]!{-(-0.25,-.0)} {2}) (
        - [d] \circ ([]!{-(-0.25,-.0)} {3}),
        - [u] \circ ([]!{-(-0.25,-.0)} {2}),
    }
$

\end{lem}

\begin{proof}
Note that a blow-up at the cusp of $F_j$ must occur since $\widetilde{C}$ has only nodes as singularities.
Hence we may assume that $\rho_{m}\colon Y_{m-1}\to Y_{m}=Y$ is a blow-up at the cusp $y_{m}$.
If $E_{m}$ is not contained in $C_{m-1}$, this is nothing but the case $(\mathrm{T1})$.
Assume that $E_{m}$ is contained in $C_{m-1}$.
Since the intersection point (say $y_{m-1}$) of $E_{m}$ and the proper transform of $D'_{j}$ is not a node, it must be blown up.
If the exceptional curve $E_{m-1}$ of the blow-up at $y_{m-1}$ is not contained in $C_{m-2}$, this is the case $(\mathrm{T2})$.
We assume that $E_{m-1}$ is contained in $C_{m-2}$.
Then $C_{m-2}$ has the triple point $y_{m-2}$.
The triple point $y_{m-2}$ must be blown up.
Moreover, the exceptional curve $E_{m-2}$ of the blow-up at $y_{m-2}$ must be contained in $C_{m-3}$.
Indeed, otherwise $\widetilde{C}$ has a $(-3)$-curve disjoint from other components and it gives a singularity of $X$ that is not $\Q$-Gorenstein smoothable, which is a contradiction.

We first assume that the proper transform of $E_{m-2}$ is not a fork in $\widetilde{C}$.
It must be a component of a T-chain after the blow-ups.
Hence, at least one of the three branches is completely separated from $E_{m-2}$.
From Lemma~\ref{strlcratbranch} $(1)$, the proper transform of $E_{m-1}$, which is a $(-2)$-curve, cannot be completely separated from $E_{m-2}$.
Assume that the proper transform of $E_{m}$ is completely separated from $E_{m-2}$.
By Lemma~\ref{strlcratbranch} $(1)$, we obtain the case $(\mathrm{T3})$.
We assume that the branch containing the proper transform of $D_{j}$ is completely separated from $E_{m-2}$.
After the blow-ups over the intersection point (say $y_{m-3}$) of $D_{j,m-3}$ and $E_{m-2}$, the proper transforms of $E_{m-1}$, $E_{m-2}$ and $E_{m}$ should form a P-admissible chain 
$$
[3,2+\beta,2]
$$ 
for some $\beta\ge 0$.
It follows by a simple computation that either $\beta=3$ and it is a T-chain or $\beta=2$ and it becomes $[4]-1-[5,2]$ after one blow-up.
The former case leads to the case $(\mathrm{T4})$. 
The latter case is excluded from the condition that $K_X$ is ample.
Indeed, let $E$ be the $(-1)$-curve obtained by the last blow-up over the point $y_{m-3}$,
and let $E'$ be the $(-2)$-curve next to $E$.
Then one see that the log discrepancy of $E'$ with respect to $X$ is is greater than $2/3$ since $E'$ becomes a curve at the end of a T-chain of the form $[2,2,\ldots]$.
While, the proper transform of $E_{m-2}$ has log discrepancy $1/3$.
Hence $K_{X}\cdot \pi(E)<0$, which is a contradiction.

Assume that the proper transform of $E_{m-2}$ is a fork of exceptional curves over some strictly lc singularity $q$.
Since $E_{m-2}$ has exactly three branches, any branch is not completely separated from $E_{m-2}$ after the blow-ups.
By Lemma~\ref{strlcratbranch}~(3), the two branches which are the proper transforms of $E_{m-1}$ and $E_{m}$ are not blown up.
Hence the singularity $q$ is of type $(2,3,6)[2]$ and there are blow-ups over the intersection point of $E_{m-2}$ and the proper transform of $D_{j}$.
By calculating the blow-ups that are certain to happen, we obtain the case $(\mathrm{R})$.
\end{proof}

\begin{lem} \label{modifyIII}
Assume that $D_{j}$ has a cycle and $F_{j}$ is of type $\mathrm{III}$.
By reordering the blow-ups if necessary, there exists some $k$ such that the graph $G_{j,k}$ coincides with one of the following:

$(\mathrm{T1})\;
 S-[\mathrm{IIIT1}]:= 
    \xygraph{
    \diamond ([]!{-(0,-.3)} {+0}) - [r]
      \circ ([]!{-(0,-.3)} {3}) 
      (- [d]   \bullet ([]!{-(-0.25,-.0)}  {1} ),
       - [r]  \circ ([]!{-(0,-.3)} {3}) 
       -[ld]  \circ ([]!{-(0,-.3)} 
       }
  $
$(\mathrm{T2})\;
    S-[\mathrm{IIIT2}]:=
    \xygraph{
    \diamond ([]!{-(0,-.3)} {+0}) - [r]
    \circ ([]!{-(-0.25,-.0)} {4}) (
        - [d] \circ ([]!{-(.3,0)} {3}) ( 
        - [r] \bullet ([]!{-(-0,-.3)} {1})
         - [r] \circ ([]!{-(0,-.3)} {2}) 
        - [d] \circ ([]!{-(.3,0)} {5}),
         - [d] \circ ([]!{-(.3,0)} {2})
    )}$
$(\mathrm{T3})\;
    S-[\mathrm{IIIT3}]:=
     \xygraph{
    \diamond ([]!{-(0,-.3)} {+0}) - [r]
    \circ ([]!{-(0,-.3)} {5}) (
        - [r] \circ ([]!{-(0,-.3)} {2}) (
         - [d] \bullet ([]!{-(-0.25,-.0)} {1}) 
         - [l] \circ ([]!{-(0,-.3)} {3}) (
        - [d] \circ ([]!{-(-0.25,-.0)} {2}),
        - [l] \circ ([]!{-(0,-.3)} {4}))
    )}
$

$(\mathrm{R})\;
    S-[\mathrm{IIIR}]_{\beta}:=
    \xygraph{
    \diamond ([]!{-(0,-.3)} {+0}) - [r]
      \circ ([]!{-(0,-.3)} {5}) - [r]
       \circ ([]!{-(0,-.3)} {2^{\beta-1}}) - [r]
      \circ ([]!{-(0,-.3)} {3}) - [r]
        \circ ([]!{-(0,-.3)} {2}) - [r]
         \circ ([]!{-(0,-.3)} {2}) - [r]
        \bullet ([]!{-(0,-.3)} {1}) - [r]
    \circ ([]!{-(0,-.3)} {4}) - [r]
    \circ ([]!{-(-0.75,-.0)} {(\beta+2)}) (
        - [d] \circ ([]!{-(-0.25,-.0)} {2}),
        - [u] \circ ([]!{-(-0.25,-.0)} {4}),
    }
$
where $\beta=0,1$.

\end{lem}

\begin{proof}
Note that the tacnode of $F_j$ must be blown up since $\widetilde{C}$ has only nodes as singularities.
We may assume that $\rho_{m}$ is a blow-up at the tacnode $y_{m}$.
If the exceptional curve $E_{m}$ is not contained in $C_{m-1}$, then it is the case $(\mathrm{T1})$.
We assume that $E_{m}$ is contained in $C_{m-1}$.
Then $C_{m-1}$ has the triple point $y_{m-1}$.
It must be blown up and the exceptional curve $E_{m-1}$ is contained in $C_{m-2}$ since $\widetilde{C}$ does not contain a $(-2)$-curve as its connected component.

We first assume that the proper transform of $E_{m-1}$ is not a fork of exceptional curves over any strictly lc rational singularity.
Hence, at least one of three branches is completely separated from $E_{m-1}$.
The proper transform of $E_{m}$, that is a $(-2)$-curve, can not be completely separated from $E_{m-1}$ after the blow-ups by Lemma~\ref{strlcratbranch}~(1).
If the proper transform of $D_j-D'_{j}$,
that is a $(-4)$-curve, is completely separated from $E_{m-1}$, then we obtain the case $(\mathrm{T2})$ from Lemma~\ref{strlcratbranch}~(1).
We assume that the branch containing the proper transform of $D'_{j}$ is completely separated from $E_{m-1}$.
After the blow-ups over the intersection point of $D'_{j,m-2}$ and $E_{m-1}$,
the proper transforms of $E_{m}$, $E_{m-1}$ and $D_j-D'_{j}$ form a P-admissible chain
\[
[2,2+\beta,4]
\]
for some $\beta\ge 0$.
It follows by a simple computation that $\beta=1$ and it is a T-chain, which is the case $(\mathrm{T3})$.

Assume that the proper transform of $E_{m-1}$ is a fork of exceptional curves over some strictly lc rational singularity $q$.
By the same argument as in the proof of Lemma~\ref{modifyII},
the singularity $q$ is of type $(2,4,4)[\beta+2]$ with $\beta=0$ or $1$. 
After the blow-ups over the intersection point of $E_{m-1}$ and the proper transform of $D'_j$, we obtain the case $(\mathrm{R})$.
\end{proof}

\begin{lem} 
\label{modifyIV}
Assume that $D_{j}$ has a cycle and $F_{j}$ is of type $\mathrm{IV}$.
By reordering the blow-ups if necessary, there exists some $k$ such that the graph $G_{j,k}$ coincides with one of the following:

$(\mathrm{T1})\;
S-[\mathrm{IVT1}]:=
 \xygraph{
    \diamond ([]!{-(0,-.3)} {+0}) - [r]
    \circ ([]!{-(0,-.3)} {3}) (
        - [d] \circ ([]!{-(.3,0)} {2}) ( 
         - [r] \bullet ([]!{-(0,-.3)} {1}) 
        - [d] \circ ([]!{-(-0.25,-.0)} {4}),
         - [d] \circ ([]!{-(-0.25,-.0)} {3})
)}
$
$(\mathrm{T2})\;
S-[\mathrm{IVT2}]:=
\xygraph{
    \diamond ([]!{-(0,-.3)} {+0}) - [r]
      \circ ([]!{-(0,-.3)} {4}) - [d]
        \bullet ([]!{-(.3,0)} {1}) - [r]
    \circ ([]!{-(-0.25,-.0)} {2}) (
        - [d] \circ ([]!{-(-0.25,-.0)} {3}),
        - [u] \circ ([]!{-(-0.25,-.0)} {3}),
}
$

$(\mathrm{R})\;
S-[\mathrm{IVR}]_{\beta}:=
 \xygraph{
    \diamond ([]!{-(0,-.3)} {+0}) - [r]
      \circ ([]!{-(0,-.3)} {4}) - [r]
       \circ ([]!{-(0,-.3)} {2^{\beta-1}}) - [r]
      \circ ([]!{-(0,-.3)} {3}) - [r]
        \circ ([]!{-(0,-.3)} {2}) - [r]
        \bullet ([]!{-(0,-.3)} {1}) - [r]
    \circ ([]!{-(0,-.3)} {3}) - [r]
    \circ ([]!{-(-0.75,-.0)} {(\beta+2)}) (
        - [d] \circ ([]!{-(-0.25,-.0)} {3}),
        - [u] \circ ([]!{-(-0.25,-.0)} {3}),
}
$

\noindent
where $\beta=0,1,2$.
\end{lem}

\begin{proof}
We may assume that $\rho_{m}$ is the blow-up at the triple point $y_{m}$ of $F_j$.
Then the exceptional curve $E_{m}$ must be contained in $C_{m-1}$ since $\widetilde{C}$ has no $(-3)$-curves as connected components.

We first assume that the proper transform of $E_{m}$ is not a fork of exceptional curves over any strictly lc rational singularity.
Then at least one of three branches is completely separated from $E_{m}$.
If the proper transform of an irreducible component of $D_j-D'_{j}$, that is a $(-3)$-curve, is completely separated from $E_{m}$, then we obtain the case $(\mathrm{T1})$ by Lemma~\ref{strlcratbranch}~(1).
We now assume that the branch containing the proper transform of $D'_{j}$ is completely separated from $E_{m}$.
Then the proper transform of $E_{m}$ and the two components of $D_j-D'_{j}$ must form a P-admissible chain
$$
[3,2+\beta,3]
$$
for some $\beta\ge 0$.
It follows by a simple computation that $\beta=0$ and it is a T-chain.
This is the case $(\mathrm{T2})$.

We assume that the proper transform of $E_{m}$ is a fork of exceptional curves over some strictly lc rational singularity $q$.
By the same argument as in the proof of Lemma~\ref{modifyII}, it follows that the singularity $q$ is of type $(3,3,3)[\beta+2]$ with $\beta=0,1$ or $2$. 
After the blow-ups over the intersection point of $E_{m}$ and the proper transform of $D'_j$, we obtain the case $(\mathrm{R})$. 
\end{proof}

\begin{lem} \label{modifyI*}
Assume that $D_{j}$ has a fork and $F_{j}$ is of type $\mathrm{I}_n^*$.
By reordering the blow-ups if necessary, there exists some $k$ such that the graph $G_{j,k}$ coincides with one of the following:

$(\mathrm{R})\;S-[\mathrm{I}_{n}^{*}\mathrm{R}]_{\beta}:=
 \xygraph{
    \diamond ([]!{-(0,-.3)} {+0}) - [r]
      \circ ([]!{-(0,-.3)} {3}) - [r]
       \circ ([]!{-(0,-.3)} {2^{\beta-1}}) - [r]
       \circ ([]!{-(0,-.3)} {3}) - [r]
        \bullet ([]!{-(0,-.3)} {1}) - [r]
    \circ ([]!{-(0,-.3)} {2}) - [r]
    \circ ([]!{-(0,-.3)} {(\beta+3)}) (
        - [d] \circ ([]!{-(.3,0)} {2}),
        - [r] \circ ([]!{-(-0.1,-.3)} {2^{n-1}})
         - [r] \circ ([]!{-(0,-.3)} {2}) (
         - [d] \circ ([]!{-(.3,0)} {2}),
        - [r] \circ ([]!{-(0,-.3)} {2})
)}$
 for some $0\le \beta\le n+3$.
 
 $(\mathrm{SR}1)$\;
$S-[\mathrm{I}_n^*\mathrm{SR1}] 
    := 
    \xygraph{
    \diamond ([]!{-(0,-.3)} {+0}) 
    - [r]  \circ ([]!{-(0,-.3)} {2}) 
    - [r]  \circ ([]!{-(0,-.3)} {2})(
        - [r] \circ ([]!{-(0,-.3)} {2}),
        - [d] \circ ([]!{-(-0.25,-.0)} {2}),
}
\quad \quad \quad $
for $n\geq0$.

$(\mathrm{SR}2)$
$S-[\mathrm{I}_n^*\mathrm{SR2}] 
    := 
    \xygraph{
    \diamond ([]!{-(0,-.3)} {+0}) 
    - [r]  \circ ([]!{-(0,-.3)} {2^{n+1}}) 
    - [r]  \circ ([]!{-(0,-.3)} {2})(
        - [r] \circ ([]!{-(0,-.3)} {2}),
        - [d] \circ ([]!{-(-0.25,-.0)} {2}),
}
\quad \quad \quad $
for $n\geq1$.
\end{lem}

\begin{proof}
We first assume that $D_{j}$ has exactly one fork $D$ of valency $3$.
Let $D_{j}-D=B_1+B_2+B_3$ denote the decomposition into three branches with $D'_{j}\subset B_{1}$.
Note that the length of $B_1$ is either $1$ or $n+1$.
From the assumption, $B_2$ and $B_3$ are chains of $(-2)$-curves.
Assume contrary that the proper transform of $D$ is not a fork of exceptional curves over any  strictly lc rational singularity, or equivalently, it becomes a part of a T-chain after the blow-ups.
In this case, at least one of the branches $B_1$, $B_2$ or $B_3$ is completely separated from $D$ after the blow-ups.
If $B_1$ is completely separated from $D$, the proper transform of $D$ and two other branches must form a  P-admissible chain
$$
[2^{l_{2}}, \beta+3,2^{l_{3}}]
$$
for some $\beta\ge 0$, where $l_{i}$ is the length of $B_{i}$.
This contradicts Lemma \ref{lem:extT_2m2}. 
Thus either $B_2$ or $B_3$ is completely separated from $D$.
We may assume $B_2$ is completely separated from $D$.
After the blow-ups, $B_2$ becomes a P-admissible chain
$$
[2^{\beta},3,2^{l_{2}-1}]
$$
for some $\beta\ge 0$, which also contradicts Lemma \ref{lem:extT_2m2}.
Hence we conclude that the proper transform of $D$ is a fork of exceptional curves over some strictly lc rational singularity $q$ after the blow-ups.
It follows that $l_2=l_3=1$ from Lemma~\ref{strlcratbranch} $(2)$.
This gives the case $(\mathrm{SR}1)$ and the case $(\mathrm{SR}2)$.

We now assume that $D_{j}$ coincides with the support of $F_j$.
For simplicity, we only consider the case $n\geq1$; the case $n=0$ follows from a similar argument.
Let $D$ be the fork of $D_j$ disjoint from $D'_{j}$.
Then the same proof as above works and so $D$ corresponds to a fork of exceptional curves over some strictly lc rational singularity $q$.
Since two branches of $D$ are of length $1$, the singularity $q$ is a strictly lc rational singularity of type $(2,2,2,2)$. 
In particular, the exceptional curves $\pi^{-1}(q)$ has two forks and so the other fork of $D_j$ must correspond to a fork of $\pi^{-1}(q)$.
After the blow-ups over the intersection point of $D'_{j}$ and the fork, we obtain the case $(\mathrm{R})$.
\end{proof}
The rest three cases are easily classified by Lemma \ref{strlcratbranch}. 
We summarize the result.

\begin{lem} \label{modifyII*}
Assume that $D_{j}$ has a fork and $F_{j}$ is of type $\mathrm{II}^*$.
By reordering the blow-ups if necessary, there exists some $k$ such that the graph $G_{j,k}$ coincides with the following graph:

$(\mathrm{SR})\;
S-[\mathrm{II}^{*}\mathrm{SR}]:= \xygraph{
    \diamond ([]!{-(0,-.3)} {+0}) 
    - [r]  \circ ([]!{-(0,-.3)} {2}) 
    - [r]  \circ ([]!{-(0,-.3)} {2}) 
    - [r]  \circ ([]!{-(0,-.3)} {2}) 
    - [r]  \circ ([]!{-(0,-.3)} {2}) 
    - [r]  \circ ([]!{-(0,-.3)} {2}) 
    - [r]  \circ ([]!{-(0,-.3)} {2})(
        - [r] \circ ([]!{-(0,-.3)} {2}),
        - [d] \circ ([]!{-(-0.25,-.0)} {2}),
}
$
\end{lem}

\begin{lem} \label{modifyIII*}
Assume that $D_{j}$ has a fork and $F_{j}$ is of type $\mathrm{III}^*$.
By reordering the blow-ups if necessary, there exists some $k$ such that the graph $G_{j,k}$ coincides with the following graph:

$(\mathrm{SR})\;
S-[\mathrm{III}^{*}\mathrm{SR}]:= \xygraph{
    \diamond ([]!{-(0,-.3)} {+0}) 
    - [r]  \circ ([]!{-(0,-.3)} {2}) 
    - [r]  \circ ([]!{-(0,-.3)} {2}) 
    - [r]  \circ ([]!{-(0,-.3)} {2}) 
    - [r]  \circ ([]!{-(0,-.3)} {2})(
        - [r] \circ ([]!{-(0,-.3)} {2}),
        - [d] \circ ([]!{-(-0.25,-.0)} {2}),
}
$
\end{lem}

\begin{lem} \label{modifyIV*}
Assume that $D_{j}$ has a fork and $F_{j}$ is of type $\mathrm{IV}^*$.
By reordering the blow-ups if necessary, there exists some $k$ such that the graph $G_{j,k}$ coincides with the following graph:

$(\mathrm{SR})\;
S-[\mathrm{IV}^{*}\mathrm{SR}]:= \xygraph{
    \diamond ([]!{-(0,-.3)} {+0}) 
    - [r]  \circ ([]!{-(0,-.3)} {2}) 
    - [r]  \circ ([]!{-(0,-.3)} {2}) 
    - [r]  \circ ([]!{-(0,-.3)} {2})(
        - [r] \circ ([]!{-(0,-.3)} {2}),
        - [d] \circ ([]!{-(-0.25,-.0)} {2}),
}
$
\end{lem}

\subsubsection{Step~2: completely separated branch} \label{subsubsec:0,1,0-step2}
In this step, we classify all the cases where $D_j$ is completely separated from $S$ after the blow-ups.
It is enough to compute the finite patterns listed in Proposition \ref{amulet}, which is the classification of initial blow-up processes.
The result is summarized in the following proposition. 

\begin{prop} \label{compsepamulet}
    Assume that $(g,n_X,l_{E_h})=(0,1,0)$ or $(g,n_X,l_{E_h})=(1,1,1)$. 
    Assume that $D_j$ is completely separated from $C_\hor=S$ after the blow-ups.
    By reordering the blow-ups if necessary, there exists some $k$ such that exactly one of the following holds:
    \begin{itemize}
        \item[$(1)$] $G_{j,k} = S-[\mathrm{I}_n \mathrm{T}]_0$ ($n\geq1$) and  
        \[G_{j,0}
=\xygraph{
    \diamond ([]!{-(0,-.3)} {+2}) - [r]
     \bullet ([]!{-(0,-.3)} {1}) - [r]
      \circ ([]!{-(0,-.3)} {2}) - [r]
       \circ ([]!{-(0,-.3)} {5})  (
        - [d] \bullet ([]!{-(.3,0)} {1}),
        - [r] \bullet ([]!{-(-0,-.3)} {1})
         - [r] \circ ([]!{-(0,-.3)} {3}) 
         - [r] \circ ([]!{-(0,-.3)} {2^{n-3}})
        - [r] \circ ([]!{-(0,-.3)} {3})
         - [lllld] \circ ([]!{-(0,-.3)})
       }
       \]
Here, if $n=1$, the subchain $1-[3, 2^{n-3},3,1]-1$ is regarded as $1$.
        \item[$(2)$] $G_{j,k} = S-[\mathrm{IIT1}]$ and $
        G_{j,0}=\xygraph{!~:{@{=}}
    \diamond ([]!{-(0,-.3)} {+2}) - [r]
     \bullet ([]!{-(0,-.3)} {1}) - [r]
      \circ ([]!{-(0,-.3)} {2}) - [r]
       \circ ([]!{-(0,-.3)} {5}) :[r]
       \bullet ([]!{-(0,-.3)} {1})
       }
       $

       \item[$(3)$] $G_{j,k} = S-[\mathrm{IIIT1}]$ and 
       $
       G_{j,0}
       =\xygraph{
    \diamond ([]!{-(0,-.3)} {+2}) - [r]
      \bullet ([]!{-(0,-.3)} {1}) - [r]
      \circ ([]!{-(0,-.3)} {2}) - [r]
       \circ ([]!{-(0,-.3)} {5}) 
        - [r] \bullet ([]!{-(-0,-.3)} {1})(
        - [d] \bullet ([]!{-(.3,0)} {2}),
         - [r] \circ ([]!{-(0,-.3)} {4}) 
}
$

        \item[$(4)$] $G_{j,k} = S-[\mathrm{IVT2}]$ and
        $
        G_{j,0}=\xygraph{
    \diamond ([]!{-(0,-.3)} {+2}) - [r]
      \bullet ([]!{-(0,-.3)} {1}) - [r]
      \circ ([]!{-(0,-.3)} {2}) - [r]
       \circ ([]!{-(0,-.3)} {5}) - [r]
        \bullet ([]!{-(0,-.3)} {1}) - [r]
    \circ ([]!{-(-0.25,-.0)} {2}) (
        - [d] \circ ([]!{-(-0.25,-.0)} {3}),
        - [u] \circ ([]!{-(-0.25,-.0)} {3}),
}
        $ 

        \item[$(5)$] $G_{j,k} = S-[\mathrm{I}_n^* \mathrm{R}]_0$ and
        \[
        G_{j,0}=
        \xygraph{
    \diamond ([]!{-(0,-.3)} {+2}) - [r]
     \bullet ([]!{-(0,-.3)} {1}) - [r]
      \circ ([]!{-(0,-.3)} {2}) - [r]
       \circ ([]!{-(0,-.3)} {5}) - [r]
       \bullet ([]!{-(0,-.3)} {1}) - [r]
    \circ ([]!{-(0,-.3)} {2}) - [r]
    \circ ([]!{-(0,-.3)} {3}) (
        - [d] \circ ([]!{-(.3,0)} {2}),
        - [r] \circ ([]!{-(-0.1,-.3)} {2^{n-1}})
         - [r] \circ ([]!{-(0,-.3)} {2}) (
         - [d] \circ ([]!{-(.3,0)} {2}),
        - [r] \circ ([]!{-(0,-.3)} {2})
)}
\]

        \item[$(6)$] $k=0$ and $G_{j,0}=S-[\mathrm{I}_n \mathrm{E}]_0 \quad\quad(n\geq0)$ 
    \end{itemize}
\end{prop}

\begin{proof}
We first remark that in the case of Proposition \ref{amulet}~(4) with $\beta=0$, $D_j$ is already completely separated from $S$ by the initial blow-up process.
This gives the case $(6)$ in the statement.
Hereafter, we consider the other cases.

Let $G_{j,k}$ be the graph associated to the initial blow-up process classified in Proposition \ref{amulet}, and let $B$ be the connected component of $D_{j,k}$ that intersects  $S$.
Let $a+1$ be the number of blow-ups on $S$ required for $B$ to be completely separated from $S$, where $a\ge 0$.

We assume that Proposition \ref{amulet}~(3) holds.
If $G_{j,k}=S-[\mathrm{I}_n\mathrm{SR1}]$ defined in Lemma \ref{modifyI}, $B$ becomes the curves corresponding to the following dual graph:
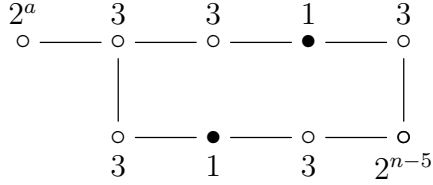
\[
\xygraph{
        \circ ([]!{-(0,-.3)} {2^a}) - [r]
    \circ ([]!{-(-0,-.3)} {3}) (
        - [d] \circ ([]!{-(-0,.3)} {3})
        - [r] \bullet ([]!{-(-0,.3)} {1})
        - [r] \circ ([]!{-(-0,.3)} {3})
        - [r] \circ ([]!{-(-0,.3)} {}),
        - [r] \circ ([]!{-(-0,-.3)} {3})
         - [r] \bullet ([]!{-(-0,-.3)} {1})
         - [r] \circ ([]!{-(-0,-.3)} {3})
          - [d] \circ ([]!{-(-0,.3)} {2^{n-5}}),
}
\]
If $a=0$, then this must be P-admissible.
However, one easily sees that $[3,3,3]$ is not P-admissible (Lemma \ref{lem:extT_I}), which is a contradiction.
Suppose $a>0$. 
If the fork in this graph does not correspond to a fork of the exceptional curves over any strictly lc singularity, then, due to Lemma \ref{strlcratbranch}, a branch consisting of a single $(-3)$-curve is completely separated and becomes a single $(-4)$-curve.
The remaining curves form the extended T-chain $[2,4,3]$ whose associated T-train is $[2,5]-1-[4]$.
In particular, we have the following train in the fiber:
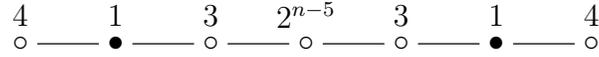
\[
      \xygraph{!~:{@{=}}
    \circ ([]!{-(0,-.3)} {4}) - [r]
     \bullet ([]!{-(0,-.3)} {1}) - [r]
      \circ ([]!{-(0,-.3)} {3}) - [r]
       \circ ([]!{-(0,-.3)} {2^{n-5}}) - [r]
      \circ ([]!{-(0,-.3)} {3}) - [r]
       \bullet ([]!{-(0,-.3)} {1}) -[r]
       \circ ([]!{-(0,-.3)} {4})
       }
\]
This leads to a contradiction as it violates the ample condition.
If the fork in the graph corresponds to a fork of the exceptional curves over some strictly lc rational singularity, then this singularity must be of type $(3,3,3)$.
However, Lemma~\ref{strlcratbranch} implies a contradiction.
In other cases, we can derive contradictions in a similar manner.

Suppose that one of (1), (2) and (4) (with $\beta\ge 1$) in Proposition \ref{amulet} holds.
From the classification, $S+B_{j,k}$ forms a chain $[\chi+m, b_1,\dots b_r]$, where $m\geq0$ and $B=[b_1,\dots b_r]$ is one of the following forms: $[2^\beta]$ ($\beta\geq1$), $[3+\beta,2^{n-2},3,2^\beta]$ ($\beta\geq0$ and $n\geq1$).
Since $B_{j,k}$ is completely separated from $S$ after $a+1$ blow-ups, the chain $L_2^aL_1B$ becomes a P-admissible chain.
It follows from Lemma~\ref{lem:extT_2m2} and \ref{lem:extT_CS} that $B$ must be of the form $[3,2^{n-1},3]$ for some $n\geq0$ and $a=1$.
This condition determines the admissible initial blow-up processes: $[\mathrm{I}_n \mathrm{T}]_0$ ($n\geq1$), $[\mathrm{IIT1}]$, $[\mathrm{IIIT1}]$, $[\mathrm{IVT1}]$, $[\mathrm{IVT2}]$ and $[\mathrm{I}_n^* \mathrm{R}]_\beta$ ($0\leq\beta\leq n+3$).
The cases $[\mathrm{IVT1}]$ and $[\mathrm{I}_n^* \mathrm{R}]_\beta$ ($\beta>0$) are excluded by using the fact that they violate the ample condition.
Hence we obtain the classification.
\end{proof}

\subsubsection{Step~3: case that $p$ is a T-singularity }
In this step, we classify possible $C$ and its resolution $\widetilde{C}$ under the assumption that $p$ is a T-singularity. 

Since $\widetilde{C}^{(p)}$ is a chain, there are the following three types for $\widetilde{C}^{(p)}$:
\begin{itemize}
    \item[(i)] There exist no branches of $S=\widetilde{C}^{(p)}$ ($J=0$).
  
    \item[(ii)] There exists only one branch of $S\subset\widetilde{C}^{(p)}$ ($J=1$).
    
    \item[(iii)] There exist exactly two branches of $S\subset\widetilde{C}^{(p)}$ ($J=2$).
\end{itemize}
Before stating the result, we introduce notation, which are subsets of the set $\cA$ of amulets. 
\begin{align*}
\cA_{[2^\beta]} :=& \{  [\mathrm{CT}]_\beta \} \cup \{ [\mathrm{I}_n \mathrm{E}]_\beta \mid \beta \le n+8 \},   \\ 
\cA_{[4]} :=& \{ [\mathrm{I}_1\mathrm{T}]_0, [\mathrm{IIT1}], [\mathrm{IVT2}] \} \cup \{ [\mathrm{I}_n^* \mathrm{R}]_0 : n\geq0 \}, \\
\cA_{[3,2^{m-1},3]}:=&
  \begin{cases}
    \{ [\mathrm{IIIT1}], [\mathrm{I}_2\mathrm{T}]_{0} \}  \cup \{ [\mathrm{I}_n^* \mathrm{R}]_1 : n\geq0 \},                 & \text{if $m=1$,} \\
    \{ [\mathrm{IVT1}], [\mathrm{I}_3\mathrm{T}]_{0} \}  \cup \{ [\mathrm{I}_n^* \mathrm{R}]_2 : n\geq0 \},       & \text{if $m=2$,} \\
     \{ [\mathrm{I}_{m+1}\mathrm{T}]_{0} \} \cup \{ [\mathrm{I}_n^* \mathrm{R}]_{m} : n\geq0 \},       & \text{if $m \geq 3$.}
  \end{cases}\\
  \cA_{[m+4,2^{m}]}:=&
  \begin{cases}
    \{ [\mathrm{I}_1\mathrm{T}]_1, [\mathrm{IIT2}], [\mathrm{IIIT3}], [\mathrm{IVR}]_{0} \}, & \text{if $m=1$,}  \\ 
   \{ [\mathrm{I}_1\mathrm{T}]_2, [\mathrm{IIT3}], [\mathrm{IIIR}]_{0} \}, & \text{if $m=2$,}   \\
   \{ [\mathrm{I}_1\mathrm{T}]_3, [\mathrm{IIT4}] \}, & \text{if $m=3$,}  \\ 
   \{ [\mathrm{I}_1\mathrm{T}]_4, [\mathrm{IIR}] \}, & \text{if $m=4$,}  \\ 
  \end{cases}\\
  \cA_{[4,2^{m},3,2]}:=&
  \begin{cases}
   \{ [\mathrm{I}_2\mathrm{T}]_1, [\mathrm{IIIT2}], [\mathrm{IVR}]_{1}\}, & \text{if $m=0$,} \\
   \{ [\mathrm{I}_3\mathrm{T}]_{1}, [\mathrm{IVR}]_2 \}, & \text{if $m=1$,}  \\
    \{ [\mathrm{I}_{m+2}\mathrm{T}]_{1} \}, & \text{if $m\geq 2 $,}  \\
  \end{cases}\\
   \cA_{[5,2^{m},3,2^2]}:=&
  \begin{cases}
   \{ [\mathrm{IIIR}]_1, [\mathrm{I}_2\mathrm{T}]_2\}, & \text{if $m=0$,} \\
   \{ [\mathrm{I}_{m+2}\mathrm{T}]_2 \}, & \text{if $m \ge 1$,} \\
  \end{cases}\\
  \cA_{[3+\beta,2^{n-2},3,2^{\beta}]}:=&\{[\mathrm{I}_n \mathrm{T}]_{\beta}\} \quad \text{if $n=1$, $\beta \ge 5$ or $n=2$, $\beta\ge 3$ or $n=3$, $\beta\ge 2$ or $n\ge 4$}
\end{align*}

\begin{thm}\label{thm:one-section-T}
Let $(g,n_X,l_{E_h})=(0,1,0)$ and $\chi=\chi(\mathcal{O}_{X})$.
Assume that $p \in X$ is a T-singularity, and there exists no branch that is completely separated from the section $S$.
Then, exactly one of the following hold:
\begin{enumerate}[label=$\bullet$]
    \item $J=0$, $\chi=4$, and $B_0=[4]$.
    \item $J=1$. 
    By reordering the blow-ups if necessary, there exists some $k$ such that $\chi$, $G_{1,k}$, $G_{2,k}$ and $B_k$ coincide with exactly one in the following list. 
    Here, the number in $B_k$ corresponding to the proper transform of $S$ is underlined.
    
    \smallskip
    
    \begin{tabular}{|l|l|l|l|} \hline
            &   $\chi$      &   $G_{1,k}=S-A$       &   $B_k$    \\
        \hline
        $(1.1)$ &   $\chi=3$    &   $A=[\mathrm{I}_n \mathrm{R}]_\beta$        &   $[\underline{4},2^{\beta-1},3,2]$        \\
        $(1.2)$ &   $\chi\geq3$    &   $A\in\cA_{[\chi+1,2^{\chi-4+n},3,2^{\chi-2}]}$        &   $[\underline{\chi},\chi+1,2^{\chi-4+n},3,2^{\chi-2}]$        \\
        $(1.3)$ &   $\chi\geq5$ &   $A\in\cA_{[2^{\chi-4}]}$ &   $[\underline{\chi},2^{\chi-4}]$ \\
        \hline
    \end{tabular}
    
    \smallskip
    
    \noindent
    In cases $(1.1)$ and $(1.3)$, the chain $B_k=B_0$ is a T-chain.
    In case $(1.2)$, the total transform of $B_k$, which is a T-train associated to $B_k$, is as follows:
    \[
        [\underline{\chi},\chi+1,2^{\chi-4},3,2^{\chi-2}]-1-[\chi+1,2^{n-2},3,2^{\chi-2}]
    \]
    
    \item $J=2$. 
    By reordering the blow-ups if necessary, there exists some $k$ such that $\chi$, $G_{1,k}$, $G_{2,k}$, $B_k$ and the total transform of $B_k$ in $\widetilde{X}$ coincide with exactly one in the following list.
    Here, the number in $B_k$ corresponding to the proper transform of $S$ is underlined. 
    
    \smallskip
    
    \begin{tabular}{|l|l|l|l|l|} \hline
                &   $\chi$  &   $G_{1,k}=S-A_1$   &   $G_{2,k}=S-A_2$       &   $B_k$     \\
        \hline
        $(2.1)$   &   $\chi=2$     &   $\cA_{[3,2^{1+m},3]}$   &   $\cA_{[4,2^{-1+n},3,2]}$   &   $[3,2^{1+m},3,\underline{2},4,2^{-1+n},3,2]$     \\
        $(2.2)$   &   $\chi=2$     &   $\cA_{[3,2^{m},3]}$   &   $\{[\mathrm{I}_n \mathrm{R}]_0\}$   &   $[3,2^m,3,\underline{4},2]$     \\
        $(2.3)$   &   $\chi=2$     &   $\cA_{[3,2^m,3]}$   &   $\{[\mathrm{I}_n \mathrm{R}]_1\}$   &   $[3,2^m,3,\underline{3},3,2]$     \\
        $(2.4)$   &   $\chi=3$     &   $\cA_{[3,2^m,3]}$   &   $\cA_{[4,2^n,3,2]}$   &   $[3,2^m,3,\underline{3},4,2^n,3,2]$     \\
        $(2.5)$   &   $\chi=3$     &   $\cA_{[3,2^m,3]}$   &   $\{[\mathrm{I}_n \mathrm{R}]_1\}$   &   $[3,2^m,3,\underline{4},3,2]$     \\
        $(2.6)$   &   $\chi=4$     &   $\cA_{[3,2^{2+m},3]}$   &   $\cA_{[4,2^{-1+n},3,2]}$   &   $[3,2^{2+m},3,\underline{4},4,2^{-1+n},3,2]$     \\
        $(2.7)$   &   $\chi=4$     &   $\cA_{[3,2^m,3]}$   &   $\cA_{[4,2^{1+n},3,2]}$   &   $[3,2^m,3,\underline{4},4,2^{1+n},3,2]$     \\
        $(2.8)$   &   $\chi=4$     &   $\cA_{[3,2^{m},3]}$   &   $\{[\mathrm{I}_n \mathrm{R}]_1\}$   &   $[3,2^m,3,\underline{5},3,2]$  \\
        $(2.9)$   &   $\chi=5$     &   $\cA_{[3,2^{-1+m},3]}$   &   $\cA_{[2^2]}$   &   $[3,2^{-1+m},3,\underline{5},2^2]$     \\
        $(2.10)$   &   $\chi=5$     &   $\cA_{[3,2^{2+m},3]}$   &   $\cA_{[4,2^{-1+n},3,2]}$   &   $[3,2^{2+m},3,\underline{5},4,2^{-1+n},3,2]$     \\
        $(2.11)$   &   $\chi\geq2$     &   $\cA_{[3,2^{\chi-2+m},3]}$   &   $\cA_{[4,2^n,3,2]}$   &   $[3,2^{\chi-2+m},3,\underline{\chi},4,2^n,3,2]$     \\
        $(2.12)$   &   $\chi \ge 2$     &   $\cA_{[3,2^{-1+m},3]}$   &   $\cA_{[5,2^{\chi-3+n},3,2^2]}$   &   $[3,2^{-1+m},3,\underline{\chi},5,2^{\chi -3 +n },3,2^2]$     \\
        $(2.13)$   &   $\chi\geq3$     &   $\cA_{[3,2^{\chi-4+m},3]}$   &   $\cA_{[2]}$    &   $[3,2^{\chi-4+m},3,\underline{\chi},2]$     \\
        $(2.14)$   &   $\chi\geq3$     &   $\cA_{[3,2^{\chi-1+m},3]}$   &   $\cA_{[4,2^{-1+n},3,2]}$   &   $[3,2^{\chi-1+m},3,\underline{\chi},4,2^{-1+n},3,2]$     \\
        $(2.15)$   &   $\chi\geq3$     &   $\cA_{[3,2^m,3]}$   &   $\cA_{[4,2^{\chi-2+n},3,2]}$   &   $[3,2^m,3,\underline{\chi},4,2^{\chi-2+n},3,2]$     \\
        $(2.16)$   &   $\chi\geq3$     &   $\cA_{[3,2^{-1+m},3]}$   &   $\cA_{[5,2^{\chi-2+n},3,2^2]}$   &   $[3,2^{-1+m},3,\underline{\chi},5,2^{\chi-2+n},3,2^2]$     \\
        \hline
    \end{tabular}
    
    \smallskip
    
    \noindent
    In the table, the variables $m,n$ are non-negative integers in each case. 
    The total transforms of $B_k$, which are T-trains associated to $B_k$, are described as follows:
    \begin{itemize}
        \item[$(2.1)$] $[3,2^{m-2},3]-1-[3,2,3,\underline{2},5,2]-1-[4,2^{n-2},3,2]$
        \item[$(2.2)$] $[3,2^{m-2},3]-1-[3,3,\underline{4},2]$ 
        \item[$(2.3)$] $[3,2^{m-2},3]-1-[3,3,\underline{4},2]-1-[5,2]$
        \item[$(2.4)$] $[3,2^{m-2},3]-1-[3,3,\underline{3},4,3,2]-1-[4,2^{n-2},3,2]$
        \item[$(2.5)$] $[3,2^{m-2},3]-1-[3,5,2]-1-[3,\underline{5},3,2]$
        \item[$(2.6)$] $[3,2^{m-2},3]-1-[3,2^2,7,2]-1-[3,2^2,\underline{5},5,2]-1-[4,2^{n-2},3,2]$
        \item[$(2.7)$] $[3,2^{m-2},3]-1-[3,3,\underline{5},3,2]-1-[3,6,2,3,2]-1-[4,2^{n-2},3,2]$
        \item[$(2.8)$] $[3,2^{m-2},3]-1-[3,3,\underline{5},3,2]$
        \item[$(2.9)$] $[3,2^{m-2},3]-1-[4,\underline{5},2^2]$
        \item[$(2.10)$] $[3,2^{m-2},3]-1-[3,2^2,3,\underline{5},5,2]-1-[4,2^{n-2},3,2]$
        \item[$(2.11)$] $[3,2^{m-2},3]-1-[3,2^{\chi-2},3,\underline{\chi+2},2]-1-[3,5,3,2]-1-[4,2^{n-2},3,2]$
        \item[$(2.12)$] $[3,2^{m-2},3]-1-[4,\underline{\chi},5,2^{\chi-3},3,2^2]-1-[5,2^{n-2},3,2^2]$
        \item[$(2.13)$] $[3,2^{m-2},3]-1-[3,2^{\chi-4},3,\underline{\chi},2]$
        \item[$(2.14)$] $[3,2^{m-2},3]-1-[3,2^{\chi-1},3,\underline{\chi+3},2]-1-[3,2,6,2]-1-[4,2^{n-2},3,2]$
        \item[$(2.15)$] $[3,2^{m-2},3]-1-[3,5,2]-1-[3,\underline{\chi+2},2^{\chi-3},3,2]-1-[3,\chi+3,2^{\chi-2},3,2]-1-[4,2^{n-2},3,2]$
        \item[$(2.16)$] $[3,2^{m-2},3]-1-[4,\underline{\chi+1},2^{\chi-4},3,2^2]-1-[4,\chi+3,2^{\chi-2},3,2^2]-1-[5,2^{n-2},3,2^2]$
    \end{itemize}
\end{enumerate}

\end{thm}

\begin{rem}\label{ignorecompsep}
    In the above, we impose the condition that no branch is completely separated from $S$.
    We can obtain the list without this assumption as follows.
    Let $J_{\operatorname{CS}}$ (resp. $J_{\operatorname{CSE}}$) denote the number of completely separated branches of the form $(1)\text{--}(5)$ (resp. $(6)$) in Proposition~\ref{compsepamulet}.
    We reorder the sequence of blow-ups $\widetilde{X}=Y_{0} \to Y_{m}=Y$ and decompose it as
    \[
    Y_{0} \to Y_{k_1} \to Y_{k_2} \to Y
    \]
    so that the morphism $Y_{k_2} \to Y$ separates all the completely separated branches from $S$, and $Y_{k_1} \to Y_{k_2}$ is the composition of the rest of the initial blow-ups.
    Let $S_{k_i}$ ($i=1,2$) be the proper transform of $S\subset Y$ on $Y_{k_i}$. 
    Then, from Proposition \ref{compsepamulet}, we have
    \[
    -S_{k_2}^2 = -S^2+2J_{\operatorname{CS}} + J_{\operatorname{CSE}}.
    \]
    Moreover, since $\chi(\mathcal{O}_{Y})=-S^2$, it follows from Proposition \ref{amulet} that
    \[
    -S_{k_1}^2 = -S^2+2J_{\operatorname{CS}} + J_{\operatorname{CSE}}+
    \sharp{\{j \in J: G_{j,k_1}=S-[\mathrm{I}_{n}\mathrm{E}]_{\beta \ge 1}\}}+c,
    \]
    where 
    \begin{align*}
    c:=2\sharp \{j \in J:{G_{j,k_1}=S-[\mathrm{I}_n \mathrm{R}]_{0}}\}
    +\sharp \{j \in J:{G_{j,k_1}=S-[\mathrm{I}_n \mathrm{R}]_{1}}\}.
    \end{align*}
    Since $l_E=l_{E_v}$ from the assumption that $l_{E_h}=0$, we have
    \begin{align*}
    -S_{k_1}^2 = &-S^2+2J_{\operatorname{CS}} + l_{E} + c    \\
               = & (\chi +2J_{\operatorname{CS}}) + c
    \end{align*}
    Therefore, by interpreting $\chi$ in the list of Theorem~\ref{thm:one-section-T} as $\chi + 2J_{\operatorname{CS}}$, the list becomes nothing but the classification of $X$ with $J_{\operatorname{CS}}$ completely separated branches of the form $(1)\text{--}(5)$ in Proposition \ref{compsepamulet}.
\end{rem}

\begin{rem}\label{rem:small_surfaces}
    Monreal-Negrete-Urz\'ua classified certain \emph{small surfaces} \cite[Theorem 3.12]{MNU}.
    In our terminology, a small surface classified in \cite[Theorem 3.12]{MNU} is a surface obtained from a klt surface satisfying Assumption \ref{ass:normal_stable_elliptic} and $(g,n_X)=(0,1)$ by resolving all Du Val singularities.
    The classification of such small surfaces can be recovered from Theorem \ref{thm:one-section-T} (together with Remark \ref{ignorecompsep}) by collecting the candidates consisting of amulets described in Proposition \ref{amulet} (1) and those in Proposition \ref{compsepamulet} (1)--(4).
\end{rem}

\begin{rem}
As observed in \cite[Remark 3.10]{MNU}, the P-admissible chain $[4,\chi,5,2^{\chi-2},3,2^2]$, which appears in (2.11) and (2.15), has two associated ample T-trains:
\begin{itemize}
    \item $[4,\chi,5,2^{\chi-3},3,2^2]-1-[6,2^2]$ 
    \item $[4,\chi+1,2^{\chi-4},3,2^2]-1-[4,\chi+3,2^{\chi-2},3,2^2]$.
\end{itemize}
The former corresponds to (2.11), while the latter corresponds to (2.15).
This is an example of a P-admissible chain that has two associated ample T-trains.
We learned from Urz\'ua that such phenomenon is related to \emph{wormhole cyclic quotient singularities} \cite{UV}.
\end{rem}

\begin{proof}[Proof of Theorem \ref{thm:one-section-T}]
We reorder the sequence of blow-ups appearing in the birational morphism $\widetilde{X}=Y_{0} \to Y_{m}=Y$ and decompose it so that the following holds:
\[
Y_{0} \to Y_{k} \to Y_{m}
\]
Here, the morphism $Y_{k} \to Y_{m}$ completes the initial blow-up process for all $D_j$.
Since there exists no branch that is completely separated from the section $S$,
$B_{k}$ is a chain of the form
\[
[a_1,\cdots,a_r,\chi+c,b_1,\cdots,b_s],
\]
where
 \begin{align*}
    c:=2\sharp \{j \in J:{G_{j,k}=S-[\mathrm{I}_n \mathrm{R}]_{0}}\}
    +\sharp \{j \in J:{G_{j,k}=S-[\mathrm{I}_n \mathrm{R}]_{1}}\}.
    \end{align*}
We note that
$[a_r,\cdots,a_1]$ and $[b_1,\cdots,b_s]$ is one of the following:
\begin{itemize}
    \item[(1)] $[2^\beta]$ $(\beta \geq 0)$
    \item[(2)] $[m+4,2^m]$ $(m \geq 0)$
    \item[(3)] $[3+\beta,2^m,3,2^\beta]$  $(\beta\geq 0, m \geq 0)$
    \item[(4)] $[3,2]$
  \end{itemize}
Since $p \in X$ is a T-singularity,
$B_k$ must be an extended T-chain.
Therefore, from the Lemma~\ref{lem:extT_step3}, we obtain the desired classification.
\end{proof}

\subsubsection{Step~4: case that $p$ is a strictly lc rational singularity}\label{subsubsec:0,1,0-step4}

In this step, we classify possible $C$ and its resolution $\widetilde{C}$ under the assumption that $p$ is a strictly lc rational singularity.
Again, thanks to Remark \ref{ignorecompsep}, we may assume that no branches of $S$ in $C$ are completely separated from $S$.
We will divide the problem into two cases:

\begin{itemize}
\item 
When the dual graph of $\widetilde{C}^{(p)}$ matches the one described in Lemma \ref{lem:smoothable_rational_strict_lc} (i).
\item 
When it corresponds to one of the cases in Lemma \ref{lem:smoothable_rational_strict_lc} (ii), (iii), or (iv).
\end{itemize}

To proceed with the classification, We first establish the following lemma.

\begin{lem}\label{lem:amulet-b}
    Let $(g,n_X,l_{E_h})=(0,1,0)$. 
    Assume that $p\in X$ is a strictly lc rational singularity, and the connected component of $D_{j,0}$ intersecting $S$ consists of a single $(-b)$-curve for some $b$.
    By reordering the blow-ups if necessary, there exists some $k$ such that exactly one of the following hold:
    \begin{enumerate}[label=$(\arabic*)$]
        \item $b=2$, $k=0$, and $G_{j,0}=S- [\mathrm{CT}]_1$. 

        \item $b=2$, $k=0$, and $G_{j,0}=S-[\mathrm{I}_{3+n}\mathrm{R}]_0$ $(n\geq 0)$
        \item $b=2$, $k=0$, and $G_{j,0}=S-[\mathrm{I}_n \mathrm{E}]_1$ $(n\geq0)$.

        \item $b=2$, $G_{j,k}=S-[\mathrm{I}_1\mathrm{T}]_1$ and $G_{j,0}$ is
        \begin{align*}
        S-[\mathrm{I}_1\mathrm{T}_1\mathrm{B2}] := 
                &\xygraph{
    \diamond ([]!{-(0,-.3)} {+3})
    - [r] \circ ([]!{-(0,-.3)} {2})
    - [r] \bullet ([]!{-(0,-.3)} {1})
     - [r] \circ ([]!{-(0,-.3)} {3})
      - [r] \circ ([]!{-(0,-.3)} {2})
    - [r] \circ ([]!{-(0,-.3)} {6}) (
        - [d] \bullet ([]!{-(.3,0)}  {1}, -[ru] \circ ),
        - [r] \circ ([]!{-(0,-.3)} {2})
}
\end{align*}
        \item $b=2$, $G_{j,k}=S-[\mathrm{I}_{2+n}\mathrm{T}]_1$ $(n\geq 0)$ and $G_{j,0}$ is 
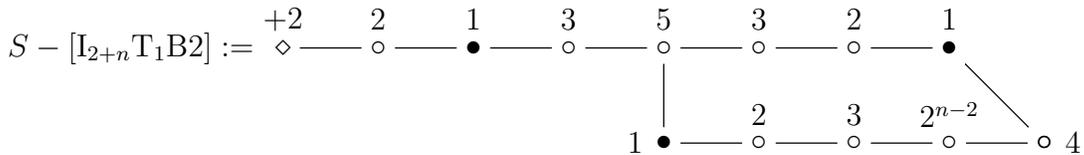
\begin{align*}
        S-[\mathrm{I}_{2+n} \mathrm{T}_1\mathrm{B2}] := 
&\xygraph{
    \diamond ([]!{-(0,-.3)} {+2}) - [r]
      \circ ([]!{-(0,-.3)} {2}) - [r] 
      \bullet ([]!{-(0,-.3)} {1}) - [r] 
       \circ ([]!{-(0,-.3)} {3}) - [r] 
        \circ ([]!{-(0,-.3)} {5})  (
        - [d] \bullet ([]!{-(.3,0)}  {1})
         - [r] \circ ([]!{-(0,-.3)}  {2})
        - [r] \circ ([]!{-(0,-.3)}  {3})
        -[r] \circ ([]!{-(0,-.3)}  {2^{n-2}})
         - [r] \circ ([]!{-(.3,0)} ),
        - [r] \circ ([]!{-(0,-.3)} {3})
        - [r] \circ ([]!{-(0,-.3)} {2})
          - [r] \bullet ([]!{-(0,-.3)} {1})
        - [dr] \circ ([]!{-(-0.3,-.0)} {4})
}
\end{align*}  
        \item $b=2$, $G_{j,k}=S-[\mathrm{IIT2}]$ and $G_{j,0}$ is
        \[ 
        S-[\mathrm{IIT2B2}] := 
    \xygraph{
    \diamond ([]!{-(0,-.3)} {+3})
    - [r]  \circ ([]!{-(0,-.3)} {2})
    - [r]  \bullet ([]!{-(0,-.3)} {1})
    - [r]  \circ ([]!{-(0,-.3)} {3})
    - [r]  \circ ([]!{-(0,-.3)} {2})
    - [r]  \circ ([]!{-(0,-.3)} {6}) 
      (- [d]   \bullet ([]!{-(-0.25,-.0)}  {1} ),
       - [r]  \circ ([]!{-(0,-.3)} {2}) 
       -[ld]  \circ ([]!{-(0,-.3)} 
       }
        \]
        \item $b=2$, $G_{j,k}=S-[\mathrm{IIIT2}]$ and $G_{j,0}$ is
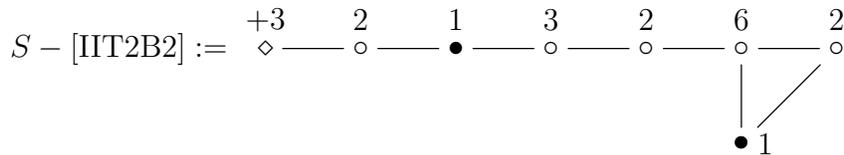
        \[
        S-[\mathrm{IIIT2B2}] := 
        \xygraph{
    \diamond ([]!{-(0,-.3)} {+2}) - [r]
    \circ ([]!{-(0,-.3)} {2}) - [r]
     \bullet ([]!{-(0,-.3)} {1}) - [r]
    \circ ([]!{-(0,-.3)} {3}) - [r]
    \circ ([]!{-(0,-.3)} {5})(
        - [r] \circ ([]!{-(0,-.3)} {3}) ( 
        - [r] \bullet ([]!{-(-0,-.3)} {1})
         - [r] \circ ([]!{-(0,-.3)} {2}) 
        - [d] \circ ([]!{-(-0.3,-.0)} {5}),
         - [d] \circ ([]!{-(-0.3,-.0)} {2})
    )}
        \]
        
        \item $b=2$, $G_{j,k}=S-[\mathrm{IIIT3}]$ and $G_{j,0}$ is
        \[
        S-[\mathrm{IIIT3B2}] := 
          \xygraph{
    \diamond ([]!{-(0,-.3)} {+3}) - [r]
    \circ ([]!{-(0,-.3)} {2}) - [r]
    \bullet ([]!{-(0,-.3)} {1}) - [d]
    \circ ([]!{-(-0.25,-.3)} {3}) 
        - [r] \circ ([]!{-(0,-.3)} {2})
        - [r] \circ ([]!{-(0,-.3)} {6})
        - [r] \circ ([]!{-(0,-.3)} {2})
         - [r] \bullet ([]!{-(0,-.3)} {1}) 
         - [r] \circ ([]!{-(-0.25,-.3)} {3}) (
        - [u] \circ ([]!{-(-0.25,-.0)} {2}),
        - [r] \circ ([]!{-(0,-.3)} {4}))
    )}
        \]

        \item $b=2$, $G_{j,k}=S-[\mathrm{I}_{3+n}\mathrm{R}]_1$ $(n\geq 0)$ and $G_{j,0}$ is
        \[
        S-[\mathrm{I}_{3+n}\mathrm{R}_1\mathrm{B2}] := 
        \xygraph{
    \diamond ([]!{-(0,-.3)} {+2}) - [r]
     \circ ([]!{-(0,-.3)} {2}) - [r]
      \bullet ([]!{-(0,-.3)} {1}) - [d]
     \circ ([]!{-(.3,0)} {5}) - [r]
      \circ ([]!{-(0,-.3)} {2}) - [r]
      \bullet ([]!{-(-0.3,-.0)} {1}) - [u]
       \circ ([]!{-(0,-.3)} {3}) - [r]
    \circ ([]!{-(-0.25,-.0)} {4}) ( 
        - [u] \circ ([]!{-(.3,0)} {3}) - [r] \bullet ([]!{-(-0,-.3)} {1}) - [r] \circ ([]!{-(-0.3,-.0)} {3}) - [d] \circ ([]!{-(-0.5,-.0)} {2^{n-2}}),
        - [d] \circ ([]!{-(.3,0)} {3})- [r] \bullet ([]!{-(-0,-.3)} {1}) - [r] \circ ([]!{-(-0.3,-.0)} {3}) - [u] \circ ([]!{-(.3,0)} {}),      
)}
        \]
        \item $b=2$, $G_{j,k}=S-[\mathrm{IVR}]_0$ and $G_{j,0}$ is
        \[
        S-[\mathrm{IVR}_0\mathrm{B2}] :=
        \xygraph{
    \diamond ([]!{-(0,-.3)} {+3}) - [r]
      \circ ([]!{-(0,-.3)} {2}) - [r]
       \bullet ([]!{-(0,-.3)} {1}) - [d]
      \circ ([]!{-(.3,0)} {3}) - [r]
        \circ ([]!{-(0,-.3)} {2}) - [r]
           \circ ([]!{-(0,-.3)} {6}) - [r]
              \circ ([]!{-(-0.25,-.0)} {2}) - [u]
        \bullet ([]!{-(0,-.3)} {1}) - [r]
    \circ ([]!{-(0,-.3)} {3}) - [r]
    \circ ([]!{-(0,-.3)} {2}) (
        - [d] \circ ([]!{-(-0.25,-.0)} {3}),
        - [r] \circ ([]!{-(0,-.3)} {3}),}
        \]
        \item $b=2$, $G_{j,k}=S-[\mathrm{IVR}]_\beta$ $(\beta=1,2)$ and $G_{j,0}$ is
        \[
        S-[\mathrm{IVR}_\beta \mathrm{B2}]  := 
         \xygraph{
    \diamond ([]!{-(0,-.3)} {+2}) - [r]
      \circ ([]!{-(0,-.3)} {2}) - [d]
      \bullet ([]!{-(-0.25,-.0)} {1}) - [d]
       \circ ([]!{-(.3,0)} {3}) - [r]
        \circ ([]!{-(0,-.3)} {5}) - [r]
         \circ ([]!{-(0,-.3)} {3}) - [r]
       \circ ([]!{-(-0.25,-.0)} {2}) - [u]
        \bullet ([]!{-(-0.25,-.0)} {1}) - [u]
      \circ ([]!{-(0,-.3)} {4}) - [r]
        \circ ([]!{-(0,-.3)} {2^{\beta-3}}) - [r]
         \circ ([]!{-(0,-.3)} {3}) - [r]
         \circ ([]!{-(0,-.3)} {2}) - [r]
        \bullet ([]!{-(-0.25,-.0)} {1}) - [d]
    \circ ([]!{-(-0.25,-.0)} {3}) - [l]
    \circ ([]!{-(0,-.3)} {(\beta+2)}) (
        - [d] \circ ([]!{-(-0.25,-.0)} {3}),
        - [l] \circ ([]!{-(.2,0)} {3}),
}
        \] 
Here if $\beta=1$, the subchain $-1-[4,2^{\beta-3},3,2]-1-$ is regarded as $-1-$.  
        \item $b=3$, $G_{j,k}=S-[\mathrm{I}_{2+n}\mathrm{T}]_0$ $(n\geq 0)$ and $G_{j,0}$ is
     \begin{align*}
        S-[\mathrm{I}_{2+n} \mathrm{T}_0\mathrm{B3}] := &[\chi+1,3]-1-[2,5,3]-1-[3,2^{n-2},3]\\
    =& \xygraph{
    \diamond ([]!{-(0,-.3)} {+1}) - [r]
      \circ ([]!{-(0,-.3)} {3}) - [r] 
      \bullet ([]!{-(0,-.3)} {1}) - [r] 
       \circ ([]!{-(0,-.3)} {2}) - [r] 
        \circ ([]!{-(0,-.3)} {5})  (
        - [d] \bullet ([]!{-(.3,0)}  {1})
         - [r] \circ ([]!{-(0,-.3)}  {3})
        - [r] \circ ([]!{-(0,-.3)}  {2^{n-2}})
        -[r] \circ ([]!{-(-0.3,-.0)}  {3}),
        - [r] \circ ([]!{-(0,-.3)} {3})
        - [r] \bullet ([]!{-(0,-.3)} {1})
         - [rd] \circ ([]!{-(0,-.3)} {})
}
 \end{align*}
Here if $n=0$, the subchain $-1-[3,2^{n-2},3]-1-$ is regarded as $-1-$.  
        \item $b=3$, $G_{j,k}=S-[\mathrm{IIIT1}]$ and $G_{j,0}$ is
        \[
        S-[\mathrm{IIIT1B3}] :=  \xygraph{
    \diamond ([]!{-(0,-.3)} {+1})
    - [r]  \circ ([]!{-(0,-.3)} {3})
    - [r]  \bullet ([]!{-(0,-.3)} {1})
    - [r]  \circ ([]!{-(0,-.3)} {2})
    - [r]  \circ ([]!{-(0,-.3)} {5}) 
      (- [d]   \bullet ([]!{-(-0.25,-.0)}  {1} ),
       - [r]  \circ ([]!{-(0,-.3)} {3}) 
       -[ld]  \circ ([]!{-(0,-.3)} 
       }
        \]

        \item $b=3$, $G_{j,k}=S-[\mathrm{IVT1}]$ and $G_{j,0}$ is
         \[
        S-[\mathrm{IVT1B3}] := 
        \xygraph{
    \diamond ([]!{-(0,-.3)} {+1}) - [r]
    \circ ([]!{-(0,-.3)} {3}) - [r]
     \bullet ([]!{-(0,-.3)} {1}) - [r]
    \circ ([]!{-(0,-.3)} {2}) - [r]
    \circ ([]!{-(0,-.3)} {5}) - [r]
     \circ ([]!{-(0,-.3)} {3})(
         ( 
        - [r] \bullet ([]!{-(-0,-.3)} {1})
         - [r] \circ ([]!{-(0,-.3)} {4}) ),
         - [d] \bullet ([]!{-(0.3,-.0)} {1},
         - [r] \circ ([]!{-(-0.3,-.0)} {4})
    )}
        \]

        \item $b=3$, $G_{j,k}=S-[\mathrm{I}_{n}^{*}\mathrm{R}]_1$ and $G_{j,0}$ is
           \[
        S-[\mathrm{I}_{n}^{*}\mathrm{R}_1 \mathrm{B3}] := 
        \xygraph{
    \diamond ([]!{-(0,-.3)} {+1}) - [r]
     \circ ([]!{-(0,-.3)} {3}) - [d]
     \bullet ([]!{-(.3,0)} {1}) - [r]
      \circ ([]!{-(0,-.3)} {2}) - [r]
       \circ ([]!{-(0,-.3)} {5}) - [r]
        \circ ([]!{-(-0.3,-.0)} {3}) - [u]
       \bullet ([]!{-(0,-.3)} {1}) - [r]
    \circ ([]!{-(0,-.3)} {2}) - [r]
    \circ ([]!{-(0,-.3)} {4}) (
        - [d] \circ ([]!{-(.3,0)} {2}),
        - [r] \circ ([]!{-(-0.1,-.3)} {2^{n-1}})
         - [r] \circ ([]!{-(0,-.3)} {2}) (
         - [d] \circ ([]!{-(.3,0)} {2}),
        - [r] \circ ([]!{-(0,-.3)} {2})
)}
\]

        \item $b=4$, $k=0$, and $G_{j,0}=S-[\mathrm{I}_1\mathrm{T}]_0$. 
        
        \item $b=4$, $k=0$, and $G_{j,0}=S-[\mathrm{IIT1}]$.

        \item $b=4$, $k=0$, and $G_{j,0}=S-[\mathrm{IVT2}]$.

        \item $b=4$, $k=0$, and $G_{j,0}=S-[\mathrm{I}_n^* \mathrm{R}]_0$ $(n\geq0)$.

        \item $b=4$, $G_{j,k}=S-[\mathrm{I}_{2+n}\mathrm{T}]_0$ $(n\geq 0)$ and $G_{j,0}$ is
        \begin{align*}
        S-[\mathrm{I}_{2+n}\mathrm{T}_0\mathrm{B4}] := 
         \xygraph{
    \diamond ([]!{-(0,-.3)} {+0}) - [r]
    \circ ([]!{-(-0.3,-.0)} {4}) ( 
        - [u] \bullet ([]!{-(.3,0)} {1}) 
        - [r] \circ ([]!{-(-0.3,-.0)} {3}) 
        - [d] \circ ([]!{-(-0.5,-.0)} {2^{n-1}})([]!{-(.5,0)} {}),
        - [d] \bullet ([]!{-(.3,0)} {1})
        - [r] \circ ([]!{-(-0.3,-.0)} {3}) - [u] \circ ([]!{-(.3,0)} {}),    
)}\end{align*}
        \item $b=4$, $G_{j,k}=S-[\mathrm{IIIT1}]$ and $G_{j,0}$ is
        \[
        S-[\mathrm{IIIT1B4}] := \xygraph{
    \diamond ([]!{-(0,-.3)} {+0})
    - [r]  \circ ([]!{-(0,-.3)} {4})
      - [r]   \bullet ([]!{-(0,-.3)}  {1} )
      ( - [d]  \bullet ([]!{-(.3,0)} {2}) ,
       - [r]  \circ ([]!{-(0,-.3)} {4}) 
       }
        \]
        
    \end{enumerate}
\end{lem}

\begin{proof}
We see $b\in\{2,3,4,6\}$ from the classification of $\Q$-Gorenstein smoothable strictly lc rational singularities (Lemma~\ref{lem:smoothable_rational_strict_lc}).

From the assumption, Proposition \ref{amulet} implies that there exists some $k$ such that, after reordering the blow-ups if necessary, $G_{j,k}$ is one of the graphs described in Proposition \ref{amulet} (1), (2), and (4).
Moreover, $S+B_{j,k}$ is one of the following:
\begin{itemize}
    \item[(i)] $[\chi,2^l]$
    \item[(ii)] $[\chi,3+\beta,2^{n-2},3,2^{\beta}]$
    \item[(iii)] $[\chi,\gamma,2^{\gamma-4}]$
    \item[(iv)] $[\chi+2,2]$
    \item[(v)] $[\chi+1,3,2]$
\end{itemize}

We divide the proof into the following three cases. 
The first case is when no blow-up occurs on $B_{j,k}$. 
The second case is when a blow-up occurs at the intersection of $S$ and $B_{j,k}$. 
The third case is when neither of these conditions hold.

We consider the first case.
Then, $B_{j,k}=B_{j,0}$ must consist of a single $(-b)$-curve.
This situation corresponds exactly to cases $(1),(2),(3),(16),(17),(18)$ and $(19)$. 

We consider the second case.
Let $a+1$ be the number of blow-ups that occur on the infinitely near to points of $S$ over $S \cap B_{j,k}$. 
After the $(a+1)$-times blow-ups there, the string $S+B_{j,k}$ becomes $SR_1^{a+1}-1-L_2^{a}L_1 B_{j,k}$. 
Let $E$ be the $(-1)$-curve in the resulting string.
Since the self-intersection number of $E$ must eventually be one of $-2,-3,-4$, or $-6$, further blow-ups occur at $E \setminus S$.
After the blow-ups, the string becomes $SR_1^{a+1}R_2R_1^{b-2}-1-L_2^{b-2}L_1L_2^{a}L_1 B_{j,k}$.

The chain $H:=L_2^{b-2}L_1L_2^{a}L_1 B_{j,k}$ is one of the following chains: 
\begin{itemize}
    \item[(i)] $[2^{b-2},3,2^{a-1},3,2^{l-1}]$ 
    \item[(ii)] $[2^{b-2},3,2^{a-1},4+c,2^{d-1},3,2^{c}]$
    \item[(iii)] $[2^{b-2},3,2^{a-1},c+1,2^{c-4}]$
    \item[(iv)] $[2^{b-2},3,2^{a-1},3]$
    \item[(v)] $[2^{b-2},3,2^{a-1},4,2]$
\end{itemize}
Here, the labels (i)--(v) correspond to those used for the chain $S+B_{j,k}$.
The chain $H$ must be P-admissible.

For case (i), it follows that $b=2$ and $l=0$ by Lemma~\ref{lem:extT_I}~(3).
Hence, on the surface $\widetilde{X}$, there exists a chain of the form 
\[
[\chi+a+2,2] -1- [3,2^{a-1},3].
\]
The log discrepancy of the $(-2)$-curve corresponding to the string $[\chi+a+2,2]$ is $1/2$.
Let $E$ denote the $(-1)$-curve connecting the two chains above.  
Then we have $\pi(E)\cdot K_{X}=0$, contradicting the ampleness of $K_{X}$.

For case (ii), we use Lemma~\ref{drill:lem:amulet-b}~(1).
The possible values are either $(b,a,\beta,d)=(2,1,1,n+1)$, yielding $H=[3,2^m,5,2^n,3,2]$, or $(3,0,0,n+1)$, giving $H=[2,3,2^{m-1},4,2^n,3]$.

If $(b,a,\beta,d)=(2,1,1,n+1)$, then the associated T-train for $[3,2^m,5,2^n,3,2]$ is \[
[3,2^{m-2},3]-1-[3,2,5,2]-1-[4,2^{n-2},3,2]. 
\]
Here, if $n=0$, then we regard $[3,2,5,2]-1-[4,2^{n-2},3,2]$ as $[3,2,5,2]$, and if $m=0$, then we regard $[3,2^{m-2},3]-1-[3,2,5,2]$ as $[3,2,5,2]$.
We prove that $m = 0$. 
Suppose, to the contrary, that $m \geq 1$. 
Then the surface $\widetilde{X}$ contains a chain of the form
$[\chi+2,2] -1- [3,2^{m-2},3]$.
The log discrepancy of the $(-2)$-curve corresponding to the chain $[\chi+2,2]$ is $1/2$.  
Let $E$ denote the $(-1)$-curve connecting the two chains above.  
Then we have $\pi(E)\cdot K_{X}=0$, contradicting the ampleness of $K_{X}$.
Therefore, we conclude that $B_{j,k}$ must be a chain of the form $[4,2^n,3,2]$ for some $n \geq 0$.
This condition determines the admissible initial blowing up process: $[\mathrm{I}_{n+2}\mathrm{T}]_1$, $[\mathrm{IIIT2}]$, $[\mathrm{IVR}]_{1}$ or $[\mathrm{IVR}]_{2}$.
Furthermore, it holds that
\[
[\chi+a+1,b]-1-H = [\chi+2,2]-1-[3,5,2^n,3,2].
\]
We obtain $(5)$, $(7)$ and $(11)$.

If $(b,a,\beta,d)=(3,0,0,n+1)$, then the associated T-train for $[2,3,2^{m-1},4,2^n,3]$ is 
\[[2,3,2^{m-2},4]-1-[2,5,3]-1-[3,2^{n-2},3].\] 
Here, if $n=0$ we regard $[2,5,3]-1-[3,2^{n-2},3]$ as $[2,5,3]$,
if $m=0$ we regard $[2,3,2^{m-2},4]-1-[2,5,3]$ as $[2,5,3]$.
We prove that $m = 0$ by contradiction.  
Suppose, to the contrary, that $m \geq 1$. Then the surface $\widetilde{X}$ contains a chain of the form
$[\chi+1,3] -1- [2,3,2^{m-2},4]$.
Let $E$ denote the $(-1)$-curve connecting the two chains above.  
Then we have $\pi(E)\cdot K_{X}=0$, contradicting the ampleness of $K_{X}$.
It follows that $B_{j,k}$ must be the form $[3,2^n,3]$ for some $n \geq 0$.
This condition determines the admissible initial blowing up process: $[\mathrm{IIIT1}]$, $[\mathrm{IVT1}]$, $[\mathrm{I}_{n+2}\mathrm{T}]_{0}$,
$[\mathrm{I}_{l}^*\mathrm{R}]_{n+1}$.
Furthermore, it holds 
\[
[\chi+a+1,b]-1-H = [\chi+1,3]-1-[2,5,2^n,3].
\]
Note that $[\mathrm{I}_{l}^*\mathrm{R}]_{n+1}$ ($n \geq 1$) is also excluded using the fact that $K_{X}$ is ample.
We obtain $(12)$, $(13)$, $(14)$ and $(15)$.

For case (iii), we use Lemma~\ref{drill:lem:amulet-b}~$(2)$.
The possible values are either $(b,a,c)=(2,m+2,5)$, yielding $H=[3,2^{m+1},6,2]$, or $(4,0,4)$, giving $H=[2,2,6]$.

If $(b,a,c)=(2,m+2,5)$, then the associated T-train $[3,2^{m+1},6,2]$ is $[3,2^{n-2},3]-1-[3,2,6,2]$.
Here, if $n=0$ we regard $[3,2^{n-2},3]-1-[3,2,6,2]$ as $[3,2,6,2]$.
We prove that $n = 0$ by contradiction.  
Suppose, to the contrary, that $n \geq 1$. 
Then the surface $\widetilde{X}$ contains a chain of the form
$[\chi+2,2] -1- [3,2^{n-2},3]$.
Let $E$ denote the $(-1)$-curve connecting the two chains above.  
Then we have $\pi(E)\cdot K_{X}=0$, contradicting the ampleness of $K_{X}$.
Therefore, it follows $B_{j,k}=[5,2]$.
This condition determines the admissible initial blowing up process: $[\mathrm{I}_{1}\mathrm{T}]_1$, $[\mathrm{IIT2}]$, $[\mathrm{IIIT3}]$ and $[\mathrm{IVR}]_{0}$.
Furthermore, it holds that
\[
[\chi+a+1,b]-1-H = [\chi+3,2]-1-[3,2,6,2].
\]
In this case, we obtain $(4)$, $(6)$, $(8)$ and $(10)$.

If $(b,a,c)=(4,0,4)$,
then it follows that $B=[2,2,6]$.
Furthermore, it holds that
\[
[\chi+a+1,b]-1-H = [\chi+1,4]-1-[2,2,6].
\]
This case is excluded using the fact that $K_{X}$ is ample.

For case (iv), we use Lemma~\ref{lem:extT_I}~$(3)$.
The possible values are $(b,a)=(2,a)$, yielding $H=[3,2^{a-1},3]$.
The chain $H$ is a T-chain.
It follows that $B_{j,k}=[2]$.
Furthermore, it holds that
\[
[\chi+a+3,b]-1-C = [\chi+3,2]-1-[4].
\]
Therefore, this case is excluded using the fact that $K_{X}$ is ample.

For case (v), we use Lemma~\ref{drill:lem:amulet-b}~$(3)$
The possible values are $(b,a)=(2,m)$, yielding $H=[3,2^{m-1},4,2]$.
The associated T-chain of $H$ is a T-train $[3,2^{m-2},3]-1-[5,2]$.
It follows that $B_{j,k}=[3,2]$.
This condition determines the admissible initial blowing up process $[\mathrm{I}_n \mathrm{R}]_1$.
Furthermore, it holds that
\[
[\chi+a+2,b]-1-H = [\chi+2,2]-1-[3,2^{m-1},4,2].
\]
Since the fact that $K_{X}$ is ample, $m=0$.
We obtain $(9)$.

We consider the third case.
In this case, a blow-up must occur at the intersection of the curves corresponding to $b_1$ and $b_2$ in the string $B=[b_1,b_2,\cdots,b_r]$.
After this blow-up,
the string $S+B_{j,k}$ becomes
$SR_2 R_1^{b_1+a-1}-1-L_{2}^{a}L_1[b_2,\cdots,b_r]$.
The chain $H=L_2^{a}L_1[b_2, \dots, b_r]$ is one of the following:
\begin{itemize}
    \item[(i)] $[2^a,3,2^{b-2}]$
    \item[(ii)] $[2^a,3,2^{d-2},3,2^{c}]$
    \item[(iii)] $[2^a,3,2^{c-5}]$
    \item[(iv)] $[2^a,3]$
\end{itemize}
Since the chain $H$ must be P-admissible, we can exclude cases (i), (iii), and (iv).  
Moreover, in case (ii), we see that the possible values are $(a,d,c) = (2, n+1, 0)$, which yields a unique possibility for $H$, namely $H = [3, 2^{n-1}, 3]$.
The chain $H$ is a T-chain.
It follows that $B_{j,k}=[3,2^{n},3]$.
This condition determines the admissible initial blowing up process: $[\mathrm{IIIT1}]_1$, $[\mathrm{IVT1}]$, $[\mathrm{I}_{n+2}\mathrm{T}]_0$ and $[\mathrm{I}_{l}^* \mathrm{R}]_{n+1}$.
Furthermore,
\[
SR_2 R_1^{2}-1-L_{2}^{2}L_1[b_2,\cdots,b_r] = [\chi,4]-1-[3,2^{n-1},3].
\]
Note that $[\mathrm{IVT1}]$ and $[\mathrm{I}_{l}^*\mathrm{R}]_{n+1}$ are excluded using the fact that $K_{X}$ is ample.
We obtain $(20)$ and $(21)$.
\end{proof}

\begin{defn}\label{defn--set-amulet-b}
We define sets of amulets $\cA_{[2]}^{\text{B}}$, $\cA_{[3]}^{\text{B}}$, $\cA_{[4]}^{\text{B}}$, and $\cA_\SR$ as follows:
\begin{align*}
&\cA_{[2]}^{\text{B}} := \{  [\mathrm{CT}]_1, [\mathrm{I}_1\mathrm{T}_1\mathrm{B2}], [\mathrm{IIT2B2}], [\mathrm{IIIT2B2}], [\mathrm{IIIT3B2}], [\mathrm{IVR}_0\mathrm{B2}], [\mathrm{IVR}_1\mathrm{B2}], [\mathrm{IVR}_2\mathrm{B2}] \} \\
&\quad\quad\quad\cup \{ [\mathrm{I}_n \mathrm{R}]_0 : n\geq3 \} \cup \{[\mathrm{I}_n \mathrm{E}]_1:n\geq 0\} \cup \{ [\mathrm{I}_n \mathrm{T}_1B2]:n\geq2 \} \cup \{[\mathrm{I}_n\mathrm{R}_1\mathrm{B2}]:n\geq3\} \\
&\cA_{[3]}^{\text{B}} := \{ [\mathrm{IIIT1B3}], [\mathrm{IVT1B3}] \} \cup \{ [\mathrm{I}_n \mathrm{T}_0\mathrm{B3}]:n \geq 2\}\cup \{ [\mathrm{I}_{n}^*\mathrm{R}_1 \mathrm{B3}]:n \ge 0\}    \\
&\cA_{[4]}^{\text{B}} := \{ [\mathrm{I}_1\mathrm{T}]_0, [\mathrm{IIT1}], [\mathrm{IVT2}], [\mathrm{IIIT1B4}] \} \cup \{ [\mathrm{I}_n^* \mathrm{R}]_0 : n\geq0 \} \cup \{ [\mathrm{I}_n \mathrm{T}_0\mathrm{B4}] : n\geq2 \}  \\
&\cA_\SR := \{[\mathrm{CSR}], [\mathrm{I}_n^*\mathrm{SR1}],[\mathrm{I}_n^*\mathrm{SR2}], [\mathrm{II}^* \mathrm{SR}], [\mathrm{III}^*\mathrm{SR}], [\mathrm{IV}^*\mathrm{SR}]\} 
\end{align*}

From Lemma \ref{lem:amulet-b}, if the connected component of $D_{j,0}$ intersecting $S$ is a $(-b)$-curve, then there exists exactly one amulet $A\in\cA_{[b]}^{\text{B}}$ such that the dual graph $G_{D_{j,0}}$ coincides with the graph $S-A$.

For $A\in\cA_{[b]}^{\text{B}}$ ($b=2,3,4$) other than $[\mathrm{I}_n \mathrm{E}]_1$, we define $\gamma(A)$ as the number assigned to $S$ in the decorated graph $S-A$, and define $\gamma([\mathrm{I}_n \mathrm{E}]_1)=0$.
For example, $\gamma( [\mathrm{CT}]_1)=0$, and $\gamma([\mathrm{I}_1\mathrm{T}_1\mathrm{B2}])=3$.
For $A\in\cA_\SR$, we define $m(A)$ to be the number of components of $A$ minus $2$.
We list $m(A)$ in the following table:
\begin{table}[H]
    {\renewcommand\arraystretch{1.5}
    \begin{tabular}{|c||c|c|c|c|c|c|} \hline
     $ A $                     &  $[\mathrm{CSR}]$   & $[\mathrm{I}_n^*\mathrm{SR1}]$ & $[\mathrm{I}_n^*\mathrm{SR2}]$  & $[\mathrm{II}^* \mathrm{SR}]$ &$[\mathrm{III}^*\mathrm{SR}]$ & $[\mathrm{IV}^*\mathrm{SR}]$\\ \hline 
    $m(A)$ & $1$ & $2$  &  $n+2$                                                        & $6$ & $4$ & $3$ \\ \hline
    \end{tabular}
  }
\end{table}
\end{defn}

We first assume that the dual graph of $\widetilde{C}^{(p)}$ corresponds to the one described in
Lemma~\ref{lem:smoothable_rational_strict_lc}~(i).
Depending on the position of $S$ within $\widetilde{C}^{(p)}$, there are four possible types:
\begin{itemize}
    \item $S$ has exactly four branches, each of which is a $(-2)$-curve.
    \item $S$ has exactly three branches:
    two are $(-2)$-curves, and the other one has a fork.
    \item $S$ has exactly two branches, both of which have forks.
    \item $S$ is located at the end of $\widetilde{C}^{(p)}$.
\end{itemize}
Applying Lemma \ref{strlcratbranch} and Lemma \ref{lem:amulet-b}, we obtain the following classification.

\begin{thm}\label{thm:one-section-halfcusp}
Let $(g,n_X,l_{E_{h}})=(0,1,0)$.
Assume that $p$ is a strictly lc rational singularity of type $(2,2,2,2)$ and there exist no branches that are completely separated from the section $S$.
Then, exactly one of the following holds:
\begin{itemize}
\item[$1.$] $(S-A_{[2]}^{\operatorname{B}},-A_{[2]}^{\operatorname{B}},-A_{[2]}^{\operatorname{B}},-A_{[2]}^{\operatorname{B}})$
\begin{itemize}
    \item $J=4$, and $G_{D_{j,0}}=S-A_{j}$ where $A_{j}\in\cA_{[2]}^{\text{B}}$ for $1\leq j\leq4$.
    \item $2 \leq\chi\leq6-\sum_{j=1}^4 \gamma(A_{j})$. 
    \item The dual graph of $B_0$ is
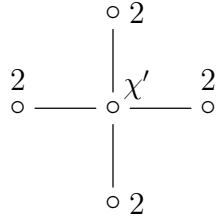
\[
\xygraph{
    \circ ([]!{-(0,-.3)} {2}) - [r]
    \circ ([]!{-(-0.25,-.25)} {\chi'}) (
         - [u] \circ ([]!{-(-0.25,-.0)} {2}),
        - [d] \circ ([]!{-(-0.25,-.0)} {2}), 
         - [r] \circ ([]!{-(0,-.3)} {2})
)}
\]
where $\chi':=\chi+\sum_{j=1}^4 \gamma(A_{j})$. 
\end{itemize}

\item[$2.$] $(S-A_{[2]}^{\operatorname{B}},-A_{[2]}^{\operatorname{B}},-A_{\SR})$ 
\begin{itemize}
    \item $J=3$, $G_{D_{1,0}}=S-A_{1}$, $G_{D_{2,0}}=S-A_{2}$, and $G_{D_{3,0}}=S-A$ where $A_{1}, A_{2} \in \cA_{[2]}^{\text{B}}$ and $A \in\cA_\SR$.
    \item $2 \leq\chi\leq 6+m(A)-\gamma(A_{1})-\gamma(A_{2})$
    \item The dual graph of $B_0$ is
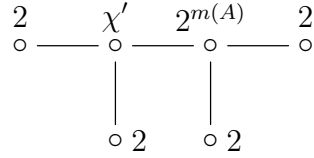
\[
\xygraph{
    \circ ([]!{-(0,-.3)} {2}) - [r]
    \circ ([]!{-(0,-.3)} {\chi'}) (
        - [d] \circ ([]!{-(-0.25,-.0)} {2}),
         - [r] \circ ([]!{-(0,-.3)} {2^{m(A)}}) (
         - [d] \circ ([]!{-(-0.25,-.0)} {2}),
        - [r] \circ ([]!{-(0,-.3)} {2})
)}
\]
where $\chi':=\chi+\gamma(A_{1})+\gamma(A_{2})$.
\end{itemize}

\item[$3.$] $(A_{\SR}-S-A_{\SR})$ 
\begin{itemize}
    \item $J=2$, $G_{D_{1,0}}=S-A_{1}$ and $G_{D_{2,0}}=S-A_{2}$ where $A_{1},A_{2}\in\cA_\SR$.
    \item $2 \leq \chi = \chi' \leq 6+m(A_{1})+m(A_{2})$.
    \item The dual graph of $B_0$ is
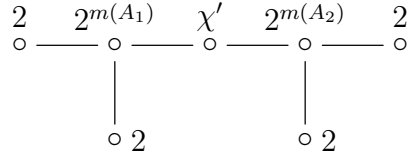
\[
\xygraph{
    \circ ([]!{-(0,-.3)} {2}) - [r]
    \circ ([]!{-(0,-.3)} {2^{m(A_{1})}}) (
        - [d] \circ ([]!{-(-0.25,-.0)} {2}),
        - [r] \circ ([]!{-(-0,-.3)} {\chi'})
         - [r] \circ ([]!{-(0,-.3)} {2^{m(A_{2})}}) (
         - [d] \circ ([]!{-(-0.25,-.0)} {2}),
        - [r] \circ ([]!{-(0,-.3)} {2})
)}
\]
\end{itemize}

\end{itemize}
\end{thm}

In the above, we impose the condition that no branch is completely separated from $S$, but it is easy to extend the above classification without this assumption; see Remark~\ref{rem:ignore_compsep_step4}.

Next, we assume that the dual graph of $\widetilde{C}^{(p)}$ corresponds to the one described in Lemma \ref{lem:smoothable_rational_strict_lc} (ii), (iii), or (iv).
Depending on the position of $S$ within $\widetilde{C}^{(p)}$, there are two possible types:
\begin{itemize}
    \item $S$ is the fork of $\widetilde{C}^{(p)}$.
    \item $S$ is a branch with respect to the fork of $\widetilde{C}^{(p)}$.
\end{itemize}
Applying Lemma \ref{strlcratbranch} and Lemma \ref{lem:amulet-b} again, we obtain the following classification.

\begin{thm}\label{thm:one-section-triangle}
Let $(g,n_X,l_{E_h})=(0,1,0)$.
Assume that the dual graph of $\widetilde{C}^{(p)}$ coincides with the one described in Lemma \ref{lem:smoothable_rational_strict_lc} $\mathrm{(ii)}$, $\mathrm{(iii)}$ or $\mathrm{(iv)}$, and there exist no branches that are completely separated from the section $S$.
Then, exactly one of the following holds:
\begin{itemize}
    \item[$1.$] $(S-A_{[3]}^{\operatorname{B}},-A_{[3]}^{\operatorname{B}},-A_{[3]}^{\operatorname{B}})$
\begin{itemize}
    \item $J=3$, and $G_{D_{j,0}}=S-A_j$ where $A_j\in\cA_{[3]}^{\text{B}}$ for $1\leq j\leq3$.
    \item $2 \leq\chi\leq4-\sum_{j=1}^3 \gamma(A_j)$.
    \item The dual graph of $B_0$ is
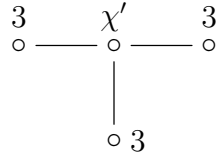
\[
\xygraph{
    \circ ([]!{-(0,-.3)} {3}) 
    - [r]  \circ ([]!{-(0,-.3)} {\chi'})(
        - [r] \circ ([]!{-(0,-.3)} {3}),
        - [d] \circ ([]!{-(-0.25,-.0)} {3}),
}
\]
where $\chi':=\chi+\sum_{j=1}^3 \gamma(A_j)$.
\end{itemize}

\item[$2.$] $(S-A_{[2]}^{\operatorname{B}},-A_{[4]}^{\operatorname{B}},-A_{[4]}^{\operatorname{B}})$
\begin{itemize}
    \item $J=3$, $G_{1,0}=S-A_1$, $G_{2,0}=S-A_2$, and $G_{3,0}=S-A_3$ where $A_1\in\cA_{[2]}^{\text{B}}$ and $A_2,A_3\in\cA_{[4]}^{\text{B}}$.
    \item $2 \leq\chi\leq3-\gamma(A_1)-\gamma(A_2)-\gamma(A_3)$.
    \item The dual graph of $B_0$ is
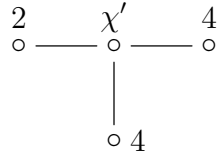
\[
\xygraph{
    \circ ([]!{-(0,-.3)} {2}) 
    - [r]  \circ ([]!{-(0,-.3)} {\chi'})(
        - [r] \circ ([]!{-(0,-.3)} {4}),
        - [d] \circ ([]!{-(-0.25,-.0)} {4}),
}
\]
where $\chi':=\chi+\gamma(A_1)+\gamma(A_2)+\gamma(A_3)$.
\end{itemize}

\item[$3.$] $(S-[\mathrm{I}_n \mathrm{SR1}])$
\begin{itemize}
    \item $J=1$, and $G_{1,0}=S-[\mathrm{I}_n \mathrm{SR1}]$ for some $n\geq3$.
    \item $\chi=\chi'=3$.
    \item The dual graph of $B_0$ is
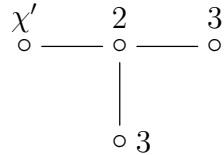
\[
\xygraph{
    \circ ([]!{-(0,-.3)} {\chi'}) 
    - [r]  \circ ([]!{-(0,-.3)} {2})(
        - [r] \circ ([]!{-(0,-.3)} {3}),
        - [d] \circ ([]!{-(-0.25,-.0)} {3}),
}
\]
\end{itemize}
\item[$4.$] $(S-[\mathrm{I}_n \mathrm{SR2}])$
\begin{itemize}  
    \item $J=1$, and $G_{1,0}=S-[\mathrm{I}_n \mathrm{SR2}]$ for some $n\geq3$.
    \item $\chi=\chi'=4$.
    \item The dual graph of $B_0$ is
\[
\xygraph{
    \circ ([]!{-(0,-.3)} {\chi'}) 
    - [r]  \circ ([]!{-(0,-.3)} {3})(
        - [r] \circ ([]!{-(0,-.3)} {2}),
        - [d] \circ ([]!{-(-0.25,-.0)} {4}),
}
\]
\end{itemize}
\end{itemize}
\end{thm}

\begin{rem}\label{rem:ignore_compsep_step4}
    A similar argument in Remark~\ref{ignorecompsep} shows that, by substituting $\chi' + 2J_{\operatorname{CS}}$ for $\chi'$ in the lists provided in Theorems~\ref{thm:one-section-halfcusp} and \ref{thm:one-section-triangle}, they give the classification of $X$ with $J_{\operatorname{CS}}$ completely separated branches of the form $(1)\text{--}(5)$ in Proposition \ref{compsepamulet}.
\end{rem}

\subsection{Case of $(g,n_X,l_{E_h})=(0,1,0)$ with an involution}\label{subsec:0,1,0,involution}

In this subsection, we classify $C$ and its resolution $\widetilde{C}$ under the following assumption.

\begin{ass}\label{ass:5-2}
    \phantom{A}
    \begin{itemize}
        \item[$(1)$] $(g,n_X,l_{E_h})=(0,1,0)$.
        \item[$(2)$] $X$ has an involution $\sigma$ which induces a fiberwise involution on $\widetilde{X}\to \PP^1$.
        \item[$(3)$] If the quotient $W:=X/\sigma$ is klt, then $W$ admits a (possibly trivial) P-resolution $W^\dagger\to W$ with $K_{W^\dagger}^2\geq8$.
        \item[$(4)$] The projection $\widetilde{X}/\sigma \to \mathbb{P}^1$ is a ruling over $\mathbb{P}^1$.
        \item[$(5)$] The singularity $p\in X$ is not a strictly lc singularity of type $(3,3,3)$. 
    \end{itemize}
\end{ass}

\begin{rem}
    \phantom{A}
    \begin{itemize}
        \item[$(1)$] 
        If $X$ is a non-standard Horikawa surface with a good involution $\sigma$, then $\sigma$ satisfies the conditions $(2)$--$(4)$ in Assumption \ref{ass:5-2} by the definition of a good involution.
        \item[$(2)$] 
        Assumption \ref{ass:5-2} $(5)$ is imposed to avoid technical difficulties.
        When $X$ is a non-standard Horikawa surface with $(g,n_X,l_{E_h})=(0,1,0)$, then Assumption \ref{ass:5-2} $(5)$ holds.
        Indeed, if $p\in X$ is a strictly lc singularity of type $(3,3,3)$, the result of the previous subsection (see Theorem \ref{thm:one-section-triangle}) implies that $\chi<4$, which contradicts the fact that $X$ is Horikawa.
        Thus, we may impose Assumption \ref{ass:5-2} for the classification of non-standard Horikawa surfaces with a good involution.
    \end{itemize}
\end{rem}


The involution $\sigma$ induces involutions on the minimal resolution $\widetilde{X}$ and on its relatively minimal model $Y_m=Y$.
The main difference from the setting in the previous subsection is the existence of a certain involution on $X$.
This leads to further refinement because many initial blow-up processes are excluded. 

We mostly follow the classification procedure in the previous subsection. 
We will deal with this classification problem in the following steps: 
\begin{enumerate}[label=\textbf{(Step~\arabic*)}]
    \item We classify the initial blow-up processes as in Proposition \ref{amulet} under Assumption \ref{ass:5-2}.
    The possible initial blow-up processes are restricted to those compatible with the involution and with Assumption \ref{ass:5-2} $(3)$, which imposes the constraint on the quotient by the involution.
    The result is given in Proposition \ref{prop:sigma-amulet}.
    \item We classify the blow-up processes as in Proposition \ref{compsepamulet} that are compatible with the involution on $Y$.
    See Section \ref{subsubsec:0,1,0-inv-step2}.
    \item Under Assumption \ref{ass:5-2}, we classify the configurations of curves in $C$ and their resolutions when $p$ is a T-singularity.
    See Section \ref{subsubsec:0,1,0-inv-step3}.
    \item Under Assumption \ref{ass:5-2}, we classify the configurations of curves in $C$ and their resolutions when $p$ is a strictly lc rational singularity.
    See Section \ref{subsubsec:0,1,0-inv-step4}.
\end{enumerate}

\subsubsection{Step 1: initial blow-up process}

In this step, we classify all possible initial blow-up processes as in Section \ref{subsubsec:0,1,0-step1}.
The overall strategy is the same as in Section \ref{subsubsec:0,1,0-step1}, but the assumption of the existence of a suitable involution on $X$ allows for a further refinement.

We begin by introducing the necessary notions, which will be used in the statement of Lemma \ref{lem:restrict_involution}, a key ingredient in this subsection.

\begin{defn}\label{defn:triple}
Let $q\in Z$ be a germ of a cyclic quotient singularity.
For a P-resolution $Z^\dagger\to Z$, let $\rho(Z^\dagger/Z)$ denote the number of irreducible components of $\pi$-exceptional curves, and let $\mu(Z^\dagger)$ be the sum of the Milnor numbers of singularities of $Z^\dagger$. 
We define $\lambda(q)$ to be the minimum of $\rho(Z^\dagger/Z)+\mu(Z^\dagger)$ over all P-resolutions $Z^\dagger \to Z$.
For a surface $Z$ with only cyclic quotient singularities, we define $\lambda(Z)$ to be the sum of $\lambda(q)$ over all cyclic quotient singularities $q\in Z$.
\end{defn}

We note that $q\in Z$ is a Wahl singularity if and only if $\lambda(q)=0$.

\begin{lem}\label{lem:restrict_involution}
    Assume that $X$ has an involution $\sigma$ and the quotient $W=X/\sigma$ is a rational surface with only cyclic quotient singularities.
    For any (possibly trivial) P-resolution $W^\dagger\to W$, it holds that 
    \[
    \lambda(W)\leq9-K_{W^\dagger}^2.
    \]
    In particular, if $K_{W^\dagger}^2\geq8$, $W$ has no singular points $q$ with $\lambda(q)\geq2$, and there exists at most one singular point $q\in W$ with $\lambda(q)=1$.
\end{lem}

\begin{proof}
    The formula in \cite[Proposition 6.1]{Pro} for the rational surface $W^\dagger$ yields the equation $\rho(W^\dagger/W)+\mu(W^\dagger)=10-K^2_{W^{\dagger}}-\rho(W)$.
    This gives the desired inequality.
\end{proof}

We introduce the \emph{$\sigma$-decorated graph}, the decorated graph enriched with the information of the involution $\sigma$.
From Assumption \ref{ass:5-2}, the involution $\sigma$ acts on $Y$ fiberwise and preserves $C$, and in particular, it preserves each $D_j$, thereby inducing an involution on $D_j$.  
Note that this involution does not necessarily lift to each $D_{j,k}$.
Let $F_{j,k}$ denote the total transform of the fiber $F_j$ on $Y_k$.

\begin{defn}
Suppose that the involution $\sigma$ lifts locally to a germ of $F_{j,k}$. 
Then we define the \emph{$\sigma$-decorated graph} $G_{j,k}^{\sigma}$ as a tuple consisting of the following data:
\begin{itemize}
    \item the decorated graph $G_{j,k}$ associated to $D_{j,k}$ (see Definition~\ref{defn:graph});
    \item the involution on $G_{j,k}$ induced by $\sigma$;
    \item for each horizontal one-dimensional $\sigma$-fixed locus, the irreducible components (possibly more than one) in the fiber that intersect the $\sigma$-fixed locus;
    \item for each isolated $\sigma$-fixed point, the irreducible components (possibly more than one) in the fiber that contain the $\sigma$-fixed point.
\end{itemize}
\end{defn}

\begin{exam}
We describe the $\sigma$-decorated graph naturally arising from the decorated graph $G_{j,k} = S - [\mathrm{IIIR}]_0$. 
From the classification of involutions on the singular fiber of type III (Proposition \ref{prop_inv_ell}), one sees that the involution $\sigma$ lifts locally to the germ of the fiber $F_{j,k}$, and we obtain the $\sigma$-decorated graph $G_{j,k}^\sigma$.
Following Notation \ref{fig:involution}, we present the data of $G_{j,k}^\sigma$ below.

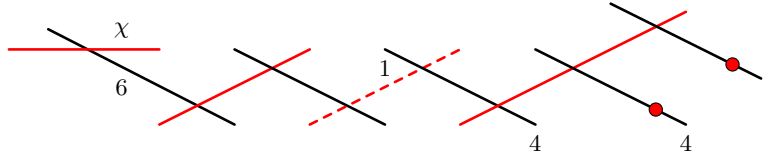
\begin{center}
\begin{tikzpicture}[line cap=round,line join=round,>=triangle 45,x=1cm,y=1cm]
\clip(-4,-1) rectangle (8,2);
\draw [line width=1pt] (-2.476872806294636,1.2686349061812408)-- (0,0);
\draw [line width=1pt,red] (-1,0)-- (1,1);
\draw [line width=1pt] (0,1)-- (2,0);
\draw [line width=1pt,red,dashed] (1,0)-- (3,1);
\draw [line width=1pt,red] (-3,1)-- (-1,1);
\draw [line width=1pt] (2,1)-- (4,0);
\draw [line width=1pt,red] (3,0)-- (6,1.5);
\draw [line width=1pt] (4,1)-- (6,0);
\draw [line width=1pt] (5,1.608733649129835)-- (7,0.6127301876375234);
\begin{scriptsize}
\draw [black] (-1.5,1.25) node {$\chi$};
\draw [black] (-1.5,0.5) node {$6$};
\draw [black] (2,0.75) node {$1$};
\draw [black] (4,-0.25) node {$4$};
\draw [black] (6,-0.25) node {$4$};
%
\draw [fill=red] (5.6,0.2) circle (2.5pt);
\draw [fill=red] (6.62,0.8) circle (2.5pt);
\end{scriptsize}
\end{tikzpicture}
\end{center}
\end{exam}

For each amulet $A$ appearing in Proposition \ref{amulet}, we define a \emph{$\sigma$-amulet} $A^{\sigma}$ in a way similar to the above example.
Let $G_{j,k} = S - A$ denote the initial blow-up process. 
If the involution $\sigma$ lifts locally to a germ of $F_{j,k}$, then the $\sigma$-decorated graph $G_{j,k}^\sigma$ is well-defined and we define $S - A^\sigma := G_{j,k}^\sigma$.
By the classification of fibers of elliptic surfaces with an involution (Proposition \ref{prop_inv_ell}) and that of initial blow-up processes (Proposition \ref{amulet}), we see that this condition holds for all amulets $A$ except $[\mathrm{I}_n \mathrm{T}]_\beta$ and $[\mathrm{IVT1}]$.

When the involution $\sigma$ does not lift to a germ of the fiber $F_{j,k}$, that is, when $A=[\mathrm{I}_n \mathrm{T}]_\beta$ or $[\mathrm{IVT1}]$, we define $S-A^\sigma$ to be the $\sigma$-decorated graph $G^\sigma_{j,k'}$ where $Y_{k'}\to Y_k$ is the minimal blow-ups of $Y_k$ to a germ of which $\sigma$ lifts.

We first consider the case $A=[\mathrm{I}_n \mathrm{T}]_\beta$.
As described in the proof of Lemma \ref{modifyI}, the graph $S-[\mathrm{I}_n \mathrm{T}]_\beta$ is obtained by blow-up $(\beta+1)$-times over a node $q$ of $D'_j$.
On the other hand, since the section $S$ is fixed, the involution must interchange two nodes $q,q'$ on $D'_j\subset Y$; see Proposition \ref{prop_inv_ell}.
Hence, the morphism $Y_{j,k'}\to Y_{j,k}$ is the composition of the $(\beta+1)$-times blow-ups over the node $q'$, and the decorated graph $G_{j,k'}$ is as follows:
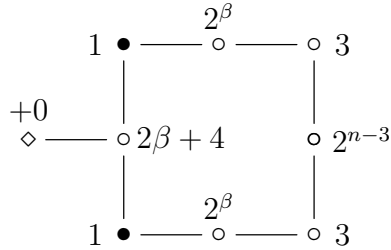
\[
\xygraph{
    \diamond ([]!{-(0,-.3)} {+0}) - [r]
    \circ ([]!{-(-0.6,-.0)} {2\beta+4}) ( 
        - [u] \bullet ([]!{-(.3,0)} {1}) 
        - [r] \circ ([]!{-(-0,-.3)} {2^\beta})
        - [r] \circ ([]!{-(-0.3,-.0)} {3})
        - [d] \circ ([]!{-(-0.5,-.0)} {2^{n-3}})([]!{-(.5,0)} {}),
        - [d] \bullet ([]!{-(.3,0)} {1})
        - [r] \circ ([]!{-(-0,-.3)} {2^\beta})
        - [r] \circ ([]!{-(-0.3,-.0)} {3}) 
        - [u] \circ ([]!{-(.3,0)} {})
    )}
\]
We see that $\sigma$ lifts locally to a germ of $F_{j,k'}$.
Since the chain $[2^\beta,3,2^{n-3},3,2^\beta]$ is not P-admissible if $\beta>0$ (Lemma \ref{lem:extT_I}~(3)), we only consider the $\sigma$-amulets $[\mathrm{I}_n \mathrm{T}]_0^\sigma$.

Let $A=[\mathrm{IVT1}]$.
By an argument similar to the case $A=[\mathrm{I}_n \mathrm{T}]_\beta$, one sees that the morphism $Y_{k'}\to Y_k$ is the blow-up at the point where $\sigma$ is not defined, and the resulting decorated graph $G_{j,k'}$ is as follows:
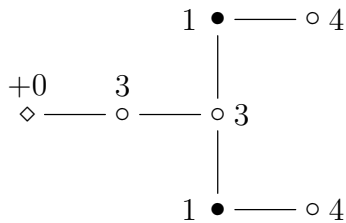
\[
\xygraph{
    \diamond ([]!{-(0,-.3)} {+0}) - [r]
      \circ ([]!{-(0,-.3)} {3}) - [r]
    \circ ([]!{-(-0.25,-.0)} {3}) (
        - [d] \bullet ([]!{-(.3,0)}{1})- [r] \circ ([]!{-(-0.25,-.0)} {4}),
        - [u] \bullet ([]!{-(.3,0)} {1}) - [r] \circ ([]!{-(-0.25,-.0)} {4})
    )}
\]

\begin{defn}
    We define the following sets of $\sigma$-amulets:
    \begin{align*}
        \cAT^\sigma &= \{ [\mathrm{IIT1}]^\sigma, [\mathrm{IIT2}]^\sigma, [\mathrm{IIT3}]^\sigma, [\mathrm{IIT4}]^\sigma, [\mathrm{IIIT1}]^\sigma, [\mathrm{IIIT2}]^\sigma, [\mathrm{IIIT3}]^\sigma, [\mathrm{IVT1}]^\sigma, [\mathrm{IVT2}]^\sigma \} \\
        &\quad \cup \{ [\mathrm{CT}]_\beta^\sigma : \beta\geq0\} \cup \{[\mathrm{I}_n \mathrm{T}]_{0}^\sigma : n\geq1 \}, \\
        \cAR^\sigma &= \{ [\mathrm{IIR}]^\sigma \} \cup \{ [\mathrm{I}_n \mathrm{R}]_\beta^\sigma : n\geq3, \beta=0,1 \} \cup \{ [\mathrm{IIIR}]_\beta^\sigma : \beta=0,1 \} \cup \{ [\mathrm{IVR}]_\beta^\sigma : \beta=0,1,2 \} \\
        &\quad \cup \{ [\mathrm{I}_n^* \mathrm{R}]_\beta^\sigma : n\geq0, 0\leq\beta\leq n+3 \}, \\
        \cASR^\sigma &= \{ [\mathrm{CSR}]^\sigma, [\mathrm{II}^* \mathrm{SR}]^\sigma, [\mathrm{III}^*\mathrm{SR}]^\sigma,[\mathrm{IV}^*\mathrm{SR}]^\sigma \} \cup \{ [\mathrm{I}_n\mathrm{SR}]^\sigma : \chi=3, n\geq3 \} \\
        &\quad \cup \{ [\mathrm{I}_n^*\mathrm{SR1}]^\sigma : n\geq0 \} \cup \{ [\mathrm{I}_n^*\mathrm{SR2}]^\sigma : n\geq1 \}, \\
        \cAE^\sigma &= \{ [\mathrm{I}_n \mathrm{E}]_\beta^\sigma : n\geq0, 0\leq\beta\leq n+8 \}, \\
        \cA^\sigma &= \cAT^\sigma \cup \cAR^\sigma \cup \cA_\SR^\sigma \cup \cAE^\sigma.
    \end{align*}
\end{defn}

The following lemma guarantees that $\lambda(W)$ is well-defined, i.e., $W$ is a surface with only cyclic quotient singularities. 

\begin{lem}
    Under Assumption~\ref{ass:5-2}, the quotient $W$ is klt.
\end{lem}

\begin{proof}

Let $\varphi\colon X \to W$ be the quotient by the involution $\sigma$. 
To show that $W$ is klt, it is sufficient to prove that the image of any lc singularity of $X$ is a klt singularity of $W$.

We first consider a singularity $q\in X$ distinct from $p\in X$.
From Section \ref{subsec:0,1,0}---especially Theorem \ref{thm:one-section-T}, \ref{thm:one-section-halfcusp}, and \ref{thm:one-section-triangle}---we see that the dual graph of the exceptional curves of $q\in X$ is a maximal connected subgraph consisting of white circles in the decorated graph described in Propostion \ref{amulet}, \ref{compsepamulet}, and \ref{lem:amulet-b}.
A direct computation shows that, for each case, the image $\varphi(q)$ is a smooth point or a cyclic quotient singularity, provided that the involution $\sigma$ on $Y$ lifts to $\widetilde{X}$.
We give a detailed computation in one case; the remaining cases are analogous.
Let $q\in X$ be the singularity of type $(2,4,4)[3]$ arising from $S-[\mathrm{IIIR}]_1$.
Since the involution $\sigma$ preserves each irreducible component of the exceptional set of $q\in X$ and there exists no isolated fixed point on the exceptional set,  
the image $\varphi(q) \in W$ is a cyclic quotient singularity of type $[2,5,2]$.

We show that $\varphi(p)\in W$ is klt.
If $p\in X$ is a T-singularity, since $\sigma$ preserves each exceptional curve of $p$, one sees that the exceptional set of $\varphi(p)$ is again a chain.
This shows that $\varphi(p)$ is klt.

Suppose $p \in X$ is a strictly lc rational singularity. 
The classification is given in Theorem \ref{thm:one-section-halfcusp} and \ref{thm:one-section-triangle}.
We give one example from each list.
Let $X$ be as in Theorem \ref{thm:one-section-T} (1).
Then, $p\in X$ is a strictly lc singularity of type $(2,2,2,2)[\chi']$, and $\sigma$ preserves each irreducible component of the exceptional set of $p\in X$.
Moreover, $\sigma$ fixes the $(-\chi')$-curve but not $(-2)$-curves, and there exists no isolated fixed point on the exceptional set.
Hence $\varphi(p) \in W$ is a cyclic quotient singularity of type $[2\chi'-4]$.
Finally, let $X$ be as in Theorem \ref{thm:one-section-triangle} (2).
In this case, $p\in X$ is a strictly lc singularity of type $(2,4,4)[\chi']$, and $\sigma$ preserves each irreducible component of the exceptional set of $p\in X$.
The involution $\sigma$ fixes the $(-\chi')$-curve but not branches, and there exists no isolated fixed point on the exceptional set.
Hence $\varphi(p) \in W$ is a cyclic quotient singularity of type $[2,2\chi'-1,2]$.
\end{proof}

We introduce subsets of $\cA^\sigma$ as follows:

\begin{align*}
        \cAT^{\sigma,0} &= \{ [\mathrm{IIT2}]^\sigma, [\mathrm{IIIT1}]^\sigma, [\mathrm{IIIT3}]^\sigma, [\mathrm{IVT1}]^\sigma \} \cup \{ [\mathrm{CT}]_\beta^\sigma : \beta\geq1\},   \\
        \cAT^{\sigma,1} &= \{ [\mathrm{IIT1}]^\sigma, [\mathrm{IVT2}]^\sigma \} \cup \{[\mathrm{I}_n \mathrm{T}]_{0}^\sigma : n\geq1 \}, \\
        \cAR^{\sigma,0} &= \{ [\mathrm{IVR}]_0^\sigma,  \} \cup \{ [\mathrm{I}_n^* \mathrm{R}]_1^\sigma : n\geq0 \},   \\
        \cAR^{\sigma,1} &= \{ [\mathrm{IIIR}]_1^\sigma \} \cup \{ [\mathrm{I}_n^* \mathrm{R}]_\beta^\sigma : n\geq0, 0\leq\beta\leq n+3, \beta\neq1 \}, \\
        \cASR^{\sigma,0} &= \{ [\mathrm{CSR}]^\sigma, [\mathrm{I}_1^*\mathrm{SR1}]^\sigma, [\mathrm{II}^* \mathrm{SR}]^\sigma, [\mathrm{IV}^*\mathrm{SR}]^\sigma,  [\mathrm{III}^*\mathrm{SR}]^\sigma \} \cup \{ [\mathrm{I}_n^*\mathrm{SR1}]^\sigma : n\geq 0 \}   \\
        &\quad \cup \{ [\mathrm{I}_n^*\mathrm{SR2}]^\sigma : n\geq1 \},  \\
        \cASR^{\sigma,1} &= \{ [\mathrm{I}_n\mathrm{SR}]^\sigma : \chi=3, n\geq3 \}, \\
        \cAE^{\sigma,0} &= \{ [\mathrm{I}_n \mathrm{E}]_\beta : n\geq0, \beta\in\{0,1,7,8\} \}, \quad\quad\quad\quad \cAE^{\sigma,1} = \cAE^\sigma\setminus\cAE^{\sigma,0}, \\
        \cA^{\sigma,i} &= \cAT^{\sigma,i} \cup \cAR^{\sigma,i} \cup \cASR^{\sigma,i} \cup \cAE^{\sigma,i} \quad (i=0,1), \quad\quad        \cA^{\sigma,\geq2} = \cA^{\sigma}\setminus(\cA^{\sigma,0}\cup\cA^{\sigma,1}).
\end{align*}

\begin{lem}\label{lem:lambda}
    Under Assumption~\ref{ass:5-2},
    the following hold.
    \begin{itemize}
        \item[$(1)$] If $G_{j,k}^\sigma=S-A$ for $A\in\cA^{\sigma,1}$, then $\lambda(W)\geq1$.
        \item[$(2)$] If $G_{j,k}^\sigma=S-A$ for $A\in\cA^{\sigma,\geq2}$, then $\lambda(W)\geq2$.
    \end{itemize}
\end{lem}

\begin{proof}
Assume that $G_{j,k}^\sigma=S-A$ for some $A\in\cA^{\sigma,1}\cup\cA^{\sigma,\geq2}$.
Define $U$ as the open subset 
\[
U=Y_k\setminus(B_k \cup \bigcup_{j'\neq j}D_{j',k}).
\]
Since the morphism $Y_{0} \to Y_{k}$ is an isomorphism over $U$, we can naturally regard $U$ as an open subset of $Y_0=\widetilde{X}$.
Let $\varphi:U \to U/\sigma$ be the quotient morphism by the involution.
Since $\sigma$ may have isolated fixed points, the quotient space $U/\sigma$ can have $A_{1}$-singularities.
Let $\nu:\widetilde{U/\sigma } \to U/\sigma$ denote the minimal resolution, and let $E$ be the exceptional divisor of $\pi\colon U\to \pi(U)$.
Let $\tau:\widetilde{U/\sigma }\to V$ denote the composition of blowing-downs of $(-1)$-curves such that the reduced divisor 
$$
E_V:=(\tau_{*}\nu^{*}\varphi_{*}E)_{\mathrm{red}}
$$
contains no $(-1)$-curves.
\[
\xymatrix{
 &  \widetilde{U/\sigma } \ar[dr]^{\nu} \ar[dl]_{\tau}  &  \\
 V & &  U/\sigma  \\
}	
\]
Then, there exists a natural morphism $V\to Z:=\pi(U)/\sigma$ that provides the minimal resolution of the open subvariety $Z\subset W$.
The divisor $E_V$ corresponds to the exceptional curves of this resolution.

We will prove the statement for the cases $A=[\mathrm{IIIR}]_{\beta}^\sigma$ ($\beta=0,1$) and $A=[\mathrm{I}_{n}^{*}\mathrm{R}]_\beta^\sigma$ ($0\leq\beta\leq n+3, \beta\neq1$). 
The remaining cases follow by similar arguments.

We first assume $A=[\mathrm{IIIR}]_{\beta}^\sigma$ ($\beta=0,1$).
By direct computation, we find that $E_V$ consists of the exceptional curves over a single cyclic quotient singularity, and the associated string is given by $[2,2\beta+3,2]$.
Since this is not a T-chain, $Z$ does not admit a trivial P-resolution.
This implies that $\lambda(W)\geq1$.
Furthermore, if $\beta=0$, we obtain $\lambda(W)\geq2$.
To see this, suppose for contradiction that $\lambda(W)=1$.
Then, there would exist a P-resolution $Z^\dagger\to Z$ with $\rho(Z^\dagger/Z)=1$ and $\mu(Z^\dagger)=0$.
However, a P-resolution $Z^\dagger\to Z$ with $\rho(Z^\dagger/Z)=1$ is unique and has two $A_1$ singularities (cf.\ Proposition~\ref{charP-resol}~(2)).
This is a contradiction.

We finally assume that $A=[\mathrm{I}_{n}^{*}\mathrm{R}]_\beta^\sigma$ ($0\leq\beta\leq n+3, \beta\neq1$).
In this case, $E_V$ consists of the exceptional curves over a single quotient singularity $q\in Z$, and the associated string is $[2\beta+2]$.
If $\beta>1$, then the P-resolution of $Z$ is unique and coincides with the minimal resolution $Z^\dagger\to Z$.
Thus, we obtain $\lambda(W)\geq\lambda(Z)=1$.
If $\beta=0$, then $q\in Z$ is an $A_1$-singularity.
It is straightforward to check that $\lambda(q\in Z)=1$.
This completes the proof.
Note that if $\beta=1$, then this singularity is a T-singularity of Milnor number zero, so $\lambda(Z)=0$.
\end{proof}

For each $\sigma$-amulet $A^{\sigma}$ in $\cA^{\sigma,0} \cup \cA^{\sigma,1}$, we provide the list of figures of the $\sigma$-decorated graph $G_{j,k}^{\sigma} = S - A^{\sigma}$.
The symbols used in the figure are as follows:
\begin{note}\label{fig:involution}
The symbols used in the figure are as follows:
\begin{itemize}
    \item 
    A number written near a line represents the negative self-intersection number of the corresponding curve.
    Any line without a number corresponds to a $(-2)$-curve.
    \item 
    Solid lines represent irreducible components of $D_{j,k}$, while dashed lines represent irreducible components of fibers that are not contained in $D_{j,k}$.
    \item 
    Red lines represent irreducible components of fibers that are fixed by the involution $\sigma$.
    \item 
    Solid red dots represent intersection points with horizontal curves that are fixed by the involution $\sigma$.
    \item 
    Hollow red dots represent isolated fixed points of the involution $\sigma$.
    \item 
    Arrows represent pairs of components exchanged by the involution.
\end{itemize}
\end{note}

\newpage

\begin{table}[H]
  \begin{tabular}{|c|c|}
    \hline $[\mathrm{CT}]^{\sigma}_{\beta}$ ($\beta$:even)  & $[\mathrm{CSR}]^\sigma$  \\ \hline
    \begin{minipage}{70mm}
    \begin{center}
    
\begin{tikzpicture}[line cap=round,line join=round,>=triangle 45,x=0.5cm,y=1cm]
\clip(-3,0) rectangle (7,1.35);
\draw [line width=1pt,red] (-3,1)-- (-1,1);
\draw [line width=1pt] (-2.476872806294636,1.2686349061812408)-- (0,0);
\draw [line width=1pt,red] (-1,0)-- (1,1);
\draw [line width=1pt] (0,1)-- (2,0);
\draw [line width=1pt,red] (1,0)-- (3,1);
\draw [line width=1pt] (4,1)-- (6,0);
\draw [line width=1pt,red] (5,0)-- (7,1);
\begin{scriptsize}
\draw [color=black] (-1.4833503667620005,1.25) node {$\chi$};
\draw [color=black] (3.499147962659239,0.4912663508701683) node {$\cdots$};
\end{scriptsize}
\end{tikzpicture}
    \end{center}
    \end{minipage} &
    \begin{minipage}{70mm}
    \begin{center}
\begin{tikzpicture}[line cap=round,line join=round,>=triangle 45,x=1cm,y=0.75cm]
\clip(-2,-2) rectangle (2,2);
\draw [line width=1pt] (-1,0.8)-- (-1,-0.8);
\draw [line width=1pt] (-1.2,0.5)-- (0,1.5);
\draw [line width=1pt] (-1.2,-0.5)-- (0,-1.5);
\draw[<->] [shift={(1,0)},line width=1pt]  plot[domain=-1.5707963267948966:1.5707963267948966,variable=\t]({1*1*cos(\t r)+0*1*sin(\t r)},{0*1*cos(\t r)+1*1*sin(\t r)});
\draw [line width=1pt,red] (-2,0)-- (-0.5,0);
\begin{scriptsize}
\draw [black] (-1.75,0.25) node {$\chi$};
\draw [fill=red] (-1,-0.25) circle (2.5pt);
\end{scriptsize}
\end{tikzpicture}
    \end{center}
    \end{minipage} \\ \hline
     $[\mathrm{I}_{n}\mathrm{E}]_{\beta}^{\sigma}$ ($\beta$:odd, $n$:even)  & $[\mathrm{I}_{n}\mathrm{E}]_{\beta}^{\sigma}$ ($\beta$:even, $n$:even) \\ \hline
    \begin{minipage}{70mm}
    \begin{center}
\begin{tikzpicture}[line cap=round,line join=round,>=triangle 45,x=0.7cm,y=1cm]
\clip(-8.258182215360351,-1.5) rectangle (2.5,1.5);
\draw [line width=1pt,red] (-8,0)-- (-7,0);
\draw [line width=1pt] (-1.5,1)-- (-1.5,-1);
\draw [line width=1pt] (1.5,1)-- (1.5,-1);
\draw [line width=1pt,dashed,red] (-0.9869053181920332,0.25018220274295794)-- (-2.254647141938089,-0.5211546900181993);
\draw [line width=1pt] (-2,-0.5)-- (-3.247457003907892,0.25018220274295794);
\draw [line width=1pt,red] (-2.995436038946327,0.2578192016811872)-- (-4.255540863754153,-0.49824369320351136);
\draw [line width=1pt] (-4.751945794739054,-0.49824369320351136)-- (-6.004413620608652,0.2425452038047287);
\draw [line width=1pt,red] (-5.760029654585316,0.2425452038047287)-- (-7,-0.5);
\draw [line width=1pt] (-6.760495152946489,-0.4919951649869175)-- (-7.9951475228779,0.2618550801154375);
\draw [line width=1pt] (-0.5,1)-- (-1.697233270487384,0.7085261899103371);
\draw [line width=1pt] (-0.5,-1)-- (-1.697233270487384,-0.7071689554196757);
\draw [line width=1pt] (0.5,1)-- (1.6950928324731982,0.6907187038055571);
\draw [line width=1pt] (0.5,-1)-- (1.6861890894208083,-0.7071689554196757);
\draw[<->] [line width=1pt] (0,-0.75)-- (0,0.75);
\begin{scriptsize}
\draw [black] (-2.25,0.5) node {$3+\beta$};
\draw [black] (-1,-0.125) node {$1$};
\draw [black] (-7.25,0.25) node {$\chi$};
\draw [fill=red] (-1.5,-0.5) circle (2pt);
\draw [fill=red] (1.5,0.5) circle (2pt);
\draw [fill=red] (1.5,-0.5) circle (2pt);
\draw [color=black] (-4.5,0) node {$\cdots$};
\draw [color=black] (0,1) node {$\cdots$};
\draw [color=black] (0,-1) node {$\cdots$};
\end{scriptsize}
\end{tikzpicture}
    \end{center}
    \end{minipage} &
    \begin{minipage}{70mm}
    \begin{center}
\begin{tikzpicture}[line cap=round,line join=round,>=triangle 45,x=0.7cm,y=1cm]
\clip(-8.258182215360351,-1.5) rectangle (2.5,1.5);
\draw [line width=1pt,red] (-8,0)-- (-7,0);
\draw [line width=1pt] (-1.5,1)-- (-1.5,-1);
\draw [line width=1pt] (1.5,1)-- (1.5,-1);
\draw [line width=1pt,dashed] (-0.9869053181920332,0.25018220274295794)-- (-2.254647141938089,-0.5211546900181993);
\draw [line width=1pt,red] (-2,-0.5)-- (-3.247457003907892,0.25018220274295794);
\draw [line width=1pt] (-2.995436038946327,0.2578192016811872)-- (-4.255540863754153,-0.49824369320351136);
\draw [line width=1pt] (-4.751945794739054,-0.49824369320351136)-- (-6.004413620608652,0.2425452038047287);
\draw [line width=1pt,red] (-5.760029654585316,0.2425452038047287)-- (-7,-0.5);
\draw [line width=1pt] (-6.760495152946489,-0.4919951649869175)-- (-7.9951475228779,0.2618550801154375);
\draw [line width=1pt] (-0.5,1)-- (-1.697233270487384,0.7085261899103371);
\draw [line width=1pt] (-0.5,-1)-- (-1.697233270487384,-0.7071689554196757);
\draw [line width=1pt] (0.5,1)-- (1.6950928324731982,0.6907187038055571);
\draw [line width=1pt] (0.5,-1)-- (1.6861890894208083,-0.7071689554196757);
\draw[<->] [line width=1pt] (0,-0.75)-- (0,0.75);
\begin{scriptsize}
\draw [black] (-2.25,0.5) node {$3+\beta$};
\draw [black] (-1,-0.125) node {$1$};
\draw [black] (-7.25,0.25) node {$\chi$};
\draw [fill=red] (-1.5,-0.5) circle (2pt);
\draw [fill=red] (1.5,0.5) circle (2pt);
\draw [fill=red] (1.5,-0.5) circle (2pt);
\draw [color=red] (-1.5,-0.05) circle (2pt);
\draw [color=black] (-4.5,0) node {$\cdots$};
\draw [color=black] (0,1) node {$\cdots$};
\draw [color=black] (0,-1) node {$\cdots$};
\end{scriptsize}
\end{tikzpicture}
    \end{center}
    \end{minipage} \\ \hline
     $[\mathrm{IIT2}]^\sigma$ & $[\mathrm{IIIT1}]^\sigma$ \\ \hline
    \begin{minipage}{70mm}
    \begin{center}
\begin{tikzpicture}[line cap=round,line join=round,>=triangle 45,x=0.75cm,y=1cm]
\clip(-5,-0.5) rectangle (7,1.5);
\draw [line width=1pt,red] (-3,1)-- (-1,1);
\draw [line width=1pt] (-2.476872806294636,1.2686349061812408)-- (0,0);
\draw [line width=1pt] (-1,0)-- (1,1);
\draw [line width=1pt,dashed] (-0.5,0)-- (-0.5,1);
\draw [color=red] (-0.5,0.25) circle (2.5pt);
\draw [color=red] (0.5,0.75) circle (2.5pt);
\begin{scriptsize}
\draw[color=black] (-1.5,0.5) node {$5$};
\draw[color=black] (-0.25,0.75) node {$1$};
\draw [color=black] (-1.5,1.25) node {$\chi$};
\end{scriptsize}
\end{tikzpicture}
    \end{center}
    \end{minipage} &
    \begin{minipage}{70mm}
    \begin{center}
\begin{tikzpicture}[line cap=round,line join=round,>=triangle 45,x=0.75cm,y=1cm]
\clip(-3,-0.5) rectangle (1,1.5);
\draw [line width=1pt,red] (-3,1)-- (-1,1);
\draw [line width=1pt] (-2.5,1.25)-- (0,0);
\draw [line width=1pt] (-1,0)-- (1,1);
\draw [line width=1pt,dashed] (-0.5,0)-- (-0.5,1);
\draw [color=red] (-0.5,0.25) circle (2.5pt);
\draw [fill=red] (0.5,0.75) circle (2.5pt);
\begin{scriptsize}
\draw [color=black] (-1.4833503667620005,1.25) node {$\chi$};
\draw[color=black] (-1.5,0.5) node {$3$};
\draw[color=black] (0.5,0.5) node {$3$};
\end{scriptsize}
\end{tikzpicture}
    \end{center}
    \end{minipage} \\ \hline
     $[\mathrm{IIIT3}]^\sigma$ & $[\mathrm{IVT1}]^\sigma$ \\ \hline
    \begin{minipage}{70mm}
    \begin{center}
\begin{tikzpicture}[line cap=round,line join=round,>=triangle 45,x=0.75cm,y=1cm]
\clip(-3.5,-1) rectangle (6.5,2);
\draw [line width=1pt] (-3,1.5)-- (0,0);
\draw [line width=1pt] (-1,0)-- (1,1);
\draw [line width=1pt,dashed] (0,1)-- (2,0);
\draw [line width=1pt,red] (1,0)-- (4,1.5);
\draw [line width=1pt] (2,1)-- (4,0);
\draw [line width=1pt] (3,1.5)-- (5,0.5);
\draw [line width=1pt,red] (-3,1)-- (-1,1);
\begin{scriptsize}
\draw [color=red] (-0.5,0.25) circle (2.5pt);
\draw [color=red] (0.5,0.75) circle (2.5pt);
\draw [fill=red] (3.507896737475771,0.2460516312621146) circle (2.5pt);
\draw [fill=red] (4.523250458123228,0.738374770938386) circle (2.5pt);
\draw[color=black] (-1.5,1.25) node {$\chi$};
\draw[color=black] (-1.5,0.5) node {$5$};
\draw[color=black] (1.25,0.75) node {$1$};
\draw[color=black] (2.75,1.25) node {$3$};
\draw[color=black] (3,0.25) node {$4$};
\end{scriptsize}
\end{tikzpicture}
    \end{center}
    \end{minipage} &
    \begin{minipage}{70mm}
    \begin{center}
\begin{tikzpicture}[line cap=round,line join=round,>=triangle 45,x=0.75cm,y=1cm]
\clip(-3.5,-1.5) rectangle (4,2);
\draw [line width=1pt] (-3,1.5)-- (0,0);
\draw [line width=1pt] (-1,0)-- (2,1.5);
\draw [line width=1pt,dashed] (0,1)-- (2,0);
\draw [line width=1pt,dashed] (1,1.5)-- (3,0.5);
\draw [line width=1pt] (1,0)-- (1.7165923059411226,0.3571001878781675);
\draw [line width=1pt] (2.2817302273700872,0.6320321496544216)-- (3,1);
\draw [line width=1pt,red] (-3,1)-- (-1,1);
\draw[<->] [shift={(3.1599851052664514,0)},line width=1pt]  plot[domain=-2.5115587440438767:0.6300339095459164,variable=\t]({1*0.6805183157207187*cos(\t r)+0*0.6805183157207187*sin(\t r)},{0*0.6805183157207187*cos(\t r)+1*0.6805183157207187*sin(\t r)});
\begin{scriptsize}
\draw [color=red] (-0.5,0.25) circle (2.5pt);
\draw [fill=red] (1,1) circle (2.5pt);
\draw [color=black] (-1.5,1.25) node {$\chi$};
\draw [color=black] (-1.5,0.5) node {$3$};
\draw [color=black] (0.25,0.4) node {$3$};
\draw [color=black] (2,-0.25) node {$1$};
\draw [color=black] (3,0.25) node {$1$};
\draw [color=black] (1,-0.25) node {$4$};
\draw [color=black] (3.25,1) node {$4$};
\end{scriptsize}
\end{tikzpicture}
    \end{center}
    \end{minipage} \\ \hline
     $[\mathrm{IVR}]^{\sigma}_0$ & $[\mathrm{I}_n^* \mathrm{R}]_\beta^{\sigma}$ ($n$:even)  \\ \hline
    \begin{minipage}{70mm}
    \begin{center}
\begin{tikzpicture}[line cap=round,line join=round,>=triangle 45,x=0.6cm,y=1cm]
\clip(-3.5,-1) rectangle (7.5,2.5);
\draw [line width=1pt] (-3,1.5)-- (0,0);
\draw [line width=1pt] (-1,0)-- (1,1);
\draw [line width=1pt,dashed] (0,1)-- (2,0);
\draw [line width=1pt,red] (1,0)-- (3,1);
\draw [line width=1pt] (2,1)-- (5,-0.5);
\draw [line width=1pt] (3,0)-- (5,1);
\draw [line width=1pt] (4,-0.5)-- (6,0.5);
\draw [line width=1pt,red] (-3,1)-- (-1,1);
\draw[<->] [shift={(6.178508935601564,1.1923969467469953)},line width=1pt]  plot[domain=-0.5468367611498977:2.5947558924398955,variable=\t]({1*0.7398515167762957*cos(\t r)+0*0.7398515167762957*sin(\t r)},{0*0.7398515167762957*cos(\t r)+1*0.7398515167762957*sin(\t r)});
\begin{scriptsize}
\draw [color=red] (-0.5,0.25) circle (2.5pt);
\draw [color=red] (0.5,0.75) circle (2.5pt);
\draw [fill=red] (3,0.5) circle (2.5pt);
\draw[color=black] (-1.5,1.25) node {$\chi$};
\draw[color=black] (-1.5,0.4) node {$5$};
\draw[color=black] (1,0.75) node {$1$};
\draw[color=black] (2,0.75) node {$3$};
\draw[color=black] (6.25,0.5) node {$3$};
\draw[color=black] (5.25,1) node {$3$};
\end{scriptsize}
\end{tikzpicture}
    \end{center}
    \end{minipage} &
    \begin{minipage}{70mm}
    \begin{center}
\begin{tikzpicture}[line cap=round,line join=round,>=triangle 45,x=0.425cm,y=1cm]
\clip(-11,-2) rectangle (4.485006294034534,2);
\draw [line width=1pt, red] (-1,1)-- (-1,-1);
\draw [line width=1pt] (-0.8,0.4)-- (-2.5,1);
\draw [line width=1pt] (-0.8,-0.6)-- (-2,-1);
\draw [line width=1pt] (-1.24,0)-- (0,0.5);
\draw [line width=1pt,red] (3,1)-- (3,-1);
\draw [line width=1pt] (2.8,0.6)-- (4,1);
\draw [line width=1pt] (2.8,-0.6)-- (4,-1);
\draw [line width=1pt,red] (-0.5,0.5)-- (0.7,0);
\draw [line width=1pt,red] (1.3,0)-- (2.5,0.5);
\draw [line width=1pt] (2,0.5)-- (3.5,0);
\draw [line width=1pt] (-3,1)-- (-5,0);
\draw [line width=1pt,dashed,red] (-4,0.85)-- (-1.5,0.85);
\draw [line width=1pt] (-4,0)-- (-6,1);
\draw [line width=1pt] (-7,1)-- (-9,0);
\draw [line width=1pt] (-8,0)-- (-11,1.5);
\draw [line width=1pt,red] (-11,1)-- (-9,1);
\begin{scriptsize}
\draw [black] (-10.75,0.75) node {$\chi$};
\draw [color=red] (-8.5,0.25) circle (2.5pt);
\draw [color=red] (-4.5,0.25) circle (2.5pt);
\draw [black] (-6.5,1) node {$\cdots$};
\draw [black] (1.075,-0.125) node {$\cdots$};
\draw [black] (-9.5,0.5) node {$3$};
\draw [black] (-3.75,0.35) node {$3$};
\draw [black] (-2.5,0.5) node {$1$};
\draw [black] (-1.85,-0.5) node {$3+\beta$};
\draw [fill=red] (-1.613562300367173,-0.8711874334557244) circle (2.5pt);
\draw [fill=red] (3.5,0.825) circle (2.5pt);
\draw [fill=red] (3.5,-0.825) circle (2.5pt);
\end{scriptsize}
\end{tikzpicture}
    \end{center}
    \end{minipage} \\ \hline
%
\end{tabular}
\end{table}

\begin{table}[H]
  \begin{tabular}{|c|c|}  
    \hline $[\mathrm{I}_{n}^{*}\mathrm{SR1}]^{\sigma}$ ($n$:odd) & $[\mathrm{I}_{n}^{*}\mathrm{SR1}]^{\sigma}$ ($n$:even) \\ \hline
    \begin{minipage}{70mm}
    \begin{center}
\begin{tikzpicture}[line cap=round,line join=round,>=triangle 45,x=0.75cm,y=1cm]
\clip(-3.045611532256289,-2) rectangle (5.6617657470576574,2);
\draw [line width=1pt,red] (-1,1)-- (-1,-1);
\draw [line width=1pt] (-0.796405425080071,0.5940678833808372)-- (-2.514305105837884,1.181896073206929);
\draw [line width=1pt] (-0.8037318293379739,-0.5928096063994598)-- (-2,-1);
\draw [line width=1pt, red] (-3,1)-- (-1.5,1);
\draw [line width=1pt] (-1.2408739500595178,0)-- (0,0.5);
\draw [line width=1pt,dashed] (3,1)-- (3,-1);
\draw [line width=1pt,dashed] (2.7935326612923945,0.6013942876387403)-- (4,1);
\draw [line width=1pt,dashed] (2.800859065550297,-0.6001360106573629)-- (4,-1);
\draw [line width=1pt,red,dashed] (-0.5,0.5)-- (0.705507447790042,0.0006291384906886806);
\draw [line width=1pt,dashed] (1.3062725969380873,0)-- (2.5,0.5);
\draw [line width=1pt,red,dashed] (2,0.5)-- (3.5,0);
\draw[<->] [shift={(4.5,0)},line width=1pt]  plot[domain=-1.5707963267948966:1.5707963267948966,variable=\t]({1*1*cos(\t r)+0*1*sin(\t r)},{0*1*cos(\t r)+1*1*sin(\t r)});
\begin{scriptsize}
\draw [color=black] (-1.4833503667620005,1.25) node {$\chi$};
\draw [fill=red] (3,-0.25) circle (2.5pt);
\draw [fill=red] (-1.613562300367173,-0.8711874334557244) circle (2.5pt);
\end{scriptsize}
\end{tikzpicture}
    \end{center}
    \end{minipage} &
    \begin{minipage}{70mm}
    \begin{center}
\begin{tikzpicture}[line cap=round,line join=round,>=triangle 45,x=0.75cm,y=0.85cm]
\clip(-3.045611532256289,-2) rectangle (5,2);
\draw [line width=1pt,red] (-1,1)-- (-1,-1);
\draw [line width=1pt] (-0.796405425080071,0.5940678833808372)-- (-2.514305105837884,1.181896073206929);
\draw [line width=1pt] (-0.8037318293379739,-0.5928096063994598)-- (-2,-1);
\draw [line width=1pt, red] (-3,1)-- (-1.5,1);
\draw [line width=1pt] (-1.2408739500595178,0)-- (0,0.5);
\draw [line width=1pt,dashed,red] (3,1)-- (3,-1); 
\draw [line width=1pt,dashed] (2.7935326612923945,0.6013942876387403)-- (4,1);
\draw [line width=1pt,dashed] (2.800859065550297,-0.6001360106573629)-- (4,-1);
\draw [line width=1pt,red,dashed] (-0.5,0.5)-- (0.705507447790042,0.0006291384906886806);
\draw [line width=1pt,red,dashed] (1.3062725969380873,0)-- (2.5,0.5);
\draw [line width=1pt,dashed] (2,0.5)-- (3.5,0);
\begin{scriptsize}
\draw [color=black] (-1.4833503667620005,1.25) node {$\chi$};
\draw [fill=red] (3.504089322370063,0.8361556731937827) circle (2.5pt);
\draw [fill=red] (3.504089322370063,0.8361556731937827) circle (2.5pt);
\draw [fill=red] (3.5217486324743974,-0.8405229158363116) circle (2.5pt);
\draw [fill=red] (-1.613562300367173,-0.8711874334557244) circle (2.5pt);
\end{scriptsize}
\end{tikzpicture}
    \end{center}
    \end{minipage} \\ \hline
     $[\mathrm{I}_{n}^{*}\mathrm{SR2}]^{\sigma}$($n$:odd) &  $[\mathrm{II}^* \mathrm{SR}]^\sigma$ \\ \hline
    \begin{minipage}{70mm}
    \begin{center}
\begin{tikzpicture}[line cap=round,line join=round,>=triangle 45,x=0.75cm,y=0.85cm]
\clip(-3.045611532256289,-2) rectangle (5.6617657470576574,2);
\draw [line width=1pt,red] (-1,1)-- (-1,-1);
\draw [line width=1pt] (-0.796405425080071,0.5940678833808372)-- (-2.514305105837884,1.181896073206929);
\draw [line width=1pt,dashed] (-0.8037318293379739,-0.5928096063994598)-- (-2,-1);
\draw [line width=1pt, red] (-3,1)-- (-1.5,1);
\draw [line width=1pt] (-1.2408739500595178,0)-- (0,0.5);
\draw [line width=1pt] (3,1)-- (3,-1);
\draw [line width=1pt] (2.7935326612923945,0.6013942876387403)-- (4,1);
\draw [line width=1pt] (2.800859065550297,-0.6001360106573629)-- (4,-1);
\draw [line width=1pt,red] (-0.5,0.5)-- (0.705507447790042,0.0006291384906886806);
\draw [line width=1pt] (1.3062725969380873,0)-- (2.5,0.5);
\draw [line width=1pt,red] (2,0.5)-- (3.5,0);
\draw[<->] [shift={(4.5,0)},line width=1pt]  plot[domain=-1.5707963267948966:1.5707963267948966,variable=\t]({1*1*cos(\t r)+0*1*sin(\t r)},{0*1*cos(\t r)+1*1*sin(\t r)});
\begin{scriptsize}
\draw [color=black] (-1.4833503667620005,1.25) node {$\chi$};
\draw [fill=red] (3,-0.25) circle (2.5pt);
\draw [fill=red] (-1.613562300367173,-0.8711874334557244) circle (2.5pt);
\end{scriptsize}
\end{tikzpicture}
    \end{center}
    \end{minipage} &
    \begin{minipage}{70mm}
    \begin{center}
\begin{tikzpicture}[line cap=round,line join=round,>=triangle 45,x=0.5cm,y=0.85cm]
\clip(-4.359562276188641,0) rectangle (13.00976638154315,1.5);
\draw [line width=1pt] (-2.9748745370407925,1.48726981236248)-- (0,0);
\draw [line width=1pt, red] (-1,0)-- (1,1);
\draw [line width=1pt] (0,1)-- (2,0);
\draw [line width=1pt ,red] (1,0)-- (3,1);
\draw [line width=1pt] (2,1)-- (4,0);
\draw [line width=1pt, red] (3,0)-- (6.0013030000667555,1.4751234286857444);
\draw [line width=1pt] (4,1)-- (6,0);
\draw [line width=1pt] (4.871689318130352,1.4751234286857444)-- (6.815110706408036,0.5398518855771103);
\draw [line width=1pt, dashed, red] (5.98915661639002,0.5762910366073168)-- (7.993309923051381,1.6451728001600416);
\draw [line width=1pt, red] (-3,1)-- (-1,1);
\begin{scriptsize}
\draw [color=black] (-1.4833503667620005,1.3023622180902614) node {$\chi$};
\draw [fill=red] (5.4450029028739655,0.27749854856301703) circle (2.5pt);
\end{scriptsize}
\end{tikzpicture}
    \end{center}
    \end{minipage} \\ \hline
     $[\mathrm{III}^* \mathrm{SR}]^\sigma$ & $[\mathrm{IV}^*\mathrm{SR}]^\sigma$ \\ \hline
    \begin{minipage}{70mm}
    \begin{center}
\begin{tikzpicture}[line cap=round,line join=round,>=triangle 45,x=0.75cm,y=0.85cm]
\clip(-2,-1.5) rectangle (6.899534288088637,1.75);
\draw [line width=1pt] (-1,1)-- (0.5,-0.5);
\draw [line width=1pt ,red] (0,-0.5)-- (1.5,1);
\draw [line width=1pt] (1,1)-- (2.5,-0.5);
\draw [line width=1pt, red] (1.9949545971809697,-0.3144062445991428)-- (3.4968674700510833,-0.3144062445991428);
\draw [line width=1pt] (3,-0.5)-- (4.5,1);
\draw [line width=1pt,red,dashed] (4,1)-- (5.5,-0.5);
\draw [line width=1pt,dashed] (5,-0.5)-- (6.5,1);
\draw [line width=1pt] (2.738643205403547,0.9783762368802792)-- (2.747547123946322,-0.4996742412203729);
\draw [line width=1pt, red] (-1.50852593950008,0.7913939474820039)-- (-0.5023831441665164,0.7913939474820039);
\begin{scriptsize}
\draw [color=black] (-1.2,1.2) node {$\chi$};
\draw [fill=red] (2.739984431439439,0.7557327149223082) circle (2.5pt);
\draw [fill=red] (6.206763099570006,0.7067630995700059) circle (2.5pt);
\end{scriptsize}
\end{tikzpicture}
    \end{center}
    \end{minipage} &
    \begin{minipage}{70mm}
    \begin{center}
\begin{tikzpicture}[line cap=round,line join=round,>=triangle 45,x=0.75cm,y=0.85cm]
\clip(-4,-1) rectangle (5,2);
\draw [line width=1pt] (-3,1.5)-- (0,0);
\draw [line width=1pt,red] (-1,0)-- (1,1);
\draw [line width=1pt] (0,1)-- (3,-0.5);
\draw [line width=1pt] (1,0)-- (3,1);
\draw [line width=1pt] (2,-0.5)-- (4,0.5);
\draw [line width=1pt,dashed] (3.2439920935869733,0.380011184692855)-- (4,0);
\draw [line width=1pt,dashed] (2,1)-- (2.7323131647256136,0.6320321496544213);
\draw[<->] [shift={(3.9944114602591707,1.0889270716418986)},line width=1pt]  plot[domain=-0.6936192170105402:2.447973436579253,variable=\t]({1*0.6429802831071187*cos(\t r)+0*0.6429802831071187*sin(\t r)},{0*0.6429802831071187*cos(\t r)+1*0.6429802831071187*sin(\t r)});
\begin{scriptsize}
\draw [color=black] (-1.4833503667620005,1.3023622180902614) node {$\chi$};
\draw [fill=red] (2,0) circle (2.5pt);
\draw [line width=1pt, red] (-3,1)-- (-1,1);
\end{scriptsize}
\end{tikzpicture}
    \end{center}
    \end{minipage} \\ \hline

    $[\mathrm{I}_{n}\mathrm{T}]^{\sigma}_0$ ($n$:even) & $[\mathrm{I}_{n}\mathrm{T}]^{\sigma}_0$ ($n$:odd) \\ \hline
    \begin{minipage}{70mm}
    \begin{center}
\begin{tikzpicture}[line cap=round,line join=round,>=triangle 45,x=1cm,y=0.75cm]
\clip(-2,-2) rectangle (2.5,2);
\draw [line width=1pt] (-1,0.8)-- (-1,-0.8);
\draw [line width=1pt,dashed] (-1.2,0.5)-- (0,1.5);
\draw [line width=1pt,dashed] (-1.2,-0.5)-- (0,-1.5);
\draw [line width=1pt] (2.0067982655938343,0.6931281411602561)-- (2.0067982655938343,-0.6891686666592434);
\draw [line width=1pt] (2.113716250729044,0.2883671974341043)-- (1.6707703123117472,1.1055260838246372);
\draw [line width=1pt] (2.0984422528525855,-0.29968172080955013)-- (1.701318308064664,-1.0862926114471658);

\draw[<->] [line width=1pt] (0.5,-1)-- (0.5,1);
\draw [line width=1pt,red] (-2,0)-- (-0.5,0);
\draw [line width=1pt] (-0.5,1.4)-- (1,1.4);
\draw [line width=1pt] (-0.5,-1.4)-- (1,-1.4);
\begin{scriptsize}
\draw [black] (-1.75,0.25) node {$\chi$};
\draw [black] (0.5,1.75) node {$3$};
\draw [black] (0.5,-1.75) node {$3$};
\draw [black] (-0.75,0.5) node {$4$};
\draw [fill=red] (-1,-0.25) circle (2pt);
\draw [color=black] (1.4,1.25) node {$\cdots$};
\draw [color=black] (1.4,-1.25) node {$\cdots$};
\draw [fill=red] (2,0.2) circle (2pt);
\draw [fill=red] (2,-0.2) circle (2pt);
\end{scriptsize}
\end{tikzpicture}
    \end{center}
    \end{minipage} &
    \begin{minipage}{70mm}
    \begin{center}
\begin{tikzpicture}[line cap=round,line join=round,>=triangle 45,x=1cm,y=0.75cm]
\clip(-2,-2) rectangle (2.5,2);
\draw [line width=1pt] (-1,0.8)-- (-1,-0.8);
\draw [line width=1pt,dashed] (-1.2,0.5)-- (0,1.5);
\draw [line width=1pt,dashed] (-1.2,-0.5)-- (0,-1.5);
\draw[<->] [line width=1pt] (0.5,-1)-- (0.5,1);
\draw [line width=1pt,red] (-2,0)-- (-0.5,0);
\draw [line width=1pt] (-0.5,1.4)-- (1,1.4);
\draw [line width=1pt] (-0.5,-1.4)-- (1,-1.4);
\draw [line width=1pt] (2.36,-0.34)-- (1.7,1.1);
\draw [line width=1pt] (2.36,0.34)-- (1.7,-1.1);
\begin{scriptsize}
\draw [black] (-1.75,0.25) node {$\chi$};
\draw [black] (0.5,1.75) node {$3$};
\draw [black] (0.5,-1.75) node {$3$};
\draw [black] (-0.75,0.5) node {$4$};
\draw [fill=red] (-1,-0.25) circle (2.5pt);
\draw [fill=red] (2.2,0) circle (2pt);
\draw [color=black] (1.4,1.3) node {$\cdots$};
\draw [color=black] (1.4,-1.3) node {$\cdots$};
\end{scriptsize}
\end{tikzpicture}
    \end{center}
    \end{minipage} \\ \hline
     $[\mathrm{I}_{n}\mathrm{SR}]^{\sigma}$($n$:even) & $[\mathrm{IIT1}]^{\sigma}$ \\ \hline
    \begin{minipage}{70mm}
    \begin{center}
\begin{tikzpicture}[line cap=round,line join=round,>=triangle 45,x=1cm,y=0.75cm]
\clip(-3.239820004969663,-2.25) rectangle (2.5,2.396178904385385);
\draw [line width=1pt] (-1.5,1)-- (-1.5,-1);
\draw [line width=1pt] (2.0067982655938343,0.6931281411602561)-- (2.0067982655938343,-0.6891686666592434);
\draw [line width=1pt] (2.113716250729044,0.2883671974341043)-- (1.6707703123117472,1.1055260838246372);
\draw [line width=1pt] (2.0984422528525855,-0.29968172080955013)-- (1.701318308064664,-1.0862926114471658);
\draw [line width=1pt] (-1.5978652332503736,0.7007651400984855)-- (-1,1.5);
\draw [line width=1pt] (-1.613139231126832,-0.7120796634739311)-- (-1,-1.5);
\draw [line width=1pt] (-0.5057743850835902,1.7928559882652724)-- (0.9987144057475729,1.785218989327043);
\draw [line width=1pt] (0.49467247582444207,1.9150479712769408)-- (1.1896393792033042,1.6248420116242281);
\draw [line width=1pt] (-0.5057743850835902,-1.8194445095171765)-- (0.9910774068093435,-1.8118075105789473);
\draw [line width=1pt] (0.5023094747626713,-1.8958144988994694)-- (1.1972763781415334,-1.6132455381849862);
\draw [line width=1pt,red] (-2.5,0)-- (-1,0);
\draw [shift={(-0.7807063468598433,1.7470339946358966)},line width=1pt, dashed]  plot[domain=-2.5959375990628:0.5456550545269931,variable=\t]({1*0.5003254410306893*cos(\t r)+0*0.5003254410306893*sin(\t r)},{0*0.5003254410306893*cos(\t r)+1*0.5003254410306893*sin(\t r)});
\draw [shift={(-0.773069347921614,-1.7468930196039985)},line width=1pt,dashed]  plot[domain=-0.525850994614057:2.615741658975736,variable=\t]({1*0.49447725185542946*cos(\t r)+0*0.49447725185542946*sin(\t r)},{0*0.49447725185542946*cos(\t r)+1*0.49447725185542946*sin(\t r)});
\draw[<->] [line width=1pt] (0.5,-1.25)-- (0.5,1.25);
\begin{scriptsize}
\draw [fill=red] (-2,0.2) node {$\chi$};
\draw [color=black] (1.4951193367324742,1.4033690424155791) node {$\cdots$};
\draw [color=black] (1.5027563356707037,-1.3994095679145664) node {$\cdots$};
\draw [fill=red] (-1.5,-0.5) circle (2pt);
\draw [fill=red] (2.0067982655938343,0.18908621123712366) circle (2pt);
\draw [fill=red] (2.0067982655938343,-0.20040073461256958) circle (2pt);
\draw [color=black] (0,1.6) node {$3$};
\draw [color=black] (0,-1.6) node {$3$};
\draw [color=black] (-1.3993032608564129,-1.3077655806558148) node {$3$};
\draw [color=black] (-1.3993032608564129,1.2964510572803691) node {$3$};
\end{scriptsize}
\end{tikzpicture}
    \end{center}
    \end{minipage} &
    \begin{minipage}{70mm}
    \begin{center}
\begin{tikzpicture}[line cap=round,line join=round,>=triangle 45,x=0.75cm,y=0.5cm]
\clip(-2,-2) rectangle (2,3);
\draw [line width=1pt,red] (-1,2)-- (1,2);
\draw [line width=1pt] (0,2.5)-- (0,-2);
\draw [shift={(1,0)},line width=1pt,dashed]  plot[domain=1.5707963267948966:4.71238898038469,variable=\t]({1*1*cos(\t r)+0*1*sin(\t r)},{0*1*cos(\t r)+1*1*sin(\t r)});
\begin{scriptsize}
\draw [black] (-0.5,1.75) node {$\chi$};
\draw [black] (-0.25,-1.5) node {$4$};
\draw [black] (0.5,1.25) node {$1$};
\draw [fill=red] (0,0) circle (2.5pt);
\draw [color=red] (0.43,-0.82) circle (2.5pt);
\end{scriptsize}
\end{tikzpicture}
    \end{center}
    \end{minipage} \\ \hline
     $[\mathrm{IIIR}]_1^{\sigma}$ & $[\mathrm{IVT2}]^\sigma$ \\ \hline
    \begin{minipage}{70mm}
    \begin{center}
\begin{tikzpicture}[line cap=round,line join=round,>=triangle 45,x=0.6cm,y=1cm]
\clip(-3.5,-1) rectangle (9,2);
\draw [line width=1pt] (-3,1.5)-- (0,0);
\draw [line width=1pt] (-1,0)-- (1,1);
\draw [line width=1pt,red] (0,1)-- (2,0);
\draw [line width=1pt] (1,0)-- (3,1);
\draw [line width=1pt,red] (-3,1)-- (-1,1);
\draw [line width=1pt,dashed,red] (2,1)-- (4,0);
\draw [line width=1pt] (3,0)-- (5,1);
\draw [line width=1pt,red] (4,1)-- (7,-0.5);
\draw [line width=1pt] (5,0)-- (7,1);
\draw [line width=1pt] (6,-0.5)-- (8,0.5);
\begin{scriptsize}
\draw [color=red] (-0.5,0.25) circle (2.5pt);
\draw [fill=red] (6.602418243400987,0.8012091217004931) circle (2.5pt);
\draw [fill=red] (7.6112256478713425,0.3056128239356714) circle (2.5pt);
\draw[color=black] (-1.5,1.25) node {$\chi$};
\draw[color=black] (-1.5,0.5) node {$5$};
\draw[color=black] (0,0.25) node {$3$};
\draw[color=black] (3,0.75) node {$1$};
\draw[color=black] (5.25,0.6) node {$3$};
\draw[color=black] (8.25,0.5) node {$4$};
\draw[color=black] (4,0.25) node {$4$};
\end{scriptsize}
\end{tikzpicture}
    \end{center}
    \end{minipage} &
    \begin{minipage}{70mm}
    \begin{center}
\begin{tikzpicture}[line cap=round,line join=round,>=triangle 45,x=0.75cm,y=1cm]
\clip(-3.5,-1) rectangle (5,2);
\draw [line width=1pt] (-3,1.5)-- (0,0);
\draw [line width=1pt,red,dashed] (-1,0)-- (1,1);
\draw [line width=1pt] (0,1)-- (3,-0.5);
\draw [line width=1pt] (1,0)-- (3,1);
\draw [line width=1pt] (2,-0.5)-- (4,0.5);
\draw[<->] [shift={(3.9944114602591707,1.0889270716418986)},line width=1pt]  plot[domain=-0.6936192170105402:2.447973436579253,variable=\t]({1*0.6429802831071187*cos(\t r)+0*0.6429802831071187*sin(\t r)},{0*0.6429802831071187*cos(\t r)+1*0.6429802831071187*sin(\t r)});
\begin{scriptsize}
\draw [color=black] (-1.4833503667620005,1.3023622180902614) node {$\chi$};
\draw [fill=red] (2,0) circle (2.5pt);
\draw [line width=1pt, red] (-3,1)-- (-1,1);
\draw[color=black] (-1.5,0.5) node {$4$};
\draw[color=black] (-0.25,0.75) node {$1$};
\draw[color=black] (3.25,1) node {$3$};
\draw[color=black] (4.25,0.5) node {$3$};
\end{scriptsize}
\end{tikzpicture}
    \end{center}
    \end{minipage} \\ \hline
\end{tabular}
\end{table}

We can now state the result of this step.

\begin{prop}\label{prop:sigma-amulet}
Suppose that Assumption \ref{ass:5-2} holds.
If $D_{j,0}$ does not intersect $S$, set $l=1$ and reordering the blow-ups so that $D_{j,1}$ intersects $S$; otherwise, set $l=0$.
Then, exactly one of the following holds:
    \begin{itemize}
        \item[$(1)$] 
        $S+B_{j,l}$ is a chain of smooth rational curves and each connected component of the effective divisor $D_{j,l}-B_{j,l}$ is a T-chain.
        By reordering the blow-ups if necessary, there exists some $k$ and a unique amulet $A\in\cAT^{\sigma,0}\cup\cAT^{\sigma,1}$ such that no blow-ups occur on $D_{j,k}-B_{j,k}$, and $\sigma$ induces the involution on $Y_k$ and the graph $G_{j,k}^\sigma$ coincides with $S-A$. 
 
        \item[$(2)$] 
        $S+B_{j,l}$ is a chain of smooth rational curves, the effective divisor $D_{j,l}-B_{j,l}$ consists of T-chains and exceptional curves of strictly lc rational singularities, and it contains at least one of the latter.
        By reordering the blow-ups if necessary, there exists some $k$ and a unique amulet $A\in\cAR^{\sigma,0}\cup\cAR^{\sigma,1}$ such that no blow-ups occur on $D_{j,k}-B_{j,k}$, and $\sigma$ induces the involution on $Y_k$ and the graph $G_{j,k}^\sigma$ coincides with $S-A$.
        
        \item[$(3)$] 
        $S+B_{j,l}$ has a fork, and $D_{j,l}=B_{j,l}$ or effective divisor $D_{j,l}-B_{j,l}$ consists of T-chains.
        By reordering the blow-ups if necessary, there exists some $k$ and a unique amulet $A\in\cASR^{\sigma,0}\cup\cASR^{\sigma,1}$ such that no blow-ups occur on the effective divisor $D_{j,k}- B_{j,k}$, and $\sigma$ induces the involution on $Y_k$ and the graph $G_{j,k}^\sigma$ coincides with $S-A$.
        
        \item[$(4)$] 
        $D_{j,0}$ contains the exceptional curves of an elliptic singularity.
        By reordering the blow-ups if necessary, there exists some $k$ and a unique amulet $A\in\cAE^{\sigma,0}\cup\cAE^{\sigma,1}$ such that no blow-ups occur on the effective divisor $D_{j,k}- B_{j,k}$, and $\sigma$ induces the involution on $Y_k$ and the graph $G_{j,k}^\sigma$ coincides with $S-A$.
    \end{itemize}
    Moreover, in each case, if $A\in\cA^{\sigma,1}$, then $\rho(W)=\lambda(W)=1$.
\end{prop}

This follows immediately from Proposition \ref{amulet}, Lemma \ref{lem:restrict_involution}, and Lemma \ref{lem:lambda}.

\subsubsection{Step 2: completely separated branch}\label{subsubsec:0,1,0-inv-step2}
In this step, we classify all the cases where $D_j$ is completely separated from $S$ after blow-ups under Assumption \ref{ass:5-2}.
By using Proposition \ref{prop:sigma-amulet} and computing $\lambda(W)$ as in the proof of Lemma~\ref{lem:lambda}, we can refine Proposition \ref{compsepamulet} as follows.

\begin{prop}\label{prop_inv_cs}
    Suppose that Assumption \ref{ass:5-2} holds and $D_j$ is completely separated from $S$ after the blow-ups.
    Then, $G_{j,0}^\sigma=S-[\mathrm{I}_n \mathrm{E}]_0^\sigma$ for some $n\geq0$.
\end{prop}

\subsubsection{Step 3: case that $p$ is a T-singularity}\label{subsubsec:0,1,0-inv-step3}
In this step, we classify possible $C$ and its resolution $\widetilde{C}$ under Assumption \ref{ass:5-2} and the assumption that $p\in X$ is a T-singularity.

\begin{thm}\label{thm:one-section-fixed-T}
    Suppose that Assumption \ref{ass:5-2} holds, the singular point $p \in X$ is a T-singularity. 
    Let $J$ be the number of brnches that are not completely separated from the section $S$.
    By reordering the blow-ups if necessary, there exists some $k$ such that exactly one of the following holds:
    \begin{itemize}
    \item $J=0$, $\chi=4$, and $B_0=[4]$.
    \item $J=1$. 
    $\chi$, $G_{1,0}$, and $B_0$ coincide with one in the following list.
    Here, the number in $B_0$ corresponding to the proper transform of $S$ is underlined.

    \smallskip
    
    \begin{tabular}{|l|l|l|l|} \hline
            &   $\chi$      &   $G_{1,0}^\sigma=S-A^\sigma$       &   $B_0$    \\
        \hline
        $(1.1)$ &   $\chi=3$    &   $A^\sigma = [\mathrm{IIT2}]^{\sigma}$       &   $[\underline{3},5,2]$        \\
        $(1.2)$ &   $\chi=3$    &   $A^\sigma = [\mathrm{IIIT3}]^{\sigma}$      &   $[\underline{3},5,2]$        \\
        $(1.3)$ &   $\chi=3$    &   $A^\sigma = [\mathrm{IVR}]^{\sigma}_{0}$      &   $[\underline{3},5,2]$        \\
        $(1.4)$ &   $\chi=5$ &   $A^\sigma= [\mathrm{CT}]_1^\sigma$ &   $[\underline{5},2]$ \\
        $(1.5)$ &   $\chi=5$ &   $A^\sigma=[\mathrm{I}_n \mathrm{E}]_1^\sigma$ &   $[\underline{5},2]$ \\
        \hline
    \end{tabular}

    \smallskip
    
    \item $J=2$. 
    By reordering the blow-ups if necessary, there exists some $k$ such that $\chi$, $G_{1,k}^\sigma$, $G_{2,k}^\sigma$, and $B_k$ coincide with exactly one in the following list.
    Here, the number in $B_k$ corresponding to the proper transform of $S$ is underlined.

    \smallskip
    
     \begin{tabular}{|l|l|l|l|l|l|} \hline
            &   $\chi$  &    $G_{1,k}^\sigma=S-A_1^\sigma$     & $G_{2,k}^\sigma=S-A_2^\sigma$  &   $B_k$    \\
        \hline
         $(2.1)$   &  $\chi=3$      &  $A_1\in\cA_{[2]}$   &   $A_2=[\mathrm{I}_{1+n}\mathrm{T}]_0^{\sigma}$   &   $[2,\underline{3},4]$     \\
         $(2.2)$   &   $\chi=3$     &   $A_1\in\cA_{[2]}$   &   $A_2=[\mathrm{IIIT1}]^{\sigma}$     &   $[2,\underline{3},3,3]$     \\
         $(2.3)$   &   $\chi=3$     &   $A_1\in\cA_{[2]}$   &   $A_2=[\mathrm{IVT1}]^{\sigma}$     &   $[2,\underline{3},3,3]$     \\
         $(2.4)$   &   $\chi=3$    &   $A_1\in\cA_{[2]}$   &   $A_2=[\mathrm{I}_{n}^{*}\mathrm{R}]^{\sigma}_{1}$    &   $[2,\underline{3},3,3]$     \\
         $(2.5)$   &   $\chi = 4$    &   $A_1\in\cA_{[2]}$   &   $A_2=[\mathrm{IIIT1}]^{\sigma}$  &   $[2,\underline{4},3,3]$     \\
         $(2.6)$   &   $\chi = 4$   &   $A_1\in\cA_{[2]}$   &   $A_2=[\mathrm{IVT1}]^{\sigma}$    &   $[2,\underline{4},3,3]$     \\
        $(2.7)$   &   $\chi\geq 4$    &   $A_1\in\cA_{[2]}$   &   $A_2=[\mathrm{I}_{n}^{*}\mathrm{R}]^{\sigma}_{\chi-3}$    &   $[2,\underline{\chi},3,2^{\chi-4},3]$     \\
        \hline
    \end{tabular}
    
    \smallskip
    
    In the table, the variable $n$ is a non-negative integer. 
    In cases $(2.1)$, $(2.5)$, $(2.6)$, and $(2.7)$, the chain $B_k$ is a T-chain, i.e., $B_k=B_0$.
    In the other cases, the T-train associated to $B_k=[2,\underline{3},3,3]$ is $[2,\underline{3},4]-1-[4]$.
    \end{itemize}

\end{thm}

\begin{rem}\label{rem--correspondence-inv-no_inv-T}
For the reader's convenience, we provide tables indicating the correspondence between the list in Theorem~\ref{thm:one-section-fixed-T} and the list in Theorem~\ref{thm:one-section-T}.
\begin{table}[htbp]
\begin{flushleft}
    \begin{tabular}{|c||c|c|c|} \hline
     Theorem~ \ref{thm:one-section-fixed-T}                       &  $(1.1)\text{--}(1.3)$   & $(1.4),(1.5)$ & $(2.1)\text{--}(2.7)$  \\ \hline 
Theorem~\ref{thm:one-section-T}   & $(1.2)$ & $(1.3)$  &  $(2.12)$    \\ \hline
    \end{tabular}
\end{flushleft}
\end{table}
\end{rem}

\begin{proof}[Proof of Theorem \ref{thm:one-section-fixed-T}]
    If $J=0$, Theorem \ref{thm:one-section-T} indicates that $\chi=4$.
    In this case, the surface $W$ is $\overline{\Sigma}_8$ which has a P-resolution $\Sigma_8\to W$.
    This gives the first case.

    Let $J=1$ or $J=2$.
    From Proposition \ref{prop:sigma-amulet}, there exist some $k$ and $A_j\in\cAT^{\sigma,0}\cup\cAT^{\sigma,1}\cup\cAR^{\sigma,0}\cup\cAR^{\sigma,1}\cup\cAE^{\sigma,0}\cup\cAE^{\sigma,1}$ ($1\leq j\leq J$) such that $\sigma$ induces the involution on $Y_k$ and $G_{j,k}^\sigma$ coincides with $S-A_j^\sigma$ for $1\leq j\leq J$.
    Moreover, we cannot use two amulets in $\cA^{\sigma,1}$, as otherwise $\lambda(W)\geq2$ which contradicts Assumption \ref{ass:5-2} by Lemma \ref{lem:restrict_involution}.
    Based on these observations, the list of Theorem \ref{thm:one-section-T} can be reduced.
    In each case, since we can determine $Y_0$ and $C_0^{(\alpha)}$ locally around $\bigcup_{j=1}^J D_{j,0}$.
    We can exclude all the cases except for the cases $(1.2)$, $(1.3)$, and $(2.12)$ in Theorem~\ref{thm:one-section-T} by computing $\lambda(W)$, which must be equal to $0$ or $1$ by Lemma \ref{lem:restrict_involution}.
    
    For example, we check that (2.8) in Theorem~\ref{thm:one-section-T} does not occur when $n=1$.
    Let us take $k$ such that $B_k=[3,3,5,2^2]$, and let $A_1\in\cA_{[3,3]}$ and $A_2\in\cA_{[2^2]}$ be amulets which give the decorated graphs $G_{1,k}$ and $G_{2,k}$.
    After the blow-ups $\widetilde{X}\to Y_k$, the divisor $B_k$ becomes the disjoint union of $[4]$ and $[4,5,2^2]$.
    We remark that, if $A_1=[\mathrm{I}_2\mathrm{T}]_0$, the involution $\sigma$ cannot be extended to $Y_k$, but it can to $\widetilde{X}$.
    Let $U$ be the complement of $\widetilde{C}\setminus B$ in $\widetilde{X}$.
    Then, $U$ is preserved by the involution $\sigma$, and $\pi\colon U\to \pi(U)$ contracts $[4]$ and $[4,5,2^2]$.
    Let $V\to Z:=\pi(U)/\sigma$ be the minimal resolution.
    Note that $V$ is the relatively minimal model of the minimal resolution of $U/\sigma$ over $Z$.
    It is straightforward to see that the exceptional divisor of $V\to Z$ is the disjoint union of $[2]$ and $[2,9,3]$.
    Since $[2,9,3]$ is not a T-chain, this implies that any P-resolution $Z^{\dagger}\to Z$ of $Z$ satisfies $\rho(Z^{\dagger}/Z)\ge 2$.
    Hence, we have $\lambda(W)\ge \lambda(Z)\ge 2$.
    The claim follows from Lemma~\ref{lem:restrict_involution}.

\end{proof}

\subsubsection{Step 4: case that $p$ is a strictly lc rational singularity}\label{subsubsec:0,1,0-inv-step4}
In this step, we classify possible $C$ and its resolution $\widetilde{C}$ under Assumption \ref{ass:5-2} and the assumption that $p\in X$ is a strictly lc rational singularity.

We follow the same strategy as the one in Section \ref{subsubsec:0,1,0-step4}.
We first establish a $\sigma$-decorated analogue of Lemma \ref{lem:amulet-b}.
Since we assume that $p\in X$ is not a strictly lc singularity of type $(3,3,3)$, we focus on $\cA_{[2]}^{\text{B}}$ and $\cA_{[4]}^{\text{B}}$.

We introduce $\sigma$-amulets for the blow-up processes classified in Lemma~\ref{lem:amulet-b}.
Let $G_{j,k}=S-A$ be an amulet in $\cA_{[2]}^B \cup \cA_{[4]}^{B}$ except for $S-[\mathrm{I}_n \mathrm{T}_1 \mathrm{B2}]$ $(n \ge 1)$.
It can be verified that the involution lifts to the germ of the fiber $F_{j,k}$ by chasing the corresponding blow-up process.
Hence, the $\sigma$-decorated graph $G_{j,k}^\sigma=S-A^{\sigma}$ is well-defined.
We define the sets of $\sigma$-decorated graphs $\cA_{[2]}^{\text{B},\sigma,0}$ and $\cA_{[4]}^{\text{B},\sigma,1}$ as follows: 
\begin{align*}
    \cA_{[2]}^{\text{B},\sigma,0} &= \{  [\mathrm{CT}]_1^\sigma \} \cup \{ [\mathrm{I}_n \mathrm{E}]_1^\sigma : n \geq0\},  \\
    \cA_{[4]}^{\text{B},\sigma,1} &= \{ [\mathrm{IIT1}]^\sigma, [\mathrm{IVT2}]^\sigma, [\mathrm{IIIT1B4}]^{\sigma} \} \cup \{ [\mathrm{I}_n \mathrm{T}]_0^\sigma : n\geq1 \} \cup \{ [\mathrm{I}_n^* \mathrm{R}]_0^\sigma : n\geq0 \}.
\end{align*}

We explicitly describe the additional data induced by the involution on the $\sigma$-decorated graph $G_{j,k}^\sigma = S - [\mathrm{IIIT1B4}]^\sigma$.
The symbols used in the figure are as described in Notation~\ref{fig:involution}.

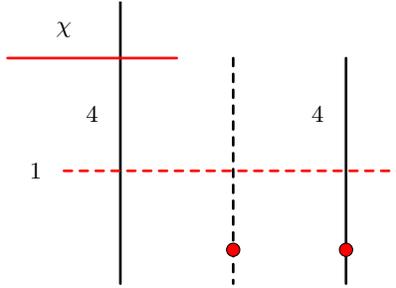
\begin{center}
\begin{tikzpicture}[line cap=round,line join=round,>=triangle 45,x=1.5cm,y=1.5cm]
\clip(-3,-1.5) rectangle (1.5,1.5);
\draw [line width=1pt] (-1.5,1.5)-- (-1.5,-1);
\draw [line width=1pt,dashed] (-0.5,1)-- (-0.5,-1);
\draw [line width=1pt] (0.5,1)-- (0.5,-1);
\draw [line width=1pt,red,dashed] (-2,0)-- (1,0);
\draw [line width=1pt,red] (-2.5,1)-- (-1,1);
\begin{scriptsize}
\draw [black] (-2,1.25) node {$\chi$};
\draw [black] (-1.75,0.5) node {$4$};
\draw [black] (0.25,0.5) node {$4$};
\draw [black] (-2.25,0) node {$1$};
\draw [fill=red] (-0.5,-0.7) circle (2.5pt);
\draw [fill=red] (0.5,-0.7) circle (2.5pt);
\end{scriptsize}
\end{tikzpicture}
\end{center}

The following is a refinement of Lemma \ref{lem:amulet-b}.

\begin{lem}\label{lem:sigma-amulet-b}
    Suppose Assumption \ref{ass:5-2} holds and the connected component of $D_{j,0}$ intersecting $S$ consists of a single $(-b)$-curve for some $b$.
    Then, exactly one of the following holds:
    \begin{enumerate}[label=$(\arabic*)$]
        \item $b=2$ and $G_{j,0}^\sigma=S-A$ for some $A\in\cA_{[2]}^{\text{B},\sigma,0}$.

        \item $b=4$, $K_{W^\dagger}^2=8$, and $G_{j,0}^\sigma=S-A$ for some $A\in\cA_{[4]}^{\text{B},\sigma,1}$.
    \end{enumerate}
\end{lem}

\begin{proof}
Note that it is necessary to satisfy the assumption $\lambda(W) \le 1$ in Assumption~\ref{ass:5-2}.
We can show $\lambda(W) \ge 2$ if $G_{j,0}^\sigma=S-A$ with $A$ not contained in $\cA_{[2]}^{B,\sigma,0} \cup \cA_{[4]}^{B,\sigma,1}$, by an argument similar to that of Lemma~\ref{lem:lambda}.

For example, we show that $G_{j,0}^{\sigma}= S-[\mathrm{IIT2B2}]^\sigma$ is not allowed.
Define $U$ as the open subset 
\[
U=Y_k\setminus(B_k \cup \bigcup_{j'\neq j}D_{j',k}).
\]
Since the morphism $Y_{0} \to Y_{k}$ is an isomorphism over $U$, we can regard $U$ as an open subset of $Y_0=\widetilde{X}$.
Then $U$ is preserved by the involution $\sigma$,
and $\pi: U \to \pi (U)$ contracts $[3,2,6,2]$.
Let $V \to Z:=\pi(U)/\sigma$ be the minimal resolution.
It is straight forward to see that the exceptional divisor of $V \to Z$ is a chain $[2^4,4,2^3]$.
Therefore, we get $\lambda(W) \ge 2$.
It does not hold the assumption~\ref{ass:5-2}.
\end{proof}

Finally, we obtain the classification by the similar argument as in Section \ref{subsubsec:0,1,0-step4}.

\begin{thm}\label{thm:one-section-fixed-lc}
Suppose that Assumption \ref{ass:5-2} holds, the singular point $p \in X$ be a strictly lc rational singularity, and there exists no branches that are completely separated from the section $S$.
Then, exactly one of the following holds:
\begin{enumerate}[label=$\arabic*.$]
\item $(S-A_{[2]},-A_{[2]},-A_{[2]},-A_{[2]})$ 
\begin{itemize}
    \item $J=4$, and $G_{j,0}^\sigma=S-A_{j}$ where $A_{j}\in\cA_{[2]}^{\text{B},\sigma,0}$ for $1\leq j\leq4$. 
    \item $3 \le \chi \le 6$
    \item The dual graph of $B_0$ is
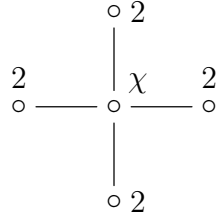
\[
\xygraph{
    \circ ([]!{-(0,-.3)} {2}) - [r]
    \circ ([]!{-(-0.25,-.25)} {\chi}) (
         - [u] \circ ([]!{-(-0.25,-.0)} {2}),
        - [d] \circ ([]!{-(-0.25,-.0)} {2}), 
         - [r] \circ ([]!{-(0,-.3)} {2})
)}
\]
\end{itemize}

\item $(S-A_{[2]},-A_{[2]},-A_{\SR})$ 
\begin{itemize}
    \item $J=3$, $G_{1,0}^\sigma=S-A_{1}$, $G_{2,0}^\sigma=S-A_{2}$, and $G_{3,0}^\sigma=S-A$ where $A_{1}, A_{2}\in\cA_{[2]}^{\text{B},\sigma,0}$ and $A \in \cASR^{\sigma,0}$. 
    \item $3 \le \chi\le 6 $. 
    \item The dual graph of $B_0$ is
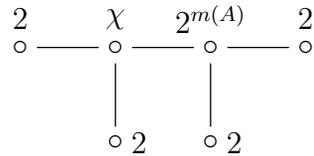
\[
\xygraph{
    \circ ([]!{-(0,-.3)} {2}) - [r]
    \circ ([]!{-(0,-.3)} {\chi}) (
        - [d] \circ ([]!{-(-0.25,-.0)} {2}),
         - [r] \circ ([]!{-(0,-.3)} {2^{m(A)}}) (
         - [d] \circ ([]!{-(-0.25,-.0)} {2}),
        - [r] \circ ([]!{-(0,-.3)} {2})
)}
\]
\end{itemize}

\item $(A_{\SR}-S-A_{\SR})$ 
\begin{itemize}
    \item $J=2$, $G_{1,0}^\sigma=S-A_{1}$ and $G_{2,0}^\sigma=S-A_{2}$ where $A_{1},A_{2}\in \cASR^{\sigma,0}$. 
    \item $2 \le \chi \le 6$.
    \item The dual graph of $B_0$ is
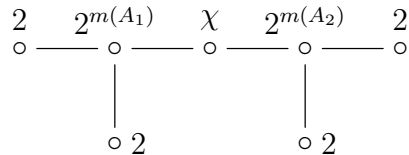
\[
\xygraph{
    \circ ([]!{-(0,-.3)} {2}) - [r]
    \circ ([]!{-(0,-.3)} {2^{m(A_{1})}}) (
        - [d] \circ ([]!{-(-0.25,-.0)} {2}),
        - [r] \circ ([]!{-(-0,-.3)} {\chi})
         - [r] \circ ([]!{-(0,-.3)} {2^{m(A_{2})}}) (
         - [d] \circ ([]!{-(-0.25,-.0)} {2}),
        - [r] \circ ([]!{-(0,-.3)} {2})
)}
\]
\end{itemize}

\item $(S-[\mathrm{I}_n\mathrm{SR}]^\sigma)$
\begin{itemize}
    \item $J=1$, and $G_{1,0}^\sigma=S-[\mathrm{I}_n\mathrm{SR}]$ for some $n\geq3$.
    \item $\chi=3$.
    \item The dual graph of $B_0$ is
\[
\xygraph{
    \circ ([]!{-(0,-.3)} {\chi=3}) 
    - [r]  \circ ([]!{-(0,-.3)} {2})(
        - [r] \circ ([]!{-(0,-.3)} {3}),
        - [d] \circ ([]!{-(-0.25,-.0)} {3}),
}
\]
\end{itemize}

\end{enumerate}
\end{thm}

\subsection{Case of $(g,n_X,l_{E_h})=(0,2,0)$}\label{subsec:0,2,0}
In this subsection, we classify $C$ and its resolution $\widetilde{C}$ under the assumption that $(g,n_X,l_{E_h})=(0,2,0)$, $\delta'_f=0$, $r=0$, and $l_{R_0}=0$.

Since $n_{X}=2$, there are two possible cases: either $C$ has a bisection, or $C$ has two sections.
In the first case, it is easy to see that $r\geq1$, which contradicts the assumption.
Hence $C$ has two sections, which are denoted by $S_1$ and $S_2$, and $C_\hor$ in \eqref{eqn:decomp} equals $S_1+S_2$.
We remark that $D_j$ may not be connected in this situation.

We introduce the following decorated graphs:
\begin{align*}
    &S_1-[\mathrm{I}_n\mathrm{E1}]_{\beta_1,\beta_2,\gamma}-S_2 := 
    \xygraph{
    \diamond ([]!{-(0,-.3)} {+1})([]!{-(.0,0.25)} {S_1}) - [r]
    \circ ([]!{-(0,-.3)} {2^{\beta_1}}) - [r]
    \bullet ([]!{-(0,-.3)} {}) - [r]
    \circ ([]!{-(0.25,-.3)} {\beta_1 +3}) ( 
        - [ur] \circ ([]!{-(0,-.3)} {2^{n-\gamma-1}}) 
        - [rd] \circ ([]!{-(-0.25,-.3)} {\beta_2 +3}) 
        ([]!{-(.5,0)} {}),
        - [rd] \circ ([]!{-(.0,0.25)} {2^{\gamma-1}})
        - [ru] \circ ([]!{-(-0.3,-.0)} {}) 
        - [r] \bullet ([]!{-(.3,0)} {})
         - [r] \circ ([]!{-(0,-.3)} {2^{\beta_2}})
          - [r] \diamond ([]!{-(0,-.3)} {+1})([]!{-(.0,0.25)} {S_2}),    
)}  \\
    &S_1-[\mathrm{I}_n\mathrm{E2}]-S_2 :=  
\xygraph{
    \diamond ([]!{(-0.3,-.0)} {+1})  
    - [r] \bullet ([]!{(-0.3,-.0)} {}) ( 
        - [u] \circ ([]!{-(.3,-.3)} {3}) (
        - [u] \circ ([]!{-(.3,0)}  {2}, -[rd] \circ ),
        - [r] \circ ([]!{-(-.3,-.3)} {2^{n-2}})) ,  
        - [r] \diamond ([]!{-(-0.3,-.0)} {+1}),
)}
  \\
    &S_1- [\mathrm{CT}]_{\beta_1,\beta_2}-S_2 := 
 \xygraph{
    \diamond ([]!{-(0,-.3)} {+0}) - [r]
    \circ ([]!{-(0,-.3)} {2^{\beta_1}})
          [r] \circ ([]!{-(0,-.3)} {2^{\beta_2}}) 
        - [r] \diamond ([]!{-(0,-.3)} {+0})
)}
\end{align*}

For $S_1-[\mathrm{I}_n\mathrm{E1}]_{\beta_1,\beta_2,\gamma}-S_2$, if $n=0$, the cycle $[\beta_1,2^{n-\gamma-1},\beta_2,2^{\gamma-1}]^{\circ}$ regard as a smooth elliptic curve with self-intersection number $-(\beta_1+\beta_2+2)$.
For $S_1-[\mathrm{I}_n\mathrm{E2}]-S_2$,
if $n=0$, the cycle $[3,2^{n-1}]^{\circ}$ regard as a amooth elliptic curve with self intersection number $-1$.

Let 
\begin{align*}
    J'':= &\sharp \{ 1\leq j\leq J : G_{j,0} = S_1-[\mathrm{I}_m \mathrm{E1}]_{0,0,\gamma}-S_2\ \text{for some $m,\gamma$} \}    \\
    + &\sharp \{ 1\leq j\leq J : G_{j,0} = S_1-[\mathrm{I}_m\mathrm{E2}]-S_2\ \text{for some $m$} \},
\end{align*}
and let $J':=J-J''$.
By reordering the branches, we may assume that $G_{j,0} = S_1-[\mathrm{I}_{n_j}\mathrm{E1}]_{0,0,\gamma}-S_2$ or $G_{j,0} = S_1-[\mathrm{I}_{n_j}\mathrm{E2}]-S_2$ for $J'+1\leq j\leq J$.
The following theorem gives the classification.

\begin{thm}\label{thm:two-sections}
    Suppose that $(g,n_X,l_{E_{h}})=(0,2,0)$, $\delta'_f=0$, $r=0$, and $l_{R_0}=0$.
    Reordering the sections $S_i$ if necessary, exactly one of the following holds:
    \begin{itemize}
        \item[$(1)$] 
        $J'=2$ and $\chi\geq5$.
        There exist $A_1, A_2$ such that
        \begin{align*}
        A_1&\in\{ [\mathrm{I}_n\mathrm{E1}]_{\chi-4,0,\gamma_1} : n\geq\chi-11\} \cup \{  [\mathrm{CT}]_{\chi-4,0} \}, \\
        A_2&\in\{ [\mathrm{I}_n\mathrm{E1}]_{0,\chi-4,\gamma_2} : n\geq\chi-11\} \cup \{  [\mathrm{CT}]_{0,\chi-4} \},
        \end{align*}
        and $G_{j,0}=S_1-A_j-S_2$ for $j=1,2$.
        
        \item[$(2)$]
        $J'=1$ and $\chi\geq5$.
        The two sections are disjoint.
        There exists a bi-amulet $A_1$ such that
        \begin{align*}
        A_1\in &\{ [\mathrm{I}_n\mathrm{E1}]_{\chi-4,\chi-4,\gamma} : \textrm{ $n\geq2\chi-14$ or $n=2\chi-15$, $\gamma\in\{\chi-10,\chi-9,\chi-6,\chi-5\}$} \} \\
        &\cup \{  [\mathrm{CT}]_{\chi-4,\chi-4} \},
         \end{align*}
        and $G_{1,0}=S_1-A_1-S_2$.
        
        \item[$(3)$]
        $J'=0$ and $\chi=4$.
 \end{itemize}

\end{thm}

\begin{proof}
Following Notation~\ref{notation_7}, we write $C=S_1+S_2+\sum_{j \in J}D_j$. 
First, we study all possible blow-up processes on $D_j$.
The two sections $S_1$ and $S_2$ either do not meet in the fiber $F_j$, or meet in $F_j$.
We deal with each case separately.

Suppose that $S_1$ and $S_2$ do not intersect in $F_j$.
There are four possibilities for $D_j$: it forms a chain, forms a disjoint union of two chains, contains a cycle, or has a fork.
The first case contradicts the assumption $r=0$.
In the second case, we see that $G_{j,0}=S_1- [\mathrm{CT}]_{\beta_1,\beta_2}-S_2$ by the assumption $l_{R_0}=0$.
In the third case, since $r=0$, the divisor $D_{j,0}$ must contain a cycle.
Hence, we find that $S_1-[\mathrm{I}_n\mathrm{E1}]_{\beta_1,\beta_2,\gamma}-S_2$.
One sees from the assumption $l_{R_0}=r=0$ that the last case does not occur.

Suppose that $S_1$ and $S_2$ intersect in $F_j$.
It follows from $r=0$ that the two sections intersect transversally and $D_j$ must contain a pushforward of the exceptional set of an elliptic singularity on $X$.
This implies that $G_{j,0}=S_1-[\mathrm{I}_n\mathrm{E2}]-S_2$.

We complete the proof by studying the configuration of $C$.
By reordering the blow-ups if necessary, we factor $\widetilde{X}\to Y$ as $\widetilde{X}\to Y_k\to Y$ so that $Y_k\to Y$ is an isomorphism away from $\bigcup_{J'+1\leq j\leq J} D_{j,k}$, and $X\to Y_k$ is an isomorphism near $\bigcup_{J'+1\leq j\leq J} D_{j,k}$. 
Let $S_{i,k}$ $(i=1,2)$ denote the proper transform of the section $S_i$ on $Y_k$.  
Then, one sees that $-S_{i,k}^2=\chi$.  
According to the above classification (along with the assumption $l_{R_0}=0$), a connected component of $B_k$ forms a chain of the form $[2^a,\chi,2^b]$. 
Lemma~\ref{lem:extT_2m2} shows that $(a,b)=(0,\chi-4)$ or $(a,b)=(\chi-4,0)$, and the desired classification follows.
\end{proof}

\subsection{Case of $(g,n_X,l_{E_h})=(0,2,0)$ with an involution}\label{subsec:0,2,0,involution}
In this subsection, we classify $C$ and its resolution $\widetilde{C}$ under the following assumption.

\begin{ass}\label{ass:two-sections}
    \phantom{A}
    \begin{itemize}
        \item[$(1)$] 
        $(g,n_X,l_{E_h})=(0,2,0)$ and $\delta'_f=r=l_{R_0}=0$.
        \item[$(2)$] $X$ has an involution $\sigma$ which induces a fiberwise involution on $\widetilde{X}\to \PP^1$.
        \item[$(3)$] The quotient $W:=X/\sigma$ admits a (possibly trivial) P-resolution $W^\dagger\to W$ with $K_{W^\dagger}^2\geq8$.
        \item[$(4)$] The projection $\widetilde{X}/\sigma \to \PP^1$ is a ruling over $\PP^1$.
    \end{itemize}
\end{ass}

From our first assumption, we can express as the sum of two sections $C_\hor=S_1+S_2$.
(For details, see the argument in the previous section.)
The involution $\sigma$ acts on $C_\hor$ in exactly one of the following two ways: 
$\sigma$ fixes both sections $S_1$ and $S_2$, or $\sigma$ interchanges $S_1$ and $S_2$.
The next proposition shows that the first case never occurs.

\begin{prop}\label{prop:interchange}
    Under Assumption \ref{ass:two-sections}, the involution $\sigma$ interchanges the two sections $S_1$ and $S_2$.
\end{prop}

\begin{proof}
    Assume contrary that $\sigma$ fixes both sections.
    According to Theorem \ref{thm:two-sections}, there are three possible configurations for $\widetilde{C}$. 
    However, each case immediately contradicts Assumption~\ref{ass:two-sections}~(3).
    For instance, in case (1) of Theorem \ref{thm:two-sections}, the quotient $W$ has two singularities $q_1$ and $q_2$ with $\lambda(q_1)=\lambda(q_2)\geq1$.
    By Lemma~\ref{lem:lambda}, this implies $K_{W^\dagger}^2<8$, contradicting Assumption \ref{ass:two-sections} (3).
    The same argument applies to all other cases, completing the proof.
\end{proof}

We use the same notation as in the previous subsection: 
\begin{align*}
    J'':= &\sharp \{ 1\leq j\leq J : G_{j,0} = S_1-[\mathrm{I}_m \mathrm{E1}]_{0,0,\gamma}-S_2\ \text{for some $m,\gamma$} \}    \\
    + &\sharp \{ 1\leq j\leq J : G_{j,0} = S_1-[\mathrm{I}_m\mathrm{E2}]-S_2\ \text{for some $m$} \},\\
    J':=&J-J''.
\end{align*}
Reordering the branches, we may assume that $G_{j,0} = S_1-[\mathrm{I}_{n_j}\mathrm{E1}]_{0,0,\gamma}-S_2$ or $G_{j,0} = S_1-[\mathrm{I}_{n_j}\mathrm{E2}]-S_2$ for $J'+1\leq j\leq J$.
The following classification can be easily derived from Theorem \ref{thm:two-sections}.

\begin{thm}\label{thm:two-sections_inv}
    Suppose that Assumption \ref{ass:two-sections} holds.
    Then, exactly one of the following holds:
    \begin{itemize}
        \item[$(1)$] 
        $J'=1$ and $\chi \geq 5$.
        There exists an amulet $A_1$ such that
       \begin{align*}
        A_1\in &\{ [\mathrm{I}_n\mathrm{E1}]_{\chi-4,\chi-4,\gamma} : \textrm{ $n\geq2\chi-14$ or $n=2\chi-15$, $\gamma\in\{\chi-10,\chi-9,\chi-6,\chi-5\}$} \} \\
        &\cup \{  [\mathrm{CT}]_{\chi-4,\chi-4} \},
        \end{align*}
        and $G_{1,0}=S_1-A_1-S_2$.
        \item[$(2)$]
        $J'=0$, $\chi=4$.
        The two sections are disjoint.
    \end{itemize}
\end{thm}

\subsection{Case of $(g,n_X,l_{E_h})=(1,1,1)$}\label{subsec:1,1,1}
In this subsection, we classify $C$ and its resolution $\widetilde{C}$ under the assumption that $(g,n_X,l_{E_h})=(1,1,1)$ and Assumption \ref{ass:normal_stable_elliptic}.

In this situation, $C_\hor=:S$ becomes a section, and all $D_j$'s are completely separated from $S$ after the blow-ups.
Using Remark \ref{ignorecompsep}, we can reduce the problem to the case where $S$ does not have any completely separated branch.
Then, according to Lemma~\ref{lem:smoothable_elliptic_singularity}, the self-intersection number of $S$ must bigger than $-10$.
The result is the following.

\begin{thm} \label{thm:elliptic_base}
    Let $(g,n_X,l_{E_h})=(1,1,1)$.
    Assume that there exists no branch that is completely separated from the section $S$.
    Then, the minimal resolution $\widetilde{X}$ is relatively minimal, and it holds that $\widetilde{C}=C=S$ and $\chi(\mathcal{O}_X)\leq 10$. 
    
\end{thm}

\subsection{Case of $(g,n_X,l_{E_h})=(1,1,1)$ with an involution}\label{subsec:1,1,1,involution}
In this subsection, we classify $C$ and its resolution $\widetilde{C}$ under Assumption \ref{ass:normal_stable_elliptic} and the following assumption.

\begin{ass}\label{ass:elliptic_base}
    \phantom{A}
    \begin{itemize}
        \item[$(1)$] $(g,n_X,l_{E_h})=(1,1,1)$. 
        \item[$(2)$] $X$ has an involution $\sigma$ which induces a fiberwise involution on $\widetilde{X}$.
    \end{itemize}
\end{ass}

As discussed in the previous subsection, we may assume that there exists no branch that is completely separated from the section $S$.
The second condition in the above assumption provides a further constraint on $-S^2+1=\chi(\mathcal{O}_X)$ as the pushforward of $S$ must be the exceptional curve of a $\Q$-Gorenstein smoothable simple elliptic singularity. 
The result is the following.

\begin{thm}
    Suppose that Assumption \ref{ass:elliptic_base} holds.
    Assume that there exists no branch that is completely separated from the section $S$.
    Then, the minimal resolution $\widetilde{X}$ is relatively minimal, and it holds that $\widetilde{C}=C=S$ and $\chi(\mathcal{O}_X)\leq 5$.
\end{thm}

\subsection{Results}\label{subsec:non-standard_Horikawa}
The aim of this subsection is to give a classification of non-standard Horikawa surfaces with only $\Q$-Gorenstein smoothable singularities.

\subsubsection{Non-standard Horikawa surfaces}
In this subsection, we classify the possible set of non-Du Val singularities on a non-standard Horikawa surface.
To state our result, we introduce a notation.
The sentence
\begin{align*}
\text{The list of non-Du Val singularities of $X$ is } \\ [2,7,3,2^2,3] + [3,3] + \mathrm{Ell}_2 + [2,5] + [3,3] + 0
\end{align*}
means the followings:
\begin{itemize}
    \item 
    $X$ has exactly one T-singularity of type $[2,7,3,2^2,3]$, two T-singularities of type $[3,3]$, one simple elliptic singularity of type $\mathrm{Ell}_2$, and one T-singularity of type $[2,5]$.
    \item
    $X$ has no non-Du Val singularities other than those listed above.
    \item 
    We regard $0$ in the list as the one that does not represent any singularity on $X$. In other words, we ignore $0$ in the list.
\end{itemize}
We introduce collections of sets of singularities as follows:
\begin{align*}
    \mathrm{Sing}(\cA_{[2^\beta]}):=&
    \begin{cases}
   \{ 0 ,\mathrm{Ell}_{\beta+1}\}\cup  \{[3+\beta,2^{n}]^{\circ}: n\ge 0 \}                 & \text{if $8 \ge \beta \ge 1$ } \\
   \{ 0\} \cup \{[3+\beta,2^{n}]^{\circ}: n\ge \beta-9 \}      & \text{if $\beta \geq 9$ }
  \end{cases}\\
   \mathrm{Sing}(\cA_{[4]}):=&\{ 0 ,[3,2,3], (2,2,2,2)[3,2^n] \}  \\
  \mathrm{Sing}(\cA_{[3,2^{m-1},3]}):=&
  \begin{cases}
   \{0, (2,2,2,2)[4,2^n] \}                 & \text{if $m=1$ } \\
    \{0,[4],(2,2,2,2)[5,2^n]\}       & \text{if $m=2$} \\
    \{0,(2,2,2,2)[m+3,2^n]\}       & \text{if $m \geq 3$ }
  \end{cases}\\
  \mathrm{Sing}(\cA_{[5,2]}):=&\{0, [4,3,2], (3,3,3)[2] \} \\
    \mathrm{Sing}(\cA_{[4,2^{m},3,2]}):=&
  \begin{cases}
   \{0,[2,5], (3,3,3)[3] \}                 & \text{if $m=0$} \\
    \{0,(3,3,3)[4]  \}       & \text{if $m=1$} \\
  \end{cases}
\end{align*}
\begin{align*}
\mathrm{Sing}(\cA_{[5,3,2^2]}):=&\{0,(4,4,2)[3] \}\\
\mathrm{Sing}(\cA_{[3+\beta,2^{n-2},3,2^{\beta}]}):=&\{0\} \quad \text{if $n=1$, $\beta \ge 5$ or $n=2$, $\beta\ge 3$ or $n=3$, $\beta\ge 2$ or $n\ge 4$}\\
\mathrm{Sing}(\cA_{[2]}^{B,2}):=
&\left\{[3,5,3,2]+[2,5], [3,2,6,2]+[2,3,4]\right\}\\
& \cup \{(3,3,3)[3]+[3,2^{n-2},3] : n \ge 0 \} \\
& \cup \{ [5,2]+(3,3,3)[4]+[3,2^{n-2},3] : n \ge 0 \} \\
& \cup \{ [3,5,3,2]+[4,2^{\beta-1},3,2]+(3,3,3)[\beta+2] : \beta  =0,1,2 \}\\
\mathrm{Sing}(\cA_{[2]}^{B,3}):=&\left\{[3,2,6,2],
[3,2,6,2]+(3,3,3)[2]\right\}\\
\mathrm{Sing}(\cA_{\operatorname{CS}}):=&\left\{[2,5] + [3,2^{n-2},3] : n \ge 0 \right\} \cup \left\{ [2,5] + (2,2,2,2)[3,2^{n}] : n \ge 0
\right\}\\
\mathrm{Sing}(\cA_{\operatorname{CSE}}):=&\left\{0,\mathrm{Ell}_1 \right\} \cup \left\{   [3,2^{n}]^{\circ} : n \ge 0 \right\}\\
\mathrm{Sing}(\cA_{\operatorname{CSE}}^{2}):=
& \{\mathrm{Ell}_1, \mathrm{Ell}_2 \}
  \cup \{[3,2^{n}]^{\circ}: n \ge 0 \} 
 \cup \{ [3,2^{n-1-\gamma},3,2^{\gamma-1}]^{\circ} : n \ge \gamma \ge 0 \} \\
\end{align*}

\begin{thm}\label{thm:classification_non-standard_Horikawa_noninv}
    Let $X$ be a non-standard Horikawa surface with only $\Q$-Gorenstein smoothable singularities with $\chi(\mathcal{O}_X)=\chi$.
    Then, $q=0$ and $X$ is one of the following.
    (Both $\sum_{i=1}^{0} \mathcal{S}_{\operatorname{CS}}^i$ and $\sum_{j=1}^{0} \mathcal{S}_{\operatorname{CSE}}^j$ are interpreted as being $0$.)
         \begin{itemize}
          \item[$(1)$] $\chi=4$.
             The  list of non-Du Val singularities of $X$ is 
             \begin{align*}
             (2,2,2,2)[4]
             + \sum_{i=1}^{4} \mathcal{S}_i
             +\sum_{j=1}^{J_{\operatorname{CSE}}} \mathcal{S}_{\operatorname{CSE}}^j,
             \end{align*}
              where $\mathcal{S}_i  \in \mathrm{Sing}(\cA_{[2]})$, 
              $\mathcal{S}_{\operatorname{CSE}}^j\in\mathrm{Sing}(\cA_{\operatorname{CSE}})$ and  $J_{\operatorname{CSE}}\geq 0$.
       \item[$(2)$]      $\chi = 4$.
              There exist $m \ge 0$ such that the list of non-Du Val singularities of $X$ is 
              \begin{align*}
              (4,4,2)[3]
              +[4,2^{m-1},3,2]
              +\sum_{j=1}^{J_{\operatorname{CSE}}} \mathcal{S}_{\operatorname{CSE}}^j,
             \end{align*}
              where
              $\mathcal{S}_{\operatorname{CSE}}^j\in\mathrm{Sing}(\cA_{\operatorname{CSE}})$ and $J_{\operatorname{CSE}}\geq0$.
             \item[$(3)$] $\chi = 5$.
             There exist integers $n,m \ge 0$ such that 
             the  list of non-Du Val singularities of $X$ is 
             \begin{align*}
             [2,5,5,3,2^2,3]
             +[3,2^{m-2},3]
             +[4,2^{n-2},3,2]
             +\mathcal{S}_1
             +\mathcal{S}_2
             +\sum_{j=1}^{J_{\operatorname{CSE}}} \mathcal{S}_{\operatorname{CSE}}^j,
             \end{align*}
             where 
             $\mathcal{S}_1  \in \mathrm{Sing}(\cA_{[4,2^{n-1},3,2]})$,
             $\mathcal{S}_2  \in \mathrm{Sing}(\cA_{[3,2^{m-1},3]})$,
             $\mathcal{S}_{\operatorname{CSE}}^j\in\mathrm{Sing}(\cA_{\operatorname{CSE}})$ and $J_{\operatorname{CSE}}\geq0$.
              \item[$(4)$]      $\chi = 5,6$.
              The list of non-Du Val singularities of $X$ is 
              \begin{align*}
              \operatorname{Ell}_{3\chi - 11} + 
             \sum_{i=1}^{\chi - 5} \mathcal{S}_{\operatorname{CS}}^i +
             \sum_{j=1}^{J_{\operatorname{CSE}}} \mathcal{S}_{\operatorname{CSE}}^j,
              \end{align*}
              where
              $\mathcal{S}_{\operatorname{CS}}^i\in\mathrm{Sing}(\cA_{\operatorname{CS}})$ and 
              $\mathcal{S}_{\operatorname{CSE}}^j\in\mathrm{Sing}(\cA_{\operatorname{CSE}})$ and $J_{\operatorname{CSE}}\geq0$.        
              \item[$(5)$] $\chi \ge 4$. 
             There exists an integer $n \ge 0$ such that 
             the  list of non-Du Val singularities of $X$ is 
             \begin{align*}
             &[3\chi-8,3\chi-7,2^{3\chi-12},3,2^{3\chi-10}]
             +[3\chi-7,2^{n-2},3,2^{3\chi-10}]
             + \mathcal{S} +
             \sum_{i=1}^{\chi-4} \mathcal{S}_{\operatorname{CS}}^i + \sum_{j=1}^{J_{\operatorname{CSE}}} \mathcal{S}_{\operatorname{CSE}}^j,
              \end{align*}
             where 
             $\mathcal{S}  \in \mathrm{Sing}(\cA_{[3\chi-7,2^{3\chi-12+n},3,2^{3\chi-10}]})$, 
             $\mathcal{S}_{\operatorname{CS}}^i\in\mathrm{Sing}(\cA_{\operatorname{CS}})$,
             $\mathcal{S}_{\operatorname{CSE}}^j\in\mathrm{Sing}(\cA_{\operatorname{CSE}})$ and $J_{\operatorname{CSE}}\geq0$.
             \item[$(6)$] $\chi\ge 4$.
             The list of non-Du Val singularities of $X$ is 
             \[ [3\chi-6,2^{3\chi-10}]
             + \mathcal{S}
             + \sum_{i=1}^{\chi-3} \mathcal{S}_{\operatorname{CS}}^i 
             + \sum_{j=1}^{J_{\operatorname{CSE}}} \mathcal{S}_{\operatorname{CSE}}^j\]
             where 
             $\mathcal{S} \in \mathrm{Sing} (\cA_{[2^{3\chi-10}]})$,
             $\mathcal{S}_{\operatorname{CS}}^i\in\mathrm{Sing}(\cA_{\operatorname{CS}})$,
             $\mathcal{S}_{\operatorname{CSE}}^j\in\mathrm{Sing}(\cA_{\operatorname{CSE}})$, and $J_{\operatorname{CSE}}\geq0$.
             \item[$(7)$] $\chi \ge 4$. 
             There exists an integer $n \ge 0$ such that 
             the  list of non-Du Val singularities of $X$ is 
             \[
             [2,3\chi-8,3,2^{3\chi-12},3]
             +[3,2^{n-2},3]
             + \mathcal{S}_1
             + \mathcal{S}_2
             +\sum_{i=1}^{\chi-4} \mathcal{S}_{\operatorname{CS}}^i 
             + \sum_{j=1}^{J_{\operatorname{CSE}}} \mathcal{S}_{\operatorname{CSE}}^j,
             \]
             where 
             $\mathcal{S}_1  \in \mathrm{Sing}(\cA_{[3,2^{3\chi+n-12},3]})$,
             $\mathcal{S}_2  \in \mathrm{Sing}(\cA_{[2]})$,
             $\mathcal{S}_{\operatorname{CS}}^i\in\mathrm{Sing}(\cA_{\operatorname{CS}})$,
             $\mathcal{S}_{\operatorname{CSE}}^j\in\mathrm{Sing}(\cA_{\operatorname{CSE}})$
              and $J_{\operatorname{CSE}}\geq0$.
              \item[$(8)$]
              $\chi \ge 4$.
             There exist $m \ge \mathrm{max}\{3\chi -13,0\} $ such that the list of non-Du Val singularities of $X$ is 
             \begin{align*}
             (2,2,2,2)[3\chi-8,2^m]
             +  \mathcal{S}_1
             +  \mathcal{S}_2
             + \sum_{i=1}^{\chi-4} \mathcal{S}_{\operatorname{CS}}^i
             +\sum_{j=1}^{J_{\operatorname{CSE}}} \mathcal{S}_{\operatorname{CSE}}^j,
             \end{align*}
              where 
              $\mathcal{S}_1, \mathcal{S}_2 \in \mathrm{Sing}(\cA_{[2]})$,
              $\mathcal{S}_{\operatorname{CS}}^i\in\mathrm{Sing}(\cA_{\operatorname{CS}})$,
              $\mathcal{S}_{\operatorname{CSE}}^j\in\mathrm{Sing}(\cA_{\operatorname{CSE}})$ and $J_{\operatorname{CSE}}\geq0$.
               \item[$(9)$]
              $\chi \ge 4$.
             There exist $m_1,m_2 \ge 0$ with $m_1+m_2 \ge 3\chi -13$ such that the list of non-Du Val singularities of $X$ is 
             \begin{align*}
             (2,2,2,2)[2^{m_1},3\chi-8,2^{m_2}]
             + \sum_{i=1}^{\chi-4} \mathcal{S}_{\operatorname{CS}}^i
             +\sum_{j=1}^{J_{\operatorname{CSE}}} \mathcal{S}_{\operatorname{CSE}}^j,
             \end{align*}
              where 
              $\mathcal{S}_{\operatorname{CS}}^i\in\mathrm{Sing}(\cA_{\operatorname{CS}})$,
              $\mathcal{S}_{\operatorname{CSE}}^j\in\mathrm{Sing}(\cA_{\operatorname{CSE}})$ and $J_{\operatorname{CSE}}\geq0$.
 \item[$(10)$]      $\chi \ge 4$.
              The list of non-Du Val singularities of $X$ is 
              \begin{align*}
              [\chi,2^{\chi-4}] + [\chi,2^{\chi-4}]+ \sum_{j=1}^{J_{\operatorname{CSE}}} \mathcal{S}_{\operatorname{CSE}}^j,
              \end{align*}
              where
              $\mathcal{S}_{\operatorname{CSE}}^j\in\mathrm{Sing}(\cA_{\operatorname{CSE}}^2)$ and $J_{\operatorname{CSE}}\geq0$.
             \item[$(11)$] $\chi \ge 5$.
             There exist integer $n, m \ge 0$ such that 
             the  list of non-Du Val singularities of $X$ is 
             \begin{align*}
            &[3,2^{3\chi-12},3,3\chi-8,2]+[3,5,3,2] \\
            +&[4,2^{n-2},3,2]
            +[3,2^{m-2},3]
            + \mathcal{S} +
             \sum_{i=1}^{\chi-5} \mathcal{S}_{\operatorname{CS}}^i + \sum_{j=1}^{J_{\operatorname{CSE}}} \mathcal{S}_{\operatorname{CSE}}^j,
             \end{align*}
             where 
             $\mathcal{S} \in  \mathrm{Sing}(\cA_{[4,2^n,3,2]})$, 
             $\mathcal{S}_{\operatorname{CS}}^i\in\mathrm{Sing}(\cA_{\operatorname{CS}})$,
             $\mathcal{S}_{\operatorname{CSE}}^j\in\mathrm{Sing}(\cA_{\operatorname{CSE}})$ and $J_{\operatorname{CSE}}\geq0$.
 \item[$(12)$] $\chi \ge 5$.
             There exist integer $n,m \ge 0$ such that 
             the  list of non-Du Val singularities of $X$ is 
             \begin{align*}
             &[4,3\chi-10,5,2^{3\chi-13},3,2^2]\\
             +&[3,2^{m-2},3]
             +[5,2^{n-2},3,2^2]
             +\mathcal{S}_1
             +\mathcal{S}_2
             +\sum_{i=1}^{\chi-5} \mathcal{S}_{\operatorname{CS}}^i 
             +\sum_{j=1}^{J_{\operatorname{CSE}}} \mathcal{S}_{\operatorname{CSE}}^j,
             \end{align*}
             where 
             $\mathcal{S}_1  \in \mathrm{Sing}(\cA_{[3,2^{m-1},3]})$,
             $\mathcal{S}_2  \in \mathrm{Sing}(\cA_{[5,2^{3\chi-13+n},3,2^2]})$,
             $\mathcal{S}_{\operatorname{CS}}^i\in\mathrm{Sing}(\cA_{\operatorname{CS}})$,
             $\mathcal{S}_{\operatorname{CSE}}^j\in\mathrm{Sing}(\cA_{\operatorname{CSE}})$ and $J_{\operatorname{CSE}}\geq0$.
             \item[$(13)$] $\chi \ge 5$.
             There exist integer $n, m \ge 0$ such that 
             the  list of non-Du Val singularities of $X$ is 
             \begin{align*}
             &[2,6,2,3]+[2,3\chi-7,3,2^{3\chi-11},3] \\
             +& [3,2^{m-2},3]
             + [4,2^{n-2},3,2]
             +\mathcal{S}_1 
             +\mathcal{S}_2 
             +\sum_{i=1}^{\chi-5} \mathcal{S}_{\operatorname{CS}}^i 
             + \sum_{j=1}^{J_{\operatorname{CSE}}} \mathcal{S}_{\operatorname{CSE}}^j,
             \end{align*}
             where 
             $\mathcal{S}_1  \in \mathrm{Sing}(\cA_{[4,2^{n-1},3,2]})$, 
             $\mathcal{S}_2 \in \mathrm{Sing}(\cA_{[3,2^{3\chi-11+m},3]})$,
             $\mathcal{S}_{\operatorname{CS}}^i\in\mathrm{Sing}(\cA_{\operatorname{CS}})$,
             $\mathcal{S}_{\operatorname{CSE}}^j\in\mathrm{Sing}(\cA_{\operatorname{CSE}})$ and $J_{\operatorname{CSE}}\geq0$.
            \item[$(14)$] $\chi \ge 5$.
             There exist integer $n,m \ge 0$ such that 
             the  list of non-Du Val singularities of $X$ is 
             \begin{align*}
             &[3,5,2]+[3,3\chi-8,2^{3\chi-13},3,2]+[3,3\chi-7,2^{3\chi-12},3,2]\\
             +& [3,2^{m-2},3]
             +[4,2^{n-2},3,2]
             +\mathcal{S}_1
             +\mathcal{S}_2
             +\sum_{i=1}^{\chi-5} \mathcal{S}_{\operatorname{CS}}^i 
             +\sum_{j=1}^{J_{\operatorname{CSE}}} \mathcal{S}_{\operatorname{CSE}}^j,
             \end{align*}
             where 
             $\mathcal{S}_1  \in \mathrm{Sing}(\cA_{[3,2^m,3]})$,
             $\mathcal{S}_2  \in \mathrm{Sing}(\cA_{[4,2^{3\chi-12+n},3,2]})$,
             $\mathcal{S}_{\operatorname{CS}}^i\in\mathrm{Sing}(\cA_{\operatorname{CS}})$,
             $\mathcal{S}_{\operatorname{CSE}}^j\in\mathrm{Sing}(\cA_{\operatorname{CSE}})$ and $J_{\operatorname{CSE}}\geq0$.
             \item[$(15)$] $\chi \ge 5$.
             There exist integer $n, m \ge 0$ such that 
             the  list of non-Du Val singularities of $X$ is 
             \begin{align*}
             & [4,3\chi-9,2^{3\chi-14},3,2^2]
             +[4,3\chi-7,2^{3\chi-12},3,2^2]  \\
             + &[3,2^{m-2},3]
             +  [5,2^{n-2},3,2^2]  
             +  \mathcal{S}_1
             + \mathcal{S}_2
             + \sum_{i=1}^{\chi-5} \mathcal{S}_{\operatorname{CS}}^i 
             + \sum_{j=1}^{J_{\operatorname{CSE}}} \mathcal{S}_{\operatorname{CSE}}^j,
             \end{align*}
             where 
             $\mathcal{S}_1  \in \mathrm{Sing}(\cA_{[3,2^{m-1},3]})$,
             $\mathcal{S}_2  \in \mathrm{Sing}(\cA_{[5,2^{3\chi-12+n},3,2^2]})$,
             $\mathcal{S}_{\operatorname{CS}}^i\in\mathrm{Sing}(\cA_{\operatorname{CS}})$,
             $\mathcal{S}_{\operatorname{CSE}}^j\in\mathrm{Sing}(\cA_{\operatorname{CSE}})$ and $J_{\operatorname{CSE}}\geq0$.
            
     \item[$(16)$] $\chi \ge 5$.
             There exist $\gamma \in \{2,3\}$ and $m \ge 3\chi -15 +\gamma$ such that the list of non-Du Val singularities of $X$ is 
             \begin{align*}
             (2,2,2,2)[3\chi-10+\gamma,2^m]
             +  \mathcal{S}_1
             +  \mathcal{S}_2
             + \sum_{i=1}^{\chi-5} \mathcal{S}_{\operatorname{CS}}^i
             +\sum_{j=1}^{J_{\operatorname{CSE}}} \mathcal{S}_{\operatorname{CSE}}^j,
             \end{align*}
              where 
              $\mathcal{S}_1  \in \mathrm{Sing}(\cA_{[2]}^{B,\gamma})$,
              $\mathcal{S}_2  \in \mathrm{Sing}(\cA_{[2]})$,
              $\mathcal{S}_{\operatorname{CS}}^i\in\mathrm{Sing}(\cA_{\operatorname{CS}})$,
              $\mathcal{S}_{\operatorname{CSE}}^j\in\mathrm{Sing}(\cA_{\operatorname{CSE}})$ and $J_{\operatorname{CSE}}\geq0$.     
       \item[$(17)$]      $\chi \ge 5$.
              There exist integer $0 \le \gamma \le n$ such that
              the list of  non-Du Val singularities of $X$ is 
              \begin{align*}
              [\chi,2^{\chi-4}]  
              + [\chi,2^{\chi-4}]
              + [\chi-1,2^{n-\gamma-1},\chi-1,2^{\gamma-1}]^{\circ} + \sum_{j=1}^{J_{\operatorname{CSE}}} \mathcal{S}_{\operatorname{CSE}}^j,
              \end{align*}
              where
              $\mathcal{S}_{\operatorname{CSE}}^j\in\mathrm{Sing}(\cA_{\operatorname{CSE}}^2)$ and $J_{\operatorname{CSE}}\geq0$.
  \item[$(18)$]     $\chi \ge 5$.
              There exist integer $0 \le \gamma_1 \le n_1$ and $0 \le \gamma_2 \le n_2$ such that the list of  non-Du Val singularities of $X$ is 
              \begin{align*}
              [\chi,2^{\chi-4}] 
              +[\chi,2^{\chi-4}]
              + [3,2^{n_1-\gamma_1-1},\chi-1,2^{\gamma_1-1}]^{\circ}
              +[\chi-1,2^{n_2-\gamma_2-1},3,2^{\gamma_2-1}]^{\circ}
              + \sum_{j=1}^{J_{\operatorname{CSE}}} \mathcal{S}_{\operatorname{CSE}}^j,
              \end{align*}
              where
              $\mathcal{S}_{\operatorname{CSE}}^j\in\mathrm{Sing}(\cA_{\operatorname{CSE}}^2)$ and $J_{\operatorname{CSE}}\geq0$.
              \item[$(19)$] $\chi \ge 6$.
             There exist $\gamma_1, \gamma_2 \in \{2,3\}$ and $m \ge 3\chi -17 +\gamma_1+\gamma_2$ such that the list of non-Du Val singularities of $X$ is 
             \begin{align*}
             (2,2,2,2)[3\chi-12+\gamma_1+\gamma_2,2^m]
             +  \mathcal{S}_1
             +  \mathcal{S}_2
             + \sum_{i=1}^{\chi-6} \mathcal{S}_{\operatorname{CS}}^i
             +\sum_{j=1}^{J_{\operatorname{CSE}}} \mathcal{S}_{\operatorname{CSE}}^j,
             \end{align*}
              where 
              $\mathcal{S}_1  \in \mathrm{Sing}(\cA_{[2]}^{B,\gamma_1})$,
              $\mathcal{S}_2  \in \mathrm{Sing}(\cA_{[2]}^{B,\gamma_2})$,
              $\mathcal{S}_{\operatorname{CS}}^i\in\mathrm{Sing}(\cA_{\operatorname{CS}})$,
              $\mathcal{S}_{\operatorname{CSE}}^j\in\mathrm{Sing}(\cA_{\operatorname{CSE}})$ and $J_{\operatorname{CSE}}\geq0$.        
\end{itemize}
\end{thm}
\begin{rem}\label{rem--correspondence-classification}
For the reader's convenience, 
we provide tables indicating which non-standard Horikawa surfaces in Theorem~\ref{thm:classification_non-standard_Horikawa_noninv}, correspond to the cases listed in Theorem~\ref{thm:one-section-T}, \ref{thm:one-section-halfcusp} and \ref{thm:two-sections}.

\begin{table}[htbp]
\begin{flushleft}
    \begin{tabular}{|c||c|c|c|c|c|c|c|c|c|} \hline
     Theorem~ \ref{thm:classification_non-standard_Horikawa_noninv}                       &  $(3)$   & $(5)$ & $(6)$  & $(7)$ &$(11)$ & $(12)$ & $(13)$ & $(14)$ & $(15)$\\ \hline 
Theorem~\ref{thm:one-section-T}   & $(2.10)$ & $(1.2)$  &  $(1.3)$    & $(2.13)$ & $(2.11)$ & $(2.12)$ & $(2.14)$ & $(2.15)$ & $(2.16)$ \\ \hline
    \end{tabular}
\end{flushleft}
\end{table}
\begin{table}[htbp]
\begin{flushleft}
    \begin{tabular}{|c||c|c|c|c|c|c|} \hline
     Theorem~ \ref{thm:classification_non-standard_Horikawa_noninv}                       &  $(1)$   & $(8)$ & $(9)$  & $(16)$ &$(19)$  \\ \hline 
Theorem~\ref{thm:one-section-halfcusp}   & $(1)$ & $(2)$  &  $(3)$    & $(2)$ & $(2)$  \\ \hline
    \end{tabular} 
\end{flushleft}
\end{table}
\begin{table}[htbp]
\begin{flushleft}
    \begin{tabular}{|c||c|c|c|} \hline
     Theorem~ \ref{thm:classification_non-standard_Horikawa_noninv}                       &  $(10)$   & $(17)$ & $(18)$  \\ \hline 
Theorem~\ref{thm:two-sections}   & $(1),(3)$ & $(2)$  &  $(1)$    \\ \hline
    \end{tabular}
\end{flushleft}
\end{table}
\end{rem}

\begin{rem}\label{rem:triangle_Horikawa}
    In general, the existence of the surfaces listed above is nontrivial.
    According to \cite[Theorem~2.3]{Artin}, the existence of a surface in the above list with only rational singularities is equivalent to the existence of an elliptic surface with a specific Euler characteristic and with singular fibers required for constructing the desired $\widetilde{C}$.
    For example, the surfaces listed in case (2) with $m=J_{\operatorname{CSE}}=0$ can be constructed from any elliptic surface with the Euler characteristic four that has at least one singular fiber of type $\mathrm{I}_3$, so they exist.
\end{rem}

\begin{proof}[Proof of Theorem \ref{thm:classification_non-standard_Horikawa_noninv}]
According to Proposition \ref{prop:non-std_Horikawa} $(2)$, non-standard Horikawa surfaces are classified into the three types:
\begin{itemize}
    \item
    $(g,n_X,l_{E_h})=(0,1,0)$, $\chi = l_{R_0} + r + 3$ (Section \ref{subsec:0,1,0})
    \item 
    $(g,n_X,l_{E_h})=(0,2,0)$, $\delta'_f=r=l_{R_0}=0$ (Section \ref{subsec:0,2,0})
    \item 
    $(g,n_X,l_{E_h})=(1,1,1)$, $\chi = l_{R_0} + r + 5\ $ (Section \ref{subsec:1,1,1})
\end{itemize}
We deal with each case separately.

In the case where $(g,n_X,l_{E_{h}}) = (0,2,0)$, we obtain cases $(10)$, $(17)$, and $(18)$. 
This result follows immediately from Theorem~\ref{thm:two-sections}.

We consider the case $(g,n_X,l_{E_h})=(0,1,0)$.
This case is studied in Section \ref{subsec:0,1,0}.
We use the same notation there.
Specifically, let $S$ denote the unique horizontal component of $C$, and let $p\in X$ denote the singular point whose exceptional locus contains $S$.
Let $J'$ be the number of branches of $S$ that is not completely separated after the blow-ups.
After reordering the branches $D_j$ ($1\leq j\leq J$), we may assume that the branches $D_j$ ($1\leq j\leq J'$) are not completely separated.
Let $J_{\operatorname{CS}}$ (resp. $J_{\operatorname{CSE}}$) denote the number of completely separated branches of the form $(1)\text{--}(5)$ (resp. $(6)$) in Proposition~\ref{compsepamulet}.

We divide into three cases whether $p\in X$ is a T-singularity, a strictly lc singularity of type $(2,2,2,2)$, or a singularity of the type in Lemma \ref{lem:smoothable_rational_strict_lc} (ii)--(iv).

Assume that $p \in X$ is a T-singularity.
Let $J_0$ denote the number of branches $D_j$ with $1 \le j \le J'$
which corresponds to an amulet belonging to $\cA_{[2^{l}]}$ with $l \ge 1$.
Then we have
\[
r + l_{R_0} = J_{\operatorname{CS}} + (J' - J_0).
\]
We will prove this.
Let $D_j$ be a branch as described above.
An initial blow-up process corresponding to an amulet in $\cA_{[2^{l}]}$ does not contribute to $r$.
Furthermore, the divisor $D_{j,0} - B_{j,0}$ does not contain any exceptional divisor of a strictly lc rational singularity of the type in Lemma~\ref{lem:smoothable_rational_strict_lc}~(i).
Therefore, such a branch $D_j$ does not contribute to $l_{R_0} + r$.
Now suppose that a branch $D_j$ with $1 \le j \le J'$ corresponds to an amulet not belonging to $\cA_{[2^l]}$ for any $l \ge 1$.
If $D_j$ corresponds to the amulet $[\mathrm{I}_n \mathrm{R}]_{\beta}$, then the initial blow-up process corresponding to $D_j$ does not contribute to $r$.
On the other hand, the divisor $D_{j,0} - B_{j,0}$ contains exactly one exceptional divisor of a strictly lc singularity  of type $(2,2,2,2)$.
Hence, this initial blow-up process contributes $1$ to $l_{R_0} + r$.
For a branch $D_j$ with $1 \le j \le J'$ corresponding to an amulet other than $[\mathrm{I}_n \mathrm{R}]_{\beta}$,
a straightforward computation shows that the initial blow-up process corresponding to $D_j$ contributes $1$ to $r$.
On the other hand, the divisor $D_{j,0} - B_{j,0}$ does not contain any exceptional divisor of a strictly lc singularity  of type $(2,2,2,2)$.
Therefore, such a branch $D_j$ also contributes $1$ to $l_{R_0} + r$.
Similarly, a branch corresponds to one of the cases (1) through (5) in Proposition~\ref{compsepamulet} also contributes $1$ to $l_{R_0}+r$.
Thus, we obtain $r+l_{R_0}=J_{\operatorname{CS}}+(J'-J_0)$.

Since we have $\chi=r+l_{R_0}+3$, we obtain
\[
\chi=J_{\operatorname{CS}}+(J'-J_0)+3.
\]
Now, since we are assuming that $p \in X$ is a $T$-singularity, note that $0 \le J_0 \le J' \le 2$.
Hence, the possible combinations of $(J', J_0)$ are $(0,0)$, $(1,1)$, $(1,0)$, $(2,0)$, $(2,1)$ and $(2,2)$.
Here, we explain only the case where $(J', J_0) = (2, 1)$.
Substituting $J'=2$ and $J_0=1$ into the equation $\chi=J_{\operatorname{CS}}+(J'-J_0)+3$, 
we obtain $J_{\operatorname{CS}}=\chi-4$.
As mentioned in Remark~\ref{ignorecompsep}, if there exists a branch that is completely separated from the section $S$, then we can apply the list given in Theorem~\ref{thm:one-section-T} by interpreting $\chi$ in that list as $\chi + 2J_{\operatorname{CS}}$.
Since $X$ is a Horikawa surface, we have $\chi \ge 4$.
In particular, $\chi+2J_{\operatorname{CS}}=3\chi-8 \ge 4$.
In the list in Theorem~\ref{thm:one-section-T}, the only possibilities satisfying $J_{0}=1$ are $(2.9)$ and $(2.13)$.
The former case is excluded.
Indeed, it must be $\chi+2J_{\mathrm{CS1}}=3\chi-8 =5$, but it is impossible.
We now consider the latter case.
According to the list in Theorem~\ref{thm:one-section-T},  
the surface $X$ has at most two $T$-singularities:  
one of type  
\[
[2,3\chi - 8,3,2^{3\chi - 12},3],
\]  
and possibly another of type  
\[
[3,2^{n-2},3] \quad \text{for some } n \ge 1.
\]
Moreover, the amulets corresponding to the branches $D_j$ with $j = 1,2$ belong to $\cA_{[3,2^{3\chi+n-12},3]}$ and $\cA_{[2]}$, respectively (after relabeling the indices if necessary).
Therefore, the list of non-Du Val singularities on $X$ is given by
\[
[2,3\chi-8,3,2^{3\chi-12},3]
+ [3,2^{n-2},3]
+ \mathcal{S}_1
+ \mathcal{S}_2
+ \sum_{i=1}^{\chi-4} \mathcal{S}_{\operatorname{CS}}^i 
+ \sum_{j=1}^{J_{\operatorname{CSE}}} \mathcal{S}_{\operatorname{CSE}}^j,
\]
where 
$\mathcal{S}_1  \in \mathrm{Sing}(\cA_{[3,2^{3\chi+n-12},3]})$,
$\mathcal{S}_2  \in \mathrm{Sing}(\cA_{[2]})$,
$\mathcal{S}_{\operatorname{CS}}^i \in \mathrm{Sing}(\cA_{\operatorname{CS}})$, and
$\mathcal{S}_{\operatorname{CSE}}^j \in \mathrm{Sing}(\cA_{\operatorname{CSE}})$.
It follows that $X$ corresponds to case $(7)$.

Assume $p \in X$ is a strictly lc rational singularity
of type $(2,2,2,2)$.
Such a non-standard Horikawa surface $X$ belongs to one of the three types (1), (2), or (3) in the classification given in Theorem~\ref{thm:one-section-halfcusp}.
Here, we focus on the case where $X$ is of type~(2) in Theorem~\ref{thm:one-section-halfcusp}.
In this case, we have $J'=3$.
Let $A_j$ denote an amulet corresponding to the branch $D_j$ for $j = 1, 2, 3$.
If necessary, we relabel the indices so that $A_1, A_2 \in \cA_{[2]}^{\mathrm{B}}$ and $A_3 \in \cA_{\mathrm{SR}}$.
Let $J_0$ be the number of branches among $D_1$ and $D_2$ such that the corresponding amulet $A_j$ belongs to $\cA_{[2]}$.
Then $0 \le J_0 \le 2$.
By an argument similar to that in the case where $p \in X$ is a $T$-singularity, we obtain 
$l_{R_0} + r = J_{\operatorname{CS}} + (2 - J_0) + 1$.
Substituting into $\chi = l_{R_0} + r + 3$, we find
$\chi = J_{\operatorname{CS}} + (2 - J_0) + 4$.
Therefore, the exceptional divisor over $p \in X$ is of the form
\[
(2,2,2,2)[3\chi - 16 + 2J_0 - l_{E_v} + \sum_{j=1}^{2} \gamma(A_j) + J_{\operatorname{CSE}},\,2^{m(A_3)}].
\]
When $J_0 = 0$, the branches $D_1$ and $D_2$ do not correspond to the amulet $[\mathrm{I}_n \mathrm{E}]_1$.
Hence, we have $J_{\operatorname{CSE}} = l_{E_v}$, and $X$ corresponds to case~$(19)$.
Similarly, when $J_0 = 1$ or $J_0 = 2$, the surface corresponds to cases~$(16)$ and~$(8)$, respectively.

Assume $p \in X$ is a strictly lc rational singularity of type (ii), (iii) or (iv) in Lemma~\ref{lem:smoothable_rational_strict_lc}.
Such a non-standard Horikawa surface $X$ is necessarily of one of the four types (1), (2), (3), or (4) as given in the classification of Theorem~\ref{thm:one-section-triangle}.
Among these, types (1), (2), and (3) are excluded.
Here, we show that type (1) is excluded.
In this case, we have $J'=3$.
Let $A_j \in \cA_{[3]}^{B}$ denote an amulet corresponding to the branch $D_j$ for $j = 1, 2, 3$.
As in the argument that $p \in X$ is a T-singularity, we obtain
$\chi=J_{\mathrm{CS1}}+6$.
In particular, $\chi \ge 6$.
Therefore, the exceptional divisor of $p \in X$ is 
\[
(3,3,3)[3 \chi-12 + \sum_{j=1,2,3} \gamma(A_j)].
\]
Since $\gamma(A_i) \ge 1$ for $i=1,2,3$,
it follows that $3\chi -12 \le 4 - \sum_{j=1,2,3}\gamma(A_i) \le 1$.
Therefore, we get $3 \chi \le 13$.
This leads to a contradiction with the fact that $\chi \ge 6$.

We consider case~(4) of Theorem~\ref{thm:one-section-triangle}, where $J' = 1$.
Arguing as in the case where $p \in X$ is a $T$-singularity, we obtain $l_{R_0} + r = 1$.
Hence, $\chi = J_{\operatorname{CS}} + 4$.
In particular, it follows that $\chi \ge 4$.
Since the self-intersection number of the section $S$ on $\widetilde{X}$ is $(-4)$, we obtain $4 = \chi + 2(\chi - 4)$.
It follows that $\chi = 4$ and $J_{\operatorname{CS}} = 0$, which corresponds to case $(2)$.

In the remainder of the proof, 
we consider the case $(g, n_X, l_{E_h}) = (1, 1, 1)$,  
which is analyzed in detail in Section~\ref{subsec:1,1,1}.  
Let $p \in X$ be the singular point whose exceptional divisor contains the section $S$.  
Let $J_{\operatorname{CS}}$ (resp. $J_{\operatorname{CSE}}$) denote the number of completely separated branches of the form $(1)\text{--}(5)$ (resp. $(6)$) in Proposition~\ref{compsepamulet}.

Arguing similarly to the case $(g, n_X, l_{E_h}) = (0, 1, 0)$, 
we obtain $r + l_{R_0} = J_{\operatorname{CS}}$.  
It follows that $\chi = J_{\operatorname{CS}} + 5$.

By Theorem~\ref{thm:elliptic_base}, 
the surface $X$ has a simple elliptic singularity of type $\operatorname{Ell}_{3\chi - 11}$.  
Since $p \in X$ is $\Q$-Gorenstein smoothable, 
we have $3\chi - 11 \le 9$.  
This yields the possible pairs $(\chi, J_{\operatorname{CS}}) = (5, 0)$ or $(6, 1)$.
Therefore, the list of non-Du Val singularities on $X$ is given by
\[
\operatorname{Ell}_{3\chi - 11} + 
\sum_{i=1}^{\chi - 5} \mathcal{S}_{\operatorname{CS}}^i +
\sum_{j=1}^{J_{\operatorname{CSE}}} \mathcal{S}_{\operatorname{CSE}}^j,
\]
which corresponds to case~(4).
\end{proof}

\subsubsection{Non-standard Horikawa surfaces with good involutions}
We conclude this section with the classification of non-standard Horikawa surafaces that admit a good involution, together with explicit constructions of the surfaces listed.

\begin{thm}\label{thm:classification_non-standard_Horikawa}
    Let $X$ be a non-standard Horikawa surface with $\chi = \chi(\mathcal{O}_X)$ and only $\Q$-Gorenstein smoothable singularities that admitting a good involution $\sigma$.
    Then, $q=0$ and the list of non-mild singularities of $X$ is one of the following.
   \begin{enumerate}[label=$(\arabic*)$]
        \item 
        $(\chi=4)$ $[2,4,3,3]$
        \item $(\chi=4)$ $[2,4,3,3]$, $[4]$, $[4]$.
        \item $(\chi = 4,\ n \ge 0)$ $[2,4,3,3], (2,2,2,2)[4,2^n]$
        \item 
        $(\chi=4, m_1\geq m_2\geq0)$ $(2,2,2,2)[2^{m_1},4,2^{m_2}]$
        \item 
        $(\chi\geq4)$ $[\chi,2^{\chi-4}]$, $[\chi,2^{\chi-4}]$
        \item 
        $(5\leq\chi\leq7)$ $[\chi,2^{\chi-4}]$, $[\chi,2^{\chi-4}]$, $\Ell_{2\chi-6}$ 
        \item 
        $(\chi\geq 5,\ m_1\geq-1,\ m_2\geq0)$ $[\chi,2^{\chi-4}]$, $[\chi,2^{\chi-4}]$, $[\chi-1,2^{m_1},\chi-1,2^{m_2}]^\circ$ 
        \item 
        $(\chi=5)$ $\Ell_4$ 
    \end{enumerate}
\end{thm}
\begin{proof}
According to Proposition \ref{prop:non-std_Horikawa} $(2)$, non-standard Horikawa surfaces are classified into the three types:
\begin{itemize}
    \item
    $(g,n_X,l_{E_h})=(0,1,0)$, $\chi = l_{R_0} + r + 3$ (Section \ref{subsec:0,1,0})
    \item 
    $(g,n_X,l_{E_h})=(0,2,0)$, $\delta'_f=r=l_{R_0}=0\ $ (Section \ref{subsec:0,2,0})
    \item 
    $(g,n_X,l_{E_h})=(1,1,1)$, $\chi = l_{R_0} + r + 5\ $ (Section \ref{subsec:1,1,1})
\end{itemize}
We deal with each case separately.

We consider the case $(g,n_X,l_{E_h})=(0,2,0)$.
In this case, the list of non-Du Val singularities on $X$ corresponds to one of $(10)$, $(17)$, or $(18)$ in Theorem~\ref{thm:classification_non-standard_Horikawa_noninv}.
The case (18) in Theorem~\ref{thm:classification_non-standard_Horikawa_noninv} corresponds to case (1) in Theorem~\ref{thm:two-sections} (see Remark~\ref{rem--correspondence-classification}), 
but this case does not admit a good involution.
Therefore, the list of non-Du Val singularities of $X$ must be one of $(10)$ or $(17)$ in Theorem~\ref{thm:classification_non-standard_Horikawa_noninv}.
If $X$ has a simple elliptic singularity of type $\operatorname{Ell}_1$ or $\operatorname{Ell}_2$, or a cusp singularity of type $[3, 2^{n - 1 - \gamma}, 3, 2^{\gamma - 1}]$ with $0 \le \gamma \le n$, then these singularities can be verified to be mild with respect to the good involution.  
Thus, we obtain cases $(5)$, $(6)$, and $(7)$.

In the case where $(g, n_X, l_{E_h}) = (1, 1, 1)$, 
we obtain case~$(8)$.  
This follows immediately from Theorem~\ref{thm:classification_non-standard_Horikawa_noninv} and Proposition~\ref{prop_inv_cs}.

In the rest of the proof, we consider the case $(g,n_X,l_{E_h})=(0,1,0)$.
This case is studied in Section \ref{subsec:0,1,0}.
We use the same notation there.
Specifically, let $S$ denote the unique horizontal component of $C$, and let $p\in X$ denote the singular point whose exceptional locus contains $S$.
Let $J_{\operatorname{CS}}$ denote the number of completely separated branches of $S$ with the form $(1)\text{--}(5)$ in Proposition~\ref{compsepamulet}.
We divide into three cases whether $p\in X$ is a T-singularity, a lc singularity of type $(2,2,2,2)$, or a singularity of the type in Lemma \ref{lem:smoothable_rational_strict_lc} (ii)--(iv).

Assume that $p \in X$ is a T-singularity.
The list of non-Du Val singularities of $X$ corresponds to one of $(3), (5), (6), (7)$, $(11)\text{--}(15)$ in Theorem~\ref{thm:classification_non-standard_Horikawa_noninv}.
As mentioned in Remark~\ref{rem--correspondence-inv-no_inv-T}, if $X$ admits a good involution, then all cases in the list of Theorem~\ref{thm:one-section-T}, except for $(1.2)$, $(1.3)$, and $(2.13)$, are excluded.
Therefore, the list of non-Du Val singularities of $X$ corresponds to either case~$(6)$ or case~$(7)$ in Theorem~\ref{thm:classification_non-standard_Horikawa_noninv} (see Remark~\ref{rem--correspondence-classification}).
Proposition~\ref{prop_inv_cs} asserts that $J_{\operatorname{CS}} = 0$.  
However, in case~$(6)$, we have $J_{\operatorname{CS}} = \chi-3 \ge 1$.
Hence, case~$(6)$ in Theorem~\ref{thm:classification_non-standard_Horikawa_noninv} is also excluded.
It follows that the list of non-Du Val singularities of $X$ corresponds to $(7)$ in Theorem~\ref{thm:classification_non-standard_Horikawa_noninv}.
In this case, we have $J_{\operatorname{CS}} = \chi-4$.
Consequently, $X$ corresponds to one of $(2.6)$, $(2.7)$, or $(2.8)$ in Theorem~\ref{thm:one-section-fixed-T} and $\chi=4$.
Thus, we obtain cases $(1)$, $(2)$, and $(3)$.

Assume $p \in X$ is a strictly lc rational singularity of type $(2,2,2,2)$.
The list of non-Du Val singularities of $X$ corresponds to one of $(1)$, $(8)$, $(9)$, $(16)$, and $(19)$ in Theorem~\ref{thm:classification_non-standard_Horikawa_noninv}.
The cases $(16)$ and $(19)$ in Theorem~\ref{thm:classification_non-standard_Horikawa_noninv} are excluded.
Indeed, if $X$ corresponds to one of these cases, then the quotient $W$ by the involution $\sigma$ satisfies $\lambda(W) \ge 2$ (see Lemma~\ref{lem:sigma-amulet-b}). 
However, this contradicts the assumption that $\sigma$ is a good involution.
Proposition~\ref{prop_inv_cs} asserts that $J_{\operatorname{CS}} = 0$.  
Therefore, it follows that the list of non-Du Val singularities of $X$ corresponds to $(1)$, $(8)$, and $(9)$ in Theorem~\ref{thm:classification_non-standard_Horikawa_noninv}.
In these cases, we have $J_{\operatorname{CS}} = \chi-4$.
Thus, we obtain case $(4)$.

Assume $p \in X$ is a strictly lc rational singularity of type (ii), (iii) or (iv) in Lemma~\ref{lem:smoothable_rational_strict_lc}.
The list of non-Du Val singularities of $X$ corresponds to case~$(2)$ in Theorem~\ref{thm:classification_non-standard_Horikawa_noninv}. 
However, such a surface does not admit a good involution. Therefore, this case is excluded.
\end{proof}

We describe an explicit construction of $X$ for each case listed in Theorem \ref{thm:classification_non-standard_Horikawa}.

\begin{const}[Non-standard Horikawa surfaces with good involutions] \label{construction}

In Construction~\ref{construction}, the labels $(1)$\text{--}$(7)$ correspond to those in Theorem~\ref{thm:classification_non-standard_Horikawa}.
\begin{itemize}
 \item[(1)]
 Let $\Gamma_{1}$ and $\Gamma_2$ be two distinct fibers of the Hirzebruch surface $\Sigma_{7}\to \mathbb{P}^1$.
 Take a reduced effective divisor $B_{\Sigma}= \Delta_0+B_{\Sigma}' \in |4\Delta_{0}+22 \Gamma|$ on $\Sigma_{7}$, satisfying the following conditions:   
\begin{itemize}  
    \item In a neighborhood of $\Gamma_1$, the branch divisor $B_{\Sigma}$ is described as one of the cases listed in Table~\ref{table_generate_[CT]}.
    
    \item In a neighborhood of $\Gamma_2$, $B_{\Sigma}$ is described as 
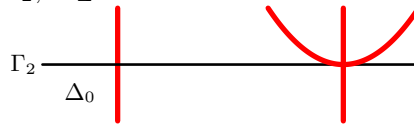
\begin{center}
\begin{tikzpicture}[line cap=round,line join=round,>=triangle 45,x=1cm,y=0.75cm]
\draw [line width=2pt, red] (-1.5,-1)-- (-1.5,1);
\draw [line width=2pt, red] (1.5,-1)-- (1.5,1);
\draw [line width=1pt] (2.5,0)-- (-2.5,0);
\draw[line width=2pt,color=black,smooth,samples=100,domain=0.5:2.5,red] plot(\x,{(\x - 1.5)^2});
\begin{scriptsize}
\draw [color=black] (-2.75,0) node {$\Gamma_2$};
\draw [color=black] (-2,-0.5) node {$\Delta_0$};
\end{scriptsize}
\end{tikzpicture}
\end{center}

The singularity of $B_{\Sigma}$ in the figure is a double point.
    \item The pair $(\Sigma_7,B_{\Sigma})$ is mild on $\Sigma_7\setminus \Gamma_1\cup \Gamma_2$. 

\end{itemize}

The existence of such $B_{\Sigma}$ is guaranteed by Remark~\ref{rem--existence-all--q-gor--smoothable}. 
We construct the Horikawa surface in Theorem \ref{thm:classification_non-standard_Horikawa} (1).
We note that Horikawa surfaces in Theorem \ref{thm:classification_non-standard_Horikawa} (1) can be constructed regardless of the specific branch type chosen from Table~\ref{table_generate_[CT]}.

First, we describe the blow-up procedure near $\Gamma_1$.  
Blow up at the intersection of the section $\Delta_0$ and the fiber $\Gamma_1$.

Next, we describe the blow-up procedure near $\Gamma_2$.
Blow up at the double point of $B_{\Sigma}$ on $\Gamma_2$.  
According to the even resolution process (Definition~\ref{defn--evenresol}), the exceptional curve $E$ is not contained in the new branch locus.  
Blow up at the intersection of $E$ and the proper transform of the fiber $\Gamma_2$.
Then, the new branch locus contains the exceptional curve and has a double point.
Blow up at this double point.
Let $\widetilde{\Sigma}$ be the surface obtained from $\Sigma_7$ via the above procedure.
The branch divisor $B_{\widetilde{\Sigma}}$ on $\widetilde{\Sigma}$ is determined according to the even resolution process.
As a result, a neighborhood of the total transform of $\Gamma_1 + \Gamma_2 + \Delta_0$ on $\widetilde{\Sigma}$ is illustrated in the following figure.

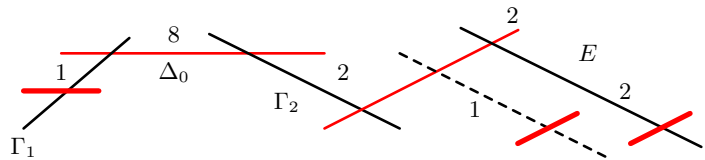
\begin{center}
\begin{tikzpicture}[line cap=round,line join=round,>=triangle 45,x=1cm,y=1cm]
\clip(-6,-1) rectangle (5,2);
\draw [line width=1pt,red] (-4.5,1)-- (-1,1);
\draw [line width=1pt] (-5,0)-- (-3.593917603691098,1.203768293279733);
\draw [line width=1pt] (0,0)-- (-2.5316725820337616,1.2568805443626);
\draw [line width=1pt,red] (-1,0)-- (1.5756748350412713,1.309992795445467);
\draw [line width=1pt,dashed] (0,1)-- (2.7825032068686335,-0.398451281053423);
\draw [line width=1pt] (0.8114485555711322,1.3365489209869006)-- (4.009986343005999,-0.24501588903625152);
\draw [line width=2pt,red] (-5,0.5)-- (-4,0.5);
\draw [line width=2pt,red] (1.57,-0.2)-- (2.36,0.2);
\draw [line width=2pt,red] (3.07,-0.2)-- (3.86,0.2);
\begin{scriptsize}
\draw [fill] (-4.5,0.75) node {$1$};
\draw [fill] (-3,1.25) node {$8$};
\draw [fill] (1.5,1.5) node {$2$};
\draw [fill] (-0.75,0.75) node {$2$};
\draw [fill] (1,0.25) node {$1$};
\draw [fill] (3,0.5) node {$2$};
\draw [fill] (-5,-0.25) node {$\Gamma_1$};
\draw [fill] (-1.5,0.35) node {$\Gamma_2$};
\draw [fill] (2.5,1) node {$E$};
\draw [fill] (-3,0.75) node {$\Delta_0$};
\end{scriptsize}
\end{tikzpicture}
\end{center}

The red line in the figure represents the branch locus $B_{\widetilde{\Sigma}}$.  
In particular, the bold red line indicates the horizontal part of $(B_{\widetilde{\Sigma}}-\Delta_0)$.  
The thin and dashed lines represent the total transform of $\Delta_0+\Gamma_1 + \Gamma_2$.  
A number written near a line denotes the minus self-intersection number of the corresponding curve.  
In the figure, $\Gamma_i$ and $E$ denote their proper transforms.

Let $[1,8,2,2,2]$ denote the chain on $\widetilde{\Sigma}$ represented by the solid lines in the figure above.
By taking the double cover of $\widetilde{\Sigma}$ branched along $B_{\widetilde{\Sigma}}$, the preimage of the chain $[1,8,2,2,2]$ becomes $[2,4,4,1,4]$.  
Contracting the chain $[2,4,4,1,4]$,
we obtain the desired Horikawa surface.

Conversely, it can also be verified that all surfaces belonging to $(1)$ in Theorem~\ref{thm:classification_non-standard_Horikawa} can be obtained through the above construction.

 \item[(2)]
 Let $\Gamma_{1}$ and $\Gamma_2$ be two distinct fibers of the Hirzebruch surface $\Sigma_{7}\to \mathbb{P}^1$.
 Take a reduced effective divisor $B_{\Sigma}= \Delta_0+B_{\Sigma}' \in |4\Delta_{0}+22 \Gamma|$ on $\Sigma_{7}$ satisfying the following conditions:   
\begin{itemize}  
    \item In a neighborhood of $\Gamma_1$, the branch divisor $B_{\Sigma}$ corresponds to  one of the cases listed in Table~\ref{table_generate_[CT]}.
    
    \item In a neighborhood of $\Gamma_2$, $B_{\Sigma}$ is described as follows:

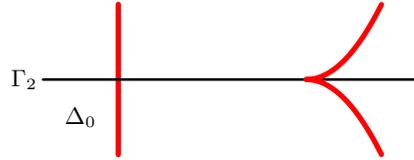
\begin{center}
\begin{tikzpicture}[line cap=round,line join=round,>=triangle 45,x=1cm,y=1cm]
\draw [line width=2pt, red] (-1.5,-1)-- (-1.5,1);
\draw [line width=1pt] (2.5,0)-- (-2.5,0);
\draw[line width=2pt,color=black,smooth,samples=100,domain=1:2,red] plot(\x,{(\x - 1)^2});
\draw[line width=2pt,color=black,smooth,samples=100,domain=1:2,red] plot(\x,{0-(\x - 1)^2});
\begin{scriptsize}
\draw [color=black] (-2.75,0) node {$\Gamma_2$};
\draw [color=black] (-2,-0.5) node {$\Delta_0$};
\end{scriptsize}
\end{tikzpicture}
\end{center}

The singularity of $B_{\Sigma}$ in the figure is a double point.

\item The pair $(\Sigma_7,B_{\Sigma})$ is mild on $\Sigma_7\setminus \Gamma_1\cup \Gamma_2$.

\end{itemize}

The existence of such $B_{\Sigma}$ is guaranteed by Remark~\ref{rem--existence-all--q-gor--smoothable}.
We proceed with the construction of the desired Horikawa surface.

First, we describe the blow-up procedure near $\Gamma_1$.  
Blow up at the intersection of the section $\Delta_0$ and the fiber $\Gamma_1$.

Next, we describe the blow-up procedure near $\Gamma_2$.
Blow up at the double point of $B_{\Sigma}$ on $\Gamma_2$.   
According to the even resolution process (Definition~\ref{defn--evenresol}), the exceptional curve $E$ is not contained in the new branch locus.
Then, blow up the intersection point of the proper transform of $\Gamma_2$ and $E$.
The new branch locus then contains the exceptional curve and has a double point.
We blow up this double point.
Finally, blow up the intersection point of the proper transform of $E$ and the resulting exceptional curve.

Let $\widetilde{\Sigma}$ be the surface obtained from $\Sigma_{7}$ by the above procedures.
The branch divisor $B_{\widetilde{\Sigma}}$ on $\widetilde{\Sigma}$ is determined according to the even resolution process.
As a result, a neighborhood of the total transform of $\Gamma_1 + \Gamma_2 + \Delta_0$ on $\widetilde{\Sigma}$ is illustrated in the following figure.

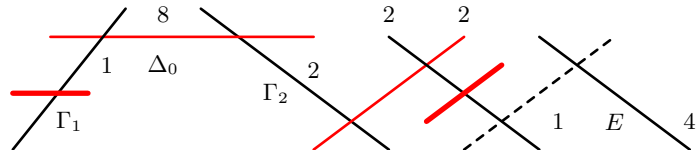
\begin{center}
\begin{tikzpicture}[line cap=round,line join=round,>=triangle 45,x=1cm,y=1.5cm]
\clip(-5.5,0) rectangle (5,1.5);
\draw [line width=1pt] (-2.5,1.25)-- (0,0);
\draw [line width=1pt,red] (-1,0)-- (1,1);
\draw [line width=1pt] (0,1)-- (2,0);
\draw [line width=1pt,dashed] (1,0)-- (3,1);
\draw [line width=1pt] (2,1)-- (4,0);
\draw [line width=1pt,red] (-4.5,1)-- (-1,1);
\draw [line width=1pt] (-3.5,1.25)-- (-5,0);
\draw [line width=2pt,red] (-5,0.5)-- (-4,0.5);
\draw [line width=2pt,red] (0.5,0.25)-- (1.5,0.75);

\begin{scriptsize}
\draw [fill] (-3.75,0.75) node {$1$};
\draw [fill] (-3,1.2) node {$8$};
\draw [fill] (0,1.2) node {$2$};
\draw [fill] (-1,0.7) node {$2$};
\draw [fill] (1,1.2) node {$2$};
\draw [fill] (2.25,0.25) node {$1$};
\draw [fill] (4,0.25) node {$4$};
\draw [fill] (-4.25,0.25) node {$\Gamma_1$};
\draw [fill] (-1.5,0.5) node {$\Gamma_2$};
\draw [fill] (3,0.25) node {$E$};
\draw [fill] (-3,0.75) node {$\Delta_0$};

\end{scriptsize}
\end{tikzpicture}
\end{center}

In the figure, $\Gamma_i$ and $E$ denote their proper transforms.
Let $[1,8,2,2,2]$ and $[4]$ denote the chains on $\widetilde{\Sigma}$ represented by the solid lines in the figure above.
By taking the double cover of $\widetilde{\Sigma}$ branched along $B_{\widetilde{\Sigma}}$, the preimages of the two chain $[1,8,2,2,2]$ and $[4]$ become $[2,4,4,1,4]$ and two chains $[4]$, respectively.  
Contracting the chain $[2,4,4,1,4]$ and the two chains $[4]$,
we obtain the desired Horikawa surface.

Conversely, it can also be verified that all surfaces belonging to $(2)$ in Theorem~\ref{thm:classification_non-standard_Horikawa} can be obtained through the above construction.

 \item[(3)]
 Let $\Gamma_{1}$ and $\Gamma_2$ be two distinct fibers of the Hirzebruch surface $\Sigma_{7}\to \mathbb{P}^1$.
 Take a reduced effective divisor
  $B_{\Sigma} = \Delta_0+B_{\Sigma}' \in |4\Delta_{0}+ 22\Gamma|$ on $\Sigma_{7}$ satisfying the following conditions:   
\begin{itemize}  
    \item In a neighborhood of $\Gamma_1$, the branch divisor $B_{\Sigma}$ corresponds to one of the cases listed in Table~\ref{table_generate_[CT]}.
    \item In a neighborhood of $\Gamma_2$, $B_{\Sigma}$ is described as one of the following:

      \begin{table}[H]
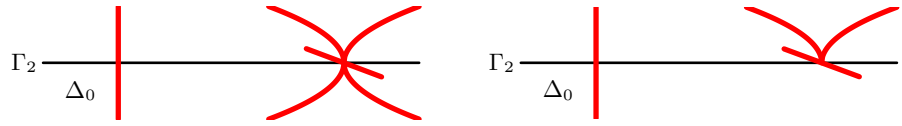
  \begin{tabular}{cc} 
    \begin{minipage}{70mm}
\begin{tikzpicture}[line cap=round,line join=round,>=triangle 45,x=1cm,y=0.75cm]
\clip(-4,-1) rectangle (3,1);
\draw [line width=1pt] (-2.5,0)--(2.5,0);
\draw [line width=2pt,color=red] (-1.5,-1)-- (-1.5,1);
\draw [line width=2pt,color=red] (1,0.25)-- (2,-0.25);
\draw[line width=2pt,color=red,smooth,samples=100,domain=1.5:2.5] plot(\x,{sqrt((\x)-1.5)});
\draw[line width=2pt,color=red,smooth,samples=100,domain=1.5:2.5] plot(\x,{0-sqrt((\x)-1.5)});
\draw[line width=2pt,color=red,smooth,samples=100,domain=0.5:1.5] plot(\x,{sqrt(-((\x)-1.5))});
\draw[line width=2pt,color=red,smooth,samples=100,domain=0.5:1.5] plot(\x,{0-sqrt(-((\x)-1.5))});
\begin{scriptsize}
\draw [color=black] (-2.75,0) node {$\Gamma_2$};
\draw [color=black] (-2,-0.5) node {$\Delta_0$};
\end{scriptsize}
\end{tikzpicture}
     \end{minipage} &
    \begin{minipage}{70mm}
\begin{tikzpicture}[line cap=round,line join=round,>=triangle 45,x=1cm,y=0.74cm]
\clip(-3,-1) rectangle (3,1);
\draw [line width=1pt] (2.5,0)-- (-2.5,0);
\draw [line width=2pt,color=red] (-1.5,-1)-- (-1.5,1);
\draw [line width=2pt,color=red] (1,0.25)-- (2,-0.25);
\draw[line width=2pt,color=red,smooth,samples=100,domain=1.5:2.5] plot(\x,{sqrt((\x)-1.5)});
\draw[line width=2pt,color=red,smooth,samples=100,domain=0.5:1.5] plot(\x,{sqrt(-((\x)-1.5))});
\begin{scriptsize}
\draw [color=black] (-2.75,0) node {$\Gamma_2$};
\draw [color=black] (-2,-0.5) node {$\Delta_0$};
\end{scriptsize}
\end{tikzpicture}
    \end{minipage} \\ 
    
\end{tabular}
\end{table}
The triple point of $B_{\Sigma}$ on $\Gamma_2$ consists of an $k$-fold node or cusp that is not tangent to $\Gamma_2$, and a smooth curve that meets $\Gamma_2$ transversally.
When $B_{\Sigma}$ contains $k$-fold node, the configuration appears as shown in the left figure,
and when $B_{\Sigma}$ contains $k$-fold cusp, it appears as shown in the right figure.   

\item The pair $(\Sigma_{7},B_{\Sigma})$ is mild on $\Sigma_{7}\setminus \Gamma_1\cup \Gamma_2$.

\end{itemize}
We will give a detailed explanation in the case where the branch locus $B_{\Sigma}$ contains a $k$-fold node in a neighborhood of $\Gamma_2$.
The existence of such $B_{\Sigma}$ is guaranteed by Remark~\ref{rem--existence-all--q-gor--smoothable}.
First, we describe the blow-up process near $\Gamma_1$.  
Blow up at the intersection of the section $\Delta_0$ and the fiber $\Gamma_1$.

Next, we describe the blow-up process near $\Gamma_2$.
Blow up at the triple point of $B_{\Sigma}$ on $\Gamma_2$.  
Let $\Gamma_2'$ be the proper transform of $\Gamma_2$, and let $E$ be the exceptional curve.
Next, blow up the intersection point of $\Gamma_2'$ and $E$.
Then, perform two further successive blow-ups at the intersection point of the total transform of $\Gamma_2'$ and the proper transform of $E$ at each step.
A neighborhood of the total transform of $\Gamma_1 + \Gamma_2 + \Delta_0$ is illustrated in the following figure.

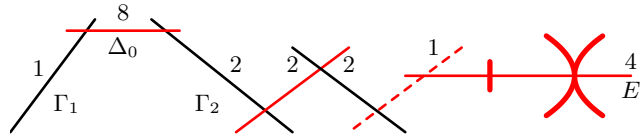
\begin{center}
\begin{tikzpicture}[line cap=round,line join=round,>=triangle 45,x=0.75cm,y=0.75cm]
\clip(-9,-1.5) rectangle (4,1.5);
\draw [line width=1pt] (-7.5,1)-- (-9,-1);
\draw [line width=1pt] (-4,-1)-- (-6.5,1);
\draw [color=red,line width=1pt] (-8,0.8)-- (-6,0.8);
\draw [color=red,line width=1pt] (-3,0.5)-- (-5,-1);
\draw [color=red,line width=1pt] (-2,0)-- (2,0);
\draw [color=red,dashed,line width=1pt] (-1.,0.5)-- (-3,-1);
\draw [line width=1pt] (-2,-1)-- (-4,0.5);
\draw [line width=2pt,red] (-0.5,0.25)-- (-0.5,-0.25);
\draw[line width=2pt,color=red,smooth,samples=100,domain=1:1.5] plot(\x,{0+sqrt((\x)-1)});
\draw[line width=2pt,color=red,smooth,samples=100,domain=1:1.5] plot(\x,{0-sqrt((\x)-1)});
\draw[line width=2pt,color=red,smooth,samples=100,domain=0.5:1] plot(\x,{0+sqrt(-((\x)-1))});
\draw[line width=2pt,color=red,smooth,samples=100,domain=0.5:1] plot(\x,{0-sqrt(-((\x)-1))});
\begin{scriptsize}
\draw [fill=red] (-8,-0.5) node {$\Gamma_1$};
\draw [fill=red] (-8.5,0.125) node {$1$};
\draw [fill=red] (-5.5,-0.5) node {$\Gamma_2$};
\draw [fill=red] (-5,0.2) node {$2$};
\draw [fill=red] (-4,0.2) node {$2$};
\draw [fill=red] (-3,0.2) node {$2$};
\draw [fill=red] (-1.5,0.5) node {$1$};
\draw [fill=red] (-7,0.5) node {$\Delta_0$};
\draw [fill=red] (-7,1.125) node {$8$};
\draw [fill=red] (2,-0.25) node {$E$};
\draw [fill=red] (2,0.25) node {$4$};
\end{scriptsize}
\end{tikzpicture}
\end{center}

In the figure, $\Gamma_i$ and $E$ denote their proper transforms.

We observe that the chain formed by the solid lines to the left of the dashed line is given by $[1,8,2^{3}]$ .
Continue the even resolution process until the branch locus becomes smooth in the neighborhood of $E$.  
Let $\widetilde{\Sigma}$ be the surface obtained from $\Sigma_{7}$ through the above procedure.
The branch divisor $B_{\widetilde{\Sigma}}$ on $\widetilde{\Sigma}$ is determined according to the even resolution process (Definition~\ref{defn--evenresol}).
As a result, a neighborhood of the total transform of $\Gamma_1 + \Gamma_2 + \Delta_0$ on $\widetilde{\Sigma}$ is illustrated in the following figure.

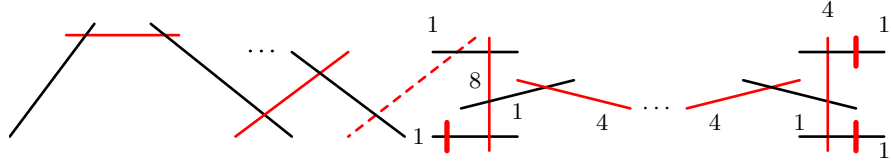
\begin{center}
\begin{tikzpicture}[line cap=round,line join=round,>=triangle 45,x=0.75cm,y=0.75cm]
\clip(-9,-2) rectangle (7,2);
\draw [color=red,line width=1pt] (-8,0.8)-- (-6,0.8);
\draw [line width=1pt] (-7.5,1)-- (-9,-1);
\draw [line width=1pt] (-4,-1)-- (-6.5,1);
\draw [color=red,line width=1pt] (-3,0.5)-- (-5,-1);
\draw [color=red,dashed,line width=1pt] (-0.75,0.75)-- (-3,-1);
\draw [line width=1pt] (-2,-1)-- (-4,0.5);
\draw [line width=1pt] (-1.5,0.5)-- (0,0.5);
\draw [line width=1pt] (-1.5,-1)-- (0,-1);
\draw [line width=1pt,red] (-0.5,0.75)-- (-0.5,-1.25);
\draw [line width=1pt] (-1,-0.5)-- (1,0);
\draw [line width=1pt,red] (0,0)-- (2,-0.5);
\draw [line width=1pt,red] (3,-0.5)-- (5,0);
\draw [line width=1pt] (4,0)-- (6,-0.5);
\draw [line width=1pt,red] (5.5,0.75)-- (5.5,-1.25);
\draw [line width=1pt] (5,0.5)-- (6.5,0.5);
\draw [line width=1pt] (5,-1)-- (6.5,-1);
\draw [line width=2pt,red] (6,0.25)-- (6,0.75);
\draw [line width=2pt,red] (6,-0.75)-- (6,-1.25);
\draw [line width=2pt,red] (-1.25,-0.75)-- (-1.25,-1.25);
\begin{scriptsize}
\draw [fill=red] (-4.5,0.5) node {$\cdots$};
\draw [fill=red] (2.5,-0.5) node {$\cdots$};
\draw [fill=red] (-0.75,0) node {$8$};
\draw [fill=red] (-1.5,1) node {$1$};
\draw [fill=red] (-1.75,-1) node {$1$};
\draw [fill=red] (0,-0.55) node {$1$};
\draw [fill=red] (1.5,-0.75) node {$4$};
\draw [fill=red] (-4.5,0.5) node {$\cdots$};
\draw [fill=red] (3.5,-0.75) node {$4$};
\draw [fill=red] (5,-0.75) node {$1$};
\draw [fill=red] (5.5,1.25) node {$4$};
\draw [fill=red] (6.5,1) node {$1$};
\draw [fill=red] (6.5,-1.25) node {$1$};
\end{scriptsize}
\end{tikzpicture}
\end{center}

In the figure above, the dual graph corresponding to the thin solid lines to the right of the dashed line is given as follows:
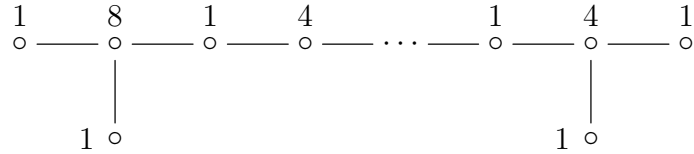
\[
\xygraph{
    \circ ([]!{-(0,-.3)} {1}) - [r]
    \circ ([]!{-(0,-.3)} {8}) (
        - [d] \circ ([]!{-(.3,0)} {1}),
        - [r] \circ ([]!{-(-0,-.3)} {1})
         - [r] \circ ([]!{-(0,-.3)} {4},
          - [r] \cdots ([]!{-(0,-.3)} {},
           - [r] \circ ([]!{-(0,-.3)} {1},
            - [r] \circ ([]!{-(0,-.3)} {4}) (
         - [d] \circ ([]!{-(.3,0)} {1}),
        - [r] \circ ([]!{-(0,-.3)} {1})
)}
\]
The length of the subchain $[8,1,4,\cdots,1,4]$ in the dual graph is $2k-1$.
We denote this entire dual graph by $(1,1,1,1)[8,1,4,\dots,1,4]$.
By abuse of notation, we also use $(1,1,1,1)[8,1,4,\dots,1,4]$ to refer to the corresponding configuration of curves on $\widetilde{\Sigma}$.

Now, take a double cover branched along the branch locus $B_{\widetilde{\Sigma}}$. 
In the double cover, the preimage of the chain $[1,8,2^3]$ coincides with $[2,4,4,1,4]$.  
Similarly, the preimage of $(1,1,1,1)[8,1,4,\dots,1,4]$ coincides with $(2,2,2,2)[4,2^{2k-2}]$.
Finally, contracting the chain $[2,4,4,1,4]$ and $(2,2,2,2)[4,2^{2k-2}]$, we obtain the desired Horikawa surface.

Conversely, it can also be verified that all surfaces belonging to $(3)$ in Theorem~\ref{thm:classification_non-standard_Horikawa} can be obtained through the above construction.

 \item[(4)] Let $B_{\Sigma}= \Delta_0+B_{\Sigma}' \in |4\Delta_{0}+ 16 \Gamma|$ be a reduced effective divisor on the Hirzebruch surface $\Sigma_{4}$ satisfying the following condition.
\begin{itemize}  
    \item [($\operatorname{SR}$,$\operatorname{SR}$)]
    In the neighborhood of two distinct fibers $\Gamma_1$, $\Gamma_2$, the branch divisor $B_{\Sigma}$ is described as one of those listed in Table~\ref{table_generate_[*SR]}, respectively.
    The pair $(\Sigma_4,B_{\Sigma})$ is mild on $\Sigma_4 \setminus \Gamma_1\cup \Gamma_2$. 
\item[($\operatorname{C}$,$\operatorname{C}$,$\operatorname{SR}$)]
     In the neighborhood of two distinct fibers $\Gamma_1$, $\Gamma_2$, $B_{\Sigma}$ is described as one of those listed in Table~\ref{table_generate_[CT]}, respectively.
    In the neighborhood of another fiber $\Gamma_3$, the branch divisor $B_{\Sigma}$ is described as one of those listed in Table~\ref{table_generate_[*SR]}.
    The pair $(\Sigma_4,B_{\Sigma})$ is mild on $\Sigma_4 \setminus \Gamma_1\cup \Gamma_2 \cup \Gamma_3$. 
\item[($\operatorname{C}$,$\operatorname{C}$,$\operatorname{C}$,$\operatorname{C}$)]
     In the neighborhood of four distinct fibers $\Gamma_i$ ($i=1,\cdots,4$), $B_{\Sigma}$ is described as one of those listed in Table~\ref{table_generate_[CT]}, respectively.
    The pair $(\Sigma_4,B_{\Sigma})$ is mild on $\Sigma_4 \setminus \bigcup_{i=1}^{4} \Gamma_i$. 
\end{itemize}
Note that the existence of such $B_{\Sigma}$ is guaranteed by Remark~\ref{rem--existence-all--q-gor--smoothable}.

Here, we consider the case where $(\operatorname{C}, \operatorname{C}, \operatorname{SR})$.
We assume that the branch locus $B_{\Sigma}$, in a neighborhood of each of $\Gamma_1$ and $\Gamma_2$, is of the form $(\mathrm{C3})$ as listed in Table~\ref{table_generate_[CT]}.
Furthermore, in the neighborhood of $\Gamma_3$, we assume that $B_{\Sigma}$ is of the form $(\mathrm{SR}3)$ in Table~\ref{table_generate_[*SR]}.
The branch locus $B_{\Sigma}$ has a $k$-fold node on $\Gamma_3$.

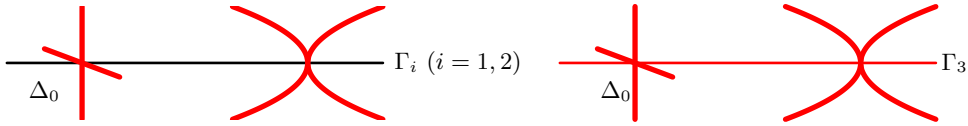
\begin{table}[H]
\centering
  \begin{tabular}{cc}
    \begin{minipage}{70mm}
\begin{tikzpicture}[line cap=round,line join=round,>=triangle 45,x=1cm,y=0.75cm]
\clip(-3,-1) rectangle (4.5,1);
\draw [line width=1pt] (-2.5,0)--(2.5,0);
\draw [line width=2pt, red] (-1.5,-1)-- (-1.5,1);
\draw [line width=2pt, red] (-2,0.25)-- (-1,-0.25);
\draw[line width=2pt,color=red,smooth,samples=100,domain=1.5:2.5] plot(\x,{sqrt((\x)-1.5)});
\draw[line width=2pt,color=red,smooth,samples=100,domain=1.5:2.5] plot(\x,{0-sqrt((\x)-1.5)});
\draw[line width=2pt,color=red,smooth,samples=100,domain=0.5:1.5] plot(\x,{sqrt(-((\x)-1.5))});
\draw[line width=2pt,color=red,smooth,samples=100,domain=0.5:1.5] plot(\x,{0-sqrt(-((\x)-1.5))});
\begin{scriptsize}
\draw [color=black] (3.5,0) node {$\Gamma_i$ $(i=1,2)$};
\draw [color=black] (-2,-0.5) node {$\Delta_0$};
\end{scriptsize}
\end{tikzpicture}
     \end{minipage} &
    \begin{minipage}{70mm}
\begin{tikzpicture}[line cap=round,line join=round,>=triangle 45,x=1cm,y=0.75cm]
\clip(-3,-1.25) rectangle (3,1.25);
\draw [line width=1pt,red] (-2.5,0)--(2.5,0);
\draw [line width=2pt,color=red] (-1.5,-1)-- (-1.5,1);
\draw [line width=2pt,color=red] (-2,0.25)-- (-1,-0.25);
\draw[line width=2pt,color=red,smooth,samples=100,domain=1.5:2.5] plot(\x,{sqrt((\x)-1.5)});
\draw[line width=2pt,color=red,smooth,samples=100,domain=1.5:2.5] plot(\x,{0-sqrt((\x)-1.5)});
\draw[line width=2pt,color=red,smooth,samples=100,domain=0.5:1.5] plot(\x,{sqrt(-((\x)-1.5))});
\draw[line width=2pt,color=red,smooth,samples=100,domain=0.5:1.5] plot(\x,{0-sqrt(-((\x)-1.5))});
\begin{scriptsize}
\draw [color=black] (2.75,0) node {$\Gamma_3$};
\draw [color=black] (-1.75,-0.5) node {$\Delta_0$};
\end{scriptsize}
\end{tikzpicture}
    \end{minipage} \\ 
\end{tabular}
\end{table}

First, we describe the procedure in the neighborhood of $\Gamma_1$ and $\Gamma_2$.
Blow up at the intersection of the section $\Delta_0$ and the fibers $\Gamma_1$ and $\Gamma_2$, respectively.

Next, we describe the procedure in the neighborhood of $\Gamma_3$.
We continue the process of even resolution (Definition~\ref{defn--evenresol}) until the germ of the branch loci in the neighborhood of the fiber $\Gamma_3$ becomes smooth.
Let $\widetilde{\Sigma}$ be the surface obtained from $\Sigma_4$ by the above procedures.
The branch divisor $B_{\widetilde{\Sigma}}$ on $\widetilde{\Sigma}$ is determined according to the even resolution process (Definition~\ref{defn--evenresol}).

Based on the above preparation, we now describe two possible constructions, depending on the choice of curves to be contracted among the total transform of $\Gamma_1 + \Gamma_2 + \Gamma_3 + \Delta_0$ on $\widetilde{\Sigma}$.

We begin by describing the first construction.
The following figure illustrates a neighborhood of the total transform of $\Gamma_1 + \Gamma_2+\Gamma_3 + \Delta_0$ on $\widetilde{\Sigma}$.

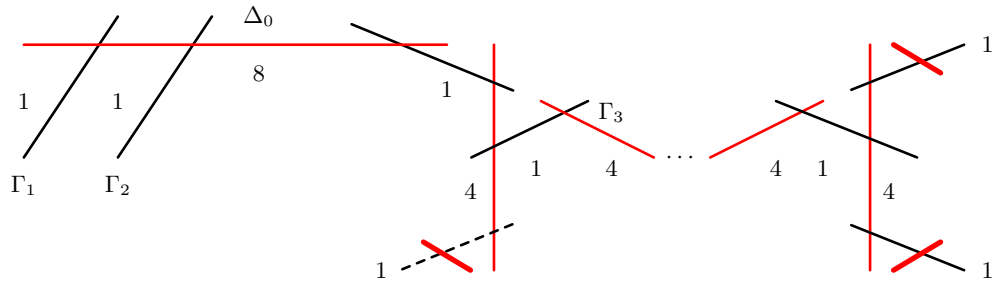
\begin{center}
\begin{tikzpicture}[line cap=round,line join=round,>=triangle 45,x=1.25cm,y=1.5cm]
\clip(-7,-1.25) rectangle (5,1.5);
\draw [line width=1pt,red] (-1,1)-- (-1,-1);
\draw [line width=1pt] (-0.796405425080071,0.5940678833808372)-- (-2.514305105837884,1.181896073206929);
\draw [line width=1pt,dashed] (-0.8037318293379739,-0.5928096063994598)-- (-2,-1);
\draw [line width=1pt] (-1.2408739500595178,0)-- (0,0.5);
\draw [line width=1pt,red] (3,1)-- (3,-1);
\draw [line width=1pt] (2.8,0.6)-- (4,1);
\draw [line width=1pt] (2.8,-0.6)-- (4,-1);
\draw [line width=1pt,red] (-0.5,0.5)-- (0.7,0);
\draw [line width=1pt,red] (1.3,0)-- (2.5,0.5);
\draw [line width=1pt] (2,0.5)-- (3.5,0);
\draw [line width=1pt] (-5,1.25)-- (-6,0);
\draw [line width=1pt] (-4,1.25)-- (-5,0);
\draw [line width=1pt, red] (-6,1)-- (-1.5,1);
\draw [line width=2pt,red] (3.25,1)-- (3.75,0.75);
\draw [line width=2pt,red] (3.25,-1)-- (3.75,-0.75);
\draw [line width=2pt,red] (-1.25,-1)-- (-1.75,-0.75);
\begin{scriptsize}
\draw [fill=black] (1,0) node {$\cdots$};
\draw [fill=black] (-6,-0.25) node {$\Gamma_1$};
\draw [fill=black] (-5,-0.25) node {$\Gamma_2$};
\draw [fill=black] (0.25,0.4) node {$\Gamma_3$};
\draw [fill=black] (-3.5,1.25) node {$\Delta_0$};
\draw [fill=black] (-3.5,0.75) node {$8$};
\draw [fill=black] (-6,0.5) node {$1$};
\draw [fill=black] (-5,0.5) node {$1$};
\draw [fill=black] (-1.5,0.6) node {$1$};
\draw [fill=black] (-1.25,-0.3) node {$4$};
\draw [fill=black] (-2.2,-1) node {$1$};

\draw [fill=black] (-0.55,-0.1) node {$1$};
\draw [fill=black] (0.25,-0.1) node {$4$};
\draw [fill=black] (2,-0.1) node {$4$};
\draw [fill=black] (2.5,-0.1) node {$1$};
\draw [fill=black] (3.2,-0.3) node {$4$};
\draw [fill=black] (4.25,1) node {$1$};
\draw [fill=black] (4.25,-1) node {$1$};
\end{scriptsize}
\end{tikzpicture}
\end{center}

\noindent
In the figure, $\Gamma_i$ and $\Delta_0$ denote their proper transforms.
The dual graph corresponding to the thin solid lines is as follows:
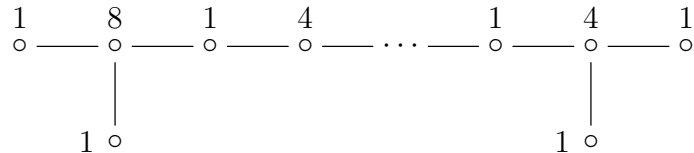
\[
\xygraph{
    \circ ([]!{-(0,-.3)} {1}) - [r]
    \circ ([]!{-(0,-.3)} {8}) (
        - [d] \circ ([]!{-(.3,0)} {1}),
        - [r] \circ ([]!{-(-0,-.3)} {1})
         - [r] \circ ([]!{-(0,-.3)} {4},
          - [r] \cdots ([]!{-(0,-.3)} {},
           - [r] \circ ([]!{-(0,-.3)} {1},
            - [r] \circ ([]!{-(0,-.3)} {4}) (
         - [d] \circ ([]!{-(.3,0)} {1}),
        - [r] \circ ([]!{-(0,-.3)} {1})
)}
\]
The length of subchain $[8,1,4,\cdots,1,4]$ in the dual graph above is $2k+5$.
We denote by $(1,1,1,1)[8,1,4,\dots,1,4]$ the dual graph above.
By abuse of notation, we also regard $(1,1,1,1)[8,1,4,\dots,1,4]$ as denoting the corresponding configuration of curves on $\widetilde{\Sigma}$.
We take a double cover of $\widetilde{\Sigma}$ branched along $B_{\widetilde{\Sigma}}$.
In the double cover, the preimages of $(1,1,1,1)[8,1,4,\dots,1,4]$ is $(2,2,2,2)[4,2^{2k+4}]$.
By contracting $(2,2,2,2)[4,2^{2k+4}]$, 
we obtain a non-standard Horikawa surface belonging to $(4)$ in Theorem~\ref{thm:classification_non-standard_Horikawa}

We now proceed to describe the second construction.
We denote by $(1,1,1,1)[8,1,4]$ the configuration of curves on $\widetilde{\Sigma}$ corresponding to the thin solid lines, as illustrated in the figure below:

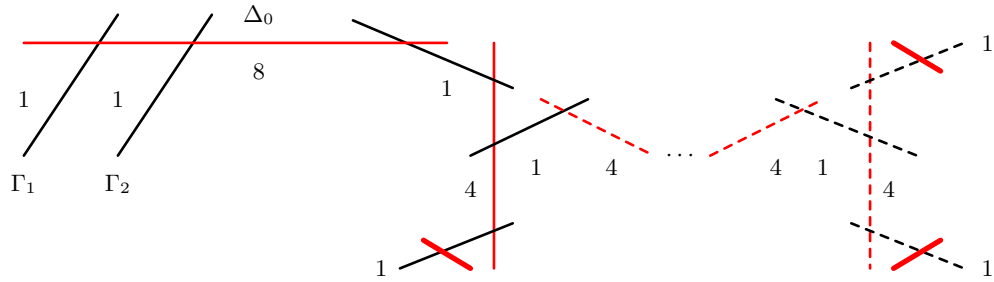
\begin{center}
\begin{tikzpicture}[line cap=round,line join=round,>=triangle 45,x=1.25cm,y=1.5cm]
\clip(-7,-1.25) rectangle (5,1.5);
\draw [line width=1pt,red] (-1,1)-- (-1,-1);
\draw [line width=1pt] (-0.8,0.6)-- (-2.5,1.2);
\draw [line width=1pt] (-0.8,-0.6)-- (-2,-1);
\draw [line width=1pt] (-1.25,0)-- (0,0.5);
\draw [line width=1pt,red,dashed] (3,1)-- (3,-1);
\draw [line width=1pt,dashed] (2.8,0.6)-- (4,1);
\draw [line width=1pt,dashed] (2.8,-0.6)-- (4,-1);
\draw [line width=1pt,red,dashed] (-0.5,0.5)-- (0.7,0);
\draw [line width=1pt,red,dashed] (1.3,0)-- (2.5,0.5);
\draw [line width=1pt,dashed] (2,0.5)-- (3.5,0);
\draw [line width=1pt] (-5,1.25)-- (-6,0);
\draw [line width=1pt] (-4,1.25)-- (-5,0);
\draw [line width=1pt, red] (-6,1)-- (-1.5,1);
\draw [line width=2pt,red] (3.25,1)-- (3.75,0.75);
\draw [line width=2pt,red] (3.25,-1)-- (3.75,-0.75);
\draw [line width=2pt,red] (-1.25,-1)-- (-1.75,-0.75);
\begin{scriptsize}
\draw [fill=black] (1,0) node {$\cdots$};
\draw [fill=black] (-6,-0.25) node {$\Gamma_1$};
\draw [fill=black] (-5,-0.25) node {$\Gamma_2$};
\draw [fill=black] (-3.5,1.25) node {$\Delta_0$};
\draw [fill=black] (-3.5,0.75) node {$8$};
\draw [fill=black] (-6,0.5) node {$1$};
\draw [fill=black] (-5,0.5) node {$1$};
\draw [fill=black] (-1.5,0.6) node {$1$};
\draw [fill=black] (-1.25,-0.3) node {$4$};
\draw [fill=black] (-2.2,-1) node {$1$};

\draw [fill=black] (-0.55,-0.1) node {$1$};
\draw [fill=black] (0.25,-0.1) node {$4$};
\draw [fill=black] (2,-0.1) node {$4$};
\draw [fill=black] (2.5,-0.1) node {$1$};
\draw [fill=black] (3.2,-0.3) node {$4$};
\draw [fill=black] (4.25,1) node {$1$};
\draw [fill=black] (4.25,-1) node {$1$};
\end{scriptsize}
\end{tikzpicture}
\end{center}

\noindent
Let $D$ be the reduced divisor on $\widetilde{\Sigma}$ defined as the sum of the curves corresponding to the dashed lines in the above figure that are disjoint from the $(1,1,1,1)[8,1,4]$.
We take a double cover of $\widetilde{\Sigma}$ branched along $B_{\widetilde{\Sigma}}$.
In the double cover, the preimages of $(1,1,1,1)[8,1,4]$ is $(2,2,2,2)[4,2,2]$.
The preimages of $D$ is as follows:
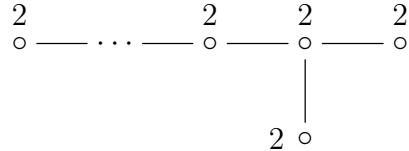
\[
\xygraph{
    \circ ([]!{-(0,-.3)} {2})
          - [r] \cdots ([]!{-(0,-.3)} {},
           - [r] \circ ([]!{-(0,-.3)} {2},
            - [r] \circ ([]!{-(0,-.3)} {2}) (
         - [d] \circ ([]!{-(.3,0)} {2}),
        - [r] \circ ([]!{-(0,-.3)} {2})
)}
\]
By contracting $(2,2,2,2)[4,2,2]$ and the preimages of $D$, 
we obtain a non-standard Horikawa surface belonging to $(4)$ in Theorem~\ref{thm:classification_non-standard_Horikawa}.

Conversely, it can also be verified that all surfaces belonging to $(4)$ in Theorem~\ref{thm:classification_non-standard_Horikawa} can be obtained through the above construction by choosing a certain branch $B_{\Sigma}$.

\item[(5)]
The construction of non-standard Horikawa surfaces corresponding to $(5)$ in Theorem~\ref{thm:classification_non-standard_Horikawa} differs depending on whether $\chi = 4$ or $\chi \ge 5$.  
We treat these cases separately.

In the case $\chi = 4$, let $\Sigma_4$ be the Hirzebruch surface, and let $B_{\Sigma} \in |4\Delta_0 + 16\Gamma|$ be a reduced divisor such that the pair $(\Sigma_4, B_{\Sigma})$ is mild.  
Take a double cover of $\Sigma_4$ branched along $B_\Sigma$.  
In the double cover, the preimage of the $(-4)$-curve consists of two disjoint $(-4)$-curves.  
By contracting these two $(-4)$-curves, we obtain the desired non-standard Horikawa surface.

Conversely, any non-standard Horikawa surface with $\chi = 4$ appearing in $(5)$ of Theorem~\ref{thm:classification_non-standard_Horikawa} can be obtained via this construction.
We refer to 
this surface as a \emph{Horikawa surface of (infinite) Lee-Park type} with $p_g=3$.

In the case $\chi \ge 5$, we consider the Hirzebruch surface $\Sigma_\chi$.
Let $B_{\Sigma}  \in |4\Delta_{0}+ 4\chi \Gamma|$ be a reduced effective divisor on the Hirzebruch surface $\Sigma_{\chi}$ satisfying the following condition:   
\begin{itemize}  
    \item $\Delta_0 \not\subset B_{\Sigma}$
    \item In the neighborhood of a certain fiber $\Gamma$, the branch divisor $B_{\Sigma}$ is described as one of those listed in Table~\ref{table_generate_[LP]}.
    \item The pair $(\Sigma_{\chi},B_{\Sigma})$ is mild on $\Sigma_{\chi}\setminus \Gamma $.  
\end{itemize}
Since the other cases can be treated in a similar way, we describe the construction in the case where the branch is locally described as $(\mathrm{LP}5)$ near the fiber $\Gamma$.
Note that the existence of such $B_{\Sigma}$ is guaranteed by Remark~\ref{rem--existence-all--q-gor--smoothable}.

We describe the procedure in the neighborhood of $\Gamma$.
Continue the process of even resolution until the germ of the branch locus $B_{\Sigma}$ in the neighborhood of the fiber $\Gamma$ becomes smooth.
Let $\widetilde{\Sigma}$ be the surface obtained from $\Sigma_{\chi}$ by the above procedures.
The branch divisor $B_{\widetilde{\Sigma}}$ on $\widetilde{\Sigma}$ is determined according to the even resolution process (Definition~\ref{defn--evenresol}).
A neighborhood of the total transform of $\Gamma + \Delta_0$ on $\widetilde{\Sigma}$ is illustrated in the following figure.

\begin{figure}[H]
    \centering
\begin{tikzpicture}[line cap=round,line join=round,>=triangle 45,x=0.4cm,y=0.25cm]
\clip(-3,-10) rectangle (20,4);
\draw [line width=1pt] (-2,2)-- (2,2);
\draw [line width=1pt] (0,3)-- (-2,-2);
\draw [dashed,line width=1pt] (-2,-4)-- (0,-9);
\draw [line width=1pt] (-1.35,1)-- (-1.3884710932081723,-7);
\draw [line width=1pt] (-2,-3)-- (2,-1);
\draw [line width=1pt] (4,-1)-- (8,-3);
\draw [dashed,line width=1pt] (6,-3)-- (10,-1);
\draw [dashed,line width=1pt] (12,-1)-- (16,-3);
\draw [dashed,line width=1pt] (14,-3)-- (18,-1);
\draw [red,line width=2pt] (16,-1)-- (18,-2);
\draw [red,line width=2pt] (15.5,-1.5)-- (17.5,-2.5);
\begin{scriptsize}
\draw  (-2,2.65) node {$\chi$};
\draw  (2,2.65) node {$\Delta_0$};
\draw  (1,-0.75) node {$2$};
\draw  (5,-0.75) node {$2$};
\draw  (8,-0.75) node {$2$};
\draw  (14,-0.75) node {$2$};
\draw  (-0.5,0.75) node {$2$};
\draw  (-0.5,-4) node {$2$};
\draw  (0.5,-8) node {$2$};
\draw  (18,-0.5) node {$1$};
\draw  (3,-1) node {$\cdots$};
\draw  (11,-1) node {$\cdots$};
\end{scriptsize}
\end{tikzpicture}
\caption{}
\label{fig:const_LPG_chain}
\end{figure}

Except for the case $\chi = 7$, the total transform of $\Delta_0$ and $\Gamma$ contains a unique chain of the form $[\chi, 2^{\chi - 4}]$.
The solid lines in Figure~\ref{fig:const_LPG_chain} denote this chain.
In the case $\chi = 7$ only, the chain $[7,2^3]$ appears in an alternative configuration, as shown below:

\begin{figure}[H]
    \centering
    \begin{tikzpicture}[line cap=round,line join=round,>=triangle 45,x=0.4cm,y=0.25cm]
\clip(-3,-10) rectangle (20,4);
\draw [line width=1pt] (-2,2)-- (2,2);
\draw [line width=1pt] (0,3)-- (-2,-2);
\draw [line width=1pt] (-2,-4)-- (0,-9);
\draw [line width=1pt] (-1.35,1)-- (-1.3884710932081723,-7);
\draw [dashed,line width=1pt] (-2,-3)-- (2,-1);
\draw [dashed,line width=1pt] (4,-1)-- (8,-3);
\draw [dashed,line width=1pt] (6,-3)-- (10,-1);
\draw [dashed,line width=1pt] (12,-1)-- (16,-3);
\draw [dashed,line width=1pt] (14,-3)-- (18,-1);
\draw [red,line width=2pt] (16,-1)-- (18,-2);
\draw [red,line width=2pt] (15.5,-1.5)-- (17.5,-2.5);
\begin{scriptsize}
\draw  (-2,2.65) node {$7$};
\draw  (2,2.65) node {$\Delta_0$};
\draw (-0.5,0.75) node {$2$};
\draw (-0.5,-4) node {$2$};
\draw (0.5,-8) node {$2$};
\draw  (18,-0.5) node {$1$};
\draw  (3,-1) node {$\cdots$};
\draw  (11,-1) node {$\cdots$};
\end{scriptsize}
\end{tikzpicture}
    \caption{}
    \label{fig:const_LP_infty_chain}
\end{figure}
Thus, in the case $\chi = 7$, the total transform of $\Delta_0$ and $\Gamma$ contains two distinct configurations of chains of the form $[7,2^3]$.
We note that this occurs only when $\chi=7$ and the chosen branch is of type $(\mathrm{LP}5)$ or $(\mathrm{LP}6)$.

Take a double cover of $\widetilde{\Sigma}$ branched along $B_{\widetilde{\Sigma}}$.

Assume that the chain $[\chi,2^{\chi-4}]$ corresponding to the Figure~\ref{fig:const_LPG_chain}.
Let $D$ be the chain consisting of all $(-1)$-curves and $(-2)$-curves contained in the total transform of $\Delta_0 + \Gamma$ that are disjoint from the chain $[\chi, 2^{\chi - 4}]$.
The chain $[\chi,2^{\chi-4}]$ pulls back under the double cover to two disjoint strings $[\chi,2^{\chi-4}]$.
The chain $D$ pulls back under the double cover to a chain consisting of $(-2)$-curves.
By contracting these two chains $[\chi,2^{\chi-4}]$ and the chain consisting of $(-2)$-curves, we obtain a non-standard Horikawa surface.

Assume that the chain $[7,2^3]$ corresponding to the Figure~\ref{fig:const_LP_infty_chain}.
Let $D'$ be the chain consisting of all $(-1)$-curves and $(-2)$-curves contained in the total transform of $\Delta_0 + \Gamma$ that are disjoint from the chain $[7,2^3]$.
The chain $[7,2^3]$ pulls back under the double cover to two disjoint strings $[7,2^3]$.
The chain $D'$ pulls back under the double cover to a chain consisting of $(-2)$-curves.
By contracting these two chains $[7,2^3]$ and the chain consisting of $(-2)$-curves, we obtain a non-standard Horikawa surface.

Conversely, the non-standard Horikawa surfaces appearing in $(5)$ of Theorem~\ref{thm:classification_non-standard_Horikawa} can be constructed by choosing a certain branch $B_{\Sigma}$ from one of the cases listed in Table~\ref{table_generate_[LP]}.  

We divide non-standard Horikawa surfaces with $\chi \ge 5$ corresponding to (5) in Theorem~\ref{thm:classification_non-standard_Horikawa} into several types according to the list of branches $B_{\Sigma}$ in Table~\ref{table_generate_[LP]}.

\begin{itemize} 

\item
$\chi \geq 5$. 
The surface obtained by choosing $(\mathrm{LP}1)$ $(\mathrm{LP}2)$, $(\mathrm{LP}3)$ or $(\mathrm{LP}4)$ in Table~\ref{table_generate_[LP]} as the branch locus $B_{\Sigma}$ is called a \emph{Horikawa surface of general Lee-Park type} $I$.
 \item
$\chi = 5$. The surface obtained by choosing $(\mathrm{LP}5)$, $(\mathrm{LP}6)$, $(\mathrm{LP}7)$ or $(\mathrm{LP}8)$ in Table~\ref{table_generate_[LP]} as the branch locus $B_{\Sigma}$ is called a \emph{Horikawa surface of special Lee-Park type $I^*$} .
Such a surface has a Du Val double cone singularity with respect to the covering involution.
\item
$\chi \geq 6$. The surface obtained by choosing $(\mathrm{LP}5)$ or $(\mathrm{LP}6)$ in Table~\ref{table_generate_[LP]} as the branch locus $B_{\Sigma}$ (in $\chi=7$, a surface by contracting the two chains $[7,2^3]$ corresponding to the solid lines in Figure~\ref{fig:const_LPG_chain}) is called a \emph{Horikawa surface of general Lee-Park type $I^*$}. 

When $\chi = 7$, the surface obtained by contracting the two chains $[7, 2^3]$ corresponding to the solid lines in Figure~\ref{fig:const_LP_infty_chain} is called a \emph{Horikawa surface of infinite Lee-Park type}.

\item
$\chi = 5,6,7$. The surface obtained by choosing $(\mathrm{LP}9)$ in Table~\ref{table_generate_[LP]} as the branch locus $B_{\Sigma}$ is called a \emph{Horikawa surface of special Lee-Park type $III^*$}. 
Such a surface has a Du Val double cone singularity with respect to the covering involution.
\item
$\chi = 5,6$. The surface obtained by choosing $(\mathrm{LP}10)$ in Table~\ref{table_generate_[LP]} as the branch locus $B_{\Sigma}$ is called a \emph{Horikawa surface of special Lee-Park type $IV^*$}.  
Such a surface has a Du Val double cone singularity with respect to the covering involution.
\end{itemize}
We add a brief remark on the terminology \emph{special Lee-Park type}. 
For non-standard Horikawa surfaces, the covering involution in this construction is a good involution.  
Among the non-standard Horikawa surfaces corresponding to cases (5), (6), and (7), those that have a double cone singularity with respect to the good involution are referred to as being of special Lee-Park type.

\item[(6),(7)]
 Let $B_{\Sigma}  \in |4\Delta_{0}+ 4\chi \Gamma|$ be a reduced effective divisor on the Hirzebruch surface $\Sigma_{\chi}$ satisfying the following condition:   
\begin{itemize}  
    \item $\Delta_0 \not\subset B_{\Sigma}$
    \item In the neighborhood of a certain fiber $\Gamma$, the branch divisor $B_{\Sigma}$ is described as one of those listed in Table~\ref{table_generate_[LPD]}.
    \item The pair $(\Sigma_{\chi},B_{\Sigma})$ is mild on $\Sigma_{\chi}\setminus \Gamma $.  
\end{itemize}
Since the other cases can be treated in a similar way, we describe the construction in the case where the branch is locally described as $(\mathrm{LP}11)$ ($k_1=k_2=1$) in Table~\ref{table_generate_[LPD]}.
Note that the existence of such $B_{\Sigma}$ is guaranteed by Remark~\ref{rem--existence-all--q-gor--smoothable}.

We describe the procedure in the neighborhood of $\Gamma$.
Blow up the surface at the ordinary quadruple point of $B_{\Sigma}$ on $\Gamma$.
Let $\Gamma'$ be the proper transform of $\Gamma$, and let $E$ be the exceptional curve.
Next, blow up the intersection point of $\Gamma'$ and $E$.
Then, perform $(\chi-5)$-times further successive blow-ups at the intersection point of the total transform of $\Gamma'$ and the proper transform of $E$ at each step.
Let $\widetilde{\Sigma}$ be the surface obtained from $\Sigma_{\chi}$ by the above procedure.
The branch divisor $B_{\widetilde{\Sigma}}$ on $\widetilde{\Sigma}$ is determined by the sense of even resolution (Definition~\ref{defn--evenresol}).

As a result, a neighborhood of the total transform of $\Gamma + \Delta_0$ on $\widetilde{\Sigma}$ is illustrated in the following figure.
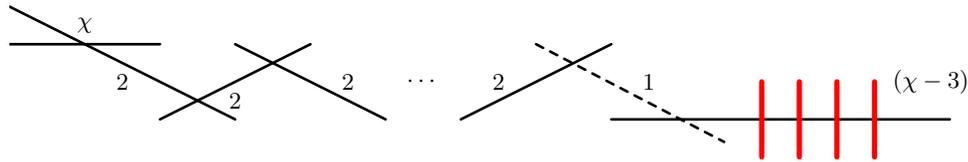
\begin{center}
\begin{tikzpicture}[line cap=round,line join=round,>=triangle 45,x=1cm,y=1cm]
\clip(-3,-1) rectangle (11,2);
\draw [line width=1pt] (-3,1.5)-- (0,0);
\draw [line width=1pt] (-1,0)-- (1,1);
\draw [line width=1pt] (0,1)-- (2,0);
\draw [line width=1pt] (-3,1)-- (-1,1);
\draw [line width=1pt] (3,0)-- (5,1);
\draw [line width=1pt,dashed] (4,1)-- (6.5,-0.3);
\draw [line width=1pt] (5,0)-- (9.5,0);
\draw [line width=2pt,red] (7,0.5)-- (7,-0.5);
\draw [line width=2pt,red] (7.5,0.5)-- (7.5,-0.5);
\draw [line width=2pt,red] (8,0.5)-- (8,-0.5);
\draw [line width=2pt,red] (8.5,0.5)-- (8.5,-0.5);
\begin{scriptsize}
\draw [black] (-2,1.25) node {$\chi$};
\draw [black] (-1.5,0.5) node {$2$};
\draw [black] (0,0.25) node {$2$};
\draw [black] (1.5,0.5) node {$2$};
\draw [black] (2.5,0.5) node {$\cdots$};
\draw [black] (3.5,0.5) node {$2$};
\draw [black] (5.5,0.5) node {$1$};
\draw [black] (9.25,0.5) node {$(\chi -3)$};
\end{scriptsize}
\end{tikzpicture}
\end{center}
Let $[\chi, 2^{\chi-4}]$ and $[\chi - 3]$ denote the chains on $\widetilde{\Sigma}$ represented by the solid lines in the figure above.
Take the double cover of $\widetilde{\Sigma}$ branched along $B_{\widetilde{\Sigma}}$.
In the double cover, the preimage of the chain $[\chi, 2^{\chi-4}]$ consists of two copies of the same chain.
The preimage of the chain $[\chi - 3]$ is an elliptic curve with self-intersection number $-(2\chi - 6)$.
By contracting the two chains $[\chi,2^{\chi-4}]$ and 
the elliptic curve, we obtained the desired Horikawa surface.

Conversely, it can also be verified that all surfaces belonging to $(6)$, $(7)$ in Theorem~\ref{thm:classification_non-standard_Horikawa} can be obtained through the above construction.

The surface obtained by choosing $(\mathrm{LP}11)$, $(\mathrm{LP}12)$ or $(\mathrm{LP}13)$ in Table~\ref{table_generate_[LPD]} as the branch locus $B$ is called a \emph{Horikawa surface of special Lee-Park type $I$}.
Such a surface has an elliptic double cone singularity with respect to the covering involution.

\begin{note}
We refer collectively to the special Lee-Park types $\mathrm{I}$, $\mathrm{I}^*$, $\mathrm{III}^*$, and $\mathrm{IV}^*$ introduced above as the \emph{special Lee-Park type}.
\end{note}

\item[(8)]
Let $C$ be an elliptic curve,
and let $\mathcal{M}$ be a line bundle on $C$ of degree $8$.
Define $\widehat{W}$ as the projective bundle $\widehat{W}:=\mathbb{P}_{C}(\O_{C}\oplus \mathcal{M}^{-1})$.
The Picard group of $\widehat{W}$ is generated by the negative section $\Delta_0$ and the pullback of the Picard group of $C$ via the projection $p\colon \widehat{W}\to C$.
The linear system $|\Delta_0+p^{*}\mathcal{M}|$ defines a contraction morphism $q \colon \widehat{W}\to W\subset \mathbb{P}^{8}$ onto an elliptic cone $W\subset \mathbb{P}^{8}$ of degree $8$.
Next, we choose a divisor $\widehat{B}'\in |3\Delta_{0}+3p^{*}\mathcal{M}|$ that is horizontal over $C$ and has only mild singularities as described in Proposition~\ref{lem:classify_mild1}.
Note that $\widehat{B}'$ and $\Delta_0$ are disjoint.
Define $\widehat{B}$ as $\widehat{B}:=\widehat{B}'+\Delta_0$.
We also take a line bundle $\widehat{L}$ on $\widehat{W}$ such that $\widehat{B}\in |2\widehat{L}|$.
Note that giving $\widehat{L}$ is equivalent to giving a line bundle $\mathcal{L}$ on $C$ satisfying $\mathcal{M}\cong \mathcal{L}^{\otimes 2}$ by a relation $\widehat{L}\cong \O_{\widehat{W}}(2\Delta_0+p^{*}(\mathcal{M}+\mathcal{L}))$,
which always exists.
Now, consider the double cover 
$\widehat{f}\colon \widehat{X}:=\mathbf{Spec}_{\widehat{W}}(\mathcal{O}_{\widehat{W}}\oplus \widehat{L}^{-1})\to \widehat{W}$
branched along $\widehat{B}$, and take a Stein factorization $\widehat{X}\to X\to W$ for the composition $q \circ \widehat{f}$.
The birational morphism $\widehat{X}\to X$ contracts the smooth elliptic curve $\widehat{f}^{-1}(\Delta_{0})$ to a simple elliptic singularity of degree $4$.
It follows easily that $X$ is Gorenstein and stable, satisfying $K_X^{2}=p_g(X)=4$.
Moreover, the elliptic pencil $X\leftarrow \widehat{X}\to \widehat{W}\to C$ is induced by $|K_X|$.
Thus, $X$ is a stable normal Gorentein Horikawa surface with a canonical pencil.
It is easy to check that the double cover $X\to W$ is induced by the bicanonical system $|2K_X|$.
\end{itemize}

We provide a list of the singularities of the branch $B_{\Sigma}$ used in the above construction.  
Each entry in the list consists of a figure of the branch on a germ of a fiber $\Gamma$ and a description of the singularities of $B_{\Sigma}$ on the corresponding fiber germ.

Assume that the horizontal part of $B_{\Sigma}$ decomposes on a fiber germ $\Gamma$ as a sum $B_{\Sigma,\mathrm{hor}}=\sum_i B_i$ of irreducible components.  
Note that each $B_i$ is an irreducible local analytic curve on the fiber germ $\Gamma$.
Following the above construction, we divide the table into four parts, Table~\ref{table_generate_[CT]}--\ref{table_generate_[LPD]}.
Although it is a slight abuse of notation, we use expressions such as $(B_i \cdot \Gamma)$ to denote the sum of local intersection numbers on the fiber in the Tables.

\newpage

\begin{table}[H]
\caption{}
\label{table_generate_[CT]}
  \begin{tabular}{|c|c|c|} 
    \hline & figure & local analytic branches of $B_{\Sigma}$ \\ \hline
 $(\mathrm{C1})$   &
    \begin{minipage}{70mm}
\begin{center}
\begin{tikzpicture}[line cap=round,line join=round,>=triangle 45,x=1cm,y=0.75cm]
\clip(-3,-1.25) rectangle (3,1.25);
\draw [line width=2pt, red] (-1.5,-1)-- (-1.5,1);
\draw [line width=2pt, red] (1.5,-1)-- (1.5,1);
\draw [line width=1pt] (2.5,0)-- (-2.5,0);
\draw[line width=2pt,color=black,smooth,samples=100,domain=-2.5:-0.5,red] plot(\x,{(\x + 1.5)^2});
\begin{scriptsize}
\draw [color=black] (-2.75,0) node {$\Gamma$};
\draw [color=black] (-1.75,-0.5) node {$\Delta_0$};
\draw [color=black] (-0.4,0.5) node {$B_2$};
\draw [color=black] (1.75,-0.25) node {$B_3$};
\end{scriptsize}
\end{tikzpicture}
\end{center}
    \end{minipage} &
    \begin{minipage}{70mm}
   $B_{\Sigma}=B_1+B_2+B_3$,\\
   $B_1 = \Delta_0$,
   $(B_2\cdot \Gamma)=2$,
   $(B_3\cdot \Gamma)=1$.
 \end{minipage} \\ \hline
 $(\mathrm{C2})$ &
\begin{minipage}{70mm}
\begin{center}
    \begin{tikzpicture}[line cap=round,line join=round,>=triangle 45,x=1cm,y=0.75cm]
\clip(-3,-1.25) rectangle (3,1.25);
\draw [line width=1pt] (2.5,0)-- (-2.5,0);
\draw [line width=2pt,color=red] (-1.5,-1)-- (-1.5,1);
\draw[line width=2pt,color=black,smooth,samples=100,domain=-2.5:-0.5,red] plot(\x,{(\x + 1.5)^2});
\begin{scriptsize}
\draw [color=black] (-2.75,0) node {$\Gamma$};
\draw [color=black] (-1.75,-0.5) node {$\Delta_0$};
\draw [color=black] (-0.4,0.5) node {$B_2$};
\end{scriptsize}
\end{tikzpicture}
\end{center}
\end{minipage} &
\begin{minipage}{70mm}
  $B_{\Sigma}=B_1+B_2$,\\
   $B_1 = \Delta_0$,
   $(B_2\cdot \Gamma)=3$.
\end{minipage} \\ \hline
$(\mathrm{C3})$ & 
\begin{minipage}{70mm}
\begin{center}
\begin{tikzpicture}[line cap=round,line join=round,>=triangle 45,x=1cm,y=0.75cm]
\clip(-3,-1.25) rectangle (3,1.25);
\draw [line width=1pt] (-2.5,0)--(2.5,0);
\draw [line width=2pt, red] (-1.5,-1)-- (-1.5,1);
\draw [line width=2pt, red] (-2,0.25)-- (-1,-0.25);
\draw[line width=2pt,color=red,smooth,samples=100,domain=1.5:2.5] plot(\x,{sqrt((\x)-1.5)});
\draw[line width=2pt,color=red,smooth,samples=100,domain=1.5:2.5] plot(\x,{0-sqrt((\x)-1.5)});
\draw[line width=2pt,color=red,smooth,samples=100,domain=0.5:1.5] plot(\x,{sqrt(-((\x)-1.5))});
\draw[line width=2pt,color=red,smooth,samples=100,domain=0.5:1.5] plot(\x,{0-sqrt(-((\x)-1.5))});
\begin{scriptsize}
\draw [color=black] (-2.75,0) node {$\Gamma$};
\draw [color=black] (-1.75,-0.5) node {$\Delta_0$};
\draw [color=black] (-2,0.5) node {$B_2$};
\draw [color=black] (1,-0.25) node {$B_3$};
\draw [color=black] (2,-0.25) node {$B_4$};
\end{scriptsize}
\end{tikzpicture}
\end{center}
\end{minipage} &
\begin{minipage}{70mm}
 $B_{\Sigma}=B_1+B_2+B_3+B_4$,\\
   $B_1 = \Delta_0$,
   $(B_2\cdot \Gamma)=1$,\\
   $B_3$ and $B_4$ form a $k$-fold node.
\end{minipage} \\ \hline
$(\mathrm{C4})$  &
\begin{minipage}{70mm}
\begin{center}
\begin{tikzpicture}[line cap=round,line join=round,>=triangle 45,x=1cm,y=0.75cm]
\clip(-3,-1.25) rectangle (3,1.25);
\draw [line width=1pt] (2.5,0)-- (-2.5,0);
\draw [line width=2pt, red] (-1.5,-1)-- (-1.5,1);
\draw [line width=2pt, red] (-2,0.25)-- (-1,-0.25);
\draw[line width=2pt,color=red,smooth,samples=100,domain=1.5:2.5] plot(\x,{sqrt((\x)-1.5)});
\draw[line width=2pt,color=red,smooth,samples=100,domain=0.5:1.5] plot(\x,{sqrt(-((\x)-1.5))});
\begin{scriptsize}
\draw [color=black] (-2.75,0) node {$\Gamma$};
\draw [color=black] (-1.75,-0.5) node {$\Delta_0$};
\draw [color=black] (-2,0.5) node {$B_2$};
\draw [color=black] (1.5,-0.25) node {$B_3$};
\end{scriptsize}
\end{tikzpicture}
\end{center}
\end{minipage} &
\begin{minipage}{70mm}
 $B_{\Sigma}=B_1+B_2+B_3$,\\
   $B_1 = \Delta_0$,
   $(B_2\cdot \Gamma)=1$,\\
   $B_3$ forms a $k$-fold cusp. 
\end{minipage} \\ \hline

$(\mathrm{C5})$  &
\begin{minipage}{70mm}
\begin{center}
\begin{tikzpicture}[line cap=round,line join=round,>=triangle 45,x=1cm,y=0.75cm]
\clip(-3,-1.25) rectangle (3,1.25);
\draw [line width=1pt,red] (-2.5,0)--(2.5,0);
\draw [line width=2pt, red] (-1.5,-1)-- (-1.5,1);
\draw [line width=2pt, red] (-0.5,-1)-- (-0.5,1);
\draw[line width=2pt,color=red,smooth,samples=100,domain=1.5:2.5] plot(\x,{sqrt((\x)-1.5)});
\draw[line width=2pt,color=red,smooth,samples=100,domain=1.5:2.5] plot(\x,{0-sqrt((\x)-1.5)});
\draw[line width=2pt,color=red,smooth,samples=100,domain=0.5:1.5] plot(\x,{sqrt(-((\x)-1.5))});
\draw[line width=2pt,color=red,smooth,samples=100,domain=0.5:1.5] plot(\x,{0-sqrt(-((\x)-1.5))});
\begin{scriptsize}
\draw [color=black] (-2.75,0) node {$\Gamma$};
\draw [color=black] (-1.75,-0.5) node {$\Delta_0$};
\draw [color=black] (-0.25,-0.25) node {$B_2$};
\draw [color=black] (1,-0.25) node {$B_3$};
\draw [color=black] (2,-0.25) node {$B_4$};
\end{scriptsize}
\end{tikzpicture}
\end{center}
\end{minipage} &
\begin{minipage}{70mm}
$B_{\Sigma}=B_1+B_2+B_3+B_4+\Gamma$,\\
$B_1 = \Delta_0$,
$(B_2\cdot \Gamma)=1$,\\
$B_3$ and $B_4$ form a $k$-fold node. 
\end{minipage} \\ \hline

$(\mathrm{C6})$  &
\begin{minipage}{70mm}
\begin{center}
\begin{tikzpicture}[line cap=round,line join=round,>=triangle 45,x=1cm,y=0.75cm]
\clip(-3,-1.25) rectangle (3,1.25);
\draw [line width=1pt,red] (2.5,0)-- (-2.5,0);
\draw [line width=2pt, red] (-1.5,-1)-- (-1.5,1);
\draw [line width=2pt, red] (-0.5,-1)-- (-0.5,1);
\draw[line width=2pt,color=red,smooth,samples=100,domain=1.5:2.5] plot(\x,{sqrt((\x)-1.5)});
\draw[line width=2pt,color=red,smooth,samples=100,domain=0.5:1.5] plot(\x,{sqrt(-((\x)-1.5))});
\begin{scriptsize}
\draw [color=black] (-2.75,0) node {$\Gamma$};
\draw [color=black] (-1.75,-0.5) node {$\Delta_0$};
\draw [color=black] (-0.25,-0.25) node {$B_2$};
\draw [color=black] (1.5,-0.25) node {$B_3$};
\end{scriptsize}
\end{tikzpicture}
\end{center}
\end{minipage} &
\begin{minipage}{70mm}
$B_{\Sigma}=B_1+B_2+B_3+\Gamma$,\\
$B_1 = \Delta_0$,
$(B_2\cdot \Gamma)=1$,\\
$B_3$ forms a $k$-fold cusp.
\end{minipage} \\ \hline
$(\mathrm{C7})$  &
\begin{minipage}{70mm}
\begin{center}
\begin{tikzpicture}[line cap=round,line join=round,>=triangle 45,x=1cm,y=0.75cm]
\clip(-3,-1.25) rectangle (3,1.25);
\draw [line width=2pt, red] (-1.5,-1)-- (-1.5,1);
\draw [line width=1pt, red] (2.5,0)-- (-2.5,0);
\draw[line width=2pt,color=black,smooth,samples=100,domain=0.5:2.5,red] plot(\x,{(\x - 1.5)^2});
\begin{scriptsize}
\draw [color=black] (-2.75,0) node {$\Gamma$};
\draw [color=black] (-1.75,-0.5) node {$\Delta_0$};
\draw [color=black] (1.5,-0.25) node {$B_2$};
\end{scriptsize}
\end{tikzpicture}
\end{center}
\end{minipage} &
\begin{minipage}{70mm}
$B_{\Sigma}=B_1+B_2+\Gamma$,\\
$B_1 = \Delta_0$,
$(B_2\cdot \Gamma)=3$.\\
\end{minipage} \\ \hline
$(\mathrm{C8})$  &
\begin{minipage}{70mm}
\begin{center}
\begin{tikzpicture}[line cap=round,line join=round,>=triangle 45,x=1cm,y=0.75cm]
\clip(-3,-1.25) rectangle (3,1.25);
\draw [line width=2pt, red] (-1.5,-1)-- (-1.5,1);
\draw [line width=2pt, red] (1.5,-1)-- (1.5,1);
\draw [line width=1pt, red] (2.5,0)-- (-2.5,0);
\draw[line width=2pt,color=black,smooth,samples=100,domain=0.5:2.5,red] plot(\x,{(\x - 1.5)^2});
\begin{scriptsize}
\draw [color=black] (-2.75,0) node {$\Gamma$};
\draw [color=black] (-1.75,-0.5) node {$\Delta_0$};
\draw [color=black] (1.25,-0.25) node {$B_2$};
\draw [color=black] (2,0.7) node {$B_3$};
\end{scriptsize}
\end{tikzpicture}
\end{center}
\end{minipage} &
\begin{minipage}{70mm}
$B_{\Sigma}=B_1+B_2+B_3+\Gamma$,\\
$B_1 = \Delta_0$,
$(B_2\cdot \Gamma)=1$,
$(B_3\cdot \Gamma)=2$.\\
\end{minipage} \\ \hline
$(\mathrm{C9})$  &
\begin{minipage}{70mm}
\begin{center}
\begin{tikzpicture}[line cap=round,line join=round,>=triangle 45,x=1cm,y=0.75cm]
\clip(-3,-1.25) rectangle (3,1.25);
\draw [line width=2pt, red] (-1.5,-1)-- (-1.5,1);
\draw [line width=1pt, red] (2.5,0)-- (-2.5,0);
\draw[line width=2pt,color=black,smooth,samples=100,domain=1:2,red] plot(\x,{(\x - 1)^2});
\draw[line width=2pt,color=black,smooth,samples=100,domain=1:2,red] plot(\x,{0-(\x - 1)^2});
\begin{scriptsize}
\draw [color=black] (-2.75,0) node {$\Gamma$};
\draw [color=black] (-1.75,-0.5) node {$\Delta_0$};
\draw [color=black] (1,-0.25) node {$B_2$};
\end{scriptsize}
\end{tikzpicture}
\end{center}
\end{minipage} &
\begin{minipage}{70mm}
$B_{\Sigma}=B_1+B_2+\Gamma$,\\
$B_1 = \Delta_0$,
$(B_2\cdot \Gamma)=3$,\\
$B_2$ has a cusp with multiplicity $2$.\\
\end{minipage} \\ \hline
$(\mathrm{C10})$  &
\begin{minipage}{70mm}
\begin{center}
\begin{tikzpicture}[line cap=round,line join=round,>=triangle 45,x=1cm,y=0.75cm]
\clip(-3,-1.25) rectangle (3,1.25);
\draw [line width=1pt,red] (-2.5,0)--(2.5,0);
\draw [line width=2pt,color=red] (-1.5,-1)-- (-1.5,1);
\draw [line width=2pt,color=red] (1,0.25)-- (2,-0.25);
\draw[line width=2pt,color=red,smooth,samples=100,domain=1.5:2.5] plot(\x,{sqrt((\x)-1.5)});
\draw[line width=2pt,color=red,smooth,samples=100,domain=1.5:2.5] plot(\x,{0-sqrt((\x)-1.5)});
\draw[line width=2pt,color=red,smooth,samples=100,domain=0.5:1.5] plot(\x,{sqrt(-((\x)-1.5))});
\draw[line width=2pt,color=red,smooth,samples=100,domain=0.5:1.5] plot(\x,{0-sqrt(-((\x)-1.5))});
\begin{scriptsize}
\draw [color=black] (-2.75,0) node {$\Gamma$};
\draw [color=black] (-1.75,-0.5) node {$\Delta_0$};
\draw [color=black] (2.25,-0.35) node {$B_2$};
\draw [color=black] (0.4,0.75) node {$B_3$};
\draw [color=black] (2.7,0.75) node {$B_4$};
\end{scriptsize}
\end{tikzpicture}
\end{center}
\end{minipage} &
\begin{minipage}{70mm}
$B_{\Sigma}=B_1+B_2+B_3+B_4+\Gamma$,\\
  $B_1 = \Delta_0$,
   $(B_2\cdot \Gamma)=1$,\\
   $B_3$ and $B_4$ form a $k$-fold node. 
\end{minipage} \\ \hline
$(\mathrm{C11})$  &
\begin{minipage}{70mm}
\begin{center}
\begin{tikzpicture}[line cap=round,line join=round,>=triangle 45,x=1cm,y=0.74cm]
\clip(-3,-1.25) rectangle (3,1.25);
\draw [line width=1pt,red] (2.5,0)-- (-2.5,0);
\draw [line width=2pt,color=red] (-1.5,-1)-- (-1.5,1);
\draw [line width=2pt,color=red] (1,0.25)-- (2,-0.25);
\draw[line width=2pt,color=red,smooth,samples=100,domain=1.5:2.5] plot(\x,{sqrt((\x)-1.5)});
\draw[line width=2pt,color=red,smooth,samples=100,domain=0.5:1.5] plot(\x,{sqrt(-((\x)-1.5))});
\begin{scriptsize}
\draw [color=black] (-2.75,0) node {$\Gamma$};
\draw [color=black] (-1.75,-0.5) node {$\Delta_0$};
\draw [color=black] (2.25,-0.35) node {$B_2$};
\draw [color=black] (0.4,0.75) node {$B_3$};
\end{scriptsize}
\end{tikzpicture}
\end{center}
\end{minipage} &
\begin{minipage}{70mm}
$B_{\Sigma}=B_1+B_2+B_3+\Gamma$,\\
$B_1 = \Delta_0$,
$(B_2\cdot \Gamma)=1$,\\
$B_3$ forms a $k$-fold cusp. 
\end{minipage} \\ \hline
\end{tabular}
\end{table}

\begin{table}[H]
\caption{}
\label{table_generate_[*SR]}
  \begin{tabular}{|c|c|c|} 
    \hline & figure & local analytic branches of $B_{\Sigma}$ \\ \hline
 $(\mathrm{SR}1)$   &
    \begin{minipage}{70mm}
\begin{center}
\begin{tikzpicture}[line cap=round,line join=round,>=triangle 45,x=1cm,y=0.75cm]
\clip(-3,-1.25) rectangle (3,1.25);
\draw [line width=1pt] (-2.5,0)--(2.5,0);
\draw [line width=2pt,color=red] (-1.5,-1)-- (-1.5,1);
\draw[line width=2pt,color=red,smooth,samples=100,domain=-2.5:-1.5] plot(\x,{sqrt(-((\x)+1.5))});
\draw[line width=2pt,color=red,smooth,samples=100,domain=-2.5:-1.5] plot(\x,{0-sqrt(-((\x)+1.5))});
\draw[line width=2pt,color=red,smooth,samples=100,domain=1.5:2.5] plot(\x,{sqrt((\x)-1.5)});
\draw[line width=2pt,color=red,smooth,samples=100,domain=1.5:2.5] plot(\x,{0-sqrt((\x)-1.5)});
\draw[line width=2pt,color=red,smooth,samples=100,domain=0.5:1.5] plot(\x,{sqrt(-((\x)-1.5))});
\draw[line width=2pt,color=red,smooth,samples=100,domain=0.5:1.5] plot(\x,{0-sqrt(-((\x)-1.5))});
\begin{scriptsize}
\draw [color=black] (-2.75,0) node {$\Gamma$};
\draw [color=black] (-1.25,-0.5) node {$\Delta_0$};
\draw [color=black] (-2,0.25) node {$B_2$};
\draw [color=black] (1,-0.25) node {$B_3$};
\draw [color=black] (2,-0.25) node {$B_4$};
\end{scriptsize}
\end{tikzpicture}
\end{center}
    \end{minipage} &
    \begin{minipage}{80mm}
$B_{\Sigma}=B_1+B_2+B_3+B_4$,\\
   $B_1 = \Delta_0$,
   $(B_2\cdot \Gamma)=1$, $(B_2 \cdot \Delta_0)=2$,\\
   $B_3$ and $B_4$ form a $k$-fold node. 
 \end{minipage} \\ \hline
 $(\mathrm{SR}2)$ &
\begin{minipage}{70mm}
\begin{center}
\begin{tikzpicture}[line cap=round,line join=round,>=triangle 45,x=1cm,y=0.75cm]
\clip(-3,-1.25) rectangle (3,1.25);
\draw [line width=1pt] (2.5,0)-- (-2.5,0);
\draw [line width=2pt,color=red] (-1.5,-1)-- (-1.5,1);
\draw[line width=2pt,color=red,smooth,samples=100,domain=1.5:2.5] plot(\x,{sqrt((\x)-1.5)});
\draw[line width=2pt,color=red,smooth,samples=100,domain=0.5:1.5] plot(\x,{sqrt(-((\x)-1.5))});
\draw[line width=2pt,color=red,smooth,samples=100,domain=-2.5:-1.5] plot(\x,{sqrt(-((\x)+1.5))});
\draw[line width=2pt,color=red,smooth,samples=100,domain=-2.5:-1.5] plot(\x,{0-sqrt(-((\x)+1.5))});
\begin{scriptsize}
\draw [color=black] (-2.75,0) node {$\Gamma$};
\draw [color=black] (-1.25,-0.5) node {$\Delta_0$};
\draw [color=black] (-2,0.25) node {$B_2$};
\draw [color=black] (1.5,-0.25) node {$B_3$};
\end{scriptsize}
\end{tikzpicture}
\end{center}
\end{minipage} &
\begin{minipage}{80mm}
$B_{\Sigma}=B_1+B_2+B_3$,\\
   $B_1 = \Delta_0$,
   $(B_2\cdot \Gamma)=1$, $(B_2 \cdot \Delta_0)=2$,\\
   $B_3$ forms a $k$-fold cusp. 
\end{minipage} \\ \hline
$(\mathrm{SR}3)$ & 
\begin{minipage}{70mm}
\begin{center}
\begin{tikzpicture}[line cap=round,line join=round,>=triangle 45,x=1cm,y=0.75cm]
\clip(-3,-1.25) rectangle (3,1.25);
\draw [line width=1pt,red] (-2.5,0)--(2.5,0);
\draw [line width=2pt,color=red] (-1.5,-1)-- (-1.5,1);
\draw [line width=2pt,color=red] (-2,0.25)-- (-1,-0.25);
\draw[line width=2pt,color=red,smooth,samples=100,domain=1.5:2.5] plot(\x,{sqrt((\x)-1.5)});
\draw[line width=2pt,color=red,smooth,samples=100,domain=1.5:2.5] plot(\x,{0-sqrt((\x)-1.5)});
\draw[line width=2pt,color=red,smooth,samples=100,domain=0.5:1.5] plot(\x,{sqrt(-((\x)-1.5))});
\draw[line width=2pt,color=red,smooth,samples=100,domain=0.5:1.5] plot(\x,{0-sqrt(-((\x)-1.5))});
\begin{scriptsize}
\draw [color=black] (-2.75,0) node {$\Gamma$};
\draw [color=black] (-1.75,-0.5) node {$\Delta_0$};
\draw [color=black] (-2,0.5) node {$B_2$};
\draw [color=black] (1,-0.25) node {$B_3$};
\draw [color=black] (2,-0.25) node {$B_4$};
\end{scriptsize}
\end{tikzpicture}
\end{center}
\end{minipage} &
\begin{minipage}{80mm}
$B_{\Sigma}=B_1+B_2+B_3+B_4+\Gamma$,\\
   $B_1 = \Delta_0$,
   $(B_2\cdot \Gamma)=1$, $(B_2 \cdot \Delta_0)=1$,\\
   $B_3$ and $B_4$ form a $k$-fold node.
\end{minipage} \\ \hline
$(\mathrm{SR}4)$  &
\begin{minipage}{70mm}
\begin{center}
\begin{tikzpicture}[line cap=round,line join=round,>=triangle 45,x=1cm,y=0.75cm]
\clip(-3,-1.25) rectangle (3,1.25);
\draw [line width=1pt,color=red] (2.5,0)-- (-2.5,0);
\draw [line width=2pt,color=red] (-1.5,-1)-- (-1.5,1);
\draw [line width=2pt,color=red] (-2,0.25)-- (-1,-0.25);
\draw[line width=2pt,color=red,smooth,samples=100,domain=1.5:2.5] plot(\x,{sqrt((\x)-1.5)});
\draw[line width=2pt,color=red,smooth,samples=100,domain=0.5:1.5] plot(\x,{sqrt(-((\x)-1.5))});
\begin{scriptsize}
\draw [color=black] (-2.75,0) node {$\Gamma$};
\draw [color=black] (-1.75,-0.5) node {$\Delta_0$};
\draw [color=black] (-2,0.5) node {$B_2$};
\draw [color=black] (1.5,-0.25) node {$B_3$};
\end{scriptsize}
\end{tikzpicture}
\end{center}
\end{minipage} &
\begin{minipage}{80mm}
$B_{\Sigma}=B_1+B_2+B_3+\Gamma$,\\
   $B_1 = \Delta_0$,
   $(B_2\cdot \Gamma)=1$, $(B_2 \cdot \Delta_0)=1$,\\
   $B_3$ forms a $k$-fold cusp. 
\end{minipage} \\ \hline

$(\mathrm{SR}5)$  &
\begin{minipage}{70mm}
\begin{center}
\begin{tikzpicture}[line cap=round,line join=round,>=triangle 45,x=1cm,y=0.75cm]
\clip(-3,-1.25) rectangle (3,1.25);
\draw [line width=1pt] (2.5,0)-- (-2.5,0);
\draw [line width=2pt,color=red] (-1.5,-1)-- (-1.5,1);
\draw[line width=2pt,smooth,samples=100,domain=-2.5:-1.5,red] plot(\x,{(\x + 1.5)^2});
\draw[line width=2pt,smooth,samples=100,domain=-2.5:-1.5,red] plot(\x,{0-(\x + 1.5)^2});
\begin{scriptsize}
\draw [color=black] (-2.75,0) node {$\Gamma$};
\draw [color=black] (-1.25,-0.5) node {$\Delta_0$};
\draw [color=black] (-2,0.6) node {$B_2$};
\end{scriptsize}
\end{tikzpicture}
\end{center}
\end{minipage} &
\begin{minipage}{80mm}
$B_{\Sigma}=B_1+B_2$,\\
   $B_1 = \Delta_0$,
   $(B_2\cdot \Gamma)=3$.\\
\end{minipage} \\ \hline

$(\mathrm{SR}6)$  &
\begin{minipage}{70mm}
\begin{center}
\begin{tikzpicture}[line cap=round,line join=round,>=triangle 45,x=1cm,y=0.75cm]
\clip(-3,-1.25) rectangle (3,1.25);
\draw [line width=1pt,color=red] (2.5,0)-- (-2.5,0);
\draw [line width=2pt,color=red] (-1.5,-1)-- (-1.5,1);
\draw [line width=2pt,color=red] (1.5,-1)-- (1.5,1);
\draw[line width=2pt,color=black,smooth,samples=100,domain=-2.5:-0.5,red] plot(\x,{(\x + 1.5)^2});
\begin{scriptsize}
\draw [color=black] (-2.75,0) node {$\Gamma$};
\draw [color=black] (-1.75,-0.5) node {$\Delta_0$};
\draw [color=black] (-2,0.65) node {$B_2$};
\draw [color=black] (1.75,-0.25) node {$B_3$};
\end{scriptsize}
\end{tikzpicture}
\end{center}
\end{minipage} &
\begin{minipage}{80mm}
$B_{\Sigma}=B_1+B_2+B_3+\Gamma$,\\
   $B_1 = \Delta_0$,
   $(B_2\cdot \Gamma)=2$, $(B_2 \cdot \Delta_0)=1$,\\
   $(B_3\cdot \Gamma)=1$. 
\end{minipage} \\ \hline
$(\mathrm{SR}7)$  &
\begin{minipage}{70mm}
\begin{center}
\begin{tikzpicture}[line cap=round,line join=round,>=triangle 45,x=1cm,y=0.75cm]
\clip(-3,-1.25) rectangle (3,1.25);
\draw [line width=1pt,color=red] (2.5,0)-- (-2.5,0);
\draw [line width=2pt,color=red] (-1.5,-1)-- (-1.5,1);
\draw[line width=2pt,color=black,smooth,samples=100,domain=-2.5:-0.5,red] plot(\x,{(\x + 1.5)^2});
\begin{scriptsize}
\draw [color=black] (-2.75,0) node {$\Gamma$};
\draw [color=black] (-1.75,-0.5) node {$\Delta_0$};
\draw [color=black] (-2,0.65) node {$B_2$};
\end{scriptsize}
\end{tikzpicture}
\end{center}
\end{minipage} &
\begin{minipage}{80mm}
$B_{\Sigma}=B_1+B_2+\Gamma$,\\
   $B_1 = \Delta_0$,
   $(B_2\cdot \Gamma)=3$.
\end{minipage} \\ \hline
\end{tabular}
\end{table}


\begin{table}[H]
\caption{}
\label{table_generate_[LP]}
  \begin{tabular}{|c|c|c|} 
    \hline  & figure & local analytic branches of $B_{\Sigma}$ \\ \hline
 $(\mathrm{LP}1)$   &
    \begin{minipage}{70mm}
\begin{center}
\begin{tikzpicture}[line cap=round,line join=round,>=triangle 45,x=1cm,y=0.75cm]
\clip(-3,-1.25) rectangle (3,1.25);
\draw [line width=1pt] (-2.5,0)--(2.5,0);
\draw [line width=1pt] (-1.5,-1)-- (-1.5,1);
\draw[line width=2pt,color=red,smooth,samples=100,domain=1.5:2] plot(\x,{sqrt((\x)-1.5)});
\draw[line width=2pt,color=red,smooth,samples=100,domain=1.5:2] plot(\x,{0-sqrt((\x)-1.5)});
\draw[line width=2pt,color=red,smooth,samples=100,domain=1:1.5] plot(\x,{sqrt(-((\x)-1.5))});
\draw[line width=2pt,color=red,smooth,samples=100,domain=1:1.5] plot(\x,{0-sqrt(-((\x)-1.5))});
\draw[line width=2pt,color=red,smooth,samples=100,domain=-0.5:0] plot(\x,{sqrt((\x)+0.5)});
\draw[line width=2pt,color=red,smooth,samples=100,domain=-1:-0.5] plot(\x,{sqrt(-((\x)+0.5))});
\draw[line width=2pt,color=red,smooth,samples=100,domain=-0.5:0] plot(\x,{0-sqrt((\x)+0.5)});
\draw[line width=2pt,color=red,smooth,samples=100,domain=-1:-0.5] plot(\x,{0-sqrt(-((\x)+0.5))});
\begin{scriptsize}
\draw [color=black] (-2.75,0) node {$\Gamma$};
\draw [color=black] (-1.75,-0.5) node {$\Delta_0$};
\draw [color=black] (-1,-0.25) node {$B_1$};
\draw [color=black] (0,-0.25) node {$B_2$};
\draw [color=black] (1,-0.25) node {$B_3$};
\draw [color=black] (2,-0.25) node {$B_4$};
\end{scriptsize}
\end{tikzpicture}
\end{center}
     \end{minipage} &
\begin{minipage}{80mm}
$B_{\Sigma}=B_1+B_2+B_3+B_4$, $(B_i \cdot \Gamma)=1$,\\
   $B_1$ and $B_2$ form a $k_1$-fold node ($k_1 \ge \chi-4$),\\
   $B_3$ and $B_4$ form a $k_2$-fold node ($k_2 \ge 1$). 
\end{minipage} \\ \hline
 $(\mathrm{LP}2)$   &
    \begin{minipage}{70mm}
\begin{center}
\begin{tikzpicture}[line cap=round,line join=round,>=triangle 45,x=1cm,y=0.75cm]
\clip(-3,-1.25) rectangle (3,1.25);
\draw [line width=1pt] (2.5,0)-- (-2.5,0);
\draw [line width=1pt] (-1.5,-1)-- (-1.5,1);
\draw[line width=2pt,color=red,smooth,samples=100,domain=1.5:2] plot(\x,{sqrt((\x)-1.5)});
\draw[line width=2pt,color=red,smooth,samples=100,domain=1.5:2] plot(\x,{0-sqrt((\x)-1.5)});
\draw[line width=2pt,color=red,smooth,samples=100,domain=1:1.5] plot(\x,{sqrt(-((\x)-1.5))});
\draw[line width=2pt,color=red,smooth,samples=100,domain=1:1.5] plot(\x,{0-sqrt(-((\x)-1.5))});
\draw[line width=2pt,color=red,smooth,samples=100,domain=-0.5:0] plot(\x,{sqrt((\x)+0.5)});
\draw[line width=2pt,color=red,smooth,samples=100,domain=-1:-0.5] plot(\x,{sqrt(-((\x)+0.5))});
\begin{scriptsize}
\draw [color=black] (-2.75,0) node {$\Gamma$};
\draw [color=black] (-1.75,-0.5) node {$\Delta_0$};
\draw [color=black] (-0.5,-0.25) node {$B_1$};
\draw [color=black] (1,-0.25) node {$B_2$};
\draw [color=black] (2,-0.25) node {$B_3$};
\end{scriptsize}
\end{tikzpicture}
\end{center}
     \end{minipage} &
    \begin{minipage}{80mm}
$B_{\Sigma}=B_1+B_2+B_3$, $(B_i \cdot \Gamma)=1$ for $i=2,3$,\\
   $B_1$ forms a $k_1$-fold cusp ($k_1 \ge \chi-4$),\\ 
   $B_2$ and $B_3$ form a $k_2$-fold node ($k_2 \ge 1$). 
    \end{minipage} \\ \hline
 $(\mathrm{LP}3)$   &
    \begin{minipage}{70mm}
\begin{center}
\begin{tikzpicture}[line cap=round,line join=round,>=triangle 45,x=1cm,y=0.75cm]
\clip(-3,-1.25) rectangle (3,1.25);
\draw [line width=1pt] (2.5,0)-- (-2.5,0);
\draw [line width=1pt] (-1.5,-1)-- (-1.5,1);
\draw[line width=2pt,color=red,smooth,samples=100,domain=-0.5:0] plot(\x,{sqrt((\x)+0.5)});
\draw[line width=2pt,color=red,smooth,samples=100,domain=-1:-0.5] plot(\x,{sqrt(-((\x)+0.5))});
\draw[line width=2pt,color=red,smooth,samples=100,domain=-0.5:0] plot(\x,{0-sqrt((\x)+0.5)});
\draw[line width=2pt,color=red,smooth,samples=100,domain=-1:-0.5] plot(\x,{0-sqrt(-((\x)+0.5))});
\draw[line width=2pt,color=red,smooth,samples=100,domain=1.5:2] plot(\x,{sqrt((\x)-1.5)});
\draw[line width=2pt,color=red,smooth,samples=100,domain=1:1.5] plot(\x,{sqrt(-((\x)-1.5))});
\begin{scriptsize}
\draw [color=black] (-2.75,0) node {$\Gamma$};
\draw [color=black] (-1.75,-0.5) node {$\Delta_0$};
\draw [color=black] (-1,-0.25) node {$B_1$};
\draw [color=black] (0,-0.25) node {$B_2$};
\draw [color=black] (1.5,-0.25) node {$B_3$};
\end{scriptsize}
\end{tikzpicture}
\end{center}
     \end{minipage} &
    \begin{minipage}{80mm}
$B_{\Sigma}=B_1+B_2+B_3$, $(B_i \cdot \Gamma)=1$ for $i=1,2$,\\
   $B_1$ and $B_2$ form a $k_1$-fold node ($k_1 \ge \chi-4$),
   \\ 
   $B_3$ forms a $k_2$-fold cusp ($k_2 \ge 1$). 
    \end{minipage} \\ \hline
 $(\mathrm{LP}4)$   &
    \begin{minipage}{70mm}
\begin{center}
\begin{tikzpicture}[line cap=round,line join=round,>=triangle 45,x=1cm,y=0.75cm]
\clip(-3,-1.25) rectangle (3,1.25);
\draw [line width=1pt] (2.5,0)-- (-2.5,0);
\draw [line width=1pt] (-1.5,-1)-- (-1.5,1);
\draw[line width=2pt,color=red,smooth,samples=100,domain=1.5:2] plot(\x,{sqrt((\x)-1.5)});
\draw[line width=2pt,color=red,smooth,samples=100,domain=1:1.5] plot(\x,{sqrt(-((\x)-1.5))});
\draw[line width=2pt,color=red,smooth,samples=100,domain=-0.5:0] plot(\x,{sqrt((\x)+0.5)});
\draw[line width=2pt,color=red,smooth,samples=100,domain=-1:-0.5] plot(\x,{sqrt(-((\x)+0.5))});
\begin{scriptsize}
\draw [color=black] (-2.75,0) node {$\Gamma$};
\draw [color=black] (-1.75,-0.5) node {$\Delta_0$};
\draw [color=black] (-0.5,-0.25) node {$B_1$};
\draw [color=black] (1.5,-0.25) node {$B_2$};
\end{scriptsize}
\end{tikzpicture}
\end{center}
     \end{minipage} &
    \begin{minipage}{80mm}
$B_{\Sigma}=B_1+B_2$,\\
   $B_1$ forms a $k_1$-fold cusp ($k_1 \ge \chi-4$),\\ 
   $B_2$ forms a $k_2$-fold node ($k_2 \ge 1$). 
    \end{minipage} \\ \hline
 $(\mathrm{LP}5)$   &
    \begin{minipage}{70mm}
    \begin{center}
\begin{tikzpicture}[line cap=round,line join=round,>=triangle 45,x=1cm,y=0.75cm]
\clip(-3,-1.25) rectangle (3,1.25);
\draw [line width=1pt] (2.5,0)-- (-2.5,0);
\draw [line width=1pt] (-1.5,-1)-- (-1.5,1);
\draw[line width=2pt,color=black,smooth,samples=100,domain=-1:1,red] plot(\x,{(\x)^2});
\draw[line width=2pt,color=black,smooth,samples=100,domain=-1:1,red] plot(\x,{0-(\x)^2});
\begin{scriptsize}
\draw [color=black] (-2.75,0) node {$\Gamma$};
\draw [color=black] (-1.75,-0.5) node {$\Delta_0$};
\draw [color=black] (0,0.5) node {$B_2$};
\draw [color=black] (0,-0.5) node {$B_1$};
\end{scriptsize}
\end{tikzpicture}
    \end{center}
     \end{minipage} &
    \begin{minipage}{80mm}
   $B_{\Sigma}=B_1+B_2$, 
   $(B_i \cdot \Gamma)=2$ ($i=1,2$),\\ 
   $(B_1\cdot B_2) \ge 2$,  $B_1$ and $B_2$ are smooth.\\
   After two times blow-ups at the singular point of $B_1+B_2$, it becomes a $k$-fold node with $k\ge \chi-5$.
    \end{minipage} \\ \hline

 $(\mathrm{LP}6)$   &
    \begin{minipage}{70mm}
    \begin{center}
\begin{tikzpicture}[line cap=round,line join=round,>=triangle 45,x=1cm,y=0.75cm]
\clip(-3,-1.25) rectangle (3,1.25);
\draw [line width=1pt] (2.5,0)-- (-2.5,0);
\draw [line width=1pt] (-1.5,-1)-- (-1.5,1);
\draw[line width=2pt,color=black,smooth,samples=100,domain=0:1,red] plot(\x,{(\x)^2});
\draw[line width=2pt,color=black,smooth,samples=100,domain=0:1,red] plot(\x,{0-(\x)^2});
\begin{scriptsize}
\draw [color=black] (-2.75,0) node {$\Gamma$};
\draw [color=black] (-1.75,-0.5) node {$\Delta_0$};
\draw [color=black] (0,-0.25) node {$B_1$};
\end{scriptsize}
\end{tikzpicture}
  \end{center}
     \end{minipage} &
    \begin{minipage}{80mm}
 $B_{\Sigma}=B_1$,
   $(B_1 \cdot \Gamma)=4$, 
   $B_1$ has a cusp with multiplicity $2$.\\
   After two times blow-ups at the singular point of $B_1$, it becomes a $k$-fold cusp with $k\ge \chi-5$. 
    \end{minipage} \\ \hline

 $(\mathrm{LP}7)$   &
    \begin{minipage}{70mm}
\begin{center}
\begin{tikzpicture}[line cap=round,line join=round,>=triangle 45,x=1cm,y=0.75cm]
\clip(-3,-1.25) rectangle (3,1.25);
\draw [line width=1pt] (-2.5,0)--(2.5,0);
\draw [line width=1pt] (-1.5,-1)-- (-1.5,1);
\draw[line width=2pt,color=black,smooth,samples=100,domain=0.75:2.25,red] plot(\x,{(\x - 1.5)^2});
\draw[line width=2pt,color=red,smooth,samples=100,domain=1.5:2.5] plot(\x,{sqrt((\x)-1.5)});
\draw[line width=2pt,color=red,smooth,samples=100,domain=1.5:2.5] plot(\x,{0-sqrt((\x)-1.5)});
\draw[line width=2pt,color=red,smooth,samples=100,domain=0.5:1.5] plot(\x,{sqrt(-((\x)-1.5))});
\draw[line width=2pt,color=red,smooth,samples=100,domain=0.5:1.5] plot(\x,{0-sqrt(-((\x)-1.5))});
\begin{scriptsize}
\draw [color=black] (-2.75,0) node {$\Gamma$};
\draw [color=black] (-1.75,-0.5) node {$\Delta_0$};
\draw [color=black] (0.5,0.35) node {$B_1$};
\draw [color=black] (1,-0.25) node {$B_2$};
\draw [color=black] (2,-0.25) node {$B_3$};
\end{scriptsize}
\end{tikzpicture}
\end{center}
     \end{minipage} &
    \begin{minipage}{80mm}
$B_{\Sigma}=B_1+B_2+B_3$,\\
$(B_1 \cdot \Gamma)=2$, $(B_i \cdot \Gamma)=1$ for $(i=2,3)$,\\
   $B_1$ is smooth, $B_2$ and $B_3$ form a $k$-fold node.
    \end{minipage} \\ \hline
     $(\mathrm{LP}8)$   &
    \begin{minipage}{70mm}
\begin{center}
\begin{tikzpicture}[line cap=round,line join=round,>=triangle 45,x=1cm,y=0.75cm]
\clip(-3,-1.25) rectangle (3,1.25);
\draw [line width=1pt] (2.5,0)-- (-2.5,0);
\draw [line width=1pt] (-1.5,-1)-- (-1.5,1);
\draw[line width=2pt,color=black,smooth,samples=100,domain=0.75:2.25,red] plot(\x,{(\x - 1.5)^2});
\draw[line width=2pt,color=red,smooth,samples=100,domain=1.5:2.5] plot(\x,{sqrt((\x)-1.5)});
\draw[line width=2pt,color=red,smooth,samples=100,domain=0.5:1.5] plot(\x,{sqrt(-((\x)-1.5))});
\begin{scriptsize}
\draw [color=black] (-2.75,0) node {$\Gamma$};
\draw [color=black] (-1.75,-0.5) node {$\Delta_0$};
\draw [color=black] (1.5,0.75) node {$B_2$};
\draw [color=black] (1.5,-0.25) node {$B_1$};
\end{scriptsize}
\end{tikzpicture}
\end{center}
     \end{minipage} &
    \begin{minipage}{80mm}
$B_{\Sigma}=B_1+B_2$, $(B_1 \cdot \Gamma)=2$, $(B_2\cdot \Gamma)=2$,\\
   $B_1$ is smooth,\\
   $B_2$ forms a $k$-fold cusp. 
\end{minipage} \\ \hline
     $(\mathrm{LP}9)$   &
    \begin{minipage}{70mm}
    \begin{center}
\begin{tikzpicture}[line cap=round,line join=round,>=triangle 45,x=1cm,y=0.75cm]
\clip(-3,-1.25) rectangle (3,1.25);
\draw [line width=1pt] (2.5,0)-- (-2.5,0);
\draw [line width=1pt] (-1.5,-1)-- (-1.5,1);
\draw[line width=2pt,color=black,smooth,samples=100,domain=0.5:1.5,red] plot(\x,{(\x - 1.5)^2});
\draw[line width=2pt,color=black,smooth,samples=100,domain=0.5:1.5,red] plot(\x,{0-(\x - 1.5)^2});
\begin{scriptsize}
\draw [color=black] (-2.75,0) node {$\Gamma$};
\draw [color=black] (-1.75,-0.5) node {$\Delta_0$};
\draw [color=black] (1.5,-0.35) node {$B_1$};
\end{scriptsize}
\end{tikzpicture}
 \end{center}
     \end{minipage} &
\begin{minipage}{80mm}
   $B_{\Sigma}=B_1$, $(B_1 \cdot \Gamma)=4$, \\
   $B_1$ has a cusp with multiplicity $3$. 
\end{minipage} \\ \hline
     $(\mathrm{LP}10)$   &
    \begin{minipage}{70mm}
    \begin{center}
\begin{tikzpicture}[line cap=round,line join=round,>=triangle 45,x=1cm,y=0.75cm]
\clip(-3,-1.25) rectangle (3,1.25);
\draw [line width=1pt] (2.5,0)-- (-2.5,0);
\draw [line width=1pt] (-1.5,-1)-- (-1.5,1);
\draw [line width=2pt,red] (1.5,0.5)-- (1.5,-0.5);
\draw[line width=2pt,color=black,smooth,samples=100,domain=0.5:1.5,red] plot(\x,{(\x - 1.5)^2}) ;
\draw[line width=2pt,color=black,smooth,samples=100,domain=0.5:1.5,red] plot(\x,{0-(\x - 1.5)^2});
\begin{scriptsize}
\draw [color=black] (-2.75,0) node {$\Gamma$};
\draw [color=black] (-1.75,-0.5) node {$\Delta_0$};
\draw [color=black] (0.5,-0.35) node {$B_1$};
\draw [color=black] (1.75,-0.35) node {$B_2$};
\end{scriptsize}
\end{tikzpicture}
 \end{center}
     \end{minipage} &
    \begin{minipage}{80mm}
   $B_{\Sigma}=B_1+B_2$, $(B_1 \cdot \Gamma)=3$, $(B_2 \cdot \Gamma)=1$,\\
   $B_1$ has a cusp with multiplicity $2$. 
    \end{minipage} \\ \hline
\end{tabular}
\end{table}


\begin{table}[H]
\caption{}
\label{table_generate_[LPD]}
  \begin{tabular}{|c|c|c|} 
    \hline  & figure & local analytic branches of $B_{\Sigma}$ \\ \hline
     $(\mathrm{LP}11)$   &
    \begin{minipage}{70mm}
    \begin{center}
\begin{tikzpicture}[line cap=round,line join=round,>=triangle 45,x=1.cm,y=1cm]
\clip(-3,-1.25) rectangle (3,1.25);
\draw [line width=1pt] (-2.5,0)--(2.5,0);
\draw [line width=1pt] (-1.5,-1)-- (-1.5,1);
\draw[line width=2pt,color=black,smooth,samples=100,domain=1:2,red,rotate around={45:(1.5,0)}] plot(\x,{(\x - 1.5)^2});
\draw[line width=2pt,color=black,smooth,samples=100,domain=1:2,red,rotate around={45:(1.5,0)}] plot(\x,{0-(\x - 1.5)^2});
\draw[line width=2pt,color=red,smooth,samples=100,domain=1.5:2.25,rotate around={45:(1.5,0)}] plot(\x,{sqrt((\x)-1.5)});
\draw[line width=2pt,color=red,smooth,samples=100,domain=1.5:2.25,rotate around={45:(1.5,0)}] plot(\x,{0-sqrt((\x)-1.5)});
\draw[line width=2pt,color=red,smooth,samples=100,domain=0.75:1.5,rotate around={45:(1.5,0)}] plot(\x,{sqrt(-((\x)-1.5))});
\draw[line width=2pt,color=red,smooth,samples=100,domain=0.75:1.5,rotate around={45:(1.5,0)}] plot(\x,{0-sqrt(-((\x)-1.5))});
\begin{scriptsize}
\draw [color=black] (-2.75,0) node {$\Gamma$};
\draw [color=black] (-1.75,-0.5) node {$\Delta_0$};
\draw [color=black] (0.5,0.35) node {$B_1$};
\draw [color=black] (1,0.75) node {$B_2$};
\draw [color=black] (1.75,0.75) node {$B_3$};
\draw [color=black] (2.25,0.25) node {$B_4$};
\end{scriptsize}
\end{tikzpicture}
\end{center}
     \end{minipage} &
    \begin{minipage}{80mm}
$B_{\Sigma}=B_1+B_2+B_3+B_4$\\
   $B_1$ and $B_2$ form a $k_1$-fold node.\\ 
   $B_3$ and $B_4$ form a $k_2$-fold node. 
    \end{minipage} \\ \hline
     $(\mathrm{LP}12)$   &
    \begin{minipage}{70mm}
    \begin{center}
\begin{tikzpicture}[line cap=round,line join=round,>=triangle 45,x=1cm,y=1cm]
\clip(-3,-1.25) rectangle (3,1.25);
\draw [line width=1pt] (2.5,0)-- (-2.5,0);
\draw [line width=1pt] (-1.5,-1)-- (-1.5,1);
\draw[line width=2pt,color=black,smooth,samples=100,domain=1:2,red,rotate around={45:(1.5,0)}] plot(\x,{(\x - 1.5)^2});
\draw[line width=2pt,color=black,smooth,samples=100,domain=1:2,red,rotate around={45:(1.5,0)}] plot(\x,{0-(\x - 1.5)^2});
\draw[line width=2pt,color=red,smooth,samples=100,domain=1.5:2.25,rotate around={45:(1.5,0)}] plot(\x,{sqrt((\x)-1.5)});
\draw[line width=2pt,color=red,smooth,samples=100,domain=0.5:1.5,rotate around={45:(1.5,0)}] plot(\x,{sqrt(-((\x)-1.5))});
\begin{scriptsize}
\draw [color=black] (-2.75,0) node {$\Gamma$};
\draw [color=black] (-1.75,-0.5) node {$\Delta_0$};
\draw [color=black] (0.5,0.35) node {$B_1$};
\draw [color=black] (1.75,0.85) node {$B_2$};
\draw [color=black] (2.25,0.25) node {$B_3$};
\end{scriptsize}
\end{tikzpicture}
\end{center}
     \end{minipage} &
    \begin{minipage}{80mm}
   $B_{\Sigma}=B_1+B_2+B_3$\\
   $B_1$  forms a $k_1$-fold cusp. \\
   $B_2$ and $B_3$ form a $k_2$-fold node.
    \end{minipage} \\ \hline
      $(\mathrm{LP}13)$   &
    \begin{minipage}{70mm}
    \begin{center}
\begin{tikzpicture}[line cap=round,line join=round,>=triangle 45,x=1cm,y=1cm]
\clip(-3,-1.25) rectangle (3,1.25);
\draw [line width=1pt] (2.5,0)-- (-2.5,0);
\draw [line width=1pt] (-1.5,-1)-- (-1.5,1);
\draw[line width=2pt,color=red,smooth,samples=100,domain=1.5:1.8,rotate around={45:(1.5,0)}] plot(\x,{sqrt((\x)-1.5)});
\draw[line width=2pt,color=red,smooth,samples=100,domain=1:1.5,rotate around={45:(1.5,0)}] plot(\x,{sqrt(-((\x)-1.5))});
\draw[line width=2pt,color=red,smooth,samples=100,domain=1.5:2,rotate around={-45:(1.5,0)}] plot(\x,{sqrt((\x)-1.5)});
\draw[line width=2pt,color=red,smooth,samples=100,domain=1.2:1.5,rotate around={-45:(1.5,0)}] plot(\x,{sqrt(-((\x)-1.5))});

\begin{scriptsize}
\draw [color=black] (-2.75,0) node {$\Gamma$};
\draw [color=black] (-1.75,-0.5) node {$\Delta_0$};
\draw [color=black] (0.5,0.35) node {$B_1$};
\draw [color=black] (2.5,0.35) node {$B_2$};
\end{scriptsize}
\end{tikzpicture}
\end{center}
     \end{minipage} &
    \begin{minipage}{80mm}
   $B_{\Sigma}=B_1+B_2$\\
   $B_1$  forms a $k_1$-fold cusp. \\
   $B_2$ forms a $k_2$-fold cusp.
    \end{minipage} \\ \hline
\end{tabular}
\end{table}

\end{const}

\section{$\mathbb{Q}$-Gorenstein deformation revisited} \label{sec--deformation-arekore-honke}

In this section, we establish several general results on $\mathbb{Q}$-Gorenstein deformations of normal surface singularities.
Our goal is to apply these results to the $\mathbb{Q}$-Gorenstein smoothing problem for stable normal Horikawa surfaces.

\subsection{$\mathbb{Q}$-Gorenstein deformation of T-singularities and divisorial sheaves}\label{sec--deformation--arekore}
We begin with a complete criterion for the extendability of divisorial sheaves over a one-parameter $\mathbb{Q}$-Gorenstein smoothing of T-singularities.
 DeVleming-Stapleton \cite[Theorem 1.6]{DVS} gave a complete criterion for the extendability of divisorial sheaves over a one-parameter $\mathbb{Q}$-Gorenstein smoothing of T-singularities for degenerations of $\mathbb{P}^2$.
The following is a generalization of their result.

\begin{prop} \label{prop--smoothability--criterion}
Let $f\colon \mathscr{X}\to C$ be a $\mathbb{Q}$-Gorenstein deformation, and let $0\in C$ be a closed point such that $X:=\mathscr{X}_0$ is a deminormal surface.   
Let $L$ be a divisorial sheaf on $X$.
Then for any T-singularity $x\in X$ smoothed by $f$, the following are equivalent:
\begin{itemize}
\item[$(1)$]
There exist an open neighborhood $U\subset X$ of $x$ and an integer $l\in \Z$ such that $L|_{U}\sim lK_U$.

\item[$(2)$]
After shrinking $\mathscr{X}$ to an open neighborhood of $x$ if necessary, 
there exists a $\Q$-Cartier divisorial sheaf $\mathscr{L}$ on $\mathscr{X}$ such that $\mathscr{L}|_{X}\cong L$.
\end{itemize}

Moreover, if $f$ is projective, $H^2(X,\mathcal{O}_{X})=0$ and $f$ induces a locally trivial deformation for any singularity of $X$ other than the T-singularities satisfying the above equivalent conditions $(1)$ and $(2)$, then there exists a $\Q$-Cartier divisorial sheaf $\mathscr{L}$ on $\mathscr{X}$ such that $\mathscr{L}|_{X}\cong L$ after taking a finite base change $C'\to C$ from a smooth curve if necessary.
\end{prop}

\begin{proof} 
First, assume (1).
By shrinking $\mathscr{X}$ around $x$, we may assume that $L\sim lK_X$.
Then the claim (2) follows by setting $\mathscr{L}:=\omega_{\mathscr{X}/C}^{[l]}$.

Now, assume (2).
We may regard $\mathscr{L}$ as a $\Q$-Cartier Weil divisor on $\mathscr{X}$ with $\mathscr{L}|_{X}\sim L$ (cf.\ \cite[Corollary 5.25]{KoMo}).
Consider the canonical cover $\pi\colon \widetilde{\mathscr{X}}\to \mathscr{X}$.
By shrinking $C$ and applying \cite[Lemma 3.16]{KSB}, we may assume that the induced map $\pi_s\colon\widetilde{\mathscr{X}}_0\to X$ is also a canonical covering.
Consequently, $\widetilde{\mathscr{X}}_0$ has only Du Val singularities.
Next, by replacing $C$ with its finite cover, we may assume that there exists a simultaneous resolution $\mu\colon \widehat{\mathscr{X}}\to \widetilde{\mathscr{X}}$ of these singularities.
Since $\mu$ is a small birational morphism,
it follows that $\mu^*\pi^*\mathscr{L}=\mu_*^{-1}\pi^*\mathscr{L}$ as Weil divisors.
Furthermore, because $\widehat{\mathscr{X}}$ is smooth, this divisor is Cartier.
By \cite[Lemma 3.1]{has-lc-trivial-fib}, we conclude that $\pi^*\mathscr{L}$ is also Cartier. 
Hence, $\pi_0^*L\sim \pi_{0}^{*}\mathscr{L}|_{X}=\pi^*\mathscr{L}|_{\widetilde{\mathscr{X}}_0}$ is also a Cartier divisor.
Applying \cite[Proposition 3.10]{KSB}, we obtain that $lK_{U}\sim L|_{U}$ for some open neighborhood $U$ of $x$ and $l\in\mathbb{Z}$.
Indeed, we consider the case where $x$ is a T-singularity of type $\frac{1}{dn^2}(1,adn-1)$.
Then, we can take an open neighborhood $U\subset X$ of $x$ such that the Weil divisor class group satisfies $\mathrm{Cl}(U)\cong\mathbb{Z}/dn^2\mathbb{Z}$. 
Let $\widetilde{U}:=\pi_0^{-1}(U)$, noting that $\widetilde{U}$ has an $A_{dn-1}$-singularity.
It is known that $\mathrm{Cl}(\widetilde{U})\cong\mathbb{Z}/dn\mathbb{Z}$ and the kernel of the natural surjective map 
\[
\mathrm{Cl}(U)\to\mathrm{Cl}(\widetilde{U})
\]
is the cyclic subgroup $\langle K_U\rangle$ of order $n$.
This confirms the claim.

Finally, we complete the proof establishing the remaining claim.
For any $n>0$, we consider the canonical closed immersion $\mathrm{Spec}(\O_{C,0}/\mathfrak{m}^{n}) \hookrightarrow C$ where $\mathfrak{m}$ is the maximal ideal.
Define $\mathscr{X}_n:=\mathscr{X}\times_C\mathrm{Spec}(\O_{C,0}/\mathfrak{m}^{n})$. 
We now prove by induction that for any $n>0$, there exists a flat coherent sheaf $\mathcal{L}_n$ on $\mathscr{X}_n$ such that $\mathcal{L}_{n}|_{\mathscr{X}_{n-1}}\cong \mathcal{L}_{n-1}$, where we set $\mathcal{L}_{0}:=L$. 
Assuming the existence of such $\mathcal{L}_n$,
we first analyze the obstruction space for extending $\mathcal{L}_n$.
Since $\mathcal{L}_n$ is flat and $\mathcal{L}_n|_{X}\cong L$ is locally free on codimension one points of $X$,
the sheaves $\mathcal{E}xt^j_{\mathscr{X}_n}(\mathcal{L}_n,L)$ for $j>0$ are supported on finitely many singularities on $X$.
Consequently, 
\begin{equation}
H^i(\mathscr{X}_n,\mathcal{E}xt^j_{\mathscr{X}_n}(\mathcal{L}_n,L))=0\label{eq--local--obstruction-1}
\end{equation} 
for all $i,j>0$.
By \cite[Tag 08L8]{Stacks}, the obstruction space to deform $\mathcal{L}_n$ to a flat coherent sheaf on $\mathscr{X}_{n+1}$ is given by
    \[
    \mathrm{Ext}^2_{\mathscr{X}_n}(\mathcal{L}_{n},L).
    \]
 On the other hand, since $\O_X(D)$ is $S_2$, we have 
 $$
 \mathcal{H}om_{\mathscr{X}_n}(\mathcal{L}_n,L)=\mathcal{H}om_{X}(L,L)=\mathcal{O}_{X}.
 $$
Therefore,
   \begin{equation}
    H^2(\mathscr{X}_n,\mathcal{H}om(\mathcal{L}_n,L))=H^2(X,\mathcal{O}_{X})=0\label{eq--local--obstruction-2}
    \end{equation}
by assumption.
Using the local-to-global spectral sequence along with \eqref{eq--local--obstruction-1} and \eqref{eq--local--obstruction-2}, 
we obtain that
   \[
   \mathrm{Ext}^2_{\mathscr{X}_n}(\mathcal{L}_{n},L)\cong H^0(\mathscr{X}_n,\mathcal{E}xt^2_{\mathscr{X}_n}(\mathcal{L}_n,L)).
   \]
This shows that all obstructions are local and supported on finitely many closed points.
By using the local-to-global spectral sequence again, the natural morphism
\[
   \mathrm{Ext}^1_{\mathscr{X}_n}(\mathcal{L}_{n},L)\to H^0(\mathscr{X}_n,\mathcal{E}xt^1_{\mathscr{X}_n}(\mathcal{L}_n,L))
   \]
   is surjective.
As a result, the natural map from the global deformation functor of $L$ to the product of the local deformation functors of $L$ at the singularities of $X$ is smooth (cf.~\cite[Lemma 1]{Manetti}). 
By the assumption of $f$, we obtain the desired sheaf $\mathcal{L}_{n+1}$.
 Then, we can apply Grothendieck's existence theorem \cite[Tag 0CTK]{Stacks} to deduce that there exists a divisorial sheaf $\widehat{\mathscr{L}}$ on $\widehat{\mathscr{X}}:=\mathscr{X}\times_{C}\mathrm{Spec}\,(\widehat{\mathcal{O}}_{C,0})$ such that $\widehat{\mathscr{L}}$ is flat over $\mathrm{Spec}\,(\widehat{\mathcal{O}}_{C,0})$ and $\widehat{\mathscr{L}}|_{\mathscr{X}_n}\cong\mathcal{L}_n$ for any $n\in\mathbb{Z}_{\ge0}$. By Popescu's theorem \cite[Tag 07GC]{Stacks}, there exist a smooth morphism $\varphi\colon B\to C$ from a variety, a coherent sheaf $\mathscr{L}_B$ on $\mathscr{X}_B:=\mathscr{X}\times_CB$ flat over $B$ and a $C$-morphism $\xi\colon \mathrm{Spec}\,(\widehat{\mathcal{O}}_{C,0})\to B$ such that $\xi^*\mathscr{L}_B\cong \widehat{\mathscr{L}}$. We may assume that $\xi$ is dominant. Take a morphism $f'\colon C'\to B$ from a smooth affine curve such that the image of $f'$ passes through a general point of $B$. Let $\mathscr{L}:=f'^*\mathscr{L}_B$. Then, we see that $\varphi\circ f'\colon C'\to C$ is dominant and $\mathscr{L}$ is a divisorial sheaf on $\mathscr{X}\times_CC'$ such that $\mathscr{L}|_{X}\cong L$, which completes the proof.
\end{proof}

Using Proposition \ref{prop--smoothability--criterion}, we deduce the following generalization of Manetti's result in a more algebro-geometric way, which is of independent interest.

\begin{cor}[cf.~{\cite[Theorem 15 (2)]{Manetti}}]
    Let $f\colon \mathscr{X}\to C$ be a $\mathbb{Q}$-Gorenstein deformation and fix a closed point $0\in C$ such that $X:=\mathscr{X}_0$ is a surface.
    Suppose that there exist a divisorial sheaf $\mathscr{L}$ and a non-zero integer $m \in\mathbb{Z}$ such that $K_\mathscr{X}\sim_C m\mathscr{L}$.
Let $x\in X$ be a T-singularity of type $\frac{1}{dn^2}(1,adn-1)$ that is smoothed by $f$.
Then, the Cartier index $n$ is coprime to $m$.
\end{cor}

\begin{proof}
    By shrinking $\mathscr{X}$ around $x$, we may assume that $\mathscr{X}$ has a unique singularity $x$ and  $\mathrm{Cl}(X)\cong\mathbb{Z}/dn^2\mathbb{Z}$.
    We know that the cyclic subgroup $\langle K_{X}\rangle$ of $\mathrm{Cl}(X)$ has order $n$.
    Applying Proposition~\ref{prop--smoothability--criterion} to $D:=\mathscr{L}|_{X}$, there exists an integer $l$ such that $\mathscr{L}|_{X}\sim lK_X$.
    Hence, we obtain that
    $$
    K_X\sim m\mathscr{L}|_{X}\sim ml K_X,
    $$
   which means that $ml \equiv 1\ (\textrm{mod $n$})$.
   Therefore, we have the assertion.
\end{proof}

\subsection{$\mathbb{Q}$-Gorenstein deformation and involutions}\label{subsec:invol-deform-8}

In this subsection, we introduce fundamental tools for constructing log $\mathbb{Q}$-Gorenstein deformations, which are essential to address the $\mathbb{Q}$-Gorenstein deformation theory of Horikawa surfaces.

First, to study $\mathbb{Q}$-Gorenstein smoothability of stable Horikawa surfaces with good involutions, it is natural to consider log $\mathbb{Q}$-Gorenstein deformations of log pairs, rather than working directly with the deformation theory of the surfaces themselves.
The following is an immediate corollary of Proposition~\ref{prop--involution}:

\begin{cor}\label{cor--too--trivial}
    Let $X$ be a stable $\mathbb{Q}$-Gorenstein smoothable normal surface, and let $f\colon \mathscr{X}\to C$ be a $\mathbb{Q}$-Gorenstein smoothing over a smooth affine curve with a closed point $0\in C$ such that $\mathscr{X}_0=X$.
    Suppose that $\mathscr{X}_c$ is a stable Horikawa surface with only Du Val singularities for any $c\in C\setminus\{0\}$.
Let $\pi\colon\mathscr{X}\to \mathscr{W}$ be the quotient morphism by the involution $\sigma_{\mathscr{X}}$ defined as in Proposition~\ref{prop--involution}, and let $\mathscr{B}$ be the branch divisor with respect to $\pi$.
Then, the pair $(\mathscr{W},\frac{1}{2}\mathscr{B}+\mathscr{W}_0)$ is lc and 
the restriction $\mathscr{X}_c\to \mathscr{W}_c$ on each fiber for $c\in C\setminus\{0\}$ is a double cover branched along 
$\mathscr{B}_c$ which is induced by the linear system $|\omega_{\mathscr{X}_c}|$.
\end{cor}

Conversely, the following proposition is useful to construct a partial $\mathbb{Q}$-Gorenstein smoothing of stable Horikawa surfaces with good involutions.

\begin{prop}\label{prop--smoothing--quotient}
Let $W$ be a projective klt surface and $B$ an effective reduced Weil divisor on $W$. 
Consider a set of T-singularities $x_1,\ldots,x_l$ on $W$.
Suppose that the following conditions hold:
\begin{itemize}
    \item[$(\mathrm{a})$] There exists a Weil divisorial sheaf $L$ on $W$ such that $L^{[2]}\sim B$,
    \item[$(\mathrm{b})$] For each $i$, there exists $l_i\in\mathbb{Z}$ such that $L\sim l_iK_W$ locally around $x_i$,
    \item[$(\mathrm{c})$] $-K_W$ is big,
    \item[$(\mathrm{d})$] $H^1(W,\mathcal{O}_W(B))=0$,
    \item[$(\mathrm{e})$] $(W,\frac{1}{2}B)$ is lc,
    \item[$(\mathrm{f})$] $K_W+\frac{1}{2}B$ is big and nef.
\end{itemize} 
    Then, there exist a flat and projective family $f\colon \mathscr{W}\to C$ over a smooth affine curve with a closed point $0\in C$ and an effective $\mathbb{Q}$-Cartier relative Mumford divisor $\mathscr{B}$ on $\mathscr{W}$ such that
    \begin{itemize}
        \item[$(1)$] $\mathscr{W}$ is a normal $\mathbb{Q}$-Gorenstein threefold,
        \item[$(2)$]  $f$ induces a smoothing of each T-singularity $x_i$ and a locally trivial deformation for any singular point of $W$ other than $x_i$,
        \item[$(3)$]  $(\mathscr{W}_0,\mathscr{B}_0)\cong(W,B)$,
        \item[$(4)$] $\mathscr{B}_c$ is reduced for any $c\in C\setminus\{0\}$, 
        \item[$(5)$] $K_{\mathscr{W}}+\frac{1}{2}\mathscr{B}$ is $f$-big and $f$-semiample,
        \item[$(6)$] There exists a divisorial sheaf $\mathscr{L}$ on $\mathscr{W}$ such that $\mathscr{L}^{[2]}\sim \mathscr{B}$ and $\mathscr{L}|_{W}\cong L$.
    \end{itemize}
    Furthermore, let $\mathscr{X}:=\mathbf{Spec}_{\mathscr{W}}(\mathcal{O}_{\mathscr{W}}\oplus \mathscr{L}^{[-1]})$ (resp.\ $X:=\mathbf{Spec}_{W}(\mathcal{O}_W\oplus L^{[-1]})$) be the double cover of $\mathscr{W}$ branched along $\mathscr{B}$ (resp.\ the double cover of $W$ branched along $B$). 
    Then, $\mathscr{X}$ is normal and $\mathbb{Q}$-Gorenstein. 
    The canonical divisor $K_{\mathscr{X}}$ is relatively big and semiample over $C$, and the central fiber satisfies $\mathscr{X}_0\cong X$.
Moreover, if we define $\overline{\mathscr{X}}$ as the relative canonical model of $\mathscr{X}$ over $C$, then the central fiber $\overline{\mathscr{X}}_0$ coincides with the canonical model of $X$, and $\overline{\mathscr{X}}$ is also $\mathbb{Q}$-Gorenstein.    
\end{prop}

\begin{proof}
First, we construct a flat and projective morphism $f\colon \mathscr{W}\to C$, where $\mathscr{W}$ is a normal threefold and $C$ is a smooth affine curve with a closed point $0\in C$, satisfying conditions (1) and (2).
Since $-K_W$ is big and $W$ has only rational singularities, it follows that $H^2(W,\mathcal{O}_W)=0$ (cf.\ Remark~\ref{rem--smoothing--quotient}~(2) below).
Moreover, because $-K_W$ is big and $W$ is lc, it follows from \cite[Proposition 3.1]{HP} and Grothendieck's algebraization theorem \cite[Tag 089A]{Stacks} that there exists a flat and projective morphism 
$$
\hat{f}\colon \widehat{\mathscr{W}}\to R:=\mathrm{Spec}\,\mathbb{C}[[t]]
$$
from a normal scheme such that
\begin{itemize}
    \item $\widehat{\mathscr{W}}_0\cong W$ for the closed point $0\in R$.
    \item For each $n$, the restriction $\widehat{\mathscr{W}}\times_{R}R_n$ is $\mathbb{Q}$-Gorenstein, where $R_n:=\mathrm{Spec}(\mathbb{C}[[t]]/(t^n))$.
    \item $\hat{f}$ induces a smoothing of each T-singularity $x_i$ and a locally trivial deformation at any singular point of $W$ other than $x_i$.    
\end{itemize}
Let $m$ be a positive integer such that $mK_W$ is Cartier.
Then, by Grothendieck's existence theorem \cite[Tag 088E]{Stacks} and \cite[Lemma 3.16]{KSB}, there exists a line bundle $\widehat{\mathscr{M}}$ on $\widehat{\mathscr{W}}$ such that 
$$
\widehat{\mathscr{M}}|_{\widehat{\mathscr{W}}\times_{R}R_n}\cong \omega^{[m]}_{\widehat{\mathscr{W}}\times_{R}R_n/R_n}.
$$
Furthermore, there exists a morphism $\omega^{\otimes m}_{\widehat{\mathscr{W}}/R}\to\widehat{\mathscr{M}}$ which  is an isomorphism in codimension one around $W$ by \cite[Tag 088E]{Stacks}.
Hence, the $S_2$-hull of $\omega^{\otimes m}_{\widehat{\mathscr{W}}/R}$ coincides with $\widehat{\mathscr{M}}$ and so $\widehat{\mathscr{W}}$ is also $\mathbb{Q}$-Gorenstein.
By \cite[Tag 07GC]{Stacks}, we have
$$
\varinjlim_{\lambda}A_\lambda=\mathbb{C}[[t]],
$$
where $\{A_\lambda\}_{\{\lambda\in\Lambda\}}$ is a cofiltered direct system of finite-type $\mathbb{C}$-algebras such that $D_{\lambda}:=\mathrm{Spec}\,(A_{\lambda})$ is smooth for any $\lambda\in\Lambda$ (cf.~\cite[Tag 04AZ]{Stacks}).
Since $\hat{f}$ is of finite type, there exist $\lambda_0\in \Lambda$ and a flat projective morphism 
$$
\hat{f}_{\lambda_0}\colon\widehat{\mathscr{W}}_{\lambda_0}\to D_{\lambda_0}
$$
from a $\mathbb{Q}$-Gorenstein normal scheme such that 
$$
\widehat{\mathscr{W}}_{\lambda_0}\times_{D_{\lambda_0}}\mathrm{Spec}\,(\mathbb{C}[[t]])\cong \widehat{\mathscr{W}},
$$
and $\hat{f}_{\lambda_0}$ induces a smoothing of each T-singularity $x_i$ 
at the closed point $0\in D_{\lambda_0}$
(cf.~\cite[1.2.13, 1.2.14]{Ols}, \cite[Tag 081C]{Stacks}).
Moreover, since $\hat{f}$ induces a locally trivial deformation at any singular point $w$ of $W$ other than $x_i$, by the local criterion for flatness, we can choose an affine open neighborhood $U_w\subset W$ such that $\widehat{W}$ has an affine open neighborhood of $w$ of the form $U_{w}\times R$. 
Hence, after possibly replacing $\lambda_0$, we may assume that $\hat{f}_{\lambda_0}$ also induces a locally trivial deformation at any singular point of $W$ other than $x_i$.
Now, choose a morphism $\varphi\colon C\to D_{\lambda_0}$ from a smooth affine curve $C$ such that $\varphi$ maps a closed point $0\in C$ to the closed point $0\in D_{\lambda_0}$ and $\varphi(C)$ meets a general closed point of $D_{\lambda_0}$.
Then define 
$$
\mathscr{W}:=\widehat{\mathscr{W}}_{\lambda_0}\times_{D_{\lambda_0}}C.
$$ 
It is straightforward to verify that the resulting morphism $f\colon \mathscr{W}\to C$ satisfies conditions (1) and (2).
    
    By the condition (b) and Proposition \ref{prop--smoothability--criterion}, there exists a $\mathbb{Q}$-Cartier divisorial sheaf $\mathscr{L}$ on $\mathscr{W}$ such that $\mathscr{L}|_{\mathscr{W}_0}\cong L$ after taking a finite base change over $C$.
    Furthermore, by the condition (d) and \cite[III, Theorem~12.11]{Ha}, the natural homomorphism 
    \[
    f_*\mathscr{L}^{[2]}\otimes \O_{C,0}/\mathfrak{m} \to H^0(W,\mathscr{L}^{[2]}|_{W})\cong H^0(W,\mathcal{O}_W(B))
    \]
    is an isomorphism. 
    Thus, there exists an effective relative Mumford divisor $\mathscr{B}\sim \mathscr{L}^{[2]}$ such that $\mathscr{B}_0=B$.
    Since $B$ is reduced, so is $\mathscr{B}_c$  for any $c\in C$ after shrinking $C$ if necessary.
    The inversion of adjunction for log canonicity \cite{kawakita} implies that $(\mathscr{W},\frac{1}{2}\mathscr{B}+\mathscr{W}_0)$ is lc and $(\mathscr{W},\mathscr{W}_0)$ is plt.
    Since $\mathscr{W}$ is a klt threefold and projective over $C$, we may run a relative $(K_\mathscr{W}+\frac{1}{2}\mathscr{B}+\mathscr{W}_0)$-MMP over $C$ and obtain a relative minimal model $h\colon \mathscr{W}\dashrightarrow\mathscr{W}'$ (cf.~\cite{Koetal}).
    Note that $K_{\mathscr{W}'}+h_*(\frac{1}{2}\mathscr{B}+\mathscr{W}_0)$ is relatively semiample over $C$ (cf.~\cite{kmm}, \cite[Corollary 4.7.8]{fujino-foundation}).
    Since $K_{W}+\frac{1}{2}B$ is nef,
    $h$ is isomorphic along the central fiber.
    Hence, we may assume that $\mathscr{W}'\cong \mathscr{W}$ after shrinking $C$ around $0$.
    Therefore, $K_{\mathscr{W}}+\frac{1}{2}\mathscr{B}$ is $f$-semiample.
    Note that the natural homormorphism
    $$
    f_*\mathcal{O}_{\mathscr{W}}(l(K_{\mathscr{W}/C}+\mathscr{L}))\otimes \O_{C,c}/\mathfrak{m}\to H^0(\mathscr{W}_c,\mathcal{O}_{\mathscr{W}_c}(l(K_{\mathscr{W}_c}+\mathscr{L}_c)))
    $$
    is an isomorphism for any sufficiently large and divisible $l\in\mathbb{Z}_{>0}$ by \cite[Lemma 4.2]{MZ} and Grauert's theorem \cite[III, Corollary~12.9]{Ha}. 
    Hence the bigness of $K_{W}+\frac{1}{2}B$ implies that $K_{\mathscr{W}}+\frac{1}{2}\mathscr{B}$ is $f$-big.
    This completes the proof of the first half of the claim.
    
    Let $\pi\colon \mathscr{X}\to \mathscr{W}$ be the double cover branch along $\mathscr{B}$. 
    Since $\mathscr{B}$ is reduced, $\mathscr{X}$ is normal.
    By the ramification formula
    \[
   K_{\mathscr{X}}+\mathscr{X}_0=\pi^*\left(K_{\mathscr{W}}+\frac{1}{2}\mathscr{B}+\mathscr{W}_0\right),
   \]
  the pair $(\mathscr{X},\mathscr{X}_0)$ is lc.
 Therefore, $K_{\mathscr{X}}$ is a $g$-big and $g$-semiample $\mathbb{Q}$-Cartier divisor, where $g\colon \mathscr{X}\to C$ is the natural morphism.
 Let 
 $$
 \overline{\mathscr{X}}:=\mathbf{Proj}_{C}(\bigoplus_{l\ge0}g_*\mathcal{O}_{\mathscr{X}}(lK_{\mathscr{X}}))
 $$
 denote the relative canonical model of $\mathscr{X}$.
 By applying \cite[Lemma 4.2]{MZ} again, 
 it follows that the central fiber $\overline{\mathscr{X}}_0$
 coincides with the canonical model of $X$.
 Hence we complete the proof.
    \end{proof}

\begin{rem}\label{rem--smoothing--quotient}
\begin{itemize}
  \item[(1)]  Let $f\colon\mathscr{W}\to C$ and $\mathscr{B}$ be as in Proposition~\ref{prop--smoothing--quotient}.
    For any $c\in C\setminus\{0\}$, the pair $(\mathscr{W}_c,\frac{1}{2}\mathscr{B}_c)$ satisfies all conditions (a)--(f) in Proposition~\ref{prop--smoothing--quotient} except that $-K_{\mathscr{W}_c}$ is big from \cite[III, Theorem 12.11]{Ha} and the proof of Proposition~\ref{prop--smoothing--quotient}.
    \item[(2)] For any pair $(W,\frac{1}{2}B)$ satisfying the all conditions (a)--(f) in Proposition~\ref{prop--smoothing--quotient},  $W$ is a rational surface.
    Indeed, let $\mu\colon\tilde{W}\to W$ be the minimal resolution.
    Then, there exists an effective $\mu$-exceptional $\mathbb{Q}$-divisor $\Delta$ such that $\lfloor\Delta\rfloor=0$ and 
    \[
    K_{\tilde{W}}+\Delta=\mu^*K_W.
    \]
    Thus, $-K_{\tilde{W}}$ is big and hence $H^0(\tilde{W},\mathcal{O}_{\tilde{W}}(2K_{\tilde{W}}))=0$.
    On the other hand, $W$ is klt and hence has only rational singularities by \cite[Theorem 5.22]{KoMo}.
    Therefore, it is easy to see that $H^1(\tilde{W},\mathcal{O}_{\tilde{W}})=0$.
    By Castelnuovo's rationality criterion \cite[V, Theorem 6.2]{Ha}, $\tilde{W}$ is rational.
\item[(3)] Let $\iota$ be the involution of the 
double cover $\mathscr{X}\to\mathscr{W}$.
Then it is straightforward to see that $\iota$ induces an involution $\overline{\iota}$ on $\overline{\mathscr{X}}$ over $C$
and the quotient $\overline{\mathscr{X}}/\overline{\iota}$
is the relative log canonical model of 
$(\mathscr{W},\frac{1}{2}\mathscr{B}+\mathscr{W}_0)$ over $C$.
It is also easy to see by \cite[Corollary 2.69]{kollar-moduli} that if $\overline{\iota}_0$ acts on $H^0(\overline{\mathscr{X}}_0,K_{\overline{\mathscr{X}}_0})$ trivially, then $\overline{\iota}_c$ is also trivial for any $c\in C$.
\end{itemize}
\end{rem}

We also put the following useful criterion for log $\mathbb{Q}$-Gorenstein smoothability of slc surface pairs.
\begin{prop}\label{prop--smoothing--slc}
Let $W$ be a projective irreducible slc surface and $B$ an effective reduced Weil divisor on $W$. 
Suppose that the codimension one part $Z$ of  $\mathrm{Sing}(W)$ is isomorphic to $\mathbb{P}^1$, and that the germ $(Z\subset W)$ is $\mathbb{Q}$-Gorenstein smoothable.
Suppose also that the following conditions hold.

\begin{itemize}
\item[$(\mathrm{a})$] 
There exists a divisorial sheaf $L$ on $W$ such that $L^{[2]}\sim B$,
\item[$(\mathrm{b})$]
$K_{W}+L$ is Cartier around $Z$,
\item[$(\mathrm{c})$]
$-K_W$ is log big and nef,
\item[$(\mathrm{d})$]
$W\setminus Z$ is klt,
\item[$(\mathrm{e})$]
$(W,\frac{1}{2}B)$ is slc,
\item[$(\mathrm{f})$]
$K_W+\frac{1}{2}B$ is big and nef,
\end{itemize}

Then, there exists a flat and projective family $f\colon \mathscr{W}\to C$ over a smooth affine curve with a closed point $0\in C$ and an effective $\mathbb{Q}$-Cartier relative Mumford divisor $\mathscr{B}$ on $\mathscr{W}$ such that
\begin{itemize}
\item[$(1)$]
$\mathscr{W}$ is a normal $\mathbb{Q}$-Gorenstein threefold,
\item[$(2)$]
$f$ induces a smoothing locally around $Z\subset W$ and a locally trivial deformation at any $w\in \mathrm{Sing}(W)\setminus Z$,
\item[$(3)$]
$(\mathscr{W}_0,\mathscr{B}_0)\cong(W,B)$,
\item[$(4)$]
$\mathscr{B}_c$ is reduced for any $c\in C\setminus\{0\}$, 
\item[$(5)$]
$K_{\mathscr{W}}+\frac{1}{2}\mathscr{B}$ and $-K_{\mathscr{W}}$ are $f$-big and $f$-semiample,
\item[$(6)$]
There exists a divisorial sheaf $\mathscr{L}$ on $\mathscr{W}$ such that $\mathscr{L}^{[2]}\sim \mathscr{B}$ and $\mathscr{L}|_{\mathscr{W}_0}\cong L$.
\end{itemize}
Furthermore, let $\mathscr{X}:=\mathbf{Spec}_{\mathscr{W}}(\mathcal{O}_{\mathscr{W}}\oplus \mathscr{L}^{[-1]})$ (resp.\ $X:=\mathbf{Spec}_{W}(\mathcal{O}_W\oplus L^{[-1]})$) be the double cover of $\mathscr{W}$ branched along $\mathscr{B}$ (resp.\ the double cover of $W$ branched along $B$). 
    Then, $\mathscr{X}$ is normal and $\mathbb{Q}$-Gorenstein. 
    The canonical divisor $K_{\mathscr{X}}$ is relatively big and semiample over $C$, and the central fiber satisfies $\mathscr{X}_0\cong X$.
Moreover, if we define $\overline{\mathscr{X}}$ as the relative canonical model of $\mathscr{X}$ over $C$, then the central fiber $\overline{\mathscr{X}}_0$ coincides with the canonical model of $X$, and $\overline{\mathscr{X}}$ is also $\mathbb{Q}$-Gorenstein. 
\end{prop}

\begin{proof}
By \cite[Theorem 1.10]{fujino--slc--vanishing}, we have \[H^1(W,\mathcal{O}_W(B))=H^1(W,\mathcal{O}_W)=H^2(W,\mathcal{O}_W)=0.\]
In addition, Theorem \ref{thm--hacking--prokhorov} implies that $H^2(W,T_W)=0$.
We now construct a flat and projective contraction $f\colon \mathscr{W}\to C$ from a klt threefold $\mathscr{W}$ to a smooth affine curve $C$ with a closed point $0\in C$, together with a $\mathbb{Q}$-Cartier divisorial sheaf $\mathscr{L}$, such that:

\begin{itemize}
\item $\mathscr{W}_0\cong W$,
\item $\mathscr{L}|_W=L$,
\item $f$ is a $\Q$-Gorenstein smoothing along $Z$,
\item $f$ is locally trivial over $W\setminus Z$. 
\end{itemize}

There is an exact sequence of sheaves:
\[
0\to \mathscr{T}^1_{QG,\mathrm{tor}}\to  \mathscr{T}^1_{QG}\to  \mathscr{T}^1_{QG,\mathrm{tf}}\to 0,
\]
where $\mathscr{T}^1_{QG,\mathrm{tor}}$ is the subsheaf of $\mathscr{T}^1_{QG}$ consisting of sections supported in dimension zero.
The quotient $\mathscr{T}^1_{QG,\mathrm{tf}}$ is a line bundle on $Z\cong\mathbb{P}^1$.
 Since the germ $(Z\subset W)$ is $\mathbb{Q}$-Gorenstein smoothable, $\mathscr{T}^1_{QG,\mathrm{tf}}$ admits a nonzero global section. 
 Hence, we have $\mathscr{T}^1_{QG,\mathrm{tf}}\cong\mathcal{O}_{\mathbb{P}^1}(m)$ for some $m\ge0$.
 In particular, this implies 
 $$
 H^1(\mathscr{T}^1_{QG})=H^1(\mathscr{T}^1_{QG,\mathrm{tf}})=0.
 $$
 By the spectral sequence \eqref{eq--T^i--spectral}, the natural map $T^2_{QG}\to H^0(\mathscr{T}^2_{QG})$ is an isomorphism and the map $T^1_{QG}\to H^0(\mathscr{T}^1_{QG})$ is surjective.
Using this and the argument in the first paragraph of the proof of Proposition \ref{prop--smoothing--quotient}, we can conatruct a flat and projective morphism $f\colon \mathscr{W}\to C$ from a threefold to a smooth affine curve with a closed point $0\in C$ such that 
$\mathscr{W}_0\cong W$,
the pair $(\mathscr{W},\mathscr{W}_0)$ is slc, and $f$ is a $\mathbb{Q}$-Gorenstein smoothing of the germ $(Z\subset W)$, while $f$ is locally trivial over $W\setminus Z$.
 Since by the assumption (d), the generic fiber of $f$ is klt, we conclude that $\mathscr{W}$ is klt.
 Now, since $K_W+L$ is Cartier around $Z$ and $f$ is locally trivial on $W\setminus Z$, the vanishing $H^2(W,\mathcal{O}_W)=0$ and the same argument as in the proof of Proposition \ref{prop--smoothability--criterion} ensure the existence of a line bundle $\mathscr{H}$ on $\mathscr{W}$ such that 
 $$
 \mathscr{H}_0\sim K_W+L.
 $$
We then define 
$$
\mathscr{L}:=\mathscr{H}\otimes \omega_{\mathscr{W}/C}^{[-1]}.
$$
By \cite[Corollary 5.25]{KoMo}, $\mathscr{L}$ is a Cohen-Macaulay sheaf, and in paticular, $\mathscr{L}|_W\cong L$.
The remaining part of the argument follows as in the proof of Proposition \ref{prop--smoothing--quotient}.
\end{proof}

\subsection{$\mathbb{Q}$-Gorenstein smoothing of rational surfaces}\label{subsec:rational-surf-deform-8}

By Proposition \ref{prop--smoothing--quotient}, the problem of finding a $\mathbb{Q}$-Gorenstein smoothing for a stable Horikawa surface $X$ with a good involution $\sigma$ reduces to constructing a log $\mathbb{Q}$-Gorenstein partial smoothing of $(W,\frac{1}{2}B)$, where $W:=X/\sigma$ and $B$ is the branch divisor.
We summarize the key properties of $W$ and $B$ under the assumption that $X$ is $\mathbb{Q}$-Gorenstein smoothable.
\begin{enumerate}
\item $h^0(-K_W)\ge9$. 
\item $-K_W$ is big.
\item $H^1(W,\mathcal{O}_W)=H^2(W,\mathcal{O}_W)=0$ (see \cite[Corollary 2.64]{kollar-moduli}).
\end{enumerate}
Properties (1) and (2) follows from the following general lemma, given that $W$ deforms to either a Hirzebruch surface or $\mathbb{P}^2$ (cf.~Corollary \ref{cor--too--trivial}).
\begin{lem}\label{lem--upper--semiconti--plurigenera}
Let $f\colon\mathscr{W}\to C$ be a projective flat morphism from a normal variety to a smooth affine curve with a closed point $0\in C$, and set $W:=\mathscr{W}_0$.
Fix an integer $r\in\mathbb{Z}$.
Then, for any general point $c\in C$, we have
\[
h^0(W,\omega_{W}^{[r]})\ge h^0(\mathscr{W}_c,\omega_{\mathscr{W}_c}^{[r]}).
\]
In particular, if $\omega_{W}^{[r]}$ is $\mathbb{Q}$-Cartier and $\omega_{\mathscr{W}_c}^{[r]}$ is big for any general $c\in C$, then so is $\omega_{W}^{[r]}$.
\end{lem}
\begin{proof}
Since $\omega_{\mathscr{W}/C}^{[r]}$ satisfies the $S_2$-condition,  the restriction $\omega_{\mathscr{W}/C}^{[r]}|_{W}$ satisfies the $S_1$-condition.
Hence, the natural morphism 
$$
\omega_{\mathscr{W}/C}^{[r]}|_{W}\to \omega_{W}^{[r]}
$$
is injective for any $r\in\mathbb{Z}$.
Furthermore, since $\omega_{\mathscr{W}/C}^{[r]}$ is flat over $C$, we have
$\omega_{\mathscr{W}/C}^{[r]}|_{\mathscr{W}_c}\cong \omega^{[r]}_{\mathscr{W}_c}$ for general $c\in C$ (cf.~\cite[Corollary 10.12]{kollar-moduli}).
Thus, by \cite[III, Theorem 12.8]{Ha}, we conclude that 
\[
h^0(W,\omega_{W}^{[r]})\ge h^0(W,\omega_{\mathscr{W}/C}^{[r]}|_W)\ge h^0(\mathscr{W}_c,\omega_{\mathscr{W}/C}^{[r]}|_{\mathscr{W}_c})=
h^0(\mathscr{W}_c,\omega_{\mathscr{W}_c}^{[r]}).
\]
The last assertion follows from the inequality above.
\end{proof}

\begin{cor}\label{cor--deformation--rational}
    Let $f\colon\mathscr{W}\to C$ be a smoothing of a normal rational surface $W=\mathscr{W}_0$ with only rational singularities. 
    Then any general fiber 
    $\mathscr{W}_c$ is also rational. 
\end{cor}

\begin{proof}
    Since $W$ is a rational surface with only rational singularities, it follows that $H^0(W,\omega_{W}^{[2]})=0$ and $H^{1}(W,\O_W)=0$.
    Thus, by Lemma~\ref{lem--upper--semiconti--plurigenera} and \cite[III, Theorem~12.8]{Ha}, we obtain $H^0(\mathscr{W}_c,\omega_{\mathscr{W}_c}^{\otimes 2})=0$ and $H^{1}(\mathscr{W}_c,\O_{\mathscr{W}_c})=0$.
    Applying Castelnuovo's rationality criterion \cite[V, Theorem 6.2]{Ha}, we conclude that $\mathscr{W}_c$ is rational.
\end{proof}

The following lemma provides a useful criterion for identifying the location of curves with negative self-intersection on a given rational surface.

\begin{lem}\label{lem--modified--manetti--lemma}
Let $W$ be a projective smooth rational surface
equipped with a ruling $p\colon W\to \mathbb{P}^1$.
Fix an integer $d\ge1$ such that $h^0(W,\mathcal{O}(-K_{W}))+\min\{d,3\}\ge 8$.
Let $C$ be a section of $p$ such that $C^{2}=-d$.
Then, any irreducible curve $D$ on $W$ with $D^{2}\le -2$ other than $C$ is vertical with respect to $p$.

Furthermore, if $h^0(W,\mathcal{O}(-K_{W}))+\min\{d,3\}\ge 9$, then any irreducible curve $D$ on $W$ with $D^2\le -1$ other than $C$ is vertical with respect to $p$.
\end{lem}

\begin{proof}
The first assertion follows from the same argument as in the proof of \cite[Proposition 9]{Manetti}.

For the second assertion, assume that
there exists a horizontal curve $D\subset W$ with $D^2\le -1$.
Consider the blow-up $\pi\colon \hat{W}\to W$ at a point on $D$.
Then, the proper transform $\hat{D}:=\pi_*^{-1}(D)$ satisfies $\hat{D}^2\le-2$.
Since the proper transform $\hat{C}:=\pi_*^{-1}(C)$ satisfies $h^0(-K_{\hat{W}})+\min\{-(\hat{C})^2,3\}\ge8$, the first assertion implies $D=C$.
Thus, the second assertion follows.
\end{proof}

The following is also useful for applying vanishing theorems, including \cite[Theorem 1.3]{Enokizono}.

\begin{lem}\label{lem--nakai--moishezon} 
    Let $X$ be a normal projective surface.
Let $B$ be a $\Q$-Cartier $\Q$-divisor on $X$.
Then, $B$ is nef (resp.\ ample) if and only if it is pseudo-effective (resp.\ big) and $B\cdot C\ge 0$ (resp.\ $B\cdot C > 0$) for any irreducible curve $C\subset X$ with $C^2<0$.
\end{lem}

\begin{proof}
This follows from the Zariski decomposition theorem (cf.\ \cite[Theorem 3.1]{Enokizono}), the Nakai-Moishezon theorem \cite[Theorem 1.37]{KoMo} and the Hodge index theorem \cite[V, Theorem 1.9]{Ha}.
\end{proof}

On the other hand, 
the following is useful for analyzing the structure of a one-parameter deformation of rational surfaces.

\begin{lem}\label{lem--lower-semi--rho}
    Let $f\colon \mathscr{W}\to C$ be a projective flat morphism over a smooth affine curve with a closed point $0\in C$.
    Suppose that the central fiber $W:=\mathscr{W}_0$ is slc and irreducible but any other fiber $\mathscr{W}_c$ is a rational surface with only klt singularities.
    Then, the function 
    $$
    C\ni c\mapsto\rho(\mathscr{W}_c)\in \Z_{>0}
    $$
    is lower semicontinuous. 
Furthermore, there exists a quasi-finite map $\pi\colon U\to C$ from a smooth curve with $0\in\pi(U)$ such that $\mathscr{W}\times_CV$ is $\mathbb{Q}$-factorial for any finite morphism $V\to U$ from a smooth curve if and only if the function $U\ni u\mapsto\rho(\mathscr{W}_u)$ is constant. 
\end{lem}

\begin{proof}
By assumption and \cite[Corollary 2.64]{kollar-moduli}, we have $H^1(W,\mathcal{O}_W)=H^2(W,\mathcal{O}_W)=0$.
To prove the assertions, we will construct a morphism $\pi\colon U\to C$ from a smooth curve $U$ with $0\in \pi(U)$ satisfying the following properties:
\begin{itemize}
    \item The fiber product $\mathscr{W}\times_C(U\setminus\pi^{-1}(0))$ is $\mathbb{Q}$-factorial, and 
    \item
    For any \'etale morphism $V\to U\setminus\pi^{-1}(0)$ and any $v\in V$, the restriction map $\mathrm{Pic}(\mathscr{W}\times_CV/V)\to \mathrm{Pic}(\mathscr{W}_v)$ is an isomorphism.
\end{itemize}
By \cite[Theorem 2.15]{KSB} and \cite[Theorem 4.28]{KoMo}, there exists a morphism $\pi\colon U\to C$ from a smooth curve with $0\in\pi(U)$ such that $\mathscr{W}\times_C(U\setminus\pi^{-1}(0))$ admits a projective simultaneous resolution
\[
h\colon \widetilde{\mathscr{W}}\to \mathscr{W}\times_C(U\setminus\pi^{-1}(0)),
\]
where the restriction of $h$ to each fiber is a minimal resolution.
We claim that after possibly replacing $U$, the restriction map 
$$
\mathrm{Pic}(\widetilde{\mathscr{W}}\times_{U\setminus\pi^{-1}(0)}V/V)\to \mathrm{Pic}(\widetilde{\mathscr{W}}_v)
$$
is an isomorphism for any \'etale morphism $V\to U\setminus\pi^{-1}(0)$ and any $v\in V$.
If the general fiber $\widetilde{\mathscr{W}}_u\cong \mathbb{P}^2$, then  we may assume that $\widetilde{\mathscr{W}}\cong \mathbb{P}^2\times (U\setminus\pi^{-1}(0))$ after possibly replacing $U$.
Similarly, if the general fiber $\widetilde{\mathscr{W}}_u\cong \Sigma_{d}$ for some $d\in\mathbb{Z}_{>0}$, then we may replace $U$ so that $\widetilde{\mathscr{W}}\cong \Sigma_d\times (U\setminus\pi^{-1}(0))$.
In both cases, the desired isomorphism of Picard groups follows immediately.
If the general fiber satisfies $\rho(\widetilde{\mathscr{W}}_u)>2$, then after shrinking $U$, we may assume that there exists a projective birational morphism $\widetilde{\mathscr{W}}\to \Sigma_d\times (U\setminus\pi^{-1}(0))$ for some $d\in\mathbb{Z}_{\ge0}$, given as a sequence of blow-ups along sections over $U\setminus\pi^{-1}(0)$, obtained  via the relative minimal model program.
In this case as well, we obtain the desired isomorphism of Picard groups.

Now, consider the equality
\[
K_{\widetilde{\mathscr{W}}}+\Delta=h^*(K_{\mathscr{W}\times_C(U\setminus\pi^{-1}(0))})+F,
\]
where $\Delta$ is an $h$-exceptional $\mathbb{Q}$-divisor.
Take a sufficiently small $\varepsilon\in\mathbb{Q}_{>0}$ such that the pair $(\widetilde{\mathscr{W}},\Delta+\varepsilon E)$ is klt, where $E$ denotes the sum of all $h$-exceptional prime divisors.
We may then run a relative log minimal model program of $(\widetilde{\mathscr{W}},\Delta+\varepsilon E)$ over $\mathscr{W}\times_C (U\setminus\pi^{-1}(0))$ (cf.~\cite{Koetal}):
\[ \widetilde{\mathscr{W}}\dashrightarrow\mathscr{W}_{1}\dashrightarrow\ldots\dashrightarrow \mathscr{W}_i\dashrightarrow\ldots\dashrightarrow \mathscr{W}_{k}=\mathscr{W}_{\mathrm{min}}.
\]
The resulting minimal model $\mathscr{W}_{\mathrm{min}}$ is klt and isomophic  to $\mathscr{W}\times_C(U\setminus\pi^{-1}(0))$ in codimension one.
After shrinking $U$, we may assume that 
\[
\mathscr{W}_{\mathrm{min}}\cong\mathscr{W}\times_C(U\setminus\pi^{-1}(0)).
\]
Since $\mathrm{dim}\,\widetilde{\mathscr{W}}=3$, each step of this MMP can be taken to be the contraction of a prime divisor that is horizontal over $U\setminus\pi^{-1}(0)$, possibly after further shrinking $U$.
The number of such steps is bounded above by $\rho(\widetilde{\mathscr{W}}_u)$ for all $u\in U\setminus\pi^{-1}(0)$.
Replacing $U$ by a finite cover if necessary, we may assume that in each step, the contracted prime divisor has irreducible fibers over $U\setminus\pi^{-1}(0)$.
It then follows that 
$$
\mathrm{Pic}(\mathscr{W}_{\mathrm{min}}\times_{U\setminus\pi^{-1}(0)}V/V)\cong \mathrm{Pic}((\mathscr{W}_{\mathrm{min}})_v)
$$
for any \'etale morphism $V\to U\setminus\pi^{-1}(0)$ and any $v\in V$.
To see this, we may assume $U\setminus\pi^{-1}(0)=V$ and let $\xi\colon \widetilde{\mathscr{W}}\to \mathscr{W}_{\mathrm{min}}$ denote the contraction.
For a line bundle $L$ on $(\mathscr{W}_{\mathrm{min}})_v$, consider the pullback $\xi_{v}^*L$.
Using the isomorphism $\mathrm{Pic}(\widetilde{\mathscr{W}}_u)\cong \mathrm{Pic}(\widetilde{\mathscr{W}}/(U\setminus\pi^{-1}(0)))$,
we obtain a line bundle $\mathscr{L}$, which is numerically trivial along the fibers of $\xi$.
Then, by applying \cite[Theorem 3.3]{KoMo} to the pair $(\widetilde{\mathscr{W}},\Delta)$ and the line bundle $\mathscr{L}$, we obtain a line bundle $\mathscr{L}'$ on $\mathscr{W}_{\mathrm{min}}$ such that $\xi^*\mathscr{L}'=\mathscr{L}$ and $\mathscr{L}'_v=L$.
 This correspondence is inverse to the natural restriction map, completing the proof of the isomorphism.

Finally, choose line bundles $L_1,\ldots, L_{\rho(W)}$ whose classes form a basis of $\mathrm{Pic}(W)\otimes\mathbb{Q}$.
By the proof of Proposition \ref{prop--smoothability--criterion}, after a finite base change of $U$, we may assume that there exist line bundles $\mathscr{L}_j$ on $\mathscr{W}\times_CU$ such that $\mathscr{L}_j|_W=L_j$.
Since $H^1(W,\mathcal{O}_W)=0$, such $\mathscr{L}_j$ are unique up to relative linear equivalence over $U$ after shrinking $U$.
Therefore, we obtain 
\begin{equation}\label{eq--pic--loc--const}
\mathrm{Pic}(W)\otimes\mathbb{Q}\cong\mathrm{Pic}(\mathscr{W}\times_CU/U)\otimes\mathbb{Q}\hookrightarrow \mathrm{Cl}(\mathscr{W}\times_CU/U)\otimes\mathbb{Q}\cong\mathrm{Pic}(\mathscr{W}_u)\otimes\mathbb{Q}.
\end{equation}
To deduce \eqref{eq--pic--loc--const}, we use the following isomorphisms: 
\begin{align*}
&\mathrm{Cl}(\mathscr{W}\times_CU/U)\otimes\mathbb{Q}\cong\mathrm{Cl}(\mathscr{W}\times_C(U\setminus\pi^{-1}(0))/U\setminus\pi^{-1}(0))\otimes\mathbb{Q}\\
&\cong\mathrm{Pic}(\mathscr{W}\times_C(U\setminus\pi^{-1}(0))/U\setminus\pi^{-1}(0))\otimes\mathbb{Q}\cong\mathrm{Pic}(\mathscr{W}_u)\otimes\mathbb{Q},
\end{align*}
which follow from the properties of $\pi\colon U\to C$.
The inclusion \eqref{eq--pic--loc--const} completes the proof of the desired statements.
Indeed, we have $\rho(W)\le \rho(\mathscr{W}_u)$, and if $\rho(W)< \rho(\mathscr{W}_u)$, then $\mathrm{Pic}(\mathscr{W}\times_CU/U)\otimes\mathbb{Q}\subsetneq \mathrm{Cl}(\mathscr{W}\times_CU/U)\otimes\mathbb{Q}$,
showing that $\mathscr{W}\times_CU$ is not $\mathbb{Q}$-factorial.
On the other hand, for any finite morphism $V\to U$, we similarly obtain
\[\mathrm{Pic}(W)\otimes\mathbb{Q}\cong\mathrm{Pic}(\mathscr{W}\times_CV/V)\otimes\mathbb{Q}\hookrightarrow \mathrm{Cl}(\mathscr{W}\times_CV/V)\otimes\mathbb{Q}\cong\mathrm{Pic}(\mathscr{W}_u)\otimes\mathbb{Q}.\]
If $\rho(W)= \rho(\mathscr{W}_u)$, then
\[
\mathrm{Pic}(\mathscr{W}\times_CV/V)\otimes\mathbb{Q}= \mathrm{Cl}(\mathscr{W}\times_CV/V)\otimes\mathbb{Q},
\]
so that $\mathscr{W}\times_CV$ is $\mathbb{Q}$-factorial.
\end{proof}

\subsection{Anti-P-resolution}\label{subsec:anti-P-8}
A P-resolution is a useful tool for describing the irreducible components of the deformation space of quotient surface singularities.
However, in some cases, a P-resolution does not provide a log $\mathbb{Q}$-Gorenstein smoothing that respects a given boundary divisor (see Remark \ref{rem--necessity-of-P-or-anti-p} below).
To address this issue, we introduce the following concept.

\begin{defn}[Anti-P-resolution]
\label{defn--anti_P}
Let $W$ be a klt surface.
Consider a projective birational morphism  $\mu\colon W^{-} \to W$ from a reduced (but not necessarily normal) surface.
We say that $\mu\colon W^{-} \to W$ is an {\it anti-P-resolution} if the following conditions hold:
 \begin{enumerate}
 \item $W^{-}$ has only $\mathbb{Q}$-Gorenstein smoothable slc singularities, and
     \item $-K_{W^{-}}$ is $\mu$-ample.
 \end{enumerate}

In this paper, we frequently consider the following situation:
 Let $B$ be an effective divisor on $W$ and let $B^{-}$ be an effective divisor on $W^{-}$ such that $\mu_{*}B^{-}=B$.
 Then, $\mu\colon (W^{-},\frac{1}{2}B^{-})\to (W,\frac{1}{2}B)$ is an {\it anti-P-resolution that is log crepant with respect to} $(W,\frac{1}{2}B)$ if $\mu$ is an anti-P-resolution and $K_{W^{-}}+\frac{1}{2}B^{-}=\mu^*\left(K_W+\frac{1}{2}B\right)$.

\end{defn}

An anti-P-resolution arises naturally when considering a log $\mathbb{Q}$-Gorenstein smoothing of a klt singularity with a given boundary divisor, as we shall see below. 

\begin{prop}\label{prop--anti-P}
    Let $f\colon\mathscr{W}\to C$ be a flat morphism from a normal variety to a smooth curve.
    Fix a point $0\in C$ where $\mathscr{W}_0$ is a klt surface.
    Suppose that there exists a relative Mumford $\mathbb{Q}$-divisor $\mathscr{B}$ over $C$ such that $(\mathscr{W}_c,\frac{1}{2}\mathscr{B}_c)$ is klt for any $c\in C\setminus\{0\}$ and $(\mathscr{W},\frac{1}{2}\mathscr{B}+\mathscr{W}_0)$ is lc.
    Then, there exists a projective small contraction $g\colon \mathscr{W}^{-} \to \mathscr{W}$ such that $-K_{\mathscr{W}^{-}}$ is a $g$-ample $\mathbb{Q}$-Cartier divisor.
    Furthermore, such $g$ is unique up to isomorphism, and the pair $(\mathscr{W}^{-},\frac{1}{2}\mathscr{B}^{-}+\mathscr{W}^{-}_0)$ is also lc.
\end{prop}

\begin{proof}
By the assumptions and \cite[Proposition 2.15]{kollar-moduli}, the pair     $(\mathscr{W},\frac{1}{2}\mathscr{B})$ is klt.
   Then, by \cite[Theorem 7.2]{fujino--some--remarks}, the required morphism $g$ exists as claimed.
   The uniqueness of $g$ follows from \cite[Lemma 6.2]{KoMo}.
\end{proof}

Note that the induced morphism $g_0\colon (\mathscr{W}^{-}_{0}, \frac{1}{2}\mathscr{B}^{-}_{0})\to (\mathscr{W}_{0},\frac{1}{2}\mathscr{B}_0)$ on central fibers is an anti-P-resolution that is log crepant with respect to $(\mathscr{W}_{0}, \frac{1}{2}\mathscr{B}_0)$. 

A strategy to construct a partial $\mathbb{Q}$-Gorenstein smoothing of a given stable Horikawa surface $X$ proceeds as follows:
Let $W$ be the quotient of $X$ by a good involution and let $B$ be the branch divisor of this quotient map.
First, we perform an anti-P-resolution $(W^{-}, \frac{1}{2}B^{-})\to (W,\frac{1}{2}B)$.
Next, using Proposition \ref{prop--smoothing--quotient}, we construct a $\Q$-Gorenstein smoothing family $\mathscr{W}^{-} \to C$ of $W^{-}$, along with a $\Q$-Cartier Weil divisor $\frac{1}{2}\mathscr{B}^{-}$ on $\mathscr{W}^{-}$ that extends $\frac{1}{2}B^{-}$.
Then, we take the relative log canonical model $(\mathscr{W}, \frac{1}{2}\mathscr{B})$ of $(\mathscr{W}^{-}, \frac{1}{2}\mathscr{B}^{-})$ and consider a double cover $\mathscr{X}\to \mathscr{W}$ branched along $\mathscr{B}$.
Through this process, we obtain a partial $\Q$-Gorenstein smoothing family $\mathscr{X}\to C$ of $X$.
However, there is a subtle ambiguity. It might be the case that the generic fiber of $\mathscr{W}$ and the generic fiber of $\mathscr{W}^{-}$ are not isomorphic.
To avoid this possibility, we note the following useful lemma: 

\begin{lem}\label{lem--existence-anti-P}
    Let $f\colon \mathscr{W}\to C$ be a family of normal surfaces over a smooth affine curve.
    Let $0\in C$ be a closed point.
    Consider a proper birational morphism $g\colon \mathscr{W}^{-} \to \mathscr{W}$ from a $\Q$-Gorenstein threefold such that $-K_{\mathscr{W}^{-}_0}$ is relatively ample over $\mathscr{W}_0$.
    If $\mathscr{W}^{-}_c$ is smooth for any $c\in C\setminus\{0\}$, then there exists an open neighborhood $0\in U\subset C$ such that $\mathscr{W}_c$ is also smooth for any $c\in U\setminus\{0\}$.
    In other words, every singularity of $\mathscr{W}_0$ is smoothed by $f$.



\end{lem}

\begin{proof}
    We may freely shrink $C$ and assume that $-K_{\mathscr{W}^{-}}$ is $g$-ample by \cite[Proposition 1.41]{KoMo}.
    In particular, $-K_{\mathscr{W}^{-}_c}$ is $g_c$-ample for any $c\in C\setminus \{0\}$.
    Let $E_c:=K_{\mathscr{W}^{-}_c}-g_c^*K_{\mathscr{W}_c}$ be the canonical cycle.
    Since $-K_{\mathscr{W}^{-}_c}$ is $g_c$-ample, it follows from the negativity lemma \cite[Lemma 3.39]{KoMo} that $E_c$ is effective and its support contains all $g_c$-exceptional curves.
   This shows that $\mathscr{W}_c$ has only terminal singularities and hence $\mathscr{W}_c$ is smooth (cf.~\cite[Theorem 4.5]{KoMo}).
\end{proof}

Next, we consider anti-P-resolutions of cone singularities of type $(d;k_1,k_2)$ as in Definition~\ref{defn:(d;k_1,k_2)_sing}.

\begin{lem}\label{lem--anti--p-resol--k-square}
Let $(W,\frac{1}{2}B)$ be an lc surface pair such that $B$ is an effective Weil divisor.
Let $x\in W$ be a closed point.
    Suppose that $(x\in W, \frac{1}{2}B)$ is a germ of a cone singularity of type $(d;k_1,k_2)$ for $d\ge 3$.
    Consider a flat morphism $f\colon \mathscr{W}\to C$ with a fixed closed point $0\in C$ and an effective Weil divisor $\mathscr{B}$ on $\mathscr{W}$ that is horizontal over $C$ such that
\begin{itemize}
\item[$(1)$]
$K_{\mathscr{W}}+\frac{1}{2}\mathscr{B}$ is $\mathbb{Q}$-Cartier but $K_{\mathscr{W}}$ is not,
\item[$(2)$]
$(\mathscr{W}_0,\mathscr{B}_0)\cong (W,B)$, and
\item[$(3)$]
$\mathscr{W}_c$ is smooth and $(\mathscr{W}_c,\frac{1}{2}\mathscr{B}_c)$ is klt for any $c\in C\setminus\{0\}$.
\end{itemize}
    Let $g \colon \mathscr{W}^{-}\to \mathscr{W}$ be the projective small birational morphism from a $\Q$-Gorenstein family of slc surfaces over $C$ such that $-K_{\mathscr{W}^{-}}$ is $g$-ample obtained as Proposition \ref{prop--anti-P}.
    Note that $g_0$ defines an anti-P-resolution that is log crepant with respect to $(W,\frac{1}{2}B)$.   
    Then, we have 
    $$
    (K_{\mathscr{W}^{-}_0}-g_{0}^*K_W)^2=-\frac{(d-2)^2}{d}.
    $$
    Moreover, 
    after shrinking $C$ around $0$, we have $\rho(\mathscr{W}^{-}/\mathscr{W})=1$ and $\mathrm{Ex}(g)\cong\mathbb{P}^1$.
\end{lem}

\begin{proof}
Since the statement is local, we may assume that $W$ is affine and that the point $x\in W$ corresponds to a germ of a $\frac{1}{d}(1,1)$-singularity.
In paticular, $W$ is \'etale locally isomorphic to an open affine neighborhood $U\subset \overline{\Sigma}_d$ of the vertex $y$. 

First, we show that it is sufficient to prove the statement assuming that $W=\overline{\Sigma}_d$ and $f$ is projective.
Consider the thickenings
$$
\mathscr{W}_n:=\mathscr{W}\times_{C}C_n \to C_n:=\mathrm{Spec}(\mathcal{O}_{C,0}/\mathfrak{m}^{n+1}),
$$
where $\mathfrak{m}$ is the maximal ideal of $\mathcal{O}_{C,0}$. 
Since $W$ and $U$ are \'etale equivalent, there exists a deformation $\mathscr{U}_n\to C_n$ of $U$ which is \'etale locally isomorphic to $\mathscr{W}_n$. 
By \cite[Proposition 3.1]{HP}, for any $n>0$, there exists a deformation $\mathscr{Y}_n\to C_n$ of $\overline{\Sigma}_d$ containing $\mathscr{U}_n$ such that $\mathscr{Y}_n=\mathscr{Y}_{n+1}\times_{C_{n+1}}C_n$ and $\mathscr{U}_n\hookrightarrow \mathscr{Y}_n$ restricts $U\subset Y$.
Next, let $\mathscr{W}^-_n:=\mathscr{W}^-\times_{C}C_n$. 
Then, We may glue $\mathscr{U}_n\setminus\{y\}$ and $\mathscr{W}^-_n\setminus\mathrm{Ex}(g)$ to construct an algebraic space $g'_n\colon \mathscr{Y}^-_{n}\to \mathscr{Y}$. 
Since $-K_{\mathscr{W}^-}$ is $g$-ample, the reflexive sheaf $\omega^{[-1]}_{\mathscr{Y}^-_n/C_n}$ is $g_n$-ample.
Therefore, $\mathscr{Y}^-_n$ is in fact a scheme. 
Let $H$ be an ample line bundle on $\overline{\Sigma}_d$. 
Since $H^2(\overline{\Sigma}_d,\mathcal{O}_{\overline{\Sigma}_d})=0$, for each $n\ge 0$, we can find a line bundle $\mathscr{H}_n$ on $\mathscr{Y}_n$ such that $\mathscr{H}_{n+1}|_{\mathscr{Y}_{n}}=\mathscr{H}_n$ and $\mathscr{H}_{n}|_{\overline{\Sigma}_d}=H$. 
Moreover, there exist integers $r,l>0$ independent of $n$, such that the divisor $-rK_{\mathscr{W}_n/C_n}+lr\mathscr{H}_n$ is ample for every $n>0$.
By \cite[089A]{Stacks}, we may then construct projective normal  $\mathbb{Q}$-Gorenstein schemes $\widehat{\mathscr{Y}}^-$ and $ \widehat{\mathscr{Y}}$, flat over $R:=\mathrm{Spec}\,(\mathbb{C}[[t]])$, together with a birational morphism $\hat{g}'\colon \widehat{\mathscr{Y}}^-\to \widehat{\mathscr{Y}}$ over $R$ such that $\hat{g}'|_{\widehat{\mathscr{Y}}^-\times_{R}C_n}=g'_n$ for each $n\ge 0$.
By construction, $\hat{g}'$ is an isomorphism over the generic point of $R$.
Furthermore, by 
\cite[Tag 07GC]{Stacks} and the argument in the first paragraph of the proof of Proposition \ref{prop--smoothing--quotient}, 
there exist a smooth curve $C'$ with a closed point $0'$, projective flat morphisms $\mathscr{Y}\to C'$ and $\mathscr{Y}^-\to C'$ from normal varieties, and a birational small morphism $g'\colon \mathscr{Y}^-\to \mathscr{Y}$ over $C'$ such that the restriction $g'_{0'}$ is isomorphic to $g'_{0}$, $g'|_{\mathscr{Y}^-\times_{C'}C'\setminus\{0'\}}$ is an isomorphism, $\mathscr{Y}^-$ is $\mathbb{Q}$-Gorenstein, and $\mathscr{Y}^-_{c'}$ is smooth for any $c'\in C'\setminus\{0'\}$.
Therefore, we may replace $W$ by $\overline{\Sigma}_d$ in our argument.

Note that the singularity of type $\frac{1}{d}(1,1)$ with $d=3$ or $d\ge 5$ admits a unique P-resolution, which is the minimal resolution (cf.~\cite[Lemma 3.14]{KSB}). 
When $d=4$, P-resolutions are trivial or the minimal resolution.
Thus, applying \cite[Theorem 3.9]{KSB}, we may assume that there exists a simultaneous small resolution $h\colon \mathscr{W}^{+}\to \mathscr{W}$ with $\rho(\mathscr{W}^{+}/ \mathscr{W})=1$.
Let $\mathscr{B}^{+}$ be the proper transform of $\mathscr{B}$ on $\mathscr{W}^{+}$.
Since 
$$
K_{\mathscr{W}^{+}}+\frac{1+\varepsilon}{2}\mathscr{B}^{+}\sim_{\mathscr{W},\Q}-\varepsilon K_{\mathscr{W}^{+}},
$$
it follows that $\mathscr{W}^{-}$ coincides with the log flip of $(\mathscr{W}^{+},\frac{1+\varepsilon}{2}\mathscr{B}^{+})$ over $\mathscr{W}$, where $\varepsilon>0$ is sufficiently small.
Thus, it follows that $\rho(\mathscr{W}^{-}/\mathscr{W})=1$.
Furthermore, it follows from \cite[Lemma 5.1.4]{kawakita--threefold} that $\mathrm{Ex}(g)\cong\mathbb{P}^1$.
On the other hand, 
since $\mathscr{W}^{+}\to C$ is a family of Hirzebruch surfaces, we compute
$$
K_{\mathscr{W}^{-}_{0}}^{2}=K_{\mathscr{W}^{-}_{c}}^{2}=K_{\mathscr{W}^{+}_{c}}^{2}=8,
$$
where $c\in C\setminus \{0\}$.
Combining this with $K_{\overline{\Sigma}_d}^2=8+\frac{(d-2)^2}{d}$, we have the first assertion.
We complete the proof.
\end{proof}

\begin{cor}\label{cor--gen--fiber--rho}
Let $(W,\frac{1}{2}B)$ be a projective lc surface pair with $K_W^2=8+\frac{(d-2)^2}{d}$ such that $B$ is an effective Weil divisor with a cone singularity of type $(d;k_1,k_2)$ at $x\in W$, and any other singularity of $W$ is a Wahl singularity.
Suppose that $K_W+\frac{1}{2}B$ is ample.
Let $\mu \colon (W^{-}, \frac{1}{2}B^{-}) \to (W,\frac{1}{2}B)$ be an anti-P-resolution that is log crepant with respect to $(W, \frac{1}{2}B)$.
Suppose that $\rho(W)=1$, $-K_{W^{-}}$ is nef and log big and $\mu$ is not an isomorphism.
Consider a projective $\Q$-Gorenstein deformation $f\colon\mathscr{W}^{-} \to C$ from a normal threefold to a smooth affine curve with a closed point $0\in C$ and an effective $\mathbb{Q}$-Cartier divisor $\mathscr{B}^{-}$ on $\mathscr{W}^{-}$ satisfying:  
\begin{itemize}
    \item[$(1)$]
    $f$ is locally trivial on $W^{-} \setminus\mathrm{Ex}(\mu)$,
    \item[$(2)$]  $(\mathscr{W}^{-}_0,\mathscr{B}^{-}_0)\cong (W^{-},B^{-})$, 
    \item[$(3)$]
    $(\mathscr{W}^{-},\frac{1}{2}\mathscr{B}^{-}+\mathscr{W}^{-}_0)$ is lc, and
    \item[$(4)$]
    $\mathscr{W}^{-}_c$ is normal for any $c\in C\setminus\{0\}$.
\end{itemize}
Then, there exists the relative log canonical model $\mathscr{W}$ of $(\mathscr{W}^{-},\frac{1}{2}\mathscr{B^{-}})$ over $C$, which satisfies $\mathscr{W}_0\cong W$. 
Let $g\colon \mathscr{W}\to C$ be the canonical morphism.
If further any singularity contained in $\mathrm{Ex}(\mu)$ is smoothed by $f$ and $K_{\mathscr{W}^{-}_c}+\frac{1}{2}\mathscr{B}^{-}_c$ is not ample for any general $c\in C$, then $d=4$ and $g$ is the $\mathbb{Q}$-Gorenstein smoothing of $\frac{1}{4}(1,1)$.
\end{cor}

\begin{proof}
Since $-K_{W^{-}}$ is log big and nef,
$$H^1(\mathcal{O}_{W^{-}}(B^{-}))=H^1(\mathcal{O}_{W^{-}}(-mK_{W^{-}}))=H^1(\mathcal{O}_{W^{-}}(m(2K_{W^{-}}+B^-)))=0$$ for any sufficiently divisible $m>0$ by \cite[Theorem 1.10]{fujino--slc--vanishing}.
Furthermore, $-K_{W^-}$ is semiample by \cite[Theorem 1.16]{fujino--slc--vanishing}.
In particular, we may assume that $-K_{\mathscr{W}^{-}_c}$ is big for any $c\in C\setminus\{0\}$ by shrinking $C$. 
On the other hand, $H^1(\mathcal{O}_{W^{-}})=0$.
Since $\mathscr{W}^{-}_c$ has only rational singularities, we may also assume that $\mathscr{W}^{-}_c$ is rational.
By Lemma \ref{lem--anti--p-resol--k-square} and \cite[Proposition 2.6]{HP}, we may assume that $\rho(\mathscr{W}^{-}_c)=2$.
Since $H^1(\mathcal{O}_{W^{-}}(m(2K_{W^{-}}+B^-)))=0$ for any sufficiently divisible $m>0$, we may assume that $K_{\mathscr{W}^{-}}+\frac{1}{2}\mathscr{B}^{-}$ is $f$-big and $f$-semiample. 
Thus, there exists the relative log canonical model $\mathscr{W}$ of $(\mathscr{W}^{-},\frac{1}{2}\mathscr{B}^{-})$ over $C$. Let $g\colon \mathscr{W}\to C$ be the canonical morphism and note that $\mathscr{W}_0\cong W$.
Locally around $x\in W$, $g$ is a smoothing by Lemma \ref{lem--existence-anti-P}.
Assume that $K_{\mathscr{W}^{-}_c}+\frac{1}{2}\mathscr{B}^{-}_c$ is not ample for any general $c\in C$. Then, $\rho(\mathscr{W}_c)=1$.
This shows that $x$ should be a Wahl singularity of type $\frac{1}{d}(1,1)$ by \cite[Proposition 2.6]{HP} and hence $d=4$.
Thus, we complete the proof.
\end{proof}

In this paper, we only use the following type of anti-P-resolutions:

\begin{defn}[Admissible anti-P-resolution] \label{def--adm--anti--P}
    Let $(0\in W,\frac{1}{2}B)$ be a germ of a cone singularity of type $(d;k_1,k_2)$ with $d\ge5$. 
    Let $\mu \colon(W^-,\frac{1}{2}B^-)\to(W,\frac{1}{2}B)$ be an anti-P-resolution that is log crepant with respect to $(W, \frac{1}{2}B)$.
    Then, we say that $\mu$ is {\em admissible} if $\mathrm{Ex}(\mu)$ is irreducible and $(K_{W^-}-\mu^*K_W)^2=-\frac{1}{d}(d-2)^2$. 
    Note that the anti-P-resolution $g_0$ appeared in Lemma~\ref{lem--anti--p-resol--k-square} is admissible.
\end{defn}

In the remainder of this subsection, we present a complete classification of admissible anti-P-resolutions (Theorem~\ref{thm--anti-P-k_1-k_2}).
Before stating the result, we introduce three distinct types:
those of normal type, of non-normal type $\mathrm{I}$, and of non-normal type $\mathrm{II}$.

Let $(0\in W,\frac{1}{2}B)$ be a germ of a cone singularity of type $(d;k_1,k_2)$ with $d\ge 5$. 
Let $W^+$ and $B^+$ be as in Definition \ref{defn:(d;k_1,k_2)_sing}.

\begin{exam}[anti-P-resolution of normal type] \label{exam--anti-P-type-C} 
We consider a sequence of blow-ups
    $$
    W^+=:W_0\xleftarrow{\psi_1} W_1\xleftarrow{\psi_2} W_2\xleftarrow{\psi_3}\cdots\xleftarrow{\psi_{d-2}} W_{d-2}
    $$
constructed as follows:
First, fix a point $w_1$ which is an intersection point of the proper transform $\widehat{B}$ of $B$ on $W^{+}$ with the divisor $\Delta_0$.
Assume that $w_1\in \widehat{B}$ is an $A_{k_1-1}$-singularity for some $k_1\ge 2d-8$.
For $1\le i\le d-2$, we define each blow-up $\psi_i$ depending on the value of $k_1-2i+1$:

\begin{itemize}
    \item 
    If $k_1-2i+1\ge 0$, then $\psi_{i}$ is the blow-up at the point of the proper transform of $\widehat{B}$ lying over $w_1$, where it has an $A_{k_1-2i+1}$-singularity.
    \item 
    If $k_1=2i-2$, then either $(i,k_1)=(d-2,2d-6)$ or $(d-3,2d-8)$.
In this case, $\psi_{i}$ is defined as the blow-up at a point on the intersection of the proper transform of $B$ and the exceptional divisor $E_{i-1}$.
    \item 
    If $k_1=2i-3$, then necessarily $(i,k_1)=(d-2,2d-7)$.
    Here, $\psi_{i}$ is the blow-up at a point on $E_{i-1}$ that does not lie on the proper transform of $B$.
    \item 
    If $k_1=2i-4$, then $(i,k_1)=(d-2,2d-8)$, and $\psi_{i}$ is the blow-up at the intersection point of the proper transform of $B$ and $E_{i-1}$.
\end{itemize}

As a result of this sequence of blow-ups, the proper transforms of $\Delta_{0}, E_{1}, \ldots, E_{d-3}$ form a T-chain $[d+1,2^{d-3}]$.
Now define a divisor $B_{i}$ on each $W_{i}$ by the log crepant condition: 
$$
K_{W_i}+\frac{1}{2}B_{i}=\psi_{i}^{*}\left(K_{W_{i-1}}+\frac{1}{2}B_{i-1}\right), \quad B_{0}:=B^{+}=\widehat{B}+2\Delta_{0}.
$$
By construction, each $B_{i}$ is effective.
We then contract the T-chain on $W_{d-2}$ to obtain a surface$W^{-}$, and define $B^{-}$ as the pushforward of $B_{d-2}$ to $W^{-}$.
Thus, the induced morphism $\mu\colon (W^{-}, \frac{1}{2}B^{-})\to (W, \frac{1}{2}B)$ gives an admissible anti-P-resolution,
which is called an {\em anti-P-resolution of normal type} in this paper.
\end{exam}

\begin{exam}[Anti-P-resolution of non-normal type $\mathrm{I}$] \label{exam--anti-P-type-B}
Assume that the proper transform $\widehat{B}$ of $B$ has  an $A_{k_1-1}$-singularity $w_1$ lying on $\Delta_{0}$, where $k_1\ge 2d-7$.
We first take a sequence of blow-ups
    $$
    W^+=:W_0\xleftarrow{\psi_1} W_1\xleftarrow{\psi_2} W_2\xleftarrow{\psi_3}\cdots\xleftarrow{\psi_{d-3}} W_{d-3}
    $$
as in Example~\ref{exam--anti-P-type-C}.
Next, we perform an additional blow-up $\psi_{d-2}\colon W_{d-2}\to W_{d-3}$ at the intersection point of $E_{d-3}$ and the proper trasform of $E_{d-4}$.
Then, the proper transforms of $\Delta_{0}, E_{1},\ldots, E_{d-2}$ form a chain 
$$
[2,1,3,2^{d-5},d+1]=R_{2}-1-L_{1}L_{2}^{d-4}L_{1}[d].
$$
Now, we give the following two types of admissible anti-P-resolutions, both of which are called {\em anti-P-resolutions of non-normal type $\mathrm{I}$}:

\begin{itemize}
\item[$(1)$]
Assume that $d\ge 6$.
We perform a sequence of further blow-ups $W_{N}\to \cdots \to W_{d-2}$,
each at one of the two points in the intersection of the $(-1)$-curve and adjacent components, resulting in the configuration: 
$$
R_{2}-1-L_{1}L_{2}^{d-4}L_{1}[d]\leftarrow \cdots \leftarrow R_{2}^{d}(R_{1}^{d-4}R_{2}^{d})^{m}R_{1}^{\alpha}-1-L_{2}^{\alpha}(L_{1}^{d}L_{2}^{d-4})^{m+1}L_{1}[d],
$$
where $m\ge 0$ and $0\le \alpha\le d-6$.
We define the following chains:
\begin{align*}    
C_{1}&:=R_{2}^{d}(R_{1}^{d-4}R_{2}^{d})^{m}R_{1}^{\alpha}=[2^{d-1}, d-2, \ldots, 2^{d-1}, \alpha+2], \\
C_{2}&:=L_{2}^{\alpha}(L_{1}^{d}L_{2}^{d-4})^{m}L_{1}[d]=[2^{\alpha}, d+2, 2^{d-5}, \ldots, d+2, 2^{d-5}, d+1], \\
C_{3}&:=L_{2}^{\alpha}L_{1}^{d}L_{2}^{d-5-\alpha}=[2^{\alpha}, d+2, 2^{d-6-\alpha}].
\end{align*}
Then, the proper transforms of $\Delta_{0}$ and the exceptional curves form a chain $C_1-1-C_3-2-C_2$.
We define a divisor $B_{i}$ on each $W_{i}$ inductively via the log crepant condition: 
$$
K_{W_i}+\frac{1}{2}B_{i}=\psi_{i}^{*}(K_{W_{i-1}}+\frac{1}{2}B_{i-1}), \quad B_{0}:=B^{+}=\widehat{B}+2\Delta_{0}.
$$
By construction, each $B_{i}$ is effective.
Note that $C_1$ and $C_2$ are conjugate in the sense that $C_1-1-C_2$ blows down to $[0]$.
We then contract the three chains $C_1$, $C_2$, $C_3$ to obtain a surface $W'$, and define $\Delta'_{1}$, $\Delta'_{2}$ as the images of the $(-1)$-curve $E_{N}$ and the $(-2)$-curve connecting $C_3$ and $C_2$, respectively.
Let $B'=B''+2\Delta'_1+2\Delta'_2$ denote the pushforward of $B_N$ to $W'$.
We note that there exists an involution $\tau$ on $\Delta'_{1}+\Delta'_{2}$ exchanging $\Delta'_1$ and $\Delta'_2$ and identifying $\mathrm{Diff}_{\Delta'_{1}}(\frac{1}{2}B')$ and $\mathrm{Diff}_{\Delta'_{2}}(\frac{1}{2}B')$ (fro the definition of the different, see \cite[\S16]{Koetal}).
Hence, by Koll\'{a}r's gluing theorem \cite[Theorem 11.38]{kollar-moduli}, $\Delta'_1$ and $\Delta'_2$ can be identified via the involution $\tau$, giving rise to a morphism $W'\to W^{-}$.
Indeed, take a divisor $H$ on $W'$ which is ample over $W$ and $H\cdot \Delta'_{1}=H\cdot \Delta'_{2}$, the existence of which follows from the negative definiteness of $\Delta'_{1}+\Delta'_{2}$.
By the Serre vanishing theorem, the restriction map 
$$
H^{0}(\O_{W'}(mH))\to H^{0}(\O_{W'}(mH)|_{\Delta'_{1}+\Delta'_{2}})
$$
is surjective for $m\gg 0$.
Then we can take a member $A\in |mH|$ satisfying $\tau^{*}A|_{\Delta'_{1}+\Delta'_{2}}=A|_{\Delta'_{1}+\Delta'_{2}}$.
Since $K_{W'}+\frac{1}{2}B'+\varepsilon A$ is ample over $W$ and lc for sufficiently small $\varepsilon>0$, we can apply the gluing theorem (\cite{kollar-mmp}, \cite[Theorem~11.38]{kollar-moduli}) to the pair $(W', \frac{1}{2}B'+\varepsilon A)$.
By construction, the canonical morphism $\mu\colon W^{-}\to W$ gives an admissible anti-P-resolution since all degenerate cusps are locally smoothable by \cite{JSte}. It is also easy to see that the exceptional locus is isomorphic to $\mathbb{P}^1$.

\item[$(2)$]
we perform another sequence of blow-ups $W_{M}\to \cdots \to W_{d-2}$ as follows:
\begin{align*}
&R_{2}-1-L_{1}L_{2}^{d-4}L_{1}[d]\leftarrow \cdots \leftarrow R_{2}^{d}(R_{1}^{d-4}R_{2}^{d})^{m}R_{1}^{\alpha}-1-L_{2}^{\alpha}(L_{1}^{d}L_{2}^{d-4})^{m+1}L_{1}[d] \\
&\leftarrow R_{2}^{d}(R_{1}^{d-4}R_{2}^{d})^{m}R_{1}^{\alpha}R_{2}-1-L_{1}L_{2}^{\alpha}(L_{1}^{d}L_{2}^{d-4})^{m+1}L_{1}[d] \\
&\leftarrow R_{2}^{d}(R_{1}^{d-4}R_{2}^{d})^{m}R_{1}^{\alpha}R_{2}-1-[3,2^{\alpha-1}, d+2, 2^{d-\alpha-6}, 3]-1-L_{1}L_{2}^{\alpha}(L_{1}^{d}L_{2}^{d-4})^{m}L_{1}[d],
\end{align*}
where $m\ge 0$ and $0\le \alpha\le d-5$.
After further blow-ups $W_{N}\to \cdots \to W_{M}$ at nodes on $(-1)$-curves, 
we obtain a chain $C_1-1-C_3-1-C_2$ consisting of the proper transform of $\Delta_{0}$ and the exceptional curves, where
we define 
\begin{align*}    
C_{1}&:=R_{2}^{d}(R_{1}^{d-4}R_{2}^{d})^{m}R_{1}^{\alpha}R_{2}\bm{R}=[2^{d-1}, d-2, \ldots, 2^{d-1}, \alpha+2,2]\bm{R}, \\
C_{2}&:=\bm{L}L_{1}L_{2}^{\alpha}(L_{1}^{d}L_{2}^{d-4})^{m}L_{1}[d]=\bm{L}[3, 2^{\alpha-1}, d+2, 2^{d-5}, \ldots, d+2, 2^{d-5}, d+1], \\
C_{3}&:=\bm{L}[3,2^{\alpha-1}, d+2, 2^{d-\alpha-6}, 3]\bm{R}, \\
\bm{L}&:=L_{i_l}\cdots L_{i_1}, \quad \bm{R}:=R_{3-i_1}\cdots R_{3-i_l}.
\end{align*}
As in case (1), by contracting the three chains $C_1$, $C_2$, $C_3$ and gluing the images of the two $(-1)$-curves via a natural involution, we obtain morphisms $W_{N}\to W'\to W^{-}$ and an admissible anti-P-resolution $\mu\colon W^{-}\to W$.
\end{itemize}

\end{exam}

\begin{exam}[Anti-P-resolution of non-normal type $\mathrm{II}$] \label{exam--anti-P-type-A} 
Let $w_i$, $i=1,2$ denote the $A_{k_i-1}$-singularities of $\widehat{B}$ lying on $\Delta_0$, and assume that $k_1>0$, $k_2>0$.
Then, we give the following two types of admissible P-resolutions, both of which are called {\em anti-P-resolutions of non-normal type $\mathrm{II}$}:

\begin{itemize}
    \item[$(1)$]
    Suppose that
$\lfloor \frac{k_1}{2} \rfloor + \lfloor \frac{k_2}{2} \rfloor \ge  d-4$.
We choose non-negative integers $\alpha_1$ and $\alpha_2$ such that $\alpha_{i}+1\le \lfloor \frac{k_i}{2} \rfloor$ and $\alpha_1+\alpha_2=d-6$.
Moreover, we assume that either $2\alpha_1+2=k_1$ and $2\alpha_2+2=k_2$ or $2\alpha_1+2<k_1$ and $2\alpha_2+2<k_2$ (this assumption is used to the gluing condition).
We consider a sequence of blow-ups
    $$
    W^+=:W_0\xleftarrow{\psi_1} W_1\xleftarrow{\psi_2} \cdots \xleftarrow{\psi_{\alpha_{1}+1}}
 W_{\alpha_1+1}\xleftarrow{\psi_{\alpha_1+2}}\cdots\xleftarrow{\psi_{d-4}} W_{d-4},
    $$    
     where $\psi_{1}\circ \cdots \circ \psi_{\alpha_1+1}$ (resp.\ $\psi_{\alpha_1+2}\circ \cdots \circ \psi_{d-3}$) is a partial resolution process at $w_1$ (resp.\ at $w_2$) as in Example~\ref{exam--anti-P-type-C}.
     Then, the proper transforms of $\Delta_0, E_1, \ldots, E_{d-4}$ forms a chain
     $$
     1-[2^{\alpha_1}, d+2, 2^{\alpha_2}]-1.
     $$
     As in Example~\ref{exam--anti-P-type-B}, by contracting the $C_3:=[2^{\alpha_1}, d+2, 2^{\alpha_2}]$ and gluing the images of the two $(-1)$-curves via a natural involution, we obtain a sequence of morphisms $W_{N}\to W'\to W^{-}$, and a morphism $\mu\colon W^{-}\to W$ which defines an admissible anti-P-resolution.

    \item[$(2)$]
    Suppose that
$\lceil \frac{k_1}{2} \rceil + \lceil \frac{k_2}{2} \rceil \ge  d-3$.
We choose non-negative integers $\alpha_1$ and $\alpha_2$ such that $\alpha_{i}+1\le \lceil \frac{k_i}{2} \rceil$ and $\alpha_1+\alpha_2=d-5$.
Similary to (1), we consider a sequence of blow-ups
    $$
    W^+=:W_0\xleftarrow{\psi_1} W_1\xleftarrow{\psi_2} \cdots \xleftarrow{\psi_{\alpha_{1}+1}}
 W_{\alpha_1+1}\xleftarrow{\psi_{\alpha_1+2}}\cdots\xleftarrow{\psi_{d-3}} W_{d-3},
    $$
    where $\psi_{1}\circ \cdots \circ \psi_{\alpha_1+1}$ (resp.\ $\psi_{\alpha_1+2}\circ \cdots \circ \psi_{d-3}$) is a partial resolution process at $w_1$ (resp.\ at $w_2$) as in Example~\ref{exam--anti-P-type-C}.
    Next, we perform two additional blow-ups $W_{d-1}\to W_{d-2}\to W_{d-3}$, centered at $E_{\alpha_1+1}\cap E_{\alpha_1}$ and $E_{d-3}\cap E_{d-4}$, where $E_i$ denotes the proper transform of the $\psi_i$-exceptional curve.
    Then, the proper transforms of $\Delta_0, E_1, \ldots, E_{d-1}$ form a chain
    $$
    [2]-1-[3,2^{\alpha_1-1}, d+2, 2^{\alpha_2-1}, 3]-1-[2].
    $$
    Following the same procedure as in Example~\ref{exam--anti-P-type-B}~(2),
   we perform further blow-ups $W_{N}\to \cdots \to W_{d-1}$ at nodes lying on $(-1)$-curves.
   As a result, we obtain a chain $C_1-1-C_3-1-C_2$ consisting of the proper transform of $\Delta_{0}$ and the exceptional divisors, where
we define 
\begin{align*}    
C_{1}&:=[2]\bm{R},\quad  C_{2}:=\bm{L}[2], \quad C_{3}:=\bm{L}[3,2^{\alpha_1-1}, d+2, 2^{\alpha_2-1}, 3]\bm{R}, \\
\bm{L}&:=L_{i_l}\cdots L_{i_1}, \quad \bm{R}:=R_{3-i_1}\cdots R_{3-i_l}.
\end{align*}
Here, by the gluing condition, we assume that
if $2\alpha_{1}+1=k_1$ and $2\alpha_{2}+1<k_2$ (resp.\ $2\alpha_{1}+1<k_1$ and $2\alpha_{2}+1=k_2$, $2\alpha_{1}+1=k_1$ and $2\alpha_{2}+1=k_2$), then we assume that $i_1=1$ (resp.\ $i_1=2$, $l=0$).
As in Example~\ref{exam--anti-P-type-B}, by contracting the three chains $C_1$, $C_2$, $C_3$ and gluing the images of the two $(-1)$-curves via a natural involution, we obtain a sequence of morphisms $W_{N}\to W'\to W^{-}$, and a morphism $\mu\colon W^{-}\to W$ which defines an admissible anti-P-resolution.

\end{itemize}
\end{exam}

\begin{thm}\label{thm--anti-P-k_1-k_2}
Let $(0\in W, \frac{1}{2}B)$ be a cone singularity of type $(d;k_1,k_2)$ with $d\ge 5$, and let $\mu\colon (W^{-}, \frac{1}{2}B^{-})\to (W, \frac{1}{2}B)$ be an admissible anti-P-resolution.
Then, the following hold:

\begin{itemize}
\item[$(1)$]
If $W^{-}$ is normal, then $\mu$ is constructed as in Example~\ref{exam--anti-P-type-C}.
In particular, $W^{-}$ is klt and either $k_1\ge 2d-8$ or $k_2\ge 2d-8$.

\item[$(2)$]
If $W^{-}$ is non-normal, then $\mu$ is constructed as in Example~\ref{exam--anti-P-type-B} or Example~\ref{exam--anti-P-type-A}.
In particular, one of the following holds:

\begin{itemize}
\item 
$k_1\ge 2d-8$ or $k_2\ge 2d-8$,
\item 
$k_1>0$, $k_2>0$ and 
$\lceil \frac{k_1}{2} \rceil + \lceil \frac{k_2}{2} \rceil \ge  d-3$,
\item 
$k_1>0$, $k_2>0$ and 
$\lfloor \frac{k_1}{2} \rfloor + \lfloor \frac{k_2}{2} \rfloor \ge  d-4$.
\end{itemize}

\end{itemize}

\end{thm}

\begin{proof}
    One sees that $W^+$, $W^-_{\mathrm{min}}$, and $W^-$ fit into the following diagram:
    \[
    W^+=:W_0\xleftarrow{\psi_1} W_1\xleftarrow{\psi_2} W_2\xleftarrow{\psi_3}\cdots\xleftarrow{\psi_N} W_N=W^{-}_{\mathrm{min}}\xrightarrow{c} (W^-)^\nu\xrightarrow{\nu}W^-.
    \]
    Here, each morphism $\psi_{i}\colon W_{i}\to W_{i-1}$ is a blow-up at a point $x_i\in W_{i-1}$ with the exceptional divisor $E_{i}\subset W_{i}$.
    Note that $N\geq1$, otherwise $W^-\cong W^+$ and $\mu$ is no longer an anti-P-resolution.
    Let $B_0:=B^+$ and let $B_i$ ($1\leq i\leq N$) be the Weil divisor on $W_i$ such that $K_{W_{i}}+\frac{1}{2}B_{i}=\psi_{i}^{*}(K_{W_{i-1}}+\frac{1}{2}B_{i-1})$.
    From now on, we classify this blow-up process.    
    We note that, for $1\leq i\leq N$, the multiplicity $m_i$ of $B_{i-1}$ at $x_i$ is at least two.
    Indeed, if $m_i<2$, then $B_{i}=(\psi_{i})_*^{-1}B_{i-1}+(m_i-2)E_{i}$ is not effective.
    This implies that $B^-$ is not effective, which is a contradiction.
    We examine the case where $W^-$ is normal in Step 1, and the non-normal case in Step 2.

{\bf Step 1.}
    We first assume that $W^-$ is normal.
    Since $\rho(W^-/W)=1$, the morphism $\nu\circ c\colon W_N\to W^-$ is the minimal resolution and contracts the proper transforms of $\Delta_0,E_1,\dots,E_{N-1}$.
    This implies that, for $1\leq i\leq N-1$, the proper transform of $E_i$ on $W_N$ is not a $(-1)$-curve, and the blow-up $W_{i+1}\to W_{i}$ must occur at a point on $E_i$, i.e., $x_{i+1}\in E_i$.

    We specify the first blow-up $\psi_1\colon W_1\to W_0$.
    Let $\widehat{B}=B^{+}-2\Delta_{0}$ denote the proper transform of $B$ on $W^{+}$.
    The point $x_1\in \widehat{B} \cap\Delta_0$ is one of the following:
    \begin{itemize}
        \item[(A)]
        $x_1\notin \widehat
        {B}$.
        \item[(B)]
        $x_1$ is a smooth point on $\widehat{B}$ such that $(\widehat{B} \cdot \Delta_0)_{x_1} =1$. 
        \item[(C)]
        $x_1$ is a (not necessarily smooth) point on $\widehat{B}$ such that $(\widehat{B} \cdot \Delta_0)_{x_1} =2$.
    \end{itemize}
    We examine each case separately.

    In the case (A), one can see from $m_i\geq2$ that the exceptional set $\mathrm{Ex}(c)=\Delta_0\cup E_1\cup\cdots\cup E_{N-1}$ is $[d+N]-1-[2^{N-1}]$ where $-1-$ represents $E_N$.
    This is a contradiction, since $[d+N]$ is not a T-chain.

    In the case (B), the coefficient of $E_1$ in $B_1$ is one.
    Note that $N\geq2$ since the proper transform of $\Delta_0$ in $W_1$ is not a T-chain.
    Let $p$ (resp. $q$) denote the intersection point of $E_1$ and the proper transform of $B^+$ (resp. $\Delta_0$).
    Since $m_1\geq2$, the point $x_2\in E_1$ must be either $p$ or $q$.
    
    If $x_2=p$, then $N=2$ since each point on $E_2$ has multiplicity of $B_2$ at most one, and $\mathrm{Ex}(c)$ becomes the chain $[d+1,2]$.
    This chain is a T-chain (of Milnor number zero) if and only if $d=4$.
    Since we assume $d\ge 5$, this case does not occur.
    
    If $x_2=q$, one can easily see from the properties $x_{i+1}\in E_i$ and $m_i\geq2$ that $\mathrm{Ex}(c)\cup E_N$ becomes $[d+N]-1-[2^{N-1}]$ ($N\geq2$) or $[d+N-1,2]-1-[3,2^{N-3}]$ ($N\geq3$).
    However, since they are not a T-chain, we see that this case does not occur.

    Finally, we consider the case (C).
    The condition in (C) implies that $x_1$ is either $w_1$ ($k_1\geq2$) or $w_2$ ($k_2\geq2$).
    We may assume that $k_1\geq2$ and $x_1=w_1$.
    Then, there exists $0\leq l'\leq\lfloor \frac{k_1}{2} \rfloor-1$ satisfying the following:
    \begin{itemize}
        \item 
        The proper transforms of $\Delta_0,E_1,\dots,E_{l'+1}$ in $W_{l'+1}$ form the chain $[d+1,2^{l'}]-1$ where $-1$ represents $E_{l'+1}$.
        \item 
        The coefficient of $E_{l'+1}$ in $B_{l'+1}$ is two.
    \end{itemize}
    Let $l$ be the largest $l'$ satisfying the above property.
    If $N=l+1$, then $[d+1,2^l]$ must be a T-chain, and hence $l=d-3$ and $\mu$ coincides with an anti-P-resolution described in Example~\ref{exam--anti-P-type-C}.
    We assume that $N>l+1$.
    Then, it follows that $2\leq m_{l+2}\leq 4$.
    We examine $\mathrm{Ex}(c)$ in each case.

    We first assume that $m_{l+2}=2$.
    In this case, $\mathrm{Ex}(c)$ becomes the chain $[d+1,2^{l+1}]$ ($N=l+2$) or $[d+1,2^l,N-l]-1-[2^{N-l-2}]$ ($N>l+2$).
    Since each connected component of $\mathrm{Ex}(c)$ must be a T-chain, we have $N=l+2$ and $l=d-4$.
    Then $\mu$ also coincides with the one described in Example~\ref{exam--anti-P-type-C}.
    
    We assume that $m_{l+2}=3$.
    This implies that $B_{l+1}-2E_{l+1}$ is smooth at $x_{l+2}\in W_{l+1}$.
    In particular, $k_1=2l+2$ or $k_1=2l+3$. 
    In the former case, the chain $\mathrm{Ex}(c)$ is of one of the following forms: $[d+1,2^{N-1}]$ 
    ($N=l+2,l+3$), $[d+1,2^{l},N-l]-1-[2^{N-l-2}]$, or $[d+1,2^l,N-l-1,2]-1-[3,2^{N-l-4}]$.
    Here $-1-$ represents $E_N$.
    Since $\mathrm{Ex}(c)$ must be a union of (possibly one) T-chains, $\mathrm{Ex}(c)$ is the chain $[d+1,2^{d-3}]$ ($N=d-2$ and $l=d-5, d-4$).
    In the latter case, if $N=l+2$, the chain $\mathrm{Ex}(c)$ is of the form $[d+1,2^{N-1}]$ ($N=d-2$ and $l=d-4$).
    If not, the proper transforms of $\Delta_0,E_1,\dots,E_{l+3}$ on $W_{l+3}$ forms the chain $[d+1,2^l,3]-1-[2]$.
    Note that the multiplicity of $E_{l+3}$, the rightmost curve in the chain, is one.
    From Lemma \ref{lem:drill_anti-P}, we see that this case does not occur.
    
    Finally, we consider the case of $m_{l+2}=4$.
    The maximality of $l$ implies that $x_{l+2}\in W_{l+1}$ is the intersection point of the proper transform of $E_{l}$ and $E_{l+1}$.
    Hence the proper transforms of $\Delta_0,E_1,\dots,E_{l+2}$ on $W_{l+2}$ form the chain $[d+1,2^{l-1},3]-1-[2]$.
    Lemma \ref{lem:drill_anti-P} indicates that $l=d-4$. 
    This falls into Example~\ref{exam--anti-P-type-C}.  

    {\bf Step 2.}
    We assume that $W^-$ is not normal.
    Note that $W^-$ is irreducible since $W_N$ is irreducible.
    Let $\Delta^\nu$ denote the conductor of the normalization $\nu\colon (W^-)^\nu\to W^-$.
    Note that $\Delta^\nu$ is the exceptional curve of $(W^-)^\nu\to W$.
    Following the arguments of the proofs of \cite[Theorem 5.3, 5.5]{Hac}, one can show that $((W^-)^\nu,\Delta^\nu)\to W^-$ is of one of the following types\footnote{\cite[Theorem 5.3, 5.5]{Hac} assumes that the surface is proper and its anticanonical divisor is ample, but the same argument works for $W^-$ from the fact that $-K_{W^-}$ is $\psi^-$-ample.}:
    \begin{itemize}
        \item[(i)] 
        $\Delta^\nu\cong\PP^1$ and $((W^-)^\nu,\Delta^\nu)$ is log terminal.
        The surface $W^-$ is obtained by folding $\Delta^\nu$ (in the sense of \cite[Notation 5.4]{Hac}).
        \item[(ii)] 
        $\Delta^\nu$ is a chain of two $\PP^1$'s. 
        $((W^-)^\nu,\Delta^\nu)$ is log canonical, and log terminal away from the node of $\Delta^\nu$.
        The surface $W^-$ is obtained by gluing a one component of $\Delta^\nu$ to the other component of $\Delta^\nu$.
    \end{itemize}
    The type (i) (resp. the type (ii)) corresponds to the type B* (resp. the type C) in \cite[Theorem 5.5]{Hac}.
    
    We show that the case (i) does not occur.
    Assuming $W^-$ is of the case (i), we derive a contradiction.
    The list of all $\Q$-Gorenstein smoothable semi-log terminal singularities in \cite[5.2]{KSB} \cite[Proposition 6.1]{Hac} indicates that, if there exists a cyclic quotient singularity $p\in(W^-)^\nu$ of the form $\frac{1}{r}(1,a)$ on $\Delta^\nu$, then there exists another cyclic quotient singularity $q\in(W^-)^\nu$ of the form $\frac{1}{r}(-1,a)$ on $\Delta^\nu$.
    Moreover, the germ of $(p\in(W^-)^\nu,\Delta^\nu)$ (resp. $(q\in(W^-)^\nu,\Delta^\nu)$) is isomorphic to $(\frac{1}{r}(1,a),(x=0))$ (resp. $(\frac{1}{r}(-1,a),(x=0))$).
    Since $\mathrm{Ex}(\mu \circ\nu)=\Delta^\nu\cong\PP^1$, the surface $(W^-)^\nu$ can be constructed from $W_N$ by contracting the proper transforms of $\Delta_0,E_1,\dots,E_{N-1}$.
    As observed in the normal case, the blow-up $W_{i+1}\to W_i$ must occur at a point on $E_i$, and the configuration of the total transform of $\Delta_0$ must be of the form $C_1-1-C_2$, where $C_1$ and $C_2$ are conjugate in the sense that $C_1-1-C_2$ blows down to $[0]$, and the $(-1)$-curve connecting $C_1$ and $C_2$ corresponds to $E_N$.
    This leads a contrtadiction since the intersection form on the total transform of $\Delta_{0}$ is negative definite.
    
    From the above discussion, the surface $W^-$ is of type (ii).
    By an argument similar to the one in the case (i), we can see the following:
    \begin{itemize}
        \item 
        $\mathrm{Ex}(\mu \circ\nu)=\Delta^\nu$.
        \item 
        The morphism $c\colon W_N\to(W^-)^\nu$ contracts the proper transforms of $\Delta_0, E_1,\ldots,E_{N-1}$ except for one, say $E_i$.
        We write $\Delta_1:=E_i$ and $\Delta_2:=E_N$.
        Note that their pushforwards $\Delta_1^\nu=c_*\Delta_1$, $\Delta_2^\nu=c_*\Delta_2$ are the components of $\Delta^\nu$.
        \item 
        If there exists a cyclic quotient singularity $p\in(W^-)^\nu$ and the germ $(p\in(W^-)^\nu,\Delta^\nu)$ is isomorphic to $(\frac{1}{r}(1,a),(x=0))$, then the germ $(q\in(W^-)^\nu,\Delta^\nu)$, where $q$ denotes the image of $p$ under the involution on $\Delta^\nu$ induced by $\nu$, is isomorphic to $(\frac{1}{r}(-1,a),(x=0))$.
        \item 
        The germ of $((W^-)^\nu,\Delta^\nu)$ at the node of $\Delta^\nu$ is isomorphic to $(\frac{1}{s}(1,b),(xy=0))$.
    \end{itemize}
    The last assertion follows from \cite[Theorem 4.15]{KoMo}.
    From these observations, one can see that the dual graph of (the support of) the total transform of $\Delta_0\subset W_0$ on $W_N$ is of the following form:
    \begin{equation}\label{eqn:nonnormal}
    \xygraph{
    \circ ([]!{-(0,-.3)}) 
    - [r]  \cdots ([]!{-(0,-.3)} {C_1}) 
    - [r]  \circ ([]!{-(0,-.3)}) 
    - [r]  \bullet ([]!{-(0,-.3)} {\Delta_1}) 
    - [r]  \circ ([]!{-(0,-.3)}) 
    - [r]  \cdots ([]!{-(0,-.3)} {C_3}) 
    - [r]  \circ ([]!{-(0,-.3)}) 
    - [r]  \bullet ([]!{-(0,-.3)} {\Delta_2}) 
    - [r]  \circ ([]!{-(0,-.3)}) 
    - [r]  \cdots ([]!{-(0,-.3)} {C_2}) 
    - [r]  \circ ([]!{-(0,-.3)}) 
    }
    \end{equation}
    Here, the white vertices correspond to the $c$-exceptional curves, and we write $C_1$ for the leftmost, $C_2$ for the rightmost, and $C_3$ for the middle connected components of the subgraph of the above graph consisting of the white vertices.
    Note that each connected component corresponds to that of $\mathrm{Ex}(c)$, and $C_i$ may be the empty set.
    If $C_1$ or $C_2$ is nonempty, then the third observation imply that the chain $C_1-1-C_2$ 
    can be obtained from $[2]-1-[2]$ by repeating the following procedures \cite[Lemma 7.3]{HP}: for any chain of the form $A-1-B$, we take $AR_1-1-L_2B$ or $AR_2-1-L_1B$.
    
    We introduce the following notation:
    \begin{itemize}
        \item 
        If $C_1$ is nonempty (which implies that $C_2$ is also nonempty), then we write $(\frac{1}{r}(1,a),(x=0))$ (resp. $(\frac{1}{r}(-1,a),(x=0))$) for the type of the singularity on $(W^-)^\nu$ obtained by contracting $C_1$ (resp. $C_2$).
        \item 
        Let $l_i$ ($i=1,2,3$) denote the length of $C_i$.
        If $C_i$ is empty, then we set $l_i=0$.
        \item 
        For $1\leq i\leq3$, we set 
        \[
        A_i:=-\sum_{E\subset C_i} a((W^-)^\nu,\Delta^\nu,E)\cdot(E\cdot(\Delta_1+\Delta_2-E)-2)
        \]
        where the sum is over all irreducible components of $C_i$, and $a((W^-)^\nu,\Delta^\nu,E)$ is the discrepancy of $((W^-)^\nu,\Delta^\nu)$ with respect to $E$, and $A_i:=0$ if $C_i$ is empty.
        
    \end{itemize}
    We claim the following equalities:
    \begin{align}\label{eqn:A_1+A_2+A_3}
        A_1+A_2+A_3&=N+(\Delta_1+\Delta_2)^2-4p_a(\Delta_1+\Delta_2)+4,   \\
        \label{eqn:A_1+A_2}
        A_1+A_2&=l_1+l_2-1+\delta,  \\
        \label{eqn:A_3}
        A_3&=l_3+(\Delta_1+\Delta_2)^2-4p_a(\Delta_1+\Delta_2)+6-\delta,
    \end{align}
    where $\delta:=0$ if $C_1$ and $C_2$ are nonempty and $\delta:=1$ otherwise.
    Note that \eqref{eqn:A_3} follows from \eqref{eqn:A_1+A_2+A_3} and \eqref{eqn:A_1+A_2}.
    From the definition of $A_i$'s and the equation $(K_{W^-/W})^2 = (K_{(W^-)^\nu/W} + \Delta^\nu)^2$, it follows that
    \[
    A_1+A_2+A_3 = K_{W^-/W}^2-(K_{W_N/W}+\Delta_1+\Delta_2)^2.
    \]
    By direct calculation, we see that
    \[
    (K_{W_N/W}+\Delta_1+\Delta_2)^2=K_{W^+/W}^2-N-(\Delta_1+\Delta_2)^2+4p_a(\Delta_1+\Delta_2)-4.
    \]
    Finally, we have $K_{W^-/W}^2=K_{W^+/W}^2$ since $\mu$ is admissible.
    These equations show \eqref{eqn:A_1+A_2+A_3}.
    
    Next, we show \eqref{eqn:A_1+A_2}. 
    Let $Y_1\to\Sigma_0$ be a composition of blow-ups of $\Sigma_0\to\PP^1$ such that $Y$ has exactly one singular fiber (as a ruled surface $Y_1\to\Sigma_0\to\PP^1$) and it is the chain of the form $C_1-1-C_2$. 
    Let $D_1$ denote the $(-1)$-curve on $Y_1$, let $Y_1\to Y_2$ be the morphism whose exceptional curves coincide with $C_1+C_2$, and let $D_2\subset Y_2$ be the pushforward of $D_1\subset Y_1$.
    Then, one can see that 
    \[
    A_1+A_2=(K_{Y_2}+D_2)^2-(K_{Y_1}+D_1)^2.
    \]
    Moreover, we have $(K_{Y_1}+D_1)^2=5-l_1-l_2$ since $K_{Y_1}^2=8-l_1-l_2$ and $D_1$ is a $(-1)$-curve.
    
    We claim that $(K_{Y_2}+D_2)^2=4$; this completes the proof of \eqref{eqn:A_1+A_2}.
    Let $S_{1,0}$ and $S_{2,0}$ be disjoint sections of $\Sigma_0\to\PP^1$.
    We assume that the first blow-up of $Y_1\to\Sigma_0$ occurs at a point on $S_{2,0}$.
    Let $S_i\subset Y_2$ ($i=1,2$) denote the pushforward of the proper transform of $S_i$, and let $F$ be a smooth fiber of $Y_2$.
    Note that $S_1+S_2+F+D_2\in|-K_{Y_2}|$, and hence $(K_{Y_2}+D_2)^2=S_1^2+S_2^2+4$.
    We claim that $S_1^2=a/r$ and $S_2^2=-a/r$.
    In fact, the first equation follows from the observation that $S_1\subset Y_2$ is locally isomorphic to the germ $((z=0)\subset\PP(r,1,a))$.
    Let $Y_2'$ be the surface obtained by blow-up $Y_2$ at $S_1\cap F$ and contracting the proper transform of $F$, and let $S_2'\subset Y_2'$ be the pushforward of the pullback of $S_2$.
    Then, $S_2^2=(S_2')^2-1$ and $S_2'\subset Y_2'$ is locally isomorphic to the germ $((z=0)\subset\PP(r,1,r-a))$.
    This shows $S_2^2=-a/r$.
    
    We classify $((W^-)^\nu,\Delta^\nu)$ by considering the two cases: when the proper transform $\Delta_0$ on $W_N$ is contained in $C_3$ and when it is not.
    Note that the former case corresponds to Example~\ref{exam--anti-P-type-A}, while the latter corresponds to Example~\ref{exam--anti-P-type-B}.

    {\bf Step 2-1.}
    We first assume that the proper transform of $\Delta_0$ is in $C_3$.
    Since the chain \eqref{eqn:nonnormal} is obtained from $\Delta_0$ by taking blow-ups and each white vertex is not a $(-1)$-curve, the curves $\Delta_1$ and $\Delta_2$ must be $(-1)$-curves.
    Hence, by reordering the blow-ups if necessary, there exist $\alpha_1,\alpha_2\in\Z_{\geq0}$ satisfying the following:
    \begin{itemize}
        \item 
        The blow-ups $\psi_i$ ($1\leq i\leq\alpha_1+1$) occur over the same point $x_1 \in\Delta_0$, and the total transform of $\Delta_0$ in $W_{\alpha_1+1}$ forms the chain $1-[2^{\alpha_1},d+1]$.
        \item 
        The blow-ups $\psi_i$ ($\alpha_1+2\leq i\leq\alpha_1+\alpha_2+2$) occur over the same point of $\Delta_0$ that differs from $x_1$, and the total transform of $\Delta_0$ in $W_{\alpha_1+\alpha_2+2}$ forms the chain $1-[2^{\alpha_1},d+2,2^{\alpha_2}]-1$.
        \end{itemize}
        
        If $N=\alpha_1+\alpha_2+2$, we conclude that $\mu$ coincides with the one described in Example~\ref{exam--anti-P-type-A}~(1).

        If $N>\alpha_1+\alpha_2+2$, the morphism $W_{\alpha_1+\alpha_2+4}\to W_{\alpha_1+\alpha_2+2}$ is a composition of two blow-ups at the nodes of the chain lying on $(-1)$-curves, and the total transform of $\Delta_0$ becomes $[2]-1-[3,2^{\alpha_1-1},d+2,2^{\alpha_2-1},3]-1-[2]$. 
        For any $n>1$, the total transform of $\Delta_0$ in $W_{\alpha_1+\alpha_2+2n}$ forms the chain $C_1^n-1-C_3^n-1-C_2^n$ and $C_1^n-1-C_2^n$ is a conjugate pair.
        In other words, for any $1<n<(N-\alpha_1-\alpha_2)/2$, the chain $C_1^{n+1}-1-C_3^{n+1}-1-C_2^{n+1}$ is either of the following two:
        \[
        C_1^nR_1-1-L_2C_3^nR_1-1-L_2C_2^n, \quad C_1^nR_2-1-L_1C_3^nR_2-1-L_1C_2^n.
        \]
        Thus, we conclude that $\mu$ coincides with the one described in Example~\ref{exam--anti-P-type-A}~(2).

    {\bf Step 2-2.}
    Finally, we assume that the proper transform of $\Delta_0$ is not contained in $C_3$.
We may then assume that it is contained in $C_2$.
In particular, both $C_1$ and $C_2$ are non-empty.
Let us write $C_{i}=[b_{i1},\ldots,b_{il_i}]$ and denote $\Delta_{i}^{2}=-m_{i}$.
Then, we have $b_{2l_2}\ge d+1$ and either $m_1=1$ or $m_2=1$.
Note that the equation~\eqref{eqn:A_3} is equivalent to $l_3-m_1-m_2+8=\sum_{j=1}^{l_3}(b_{3j}-2)$.

First, assume that $m_1=1$ and $m_2\ge 2$.
Let $C'_{2}$ be the chain of rational curves (with no $(-1)$-curves) obtained by blowing down 
$$
C_{3}-m_2-C_2R_{1}^{-d}=[b_{31},\ldots,b_{3l_3},m_2,b_{21},\ldots,b_{2l_2-1},b_{2l_2}-d].
$$
Since the chain $C_1-1-C_3-m_2-C_2$ contracts to $\Delta_0$, 
the chain $C_1-1-C_3-m_2-C_2R_{1}^{-d}$ must contract to $[0]$.
Thus, $C_1$ and $C_2'$ are conjugate, and we have $C_2=C'_{2}$.
Suppose the blow-down $C_{3}-m_2-C_2R_{1}^{-d}\to \cdots \to C_2$ does not contract $\Delta_2$.
Then we can write $C_2=C_3-m_2-C_2^{(1)}$, where $C_2^{(1)}$ is the reslt of blowing down $C_2R_{1}^{-d}$.
Hence, $C_3-m_2-C_2^{(1)}R_{1}^{-d}=C_2R_{1}^{-d}$ blows down to $C_{2}^{(1)}$, and so we may replace $C_2$ by $C_{2}^{(1)}$.
Repeating this, we obtain a sequence $C_{2}^{(i)}$ for $i=0,\ldots,m$ such that:
\begin{itemize}
\item 
$C_2^{(0)}:=C_2$,
\item 
$C_2^{(i-1)}R_1^{-d}$ blows down to $C_{2}^{(i)}$, and
\item 
$C_3-m_2-C_2^{(i)}=C_{2}^{(i-1)}$.
\end{itemize}
Eventually, the blow-down process $C_{3}-m_2-C_2^{(m)}R_{1}^{-d}\to \cdots \to C_2^{(m)}$ contracts $\Delta_2$.
then, we can express: 
$$
C_3-m_2-C_2^{(m)}=[\cdots, e,2^{\beta}]-2-[2^{\alpha},d+1]
$$
for some $\alpha, \beta\ge 0$,
and 
$C_3-m_2-C_2^{(m)}R_{1}^{-d}$ blows down to $[\cdots, e-1]=C_2=[2^{\alpha},d+1]$.
This implies $e=d+2$ and $C_3=[2^{\alpha},d+2,2^{\beta}]$,
so in particular, $l_3=\alpha+\beta+1$.
From \eqref{eqn:A_3}, it follows that $\beta=d-\alpha-6$.
Now we consider $C_2$.
Next, observe that for each $i$, we can write: 
$$
[2^{\alpha},d+2, 2^{d-6-\alpha}]-2-C_{2}^{(i)}=C_{2}^{(i-1)}=C_{2}^{(i)}R_{1}-[2^{k_i},d+1]
$$ for some $k_i\ge 0$.
For $i=m$, we have
$[2^{\alpha},d+2, 2^{d-5},d+1]=[2^{\alpha},d+2,2^{k_m},d+1]$.
Thus, $k_m=d-5$.
Inductively, $k_{i}=d-5$ for each $i$, so the morphism $\mu$ is an anti-P-resolution constructed as in Example~\ref{exam--anti-P-type-B}~(1).

Next, we assume that $m_1\ge 2$ and $m_2=1$.
Let $C_{2}^{(1)}$ be the chain obtained by blowing down $C_2R_1^{-d}$.
Since $C_1-m_1-C_3-1-C_2R_1^{-d}$ blows down to $[0]$, the two chains $C_1-m_1-C_3$ and $C_2^{(1)}$ must be conjugate.
However, since $C_1$ and $C_2$ are conjugate and the blow-down $C_1-1-C_2\to \cdots \to [0]$ factors through $[2^d]-1-[d+1]$, the chain $L_{2}^{-d}C_1-1-C_2R_1^{-d}$ also contracts to $[0]$, where
$$
L_2^{-d}[2^{d-1},b_1,b_2,\cdots]:=[b_2,\cdots].
$$
Thus, $L_2^{-d}C_1$ and $C_{2}^{(1)}$ are conjugate.
But then $L_{2}^{-d}C_1=C_1-m_1-C_3$, which contradicts  the number of rational curves.

Finally, assume that $m_1=m_2=1$.
Let $C_{2}^{(1)}$ be the chain obtained by blowing down  $C_2R_1^{-d}$,
and let $E$ be the rightmost component of $C_2^{(1)}$.
Since $C_1-1-C_3-1-C_2R_{1}^{-d}$ contracts to $[0]$, the chain $C_3-1-C_2^{(1)}$ must contract to $C_2$.
If this does not contract $E$, then $-E^{2}\ge d+2$ and $C_3-1-C_2^{(1)}R_1^{-d}$ contracts to $C_{2}^{(1)}$.
Hence, we can replace $C_2$ by $C_{2}^{(1)}$.
Iterating this process, we get chains $C_{2}^{(i)}$, $i=1,\ldots,m$ such that
\begin{itemize}
    \item
    $C_{2}^{(0)}:=C_2$,
    \item 
    $C_{2}^{(i-1)}R_1^{-d}$ contracts to $C_{2}^{(i)}$,
    \item 
    $C_3-1-C_2^{(i)}$ contrtacts to $C_2^{(i-1)}$,  
\end{itemize}
 and eventually the last blow-down $C_3-1-C_2^{(m)} \to \cdots \to C_2^{(m-1)}$ contracts the rightmost component of $C_2^{(m)}$.
Hence, we can write 
$$
C_3-1-C_2^{(m)}=C_{2}^{(m-1)}R_{1}R_2^{\beta}R_1\bm{R}-1-\bm{L}[2]
$$
for some $\beta\ge 0$ and operations $\bm{L}=L_{i_l}\cdots L_{i_1}$, $\bm{R}=R_{3-i_1}\cdots R_{3-i_l}$.
Since we can write $C_{2}^{(m-1)}=C_2^{(m)}R_1-[2^{\alpha-1},d+1]$ for some $\alpha \ge 0$, 
we have $C_3=C_2^{(m)}R_1-[2^{\alpha-1}, d+1]R_1R_{2}^{\beta}R_1\bm{R}=\bm{L}[3,2^{\alpha-1},d+2,2^{\beta-1},3]\bm{R}$.
From \eqref{eqn:A_3}, we obtain $\beta=d-\alpha-5$.
Now we show that $C_{2}^{(i-1)}=C_2^{(i)}R_1-[2^{d-5}, d+1]$ for each $i=1,\ldots,m-1$.
Since the chain
$$
C_3-1-C_{2}^{(m-1)}=\bm{L}[3,2^{\alpha-1},d+2,2^{d-\alpha-6},3]\bm{R}-1-\bm{L}[3,2^{\alpha-1},d+1]
$$
blows down to $C_{2}^{(m-2)}$, the claim holds for $i=m-1$.
Similar argument shows the claim by descending induction.
Hence, $\mu$ is an anti-P-resolution as in Example~\ref{exam--anti-P-type-B}~(2).
\end{proof}

We introduce the following useful notions for anti-P-resolutions.

\begin{defn}[Length of anti-P-resolution] \label{defn--length} 
Let $(x\in W,\frac{1}{2}B)$ be a germ of a cone singularity of type $(d;k_1,k_2)$, and let $\mu\colon W^{-}\to W$ be an admissible anti-P-resolution.
Define $\mathcal{B}$ as the set of effective divisors $D$ on $W$ satisfying the following conditions:
\begin{itemize}
    \item $D\sim B$,
    \item $(x\in W,\frac{1}{2}D)$ is a germ of a cone singularity of type $(d;k'_1,k'_2)$ for some $k'_1, k'_2$,
    \item There exists an effective divisor $D^{-}$ on $W^{-}$ such that   
    $\mu\colon (W^{-},\frac{1}{2}D^{-})\to (W,\frac{1}{2}D)$ is log crepant.
\end{itemize}
Define $l_{i}:=\min_{D\in\mathcal{B}}k'_i$.
Then, the values $l_1$ and $l_2$ are achieved simultaneously by the same divisor $D$.
We call $(l_1,l_2)$ the {\it length} of $\mu$.
Clearly, we have $l_i\le k_i$.
\end{defn}

\begin{rem}\label{rem--antiP-and-N-resol}
 Tevelev and Urz\'ua \cite{TU} introduced a notion of N-resolution, which is also obtained by taking the anti flip of a certain $\mathbb{Q}$-Gorenstein deformation over a smooth curve for an M-resolution (see \cite{BC}) with a certain condition.
 We note that their N-resolution does not appear in the classification of admissible anti-P-resolutions in Theorem \ref{thm--anti-P-k_1-k_2}.
 Note that an N-resolution uniquely exists for an arbitrary M-resolution.
 Considering anti-P-resolutions of normal type of length $(2d-8,0)$ and $(2d-9,0)$ for the cone singularity $\frac{1}{d}(1,1)$, it is easy to see that there are several anti-P-resolutions of normal type that are not isomorphic to each other.
 In this respect, anti-P-resolutions are different from N-resolutions.
\end{rem}

\begin{defn}[Good resolution]\label{defn--good-resolution}
    Let $(x\in W, \frac{1}{2}B)$ be a germ of a cone singularity of type $(d;k_1,k_2)$ with $d\ge 5$, and let $\mu\colon (W^{-}, \frac{1}{2}B^{-}) \to (W, \frac{1}{2}B)$ be an admissible anti-P-resolution of length $(k_1,k_2)$. 
    Let $\pi\colon W^{+}\to W$ be the minimal resolution, and let $B^{+}$ be an effective divisor on $W^{+}$ such that $\pi \colon (W^{+}, \frac{1}{2}B^{+})\to (W, \frac{1}{2}B)$ is log crepant.
    Then, we define a birational morphism $\psi\colon W_{\mathrm{good}}\to W^{+}$
with an effective divisor $B_{\mathrm{good}}$ on $W_{\mathrm{good}}$
    as follows:
    Let $\Delta_0$ be the exceptional divisor of $\pi$.
    Then, $B_0:=B^{+}-2\Delta_{0}$ coincides with the proper transform of $B$, which is the branch divisor of the normalized base change of the double covering $X\to W$ by $\pi$.
    Let 
    $$
    W_{N}\xrightarrow{\psi_N}\cdots \xrightarrow{\psi_{2}}W_1\xrightarrow{\psi_1} W_0:=W^{+}
    $$
    be a sequence of blow-ups that is a minimal log resolution of $(W,B)$.
    More precisely, 
    we first take successive blow-ups at singular points of $B_0$
    so that the proper transform of $B_0$ becomes smooth.
    If the proper transform of $B_0$ is tangent to the total transform of $\Delta_{0}$ at some point, we further take blow-ups twice at this point.
    Define $W_{\mathrm{good}}:=W_{N}$ and $\psi:=\psi_{1}\circ \cdots \circ \psi_{N}$.
    Let $B_{\mathrm{good}}\in |2L_{\mathrm{good}}|$ be the branch divisor of the the normalized base change of the double covering $X\to W$ by $\pi$, which can be computed by the even resolution of $B_{0}$.
    Then, $B_{\mathrm{good}}$ is smooth and has a decomposition $B_{\mathrm{good}}=B'_{\mathrm{good}}+E$,
    where $B'_{\mathrm{good}}$ is the proper transform of $B$ and $E$ is a sum of disjoint $(-2)$-curves appeared in the blow-up at a smooth point of the proper transform of $B_0$.
    We call $\psi$ the {\em good resolution} of $(W,B)$ with respect to $\mu$.
    
    If $D$ is an effective divisor on $W$ belonging to $\mathcal{B}$ defined in Definition~\ref{defn--length}, we can similarly define the associated effective divisor $D_{\mathrm{good}}$ and $D'_{\mathrm{good}}$ on $W_{\mathrm{good}}$ as
    the even resolution of the proper transform $D_0$ of $D$ on $W_0$ and the proper transform of $D$ on $W_{\mathrm{good}}$, respectively.
    In this case, we have $D_{\mathrm{good}}\sim B_{\mathrm{good}}$ and $D_{\mathrm{good}}=D'_{\mathrm{good}}+E$.

\end{defn}

\section{Deformations of normal stable Horikawa surfaces}\label{sec:deformation}

The aim of this section is to prove Theorems \ref{intro-qGsmHor}, \ref{thm--main--ii}, \ref{thm--stratification--intro} and \ref{thm--connectedness}. 

Let $X$ be a normal stable Horikawa surface that is $\Q$-Gorenstein smoothable.
By Proposition~\ref{prop:non-std_Horikawa}~(2) and Theorem~\ref{thm:standard_Horikawa}, we have 
 $q(X)=0$ and $\chi(\O_X)=p_g(X)+1\ge 4$.
 
Assume first that $X$ is standard, that is, it is Gorenstein and its canonical linear system is not composed with a pencil.
In this case, Proposition~\ref{prop:standard_Horikawa} directly implies Theorem~\ref{intro-qGsmHor}~(1).

Now, assume that $X$ is non-standard.
Non-standard Horikawa surfaces with only $\Q$-Gorenstein smoothable singularities
have already been classified in Theorem~\ref{thm:classification_non-standard_Horikawa}.
The correspondence between this classification  and Theorem~\ref{intro-qGsmHor} is as follows:

\begin{itemize}
\item
Theorem~\ref{thm:classification_non-standard_Horikawa}~(8) corresponds to Theorem~\ref{intro-qGsmHor}~(2).
\item
Theorem~\ref{thm:classification_non-standard_Horikawa}~(5), (6) and (7) correspond to Theorem~\ref{intro-qGsmHor}~(3).
    \item 
    Theorem~\ref{thm:classification_non-standard_Horikawa}~(1) corresponds to Theorem~\ref{intro-qGsmHor}~(4).
    \item 
    Theorem~\ref{thm:classification_non-standard_Horikawa}~(4) corresponds to Theorem~\ref{intro-qGsmHor}~(5).
    \item
    Theorem~\ref{thm:classification_non-standard_Horikawa}~(2) corresponds to Theorem~\ref{intro-qGsmHor}~(6).
     \item 
    Theorem~\ref{thm:classification_non-standard_Horikawa}~(3) corresponds to Theorem~\ref{intro-qGsmHor}~(7).
\end{itemize}
It is easy to see that a good involution of any stable Horikawa surface $X$ listed in Theorem \ref{intro-qGsmHor} acts on $H^0(X,\omega_X)$ trivially since $H^0(W,\omega_W)=0$.

In this section, we use the following notations unless otherwise stated:
\begin{note} \label{note--section 9}
Let $X$ be a Horikawa surface with a good involution $\sigma$ which is standard or non-standard classified in Theorem~\ref{thm:classification_non-standard_Horikawa}.
Let $W$ be the quotient of $X$ by $\sigma$, and let $B\in |L^{[2]}|$ be the branch divisor of $X\to W$.
If $(W, \frac{1}{2}B)$ has a cone singularity of type $(d;k_1,k_2)$ (that is, $X$ is standard of type $(d)'$ or special Lee-Park type with $p_g(X)=d+2$), then let $(W^{+}, \frac{1}{2}B^{+})\to (W, \frac{1}{2}B)$ be the minimal resolution which is log crepant.
In this case, the multiplicity $m(X):=\mathrm{mult}_{E}(B^{+})$ of $B^{+}$ along the exceptional divisor $E$ satisfies
$$
\left\lceil \frac{2d-4}{d} \right\rceil \le m(X) \le 2.
$$
We call $m(X)$ the {\em multiplicity} of $X$ for simplicity.

If $\mathscr{X}\to C$ is a $\Q$-Gorenstein smoothing of a Horikawa surface $X$, let $\mathscr{W}$ be the quotient of $\mathscr{X}$ by the good involution $\sigma$ on $\mathscr{X}$,
and its branch divisor is denoted by $\mathscr{B}\in |\mathscr{L}^{[2]}|$.
Note that $W=\mathscr{W}_0$.
\end{note}

\subsection{Gorenstein Horikawa surfaces with a canonical pencil}\label{subsec:eliip-cone--smoothable}


First, we deal with Gorenstein Horikawa surfaces with a canonical pencil that correspond to Theorem \ref{intro-thm:involution}.
To show Theorem \ref{intro-thm:involution}, we have to show that such Horikawa surfaces admit good involutions.
The following proposition asserts the existence of good involutions and their smoothability.

\begin{prop}\label{prop--canonical--pencil--invol}
    Let $X$ be a stable normal Gorenstein Horikawa surface with a canonical pencil.
    Then, the following hold:

 \begin{itemize}
        \item[$(1)$] 
        $p_{g}(X)=4$ and $X$ can be constructed as described in Construction~\ref{construction}~(8). 
         In particular, $X$ has a good involution that trivially acts on $H^0(X,\omega_X^{[2]})$.
        \item[$(2)$] 
        $X$ is smoothable.
        More precisely, it has a partial smoothing to canonical Horikawa surfaces of types $(2')$ and $(0)$.
 \end{itemize}
    
\end{prop}

\begin{proof}
(1): 
Let $X$ be a Gorenstein stable normal Horikawa surface such that $|K_X|$ is composed with a pencil.
By Theorem~\ref{normalstable}, $X$ belongs to case (1.4) in Theorem~\ref{normalstable}.
In particular, the minimal resolution $\widetilde{X}$ of $X$ admits an elliptic fibration $\widetilde{X}\to C$.
Since $X$ is Gorenstein, 
we have $\widetilde{X}=Y$ and $K_Y+\Delta=\pi^{*}K_X$ as in the notation in Theorem~\ref{normalstable} 
Moreover, since $X$ is a Horikawa surface, it satisfies $K_X^{2}=p_g(X)=4$.
This implies that the horizontal part of $\Delta$ consists of a section $S$ with $S^2=-4$.

Now, we claim that the bicanonical linear system $|2K_X|$ defines a double cover $X\to W\subset \mathbb{P}^{8}$ onto an elliptic cone of degree $8$.
To prove this, we first show that $|2(K_{Y}+\Delta)|$ is basepoint-free.
For each fiber $F$, it follows easily that the restricted complete linear system $|2(K_Y+\Delta)|_{F}|=|2S|_{F}|$ is basepoint-free.
Thus, it suffices to prove that the restriction map 
$$
H^{0}(2(K_Y+\Delta))\to H^{0}(2(K_Y+\Delta)|_{F})
$$
is surjective.
Equivalently, it suffices to show that the induced map 
$$
H^{1}(2(K_Y+\Delta)-F)\to H^{1}(2(K_Y+\Delta))
$$
is injective.
Since $(K_Y+\Delta)|_{S}=K_{S}=0$,
it follows that 
$$
H^{1}((2(K_Y+\Delta)-F)|_{S})\cong  H^{1}(2(K_Y+\Delta)|_{S})\cong \mathbb{C}.
$$
Thus, it suffices to show that 
$$
H^{1}(2(K_Y+\Delta)-F-S)\cong H^{1}(2(K_Y+\Delta)-S).
$$
It is easy to check that both $K_{Y}+\Delta$ and $K_Y+\Delta-F$ are big, $\Z$-positive, and numerically trivial over $\Delta-S$.
Thus, applying \cite[Theorem~4.1]{Enokizono}, we obtain 
\begin{align*}
H^{1}(2(K_Y+\Delta)-F-S)&=H^{1}((2(K_Y+\Delta)-F-S)|_{\Delta-S})=H^{1}(\O_{\Delta-S}), \\
H^{1}(2(K_Y+\Delta)-S)&=H^{1}((2(K_Y+\Delta)-S)|_{\Delta-S})=H^{1}(\O_{\Delta-S}).
\end{align*}
Therefore, we conclude that $|2(K_Y+\Delta)|$ is basepoint-free.

Since $(K_Y+\Delta)\cdot F=1$ and $K_{Y}+\Delta$ is numerically trivial over $\Delta$, the linear system $|2(K_Y+\Delta)|$ defines a generically finite morphism $Y\to W\subset \mathbb{P}^{8}$ of degree $2$ onto its image, and this morphism factors through $X$. 
Since $K_X$ is ample, it follows that $X\to W$ is a double cover.
Let $\widehat{X}\to X$ be the minimal resolution of the simple elliptic singularity corresponding to $S$.
The covering involution $\sigma\colon X\to X$ induces an involution of $\widehat{X}$ over $C$ that fixes $S$.
By the classification of singular fibers of elliptic fibrations with an involution (Proposition~\ref{prop_inv_ell}), it follows that the quotient $\widehat{W}:=\widehat{X}/\sigma$ is a $\mathbb{P}^1$-bundle over $C$ and the image of $S$ on $\widehat{W}$ is a section $\Delta_{0}$ with $\Delta_{0}^{2}=-8$.
It is easy to check that $\widehat{W}\cong \mathbb{P}_{C}(\O_{C}\oplus \mathcal{M}^{-1})$ for some line bundle $\mathcal{M}$ of degree $8$, and the branch divisor on $\widehat{W}$ equals $\Delta_0+\widehat{B}'$, where $\widehat{B}'$ is a divisor as described in Construction~\ref{construction}~(8).
In particular, $W$ is an elliptic cone of degree $8$. 
Note that the double cover structure $\widehat{X}\to \widehat{W}$ determines a line bundle $\mathcal{L}$ on $C$ with $\mathcal{L}^{\otimes 2}\cong \mathcal{M}$ as described in Construction~\ref{construction}~(8).

(2):
We use the notations introduced in Construction~\ref{construction}~(8).
The linear system $|\mathcal{L}|$ defines an embedding $C\subset \mathbb{P}^3$.
The pencil of quadrics in $\mathbb{P}^3$ containing $C$ includes a quadric $Q$ of rank $3$, which is isomorphic to $\overline{\Sigma}_{2}$.
Consider the second Veronese embedding $\mathbb{P}^3\hookrightarrow \mathbb{P}^{9}$ and take a hyperplane $H\subset \mathbb{P}^9$ such that $C=Q\cap H$.
By taking projective cones, we obtain 
$$
\mathrm{Cone}(C)=\mathrm{Cone}(Q)\cap \mathrm{Cone}(H)
$$
in $\mathbb{P}^{10}$.
Since $\O_{\mathbb{P}^{9}}(1)|_{C}\cong \mathcal{M}$, there exists a natural isomorphism $\mathrm{Cone}(C)\cong W$.
Moreover, since the cone over a hyperplane is also a hyperplane, $\mathrm{Cone}(H)$ is a hyperplane in $\mathbb{P}^{10}$.
Now, consider a one-parameter family of hyperplanes $\{H_{t}\}_{t\in T}$ in $\mathbb{P}^{10}$ such that $H_0\cong \mathrm{Cone}(H)$ and $\mathrm{Cone}(Q)\cap H_{t}\cong Q$ for each $t\neq 0$.
This induces a partial smoothing family $g\colon \mathscr{W}\to T$ from $W$ to $Q$, where 
$$
\mathscr{W}:=\{(x,t)\in \mathrm{Cone}(Q)\times T\ |\ x\in H_t\}.
$$

Next, we show that there exist a $\Q$-Cartier divisorial sheaf $\mathscr{L}$ on $\mathscr{W}$ and an effective divisor $\mathscr{B}\in |\mathscr{L}^{[2]}|$ such that $\mathscr{L}|_{W}\cong L$ and $\mathscr{B}|_{W}=B$.
Since $\O_{\mathbb{P}^{9}}(1)|_{Q}\cong \O_{Q}(-K_{Q})$, the cone $\mathrm{Cone}(Q)$ is obtained by the contraction 
$$
q\colon P:=\mathbb{P}_{Q}(\O_{Q}\oplus \O_{Q}(K_Q))\to \mathrm{Cone}(Q),
$$
which contracts the tautological section $Q_0$.
Let $p\colon P\to Q$ be the natural projection
and define $D:=q_{*}p^{*}\O_{\mathbb{P}^{3}}(1)$.
Then, it follows that $K_{\mathrm{Cone}(Q)}\sim -4D$ since $K_{\mathrm{Cone}(Q)}=p_{*}K_{P}$, $K_{P}\sim -2Q_{0}+2p^{*}K_{Q}$ and
$-K_{Q}\sim \mathcal{O}_{\mathbb{P}^{3}}(2)|_{Q}$.
By the adjunction formula, we obtain
$$
K_{\mathscr{W}}=(K_{\mathrm{Cone}(Q)\times T}+\mathscr{W})|_{\mathscr{W}}\sim_{T} r^{*}(K_{\mathrm{Cone}(Q)}+q_{*}p^{*}\O_{\mathbb{P}^9}(1))\sim -2r^{*}D,
$$
where $r\colon \mathscr{W}\to \mathrm{Cone}(Q)$ is the first projection.
In particular, $2r^{*}D$ is Cartier.
Define $\mathscr{L}:=\O_{\mathscr{W}}(3r^{*}D)$.
Since $\mathscr{W}$ is klt and $3r^{*}D$ is $\Q$-Cartier, $\mathscr{L}$ is Cohen-Macaulay by \cite[Corollary~5.25]{KoMo}.
As Weil divisors, we have
$$
L=(q|_{\widehat{W}})_{*}\widehat{L}=(q|_{\widehat{W}})_{*}(p|_{\widehat{W}})^{*}\O_{\mathbb{P}^{3}}(3)=(q_{*}p^{*}\O_{\mathbb{P}^3}(3))|_{W}=3D|_{W}.
$$
Thus, $\mathscr{L}|_{W}\cong L$.

Since $-K_W$ is ample and Cartier, 
 \cite[Theorem 1.3]{Enokizono} implies  $H^1(-3K_W)=0$. 
Thus, there exists an effective divisor $\mathscr{B}\in |-3K_{\mathscr{W}}|=|\mathscr{L}^{[2]}|$ satisfying $\mathscr{B}|_{W}=B$.
Shrinking $T$ around $0$ and perturbing the section of $g_*\mathcal{O}_{\mathscr{W}}(-3K_{\mathscr{W}})$ defining $\mathscr{B}$, we may assume that a general fiber of $\mathscr{B}$ over $T$ is smooth and does not contain the vertex of $Q$.
Then, by defining $\mathscr{X}:=\mathbf{Spec}_{\mathscr{W}}(\mathcal{O}_{\mathscr{W}}\oplus\mathscr{L}^{[-1]})$ as the double cover branched along $\mathscr{B}$, we obtain that $\mathscr{X}_{0}=X$ and general fibers of $\mathscr{X}$ are canonical Horikawa surfaces of type $(2)'$.
Thus, we conclude that $X$ is smoothable to smooth Horikawa surfaces of type $(0)$.
\end{proof}

This proposition together with Proposition \ref{prop:standard_Horikawa} shows Theorem \ref{intro-thm:involution}.

\subsection{Non-standard Horikawa surfaces with $p_g=3$}\label{subsec:p_g=3--smoothable}

Let $X$ be a normal stable Horikawa surface that belongs to the classification in Theorem~\ref{thm:classification_non-standard_Horikawa} with $\chi(\O_X)=4$.
By Propositrion~\ref{prop:non-std_Horikawa}~(2), $q(X)=0$ and $p_g(X)=3$.
\begin{prop} \label{prop--smoothing-p_g=3}
Let $X$ be a non-standard Horikawa surface described in  Theorem~\ref{thm:classification_non-standard_Horikawa} with $p_g(X)=3$.
Then, $X$ is $\Q$-Gorenstein smoothable.
\end{prop}

\begin{proof}
Note that $X$, $W$ and $B$ are explicitly described as in Construction~\ref{construction}.
In particular, $\rho(W)=1$.
It is easy to see that $3L+4K_W\sim 0$, and both $-K_W$ and $L$ are ample.
Thus, applying the Kodaira vanishing theorem, we obtain
$$
H^1(W,\mathcal{O}_W(B))=H^1(W,\mathcal{O}_W)=H^2(W,\mathcal{O}_W)=0.
$$
Furthermore, $W$ has only T-singularities of types $\frac{1}{4}(1,1)$ and $\frac{1}{25}(1,4)$.
It is easily seen that there exist integers $l_1$ and $l_2$ such that $l_1K_W\sim L$ and $l_2K_W\sim L$ locally around $\frac{1}{4}(1,1)$ and $\frac{1}{25}(1,4)$, respectively.
By Proposition \ref{prop--smoothing--quotient}, there exists a projective flat morphism $f\colon \mathscr{W}\to C$ from a normal $\mathbb{Q}$-Gorenstein threefold to a smooth curve $C$ with a closed point $0\in C$, along with an effective relative Mumford divisor $\mathscr{B}$, such that there exists a divisorial sheaf $\mathscr{L}$ satisfying
\begin{itemize}
    \item 
$\mathscr{L}_0\cong L$ and $(\mathscr{W}_0,\mathscr{B}_0)\cong (W,B)$,
    \item $\mathscr{L}^{[2]}\cong\mathcal{O}_{\mathscr{W}}(\mathscr{B})$, 
    \item
    $\mathscr{W}_c$ is smooth for any $c\in C\setminus\{0\}$.
\end{itemize}
Since $W$ has only Wahl singularities, we conclude that $\mathscr{W}_c\cong\mathbb{P}^2$.
Thus, applying Proposition \ref{prop--smoothing--quotient} again, we see that $X$ is $\mathbb{Q}$-Gorenstein smoothable.
\end{proof}

\begin{rem}\label{rem--existence--toric}
Note that non-standard Horikawa surfaces $X$ described in  Theorem~\ref{thm:classification_non-standard_Horikawa} with $p_g(X)=3$ form a non-empty irreducible family.
It is easy to check this fact for surfaces as Theorem \ref{thm:classification_non-standard_Horikawa} as follows.
Note that $W$ is all isomorphic to $\mathbb{P}(1,4,25)$ and $B\sim\mathcal{O}_{\mathbb{P}(1,4,25)}(80)$.
Let $x,y,z$ be the canonical weighted homogeneous coordinates of $\mathbb{P}(1,4,25)$.
Take the minimal resolution $\pi\colon\widetilde{W}\to W$ and note that $\widetilde{W}$ is toric again.
Let $E$ be the exceptional divisor of $\frac{1}{4}(1,1)$.
It is easy to see that $\lfloor\pi^*\mathcal{O}(80)\rfloor-E-D$ is nef but $(\lfloor\pi^*\mathcal{O}(80)\rfloor-E)\cdot D<0$, where $D$ is the proper transform of $(x=0)$.
By \cite[Theorem 6.3.12]{CLS}, we see that $\lfloor\pi^*\mathcal{O}(80)\rfloor-E-D$ is globally generated.
For any general $D'\in|\lfloor\pi^*\mathcal{O}(80)\rfloor-E-D|$, we see that the double cover of $W$ branch along $\pi_*(D+D')$ is a stable Horikawa surface that belongs to Theorem~\ref{thm:classification_non-standard_Horikawa} with $p_g(X)=3$.
For other cases, a similar argument shows the assertion.
\end{rem}

\subsection{Horikawa surfaces of general Lee-Park type}\label{subsec:LPgeneral--smoothable}

In this subsection, we consider Horikawa surfaces of general Lee-Park type described in Theorem~\ref{thm:classification_non-standard_Horikawa} (5). 

First, we deal with the case where $\rho(W)=2$.
For Horikawa surfaces of general Lee-Park type $\mathrm{I}$, we describe fundamental properties.

\begin{prop}\label{prop--construction--lee--park--i} 
Let $X$ be a Horikawa surface of general Lee-Park type $\mathrm{I}$ with $p_g\ge 4$.
Then, the following hold:

\begin{itemize}
    \item[$(1)$]
    $-K_W$ is ample and $h^{0}(-K_W)\ge 9$.
    \item[$(2)$]
    $|2L|$ is basepoint-free and $H^{1}(2L)=0$.
\end{itemize}

\end{prop}

\begin{proof}
We use the notation introduced in Definition~\ref{defn--evenresol}.
Then, the partial even resolution process 
$$
\psi=\psi_{1}\circ \cdots \circ \psi_{\chi-3}\colon \widehat{W}=W_{\chi-3} \to \cdots \to W_{1}\to W_{0}=\Sigma_{\chi}
$$
is described as in the diagram. 
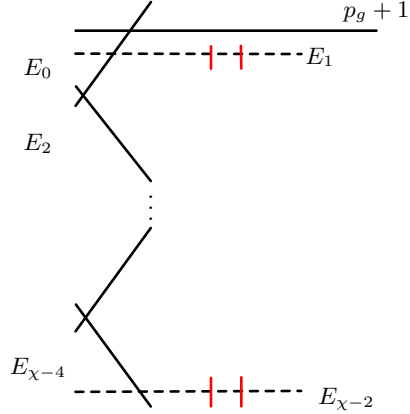
\begin{figure}[H]
\begin{tikzpicture}[line cap=round,line join=round,>=triangle 45,x=1cm,y=1cm]
\clip(-3,-4) rectangle (3,3);
\draw [line width=1pt] (-2,2)-- (2,2);
\draw [line width=1pt] (-0.9971203549676338,2.3813912622704274)-- (-2,1);
\draw [line width=1pt] (-1.993987670038532,1.2578064625260739)-- (-1,0);
\draw [line width=1pt] (-1,-0.625)-- (-2,-2);
\draw [line width=1pt] (-1.995086292612181,-1.6430148558396644)-- (-1,-3);
\draw [line width=1pt,dashed] (-2.0007087912525727,1.695060673978309)-- (1.0025873711331752,1.6893940774455056);
\draw [line width=1pt,dashed] (-2.001321734541306,-2.8052652720107796)-- (0.9988115203214603,-2.7928451752996906);
\draw [line width=1pt,red] (-0.19873109382112394,1.7913928150359655)-- (-0.19951474350404946,1.4937087696580595);
\draw [line width=1pt,red] (0.19793066347510693,1.7970594115687688)-- (0.20359726000791026,1.4967297953301926);
\draw [line width=1pt,red] (-0.19951474350404946,-2.6193046666941124)-- (-0.19951474350404946,-2.998347081416764);
\draw [line width=1pt,red] (0.20372186790302457,-2.611239934465971)-- (0.20372186790302457,-2.998347081416764);
\begin{scriptsize}
\draw [fill] (2,2.25) node {$p_g+1$};
\draw [fill] (-2.5,1.5) node {$E_0$};
\draw [fill] (1.25,1.65) node {$E_1$};
\draw [fill] (-2.5,0.5) node {$E_2$};
\draw [fill] (-2.5,-1.5) node {$ $};
\draw [fill] (-2.5,-2.5) node {$E_{\chi -4}$};
\draw [fill] (1.6,-2.9) node {$E_{\chi -2}$};
\draw [color=black] (-1,-0.25) node {$\vdots$};
\end{scriptsize}
\end{tikzpicture}
\end{figure}
Define $D:=\psi_{*}E_1$ and $E:=\psi_{*}E_{\chi-3}$.
By construction, $D$ and $E$ intersect at the singularity of $W$.
Since the divisors
 $\Gamma_1+\sum_{i=1}^{\chi-2}E_i$ and $\Gamma_1+\sum_{i=2}^{\chi-3}E_i$ are both negative definite, it follows that $D^2<0$ and $E^2<0$.
 In particular, each of $D$ and $E$ generates an extremal ray in the numerical class group of $w$.
Since $\rho(W)=2$, the curves $D$ and $E$ generate the effective cone $\mathrm{Eff}(W)$.
Pursuing the even resolution process, we have 
\begin{align*} 
-K_{\widehat{W}}&=\psi^{*}(2\Delta_0+(\chi+2)\Gamma)-\left(E_1+\sum_{i=2}^{\chi-3}(i-1)E_i  \right). 
\end{align*}
Hence, we conclude that
$$
h^{0}(-K_W)\ge h^{0}(-K_{\widehat{W}})\ge h^{0}(2\Delta_0+(\chi+1)\Gamma)\ge 9.
$$
Furthermore, it is straightforward to verify that $-K_{\widehat{W}}$ is big.
Thus, $-K_W=\delta_{*}(-K_{\widehat{W}})$ is also big.
Since $W$ is klt, we have
$$
-K_{W}\cdot D=-\delta^{*}K_{W}\cdot E_1=-K_{\widehat{W}}\cdot E_{1}+a \Gamma_1\cdot E_{1} = 1+a>0,
$$
where $a$ is the discrepancy of $\Gamma_1$ with respect to $W$.
Similarly, we obtain $-K_{W}\cdot E>0$.
Thus, we conclude that $-K_{W}$ is ample.

Since $K_{W}+L$ is ample, it follows that
 $L$ is also ample.
Moreover, by the definition of $D$ and $E$, we have $L\cdot D=L\cdot E=1$.
Since $L$ is Cartier, we obtain $L^{2}\ge 1$.
Applying the Kawamata-Viehweg vanishing theorem to the ample divisor $2L-K_{W}$, we deduce that $H^{1}(2L)=0$.
By the Reider-type theorem (cf.\ \cite[Corollary 5.11, see also 5.1 and 5.9]{Enokizono}), we conclude that the linear system $|2L|$ is basepoint-free.
\end{proof}

\begin{prop}\label{prop--smoothing--lee--park--i}
    Let $X$ be a Horikawa surface of general Lee-Park type I.
    Then, $X$ is $\mathbb{Q}$-Gorenstein smoothable.
\end{prop}
\begin{proof}
First, we verify that $(W,\frac{1}{2}B)$ satisfies all the conditions of Proposition \ref{prop--smoothing--quotient}.
The conditions (a), (e) and (f) are immediately satisfied since $B\in |2L|$ is the branch divisor of the double cover $X\to W$ and $K_X$ is ample.
Moreover, since $L$ is Cartier, condition (b) at the singularity of $W$ holds trivially.
Conditions (c) and (d) follow from Proposition~\ref{prop--construction--lee--park--i}.
Thus, we conclude that $(W,\frac{1}{2}B)$ satisfies all the conditions of Proposition \ref{prop--smoothing--quotient}.
Therefore, there exists a $\mathbb{Q}$-Gorenstein family $f\colon \mathscr{W}\to C$ with an effective $\mathbb{Q}$-Cartier Weil divisor $\mathscr{B}$ on $\mathscr{W}$ such that for some  closed point $0\in C$, 
 we have $(\mathscr{W}_0,\mathscr{B}_0)\cong (W,B)$.
 Furthermore, for any $c\in C\setminus\{0\}$, the fiber $\mathscr{W}_c$ is a smooth rational surface with $\rho(\mathscr{W}_c)=2$, since $W$ has a unique singularity, which is Wahl.
Since $-K_W$ is ample, it follows from \cite[Proposition 1.41]{KoMo} that $-K_{\mathscr{W}_c}$ is also ample for any general $c\in C$.
Thus, by the classification of smooth del Pezzo surfaces, we conclude that $\mathscr{W}_c\cong\Sigma_1$ or $\Sigma_0$. 
Moreover, we see that $\mathscr{B}_c$ is basepoint-free, and since $H^1(W,\mathcal{O}_W(B))=0$, we may assume that $\mathscr{B}_c$ is smooth.
By taking the double cover $\mathscr{X}$ of $\mathscr{W}$ along $\mathscr{B}$, we obtain a partial $\mathbb{Q}$-Gorenstein smoothing of $X$ to $\mathscr{X}_c$, which is obtained as the double cover of $\mathscr{W}_c$ branched along $\mathscr{B}_c$.
Since $\mathscr{X}_c$ is smooth, it follows that $X$ is $\mathbb{Q}$-Gorenstein smoothable.
\end{proof}

Using Proposition \ref{prop--construction--lee--park--i} again, we can describe the moduli boundary that consists of Horikawa surfaces of general Lee-Park type as follows.

\begin{rem}\label{rem-existence--LPI}
    Each $l\in\mathbb{Z}_{\ge4}$, Horikawa surfaces of general Lee-Park type $\mathrm{I}$ with geometric genus $l$ form a non-empty irreducible family.
    Indeed, it is easy to see by using the notion of even resolution (cf.~Definition \ref{defn--evenresol}) that the surfaces $\widehat{W}$ as in the proof of Proposition~\ref{prop--construction--lee--park--i} form an irreducible family.
    Thus, $W$ also forms a non-empty irreducible family.
    Furthermore, $H^1(W,\mathcal{O}_W(B))=0$ and $B$ is basepoint-free, we obtain the assertion.
\end{rem}

Next, we consider $\Q$-Gorenstein deformations of 
 Horikawa surfaces of general Lee-Park type $\mathrm{I}^{*}$.
 In this case, $p_g\ge 5$ automatically holds (see Construction \ref{construction}).
The fundamental properties of these surfaces are as follows:

\begin{prop}\label{prop--construction--lee--park--i*}
Let $X$ be a Horikawa surface of general Lee-Park type $\mathrm{I}^{*}$ with $p_g\ge 5$.
Then, the following hold:

\begin{itemize}
    \item[$(1)$]
    $-K_W$ is big but not nef and $h^{0}(-K_W)\ge 9$. 
    \item[$(2)$]
    $|2L|$ is basepoint-free and $H^{1}(2L)=0$.
\end{itemize}

\end{prop}

\begin{proof}
We use the notation introduced in Definition~\ref{defn--evenresol}. 
Then, the partial even resolution process 
$$
\psi=\psi_{1}\circ \cdots \circ \psi_{\chi-3}\colon \widehat{W}=W_{\chi-3} \to \cdots \to W_{1}\to W_{0}=\Sigma_{\chi}
$$
proceeds as in the diagram.

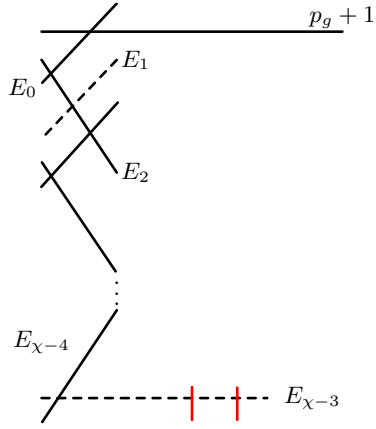
\begin{figure}[H]
    \begin{tikzpicture}[line cap=round,line join=round,>=triangle 45,x=1cm,y=0.75cm]
\clip(-3,-5) rectangle (3,3);
\draw [line width=1pt] (-2,2.5)-- (2,2.5);
\draw [line width=1pt] (-1,3)-- (-1.9909535324951024,1.592961574819557);
\draw [line width=1pt] (-2,2)-- (-1,0);
\draw [line width=1pt,dashed] (-1,2)-- (-2,0.6);
\draw [line width=1pt] (-0.9869053181920332,1.2429920647127655)-- (-2.002626176976524,-0.25385972718017363);
\draw [line width=1pt] (-1.9949891780382947,0.18144921229889535)-- (-1.0151237458887232,-1.7521886543562049);
\draw [line width=1pt] (-0.9997696351516704,-2.4415756974656944)-- (-1.9909535324951024,-4.423943492152565);
\draw [line width=1pt,dashed] (-2,-4)-- (1,-4);
\draw [line width=1pt,red] (0,-3.8044936953196897)-- (0,-4.3881345065590684);
\draw [line width=1pt,red] (0.6018638968205038,-3.815717557074293)-- (0.6018638968205038,-4.3881345065590684);
\begin{scriptsize}
\draw [fill] (2,2.75) node {$p_g+1$};
\draw [fill] (-2.25,1.5) node {$E_0$};
\draw [fill] (-0.75,2) node {$E_1$};
\draw [fill] (-0.75,0) node {$E_2$};
\draw [fill] (-2.5,-1.5) node {$ $};
\draw [fill] (-2,-3) node {$E_{\chi -4}$};
\draw [fill] (1.6,-4) node {$E_{\chi -3}$};
\draw [fill] (-1,-2) node {$\vdots$};
\end{scriptsize}
\end{tikzpicture}
\end{figure}
Define $D:=\psi_{*}E_1$ and $E:=\psi_{*}E_{\chi-3}$.
By construction, $D$ and $E$ intersect at the singularity of $W$.
Similar to the proof of Proposition~\ref{prop--construction--lee--park--i},
we see that $D$ and $E$ generate $\mathrm{Eff}(W)$ and satisfy $D^2<0$ and $E^2<0$.
Pursuing the even resolution process, we obtain
\begin{align*} 
-K_{\widehat{W}}&=\psi^{*}(2\Delta_0+(\chi+2)\Gamma)-\sum_{i=1}^{\chi-3}iE_i. 
\end{align*}
Hence, we conclude that
$$
h^{0}(-K_W)\ge h^{0}(-K_{\widehat{W}})\ge h^{0}(2\Delta_0+(\chi+1)\Gamma)\ge 9,
$$
and that both $-K_{\widehat{W}}$ and $-K_W=\delta_{*}(-K_{\widehat{W}})$ are big.
Similar to the the proof of Proposition~\ref{prop--construction--lee--park--i}, we have $-K_{W}\cdot E>0$.
However, since $E_1^{2}=-2$, we compute
$$
-K_{W}\cdot D=-\delta^{*}K_{W}\cdot E_1=-K_{\widehat{W}}\cdot E_{1}+(A_W(E_2)-1) E_2\cdot E_{1} =A_W(E_2)-1\le 0.
$$
Thus, we conclude that $-K_{W}$ is not ample.

Let $\gamma\colon W\to W'$ be the Grauert contraction of $D$ to a compact complex surface $W'$ (cf.~\cite[III, Theorem 2.1]{BHPV}),
and define $L':=\gamma_{*}L$, $B':=\gamma_{*}B$ and $E':=\gamma_{*}E$.
Since the Stein factorization of $X\to W\to W'$ induces a double cover $X'\to W'$ branched along $B'\in |2L'|$ which does not branch at the singularity of $W'$, the divisor $L'$ is Cartier.
Moreover, since $L\cdot D=0$, we have $L=\gamma^{*}L'$.
Since $L'\cdot E'=L\cdot E=1$ and $\rho(W')=1$, it follows from the Nakai-Moishezon criterion that $L'$ is ample (in particular, $W'$ is projective).
Since $-K_{W'}\cdot E'=-K_{W}\cdot E>0$, we obtain
$$
(2L'-K_{W'})\cdot E'>2,\quad (2L'-K_{W'})^{2}>4L'^{2}\ge 4.
$$
Thus, by applying the Reider-type theorem (cf.\ \cite[Corollary 5.11]{Enokizono}), we conclude that the linear system $|2L'|$ is basepoint-free.
Since $L=\gamma^{*}L'$, it follows that $|2L|$ is basepoint-free, which induces the contraction $\gamma$.

Furthermore, by applying \cite[Theorem 1.3]{Enokizono} to $2L'-K_{W'}$, we deduce that $$H^1(\gamma_{*}\O_{W}(2L))=H^{1}(2L')=0.$$
It is not hard to see that $W'$ has only rational singularities. Since $B$ and $D$ are disjoint, we obtain $R^{1}\gamma_{*}\O_{W}(2L)=0$.
Therefore, by the Leray spectral sequence, we conclude that $H^{1}(2L)=0$.
\end{proof}

By Proposition~\ref{prop--construction--lee--park--i*}, we see that Horikawa surfaces of general Lee-Park type $\mathrm{I}^*$ have different properties from surfaces of general Lee-Park type $\mathrm{I}$.
We note that it is meaningful to divide $\Q$-Gorenstein smoothable Horikawa surfaces with only two Wahl singularities into two classes because their behaviors are quite different, as we will show in Theorem \ref{thm--stratification} below.
    For example, $\Q$-Gorenstein smoothable Horikawa surfaces of special Lee-Park type are never partially $\Q$-Gorenstein smoothed to Horikawa surfaces of general Lee-Park type $\mathrm{I}^*$ when the geometric genus is more than five.
    
By Proposition~\ref{prop--construction--lee--park--i*}, we see that $(W,\frac{1}{2}B)$ satisfies the assumption of Proposition \ref{prop--smoothing--quotient} as the proof of Proposition \ref{prop--smoothing--lee--park--i}.

\begin{prop}\label{prop--Lee-Park-from-I*-to-I}
    Let $X$ be a Horikawa surface of general Lee-Park type $\mathrm{I}^*$.
    Then, $X$ is partially $\Q$-Gorenstein smoothable to Horikawa surfaces of general Lee-Park type $\mathrm{I}$.
    In particular, $X$ is $\Q$-Gorenstein smoothable.
\end{prop}

\begin{proof}

First, we construct a $\mathbb{Q}$-Gorenstein family $\mathscr{W}\to C$ over a germ of a smooth curve $0\in C$ such that the central fiber $\mathscr{W}_0$ is the quotient of $X$ by the good involution and general fibers $\mathscr{c}$ are the quotients of Horikawa surfaces of general Lee-Park type $\mathrm{I}$ by the good involutions.
    By deforming the centers of the blow-ups over $\Sigma_{\chi}$ described in the proofs of Propositions~\ref{prop--construction--lee--park--i} and \ref{prop--construction--lee--park--i*}, we can construct a blow-up sequence 
    $$
    \widehat{\mathscr{W}}=\mathscr{W}_{\chi-3}\to \cdots \to \mathscr{W}_1\to \mathscr{W}_0=\Sigma_{\chi}\times C,
    $$
    where the central fiber (resp.\ general fiber) of the map over $C$ coincides with the blow-up sequence described in the proof of Proposition~\ref{prop--construction--lee--park--i*} (resp.\ Proposition~\ref{prop--construction--lee--park--i}).
    Let $\mathscr{E}_i$ be the proper transform of the exceptional divisor of $\mathscr{W}_i\to \mathscr{W}_{i-1}$ on $\widehat{\mathscr{W}}$.
    Then, we note that $\mathscr{E}_{1}|_{\widehat{\mathscr{W}}_{0}}=\sum_{i=1}^{\chi-3}E_{\chi-3}$ and $\mathscr{E}_{i}|_{\widehat{\mathscr{W}}_{c}}=E_{i}$ for $i\neq 1$ or $c\neq 0$.
    Now, we define the divisor on $\widehat{\mathscr{W}}$ as
    $$
    \mathscr{D}:=(\chi-2)^{2}\Gamma_C+(\chi-3)\Delta_{C,0}+(\chi-4)\Gamma_{C,1}+\sum_{i=2}^{\chi-4}(\chi-3-i)\mathscr{E}_i,
    $$
    where $\Gamma_C$ (resp.\ $\Gamma_{C,1}$) denotes a general fiber (the proper transform of the fixed fiber) of $\widehat{W}\to \Sigma_{\chi}\times C\to \mathbb{P}^1\times C\to \mathbb{P}^1$ and $\Delta_
{C,0}$ denotes the proper transform of $\Delta_0\times C\subset \Sigma_0\times C$. 
    By construction, for each $c\in C$, there exists a contraction $\delta_c\colon \widehat{\mathscr{W}}_c\to W_c$ to the quotient of a Horikawa surface of general Lee-Park type by the good involution,
    which satisfies 
    $$
    (\chi-2)^{2}\delta_c^{*}\delta_{c*}(\Gamma_C|_{\widehat{\mathscr{W}}_c})=\mathscr{D}|_{\widehat{\mathscr{W}}_c}=:\mathscr{D}_c.
    $$ 
    Since $\delta_{c*}(\Gamma_C|_{\widehat{\mathscr{W}}_c})$ is ample  and $W_c$ has only rational singularities, we see that $H^1(l\mathscr{D}_c)=0$ for any sufficiently large and divisible $l\in\mathbb{Z}_{>0}$.
    This implies that $\mathscr{D}$ is relatively semiample over $C$ and 
    $$
\mathbf{Proj}_C(\bigoplus_{l\ge0}\widehat{f}_*\mathcal{O}_{\widehat{\mathscr{W}}}(l\mathscr{D}))_c\cong \mathrm{Proj}(\bigoplus_{l\ge0}H^{0}(\mathcal{O}_{\widehat{\mathscr{W}}_{c}}(l\mathscr{D}_{c})))\cong W_c,
    $$
    where $\widehat{f}\colon \widehat{\mathscr{W}}\to C$ is the natural morphism.
    Thus, 
    $$
\mathscr{W}:=\mathbf{Proj}_C(\bigoplus_{l\ge0}\widehat{f}_*\mathcal{O}_{\widehat{\mathscr{W}}}(l\mathscr{D}))
$$
with the canonical morphism $f\colon \mathscr{W}\to C$
is the desired family over $C$.
    By construction, $f$ is locally trivial around any singularity of $W$.
    In particular, $\mathscr{W}$ is a $\mathbb{Q}$-Gorenstein family over $C$.
    By Proposition \ref{prop--smoothing--quotient} and Proposition~\ref{prop--construction--lee--park--i*}, we may assume that there exist a $\mathbb{Q}$-Cartier divisorial sheaf $\mathscr{L}$ on $\mathscr{W}$ and an effective relative Mumford divisor $\mathscr{B}$ on $\mathscr{W}$ such that $\mathscr{L}_0=L$, $\mathscr{B}_0=B$ and $\mathscr{B}\sim_C\mathscr{L}^{[2]}$.
    Taking the double cover $\mathscr{X}\to \mathscr{W}$ branched along $\mathscr{B}$,
    we obtain a partial $\mathbb{Q}$-Gorenstein smoothing of $X$ to Horikawa surfaces of general Lee-Park type $\mathrm{I}$.
The last assertion follows from Proposition \ref{prop--smoothing--lee--park--i}.
\end{proof}


On the other hand, we remark the $\Q$-Gorenstein smoothability of Horikawa surfaces of general Lee-Park type $(\infty)$ first obtained by \cite{MNU} using results by \cite{DVS}.

\begin{prop}[{\cite[Example 4.4]{MNU}}]\label{prop--smoothing--urzua-type}
    Let $X$ be a Horikawa surface of infinite Lee-Park type.
    Then, $X$ is $\mathbb{Q}$-Gorenstein smoothable to Horikawa surfaces of type $(\infty)$.
\end{prop}

\begin{proof}
    This follows from the same argument as that of Proposition \ref{prop--smoothing-p_g=3}.
\end{proof}

\begin{rem}
    Lee and Park \cite{LP} first find Horikawa surfaces of general Lee-Park type $\mathrm{I}$ and show $\mathbb{Q}$-Gorenstein smoothability of them.
    Monreal, Negrete and Urz\'ua (cf.~\cite[Theorem 1.5]{MNU}) show that if $X$ is a stable non-Gorenstein klt $\mathbb{Q}$-Gorenstein smoothable Horikawa surface with $p_g(X)\ge10$, then $X$ should be of general Lee-Park type $\mathrm{I}$ or $\mathrm{I}^*$. 
\end{rem}

\subsection{Standard Horikawa surfaces}\label{subsec:standard-smoothable}


In this subsection, we study the smoothability of standard Horikawa surfaces as stated in Theorem \ref{intro-qGsmHor}~(1).
Let $X$ be a standard Horikawa surface,
and we adopt the notation from Theorem~\ref{thm:standard_Horikawa} without further explanation.
If $W=\varphi_{K_X}(X)$ is smooth, then we can construct a smoothing from $X$ to Horikawa surfaces with only Du Val singularities by partially smoothing the branch divisor $B$ in the associated linear system on $W$.
Note also that stable Horikawa surfaces with only Du Val singularities are smoothable by \cite{horikawa}.
Therefore, it suffices to consider the case where $W$ is singular.
In this case, $X$ is of type $(d)'$, and we have $W\cong \overline{\Sigma}_d$, where $d=p_g(X)-2$.
We begin by considering the case of small geometric genus.
It is straightforward to see the following:
\begin{lem}\label{lem--standard-for-low-p_g-smoothable}
Let $X$ be a standard Horikawa surface of type $(p_g-2)'$ with $4\le p_g\le 6$.
Suppose that $m(X)>\lceil\frac{2p_g-8}{p_g-2}\rceil$.
Then, $X$ is partially $\mathbb{Q}$-Gorenstein smoothable to standard Horikawa surfaces with multiplicity $m(X)-1$.
In particular, $X$ is smoothable.
\end{lem}
Note that the last assertion follows from the fact that stable Horikawa surfaces with only Du Val singularities are smoothable \cite{horikawa}.

Next, we consider the case where $p_g(X)\ge7$. 
By Proposition~\ref{prop:classification_cone_sing}, $(W, \frac{1}{2}B)$ has a cone singularity of type $(d;k_1,k_2)$ for some $k_1$, $k_2$.
We need to distinguish the following types.

\begin{defn}[Double Fano type, double non-Fano type, supersingular]\label{defn--of--fano}
Let $X$ be a standard Horikawa surface of type $(p_g-2)'$ with $p_g\ge 7$.
We say that $X$ is {\it supersingular} if $B^{+}-2\Delta_0$ contains a fiber $\Gamma$ of the ruling $W^{+}=\Sigma_{p_g-2}\to \mathbb{P}^1$ such that the residual part $B^{+}-2\Delta_0-\Gamma$ intersects $\Gamma$ at the point $\Gamma\cap \Delta_{0}$ with multiplicity greater than $2$.
We say that $X$ is of {\it double non-Fano type} if either
\begin{itemize}
   \item 
$B^{+}-2\Delta_0$ contains a fiber $\Gamma$ such that the residual part $B^{+}-2\Delta_0-\Gamma$ intersects $\Gamma$ at $F\cap \Delta_0$ with multiplicity $2$, or
   \item  
there exists a fiber $\Gamma$ not contained in $B^{+}$ such that $B^{+}-2\Delta_0$ intersects $\Gamma$ at $\Gamma\cap \Delta_0$ with multiplicity greater than $\mathrm{mult}_{\Gamma\cap \Delta_0}(B^{+}-2\Delta_0)$.
\end{itemize}
We say that $X$ is {\it of double Fano type} if it is neither supersingular nor of double non-Fano type.
\end{defn}

\begin{rem} \label{rem--houshin}
Our approach to detecting $\mathbb{Q}$-Gorenstein smoothability of standard Horikawa surfaces is as follows:
First, we take an admissible anti $P$-resolution 
$$
\mu\colon \left(W^-,\frac{1}{2}B^-\right)\to \left(W,\frac{1}{2}B\right).
$$
Second, we verify that there is no local-to-global obstruction to $\Q$-Gorenstein deformation of $W^-$, and that $H^1(B^-)=0$.
This method works well when $X$ is not supersingular. However, for supersingular Horikawa surfaces of geometric genus $10$, we have $H^1(B^-)\ne0$ .
\end{rem}
In this subsection, we study standard Horikawa surfaces of double Fano type and of double non-Fano type before turning to supersingular case. 
We begin by introducing the following definition.

\begin{defn}
Let $X$ be a normal standard Horikawa surface of type $(p_g-2)'$ with $p_g\ge 7$.
Let $\mu\colon (W^-,\frac{1}{2}B^-)\to (W,\frac{1}{2}B)$ be an admissible anti $P$-resolution of length $(k_1,k_2)$.
We say that $\mu$ is {\em supersingular} if
there exists a fiber $\Gamma$ of the ruling $W^{+}=\Sigma_{p_g-2}\to \mathbb{P}^1$ such that $\Gamma\subset B^{+}$, the residual part $B^{+}-2\Delta_{0}-\Gamma$ intersects $\Gamma$ at $\Gamma\cap \Delta_0$ with multiplicity greater than $2$,  and the corresponding length $k_i$ is greater than $4$.
In particular, $X$ is supersingular.
\end{defn}

We note
the following:

\begin{lem}\label{lem--existence--smoothable--cusp}
If we fix intergers $d$ and $k_1$ with $d\ge 5$ and $k_1\le 2d-8$, then there exist both standard Horikawa surfaces of double non-Fano type and of double Fano type with a singularity of type $(d;k_1,0)$.
\end{lem}

\begin{proof}
   We will show that there exists an effective Weil divisor $B\in |(4d+4)\overline{\Gamma}|$ on $\overline{\Sigma}_d$ such that the pair $(\overline{\Sigma}_d,\frac{1}{2}B)$ defines a singularity of type $(d;2d-8,0)$ and the associated double cover is a standard Horikawa surface of type $(d)'$.
    To construct this, we consider a projective birational morphism $\psi\colon \widetilde{W}\to \Sigma_d$ between smooth surfaces, which can be decomposed as a sequence of blow-ups $\psi=\psi_1\circ\ldots\circ \psi_{d-4}$.
    The sequence is defined as follows:
    \begin{itemize}
        \item $\psi_1$ is the blow-up at the intersection point of $\Delta_0$ and a fiber $\Gamma_0$ of the canonical ruling $p\colon \Sigma_d\to \mathbb{P}^1$.
        \item $\psi_2$ is the blow-up at a point on the exceptional divisor of $\psi_{1}$
        that does not lie on the proper transform of $\Delta_0$.
        \item $\psi_i$ for $i=3,\ldots,d-4$ is the blow-up at a point on the exceptional divisor of $\psi_{i-1}$
        that does not lie on the proper transform of the exceptional divisor of $\psi_{i-2}$.
    \end{itemize}
    Let $E_i$ denote the proper transform of the $\psi_i$-exceptional divisor, and set $E_0$ as the proper transform of $\Gamma_0$.
    Note that $E_{0}^{2}=-1$ or $E_{0}^{2}=-2$ since the proper transform $B'$ on $\Sigma_d$ of $B$ satisfies $B'\cdot \Gamma=4$.
    Define 
    $$
    \widetilde{L}:=\pi^*(2\Delta_0+(2d+2)\Gamma)-\sum_{i=1}^{d-4}iE_i.
    $$ 
    Now, we claim that $|2\widetilde{L}|$ is basepoint-free.
    By construction, it is straightforward to see that both $\widetilde{L}$ and $-K_{\widetilde{W}}$ are big, and that $h^0(-K_{\widetilde{W}})\ge 9$.
    Since 
    $$
    D:=2\widetilde{L}-K_{\widetilde{W}}=\pi^{*}(6\Delta_0+(5d+6)\Gamma)-3\sum_{i=1}^{d-4}iE_i,
    $$
    we have
    $$
    D\cdot \psi^{-1}_*\Delta_0=3-d,\quad D\cdot E_0=6+3E_{0}^2, \quad D\cdot E_{i}=0 \ (i=1,\ldots, d-5), \quad D\cdot E_{d-4}=3.
    $$
    Thus, by using Lemma \ref{lem--modified--manetti--lemma}, we can check that
    $D$ is $\Z$-positive in the sense of \cite[Definition 3.2]{Enokizono}.
    Moreover, $D-\psi^{-1}_*\Delta_0$ and $D-E_i$ for $0\le i\le d-4$ are also big and $\Z$-positive.
    Applying \cite[Theorem 1.3]{Enokizono}, we conclude that 
    \[
   H^{1}(2\widetilde{L})=H^1(2\widetilde{L}-\psi^{-1}_*\Delta_0)=H^1(2\widetilde{L}-E_i)=0
    \]
    for any $0\le i\le d-4$. 
    On the other hand, it is straightforward to see that the base locus of $|2\widetilde{L}|$ is contained in $\psi^{-1}_*\Delta_0\cup\bigcup_{i=0}^{d-4}E_i$.
    Hence, by using the vanishing of the above cohomology groups, we see that $|2\widetilde{L}|$ is basepoint-free.
    Thus, taking a general member $\widetilde{B}\in |2\widetilde{L}|$ and defining $B:=\pi_*\widetilde{B}$, we obtain the desired result.
    Note that the case where $E_0^2=-1$ (resp.\ $E_0^2=-2$) corresponds to the double Fano type (resp.\ double non-Fano type).
\end{proof}

Surfaces as in Lemma \ref{lem--existence--smoothable--cusp} are indeed smoothable as we will show in Proposition \ref{prop--smoothing--c-type} below but there are non-smoothable standard stable Horikawa surfaces as we explain below.
This means that the smoothability for standard Horikawa surfaces is more complicated than for surfaces of general Lee-Park type.

\begin{rem}\label{rem--necessity-of-P-or-anti-p}
    It is straightforward to construct a standard stable normal Horikawa surface with a double cone singularity that is simple elliptic for each geometric genus $p_g\ge 4$.
    Let $E$ be the exceptional divisor of the minimal resolution of the singularity.
    A direct computation shows that $E^2=2p_g-4$.
    By Lemma \ref{lem:smoothable_elliptic_singularity}, if $p_g \ge 7$, the singularity is not smoothable.
    Consequently, there exists a standard stable normal Horikawa surface that is not smoothable.
    Furthermore, the ample condition of $B-K_W$ does not necessarily imply the existence of log $\mathbb{Q}$-Gorenstein deformation. 
    This means that we cannot apply Proposition \ref{prop--smoothing--quotient} directly without the assumption on T-singularities.
\end{rem}

The following is a fundamental lemma for investigating several properties of standard Horikawa surfaces and their non-supersingular anti-P-resolutions.
\begin{lem}\label{lem--vanishing--for--standard}
    Let $B\in |(4d+4)\overline{\Gamma}|$ be an effective Weil divisor on $\overline{\Sigma}_{d}$ for some $d\ge 5$ such that $(\overline{\Sigma}_{d}, \frac{1}{2}B)$ defines a cone singularity of type $(d;k'_1,k'_2)$ at the cone point.
    Suppose that there exists an admissible anti-P-resolution $\mu\colon (W^{-},\frac{1}{2}B^{-})\to (\overline{\Sigma}_{d},\frac{1}{2}B)$ of length $(k_1,k_2)$ such that $k_1+k_2\le 2d-2$ for some $0\le k_i\le k'_i$ for $i=1,2$. 
   Then, the following hold:
   \begin{itemize}
       \item[$(1)$]
       Let $W^{-}_{\min}$ be the minimal resolution of the normalization of $W^{-}$.
       Then, $-K_{W^{-}_{\min}}$ is big and
    \begin{equation}\label{eq--antigenus-of-anti-P-for-standard}
  \min\{h^0(-K_{W^{-}_{\min}}),h^0(-K_{W_{\mathrm{good}}})\}\ge 6.
    \end{equation}
    In particular, $W_{\mathrm{good}}$ and $W^{-}_{\min}$ satisfy the assumption of Lemma \ref{lem--modified--manetti--lemma}.

    \item[$(2)$]
    Let $B'_{\mathrm{good}}$ be the proper transform of $B$ on the good resolution $W_{\mathrm{good}}$ defined in Definition~\ref{defn--good-resolution}.
    Suppose that $\mu$ is not supersingular.
    Then, $\mathcal{O}_{W_{\mathrm{good}}}(B'_{\mathrm{good}})$ is big and globally generated, and 
    $$  H^1(W_{\mathrm{good}},\mathcal{O}_{W_{\mathrm{good}}}(B'_{\mathrm{good}}))=0.
   $$  
   \end{itemize}

    \end{lem}

\begin{proof}
(1): We use the notation as in Definition \ref{defn--good-resolution} freely.
Let $D_i$ denote the reduced total transform of $\Delta_{0}$ on $W_{i}$.
If $\psi_{i}$ is the blow-up at a node in $D_{i-1}$, then we have $K_{W_{i}}+D_{i}=\psi_{i}^{*}(K_{W_{i-1}}+D_{i-1})$, and hence 
$$
h^{0}(-(K_{W_{i}}+D_{i}))=h^{0}(-(K_{W_{i-1}}+D_{i-1})).
$$
If $\psi_{i}$ is the blow-up at a smooth point of $D_{i-1}$, then we have $K_{W_{i}}+D_{i}=\psi_{i}^{*}(K_{W_{i-1}}+D_{i-1})+E_i$, where $E_i$ is the exceptional divisor of $\psi_{i}$.
In particular, we have 
$$
h^{0}(-(K_{W_{i}}+D_{i}))\ge h^{0}(-(K_{W_{i-1}}+D_{i-1}))-1.
$$
Hence, by the construction of the good resolution, it suffices to show that $-(K_{W_{N}}+D_{N})$ is big and $h^{0}(-(K_{W_{N}}+D_{N}))\ge 6$.
Note that 
$$
h^{0}(-(K_{W_{0}}+D_{0}))=h^{0}(\Delta_{0}+(d+2)\Gamma)=d+6.
$$
By the definition of the good resolution, the number of blow-ups at a smooth point of $D_i$ for some $i$ equals $\lceil \frac{k_1}{2} \rceil +\lceil \frac{k_2}{2} \rceil$.
Thus, we conclude that
$$
h^{0}(-(K_{W_{N}}+D_{N}))\ge h^{0}(-(K_{W_{0}}+D_{0}))-\left( \left\lceil \frac{k_1}{2} \right\rceil +\left\lceil \frac{k_2}{2} \right\rceil \right)\ge  6.
$$
The statement of bigness follows from the similar argument as above applying to $m(K_{W_{i}}+D_{i})$ for $m\gg 0$.

(2): Using (1), it is straightforward to show the assertion by a similar argument as in the proof of Lemma \ref{lem--existence--smoothable--cusp}.
\end{proof}
We remark about the above proof.
\begin{rem}\label{rem--supersingular--anti--p-resol}
In the proof of Lemma \ref{lem--vanishing--for--standard}, we can use the similar argument to Lemma \ref{lem--existence--smoothable--cusp} because  $E_0^2=-1$ or $E_0^2=-2$ when $\mu$ is not supersingular, where $E_0$ is defined in the proof of Lemma \ref{lem--existence--smoothable--cusp}.
This condition is the key to show the ampleness of $B^-$.
If $\mu$ is supersingular, then $E_0^2=-3$ or $E_0^2=-4$.
In this case, $B^-$ is not ample and $H^1(B^-)=0$ is not true in general.
This fact can be easily verified by the argument of the proof of Proposition \ref{prop--super--singular--typeC} below.
\end{rem}

 We put the following two corollaries of Lemma \ref{lem--vanishing--for--standard}.

\begin{cor}\label{cor--properties--anti-P-standard}
    Let $X$ be a standard Horikawa surface of type $(p_g-2)'$ with $p_g\ge 7$.
Let $\mu\colon(W^-,\frac{1}{2}B^-)\to (W,\frac{1}{2}B)$ be an admissible anti-P-resolution such that the normalization $\widetilde{W}$ of $W^-$ satisfies $h^0(\widetilde{W},\mathcal{O}_{\widetilde{W}}(-K_{\widetilde{W}}))\ge6$. 
 Suppose that $\mu$ is not supersingular.
Then, the following hold:
\begin{itemize}
    \item[$(1)$] $B^-$ is ample.
   
    \item[$(2)$] If $X$ is of double Fano type, then $-K_{W^-}$ is log big and nef with respect to $W$. 
    If $W^-$ is further klt, then $-K_{W^-}$ is ample.
\end{itemize}

\end{cor}

\begin{proof}
(1): Note that the morphism $\mu$ belongs to one of the three types of anti-P-resolutions classified in Theorem \ref{thm--anti-P-k_1-k_2}. 
Let $\xi\colon W^-_{\mathrm{min}}\to W^-$ be the minimal resolution.
Then, there exists a birational morphism $\beta\colon  W^-_{\mathrm{min}}\to \Sigma_{p_g-2}$ such that the exceptional locus $\mathrm{Exc}(\beta)$ is contained in either one or two fibers of the composition $p\circ\beta$, where $p\colon \Sigma_{p_g-2}\to\mathbb{P}^1$ is the canonical projection.
Let $F$ be one of the these fibers. 
    Let $E$ be the $\mu$-exceptional locus on $W^-$, and let $D$ be the unique irreducible divisor on $W^-$ that is not $\mu$-exceptional but whose proper transform $\xi^{-1}_*D$ is contained in $F$.
    Since $D$ is not contained in the support of $B^-$ but intersects it, we have $D\cdot B^- >0$.
    On the other hand, by definition of $\mu$, we have $(K_{W^-}+\frac{1}{2}B^-)\cdot E=0$ and $K_{W^-}\cdot E<0$, which implies $B^-\cdot E>0$.
    Furthermore, according to \cite[(3.9.2)]{Wahl}, we have $h^0(W^-_{\min},\mathcal{O}_{W^-_{\min}}(-K_{W^-_{\min}}))\ge6$.
    Therefore, by Lemmas \ref{lem--modified--manetti--lemma} and \ref{lem--nakai--moishezon}, $B^{-}$ is ample, and the first assertion follows.

(2): The last two assertions follow from a similar argument to that of the first assertion:
If $X$ is of double Fano type, since $-K_{W^-}$ is $\mu$-ample, it suffices to show that $K_{W^-}\cdot D\le 0$.
Here, $D$ is the proper transform of a $(-1)$-curve $\hat{D}$ on $W^-_{\mathrm{min}}$, which intersects transversally and exactly once  with the exceptional divisor of $W^-_{\mathrm{min}}\to  W^-$.
Thus, $K_{W^-}\cdot D\le 0$. 
Moreover, if $W^-$ is klt, then $K_{W^-}\cdot D< 0$.
Therefore, in the double Fano case, 
the last assertion follows from Lemmas \ref{lem--modified--manetti--lemma} and \ref{lem--nakai--moishezon}.
This completes the proof.
\end{proof}

\begin{rem}
    Corollary \ref{cor--properties--anti-P-standard} (2) does not hold in general for the double non-Fano case and (1) does not hold in general for the supersingular case. 
    If $X$ is of double non-Fano type, then by the proof of Theorem \ref{thm--anti-P-k_1-k_2},  there are cases where there exists a $(-2)$-curve $C$ on $W^-_{\min}$ such that $\xi_*C$ passes through a singularity of $W^-$.
In such a case, $\xi_*C\cdot K_{W^-}>0$ and hence $-K_{W^-}$ is not nef.
\end{rem}

\begin{cor}\label{cor--partial--smoothing--from--non-Fano-to-Fano}
Fix an arbitrary pair $(k_1,k_2)\in\mathbb{Z}_{\ge0}^{\times2}$ such that $k_1+k_2\le 2d-2$ for $d\ge 5$. 
Then, every standard Horikawa surface $X$  of double non-Fano type with a double cone singularity of type $(d; k_1,k_2)$ can be partially smoothed to standard Horikawa surfaces of double Fano type with a double cone singularity of the same type.
\end{cor}

\begin{proof}
Let $B_X$ be the branch divisor of the quotient $X\to W=\overline{\Sigma}_{d}$.
Let $Y$ be a general standard Horikawa surface of double Fano type with a double cone singularity of the same type $(d;k_1,k_2)$, and let $B_Y$ be the branch divisor on $W$.
Let $\pi_X\colon W_X\to W$ (resp.\ $\pi_Y\colon W_Y\to W$) be the composition of the good resolution of $(W, \frac{1}{2}B_{X})$ (resp.\ $(W, \frac{1}{2}B_{Y})$).
By deforming the centers of the blow-ups over $\Sigma_{d}$ described in the proof of Lemma~\ref{lem--existence--smoothable--cusp},
we can construct a blow-up sequence
$$
\mathscr{W}_{\mathrm{good}}=\mathscr{W}_{N}\to \cdots \to \mathscr{W}_{1}\to \mathscr{W}_{0}:=\Sigma_{d}\times C
$$
which are flat over a smooth affine $C$ with two closed points $0,c\in C$ such that $\mathscr{W}_{\mathrm{good},0}=W_{X}$ and $\mathscr{W}_{\mathrm{good},c}=W_{Y}$.
Let $B'_{X}$ (resp.\ $B'_{Y}$) be the proper transform of $B_{X}$ on $W_{X}$ (resp.\ $B_Y$ on $W_Y$).
It follows from Lemma~\ref{lem--vanishing--for--standard} that $H^{1}(\O_{W_X}(B'_{X}))=H^{1}(\O_{W_Y}(B'_{Y}))=0$.
Thus, there exists a relative effective Cartier divisor $\mathscr{B}'_{\mathrm{good}}$ on $\mathscr{W}_{\mathrm{good}}$ such that $\mathscr{B}'_{\mathrm{good,0}}=B'_{X}$ and $\mathscr{B}'_{\mathrm{good,c}}\cong B'_Y$.
Moreover, there exist a line bundle $\mathscr{L}_{\mathrm{good}}$ and a relative effective Cartier divisor $\mathscr{B}_{\mathrm{good}}$ on $\mathscr{W}_{\mathrm{good}}$ such that $\mathscr{B}_{\mathrm{good}}\sim_{C} 2\mathscr{L}_{\mathrm{good}}$ and $\mathscr{B}_{\mathrm{good}}-\mathscr{B}'_{\mathrm{good}}$ is the families of disjoint $(-2)$-curves $E$ appeared in Definition~\ref{defn--good-resolution}.
The divisor $K_{\mathscr{W}_{\mathrm{good}}}+\mathscr{L}_{\mathrm{good}}$ is semiample over $C$ and the corresponding contraction $\pi\colon \mathscr{W}_{\mathrm{good}}\to \mathscr{W}$ coincides with the composition of $\mathscr{W}_{\mathrm{good}}\to \Sigma_{d}\times C$
and $\Sigma_{d}\times C\to \overline{\Sigma}_{d}\times C=:\mathscr{W}$.
The double cover $\mathscr{X}\to \mathscr{W}$ branched along $\pi_{*}\mathscr{B}_{\mathrm{good}}$
defines the desired partial smoothing from $X$ to $Y$.
We complete the proof.
\end{proof}

Now, we can show that there exists a $\Q$-Gorenstein smoothable standard Horikawa surface.

\begin{prop}\label{prop--smoothing--c-type}
   Let $X$ be a  standard Horikawa surface with a double cone singularity of type $(p_g-2; k_1,k_2)$ such that $k_1\ge 2p_g-12$ or $k_2\ge 2p_g-12$.
   Suppose that $(W,\frac{1}{2}B)$ admits a non-supersingular admissible anti-P-resolution of normal type.
    Then, $X$ is $\mathbb{Q}$-Gorenstein smoothable.

    In particular, there exists a $\Q$-Gorenstein smoothable standard Horikawa surface with a double cone singularity both of double Fano type and of double non-Fano type.
\end{prop}

\begin{proof}
By Lemma \ref{lem--vanishing--for--standard},
we may assume that either $k_1=2p_g-12$ or $k_2= 2p_g-12$.
Moreover, by Corollary \ref{cor--partial--smoothing--from--non-Fano-to-Fano}, we may further assume that $X$ is of double Fano type.  
Then, we can take an admissible anti-P-resolution 
$\mu\colon (W^{-}, \frac{1}{2}B^{-}) \to (W, \frac{1}{2}B)$ of normal type  and of length $(2p_g-12, 0)$.
Let $E_{1}$ denote the proper transform on $W^{-}$ of the fiber on which  $\Sigma_{p_g-2}$ fails to be isomorphic to $W^{-}$.
Let $E_2$ be the exceptional curve of $\mu$.
By Lemmas \ref{lem--vanishing--for--standard} and \ref{lem--modified--manetti--lemma}, 
any irreducible curve $C$ on $W^{-}$ other than $E_1$ and $E_2$ has non-negative self-intersection number.
Moreover, it is straightforward to verify that $K_{W^{-}}\cdot E_{i} <0$ for $i=1,2$.
Since $-K_{W^{-}}$ is big, Lemma \ref{lem--nakai--moishezon} implies that it is in fact ample.
Let $L^{-}$ be the Weil divisor on $W^{-}$ such that $2L^{-}\sim B^{-}$.
Then $K_{W^{-}}+L^{-}$ is a big and nef Cartier divisor since it is the pullback of an ample Cartier divisor $K_{W}+\mu_*L^-$ on $W$.
Hence, all the assumptions of Proposition \ref{prop--smoothing--slc} are satisfied for the pair $(W^{-},B^{-})$.
Therefore, there exists a projective flat morphism 
$$
f\colon \mathscr{W}^{-} \to C
$$
from a normal threefold to a smooth affine curve with a closed point $0\in C$, along with an effective $\mathbb{Q}$-Cartier relative Mumford divisor $\mathscr{B}^{-}$ on $\mathscr{W}^{-}$, such that $K_{\mathscr{W}^{-}}$ is $\mathbb{Q}$-Cartier, $(\mathscr{W}^{-}_0,\mathscr{B}^{-}_0)\cong (W^{-},B^{-})$ and $\mathscr{W}^{-}_c$ is smooth for any $c\in C\setminus\{0\}$.
We compute that $K_{W^{-}}^2=8$, and since $-K_{W^{-}}$ is ample, the general fiber $\mathscr{W}^{-}_c$ is isomorphic to either $\Sigma_0$ or $\Sigma_1$.
Furthermore, since $H^1(W^{-},\mathcal{O}_{W^{-}}(B^{-}))=0$ and $\O_{\mathscr{W}^{-}_c}(\mathscr{B}^{-}_c)$ is globally generated, we may assume that $\mathscr{B}^{-}_c$ is smooth.
By Proposition \ref{prop--smoothing--quotient} and Corollary \ref{cor--gen--fiber--rho}, it follows that $X$ admits a $\mathbb{Q}$-Gorenstein smoothing.
Therefore, we obtain the first assertion.
The last assertion follows from Lemma \ref{lem--existence--smoothable--cusp}.
\end{proof}

In Corollary \ref{cor--partial--smoothing--from--non-Fano-to-Fano}, we show that we can partially smooth standard Horikawa surfaces of double non-Fano type to surfaces of double Fano type. 
However, it is still unclear whether the smoothability would be preserved or not.
To handle this problem, we also discuss partial $\Q$-Gorenstein smoothability of anti-P-resolutions as follows.

\begin{prop}\label{prop--partial--smoothing--from--non-Fano-to-Fano-with-anti-P} 
Let $X$ be a standard Horikawa surface with a double cone singularity of type $(d;k_1,k_2)$, where $k_1+k_2\le 2d-2$ and $d\ge5$.
Suppose that $X$ is not of double Fano type.
Let $\mu\colon (W^{-}, \frac{1}{2}B^{-})\to (W, \frac{1}{2}B)$ be any admissible anti-P-resolution such that $\mathrm{Ex}(\mu)\cong\mathbb{P}^1$ and $\mu$ is not supersingular.
Then, there exists a projective flat morphism $g\colon \mathscr{W}^{-} \to C$ from a $\mathbb{Q}$-Gorenstein slc scheme to a smooth curve with an effective $\mathbb{Q}$-Cartier relative Mumford divisor $\mathscr{B}^{-}$ on $\mathscr{W}^{-}$ and a closed point $0\in C$ such that the following hold:
\begin{itemize}
\item[$(1)$] $g$ is locally trivial as a deformation.
\item[$(2)$] $(\mathscr{W}^{-}_0,\mathscr{B}^{-}_0)\cong (W^{-},B^{-})$.
\item[$(3)$]
For any general $c\in C$, let $\mathscr{W}_c$ denote the log canonical model of $(\mathscr{W}^{-}_c,\frac{1}{2}\mathscr{B}^{-}_c)$.
Then, the canonical morphism $\eta\colon (\mathscr{W}^{-}_c,\frac{1}{2}\mathscr{B}^{-}_c)\to (\mathscr{W}_c,\frac{1}{2}\mathscr{B}_c)$ is an admissible P-resolution and the double cover $\mathscr{X}_{c}\to \mathscr{W}_{c}$ branched along $\mathscr{B}_c$ is a stable normal standard Horikawa surface of double Fano type, where $\mathscr{B}_c:=\eta_*\mathscr{B}^{-}_c$.
\item[$(4)$]
If the germ $(\mathrm{Ex}(\mu)\subset W^{-})$ is $\mathbb{Q}$-Gorenstein smoothable, then so is $(\mathrm{Ex}(\eta)\subset \mathscr{W}^{-}_c)$.
\end{itemize}

\end{prop}

Before giving the proof of Proposition \ref{prop--partial--smoothing--from--non-Fano-to-Fano-with-anti-P}, we deal with the following lemma, which is also used to discuss smoothability of supersingular standard Horikawa surfaces later.

\begin{lem}\label{lem--gluing--tool}
Let $f\colon (\mathscr{W},\frac{1}{2}\mathscr{B})\to C$ be a log $\mathbb{Q}$-Gorenstein deformation such that $\mathscr{W}_c\cong\overline{\Sigma}_d$ and $\mathscr{B}_c\in|(4d+4)\overline{\Gamma}|$ for any $c\in C$. Fix a closed point $0\in C$ and suppose that there exists an admissible anti-P-resolution $\mu\colon (W^-,\frac{1}{2}B^-)\to (W,\frac{1}{2}B)$ of non-normal type such that $\mathrm{Ex}(\mu)\cong\mathbb{P}^1$. Let $\widetilde{W}$ be the normalization of $W^-$ and $\widetilde{B}$ the proper transform of $B^-$.
Suppose there exists a projective flat morphism $h\colon\widetilde{\mathscr{W}} \to C$ with a $\mathbb{Q}$-Cartier effective relative Mumford divisor $\widetilde{\mathscr{B}}$ and a birational contraction $\xi\colon \widetilde{\mathscr{W}} \to \mathscr{W}$ such that:

\begin{itemize}
 \item   
$\widetilde{\mathscr{W}}$ is $\mathbb{Q}$-Gorenstein,
\item 
$\xi_*\widetilde{\mathscr{B}}=\mathscr{B}$,
\item 
$K_{\widetilde{\mathscr{W}}}+\frac{1}{2}\widetilde{\mathscr{B}}+\mathrm{Ex}(\xi)=\xi^*\left(K_{\mathscr{W}}+\frac{1}{2}\mathscr{B}\right)$,
\item 
$(\widetilde{\mathscr{W}}_0,\widetilde{\mathscr{B}}_0)\cong (\widetilde{W},\widetilde{B})$, and
\item
$h$ is trivial \'etale locally around the vertex of $\overline{\Sigma}_d$.
\end{itemize}
Then, by shrinking $C$, there exists a projective flat morphism $g\colon \mathscr{W}^{-} \to C$ from a $\mathbb{Q}$-Gorenstein slc scheme to a smooth curve with an effective $\mathbb{Q}$-Cartier relative Mumford divisor $\mathscr{B}^{-}$ on $\mathscr{W}^{-}$ and a closed point $0\in C$ such that the following hold:
\begin{itemize}
\item[$(1)$] there exists a morphism $\nu\colon\widetilde{\mathscr{W}}\to \mathscr{W}^{-}$ that is the normalization and $\widetilde{\mathscr{B}}=\nu^{-1}_*\mathscr{B}^-$,
\item[$(2)$] there exists a morphism $\eta\colon(\mathscr{W}^{-},\frac{1}{2}\mathscr{B}^-)\to (\mathscr{W},\frac{1}{2}\mathscr{B})$ that is log crepant and satisfies that $\eta_c$ is an anti-P-resolution,
\item[$(3)$] $(\mathscr{W}^{-}_0,\frac{1}{2}\mathscr{B}^-_0)\cong (W^-,\frac{1}{2}B^-)$ and $\eta_0=\mu$ via this identification, and
\item[$(4)$] $g$ is locally trivial.
\end{itemize}
Furthermore, if the germ $(\mathrm{Ex}(\mu)\subset W^-)$ is $\Q$-Gorenstein smoothable, then $(\mathrm{Ex}(\eta_c)\subset \mathscr{W}_c^-)$ is also $\Q$-Gorenstein smoothable for any $c\in C$.
\end{lem}

\begin{proof}
It is easy to see that there exists a divisorial sheaf $\widetilde{\mathscr{L}}$ such that $\widetilde{\mathscr{L}}^{[2]}\cong\mathcal{O}_{\widetilde{\mathscr{W}}}(\widetilde{\mathscr{B}})$.
By the proof of Theorem \ref{thm--anti-P-k_1-k_2}, there exists $r\in\mathbb{Z}_{>0}$ such that $\widetilde{W}$ has a singularity of the form  $(xy=0)\subset \frac{1}{r}(1,-1,a)$ and 
$\widetilde{\mathscr{L}}^{[r]}$ is Cartier and $\xi$-ample.
Take a relative very ample Cartier divisor $\mathscr{H}$ on $\overline{\mathscr{W}}$.
Then, using \cite[Theorem 1.10]{fujino--slc--vanishing}, we obtain 
\begin{equation*} \label{eq--vanishing-to-gluing} 
H^1(r\widetilde{\mathscr{L}}_0+u\xi_0^*\mathscr{H}_0)=H^1(r\widetilde{\mathscr{L}}_0 +u\xi_0^*\mathscr{H}_0-\mathrm{Ex}(\xi_0))=0
\end{equation*}
 for any sufficiently large integer $u>0$.
Thus, there exists an effective relative Cartier divisor $\mathscr{D}$ such that 
$$
\mathscr{D}\sim_{C} r\widetilde{\mathscr{L}}+u\xi^*\mathscr{H}
$$ 
and $\mathscr{D}_0$ intersects transversally with each irreducible component of $\mathrm{Ex}(\xi_0)$.
We let $\tau_0\colon \mathrm{Ex}(\xi_0)^\nu\to \mathrm{Ex}(\xi_0)^\nu$ be the canonical involution on the normalization of $\mathrm{Ex}(\xi_0)$.
Using the vanishing result above, we may choose $\mathscr{D}$ so that its restriction $(\mathscr{D}_0)|_{\mathrm{Ex}(\xi_0)^\nu}$ is preserved by $\tau_0$.
Note that the pair $(\widetilde{\mathscr{W}},\frac{1}{2}\widetilde{\mathscr{B}}+t\mathscr{D}+\mathrm{Ex}(\xi)+\widetilde{\mathscr{W}}_0)$ is log canonical around $\mathrm{Ex}(\xi)$ (cf.~\cite{kawakita}) and the divisor 
$$
K_{\widetilde{\mathscr{W}}}+\frac{1}{2}\widetilde{\mathscr{B}}+t\mathscr{D}+\mathrm{Ex}(\xi)+\widetilde{\mathscr{W}}_0
$$
is $\xi$-ample for any $t\in \mathbb{Q}\cap (0,1]$.
After taking a finite cover of $C$ if necessary, we may assume that $\mathrm{Ex}(\xi)$ is an slc surface consists of two $\mathbb{P}^1$-bundles over $C$, glued along a section.
Let $\mathrm{Ex}(\xi)^\nu$ be the normalization and denote by $S$ the conductor.
It is not hard to see that there exists an involution $\tau\colon\mathrm{Ex}(\xi)^\nu\to \mathrm{Ex}(\xi)^\nu$ over $C$, which preserves $S$ and exchanges the components of the pullback of $\mathscr{D}|_{\mathrm{Ex}(\xi)}$. 
Moreover, $\tau$ preserves the different $\mathrm{Diff}_{\mathrm{Ex}(\xi)^\nu}(\frac{1}{2}\widetilde{\mathscr{B}}+t\mathscr{D})$.
By Koll\'ar's gluing theorem \cite{kollar-mmp}, \cite[Theorem 11.38]{kollar-moduli}, we see that there exists a finite morphism $\nu \colon \widetilde{\mathscr{W}} \to \mathscr{W}^{-}$ to an slc scheme that is projective and flat over $C$ such that $\nu$ is the normalization and
\[
K_{\widetilde{\mathscr{W}}}+\frac{1}{2}\widetilde{\mathscr{B}}+t\mathscr{D}+\mathrm{Ex}(\xi)+\widetilde{\mathscr{W}}_0=\nu^*\left(K_{\mathscr{W}^{-}}+\frac{1}{2}\nu_*\widetilde{\mathscr{B}}+t\nu_*\mathscr{D}+\mathscr{W}^{-}_0\right).
\]
Again by \cite[Theorem 11.38]{kollar-moduli}, the central fiber $\mathscr{W}^{-}_0$ is isomorphic to $W^{-}$, and the morphism $g\colon \mathscr{W}^{-}\to C$ is a locally trivial deformation.
Let $\mathscr{B}^{-}:=\nu_*\widetilde{\mathscr{B}}$.
Then, it is straightforward to verify that the pair $(\mathscr{W}^{-},\mathscr{B}^{-})$ satisfies the required properties.
This completes the proof of the first assertion.

Since $g$ is locally trivial, the sheaf $\mathscr{T}^1_{\mathrm{QG}}(W^{-})$ deforms to $\mathscr{T}^1_{\mathrm{QG}}(\mathscr{W}^{-}_c)$ for general $c\in C$.
We may assume that the non-normal locus $\mathrm{Ex}(g_c)$ of $\mathscr{W}^{-}_c$ is isomorphic to $\mathbb{P}^1$ for any general $c\in C$.
Let $\mathscr{T}^1_{\mathrm{QG}}(\mathscr{W}^{-}_c)_{\mathrm{tor}}$ be the torsion part of $\mathscr{T}^1_{\mathrm{QG}}(\mathscr{W}^{-}_c)$ supported on $\mathrm{Ex}(g_c)$ and set around $\mathrm{Ex}(g_c)$
$$
\mathscr{T}^1_{\mathrm{QG}}(\mathscr{W}^{-}_c)_{\mathrm{tf}}:=\mathscr{T}^1_{\mathrm{QG}}(\mathscr{W}^{-}_c)/\mathscr{T}^1_{\mathrm{QG}}(\mathscr{W}^{-}_c)_{\mathrm{tor}}.
$$
Since $\mathscr{T}^1_{\mathrm{QG}}(W^{-})$ deforms to $\mathscr{T}^1_{\mathrm{QG}}(\mathscr{W}^{-}_c)$ for any general $c\in C$, 
the torsion-free and torsion parts satisfy:
$$
\mathscr{T}^1_{\mathrm{QG}}(\mathscr{W}^{-}_c)_{\mathrm{tf}}\cong 
\mathscr{T}^1_{\mathrm{QG}}(W^{-})_{\mathrm{tf}},
\quad
\mathscr{T}^1_{\mathrm{QG}}(\mathscr{W}^{-}_c)_{\mathrm{tor}}\cong
\mathscr{T}^1_{\mathrm{QG}}(W^{-})_{\mathrm{tor}}.
$$
From this, we conclude that the germ $(\mathrm{Ex}(g_c)\subset \mathscr{W}^{-}_c)$ is also $\mathbb{Q}$-Gorenstein smoothable. 
We complete the proof of the lemma.
\end{proof}

\begin{proof}[Proof of Proposition \ref{prop--partial--smoothing--from--non-Fano-to-Fano-with-anti-P}]
Let $f\colon (\mathscr{W},\frac{1}{2}\mathscr{B})\to C$ be as in the proof of Corollary \ref{cor--partial--smoothing--from--non-Fano-to-Fano}.
If the anti-P-resolution $\mu$ is of normal type, then the claim follows immediately.
Indeed, by performing two successive blow-ups $\mathscr{W}''\to \mathscr{W}'\to \mathscr{W}_{\mathrm{good}}$ along sections over $C$ if the length of $\mu$ is even or setting $\mathscr{W}:=\mathscr{W}_{\mathrm{good}}$ otherwise, and then contracting the resulting family of T-chains $\mathscr{W}''\to \mathscr{W}^{-}$, we obtain the desired family $g\colon \mathscr{W}^{-}\to C$.
In this case, we can define $\mathscr{B}^{-}$ as the effective divisor such that the natural map $(\mathscr{W}^{-}, \frac{1}{2}\mathscr{B}^{-})\to (\mathscr{W}, \frac{1}{2}\mathscr{B})$ is log crepant.

Thus, we may assume that $\mu$ is of non-normal type.
Let $\widetilde{W}$ be the normalization of $W^{-}$, and let $\widetilde{B}$ be the proper transform of $B^{-}$.
By the construction of $f$ and an argument similar to that used in the proof of Proposition \ref{prop--Lee-Park-from-I*-to-I}, there exists a projective flat morphism $h\colon\widetilde{\mathscr{W}} \to C$ with a $\mathbb{Q}$-Cartier effective relative Mumford divisor $\widetilde{\mathscr{B}}$ and a birational contraction $\xi\colon \widetilde{\mathscr{W}} \to \mathscr{W}$ satisfying the assumption of Lemma \ref{lem--gluing--tool}.
We complete the proof. 
\end{proof}

   Finally, we can detect smoothability of non-supersingular standard Horikawa surfaces.

\begin{thm}\label{thm--complete--classif-of-sm-st-horikawa}
Let $X$ be a standard Horikawa surface with a double cone singularity of type $(p_g-2;k'_1,k'_2)$, where $p_g \ge 7$.
Suppose that $X$ is not supersingular.
Then the following are equivalent.
\begin{itemize}
    \item[$(1)$] $X$ is $\mathbb{Q}$-Gorenstein smoothable.
    \item[$(2)$] The elliptic double cone singularity is equivariantly smoothable with respect to the good involution on $X$.
    \item[$(3)$] 
    One of the following hold.
    \begin{itemize}
    \item[$(a)$] $\max\{k_1',k_2'\}\ge 2p_g-12$, or
    \item[$(b)$] 
   $W$ admits an admisible anti-P-resolution $\mu\colon W^{-} \to W$ log crepant with respect to $(W,\frac{1}{2}B)$ of non-normal type $\mathrm{II}$ and the germ $\mathrm{Ex}(\mu)\subset W^{-}$ is $\mathbb{Q}$-Gorenstein smoothable.
    \end{itemize}
\end{itemize}
\end{thm}

\begin{proof}
The implication $(1)\Rightarrow (2)$ is trivial by Corollary \ref{cor--too--trivial}.
We first consider the implication $(2)\Rightarrow (3)$.
Assume that the germ of $x\in X$ is equivariantly $\mathbb{Q}$-Gorenstein smoothable, but $(a)$ does not hold.
By Lemma~\ref{lem--anti--p-resol--k-square}, there exists an admissible anti-P-resolution $\mu\colon W^{-} \to W$ which is log crepant with respect to $(W,\frac{1}{2}B)$ such that $\mathrm{Ex}(\mu)\cong\mathbb{P}^1$ and the germ $(\mathrm{Ex}(\mu)\subset W^{-})$ is $\mathbb{Q}$-Gorenstein smoothable.
By Theorem \ref{thm--anti-P-k_1-k_2}, this implies that $\mu$ has length $(k_1,k_2)$ such that $k_1k_2>0$.
Thus, Theorem \ref{thm--anti-P-k_1-k_2} implies that $\mu$ is of non-normal type $\mathrm{II}$. 
Therefore, $(b)$ holds.

Finally, we turn to the implication $(3)\Rightarrow (1)$.
First, Suppose that $(a)$ holds.
By Proposition \ref{prop--partial--smoothing--from--non-Fano-to-Fano-with-anti-P} and Lemma \ref{lem--vanishing--for--standard}, the surface $X$ can be smoothed to standard Horikawa surfaces of double Fano type without changing the germ of the anti-P-resolution.
Thus, we may assume that $k'_1=k_1$, $k'_2=k_2$ and $X$ is of double Fano type.
By Proposition \ref{prop--smoothing--c-type}, this implies $(1)$, as required.
On the other hand, suppose that $(b)$ holds.
Since $\mu$ is of non-normal type $\mathrm{II}$, we see that $$2p_g-12\le k_1+k_2\le 2p_g-10$$ by Theorem \ref{thm--anti-P-k_1-k_2}.
Thus, we may assume that $k'_1=k_1$, $k'_2=k_2$ and $X$ is of double Fano type again by Proposition \ref{prop--partial--smoothing--from--non-Fano-to-Fano-with-anti-P}.
In this case, it is not difficult to see from Lemma \ref{lem--vanishing--for--standard} that $-K_{W^{-}}$ is log big and nef and $B^{-}$ is nef.
Thus, the pair $(W^{-},\frac{1}{2}B^{-})$ satisfies all the conditions of Proposition \ref{prop--smoothing--slc}.
Applying Corollary \ref{cor--gen--fiber--rho} and the latter part of Proposition \ref{prop--smoothing--slc}, we conclude that $X$ is $\mathbb{Q}$-Gorenstein smoothable.
This completes the proof.
\end{proof}

In the rest of this subsection, we investigate the smoothability of supersingular standard Horikawa surfaces.
As shown in the following theorem, some admissible anti-P-resolutions are not $\Q$-Gorenstein smoothable.

\begin{thm}\label{thm--incomplete}
    Let $(x\in W,\frac{1}{2}B)$ be a cone singularity of type $(d;k_1,k_2)$ with $d\ge5$.
    Let $\mu\colon (W^-,\frac{1}{2}B^-)\to (W,\frac{1}{2}B)$ be an admissible anti-P-resolution of non-normal type $\mathrm{II}$ and of length $(l_1,l_2)$, where $0\le l_1\le k_1$ and $0\le l_2\le k_2$. 
    Assume that the germ $(\mathrm{Ex}(\mu)\subset W^-)$ is $\Q$-Gorenstein smoothable.
    Then, we have $l_1\le2$ or $l_2\le2$, and $\lceil\frac{l_1}{2}\rceil+\lceil\frac{l_2}{2}\rceil\ge d-3$.
\end{thm}

\begin{proof}
We may assume by the proof of Lemma \ref{lem--anti--p-resol--k-square} that $W$ is an open subset $\overline{\Sigma}_d$ and $\mathrm{Ex}(\mu)\cong\mathbb{P}^1$.
Let $\tilde{\mu}\colon \widetilde{W}^-\to \overline{\Sigma}_d$ be the birational morphism obtained by gluing $\overline{\Sigma}_d$ and $W^-$.
\begin{claim*}
    We may assume that there exists an effective Weil divisor $\widetilde{B}$ on $\overline{\Sigma}_d$ such that the proper transform of $\widetilde{B}$ on $\Sigma_{d}$ is not tangent to the proper transform of any fiber of the ruling $\Sigma_d\to\mathbb{P}^1$ over $\Delta_0$ and $\tilde{\mu}$ is an admissible anti-P-resolution with respect to $(\overline{\Sigma}_d,\frac{1}{2}\widetilde{B})$.
\end{claim*}

\begin{proof}[Proof of Claim]
    This follows from Lemma \ref{lem--vanishing--for--standard} (1) and arguments similar to those in Lemma \ref{lem--vanishing--for--standard} (2), Proposition \ref{prop--partial--smoothing--from--non-Fano-to-Fano-with-anti-P} and Lemma \ref{lem--gluing--tool}.
    Indeed, it is easy to see that Lemma \ref{lem--vanishing--for--standard} (1) also holds in this case.
    Let $(\overline{\Sigma}_d)_{\mathrm{good}}$ be the good resolution of $\Sigma_d$ with respect to $\widetilde{\mu}$.
    Let $\eta \colon \Sigma_d\to\overline{\Sigma}_d$ be the minimal resolution.
    For any sufficiently large $m\in\mathbb{Z}_{>0}$,
    set $D:=\eta_*(2m(\Delta_0+d\Gamma)+4\Gamma)$.
    Let $D'_{\mathrm{good}}$ be the associated Cartier divisor on $(\overline{\Sigma}_d)_{\mathrm{good}}$ as in Definition~\ref{defn--good-resolution}.
    Then $\mathcal{O}_{(\overline{\Sigma}_d)_{\mathrm{good}}}(D'_{\mathrm{good}})$ is big and globally generated, and 
    $$  H^1(\mathcal{O}_{(\overline{\Sigma}_d)_{\mathrm{good}}}(D'_{\mathrm{good}}))=0
   $$
   by the proof of Lemma \ref{lem--vanishing--for--standard} (2) for any sufficiently divisible $m>0$. 
   Thus, the similar argument to the proofs of Proposition \ref{prop--partial--smoothing--from--non-Fano-to-Fano-with-anti-P} and Lemma \ref{lem--gluing--tool} works, and we can deform $\tilde{\mu}$ to an admissible anti-P-resolution $\tilde{\mu}'\colon\widetilde{W}'^-\to\overline{\Sigma}_d$ with respect to $(\overline{\Sigma}_d,\frac{1}{2}\widetilde{B})$ for general $\widetilde{B}\in |D|$ of the same length and the same type as $\tilde{\mu}$.
   By the argument of the proof of Lemma \ref{lem--gluing--tool}, $(\mathrm{Ex}(\tilde{\mu}')\subset \widetilde{W}'^-)$ is again $\Q$-Gorenstein smoothable.
   Then, we may replace $\tilde{\mu}$ with $\tilde{\mu}'$ to show the assertion of Theorem \ref{thm--incomplete}.
\end{proof}

Returning to the proof of the theorem, replace $\mu$ with $\tilde{\mu}$.
By Theorem \ref{thm--hacking--prokhorov}, $W^-$ is globally $\Q$-Gorenstein smoothable.
    Moreover, $-K_{W^-}$ is log big and nef by the same argument of the proof of Corollary \ref{cor--properties--anti-P-standard}, so $W^{-}$ is $\mathbb{Q}$-Gorenstein smoothable to $\Sigma_0$ or $\Sigma_1$.
    Let $\pi\colon\mathscr{W}^-\to C$ be such a  smoothing.
    Note that $\mathcal{W}^-$ has only canonical singularities by \cite{kawakita}.
     Assume that $l_1>2$ and $l_2>2$, or $\lceil\frac{l_1}{2}\rceil+\lceil\frac{l_2}{2}\rceil< d-3$.
     Note that in the latter case, both $l_1$ and $l_2$ are even, and $l_1+l_2=2d-8$ (see Example \ref{exam--anti-P-type-A}).
    Let $\Gamma_1$ and $\Gamma_2$ be the two fibers where $(\overline{\Sigma}_d)_{\mathrm{good}}\to \Sigma_{d}$ is not an isomorphism.
    Let $D_1$ and $D_2$ be their proper transforms.
    Since the proper transform of $\widetilde{B}$ does not tangent to the proper transform of any fiber on $\Delta_0$, we can check that the proper transform of $\widetilde{B}$ does not meet both $D_1$ and $D_2$ and that the proper transforms of $D_i$ intersects the slc center of $W^-$ on $(\overline{\Sigma}_d)_{\mathrm{good}}$ at one point transversally
    by the same argument as in the proof of Corollary~\ref{cor--properties--anti-P-standard}.
    Hence, we have 
    $$
    K_{W^-}\cdot D_1=K_{W^-}\cdot D_2=0.
    $$
By applying \cite[Theorems 1.10, 1.15]{fujino--slc--vanishing} to $K_{W^-}$, we conclude that $-K_{\mathscr{W}^-}$ is $\pi$-big and $\pi$-semiample.
Define
$$
\mathscr{Y}:=\mathbf{Proj}_{C}(\oplus_{m\ge0}\pi_*\mathcal{O}_{\mathscr{W}^-}(-mK_{\mathscr{W}^-})).
$$
Since the canonical morphism $g\colon \mathscr{W}^-\to\mathscr{Y}$ is log crepant, $\mathscr{Y}$ has only canonical singularities.
By Lemma \ref{lem--lower-semi--rho} and shrinking $C$, we may assume that $\rho(\mathscr{W}^-/C)=2$ and hence $\rho(\mathscr{Y}/C)=1$.
Since $-K_{\mathscr{W}^-}$ is also relatively log big and nef over $\mathscr{Y}$, $R^1g_*\mathcal{O}_{\mathscr{W}^-}=0$ by \cite[Theorem 1.10]{fujino--slc--vanishing}.
However, since $g$ is small and contracts two curves, this contradicts \cite[12.1.5]{KM2} and \cite[Lemma 5.1.4]{kawakita--threefold}.
Therefore, we conclude the proof.
\end{proof}

Now, we state the following main result on smoothability of supersingular standard Horikawa surfaces.

\begin{thm}\label{thm--supersingular}
    Let $X$ be a supersingular standard Horikawa surface.
    Then the following hold.
\begin{itemize}
\item[$(1)$] 
$X$ is smoothable if and only if $X$ has only equivariantly smoothable double cone singularities with respect to the good involution.
\item[$(2)$] 
Suppose that $p_g(X)\ge 11$.
If the associated pair $(W,\frac{1}{2}B)$ admits an admissible anti-P-resolution $\mu\colon (W^-,\frac{1}{2}B^-)\to (W,\frac{1}{2}B)$ such that $(\mathrm{Ex}(\mu)\subset W^-)$ is $\Q$-Gorenstein smoothable, then $\mu$ is not supersingular.
\item[$(3)$] 
If $p_g(X)\le 9$, then $X$ is partially smoothable to standard Horikawa surfaces of double Fano type with the double cone singularity of type $(p_g(X)-2;2p_g(X)-12,0)$.
In particular, $X$ is smoothable.
\end{itemize}
\end{thm}

We consider the following concept of supersingular standard Horikawa surfaces.
\begin{defn}[Supersingular index]\label{defn--supersingular--index}
    Let $X$ be a standard Horikawa surface of type $(p_g-2)'$ with $p_g\ge 7$.
    Let $\Gamma$ be a fiber of the ruling $W^{+}=\Sigma_{p_g-2}\to \mathbb{P}^1$ contained in $B^{+}$.
We say that $X$ has {\it supersingular index} $s(X;\Gamma)$ {\it along} $\Gamma$ if $B^+-2\Delta_0-\Gamma$ intersect $\Gamma$ at $\Gamma\cap \Delta_0$ with multiplicity $s(X;\Gamma)$.
We define the {\it supersingular index} of $X$ as $s(X):=\max_{\Gamma}s(X;\Gamma)$.
Note that $0\le s(X)\le 4$ and if $0\le s(X)\le1$ (resp.~$s(X)=2$), then $X$ is of double Fano (resp.~of double non-Fano) type.
\end{defn}

We will deal with some special cases as follows.
\begin{prop}\label{prop--super--singular--typeC}
    Let $X$ be a standard normal Horikawa surface as Definition \ref{defn--supersingular--index}.
    Then, the following hold.
    \begin{itemize}
        \item[$(1)$] Suppose that $7\le p_g(X)\le9$ and $s(X)=p_g(X)-6$. Then $X$ is partially smoothable to standard Horikawa surfaces of double Fano type whose double cone singularity is of type $(p_g-2;2p_g-12,0)$.
        In particular, $X$ is smoothable.
        \item[$(2)$] Suppose that $p_g(X)=10$ and $s(X)=4$.
        Then, $X$ is $\mathbb{Q}$-Gorenstein smoothable to smooth Horikawa surfaces of type $(6)$.
    \end{itemize}
\end{prop}

\begin{proof}
We first treat case (1).
If $p_g(X)\le 8$, then the assertion follows from Proposition \ref{prop--smoothing--c-type}.
Therefore, we may assume that $p_g(X)=9$.
In this case, we show that $X$ is smoothable to standard Horikawa surfaces of double non-Fano type that are smoothable.
Let $\Gamma$ be a fiber of $\Sigma_d$ such that $s(X;\Gamma)=s(X)$.
There exists an admissible anti-P-resolution $\mu\colon (W^-,\frac{1}{2}B^-)\to (\overline{\Sigma}_d,\frac{1}{2}B)$ such that $W^-$ has only klt singularities and the proper transform $\Gamma'$ of $\Gamma$ on the minimal resolution $h\colon W^{-}_{\min}\to W^-$ is a $(-3)$-curve.
In this case, it is easy to see that $-K_{W^-}$ is big.
To prove the assertion, it suffices to show that $H^1(\mathcal{O}_{W^-}(B^-))=0$ by \cite[Proposition 3.1]{HP}, Proposition \ref{prop--smoothing--quotient} and Corollary \ref{cor--gen--fiber--rho}.
By an explicit calculation, the divisor $B^--K_{W^-}$ admits a $\mathbb{Z}$-Zariski decomposition with the $\mathbb{Z}$-positive part $B^--K_{W^-}-h_*\Gamma'$ and the $\mathbb{Z}$-negative part $h_*\Gamma'$.
Then, \cite[Theorem 1.3]{Enokizono} implies
\[
H^1(\mathcal{O}_{W^-}(B^-))\cong H^1(\mathcal{O}_{h_*\Gamma'}(\lfloor B^-|_{h_*\Gamma'}\rfloor))\cong H^1(\mathcal{O}_{\mathbb{P}^1}(-1))=0.
\]
By the same argument of Proposition \ref{prop--Lee-Park-from-I*-to-I}, we obtain a projective flat morphism $f\colon \mathscr{W}^-\to C$ from a normal klt threefold to a smooth curve with a closed point $0\in C$ such that
\begin{itemize}
    \item 
    $\mathscr{W}^-_0\cong W^-$.
    \item 
    $\mathscr{W}^-_c$ is a non-supersingular admissible anti-P-resolution of $(\overline{\Sigma}_7,\frac{1}{2}B')$ for some $B'\in |32\overline{\Gamma}|$.
    \item 
    The double cover of $\overline{\Sigma}_7$ branched along $B'$ is a standard Horikawa surface of double non-Fano type.    
\end{itemize}
The vanishing $H^1(\mathcal{O}_{W^-}(B^-))=0$, together with the same argument of the proof of Proposition \ref{prop--Lee-Park-from-I*-to-I}, implies the existence of a $\Q$-Cartier relative Mumford divisor $\mathscr{B}^-$ such that $\mathscr{B}^-_0=B^-$, $K_{\mathscr{W}^-}+\frac{1}{2}\mathscr{B}^-$ is $f$-big and $f$-semiample.
Let $g\colon \mathscr{W}^-\to\mathscr{W}$ be the relative canonical model of $(\mathscr{W}^-,\frac{1}{2}\mathscr{B}^-)$ over $C$.
Then $(\mathscr{W},\frac{1}{2}g_*\mathscr{B}^-)$ gives a log $\Q$-Gorenstein deformation of $(\overline{\Sigma}_7,\frac{1}{2}B)$ to the associated pairs of non-supersingular standard Horikawa surfaces of double non-Fano type with the double cone singularity of type $(7;6,0)$.
Hence, the assertion (1) follows from Corollary \ref{cor--partial--smoothing--from--non-Fano-to-Fano} and Proposition \ref{prop--smoothing--c-type}.

   We now treat the case (2).
    Let $\Gamma$ be a fiber of $\Sigma_d$ such that $s(X;\Gamma)=4$.
    Then $(\overline{\Sigma}_d,\frac{1}{2}B)$ admits an admissible anti-P-resolution $\mu\colon (W^-,\frac{1}{2}B^-)\to (\overline{\Sigma}_d,\frac{1}{2}B)$ such that $W^-$ is normal and the proper transform $\Gamma'$ of $\Gamma$ is a $(-6)$-curve.
    Note that $B^--\Gamma'$ is effective and disjoint from $\Gamma'$.
    Contracting $\Gamma'$ gives a birational morphism $\eta'\colon W^-\to W'$, and 
    we define $B':=\eta'_*(B^--\Gamma')$.
    Since $\rho(W')=1$, it follows from Proposition \ref{prop--smoothing--quotient} that there exists a projective flat morphism $f'\colon \mathscr{W}'\to C$ from a normal klt threefold to a smooth curve with a closed point $0\in C$, and a $\Q$-Cartier relative Mumford divisor $\mathscr{B}'$, such that 
    \begin{itemize}
        \item 
        $(\mathscr{W}'_0,\mathscr{B}'_0)\cong (W',B')$.
        \item 
        $f'$ is a $\Q$-Gorenstein smoothing of the singularity of type $\frac{1}{81}(1,8)$.
        \item 
        $f'$ is locally trivial around $\frac{1}{6}(1,1)$.
    \end{itemize}
    We may assume that there exists a projective birational morphism $\pi\colon \mathscr{W}^-\to\mathscr{W}'$ which is an isomorphism outside the section of the singularity of type $\frac{1}{6}(1,1)$, and provides a simultaneous resolution near that section.
By \cite[Proposition 2.6]{HP} and Lemma \ref{lem--lower-semi--rho}, we have $\rho(\mathscr{W}^-_{c})=2$ for general $c\in C$, hence $\mathscr{W}^-_{c}\cong \Sigma_6$.
Set $\mathscr{B}^-:=\pi^{-1}_*\mathscr{B}'+\mathscr{E}$, where $\mathscr{E}$ is the $\pi$-exceptional divisor.
By Proposition \ref{prop--smoothability--criterion}, we may assume that there exists a divisorial sheaf $\mathscr{L}$ such that $\mathscr{L}^{[2]}\sim_C\mathscr{B}^-$. 
Taking the relative canonical model of $(\mathscr{W}^-,\frac{1}{2}\mathscr{B}^-)$, we obtain a birational morphism $\pi^-\colon \mathscr{W}^-\to\mathscr{W}$.
By Corollary \ref{cor--gen--fiber--rho}, we have $\mathscr{W}^-\setminus\mathscr{W}^-_0\cong \mathscr{W}\setminus\mathscr{W}_0$.
Thus, the pair $(\mathscr{W},\frac{1}{2}\pi^{-}_*\mathscr{B}^-)$ defines a partial $\Q$-Gorenstein smoothing of $X$ to standard Horikawa surfaces of type $(6)$.
This completes the proof.
\end{proof}

To state the next lemma, we define the following notion.

\begin{defn}
    Let $d$ be a positive integer, and let $\Gamma$ be a fiber of the ruling $\Sigma_d\to \mathbb{P}^1$.
For each $m\in\mathbb{Z}_{>0}$, we define a birational birational morphism $h_m\colon W_m\to \Sigma_d$, called the {\it $m$-th blow-up along $\Gamma$},inductively as follows:
\begin{itemize}
\item 
For $m=1$, define $h_1\colon W_1\to \Sigma_d$ as the blow-up at the point $\Delta_0\cap F$.

\item 
For $m>1$, assuming $h_{m-1}$ has been defined,
let $\psi_{m}\colon W_{m}\to W_{m-1}$ be the blow-up at the intersection point of the proper transform of $\Gamma$ and the other irreducible components of $h_{m-1}^*\Gamma$.
Then set $h_m:=h_{m-1}\circ \psi_m$
\end{itemize}

Let $B_0$ be a divisor on $\Sigma_d$.
We define a sequence of divisors $B_m$ on $W_m$ for $m\ge 1$ inductively as $B_m:=\psi_{m}^{*}B_{m-1}-2E_{m}$,
where $E_{m}$ is the exceptional curve of $\psi_{m}\colon W_{m}\to W_{m-1}$ ($\psi_{1}:=h_1$ for $m=1$).
We refer to $B_m$ as the {\it associated divisor} to $B_{0}$.
\end{defn}

\begin{lem}\label{lem--supersingular--basepointfree}
Let $d$ be a positive integer. 
Let $\Gamma$ be a fiber of $\Sigma_d$, and set $B_0:= 4\Delta_0+(4d+4)\Gamma$ on ${\Sigma}_d$.
Then, the following hold:
\begin{itemize}
    \item[$(1)$]
    Let $1\le m\le 4$ and suppose $5\le d\le 8$.
    Consider the $m$-th blow-up $h_m\colon W_m\to \Sigma_d$ along $\Gamma$, and let $B_m$ be the associated divisor to $B_{0}$.
    Then the divisor $B_m-\Gamma_m$ is globally generated, where $\Gamma_m$ denotes the proper transform of $\Gamma$ on $W_{m}$.
    \item[$(2)$]
    Suppose that $d=8$.
    Let $\Gamma'$ be a fiber of $\Sigma_8$ distinct from $\Gamma$.
    Consider the fourth blow-up $h_4\colon W_4 \to \Sigma_8$ along $\Gamma$, and $B_4$ the associated divisor.
    For $5\le i\le 8$, define blow-ups
    $\psi_{i}\colon W_{i}\to W_{i-1}$ with the exceptional divisor $E_i$ and divisors $B_{i}$ on $W_{i}$ as follows:
\begin{itemize}
    \item 
    Define $\psi_{5}$ as the blow-up at the point $\Delta_{0}\cap \Gamma'$, and set $B_{5}:=\psi_{4}^{*}B_4-2E_{5}$.
    \item 
    For $6\le i\le 8$, we inductively define $\psi_{i}$ as the blow-up at a general point on $E_{i-1}$, and set $B_{i}:=\psi_{i}^{*}B_{i-1}-2E_{i}$.   
\end{itemize}
    Then the divisor $B_8-\Gamma_{8}$ is globally generated, where $\Gamma_8$ is the proper transform of $\Gamma$ on $W_8$. 
\end{itemize}
\end{lem}

\begin{proof}
    We focus on proving (2), as it is more involved; 
    the case (1) is simpler and will be omitted.
    Note that $h^0(-K_{W_8})\ge 6$ and $W_8$ satisfies the assumption of Lemma~\ref{lem--modified--manetti--lemma}.
    Therefore, if $D\subset W_8$ is an irreducible curve with $D^2<0$, then $D\cong\mathbb{P}^1$ and $D$ is either the proper transform of $\Delta_0$ or an irreducible component of a fiber of the composed morphism $W_8 \to\Sigma_8\to \mathbb{P}^1$.
    By using Lemma~\ref{lem:Z-poschar} and Remark \ref{rem:Z-poschar}, it is straightforward to check that $B_8-\Gamma_8-D-K_{W_8}$ is big and $\mathbb{Z}$-positive.
    Applying \cite[Theorem 1.3]{Enokizono} to $B_8-\Gamma_8-D$, we conclude that the restriction map 
    \[
    H^0(\mathcal{O}_{W_8}(B_8-\Gamma_8))\to H^0(\mathcal{O}_D((B_8-\Gamma_8)|_D))
    \]
    is surjective for all such curves $D$.
    This implies that the linear system $|B_8-\Gamma_8|$ has no base points along $D$.
Finally, examining the structure of $B_8-\Gamma_8$, we see that $|B_8-\Gamma_8|$ is basepoint free.
\end{proof}

We obtain immediate corollaries as follows.

\begin{cor}\label{cor--supersingular--part--smoothing}
 Let $X$ be a supersingular standard Horikawa surface with geometric genus $p_g$ and supersingular index $s$. 
 If $p_g\le10$,
  then $X$ is partially smoothable to standard Horikawa surfaces with supersingular index $s-1$ as in Definition \ref{defn--supersingular--index}.
  
  In particular, if $p_g\le 9$ and $s\ge p_g-6$, then $X$ is partially smoothable to standard Horikawa surfaces of double Fano type whose double cone singularity is of type $(p_g-2;2p_g-12,0)$.
\end{cor}

\begin{proof}
This immediately follows from Lemma \ref{lem--supersingular--basepointfree} (1) and Proposition \ref{prop--super--singular--typeC}.
\end{proof}

\begin{cor}\label{cor--p_g=10--connected}
    There exists a standard supersingular Horikawa surface $X$ with geometric genus $10$ and supersingular index $4$ satisfying the following condition:
    There exist two fibers $\Gamma$ and $\Gamma'$ of the ruling $W^{+}=\Sigma_8\to \mathbb{P}^1$ such that 
    \begin{itemize}
    \item 
    $\Gamma$ is contained in $B^{+}$ with $s(X;\Gamma)=4$.
    \item 
    $\Gamma'$ is not contained in $B^{+}$.
    \item 
    $B^+-2\Delta_0-\Gamma$ intersects $\Delta_0$ at the two points $\Gamma \cap \Delta_0$ and $\Gamma' \cap \Delta_0$.
    \item 
    There exists an admissible anti-P-resolution $\mu\colon (W^-,\frac{1}{2}B^-)\to (\overline{\Sigma}_8,\frac{1}{2}B)$ such that $W^-$ is klt, $-K_{W^-}$ is ample, and $\mu$ is not an isomorphism around the proper transform of $\Gamma'$.
    \end{itemize}
    In particular, there exists a stable normal Gorenstein Horikawa surface of geometric genus $10$ that is smoothable to both Horikawa surfaces of type~$(0)$ and type~$(6)$.
\end{cor}

\begin{proof}
    We can construct a divisor $B_0=B^{+}-2\Delta_0$
    with the required properties using Lemma \ref{lem--supersingular--basepointfree}.
    The first part of the corollary then follows from
    the argument given in the proof of Lemma \ref{lem--existence--smoothable--cusp}.
    The final claim follows from Proposition \ref{prop--smoothing--c-type} and Proposition \ref{prop--super--singular--typeC} (see also Example~\ref{exam:connectingGieseker}).
\end{proof}

\begin{proof}[Proof of Theorem \ref{thm--supersingular}]
    If $p_g(X)\le9$, then by Proposition \ref{prop--super--singular--typeC} and Corollary \ref{cor--supersingular--part--smoothing}, the surface $X$ is partially smoothable to standard Horikawa surfaces of double Fano type, which are in turn smoothable.
    This establishes claim (3), so we may assume that $p_g(X)\ge10$.

    We now consider the case $p_g(X)=10$.
    If $s(X)=4$, then the claim is a consequence of Proposition \ref{prop--super--singular--typeC}.
    Therefore, we focus on the case $s(X)=3$.
    If all singularities of $X$ are equivariantly smoothable, then there exists an admissible anti-P-resolution $\mu\colon(W^-,\frac{1}{2}B^-)\to (\overline{\Sigma}_d,\frac{1}{2}B)$ of length $(k_1,k_2)$ for some $k_1\ge0$ and $k_2\ge0$.
    By assumption, one of the integers $k_1$ or $k_2$ corresponds to a fiber that is contained in $B$.
   Suppose that this is the case for $k_2$.
   Then, $k_2$ must be even and at most six.
   Theorems \ref{thm--anti-P-k_1-k_2} and \ref{thm--incomplete} imply that in fact $k_2\le2$.
   In this case, it follows easily by arguments similar to those in the proof of Theorem \ref{thm--complete--classif-of-sm-st-horikawa} that $X$ is $\Q$-Gorenstein smoothable.

   Finally, we consider the case where $p_g(X)\ge11$.
   As seen in the proof of Theorem \ref{thm--complete--classif-of-sm-st-horikawa}, it suffices to prove statement (2) to deduce statement (1).
   Let us consider an arbitrary admissible anti-P-resolution $\mu\colon(W^-,\frac{1}{2}B^-)\to (\overline{\Sigma}_{p_g-2},\frac{1}{2}B)$ of length $(k_1,k_2)$.
   Suppose that $k_2$ corresponds to a fiber contained in $B$.
   Then, $k_2$ must be even and at most eight.
   Applying Theorems \ref{thm--anti-P-k_1-k_2} and \ref{thm--incomplete}, we find that $k_2\le2$,
   completing the proof of statement (2), and thus the theorem.
\end{proof}

We put the following corollary.

\begin{cor}\label{cor--properties--anti-P-standard-ii}
Let $X$ be a standard Horikawa surface of type $(p_g-2)'$ with $p_g\ge 11$.
    If $p_g=4k+2$ for some integer $k\ge 2$, then $X$ does not admit any $\mathbb{Q}$-Gorenstein smoothing to smooth Horikawa surfaces of type $(2k+2)$. 
\end{cor}

\begin{proof}
    Assume contrary that $X$ admits a $\Q$-Gorenstein smoothing to smooth Horikawa surfaces of type $(2k+2)$.
  Applying Proposition~\ref{prop--anti-P} to $(\mathscr{W}, \mathscr{B})$, we obtain an admissible anti-P-resolution $\mu\colon W^-\to W$ such that
  $W^{-}$ is $\Q$-Gorenstein smoothable to $\Sigma_{2k+2}$.
Note that $\mu$ is not supersingular by Theorem \ref{thm--supersingular}.
  By Lemma \ref{lem--upper--semiconti--plurigenera}, we have the inequality
  $$
  h^0(\widetilde{W},\mathcal{O}_{\widetilde{W}}(-K_{\widetilde{W}}))\ge h^0(W^-,\mathcal{O}_{W^-}(-K_{W^-}))\ge9.
  $$
  Then, combining Corollary \ref{cor--too--trivial}, Proposition \ref{prop--anti-P} and \cite[Proposition 1.41]{KoMo}, as well as the fact that the branch divisor on $\Sigma_{2k+2}$ associated with any Horikawa surface of type $(2k+2)$ is not ample (cf.~\cite[Lemma 1.4]{horikawa}),
  we conclude that $B^-$ is not ample.
This contradicts Corollary \ref{cor--properties--anti-P-standard}. 
Hence, the corollary follows. 
\end{proof}

\subsection{Horikawa surfaces of special Lee-Park type with $p_g\le 6$}\label{subsec:LPS-p_g<7-smoothable}

In this subsection, we consider Horikawa surfaces of special Lee-Park type for geometric genus at most $6$.


Let $m(X)$ be the multiplicity of $X$ introduced in Notation~\ref{note--section 9}.
Note that a Horikawa surface of general Lee-Park type $\mathrm{I}^*$ with geometric genus $4$ is nothing but of special Lee-Park type with multiplicity $0$.

\begin{lem}\label{lem--existence--low--special}
For any integers $4\le l\le6$ and $\lceil\frac{2l-8}{l-2}\rceil\le m\le 2$,
 there exists at least one Horikawa surface $X$ of special Lee-Park type with $p_g(X)=l$ and $m(X)=m$. 
\end{lem}

\begin{lem}\label{lem--Lee-Park--vanish--index}
 Let $X$ be a Horikawa surface of special Lee-Park type with $4\le p_g\le 6$, and let $\eta\colon W^{+}\to W$ and $E$ be as in Notation~\ref{note--section 9}.
 Then, the following hold:

\begin{itemize}
    \item[$(1)$] Both $W^{+}$ and $W$ are toric surfaces.
    \item[$(2)$] There exists a Weil divisor $L^{+}$ on $W^{+}$ such that $2L^{+}\sim B^{+}$.
    \item[$(3)$] $|B^{+}-m(X)E|$ is basepoint-free and $H^1(W^{+},\mathcal{O}_{W^{+}}(B^{+}-m(X)E))=0$.
\end{itemize}
\end{lem}

\begin{proof}
     From the construction of $W$, we see that $W^{+}$ is a toric surface with a unique Wahl singularity of type $\frac{1}{(p_g-1)^2}(1,p_g-2)$. 
    On the other hand, since 
    \[
    K_{W^{+}}+\frac{p_g-4}{p_g-2}E=\eta^*K_{W},
    \]
    it follows that $-K_{W^{+}}$ is big.
    By Theorem \ref{thm--hacking--prokhorov}, \cite[Proposition 2.6]{HP} and Corollary \ref{cor--deformation--rational}, it follows that $W^{+}$ is $\Q$-Gorenstein smoothable to some Hirzebruch surfaces.  
    Furthermore, by Lemma \ref{lem--upper--semiconti--plurigenera}, we obtain $h^0(-K_{W^{+}})\ge 9$.
    Let $\pi\colon W_{\mathrm{min}}\to W^{+}$ be the minimal resolution, and let $\psi \colon W_{\mathrm{min}}\to \Sigma_{p_g+1}$ be a birational morphism constructed as in Construction~\ref{construction}~(6).
    It is easy to see that there exists a Weil divisor $L^{+}$ on $W^{+}$ such that $2L^{+}\sim B^{+}$ and $B^+-2E$ is big.
    By Lemma \ref{lem--modified--manetti--lemma}, the surface $W_{\min}$ has a unique $(-1)$-curve, denoted by $C_1$.
    Let $C_2$ be the proper transform on $W_{\min}$ of the curve $E$.
    Note that neither $C_1$ or $C_2$ are $\pi$-exceptional.
    Applying Lemma \ref{lem--modified--manetti--lemma} again, we see that if $C\subset W^{+}$ is an irreducible and reduced curve with $C^2<0$, then $C=\pi_*C_1$ or $C=E=\pi_*C_2$.
    It is straightforward to check that $(B^{+}-jE)\cdot \pi_*C_i\ge 0$ for any $i=1,2$ and $\lceil\frac{2p_g-8}{p_g-2}\rceil\le j\le 2$.
    Therefore, by Lemma \ref{lem--nakai--moishezon}, the divisor $B^{+}-jE$ is nef.
    Since $W^{+}$ is toric and $B^{+}-jE$ is Cartier, it follows from  \cite[Theorems 6.1.7, 6.3.12 and 9.2.3]{CLS} that 
    $$
    H^1(W^{+},\mathcal{O}_{W^{+}}(B^{+}-jE))=0
    $$
    and the linear system $|B^{+}-jE|$ is basepoint-free.
\end{proof}

We can give a proof of Lemma \ref{lem--existence--low--special} by Lemma \ref{lem--Lee-Park--vanish--index} as in Remark \ref{rem--existence--toric} and we omit it.

\begin{cor}
Let $X$ be a Horikawa surface of special Lee-Park type with $4\le p_g\le 6$ and multiplicity $m(X)>\lceil\frac{2p_g-8}{p_g-2}\rceil$.
Then, $X$ is partially $\Q$-Gorenstein smoothable to Horikawa surfaces of special Lee-Park type with multiplicity $m(X)-1$.\label{cor--m(X)--down--special--type}
\end{cor}

\begin{proof}
This immediately follows from Lemma \ref{lem--Lee-Park--vanish--index}.
\end{proof}

\begin{prop}\label{prop--smoothability-of-special--low--p_g}
Let $X$ be a Horikawa surface of special Lee-Park type with $4\le p_g\le 6$.
Then, $X$ is partially $\Q$-Gorenstein smoothable to standard Horikawa surfaces with a double cone singularity with the same multiplicity $m(X)$.
In particular, $X$ is $\Q$-Gorenstein smoothable.
\end{prop}
\begin{proof}
   It is easy to check that both $-K_W$ and $B$ are ample and $\rho(W)=1$.
   By Theorem~\ref{thm--hacking--prokhorov}, there exists a projective flat morphism $f\colon \mathscr{W}\to C$ from a $\Q$-Gorenstein normal threefold to a smooth curve with a closed point $0\in C$ such that $\mathscr{W}_0\cong W$ and $f$ is locally trivial around the singularity of type $\frac{1}{p_g-2}(1,1)$ and a smoothing around the singularity of type $\frac{1}{(p_g-1)^2}(1,p_g-2)$. 
  By taking a finite cover of $C$ if necessary, we may take $g\colon \overline{\mathscr{W}}\to \mathscr{W}$,   a simultaneous resolution of the singularity of type $\frac{1}{p_g-2}(1,1)$ over $C$ (cf.~\cite[Theorem 2.22]{KSB}). 
  Let $\mathscr{E}$ be the exceptional divisor of $g$, and define $\overline{W}:=\overline{\mathscr{W}}_0$ and $E:=\mathscr{E}_0$.
  Note that there exists an effective Weil divisor $\overline{B}$ on $\overline{W}$ such that $(g_0)_*\overline{B}=B$ and 
  \[
  K_{\overline{W}}+\frac{1}{2}\overline{B}=(g_0)^*\left(K_W+\frac{1}{2}B\right).
  \]
  Let $\overline{L}$ be the Weil divisor on $\overline{W}$ such that $2\overline{L}\sim\overline{B}$.
  By construction, $\overline{L}$ is Cartier at the singularity of type $\frac{1}{(p_g-1)^2}(1,p_g-2)$.
  By Proposition \ref{prop--smoothability--criterion} and Lemma \ref{lem--Lee-Park--vanish--index}, we may assume that there exist a divisorial sheaf $\overline{\mathscr{L}}$ and an effective relative Mumford divisor $\overline{\mathscr{B}}'$ on $\overline{\mathscr{W}}$ such that $\overline{\mathscr{L}}_0\cong \overline{L}$ and $\overline{\mathscr{B}}'_0=\overline{B}-m(X)E$.
  By \cite[2.3]{kollar-moduli}, we note that $(\overline{\mathscr{W}},\frac{1}{2}(\overline{\mathscr{B}}'+m(X)\mathscr{E})+\overline{\mathscr{W}}_0)$ is lc.
  Since $\overline{W}$ is toric and $K_{\overline{W}}+\frac{1}{2}\overline{B}$ is nef, we have
  \[
  H^1(\overline{W},\mathcal{O}_{\overline{W}}(d(2K_{\overline{W}}+\overline{B})))=0
  \]
  for any sufficiently large and divisible $d$.
  By taking the relative ample model of $K_{\overline{\mathscr{W}}}+\frac{1}{2}(\overline{\mathscr{B}}'+m(X)\mathscr{E})$ over $C$, we get a birational morphism $\mu\colon \overline{\mathscr{W}}\to \mathscr{W}'$ to a normal variety which is flat over $C$.
  Since $(2K_{\overline{W}}+\overline{B})\cdot E=0$, $\mu$ contracts $\mathscr{E}$.
  By \cite[Proposition 2.6]{HP}, we see that $\rho(\overline{\mathscr{W}}_c)=2$ for any $c\in C$.
  Thus, $\rho(\mathscr{W}_c')=1$ and $\mathscr{W}= \mathscr{W}'$. 
  Set $\mathscr{B}:=\mu_*\overline{\mathscr{B}}'$.
  Then, it is easy to see that $(\mathscr{W}_c,\frac{1}{2}\mathscr{B}_c)$ is log canonical and by using the data of $(\mu_*\overline{\mathscr{L}})_c$, the double cover of $\mathscr{W}_c$ branch along $\mathscr{B}_c$ is a standard Horikawa surface $Y_c$ with one double cone singularity for any general $c\in C$.
  Moreover, it is easy to see that $m(Y_c)=m(X)$ by construction.
  Thus, we obtain the first assertion.
  The last assertion follows from Lemma \ref{lem--standard-for-low-p_g-smoothable}.
\end{proof}

\subsection{Horikawa surfaces of special Lee-Park type with $p_g\ge 7$}\label{subsec:LPS-p_g>6--smoothable}

The aim in this subsection is to show the following theorem for Horikawa surfaces of special Lee-Park type for any geometric genus and to conclude Theorem \ref{intro-thm:Q2}.

\begin{thm}\label{thm--smoothability-of-special--cusp}
Let $X$ be a Horikawa surface of special Lee-Park type with $p_g\ge 7$.
Then, the double cone singularity of $X$ is equivariantly smoothable if and only if $X$ is $\mathbb{Q}$-Gorenstein smoothable.
\end{thm}

Before presenting the proof of Theorem \ref{thm--smoothability-of-special--cusp}, we show the following lemma, which is analogous to Lemma \ref{lem--vanishing--for--standard}.

\begin{lem}\label{lem--vanishing--for--special}
Let $X$ be an arbitrary Horikawa surface of special Lee-Park type with $p_g\ge 7$.
Suppose that $(W,\frac{1}{2}B)$ has a singularity $w_1$ of type $(p_g-2;k'_1,k'_2)$.
Suppose further that there exists an admissible anti-P-resolution $\mu\colon (W^{-},\frac{1}{2}B^{-})\to (W,\frac{1}{2}B)$ of length $(k_1,k_2)$  such that $k_1+k_2\le 2p_g-6$ for some $0\le k_i\le k'_i$ for $i=1,2$.
   Take the good resolution of $h\colon W_{\mathrm{good}}\to W^{+}$ with respect to $\mu$, where $W^{+}\to W$ is the minimal resolution of the singularity $w_1$.
   
   Then, the following hold:
   \begin{itemize}
       \item[$(1)$]
       Let $W^{-}_{\min}$ be the minimal resolution of the normalization of $W^{-}$.
       Then, $-K_{W^{-}_{\min}}$ is big and
    \begin{equation}\label{eq--antigenus-of-anti-P-for-special}
  \min\{h^0(-K_{W^{-}_{\min}}),h^0(-K_{W_{\mathrm{good}}})\}\ge 6.
    \end{equation}
    In particular, $W_{\mathrm{good}}$ and $W^{-}_{\min}$ satisfy the assumption of Lemma \ref{lem--modified--manetti--lemma}.

    \item[$(2)$]
    Let $B'_{\mathrm{good}}$ be the proper transform of $B$ on the good resolution $W_{\mathrm{good}}$ defined in Definition~\ref{defn--good-resolution}.
    Then, $\mathcal{O}_{W_{\mathrm{good}}}(B'_{\mathrm{good}})$ is big and globally generated, and 
    $$  H^1(W_{\mathrm{good}},\mathcal{O}_{W_{\mathrm{good}}}(B'_{\mathrm{good}}))=0.
   $$  
\end{itemize}
   \end{lem}

\begin{proof}
Note that $h^0(-K_{W^+})\ge p_g+4$ by Lemma \ref{lem--upper--semiconti--plurigenera} and the fact that $W^+$ is $\Q$-Gorenstein smoothable to $\Sigma_{p_g-2}$.
    Applying the same argument to $W^+$ as those of the proofs of Proposition \ref{prop--construction--lee--park--i} and Lemma \ref{lem--vanishing--for--standard}, we obtain the proof.
\end{proof}

We state the fundamental property of Horikawa surfaces of special Lee-Park type.

\begin{lem}\label{lem--property-of-special}
    Let $X$ be a Horikawa surface of special Lee-Park type.
    Then, $-K_W$ and $L$ are ample, $\rho(W)=1$, the divisor $B$ does not pass through the singularity of type $\frac{1}{(p_g-1)^2}(1,p_g-2)$, and $K_W+L$ is Cartier at the cone singularity.
\end{lem}

It is easy to see that Lemma \ref{lem--property-of-special} holds and we omit the proof.

\begin{cor}\label{cor--properties--anti-P-special}
    Let $X$ be a Horikawa surface of special Lee-Park type 
    with geometric genus $p_g \ge 7$.
Let $\mu\colon(W^{-},\frac{1}{2}B^{-})\to (W,\frac{1}{2}B)$ be an admissible anti-P-resolution that is $\mathbb{Q}$-Gorenstein smoothable to Hirzebruch surfaces.
Then, the following hold:

\begin{itemize}
\item[$(1)$]
$B^{-}$ is ample.
\item[$(2)$]
If $p_g=4k+2$ for some $k\in\mathbb{Z}_{\ge2}$, then $X$ does not admit any $\mathbb{Q}$-Gorenstein smoothing to smooth Horikawa surfaces of type $(2k+2)$. 
\item[$(3)$]
$-K_{W^{-}}$ is log big and nef.
If $W^{-}$ is klt, then $-K_{W^-}$ is ample.    
\end{itemize}
\end{cor}

\begin{proof}
    It follows from the same argument as the proof of Corollaries \ref{cor--properties--anti-P-standard} and \ref{cor--properties--anti-P-standard-ii} by using Lemmas \ref{lem--vanishing--for--special} and \ref{lem--property-of-special}.
\end{proof}

Applying the above results, we see that there exists at least one $\Q$-Gorenstein smoothable Horikawa surface of special Lee-Park type as follows.

\begin{prop}\label{prop--smoothability-of-special--cusp--k_1=0}
Let $X$ be a Horikawa surface of special Lee-Park type with $p_g\ge 7$.
Suppose that the double cone singularity of $X$ is of type $(p_g-2;k'_1,k'_2)$ for some $k'_1\ge2p_g-12$ and $k'_2\ge 0$. 
 Then, $X$ is partially $\Q$-Gorenstein smoothable to standard Horikawa surfaces with a double cone singularity of type $(p_g-2;2p_g-12,0)$.
In particular, $X$ is $\Q$-Gorenstein smoothable.
\end{prop}

\begin{proof}
Let $w_1$ and $w_2$ denote the singularities on $W$ of type $\frac{1}{p_g-2}(1,1)$ and of type $\frac{1}{(p_g-1)^{2}}(1,p_g-2)$, respectively.
Note that $W$ admits an anti-P-resolution $\mu\colon W^-\to W$ of normal type in this case.
By Lemma \ref{lem--vanishing--for--special}, $X$ has a $\Q$-Gorenstein partial smoothing to stable normal Horikawa surfaces of special Lee-Park type with a double cone singularity of type $(p_g-2;2p_g-12,0)$.
Thus, we may assume that $k'_1=2p_g-12$ and $k'_2=0$.
Since $-K_W$ is ample by Lemma \ref{lem--property-of-special}, there exists a projective flat morphism $f\colon \mathscr{W}\to C$ from a $\mathbb{Q}$-Gorenstein normal threefold to a smooth curve $C$ with a closed point $0\in C$ such that $\mathscr{W}_0\cong W$ and $f$ is a $\mathbb{Q}$-Gorenstein smoothing around $w_2$ but locally trivial around $w_1$ by Theorem \ref{thm--hacking--prokhorov}.
Furthermore, we may assume that there exists a $\mathbb{Q}$-Cartier divisorial sheaf $\mathscr{L}$ on $\mathscr{W}$ such that $\mathscr{L}_0\cong L$ by Proposition \ref{prop--smoothability--criterion}. 
By \cite[Theorem 2.22]{KSB}, we see that there exists a projective birational morphism $g\colon \mathscr{W}_{\mathrm{good}}\to \mathscr{W}$ such that for any $c\in C$, $g_c$ is the composition of the minimal resolution at $w_1$ and the good resolution with respect to $\mu$ at $w_1$. 
Note that $g_0$ is isomorphic locally around $w_2$.
By Lemma \ref{lem--vanishing--for--special}, by shrinking $C$ if necessary, we may assume that there exists an effective relative Mumford divisor $\mathscr{B}_{\mathrm{good}}\sim_C\mathscr{L}_{\mathrm{good}}^{[2]}$ such that for any $c\in C$, $\mathscr{B}_{\mathrm{good}, c}$ is smooth, where $\mathscr{L}_{\mathrm{good}}$ is the line bundle such that $(\mathscr{L}_{\mathrm{good}})_c$ is the  line bundle associated to the good resolution.
Define $\mathscr{B}:=g_{*}\mathscr{B}_{\mathrm{good}}$.
Then, $(\mathscr{W}_c,\frac{1}{2}\mathscr{B}_c)$ defines a singularity of type $(p_g-2;2p_g-12,0)$ for any $c\in C$. 
Note also that $\mathscr{L}^{[2]}\sim_C \mathscr{B}$.
By taking the double cover of $\mathscr{W}$ branched along $\mathscr{B}$, we obtain a partial $\mathbb{Q}$-Gorenstein smoothing $f\colon \mathscr{X}\to C$ of $X$ to stable normal Gorenstein Horikawa surfaces with a unique double cone singularity of type $(p_g-2;2p_g-12,0)$.
It follows that $\mathscr{X}_c$ is standard for any general $c\in C$.
Note that locally around the singularity, $\mu$ is also an anti-P-resolution log crepant with respect to $(\mathscr{W}_c,\frac{1}{2}\mathscr{B}_c)$.

Finally, the last assertion follows from Proposition \ref{prop--smoothing--c-type} and the first assertion.
\end{proof}

Note that there exists a Horikawa surface as in Proposition \ref{prop--smoothability-of-special--cusp--k_1=0}.
This fact can be shown as Lemma \ref{lem--existence--smoothable--cusp}.
Next, we deal with the case where $W$ admits a non-normal anti-P-resolution.
Now, we want to show an analogous result to Proposition \ref{prop--partial--smoothing--from--non-Fano-to-Fano-with-anti-P} for the construction of partial $\Q$-Gorenstein smoothing of Horikawa surfaces of special Lee-Park type to standard Horikawa surfaces of double Fano type that preserves the $\Q$-Gorenstein smoothability.
To address this, we first prove the following.

\begin{prop}\label{prop--partial--smoothing--result--for--special}
    Let $X$ be a Horikawa surface of special Lee-Park type with $p_g\ge7$. 
    Then, $X$ does not admit any partial $\Q$-Gorenstein smoothing either to standard Horikawa surfaces of double non-Fano type or to supersingular standard Horikawa surfaces.
\end{prop}

\begin{proof}
Assume the contrary and 
then there exists a log $\Q$-Gorenstein deformation $f\colon (\mathscr{W},\frac{1}{2}\mathscr{B})\to C$ over a curve with a closed point $0\in C$ such that $(\mathscr{W}_0,\frac{1}{2}\mathscr{B}_0)$ is isomorphic to the quotient $(W,\frac{1}{2}B)$ of $X$ and $(\mathscr{W}_c,\frac{1}{2}\mathscr{B}_c)$ defines standard Horikawa surfaces with double cone singularities for any $c\in C\setminus\{0\}$.
Let $\xi\colon \mathscr{W}^+\to\mathscr{W}$ be the simultaneous resolution of the horizontal cone singularity and $\mathscr{B}^+$ the effective relative Mumford divisor such that $(\mathscr{W}^+,\frac{1}{2}\mathscr{B}^+)$ and $(\mathscr{W},\frac{1}{2}\mathscr{B})$ are log crepant.

First, assume that for any $c\in C\setminus\{0\}$, the support of $\mathscr{B}^+_c$ contains a fiber of the canonical ruling $\mathscr{W}^+_0\cong\Sigma_{p_g-2}\to\mathbb{P}^1$.
It is easy to see that there exists a relative Cartier divisor $\mathscr{F}^\circ\subset \mathscr{W}^+\times_C(C\setminus\{0\})$ such that $\mathscr{F}^\circ_c$ is a fiber for any $c\in C\setminus\{0\}$ and $\mathscr{F}^\circ\subset \mathrm{Supp}(\mathscr{B}^+)$.
Let $\mathscr{F}$ be the Zariski closure of $\mathscr{F}^\circ$.
We note that $\rho(\mathscr{W}^+_0)=2$ and $\mathrm{Eff}(\mathscr{W}^+_0)$ is spanned by $\mathrm{Ex}(\xi_0)$ and the proper transform $D$ of the unique $(-1)$-curve of the minimal resolution of $W$.
Note also that $(\mathrm{Ex}(\xi_0))^2<0$ and $D^2<0$.
We may assume that $\mathscr{W}^+$ is $\Q$-factorial by Lemma \ref{lem--lower-semi--rho}.
Hence, $\mathscr{F}_0^2=0$ and it is not hard to see that $D\subset \mathrm{Supp}(\mathscr{F}_0)$.
This means that $\mathscr{B}^+_0$ passes through the unique Wahl singularity.
This contradicts the log canonicity of $(\mathscr{W}^+_0,\frac{1}{2}\mathscr{B}^+_0)$.
 
 Next, we deal with the case where $(\mathscr{W}_c,\frac{1}{2}\mathscr{B}_c)$ defines a standard Horikawa surface of double non-Fano type for any $c\in C\setminus\{0\}$.
 By the previous paragraph, we may assume that the support of $\mathscr{B}^+_c$ contains no fiber for any $c\in C\setminus\{0\}$.
 We may further assume that there exists a section $\mathcal{S}_1\subset \mathrm{Ex}(\xi)$ over $C$ such that $\mathscr{B}^+_c$ and the fiber of the ruling intersect at $(\mathcal{S}_1)_c$ with multiplicity more than two for any $c\in C\setminus\{0\}$.
 Let $\mathscr{F}^\circ\subset \mathscr{W}^+\times_C(C\setminus\{0\})$ be the relative Cartier divisor such that $\mathscr{F}^\circ_c$ passes through $(\mathcal{S}_1)_c$ for any $c\in C\setminus\{0\}$.
 Take the horizontal blow up $\beta_1\colon \mathscr{W}_1\to\mathscr{W}^+$ along $\mathcal{S}_1$.
 Note that there exists a section $\mathcal{S}_2\subset\mathrm{Ex}(\beta_1)$ such that the proper transforms of $\mathscr{F}^\circ$ and $\mathscr{B}^+-2\mathrm{Ex}(\xi)$ intersect at $\mathcal{S}_2\times_C(C\setminus\{0\})$.
 Take the horizontal blow up $\beta_2\colon \mathscr{W}_2\to\mathscr{W}_2$ along $\mathcal{S}_2$.
By construction, it is not hard to see that there exists an effective relative Mumford divisor $\mathscr{B}_2$ such that $(\mathscr{W}_2,\frac{1}{2}\mathscr{B}_2)$ and $(\mathscr{W},\frac{1}{2}\mathscr{B})$ are log crepant.
Note that $(\beta_2)_*^{-1}\mathrm{Ex}(\beta_1)_c$ is a $(-2)$-curve for any $c\in C$.
By contracting $(\beta_2)_*^{-1}\mathrm{Ex}(\beta_1)$, we obtain a new family of surfaces $\mathscr{W}'\to C$.
We may assume that $\mathscr{W}'$ is $\Q$-factorial by Lemma \ref{lem--lower-semi--rho}. 
As Corollary \ref{cor--properties--anti-P-special}, we see that $\mathscr{W}'_0$ is klt and $-K_{\mathscr{W}'_0}$ is ample.
On the other hand, $\mathscr{W}'_c$ contains a $(-2)$-curve and $-K_{\mathscr{W}'_c}$ is not ample.
Note that $\mathscr{W}'$ is $\Q$-factorial by Lemma \ref{lem--lower-semi--rho}.
This contradicts \cite[Proposition 1.41]{KoMo}.
    Therefore, we obtain the assertion.
\end{proof}

Now, we can show the analog of Proposition \ref{prop--partial--smoothing--from--non-Fano-to-Fano-with-anti-P} for Horikawa surfaces of special Lee-Park type.

\begin{prop}\label{prop--smoothability-of-special--cusp--k_1k_2not0}
Let $X$ be a Horikawa surface of special Lee-Park type with $p_g\ge 7$.
Suppose that there exists an admissible anti-P-resolution $\mu\colon (W^{-}, \frac{1}{2}B^{-}) \to (W, \frac{1}{2}B)$ of non-normal type and of length $(k_1,k_2)$ for some $k_1\ge0$ and $k_2\ge0$ with $k_1+k_2\le 2p_g-6$.
Suppose further that the germ $(\mathrm{Ex}(\mu)\subset W^{-})$ is $\mathbb{Q}$-Gorenstein smoothable.
 Then, $X$ is partially $\Q$-Gorenstein smoothable to standard Horikawa surfaces of double Fano type with a double cone singularity that is equivariantly $\Q$-Gorenstein smoothable.
 
 In particular, $X$ is $\Q$-Gorenstein smoothable.
\end{prop}

\begin{proof}
    By Lemmas \ref{lem--vanishing--for--special}, \ref{lem--modified--manetti--lemma} and \ref{lem--nakai--moishezon}, we see that $-K_{W^{-}}$ is log big and nef and $B^{-}$ is nef.
    It is easy to see that there exists a divisorial sheaf $L^{-}$ such that $(L^{-})^{[2]}\cong\mathcal{O}_{W^{-}}(B^{-})$ and $K_{W^{-}}+L^{-}$ is Cartier along $\mathrm{Ex}(\mu)$.
    Note that $B^{-}$ is Cartier locally around the singularity $w$ of type $\frac{1}{(p_g-1)^2}(1,p_g-2)$.
    By Theorem \ref{thm--hacking--prokhorov}, we see that there exists a projective flat morphism $f\colon \mathscr{W}^{-} \to C$ from a $\mathbb{Q}$-Gorenstein slc scheme to a smooth curve with a closed point $0$ such that $\mathscr{W}^{-}_0\cong W^{-}$ such that $f$ is smoothing around $w$ but locally trivial for any other singularities.
    By Proposition \ref{prop--smoothability--criterion}, we may assume that there exists a $\mathbb{Q}$-Cartier divisorial sheaf $\mathscr{L}^{-}$ such that $\mathscr{L}^{-}_0\cong L^{-}$.
    Note that 
    \[
   H^1(\mathcal{O}_{W^{-}}(B^{-}))= H^1(\mathcal{O}_{W^{-}}(-dK_{W^{-}}))=H^1(\mathcal{O}_{W^{-}}(d(2K_{W^{-}}+B^{-})))=0
    \]
    for any sufficiently large and divisible $d\in\mathbb{Z}_{>0}$ by \cite[Theorem 1.10]{fujino--slc--vanishing}.
It shows that there exists an effective relative Mumford divisor $\mathscr{B}^{-} \sim (\mathscr{L}^{-})^{[2]}$ such that $\mathscr{B}^{-}_0=B^{-}$.
Furthermore, $K_{\mathscr{W}^{-}}+\frac{1}{2}\mathscr{B}^{-}$ is relatively big and semiample over $C$ and $-K_{\mathscr{W}^{-}}$ is relatively big over $C$.
Let $\pi\colon \mathscr{W}^{-}\to\mathscr{W}$ be the relative ample model of $K_{\mathscr{W}^{-}}+\frac{1}{2}\mathscr{B}^{-}$ and set $\mathscr{B}:=\pi_*\mathscr{B^{-}}$.
Note that $\mathscr{B}\sim\pi_*(\mathscr{L}^{-})^{[2]}$.
It is easy to see that $\mathscr{W}_0\cong W$.
Since $W$ is normal, we see that $\mathscr{W}_c$ is also normal for any general $c\in C$. 
This shows that $\pi_c\colon \mathscr{W}^{-}_c\to\mathscr{W}_c$ is not an isomorphism for any general $c\in C$ and the non-normal locus should be contracted. 

Next, we show that $\mathscr{W}^{-}_c$ is $\Q$-Gorenstein smoothable to $\Sigma_2$, $\Sigma_1$ or $\Sigma_0$ for any very general $c\in C$.
We note that $-K_{\mathscr{W}^{-}_c}$ is big and nef for any very general $c\in C$.
Moreover, by construction, we see that $\mathscr{W}^{-}_c$ has only $\mathbb{Q}$-Gorenstein smoothable slc singularities such that $K_{\mathscr{W}^{-}_c}^2=8$.
It is also easy to see that non-normal locus of $\mathscr{W}^{-}_c$ is isomorphic to $\mathbb{P}^1$ for any general $c\in C$.
As the last paragraph of the proof of Proposition \ref{prop--partial--smoothing--from--non-Fano-to-Fano-with-anti-P}, we see that by setting $N_c\subset \mathscr{W}^{-}_c$ as the non-normal locus, the germ $(N_c\subset \mathscr{W}^{-}_c)$ is also $\Q$-Gorenstein smoothable.
Since $-K_{\mathscr{W}^{-}_c}$ is big and nef and $(K_{\mathscr{W}^{-}_c})^2=8$, $\mathscr{W}^{-}_c$ is $\mathbb{Q}$-Gorenstein smoothable to $\Sigma_2$, $\Sigma_1$ or $\Sigma_0$ by \cite{D} and Theorem \ref{thm--hacking--prokhorov}.
This shows that $\rho(\mathscr{W}^{-}_c)\le 2$ by Lemma \ref{lem--lower-semi--rho}.
Since $\pi_c$ is not isomorphic, $\rho(\mathscr{W}_c)=1$, $\mathrm{Ex}(\pi_c)=N_c$ and $\mathscr{W}_c$ has at most one singularity.
Thus, we may assume that $\mathscr{W}$ is $\mathbb{Q}$-factorial by Lemma \ref{lem--lower-semi--rho}.
Since the singularity of type $\frac{1}{l}(1,1)$ has the minimal resolution as only one P-resolution for $l\ge5$, we see by \cite[Theorem 3.9]{KSB} that an arbitrary $\Q$-Gorenstein smoothing of this singularity should be trivial (see Remark \ref{rem--deformation-space-of-1/l(1,1)} below).  
Hence, it is easy to see that $\mathscr{W}_c\cong\overline{\Sigma}_{p_g-2}$.
Therefore, the double cover of $\mathscr{W}_c$ branched along $\mathscr{B}_c$ is a standard Horikawa surface $Y_c$ for any general $c\in C$.
Furthermore, $\pi_c$ is an anti-P-resolution log crepant with respect to $(\mathscr{W}_c,\frac{1}{2}\mathscr{B}_c)$ such that $(\mathrm{Ex}(\pi_c)\subset \mathscr{W}^{-}_c)$ is $\Q$-Gorenstein smoothable (cf.~the last assertion of Lemma \ref{lem--gluing--tool}).
On the other hand, $Y_c$ is of double Fano type by Proposition \ref{prop--partial--smoothing--result--for--special}.
By Theorem \ref{thm--complete--classif-of-sm-st-horikawa}, we see that $Y_c$ is $\mathbb{Q}$-Gorenstein smoothable.
We complete the proof.
\end{proof}

\begin{rem}\label{rem--deformation-space-of-1/l(1,1)}
    As in the proof above, we can completely describe the deformation space of the singularity of type $\frac{1}{l}(1,1)$ for any $l\ge2$.
The case where $l=2$ is well-known (cf.~\cite[Theorem 4.61]{KoMo}).
In the case where $l\ne2$ or $4$, the singularity admits the minimal resolution as only one P-resolution.
In the case where $l=4$, the singularity admits the minimal resolution and the identity morphism as P-resolutions.
Therefore, if $l\ne 4$, the reduced structure of the semiversal deformation space of the singularity of type $\frac{1}{l}(1,1)$ has one irreducible component and if $l=4$, this has exactly two irreducible components by \cite[Theorem 3.9]{KSB}.
For more details and structures of the semiversal deformation spaces, see \cite[Section 8]{Pink}.
By this, we conclude that an arbitrary nontrivial deformation of the singularity of type $\frac{1}{l}(1,1)$ should be smooth.
\end{rem}

\begin{proof}[Proof of Theorem \ref{thm--smoothability-of-special--cusp}]
Since the if part follows from Corollary \ref{cor--too--trivial}, it suffices to deal with the only if part.
    Suppose that the double cone singularity is equivariantly $\Q$-Gorenstein smoothable.
    Then, we see $(W,\frac{1}{2}B)$ admits an admissible anti-P-resolution $\mu$ by Proposition \ref{prop--anti-P}.
    By Proposition \ref{prop--smoothability-of-special--cusp--k_1=0} and Theorem \ref{thm--anti-P-k_1-k_2}, we may assume that $\mu$ is of non-normal type.
    In the case where $\mu$ is of non-normal type $\mathrm{II}$, $\mu$ has length $(k_1,k_2)$ such that $2p_g-12\le k_1+k_2\le 2p_g-10$ by Theorem \ref{thm--anti-P-k_1-k_2} and hence $X$ is $\Q$-Gorenstein smoothable by Proposition \ref{prop--smoothability-of-special--cusp--k_1k_2not0}.
    Consider the case where $\mu$ is of non-normal type $\mathrm{I}$.
    Theorem \ref{thm--anti-P-k_1-k_2} and Proposition \ref{prop--smoothability-of-special--cusp--k_1k_2not0} show that if $\mu$ is of length $(k_1,0)$, then $k_1\ge2p_g-12$.
     In this case, $W$ admits an anti-P-resolution of normal type and $X$ is $\Q$-Gorenstein smoothable by Proposition \ref{prop--smoothability-of-special--cusp--k_1=0}.
    We complete the proof.
\end{proof}

\subsection{Conclusion on $\Q$-Gorenstein smoothability}\label{subsec:conclusion--equiv--smoothability}

Now, we are ready to answer Question \ref{question} (2). 

\begin{proof}[Proof of Theorems \ref{thm--main--ii} and \ref{intro-thm:Q2}]
These follow from Lemma \ref{lem--standard-for-low-p_g-smoothable}, Propositions \ref{prop--canonical--pencil--invol}, \ref{prop--smoothing-p_g=3}, \ref{prop--smoothing--lee--park--i}, \ref{prop--Lee-Park-from-I*-to-I}, \ref{prop--smoothing--urzua-type}, \ref{prop--smoothability-of-special--low--p_g} and Theorems \ref{intro-qGsmHor}, \ref{thm--complete--classif-of-sm-st-horikawa} and \ref{thm--smoothability-of-special--cusp}.
\end{proof}

We propose the following question, which remains open and, to our knowledge, has not been addressed in the existing literature. 
\begin{ques} \label{ques--cusp} \phantom{A} \begin{itemize} \item[$(1)$] Under what conditions on $d, k_1, k_2$ is any double cone singularity of type $(d; k_1, k_2)$ equivariantly smoothable? \item[$(2)$] When is an admissible anti-$P$-resolution of non-normal type $\mathrm{II}$ $\mathbb{Q}$-Gorenstein smoothable around the germ of the exceptional locus? \end{itemize} \end{ques}
We note that (1) and (2) are equivalent. 
A complete answer to Question \ref{ques--cusp} would significantly enhance the classification of $\Q$-Gorenstein smoothable normal stable Horikawa surfaces by providing a definitive solution to the global $\mathbb{Q}$-Gorenstein smoothing problem.

We remark on the equivariant smoothability of the elliptic double cone singularities.

\begin{rem}\label{rem--smoo}
Although there is no complete criterion for the equivariant cusp smoothability, we give a partial answer of Question~\ref{ques--cusp}~(1) as follows:

\begin{itemize}
\item If $\max\{k_1,k_2\}\ge 2d-8$, then any double cone singularity of type $(d;k_1,k_2)$ is equivariantly smoothable (due to the proof of Proposition~\ref{prop--smoothing--c-type}).
\item If $\max\{k_1,k_2\}< 2d-8$ and $\lceil\frac{k_1}{2}\rceil+\lceil\frac{k_2}{2}\rceil\le d-4$, then any double cone singularity of type $(d;k_1,k_2)$ is not equivariantly smoothable (Theorem~\ref{thm--incomplete}).
\end{itemize}

Note that there exists a smoothable cusp singularity that is not equivariantly smoothable, combining Theorem \ref{thm--incomplete} with observations in Section \ref{app:cusp}.
\end{rem}

\begin{rem}\label{rem--equiv--smooth--detail}
Before discussing the difficulty of equivariant cusp smoothability, we briefly review the known results on the smoothability of cusp singularities.
For a cusp singularity, smoothability is equivalent to the condition that the exceptional divisor of its dual cusp, defined combinatorially, can be realized as an anticanonical divisor on a smooth rational surface.
Looijenga \cite{Loo} studied the deformation theory for Inoue-Hirzebruch surfaces and observed that the former condition implies the latter.
The converse implication, known as the {\em Looijenga conjecture}, was later proved by Gross-Hacking-Keel \cite{GHK} (see also \cite{Engel}).

In contrast, for {\em equivariant} cusp smoothing, such a criterion is not yet known in general, except in a few special cases \cite{Simonetti,J}.
Simonetti \cite{Simonetti} and Jiang \cite{J} studied the equivariant smoothability of cusp singularities under the assumption that the finite group acts freely on the germ outside the singularity. 
They showed that the equivariant smoothability in this setting follows if the exceptional divisor of the dual cusp can be realized as an anticanonical divisor on a rational surface with a compatible finite group action, along with some additional conditions.
Moreover, Jiang proved the converse implication in the special case where the quotient of the cusp singularity is again a cusp, by analyzing Inoue-Hirzebruch surfaces.
However, when the group action is not free outside the cusp singularity, no general criterion for equivariant smoothability is currently known.

Question \ref{ques--cusp} adresses this more subtle case, where the finite group does not act freely on the germ outside the cusp singularity.
\end{rem}

\subsection{Comparison of two moduli spaces}\label{subsec:comparison}
Before stating the result of this subsection, we explain the notion of seminormalization.

\begin{defn}\label{defn--seminormal}
    Let $A$ be a Noetherian reduced ring with its integral closure $\bar{A}$.
    Suppose that $\bar{A}$ is finite over $A$.
    We say that $A^{\mathrm{sn}}$ is the {\it seminormalization} of $A$ if $A^{\mathrm{sn}}$ is the largest among subrings of $\bar{A}$ satisfying the following.
    \begin{enumerate}
        \item $B$ contains $A$ and the induced morphism $\mathrm{Spec}(B)\to \mathrm{Spec}(A)$ is a homeomorphism, and
        \item for any prime ideal $\mathfrak{p}$ of $A$, $(B\otimes_{A}\kappa(\mathfrak{p}))_{\mathrm{red}}\cong \kappa(\mathfrak{p})$, where $\kappa(\mathfrak{p})$ is the residue field of $A$ at $\mathfrak{p}$.
    \end{enumerate}
    We say that $A$ is {\it seminormal} if $A=A^{\mathrm{sn}}$.
    This is also equivalent to that for any $x, y\in A$ such that $x^2=y^3$, there exists a unique element $a\in A$ such that $a^3=x$ and $a^2=y$ (cf.~\cite[Tags 0EUL, 0EUR]{Stacks}).
    We note that $A$ is seminormal if and only if $A_{\mathfrak{m}}$ is seminormal for any maximal ideal $\mathfrak{m}$ of $A$ (cf.~\cite[Corollary 2.7]{GT}).

    For any scheme $X$ of finite type over $\mathbb{C}$, we can always take the {\it seminormalization} $X^{\mathrm{sn}}\to X$ such that for any affine open subset $U=\mathrm{Spec}(A)$ of $X$, $X^{\mathrm{sn}}\times_XU=\mathrm{Spec}((A_{\mathrm{red}})^{\mathrm{sn}})$.
    If $X\cong X^{\mathrm{sn}}$, we say that $X$ is {\it seminormal}.
\end{defn}

It is easy to see that for any $\mathbb{C}$-algebra $A$ of finite type, $A^{\mathrm{sn}}$ is the largest among subrings $B$ of $\bar{A}$ such that
    $B$ contains $A$ and the induced morphism $\mathrm{Spec}(B)\to \mathrm{Spec}(A)$ is a homeomorphism.
    Indeed, suppose that $B$ is a subring of $\bar{A}$ such that
    $B$ contains $A$ and the induced morphism $\mathrm{Spec}(B)\to \mathrm{Spec}(A)$ is a homeomorphism.
    It suffices to show that $B$ is a subring of $A^{\mathrm{sn}}$, that is, for any prime ideal $\mathfrak{p}$, $(B\otimes_{A}\kappa(\mathfrak{p}))_{\mathrm{red}}\cong \kappa(\mathfrak{p})$.
    Now, this is the case since $\mathrm{Spec}((B/\mathfrak{p}B)_{\mathrm{red}})\to \mathrm{Spec}(A/\mathfrak{p})$ is a homeomorphism and hence they are generically isomorphic.

For any morphism $f\colon X\to Y$ of schemes, there exists the unique induced morphism $f^{\mathrm{sn}}\colon X^{\mathrm{sn}}\to Y^{\mathrm{sn}}$ by \cite[Tag 0EUS]{Stacks}.
Furthermore, for any faithfully flat morphism $f\colon X\to Y$, $Y$ is seminormal if so is $X$ (cf.~\cite[Corollary 1.7]{GT}).
Conversely, if $f\colon X\to Y$ is flat with all fibers geometrically reduced, then we see by \cite[Corollary 4.5]{GT} that $X^{\mathrm{sn}}\cong Y^{\mathrm{sn}}\times _YX$.
By using these facts, we can define the seminormalization and seminormality of Deligne-Mumford stacks of finite type over $\mathbb{C}$ in the same way.

    Let $X$ be an affine scheme and $G$ a finite group acting on $X$.
    Then, it is easy to see that $(X/G)^{\mathrm{sn}}\cong X^{\mathrm{sn}}/G$.
    This fact and \cite[Theorem 11.3.1]{Ols} show that for any separated Deligne-Mumford stack $\mathfrak{X}$ of finite type over $\mathbb{C}$, if we let $X$ and $Y$ be the coarse moduli spaces of $\mathfrak{X}$ and $\mathfrak{X}^{\mathrm{sn}}$ respectively, then $Y\cong X^{\mathrm{sn}}$.

Fix a positive integer $p_g\ge 3$.
Let $\mathcal{M}^{\mathrm{Gie}}_{2p_g-4,p_g}$ be the moduli stack of stable Horikawa surfaces with only Du Val singularities and with geometric genus $p_g$. 
Then, the universal family $\mathcal{U}$ has a canonical involution $\iota\colon \mathcal{U}\to\mathcal{U}$ over $\mathcal{M}^{\mathrm{Gie}}_{2p_g-4,p_g}$.
We mention the quotient of $\mathcal{U}$ by $\iota$ as follows.

\begin{lem}\label{lem--universal--involution}
Let $\pi\colon \mathscr{X}\to S$ be a $\mathbb{Q}$-Gorenstein family of stable Horikawa surfaces over a seminormal scheme $S$ of finite type over $\mathbb{C}$.
Suppose that there exists an open dense subset $S^\circ\subset S$ such that $\mathscr{X}_{\bar{s}}$ is standard for any geometric point $\bar{s}\in S$.
Then, $\mathscr{X}$ admits a unique involution
over $S$ the restriction of which to any geometric fiber $\mathscr{X}_{\bar{s}}$ coincides with the good involution.
\end{lem}

To begin with, we state the following useful criterion for properness of morphisms of schemes of finite type over $\mathbb{C}$.

\begin{lem}\label{lem--curve-criterion-for-properness}
    Let $f\colon X\to Y$ be a separated morphism of schemes of finite type over $\mathbb{C}$.
    Then $f$ is proper if and only if $f$ satisfies the following condition.
    \begin{itemize}
        \item[(*)] Let $Y^\circ\subset Y$ be a Zariski open dense subset such that $f^{-1}(Y^\circ)$ is Zariski dense in $X$.
        For any morphism $\varphi\colon C\to Y$ from a smooth affine curve with a closed point $0\in C$ and $\psi^\circ\colon C\setminus\{0\}\to f^{-1}(Y^\circ)$ such that $f\circ \psi^\circ=\varphi|_{C\setminus\{0\}}$, there exists an extension $\psi\colon C\to X$ of $\psi^\circ$.
    \end{itemize}
\end{lem}

\begin{proof}
First note that if $f$ is proper, then (*) holds by \cite[II, Theorem 4.7]{Ha}.
It suffices to deal with the converse direction.
Suppose that (*) holds.
    By Chow's lemma \cite[Tag 0200]{Stacks}, there exists a proper surjective morphism $\pi\colon X'\to X$ such that $X'$ is quasi-projective over $Y$ and there exists an open dense subset $U\subset X$ such that $\pi^{-1}(U)\to U$ is an isomorphism.
    Here, by shrinking $U$, we may assume that $U\subset f^{-1}(Y^\circ)$.
    Let $Z$ be a scheme projective over $Y$ such that there exists a dominant open immersion $\iota\colon X'\hookrightarrow Z$ over $Y$.
    To show that $f$ is proper, it is enough to show that $\iota$ is surjective.
    Then, we see by \cite[\S6, Lemma]{M} that for any closed point $z\in Z$, there exists a morphism $\theta\colon C\to Z$ from a smooth affine curve $C$ with a closed point $0\in C$ such that $\theta(0)=z$ but $(\pi\circ \theta)^{-1}(U)=C\setminus\{0\}$. 
    By (*) and \cite[II, Theorem 4.7]{Ha}, $z\in X'$.
    Thus, $f$ is proper.
\end{proof}

\begin{proof}[Proof of Lemma \ref{lem--universal--involution}]
First, we deal with the case when $S$ is normal and all geometric fibers are standard.
We note that $\mathscr{X}$ is normal in this case.
By \cite[Corollary 2.69]{kollar-moduli}, $\pi_*\omega_{\mathscr{X}/S}$ is a locally free sheaf and for any point $s\in S$, $\pi_*\omega_{\mathscr{X}/S}\otimes_{\mathcal{O}_{\mathscr{X}}}\kappa(s)\cong H^0(\mathscr{X}_{s},\mathcal{O}_{\mathscr{X}_{\bar{s}}}(K_{\mathscr{X}_{s}}))$.
Since $\mathcal{O}_{\mathscr{X}_{\bar{s}}}(K_{\mathscr{X}_{\bar{s}}})$ is globally generated and defines a morphism $\mathscr{X}_{\bar{s}}\to \mathbb{P}^{p_g(\mathscr{X}_{\bar{s}})-1}$ for any geometric point $\bar{s}\in S$, we see that $\pi_*\omega_{\mathscr{X}/S}$ defines $f\colon \mathscr{X}\to \mathbb{P}_S(\pi_*\omega_{\mathscr{X}/S})$.
Let $\mathscr{Y}$ be the scheme-theoretic image of $f$ and let $\mathscr{Y}^\nu$ the normalization of $\mathscr{Y}$.
Let $g\colon \mathscr{X}\to\mathscr{Y}^{\nu}$ be the natural morphism.
Since $f_{s}\colon \mathscr{X}_s\to \mathscr{Y}_s$ is a double cover for any general closed point $s\in S$, $g$ is also a finite morphism of degree two.
Thus, $\mathscr{X}$ is a double cover of $\mathscr{Y}^\nu$ of degree two and has an involution $\sigma$ over $\mathscr{Y}^\nu$ since $\mathscr{X}$ is normal.
By construction, $\sigma$ acts trivially on $\pi_*\omega_{\mathscr{X}/S}$ and hence $\sigma_{\bar{s}}$ is the good involution for any geometric point $\bar{s}\in S$.

Next, we deal with the general case.
By shrinking $S^\circ$, if necessary, we may assume that $S^\circ$ is normal.
Since the KSBA moduli is separated by Theorem \ref{thm-ksba}, we note that $\mathrm{Aut}_{S}(\mathscr{X})$ is a scheme finite over $S$.
By what we have shown in the first paragraph of this proof, there exists a section $\mu^\circ\colon S^\circ\to \mathrm{Aut}_{S}(\mathscr{X})$ corresponding to the involution $\sigma^\circ$ on $\mathscr{X}\times_SS^\circ$ constructed as in the first paragraph such that $\sigma_{\bar{s}}$ coincides with the good involution of $\mathscr{X}_{\bar{s}}$ for any geometric point $\bar{s}\in S^\circ$.
Let $V$ be the closure of $\mu^\circ(S^\circ)$ in $\mathrm{Aut}_{S}(\mathscr{X})$ and let $h\colon V\to S$ be the canonical morphism.
By Proposition \ref{prop--involution}, we obtain the following.
\begin{itemize}
    \item Let $\varphi\colon C\to S$ be a morphism from a smooth affine curve $C$ with a closed point $0\in C$ that maps $C\setminus \{0\}$ to $S^\circ$.
    Then, $h$ induces the following isomorphism $V\times_SC\to C$.
\end{itemize}
Thus, $h$ is proper by Lemma \ref{lem--curve-criterion-for-properness}.
Furthermore, $h$ is injective by Proposition \ref{prop:good_involution_uniqueness}. 
Therefore, $h$ is a homeomorphism and since $S$ is seminormal, $h$ is an isomorphism.
The involution $\sigma$ corresponding to the section $h^{-1}$ is the desired one.
\end{proof}

On the other hand, we explain the following moduli scheme of relative Mumford divisors.

    \begin{thm}[{\cite[Theorem 4.76]{kollar-moduli}}]\label{thm--kollar--mumford-divisor}
Let $f\colon \mathscr{X}\to S$ be a flat projective morphism of relative dimension $d$.
Suppose that $S$ is reduced and $\mathscr{X}\subset \mathbb{P}_S^{n}$ for some $n\in\mathbb{Z}_{>0}$.
Then, there exists a scheme $\mathbf{MDiv}_{\mathscr{X}/S}$ that is separated over $S$ such that the set of $S$-morphisms $T\to \mathbf{MDiv}_{\mathscr{X}/S}$ is functorially bijective to the set of relative Mumford divisors on $\mathscr{X}_T$ over $T$.
    \end{thm}

Applying Theorem \ref{thm--kollar--mumford-divisor}, we obtain the following.
\begin{lem}\label{lem--family-wise--quotient}
    Let $d\in\mathbb{Z}_{\ge0}$ and $f\colon \mathscr{X}\to S$ be a $\mathbb{Q}$-Gorenstein family over a seminormal scheme $S$ of finite type over $\mathbb{C}$.
    Suppose that $f$ is projective, every geometric fiber is lc, irreducible and of dimension $d$. 
    Suppose further that there exists an involution $\sigma$ on $\mathscr{X}$ over $S$.
    Then, there exist a flat projective morphism $g\colon \mathscr{Y}\to S$, $\pi\colon\mathscr{X}\to\mathscr{Y}$, which is the quotient morphism by $\sigma$ and a unique relative Mumford divisor $\mathscr{B}$ on $\mathscr{Y}$ such that $(\mathscr{Y},\frac{1}{2}\mathscr{B})$ is a log $\mathbb{Q}$-Gorenstein family over $S$ and $K_{\mathscr{X}/S}=\pi^*(K_{\mathscr{Y}/S}+\frac{1}{2}\mathscr{B})$. 
\end{lem}

\begin{proof}
    Since $f$ is projective, we see that there exists a quotient morphism $\pi\colon\mathscr{X}\to\mathscr{Y}$ with respect to $\sigma$ and $\mathscr{Y}$ is a projective scheme over $S$ such that any geometric fiber is normal and of dimension $d$.
    Let $g\colon\mathscr{Y}\to S$ be the structure morphism.
    It follows that $\mathscr{Y}$ is projective and flat over $S$ from the fact that $\mathcal{O}_{\mathscr{Y}}$ is a direct summand of $\pi_*(\mathcal{O}_{\mathscr{X}})$.
    We note that for any geometric point $\bar{s}\in S$, $\mathscr{Y}_{\bar{s}}$ is the quotient $\mathscr{X}_{\bar{s}}/\sigma_{\bar{s}}$ and hence normal.
    Thus, we can take an open subset $U\subset \mathscr{Y}$ such that $g|_{U}$ is smooth and $\mathrm{codim}_{\mathscr{Y}_{s}}(\mathscr{Y}_s\setminus U)\ge2$ for any $s\in S$.
    If $S$ is smooth, then let $D$ be the branch divisor of $\pi|_{\pi^{-1}(U)}$.
    Since $D$ contains no fiber of $g|_U$, we see that $D$ is a relative Cartier divisor (cf.~\cite[Definition-Lemma 4.20]{kollar-moduli}).
    Furthermore, it is easy to see that $K_{\pi^{-1}(U)/S}=\pi^*(K_{U/S}+\frac{1}{2}D)$.
    We note that $\mathscr{Y}$ and $\mathscr{X}$ are normal in this case.
    Therefore, it is easy to see that the closure $\mathscr{B}$ of $D$ in $\mathscr{Y}$ is the desired relative Mumford divisor. 

    Next, we deal with the case when $S$ is seminormal.
    By shrinking $S$ if necessary, we may assume that there exists $n\in\mathbb{Z}_{>0}$ such that $\mathscr{Y}\subset \mathbb{P}^n_S$.
    Take a Zariski open dense subset $S^{\mathrm{reg}}\subset S$ such that $S^{\mathrm{reg}}$ is smooth.
    Then, by what we have shown in the previous paragraph, there exists a relative Mumford divisor $\mathscr{B}^{\mathrm{reg}}$ on $\mathscr{Y}\times_SS^{\mathrm{reg}}$ as the assertion.
    $\mathscr{B}^{\mathrm{reg}}$ corresponds to a section $\mu^\circ\colon S^{\mathrm{reg}}\to \mathbf{MDiv}_{\mathscr{Y}/S}$.
    Let $V$ be the closure of the image of $\mu^\circ$ in $\mathbf{MDiv}_{\mathscr{Y}/S}$.

    Now, we claim that $V$ is homeomorphic to $S$.
    To show this, we have to show that the canonical morphism $h\colon V\to S$ is proper and bijective
    .
    For any closed point $s\in S$, we can take a morphism $\varphi\colon C\to S$ from a smooth affine curve $C$ with a closed point $0\in C$ such that $\varphi(0)=s$ and $\varphi(C\setminus\{0\})\subset S^{\mathrm{reg}}$ by \cite[\S6, Lemma]{M}.
    For such $C$, there exists a relative Mumford divisor $\mathscr{B}_C$ on $\mathscr{Y}_C$ that is the branch divisor of $\mathscr{X}_C\to \mathscr{Y}_C$ by what we have shown in the first paragraph of this proof.
    Then, $h|_{V\times_SC}\colon V\times_SC\to C$ has a section induced by $\mathscr{B}_C$ and is surjective.
We note that $h$ is injective.
Indeed, for any closed point $v\in V$, there exists a morphism $\psi\colon C\to V$ from a smooth curve with a closed point $0\in C$ such that $\psi(C\setminus\{0\})\subset \mu^\circ(S^{\mathrm{reg}})$ and $\psi(0)=v$.
If we let $\varphi=h\circ\psi$ and $s=\psi(v)$, then there exists a relative Mumford divisor $\mathscr{B}_C$ on $\mathscr{Y}_C$ that is the branch divisor of $\mathscr{X}_C\to \mathscr{Y}_C$ such that $(\mathscr{B}_C)_s$ corresponds to $v\in \mathbf{MDiv}_{\mathscr{Y}_s/\mathrm{Spec}(\C)}$.
Here, $(\mathscr{B}_C)_s$ is defined as the branch divisor of $\mathscr{X}_s\to \mathscr{Y}_s$ and hence $h^{-1}(s)=\{v\}$. 
 Therefore, $h$ is injective and we can see by Lemma \ref{lem--curve-criterion-for-properness} that $h$ is proper. 
    Since $h$ is proper, we also have that $h$ is surjective and hence a homeomorphism.
    Thus, it follows that $h$ is an isomorphism from the assumption that $S$ is seminormal.
Take the corresponding relative Mumford divisor $\mathscr{B}$ on $\mathscr{Y}$ to $h^{-1}$.
Note that by construction, $\mathscr{B}_{\bar{s}}$ is the branch divisor with respect to $\pi_{\bar{s}}$ for any geometric point $\bar{s}\in S$.
By what we have shown in the first paragraph and \cite[Corollary 4.35]{kollar-moduli}, we see that $(\mathscr{Y},\frac{1}{2}\mathscr{B})$ is log $\mathbb{Q}$-Gorenstein.
It is not hard to check that $K_{\mathscr{X}/S}=\pi^*(K_{\mathscr{Y}/S}+\frac{1}{2}\mathscr{B})$.
Thus, we complete the proof.
\end{proof}
Take the seminormalizations $(\mathcal{M}^{\mathrm{Gie}}_{2p_g-4,p_g})^{\mathrm{sn}}$ of $\mathcal{M}^{\mathrm{Gie}}_{2p_g-4,p_g}$ and $\mathcal{U}^{\mathrm{sn}}$ of $\mathcal{U}$ respectively.
Note that $$\mathcal{U}\times_{\mathcal{M}^{\mathrm{Gie}}_{2p_g-4,p_g}}(\mathcal{M}^{\mathrm{Gie}}_{2p_g-4,p_g})^{\mathrm{sn}}=\mathcal{U}^{\mathrm{sn}}.$$
Let $\sigma^{\mathrm{sn}}$ be the induced involution on $\mathcal{U}^{\mathrm{sn}}$ (cf.~Lemma \ref{lem--universal--involution}).
By Lemma \ref{lem--family-wise--quotient}, we can take the quotient $\pi\colon \mathcal{U}^{\mathrm{sn}}\to\mathcal{Y}$ by $\sigma^{\mathrm{sn}}$ and the unique relative Mumford divisor $\mathcal{D}$ such that $(\mathcal{Y},\frac{1}{2}\mathcal{D})$ is log $\Q$-Gorenstein family over $(\mathcal{M}^{\mathrm{Gie}}_{2p_g-4,p_g})^{\mathrm{sn}}$ such that $$K_{\mathcal{U}^{\mathrm{sn}}/(\mathcal{M}^{\mathrm{Gie}}_{2p_g-4,p_g})^{\mathrm{sn}}}=\pi^*\left(K_{\mathcal{Y}/(\mathcal{M}^{\mathrm{Gie}}_{2p_g-4,p_g})^{\mathrm{sn}}}+\frac{1}{2}\mathcal{D}\right).$$
$(\mathcal{Y},\frac{1}{2}\mathcal{D})$ induces the following natural morphism 
\[
\alpha\colon (\mathcal{M}^{\mathrm{Gie}}_{2p_g-4,p_g})^{\mathrm{sn}}\to \mathcal{M}^{\mathrm{KSBA}}_{2,2p_g-4,\frac{1}{2}}.
\]
Let $\mathcal{M}^{\mathrm{quot}}_{p_g}$ be the scheme-theoretic image of $\alpha$ and take the open locus $\mathcal{M}^{\mathrm{nq}}_{p_g}$ of $\mathcal{M}^{\mathrm{quot}}_{p_g}$ such that if we let $(W,\frac{1}{2}B)$ denote the log pair corresponding to a geometric point $\bar{s}\in \mathcal{M}_{p_g}^{\mathrm{quot}}$, $W$ is normal and $B$ is reduced if and only if $\bar{s}\in \mathcal{M}^{\mathrm{nq}}_{p_g}$.

\begin{thm}\label{prop--moduli--seminormaliazation}
    Let $(\mathcal{M}^{\mathrm{nq}}_{p_g})^{\mathrm{sn}}$ and $(\overline{\mathcal{M}}^{\mathrm{Gie}}_{2p_g-4,p_g})^{\mathrm{sn}}$ be the seminormalizations of $\mathcal{M}^{\mathrm{nq}}_{p_g}$ and of the Zariski closure  $\overline{\mathcal{M}}^{\mathrm{Gie}}_{2p_g-4,p_g}$ of $\mathcal{M}^{\mathrm{Gie}}_{2p_g-4,p_g}$ in $\mathcal{M}^{\mathrm{KSBA,normal}}_{2,2p_g-4,0}$,  respectively.

    Then there exists a proper quasi-finite surjective morphism 
    \[
    \overline{\alpha}\colon (\overline{\mathcal{M}}^{\mathrm{Gie}}_{2p_g-4,p_g})^{\mathrm{sn}} \to (\mathcal{M}^{\mathrm{nq}}_{p_g})^{\mathrm{sn}}
    \]
    that is an extension of the induced morphism $\alpha^{\mathrm{sn}}\colon (\mathcal{M}^{\mathrm{Gie}}_{2p_g-4,p_g})^{\mathrm{sn}} \to (\mathcal{M}^{\mathrm{nq}}_{p_g})^{\mathrm{sn}}$ by $\alpha$.
    Let $\overline{\alpha}'\colon (\overline{M}^{\mathrm{Gie}}_{2p_g-4,p_g})^{\mathrm{sn}} \to (M^{\mathrm{nq}}_{p_g})^{\mathrm{sn}}$ be the morphism of the coarse moduli spaces induced by $\overline{\alpha}$.
    Then $\overline{\alpha}'$ satisfies the following.
    \begin{itemize} 
   \item[$(1)$] If we further let $(M^{\mathrm{rat}}_{p_g})^{\mathrm{sn}}$ be the open locus of $(M^{\mathrm{nq}}_{p_g})^{\mathrm{sn}}$ parametrizing $(Y,\frac{1}{2}B)$ such that $Y$ has only rational singularities, then $\overline{\alpha}'|_{\overline{\alpha}'^{-1}((M^{\mathrm{rat}}_{p_g})^{\mathrm{sn}})}$ is an isomorphism, and
   \item[$(2)$] If $p_g=4$, the cardinality of $\overline{\alpha}'^{-1}(s)$ is four for any $s\in (M^{\mathrm{nq}}_{4})^{\mathrm{sn}}\setminus(M^{\mathrm{rat}}_{4})^{\mathrm{sn}}$. In particular, $\overline{\alpha}'$ is a finite morphism.
    \end{itemize}
\end{thm}

\begin{proof}
By applying Lemma \ref{lem--universal--involution}, we see that if $\mathcal{U}\to (\overline{\mathcal{M}}^{\mathrm{Gie}}_{2p_g-4,p_g})^{\mathrm{sn}}$ is the universal $\Q$-Gorenstein family of stable Horikawa surfaces, then $\mathcal{U}$ admits an involution $\sigma$ such that for any geometric point $\bar{s}\in (\overline{\mathcal{M}}^{\mathrm{Gie}}_{2p_g-4,p_g})^{\mathrm{sn}}$, the restriction of $\sigma$ to $\mathcal{U}_{\bar{s}}$ is the good involution.
Then, we construct $\overline{\alpha}$ by using Lemma \ref{lem--family-wise--quotient} in the same way as $\alpha$.

    On the other hand, it is easy to check that $\overline{\alpha}$ is proper by Lemma \ref{lem--curve-criterion-for-properness} and Corollary \ref{cor--too--trivial}.
    Therefore, $\overline{\alpha}$ is a proper surjective morphism between seminormal Deligne-Mumford stacks.
    
    Next, we show that the restriction $\overline{\alpha}'|_{\overline{\alpha}'^{-1}((M^{\mathrm{rat}}_{p_g})^{\mathrm{sn}})}$ is injective.
    In other words, we show the following:
    Lett $X$ be a $\mathbb{Q}$-Gorenstein smoothable normal stable Horikawa surface such that the quotient $W$ has only rational singularities.
    Then for any two divisors $L,L'\in\mathrm{Cl}(W)$ with $2L\sim 2L'$ and any effective divisor $B\in |2L|=|2L'|$, if the two double covers corresponding to the branch data $B\in |2L|$ and $B'\in |2L'|$ are both  $\mathbb{Q}$-Gorenstein smoothable stable normal Horikawa surfaces, then we must have $L\sim L'$.
    To prove this, it suffices to show that if $W$ is a klt rational surface that arises as the quotient of a normal stable Horikawa surface by its good involution, then $\mathrm{Cl}(W)$ is torsion-free.
    
    If $W$ is a smooth rational surface, then $\mathrm{Pic}(W)$ is torsion-free, and the statement is trivial.
    In general, let $\beta\colon \tilde{W}\to W$ be a resolution of singularities, and suppose $D$ is a Weil divisor on $W$ such that  $nD\sim 0$ for some $n\in\mathbb{Z}_{>0}$.
    Then $\beta^*D$ is a Cartier divisor such that $n\beta^*D\sim 0$.
    This can be verified using our classification.
    For instance, if $W$ is the quotient of a $\mathbb{Q}$-Gorenstein smoothable normal stable Horikawa surface of special Lee-Park type, then $$
    \widetilde{\Gamma}\cdot(\lceil \beta^*D\rceil-\beta^*D)\in\mathbb{Z},
    $$
    where $\widetilde{\Gamma}$ is the strict transform of a general fiber $\Gamma$.
    This implies that the proper transform $\tilde{\Delta}_{0}$ of $\Delta_0\subset \Sigma_{p_g-2}$ is not contained in the support of $\lceil \beta^*D\rceil-\beta^*D$.
    In particular, since $\tilde{\Delta}_0\cdot (\lceil \beta^*D\rceil-\beta^*D)\in \mathbb{Z}$, no curve intersecting $\tilde{\Delta}_0$ lies in the support of $\lceil \beta^*D\rceil-\beta^*D$.
    By iterating such arguments, we conclude that $\mathrm{Supp}(\lceil \beta^*D\rceil-\beta^*D)$ does not contain any irreducible component of the T-chain and the $-(p_g-2)$-curve.
    Hence, $\lceil \beta^*D\rceil=\beta^*D$, and $\beta^{*}D$ is integral.
    For the other cases in our classification, we similarly see that $\beta^*D$ is a Cartier divisor. 
    Since $\widetilde{W}$ is smooth and rational, $\mathrm{Pic}(\tilde{W})$ is torsion-free.
    Thus, $n\beta^{*}D\sim 0$ implies that $\beta^{*}D\sim 0$.
    Hence, we conclude that $D=\beta_{*}\beta^{*}D\sim 0$.

Finally, we show that for any closed point $s\in (\mathcal{M}^{\mathrm{nq}}_{4})^{\mathrm{sn}}\setminus(\mathcal{M}^{\mathrm{rat}}_{4})^{\mathrm{sn}}$, there are only four objects of $\overline{\alpha}^{-1}(s)$.
We first note that $s$ corresponds to a log pair $(W,\frac{1}{2}B)$ such that $W$ is an elliptic cone of degree eight by the classification of $\mathbb{Q}$-Gorenstein smoothable stable Horikiwa surfaces.
If a $\mathbb{Q}$-Gorenstein smoothable stable normal Horikawa surface $X$ corresponds to a point of $\overline{\alpha}^{-1}(s)$, then there exists a Weil divisor $L'$ on $W$ such that $B\in |2L'|$ and $X$ is the corresponding double cover.
Let $C\subset \mathbb{P}^8$ correspond to the elliptic cone $W$.
It is easy to see that the linear equivalence class of $L'$ such that $B\in |2L'|$ corresponds to the linear equivalence class of a line bundle $\mathcal{L}$ on the corresponding elliptic curve $C$ of degree four such that $|2\mathcal{L}|$ is corresponding to the embedding $C\subset \mathbb{P}^8$ bijectively.
We have already shown that all $X$ for such $L'$ are $\mathbb{Q}$-Gorenstein smoothable in Proposition \ref{prop--canonical--pencil--invol}.
Thus, we obtain the assertion.
\end{proof}

\begin{rem}
 We remark that the restriction $\overline{\alpha}|_{\overline{\alpha}^{-1}((\mathcal{M}^{\mathrm{rat}}_{p_g})^{\mathrm{sn}})}$ is not an isomorphism of stacks onto its image.  
 Indeed, take any point $(W,\frac{1}{2}B)\in (\mathcal{M}^{\mathrm{rat}}_{p_g})^{\mathrm{sn}}(\mathbb{C})$ and consider the double cover $X$ of $W$ branched along $B$.
 Assume that $X$ is smooth and of type $(0)$, $(1)$ or $(\infty)$, and that it corresponds to a general point $x$ of $\mathcal{M}^{\mathrm{Gie}}_{2p_g-4,p_g}$.
 According to the proof of \cite[Theorem 2.1]{horikawa}, the semiversal deformation space of $X$ is smooth.
Therefore, the Gieseker moduli stack $\mathcal{M}^{\mathrm{Gie}}_{2p_g-4,p_g}$ is smooth at $x$, and locally around $x$ it agrees with its seminormalization.
 The automorphism group $\mathrm{Aut}(X)$ contains the good involution $\sigma$. 
 Let $G$ be the stabilizer group of the seminormalization $(\overline{\mathcal{M}}^{\mathrm{Gie}}_{2p_g-4,p_g})^{\mathrm{sn}}$ at $x$, and $H$ the stabilizer group of $(\mathcal{M}^{\mathrm{rat}}_{p_g})^{\mathrm{sn}}$ at the point corresponding to $(W,\frac{1}{2}B)$.
 Then $G$ is isomorphic to $\mathrm{Aut}(X)$, and the natural map $G\to H$ induced by $\overline{\alpha}$ is not injective since $\sigma$ is contained in the kernel.
This implies that $\overline{\alpha}$ cannot be an isomorphism.
\end{rem}

We obtain the following immediate corollary, which reduces the problem of stratification of the moduli of $\Q$-Gorenstein smoothable stable normal Horikawa surfaces to stratification of the moduli of the pairs associated to Horikawa surfaces.

\begin{cor}\label{cor--reduction--to--quot}
Let $X$ be a $\mathbb{Q}$-Gorenstein smoothable stable normal Horikawa surface.
Let $\mathscr{F}$ be a locally closed substack of $(\overline{\mathcal{M}}^{\mathrm{Gie}}_{2p_g-4,p_g})^{\mathrm{sn}}$.
Then, the geometric point $\bar{x}\in(\overline{\mathcal{M}}^{\mathrm{Gie}}_{2p_g-4,p_g})^{\mathrm{sn}}$ corresponding to $X$ belongs to the Zariski closure of $(\overline{\mathcal{M}}^{\mathrm{Gie}}_{2p_g-4,p_g})^{\mathrm{sn}}$ if and only if $\overline{\alpha}(\bar{x})$ belongs to the Zariski closure of $\overline{\alpha}(\mathscr{F})$.
\end{cor}

\subsection{Stratification of the moduli space}\label{subsec:stratification}
Let $(\mathcal{M}^{\mathrm{rat}}_{p_g})^{\mathrm{sn}}$ be the open locus of $(\mathcal{M}^{\mathrm{nq}}_{p_g})^{\mathrm{sn}}$ parametrizing $(W,\frac{1}{2}B)$ such that $W$ has only rational singularities.
We put the following lemma, which claims under some condition, partial $\mathbb{Q}$-Gorenstein smoothings of $\mathbb{Q}$-Gorenstein smoothable Horikawa surfaces are again $\mathbb{Q}$-Gorenstein smoothable.

\begin{lem}\label{lem--part--smoothing--again--smoothable}
Let $f\colon \mathscr{W}\to C$ be a projective flat morphism from a $\Q$-Gorenstein normal threefold to a smooth curve $C$ with a closed point $0\in C$.
Let $\mathscr{L}$ be a $\Q$-Cartier divisorial sheaf on $\mathscr{W}$, and let $\mathscr{B}$ be an effective relative Mumford divisor $\mathscr{B}$ such that $\mathscr{L}^{[2]}\sim_C\mathscr{B}$.
Assume the following conditions:
\begin{itemize}  
    \item 
    $\mathscr{W}_0$ is slc and $\mathscr{W}_c$ has only T-singularities for any $c\in C\setminus\{0\}$.
    \item 
    $\mathscr{B}_0$ and $-K_{\mathscr{W}_0}$ are log big and nef.
    \item 
    For each singularity of $\mathscr{W}_0$, there exists an integer $l$ such that $\mathscr{L}_{0}\sim lK_{\mathscr{W}_0}$ locally around the singularity.
    \item 
    $(\mathscr{W}_0,\frac{1}{2}\mathscr{B}_0)$ is an admissible anti-P-resolution of a pair $(W,\frac{1}{2}B)\in\mathcal{M}^{\mathrm{rat}}_{p_g}(\mathbb{C})$
    \item 
    $(\mathscr{W}_c,\frac{1}{2}\mathscr{B}_c)$ is stable for any $c\in C\setminus\{0\}$.
\end{itemize}
Then we have $(\mathscr{W}_c,\frac{1}{2}\mathscr{B}_c)\in \mathcal{M}^{\mathrm{rat}}_{p_g}(\mathbb{C})$ for general $c\in C$.
\end{lem}

\begin{proof}
Let $w\in \mathscr{W}_{0}$ be a singularity.
Suppose that $\mathscr{L}_{0}\sim lK_{\mathscr{W}_0}$ locally around $w$ for some $l$.
Then, by applying \cite[Proposition 2.17]{kollar-moduli} to the sheaf $\mathscr{L}[\otimes]\omega_{\mathscr{W}/C}^{[-l+1]}$, we obtain that 
\[
\mathscr{L}\sim_C\omega_{\mathscr{W}/C}^{[l]}
\]
in a neighborhood of $w\in \mathscr{W}$. 
This and the properness of $f$ imply that for any closed point $w'\in \mathscr{W}_c$ lying over a general point $c\in C$, there exists an integer $l$ such that $\mathscr{L}_c\sim lK_{\mathscr{W}_c}$ locally around $w'$.
Furthermore, since 
\[
H^1(\mathscr{W}_0,\mathcal{O}_{\mathscr{W}_0}(\mathscr{B}_0))=H^2(\mathscr{W}_0,\mathcal{O}_{\mathscr{W}_0})=H^1(\mathscr{W}_0,\mathcal{O}_{\mathscr{W}_0})=H^1(\mathscr{W}_0,\mathcal{O}_{\mathscr{W}_0}(l(2K_{\mathscr{W}_0}+\mathscr{B}_0)))=0
\]
for any sufficiently large and divisible $l\in\mathbb{Z}_{>0}$ by \cite[Theorem 1.10]{fujino--slc--vanishing},
it follows from \cite[III, Theorem 12.8]{Ha} that the same vanishing holds for general fibers: 
\[
H^1(\mathscr{W}_c,\mathcal{O}_{\mathscr{W}_c}(\mathscr{B}_c))=H^2(\mathscr{W}_c,\mathcal{O}_{\mathscr{W}_c})=H^1(\mathscr{W}_c,\mathcal{O}_{\mathscr{W}_c})=H^1(\mathscr{W}_c,\mathcal{O}_{\mathscr{W}_c}(l(2K_{\mathscr{W}_c}+\mathscr{B}_c)))=0.
\]
Additionally, by \cite[Theorem 1.10]{fujino--slc--vanishing}, we have 
$$
H^1(\mathscr{W}_0,\mathcal{O}_{\mathscr{W}_0}(-lK_{\mathscr{W}_0}))=0
$$
for any sufficiently divisible and large integer $l$.
Therefore, using \cite[Theorem 1.15]{fujino--slc--vanishing}, we conclude that $-K_{\mathscr{W}_c}$ is semiample and big for general $c\in C$.
Thus, it follows from Proposition \ref{prop--smoothing--quotient} that for each general $c\in C$, there exists a projective flat morphism $f'\colon \mathscr{W}'\to C'$ from a $\Q$-Gorenstein normal threefold to a smooth curve $C'$ with a closed point $0'\in C'$, together with a $\Q$-Cartier divisorial sheaf $\mathscr{L}'$ and an effective relative Mumford divisor $\mathscr{B}'$, such that 
\begin{itemize}
    \item 
    $\mathscr{L}'^{[2]}\sim_{C'}\mathscr{B}'$.
    \item 
    $(\mathscr{W}'_0,\mathscr{B}'_0)\cong (\mathscr{W}_c,\mathscr{B}_c)$.
    \item 
    $\mathscr{W}'_{c'}$ is smooth for any general $c'\in C'$.   
\end{itemize}
Since $-K_{\mathscr{W}'_{c'}}$ is big and nef and satisfies $K_{\mathscr{W}'_{c'}}^2=8$ (resp.~$K_{\mathscr{W}'_{c'}}^2=9$), the classification of smooth weak del Pezzo surfaces \cite{D} implies that $\mathscr{W}'_{c'}\cong \Sigma_0,\Sigma_1$ or $\Sigma_2$ (resp.~$\mathscr{W}'_{c'}\cong \mathbb{P}^2$).
Thus, we may assume that $\mathscr{B}'_{c'}$ is smooth.
Taking a double cover $Y\to \mathscr{W}'_{c'}$ branched along $\mathscr{B}'_{c'}$, we obtain a smooth stable surface $Y$ with $K_Y^2=2p_g(Y)-4$ by \cite[Corollary 2.69]{kollar-moduli}.
Thus, $Y$ is a smooth Horikawa surface.
Furthermore, by Remark \ref{rem--smoothing--quotient} (3), the involution on $Y$ induced by the double cover acts trivially on $H^0(Y,K_Y)$.
Therefore, the pair $(\mathscr{W}'_{c'},\frac{1}{2}\mathscr{B}'_{c'})$ is the quotient of $Y$ by a good involution.
This completes the proof of the claim.
\end{proof}

In this subsection, we explain the concrete description of the moduli space parametrizing $\mathbb{Q}$-Gorenstein smoothable stable Horikawa surfaces for fixed geometric genus.
By Theorem \ref{prop--moduli--seminormaliazation}, we focus on $(W,\frac{1}{2}B)$ rather than $X$.
We collect the following notions.

\begin{note}\label{note--stratification}
    \begin{enumerate}
        \item  Let $\mathscr{F}^3_{2[4]}$ be the locally closed substack of $(\mathcal{M}^{\mathrm{nq}}_{3})^{\mathrm{sn}}$ parameterizing pairs associated to stable Horikawa surfaces with two singularities of type $\frac{1}{4}(1,1)$ with geometric genus $3$.
        
        \item  Let $\mathscr{F}^3_{(2,2,2,2)}$ be the locally closed substack of $(\mathcal{M}^{\mathrm{nq}}_{3})^{\mathrm{sn}}$ parameterizing pairs associated to stable Horikawa surfaces with one singularity of type $(2,2,2,2)$ with geometric genus $3$ (see Lemma-Definition \ref{lem:smoothable_rational_strict_lc}).
        
        \item Let $\mathscr{F}^3_{[2433]}$ be the locally closed substack of $(\mathcal{M}^{\mathrm{nq}}_{3})^{\mathrm{sn}}$ parameterizing pairs associated to stable Horikawa surfaces with one singularity of type $\frac{1}{50}(1,29)$ with geometric genus $3$.
        
        \item Let $\mathscr{F}^3_{2[4],[2433]}$ be the locally closed substack of $(\mathcal{M}^{\mathrm{nq}}_{3})^{\mathrm{sn}}$ parameterizing pairs associated to stable Horikawa surfaces with two singularities of type $\frac{1}{4}(1,1)$ and one singularity of type $\frac{1}{50}(1,29)$ with geometric genus $3$.
        
        \item Let $\mathscr{F}^3_{(2,2,2,2),[2433]}$ be the locally closed substack of $(\mathcal{M}^{\mathrm{nq}}_{3})^{\mathrm{sn}}$ parametrizing pairs associated to stable Horikawa surfaces with one strict lc rational singularity of type $(2,2,2,2)$ and one singularity of type $\frac{1}{50}(1,29)$ with geometric genus $3$.
        
        \item Let $p_g\in\mathbb{Z}_{\ge3}$. Let $\mathscr{F}^{p_g}_{\mathrm{dV},(k)}$ for $k\in\mathbb{Z}_{\ge0}\cup\{\infty\}$ be the locally closed substack of $(\mathcal{M}^{\mathrm{nq}}_{p_g})^{\mathrm{sn}}$ parameterizing pairs associated to stable Horikawa surfaces with canonical singularities partially smoothable to smooth Horikawa surfaces of type $(k)$ with geometric genus $p_g$.
        
        \item Let $\mathscr{F}^4_{\mathrm{cone}}$ be the locally closed substack of $(\mathcal{M}^{\mathrm{nq}}_{4})^{\mathrm{sn}}$ parameterizing pairs associated to stable Gorenstein Horikawa surfaces with canonical pencil.

        \item Let $p_g\in\mathbb{Z}_{\ge4}$. Let $\mathscr{F}^{p_g}_{\mathrm{LPI}}$ (resp.~$\mathscr{F}^{p_g}_{\mathrm{LPI}^*}$) be the locally closed substack of $(\mathcal{M}^{\mathrm{nq}}_{p_g})^{\mathrm{sn}}$ parameterizing pairs associated to stable Horikawa surfaces of general Lee-Park type I (resp.~I$^*$) with geometric genus $p_g$.

        \item Let $\mathscr{F}^6_{\mathrm{LP}\infty}$ be the locally substack of $(\mathcal{M}^{\mathrm{nq}}_{6})^{\mathrm{sn}}$ parameterizing pairs associated to stable Horikawa surfaces of infinite Lee-Park type with geometric genus six.
        
        \item Let $p_g$ be a positive integer such that $4\le p_g\le 6$. Let $\mathscr{F}^{p_g,m}_{\mathrm{std}}$ be the locally closed substack of $(\mathcal{M}^{\mathrm{nq}}_{p_g})^{\mathrm{sn}}$ parameterizing pairs associated to standard Horikawa surfaces with a double cone singularity with geometric genus $p_g$ and multiplicity $m$.
        
        \item Let $p_g$ be a positive integer such that $4\le p_g\le 6$. Let $\mathscr{F}^{p_g,m}_{\mathrm{LPS}}$ be the locally closed substack of $(\mathcal{M}^{\mathrm{nq}}_{p_g})^{\mathrm{sn}}$ parameterizing pairs associated to stable Horikawa surfaces of special Lee-Park type with a double cone singularity with geometric genus $p_g$ and multiplicity $m$.

\item Let $p_g\in\mathbb{Z}_{\ge7}$. Let $\mathscr{F}^{p_g}_{\mathrm{stdF}}$ be the locally closed substack of $(\mathcal{M}^{\mathrm{nq}}_{p_g})^{\mathrm{sn}}$ parametrizing pairs associated to standard Horikawa surfaces of double Fano type with geometric genus $p_g$.
        
         \item Let $p_g\in\mathbb{Z}_{\ge7}$. Assume that $p_g\ne10$. Let $\mathscr{F}^{p_g}_{\mathrm{\mathrm{stdnF}}}$ be the locally closed substack of $(\mathcal{M}^{\mathrm{nq}}_{p_g})^{\mathrm{sn}}$ parametrizing pairs associated to standard Horikawa surfaces of double non-Fano type and of supersingular type with geometric genus $p_g$. If $p_g=10$, let $\mathscr{F}^{10}_{\mathrm{stdnF}}$ be the locally closed substack of $(\mathcal{M}^{\mathrm{nq}}_{10})^{\mathrm{sn}}$ parameterizing pairs associated to standard Horikawa surfaces of non-Fano type with geometric genus ten.
         
        \item Let $p_g\in\mathbb{Z}_{\ge7}$. Let $\mathscr{F}^{p_g}_{\mathrm{LPS}}$ (resp.~$\mathscr{F}^{p_g}_{\mathrm{LPI}^*}$) be the locally closed substack of $(\mathcal{M}^{\mathrm{nq}}_{p_g})^{\mathrm{sn}}$ parameterizing pairs associated to stable Horikawa surfaces of special Lee-Park type with geometric genus $p_g$.

        \item Set $\overline{\alpha}$ as in Theorem \ref{prop--moduli--seminormaliazation}.
        We set the reduced structure of the locally closed substacks in $(\overline{\mathcal{M}}^{\mathrm{Gie}}_{2p_g-4,p_g})^{\mathrm{sn}}$ corresponding to the above locally closed substacks as follows.
        For example, we set $\mathscr{G}^{p_g}_{\mathrm{stdnF}}$ as the reduced structure of $\overline{\alpha}^{-1}(\mathscr{F}^{p_g}_{\mathrm{stdnF}})$.
        We set the other locally closed substacks in the same way. 
    \end{enumerate}
\end{note}

\begin{rem}\label{rem--existence-all--q-gor--smoothable}
As in Remark \ref{rem-existence--LPI}, combining Lemmas \ref{lem--existence--smoothable--cusp}, \ref{lem--vanishing--for--standard}, \ref{lem--existence--low--special}, \ref{lem--vanishing--for--special}, Propositions \ref{prop--canonical--pencil--invol}, \ref{prop--construction--lee--park--i*}, \ref{prop--smoothing--c-type} and similar arguments to Remark \ref{rem--existence--toric}, it is easy to see that all the above substacks are non-empty. 
Furthermore, it is also easy to see that all substacks are irreducible except $\mathscr{F}^{p_g}_{{\mathrm{stdF}}}$, $\mathscr{F}^{p_g}_{{\mathrm{stdnF}}}$ or $\mathscr{F}^{p_g}_{\mathrm{LPS}}$ for $p_g\ge7$ by a similar argument to Remark \ref{rem-existence--LPI}.
\end{rem}

We state the following general statements to see the stratification of our moduli stack $(\mathcal{M}^{\mathrm{nq}}_{p_g})^{\mathrm{sn}}$.

We recall the result of Horikawa as follows.
Note that $M^{\mathrm{Gie}}_{2p_g-4,p_g}$ is canonically homeomorphic to the coarse moduli space of the open substack $\bigcup_{k\in\mathbb{Z}_{\ge0}\cup\{\infty\}}\mathscr{F}^{p_g}_{\mathrm{dV},(k)}\cup \bigcup_{m<2}\mathscr{F}^{p_g,m}_{\mathrm{std}}$ of $(\mathcal{M}^{\mathrm{rat}}_{p_g})^{\mathrm{sn}}$ by Theorem \ref{prop--moduli--seminormaliazation}.

\begin{thm}[{\cite{horikawa}}]\label{thm--horikawa--smooth--stratification}
Let $p_g$ be a positive integer.
Then, all the members of $M^{\mathrm{Gie}}_{2p_g-4,p_g}$ are smoothable.

If $p_g$ is not divisible by $8$, then $M^{\mathrm{Gie}}_{2p_g-4,p_g}$ is irreducible.
If $p_g$ is divisible by $8$ and $p_g>8$ (resp.~$p_g=8$), $M^{\mathrm{Gie}}_{2p_g-4,p_g}$ has two connected components.
In these cases, one component consists of stable Horikawa surfaces with only Du Val singularities of Type $(0), (2),\ldots,(\frac{p_g}{2}-1)$ (resp.~of Type $(0)$ and $(2)$), while the other component consists of Type $(\frac{p_g}{2}+1)$ (resp.~of Type $(\infty)$ and $(4')$).
\end{thm}

Theorem \ref{thm--horikawa--smooth--stratification} shows that the deformation problem of Horikawa surfaces with only Du Val singularities is completely solved.
In this subsection, we aim to adress the partial $\Q$-Gorenstein smoothing problem for all normal stable Horikawa surfaces except supersingular standard Horikawa surfaces.  
In particular, we completely solve the partial $\Q$-Gorenstein smoothing problem for the case where the geometric genus is less than $7$.
To achieve this, we first identify the strata whose closures contain stable Horikawa surfaces of special Lee-Park type.
Although the following result was originally obtained by topological methods in \cite[Theorem 4.5]{MNU},
we provide a more algebraic proof.

\begin{prop}\label{prop--no--defo--spin}
$(\mathscr{F}^6_{\mathrm{LP}\infty}\cup \mathscr{F}^{6,1}_{\mathrm{LPS}})\subset\overline{\mathscr{F}}^{6,1}_{\mathrm{std}}\setminus\overline{\mathscr{F}}^6_{\mathrm{dV},(0)}$.
 \end{prop}

\begin{proof}
We begin the proof by showing that 
the intersection $\mathscr{F}^{6,1}_{\mathrm{LPS}}\cap\overline{\mathscr{F}}^6_{\mathrm{dV},(0)}$ is empty. 
Let $X$ be a Horikawa surface of special Lee-Park type with geometric genus $6$ and multiplicity $1$.
By Proposition \ref{prop--smoothability-of-special--low--p_g}, the surface $X$ is $\mathbb{Q}$-Gorenstein smoothable.
Furthermore, Proposition \ref{prop:non-std_Horikawa} shows that $H^1(X,\mathcal{O}_X)=0$.
In this situation, we can construct a birational morphism $\widetilde{W}\to\Sigma_{7}$ which is an isomorphism near the $(-7)$-curve $\Delta_0$, where $\delta\colon \widetilde{W}\to W$ is the minimal resolution. 
    Let $E$ be the unique $(-1)$-curve on $\widehat{W}$ (cf.~Lemma \ref{lem--modified--manetti--lemma}).
   An explicit calculation shows that
    \[
    K_W+L\sim 20\delta_*E.
    \]
    Let $\varphi\colon X\to W$ be the canonical morphism.
    Then we have $K_X=\varphi^*(K_W+L)$.
    The surface $W$ has exactly two singularities:
    one of type $\frac{1}{25}(1,4)$ and another of type $\frac{1}{4}(1,1)$.
    Suppose, for contradiction, that there exists a $\mathbb{Q}$-Gorenstein family $f\colon \mathscr{X}\to C$ over a smooth curve $C$ such that the central fiber $\mathscr{X}_0\cong X$ and the general fiber $\mathscr{X}_c$ is a stable Horikawa surface of type $(0)$.
    Let $M:=10\pi^*\delta_*E$, so that $M^{[2]}\sim K_X$.
    Since $M$ is either Cartier or $M\sim 3K_X$ locally around $z$ for each point $z\in \mathscr{X}_0$, there exist an affine open neighborhood $z\in U\subset \mathscr{X}$ and a divisorial sheaf $\mathscr{M}_U$ on $U$ such that $\mathscr{M}_U|_{\mathscr{X}_0}\cong M|_{U\cap \mathscr{X}_0}$.
    By shrinking $C$ if necessary, we may choose finitely many such pairs $(U_i, \mathscr{M}_{U_i})$ covering $\mathscr{X}$.
    Let $\mathfrak{m}$ be the maximal ideal of $\mathcal{O}_{C,0}$.
    For each $n>0$, we consider the closed immersion 
    $$
    C_n:=\mathrm{Spec}(\mathcal{O}_{C,0}/\mathfrak{m}^n)\hookrightarrow C,
    $$
    and define $\mathscr{X}_n:=\mathscr{X}\times_CC_n$. 
    We claim that for every $n>0$, there exists
    a coherent sheaf $\mathscr{M}_n$ on $\mathscr{X}_n$ flat over $C_n$, such that $\mathscr{M}_{n}|_{\mathscr{X}_{n-1}}=\mathscr{M}_{n-1}$ and $\mathscr{M}_{n}|_{U_i\cap \mathscr{X}_{n}}\cong \mathscr{M}_{U_i}|_{U_i\cap \mathscr{X}_{n}}$, where $\mathscr{M}_{0}:=M$.
    We prove this claim by induction on $n$.
    Assuming $\mathscr{M}_{n-1}$ exists,
    the obstruction to lifting it to $\mathscr{M}_n$ lies in $H^2(X,\mathcal{O}_X)$.
    This obstruction, denoted by $o_n$, is constructed using the local deformation data $\{(U_i,\mathscr{M}_{U_i})\}$.
    Since $H^1(X,\mathcal{O}_X)=0$, we have an isomorphism 
    \[\omega_{\mathscr{X}_{n-1}/C_{n-1}}\cong \mathscr{M}_{n-1}^{[2]}.
    \]
    To see this, consider the sheaf
    \[
   \mathscr{F}_{n-1}:= \mathbf{Hom}( \mathscr{M}_{n-1}^{[2]},\omega_{\mathscr{X}_{n-1}/C_{n-1}}).
    \]
    Since both sheaves are locally isomorphic, $\mathscr{F}_{n-1}$ is a line bundle.
    Because $H^1(X,\mathcal{O}_X)=0$, we find $\mathscr{F}_{n-1}\cong\mathcal{O}_{\mathscr{X}_{n-1}}$, hence the desired isomorphism follows.
    It then follows that the obstruction to lifting 
    $\omega_{\mathscr{X}_{n-1}/C_{n-1}}$ to some flat sheaf is $2o_n$.
    Since the canonical sheaf can be lifted 
    to $\omega_{\mathscr{X}_{n}/C_n}$, we have $2o_n=0$.
    Thus, $o_n=0$.
    This allows us to construct the desired lift $\mathscr{M}_n$.
    By \cite[Tag 0CTK]{Stacks} and Popescu's theorem \cite[Tag 07GC]{Stacks}, we may assume that there exists a divisorial sheaf $\mathscr{M}$ on the entire family $\mathscr{X}$ such that $\mathscr{M}|_{X}\cong M$.
    Moreover, by construction, $\mathscr{M}^{[2]}\sim_CK_{\mathscr{X}}$.
    However, this contradicts the result of \cite[Lemma 4.2]{horikawa}, which implies that $X$ can not admit any $\mathbb{Q}$-Gorenstein smoothing to smooth Horikawa surfaces of type $(0)$.

On the other hand, Proposition \ref{prop--smoothability-of-special--low--p_g} shows that $\mathscr{F}^{6,1}_{\mathrm{LPS}}\subset\overline{\mathscr{F}}^{6,1}_{\mathrm{std}}$.
The argument for $\mathscr{F}^6_{\mathrm{LP}\infty}$ proceeds in a similar manner.
This completes the proof.
\end{proof}

\begin{prop}\label{prop--p_g=5-special-case}
    Let $X$ be a stable Horikawa surface of special Lee-Park type. 
    Suppose that $p_g(X)\ge5$. 
        Then, $X$ admits a $\mathbb{Q}$-Gorenstein partial smoothing to stable Horikawa surfaces of general Lee-Park type $\mathrm{I}^*$ if and only if $p_g(X)=5$ and $X$ is a stable Horikawa surface of special Lee-Park type with multiplicity two.
\end{prop}

\begin{proof}
   First, we show that a stable Horikawa surface $X$ of special Lee-Park type with $p_g(X)=5$ and multiplicity two admits a $\mathbb{Q}$-Gorenstein partial smoothing to stable Horikawa surfaces of general Lee-Park type $\mathrm{I}^*$. 
   We may assume that the double cone singularity on $X$ is a simple elliptic singularity.
   Let $\widetilde{W}\to W$ be the minimal resolution.
   Then $h^0(-K_{\widetilde{W}})\ge 6$.
    By Construction~\ref{construction}~(6) and Lemma \ref{lem--nakai--moishezon}, the surface $\widetilde{W}$ has a unique $(-1)$-curve $E_3$ and a $(-3)$-curve $E_1$.
    Let $\widetilde{B}$ be the divisor on $\widetilde{W}$ such that $(\widetilde{W}, \frac{1}{2}\widetilde{B})$ is log crepant with respect to $(W, \frac{1}{2}B)$.
   Then $\widetilde{B}-2E_1$ is effective and disjoint from $E_3$.
   Let $h\colon\widehat{W}\to \widetilde{W}$ be the blow-up at the point $E_3\cap E_1$, and and let $\widehat{B}$ be the proper transform of $\widetilde{B}$. 
   Denote by $E_4$ be the exceptional curve of $h$.
   We exhibit the diagram of curves with negative self intersection in $\widehat{W}$ as below.

   
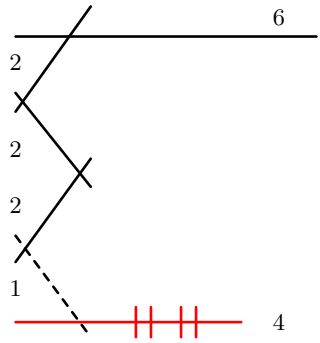
\begin{figure}[H]
       \begin{tikzpicture}[line cap=round,line join=round,>=triangle 45,x=1cm,y=1cm]
\clip(-3,-2.5) rectangle (3,3);
\draw [line width=1pt] (-2,2)-- (2,2);
\draw [line width=1pt] (-1,2.4)-- (-2,1);
\draw [line width=1pt] (-2,1.25)-- (-1,0);
\draw [line width=1pt] (-1,0.375)-- (-2,-1);
\draw [line width=1pt,dashed] (-2,-0.65)-- (-1,-2);
\draw [line width=1pt,red] (-2,-1.8)-- (1,-1.8);
\draw [line width=1pt,red] (-0.4,-1.6)-- (-0.4,-2); 
\draw [line width=1pt,red] (-0.2,-1.6)-- (-0.2,-2); 
\draw [line width=1pt,red] (0.2,-1.6)-- (0.2,-2); 
\draw [line width=1pt,red] (0.4,-1.6)-- (0.4,-2); 
\begin{scriptsize}
\draw [fill] (1.5,2.25) node {$6$};
\draw [fill] (-2,1.65) node {$2$}; 
\draw [fill] (-2,0.5) node {$2$}; 
\draw [fill] (-2,-0.25) node {$2$}; 
\draw [fill] (-2,-1.35) node {$1$}; 
\draw [fill] (1.5,-1.8) node {$4$};
\end{scriptsize}
\end{tikzpicture}
   \end{figure}
   The surface $\widehat{W}$ contains a chain $[6,2,2,2]$ and a $(-4)$-curve.
   The divisor $\widehat{B}$ does not intersect any component of the chain $[6,2,2,2]$. 
   Contracting the chain $[6,2,2,2]$ and the $(-4)$-curve, yields a projective birational morphism $p\colon \widehat{W}\to Z$, where $Z$ is a klt surface.
   Since $\rho(Z)=1$, both $-K_Z$ and $p_*\widehat{B}$ are ample. 
   Explicit calculation shows that 
   $$
   -K_Z\sim 26 p_*E_4,\quad p_*\widehat{B}\sim84 p_*E_4.
   $$
   It follows that  
   $p_{*}\widehat{B} \sim 0$ locally around the singularity of type $\frac{1}{4}(1,1)$.
   Applying \cite[Proposition 3.1]{HP}, the surface $Z$ admits a one-parameter partial $\mathbb{Q}$-Gorenstein smoothing $\mathscr{Z}\to C$ in which the singularity $\frac{1}{4}(1,1)$ is smoothed, while $\mathscr{Z}\to C$ remains locally trivial around the singularity $\frac{1}{21}(1,4)$.
   Since $-K_{Z}$ is ample, \cite[Proposition 1.41]{KoMo} implies that $-K_{\mathscr{Z}_c}$ is also ample.
   By the Kawamata--Viehweg vanishing theorem \cite[Theorem 2.70]{KoMo}, we have $H^1(Z,\mathcal{O}_{Z}(p_*\widehat{B}))=0$.
   Moreover, since $p_*\widehat{B}$ is Cartier around the singularity $\frac{1}{4}(1,1)$, we may assume the existence of a $\mathbb{Q}$-Cartier relative Mumford divisor $\mathscr{B}_{\mathscr{Z}}$ on $\mathscr{Z}$ such that $\mathscr{B}_{\mathscr{Z}}|_{\mathscr{Z}_0}=p_*\widehat{B}$.
   As $\mathscr{Z}$ is locally trivial around the singularity $\frac{1}{21}(1,4)$ and $p_*\widehat{B}$ avoids the singularity, we may further assume that $\mathscr{B}_{\mathscr{Z}}$ never intersects the singularity $\frac{1}{21}(1,4)$ in any fiber over $C$.
We now perform a horizontal partial blow-up of the singularity $\frac{1}{21}(1,4)$ to $\frac{1}{16}(1,3)$  (see the proof of Proposition \ref{prop--Lee-Park-from-I*-to-I}), yielding a projective birational morphism $\pi\colon \overline{\mathscr{W}}\to \mathscr{Z}$ of normal threefolds.
The central fiber $\overline{W}:=\overline{\mathscr{W}}_0$ is a normal surface with only Wahl singularities $\frac{1}{16}(1,3)$ and $\frac{1}{4}(1,1)$.
Let $\overline{\mathscr{B}}:=\pi^*\mathscr{B}_{\mathscr{Z}}$ and $\overline{B}:=\overline{\mathscr{B}}_0$.
Let $f\colon \overline{W}\to W$ be the canonical birational morphism.
Then, 
\[
K_{\overline{W}}+\frac{1}{2}\overline{B}=f^*\left(K_W+\frac{1}{2}B\right).
\] 
Explicit computation shows that there exists a divisorial sheaf $L_{\overline{W}}$ such that $L_{\overline{W}}^{[2]}\cong\mathcal{O}_{\overline{W}}(\overline{B})$ and $L_{\overline{W}}\sim lK_{\overline{W}}$ for some $l\in\mathbb{Z}$ locally around each singularity. 
By Proposition \ref{prop--smoothability--criterion}, this sheaf lifts to a divisorial sheaf $\mathscr{L}$ on $\overline{\mathscr{W}}$ with $\mathscr{L}|_{\overline{W}}=L_{\overline{W}}$.
Since $\overline{\mathscr{W}}_c$ and $\mathscr{Z}_c$ has the same minimal resolution and $-K_{\mathscr{Z}_c}$ is ample with only klt singularities, $-K_{\overline{\mathscr{W}}_c}$ is big.
Thus, $(\overline{\mathscr{W}}_c, \overline{\mathscr{B}}_c)$ satisfies the assumptions of Proposition \ref{prop--smoothing--quotient}, and $\rho(\overline{\mathscr{W}}_c)=\rho(W)=2$ by \cite[Proposition 2.6]{HP}.
Take the relative log canonical model $\eta\colon \overline{\mathscr{W}}\to\mathscr{W}$ with respect to  $K_{\overline{\mathscr{W}}}+\frac{1}{2}\overline{\mathscr{B}}$ over $C$.
By \cite[Lemma 4.2]{MZ}, we have $\mathscr{W}_0\cong W$.

From Lemma \ref{lem--lower-semi--rho}, 
\[
1=\rho(W)\le \rho(\mathscr{W}_c)\le \rho(\overline{\mathscr{W}}_c)=2.
\]
Assuming $\rho(W)= \rho(\mathscr{W}_c)=1$,
we deduce that $\mathscr{W}$ is $\mathbb{Q}$-factorial and hence $\mathbb{Q}$-Gorenstein by Lemma \ref{lem--lower-semi--rho}.
But this contradicts the existence of non-$\mathbb{Q}$-Gorenstein smoothable singularity on $W$.
Therefore, $\rho(\mathscr{W}_c)=2$, and so $\mathscr{W}_c\cong\overline{\mathscr{W}}_c$ by Proposition \ref{prop--anti-P}.
By Remark \ref{rem--smoothing--quotient}, $X$ admits a partial $\mathbb{Q}$-Gorenstein smoothing to surfaces $\mathscr{X}_c$ defined as the double cover of $\overline{\mathscr{W}}_c$ branched along $\overline{\mathscr{B}}_c$.
Since $\overline{\mathscr{W}}_c$ has only one Wahl singularity, and $-K_{\overline{\mathscr{W}}_c}$ is big, we conclude from Proposition \ref{prop--smoothing--quotient} and Remark \ref{rem--smoothing--quotient} (1) that $\mathscr{X}_c$ is partially $\mathbb{Q}$-Gorenstein smoothable to normal standard Horikawa surfaces.
Since $p_g(\mathscr{X}_c)=5$, these surfaces are $\mathbb{Q}$-Gorenstein smoothable.
and so $\mathscr{X}_c$ is also $\mathbb{Q}$-Gorenstein smoothable as well.
By Theorem \ref{intro-thm:Q2} and the proof of Lemma \ref{lem--part--smoothing--again--smoothable}, as $\overline{\mathscr{W}}_c$ has one $(-2)$-curve and one Wahl singularity, we conclude that $\mathscr{X}_c$ is a stable Horikawa surface of general Lee-Park type $\mathrm{I}^*$.

Next, we consider the remaining cases.
If $p_g(X)\ge7$, and $X$ admits a $\mathbb{Q}$-Gorenstein partial smoothing to stable Horikawa surfaces of general Lee-Park type $\mathrm{I}^*$,
then Corollary \ref{cor--properties--anti-P-special} implies that $W$ admits an admissible anti-P-resolution $W^-$ such that $-K_{W^-}$ is nef and log big.
Moreover, there exists a log $\mathbb{Q}$-Gorenstein family $(\mathscr{W}^-,\frac{1}{2}\mathscr{B}^-)\to C$ with central fiber $(\mathscr{W}_0^-,\mathscr{B}^-_0)\cong (W^-,B^-)$, such that the double cover of $\mathscr{W}_c^-$ branched along $\mathscr{B}^-_c$ is a Horikawa surface of general Lee-Park type $\mathrm{I}^*$ for all $c\in C\setminus\{0\}$.
By Lemma \ref{lem--lower-semi--rho}, $\mathscr{B}^-$ is $\mathbb{Q}$-Cartier and $-K_{\mathscr{W}_c^-}$ is not nef.
However, \cite[Theorems 1.10 and 1.15]{fujino--slc--vanishing} imply that $-K_{\mathscr{W}^-}$ is relatively semiample after shrinking $C$, giving a contradiction.

Finally, for $p_g=6$, the claim follows from Proposition \ref{prop--no--defo--spin}. 
This completes the proof.
\end{proof}

\begin{prop}\label{prop--LPI-to-LPS}
 Let $X$ be a stable Horikawa surface of special Lee-Park type with only equivariantly smoothable singularities. 
    Suppose that $p_g(X)\ne 6$.
    Then, $X$ admits a $\mathbb{Q}$-Gorenstein partial smoothing to stable Horikawa surfaces of general Lee-Park type $\mathrm{I}$.  
\end{prop}

\begin{proof}
    We first assume that $p_g(X)\ge7$.
    In this case, the surface $W$ admits an admissible anti-P-resolution $\mu\colon (W^-,\frac{1}{2}B^-)\to (W,\frac{1}{2}B)$ which is log $\mathbb{Q}$-Gorenstein smoothable to either $\Sigma_0$ or $\Sigma_1$ by Theorem \ref{thm--smoothability-of-special--cusp}.
By Corollary \ref{cor--properties--anti-P-special}, the divisor $-K_{W^-}$ is nef and log big, and $K_{W^-}+\frac{1}{2}B^-$ is big and semiample.
Therefore, by Proposition \ref{prop--smoothing--slc},    there exists a morphism $(\mathscr{W}^-,\frac{1}{2}\mathscr{B}^-)\to C$ from a normal $\mathbb{Q}$-Gorenstein threefold to a smooth curve $C$ with a closed point $0\in C$, along with a $\mathbb{Q}$-Cartier effective relative Mumford divisor $\mathscr{B}^-$, such that
\begin{itemize}
    \item 
    $(W^-,B^-)\cong (\mathscr{W}^-_0,\mathscr{B}_0^-)$, and
    \item 
    $\mathscr{W}^-\to C$ smooths all singularities of $W^-$ except for the Wahl singularity of type $\frac{1}{(p_g-1)^2}(1,p_g-2)$, around which the family is locally trivial.
\end{itemize}
Moreover, by Corollary \ref{cor--gen--fiber--rho}, 
the divisor $K_{\mathscr{W}^-_c}+\frac{1}{2}\mathscr{B}^-_c$ is ample for general $c\in C$.
By Lemma \ref{lem--part--smoothing--again--smoothable},  the double cover of $\mathscr{W}^-_c$ branched along $\mathscr{B}^-_c$ is a $\Q$-Gorenstein smoothable normal Horikawa surface for any general $c\in C$.
Since $-K_{W^-}$ is nef, it follows that $-K_{\mathscr{W}^-_c}$ is also nef for very general $c\in C$. 
Theregore, by Theorem \ref{intro-thm:Q2}, the double cover of $\mathscr{W}^-_c$ branched along $\mathscr{B}^-_c$ is a Horikawa surface of general Lee-Park type $\mathrm{I}$.
This establishes the claim in the case $p_g\ge7$.

The remaining case can be proved by a similar argument, using Corollary \ref{cor--gen--fiber--rho}.
\end{proof}

In the case where geometric genus is at most six, various sporadic phenomena arise.
Before proving Theorem \ref{thm--stratification--intro}, we treat some difficult cases separately.
We begin with the following lemma, which adresses a partial smoothing problem for Horikawa surfaces in the case $p_g=3$.

\begin{lem}\label{lem--p_g=3--part--sm}
   \[
   \mathscr{F}^3_{(2,2,2,2),[2433]}\subset \overline{\mathscr{F}}^3_{(2,2,2,2)}.
   \]
\end{lem}

\begin{proof}
    Let $X$ be a stable normal Horikawa surface with geometric genus three with one strictly lc singularities of type $(2,2,2,2)$ and one singularity of type $\frac{1}{50}(1,29)$.
    Then the associated surface $W$ is a projective klt surface with only Wahl singularities with $\rho(W)=1$.
    Since $-K_W$ is ample, we can construct a projective flat morphism $f\colon \mathscr{W}\to C$ such that 
    \begin{itemize}
    \item 
    $f$ is $\mathbb{Q}$-Gorenstein,
    \item 
    the central fiber satisfies $\mathscr{W}_0\cong W$, and
    \item 
    $f$ gives a $\mathbb{Q}$-Gorenstein smoothing locally around the singularity of type $\frac{1}{25}(1,4)$,
    while being locally trivial around the singularity of type $\frac{1}{4}(1,1)$.
    \end{itemize}
    
    Note that $W\cong \mathbb{P}(1,4,25)$.
    Let $\xi\colon \widehat{W}\to W$ be the minimal resolution of the singularity $\frac{1}{4}(1,1)$, and let $E$ denote the exceptional divisor.
    Note that the divisor $\xi^*B-K_{\widehat{W}}-E$ is big and nef.
    By the Kawamata-Viehweg vanishing theorem, we have
    $H^1(\widehat{W},\mathcal{O}_{\widehat{W}}(\xi^*B-E))=0$.
    Let $\pi\colon \widehat{\mathscr{W}}\to \mathscr{W}$ be the partial simultaneous minimal resolution of the horizontal singularity of type $\frac{1}{4}(1,1)$, and let $\mathscr{E}$ be the exceptional divisor.
    Applying Proposition \ref{prop--smoothability--criterion} to the divisor $\xi^*B-E$ and the vanishing above, we may assume the existence of an effective $\mathbb{Q}$-Cartier relative Mumford divisor $\widehat{\mathscr{B}}$ on $\widehat{\mathscr{W}}$ such that $\widehat{\mathscr{B}}_0=\xi^*B-E$.
    We then define 
    $$
    \mathscr{B}=\pi_*(\widehat{\mathscr{B}}+\mathscr{E}).
    $$
    By again applying Proposition \ref{prop--smoothability--criterion} to the Weil divisor $L$ on $W$ satisfying $L^{[2]}=B$, we obtain a divisorial sheaf $\mathscr{L}$ on $\mathscr{W}$ such that $\mathscr{L}_0=L$ after replacing a finite base change over $C$. 
    Since $H^1(W,\mathcal{O}_W)=0$, it follows that $\mathscr{L}^{[2]}\sim_{C}\mathcal{O}_{\mathscr{W}}(\mathscr{B})$.
    Then, by Lemma \ref{lem--part--smoothing--again--smoothable}, the double cover $\mathscr{X}\to \mathscr{W}$ branched along $\mathscr{B}$ defines a partial $\mathbb{Q}$-Gorenstein smoothing from $X$ to stable normal Horikawa surfaces with only one strictly lc singularity of type $(2,2,2,2)$.
\end{proof}

\begin{lem}\label{lem--cone--part--smooth}
    \[\mathscr{F}^4_{\mathrm{cone}}\subset\overline{\mathscr{F}}^{4,0}_{\mathrm{std}}\setminus(\overline{\mathscr{F}}^{4,1}_{\mathrm{std}}\cup\overline{\mathscr{F}}^{4,2}_{\mathrm{std}}\cup\overline{\mathscr{F}}^{4}_{\mathrm{LPI}}).\]
\end{lem}

\begin{proof}
    In Proposition \ref{prop--canonical--pencil--invol}, we showed that $(W,\frac{1}{2}B)\in \mathscr{F}^4_{\mathrm{cone}}(\mathbb{C})$ is log $\mathbb{Q}$-Gorenstein smoothable to $(\overline{\Sigma}_2,\frac{1}{2}B_{\overline{\Sigma}_2})$, where the double cover of $\overline{\Sigma}_2$ branched along $B_{\overline{\Sigma}_2}$ is a standard Horkawa surface.
    We observe that the singularity of $\overline{\Sigma}_2$ degenerates to the cone singularity of $W$.
    On the other hand, since $(W,\frac{1}{2}B)$ is log canonical, $B$ does not pass through the singularity.
    It follows that $B_{\overline{\Sigma}_2}$ also avoids the singular point of $\overline{\Sigma}_2$.
    Therefore, the corresponding point of $(W,\frac{1}{2}B)$ lies in 
    $$  
    \overline{\mathscr{F}}^{4,0}_{\mathrm{std}}\setminus (\overline{\mathscr{F}}^{4,1}_{\mathrm{std}} \cup \overline{\mathscr{F}}^{4,2}_{\mathrm{std}}).    
    $$
    Now, assume that the point $[(W,\frac{1}{2}B)]$ lies in $\overline{\mathscr{F}}^{4}_{\mathrm{LPI}}$.
In this case, $W$ would deform to non-Gorenstein surfaces.
This contradicts the fact that $W$ is Gorenstein (cf.~\cite[Theorem 9.1.6]{Ishii}).
This completes the proof.
\end{proof}

We now turn to the case where $p_g(X)=5$.
As shown in Proposition \ref{prop--p_g=5-special-case}, a special phenomenon occurs for Horikawa surfaces of special Lee-Park type with multiplicity two when the geometric genus is five. 
On the other hand, Horikawa surfaces of special Lee-Park type with multiplicity one and geometric genus five behave similarly to those in other geometric genus cases.

\begin{lem}\label{lem--p_g=5--ordinary}
    \[
    \mathscr{F}^{5,1}_{\mathrm{LPS}}\subset\overline{\mathscr{F}}^{5,1}_{\mathrm{std}}\cap\overline{\mathscr{F}}^{5}_{\mathrm{LPI}}\setminus\overline{\mathscr{F}}^{5}_{\mathrm{LPI}^*}.
    \]
\end{lem}
\begin{proof}
Let $(W,\frac{1}{2}B)\in\mathscr{F}^{5,1}_{\mathrm{LPS}}(\mathbb{C})$ be an arbitrary object.
    Then, $(W,\frac{1}{2}B)\not\in\overline{\mathscr{F}}^{5}_{\mathrm{LPI}^*}(\mathbb{C})$ follows from a similar argument to the proof of Corollary \ref{cor--properties--anti-P-special} by considering anti-P-resolutions of the singularity of type $\frac{1}{3}(1,1)$.
    The remaining part follows from Lemma \ref{lem--part--smoothing--again--smoothable} and Theorem \ref{intro-thm:Q2} by taking an anti-P-resolution similar to Remark \ref{rem--anti--p}.
\end{proof}

Finally, we consider $\Q$-Gorenstein smoothable Horikawa surfaces with geometric genus six.

\begin{prop}
    Let $X$ be a stable normal Horikawa surface of special Lee-Park type with $p_g(X)=6$ and multiplicity two.
    Then, $X$ is partially $\mathbb{Q}$-Gorenstein smoothable to normal stable Horikawa surfaces of general Lee-Park type $\mathrm{I}$.
\end{prop}

\begin{proof}
It is straightforward to verify that $(W,\frac{1}{2}B)$ admits an anti-P-resolution 
$$
\mu\colon \left(W^-,\frac{1}{2}B^-\right)\to \left(W,\frac{1}{2}B\right)
$$
such that both $B^-$ and $K_{W^-}$ are ample and $W^-$ has only T-singularities as Proposition \ref{prop--smoothability-of-special--cusp--k_1=0}.
Then, by Proposition \ref{prop--smoothing--quotient}, there exists a projective flat morphism $f\colon \mathscr{W}^-\to C$ from a $\Q$-Gorenstein klt threefold to a smooth curve with a closed point $0\in C$, an effective relative Mumford divisor $\mathscr{B}^-$, and a $\mathbb{Q}$-Cartier divisorial sheaf $\mathscr{L}$ such that 
\begin{itemize}
    \item 
    $\mathscr{L}^{[2]}\sim_C \mathscr{B}^-$,
    \item 
    $(\mathscr{W}^-_0,\mathscr{B}^-_0)\cong (W^-,B^-)$, 
    \item 
    $f$ is locally trivial at the Wahl singularity of type $\frac{1}{25}(1,4)$, and
    \item 
    $f$ is a $\mathbb{Q}$-Gorenstein smoothing around every other singularity.  
\end{itemize}
Let $(\mathscr{W},\frac{1}{2}\mathscr{B}+\mathscr{W}_0)$ be the relative log canonical model of $(\mathscr{W}^-,\frac{1}{2}\mathscr{B}^-+\mathscr{W}^-_0)$ over $C$. 
Note that $\mathscr{W}_0\cong W$ and that $-K_{\mathscr{W}^-}$ is relatively ample over $\mathscr{W}$.
This shows that $K_{\mathscr{W}}$ is not $\Q$-Cartier. 
Then, Corollary \ref{cor--gen--fiber--rho} shows that 
$\mathscr{W}^-\setminus\mathscr{W}^-_0\cong\mathscr{W}\setminus\mathscr{W}_0$.
By applying Theorem \ref{intro-thm:involution}, Proposition \ref{prop--construction--lee--park--i*}, Lemma \ref{lem--part--smoothing--again--smoothable}, and using the fact that $-K_{W^-}$ is nef, we conclude that the double cover of $\mathscr{W}_c$ branched along $\mathscr{B}_c$ defines a Horikawa surface of general Lee-Park type $\mathrm{I}$.
\end{proof}

\begin{prop}
    $\mathscr{F}^{6,2}_{\mathrm{std}}\subset \overline{\mathscr{F}}^6_{\mathrm{dV},(0)}$.
\end{prop}

\begin{proof}
    Let $(W,\frac{1}{2}B)\in \mathscr{F}^{6,2}_{\mathrm{std}}(\mathbb{C})$ be an object.
    Note that $W\cong\overline{\Sigma}_4$.
    To prove the assertion, we may replace $(W,\frac{1}{2}B)$ with a general member of the same stratum.
    In particular, we may assume that 
    \begin{itemize}
        \item 
        $B^+-2\Delta_0$ intersects $\Delta_0$ transversally, and
        \item
        $B^+$ contains no fiber of the ruling $W^{+}=\Sigma_4\to\mathbb{P}^1$.
    \end{itemize}    
    It is easy to check that there exists an admissible anti-P-resolution $\mu\colon(W^-,\frac{1}{2}B^-)\to (W,\frac{1}{2}B)$ such that both $-K_{W^-}$ and $B^-$ are ample as in Corollary \ref{cor--properties--anti-P-standard}.
    By Proposition \ref{prop--smoothing--quotient}, there exists a flat projective morphism $f\colon\mathscr{W}^-\to C$ from a klt threefold to a smooth curve $C$ with a closed point $0\in C$, an effective relative Mumford divisor $\mathscr{B}^-$ and a $\mathbb{Q}$-Cartier divisorial sheaf $\mathscr{L}$ on $\mathscr{W}^-$ such that
    \begin{itemize}
        \item 
        $\mathscr{B}^-\sim_C\mathscr{L}^{[2]}$,
        \item 
         $(\mathscr{W}^-_0,\mathscr{B}^-_0)\cong (W^-,B^-)$, and
         \item 
         $\mathscr{W}^-_c$ is smooth for any $c\in C\setminus\{0\}$.    
    \end{itemize}
     
    Since $-K_{\mathscr{W}^-_c}$ is ample by \cite[Proposition 1.41]{KoMo}, we conclude from \cite[Proposition 2.6]{HP} that $\mathscr{W}^-_c\cong \Sigma_0$.
    Thus, we may assume that $\mathscr{B}^-_c$ is smooth and that $K_{\mathscr{W}^-_c}+\frac{1}{2}\mathscr{B}_{c}^-$ is ample for any $c\in C\setminus\{0\}$.
    Again by Proposition \ref{prop--smoothing--quotient}, it follows that $(\mathscr{W}^-_c,\frac{1}{2}\mathscr{B}^-_c)$ belongs to $\mathscr{F}^6_{\mathrm{dV},(0)}$ for any $c\in C\setminus\{0\}$.
    This completes the proof.
\end{proof}

Now, we state the main result of this subsection.

\begin{thm}\label{thm--stratification}
    The following are the Hasse diagrams of strata in the moduli stacks $(\overline{\mathcal{M}}^{\mathrm{Gie}}_{2p_g-4,p_g})^{\mathrm{sn}}$ of $\Q$-Gorenstein smoothable normal stable Horikawa surfaces of every geometric genus $p_g$.
    
    A directed path in the diagram indicates that the moduli stack represented by the lower node is contained in the Zariski closure of the moduli stack represented by the upper node.
If there is no connected line between two nodes, then the corresponding moduli stack in one node is disjoint from the closure of the other.
    
   \newpage
    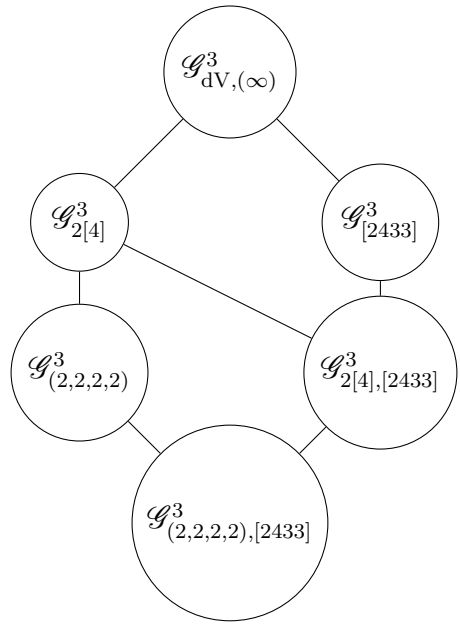
\begin{figure}[htbp]
  \centering
  \begin{tikzpicture}[scale=1, transform shape,
    every node/.style={draw, circle, minimum size=1cm},
    every edge/.style={draw, thick}]
    
    \node (empty) at (0,0) {$\mathscr{G}^3_{(2,2,2,2),[2433]}$};
    \node (1)     at (-2,2) {$\mathscr{G}^3_{(2,2,2,2)}$};
    \node (2)     at ( 2,2) {$\mathscr{G}^3_{\mathrm{2[4],[2433]}}$};
    \node (12)    at (-2,4) {$\mathscr{G}^3_{2[4]}$};
    \node (13)    at ( 2,4) {$\mathscr{G}^3_{[2433]}$};
    \node (123)   at (0,6)  {$\mathscr{G}^3_{\mathrm{dV},(\infty)}$};

    \draw (empty) -- (1);
    \draw (empty) -- (2);
    \draw (1) -- (12);
    \draw (2) -- (12);
    \draw (2) -- (13);
    \draw (12) -- (123);
    \draw (13) -- (123);
  \end{tikzpicture}
  \caption{$p_g=3$ Case}
  \label{fig:hasse-p_g=3}
\end{figure}

  \begin{figure}[htbp]
  \centering
  \begin{tikzpicture}[scale=1, transform shape,
    every node/.style={draw, circle, minimum size=1cm},
    every edge/.style={draw, thick}]
    
    \node (LPS2)     at ( -2,0) {$\mathscr{G}^{4,2}_{\mathrm{LPS}}$};
    \node (LPS1)     at ( -2,2) {$\mathscr{G}^{4,1}_{\mathrm{LPS}}$};
    \node (LPS0)     at ( -2,4) {$\mathscr{G}^{4,0}_{\mathrm{LPI}}$};
    \node (2'2)   at (0,2)  {$\mathscr{G}^{4,2}_{\mathrm{std}}$};
    \node (LPI)     at ( -2,6) {$\mathscr{G}^4_{\mathrm{LPI}}$};
    \node (cone)    at (2,4) {$\mathscr{G}^4_{\mathrm{cone}}$};
    \node (2'1)   at (0,4)  {$\mathscr{G}^{4,1}_{\mathrm{std}}$};
    \node (2'0)   at (0,6)  {$\mathscr{G}^{4,0}_{\mathrm{std}}$};
    \node (0)   at (0,8)  {$\mathscr{G}^4_{\mathrm{dV},(0)}$};
\draw (LPS2) -- (LPS1);
\draw (LPS2) -- (2'2);
\draw (LPS1) -- (LPS0);
\draw (LPS1) -- (2'1);
\draw (LPS0) -- (LPI);
\draw (LPS0) -- (2'0);
\draw (cone) -- (2'0);
\draw (2'0)  --  (0) ;
\draw (LPI) -- (0);
\draw (2'2) -- (2'1);
\draw (2'1) -- (2'0);

  \end{tikzpicture}
  \caption{$p_g=4$ Case}
  \label{fig:hasse-p_g=4}
\end{figure}
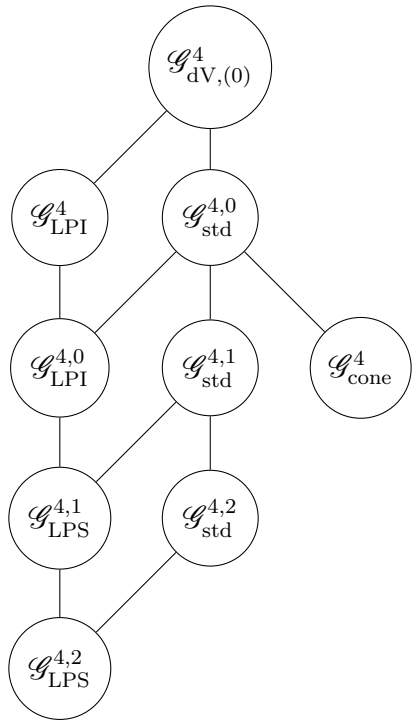
\newpage
\begin{figure}[htbp]
  \centering
  \begin{tikzpicture}[scale=1, transform shape,
    every node/.style={draw, circle, minimum size=1cm},
    every edge/.style={draw, thick}]
    
    \node (LPS2)     at ( 0,0) {$\mathscr{G}^{5,2}_{\mathrm{LPS}}$};
    \node (LPS1)     at ( 0,2) {$\mathscr{G}^{5,1}_{\mathrm{LPS}}$};
    \node (LPI*)     at ( -2,2) {$\mathscr{G}^5_{\mathrm{LPI}^*}$};
    \node (LPI)     at ( 0,4) {$\mathscr{G}^5_{\mathrm{LPI}}$};
    \node (2'2)   at (2,2)  {$\mathscr{G}^{5,2}_{\mathrm{std}}$};
    \node (2'1)   at (2,4)  {$\mathscr{G}^{5,1}_{\mathrm{std}}$};
    \node (0)   at (2,6)  {$\mathscr{G}^5_{\mathrm{dV},(1)}$};
\draw (LPS2) -- (LPS1);
\draw (LPS2) -- (2'2);
\draw (LPS1) -- (2'1);
\draw (LPI) -- (0);
\draw (2'2) -- (2'1);
\draw (2'1) -- (0);
\draw (LPS1) -- (LPI);
\draw (LPI*) -- (LPI);
\draw (LPS2) -- (LPI*);

  \end{tikzpicture}
  \caption{$p_g=5$ Case}
  \label{fig:hasse-p_g=5}
\end{figure}
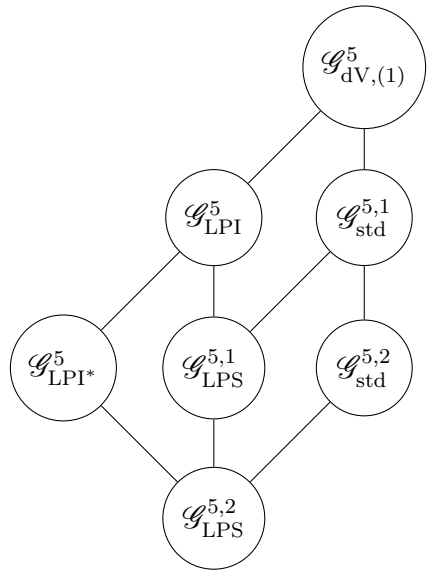
\begin{figure}[htbp]
  \centering
  \begin{tikzpicture}[scale=1, transform shape,
    every node/.style={draw, circle, minimum size=1cm},
    every edge/.style={draw, thick}]
    
    \node (LPS2)     at (-2,0) {$\mathscr{G}^{6,2}_{\mathrm{LPS}}$};
    \node (LPS1)     at ( 2,4) {$\mathscr{G}^{6,1}_{\mathrm{LPS}}$};
    \node (MNU)     at ( 2,6) {$\mathscr{G}^6_{\mathrm{LP}\infty}$};
    \node (2'2)   at (2,2)  {$\mathscr{G}^{6,2}_{\mathrm{std}}$};
    \node (2'1)   at (4,6)  {$\mathscr{G}^{6,1}_{\mathrm{std}}$};
    \node (infty)   at (4,8)  {$\mathscr{G}^6_{\mathrm{dV},(\infty)}$};
    \node (0)   at (-2,8)  {$\mathscr{G}^6_{\mathrm{dV},(0)}$};
    \node (LPI)   at (-2,6)  {$\mathscr{G}^6_{\mathrm{LPI}}$};
    \node (LPI*)   at (-4,4)  {$\mathscr{G}^6_{\mathrm{LPI}*}$};
\draw (LPS2) -- (LPS1);
\draw (LPS2) -- (2'2);
\draw (LPS1) -- (2'1);
\draw (MNU) -- (infty);
\draw (2'2) -- (2'1);
\draw (2'1) -- (infty);
\draw (LPS1) -- (MNU);
\draw (0) -- (LPI);
\draw (LPI*) -- (LPI);
\draw (LPS2) -- (LPI);
    \draw (2'2) -- (0);
  \end{tikzpicture}
  \caption{$p_g=6$ Case}
  \label{fig:hasse-p_g=6}
\end{figure}
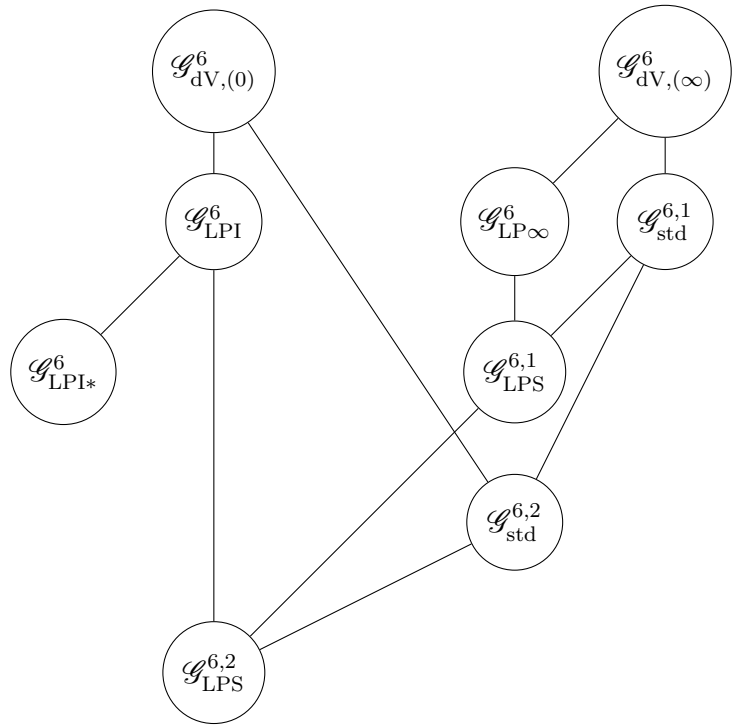

\newpage
\begin{figure}[htbp]
  \centering
  \begin{tikzpicture}[scale=1, transform shape,
    every node/.style={draw, circle, minimum size=1cm},
    every edge/.style={draw, thick}]
    
    \node (LPS)     at ( 1,2) {$\mathscr{G}^{p_g}_{\mathrm{LPS}}$};
    \node (LPI*)     at ( -2,2) {$\mathscr{G}^{p_g}_{\mathrm{LPI}^*}$};
    \node (LPI)     at ( 0,4) {$\mathscr{G}^{p_g}_{\mathrm{LPI}}$};
    \node (stdnf)   at (4,2)  {$\mathscr{G}^{p_g}_{\mathrm{stdnF}}$};
    \node (stdf)   at (2,4)  {$\mathscr{G}^{p_g}_{\mathrm{stdF}}$};
    \node (0)   at (1,6)  {$\mathscr{G}^{p_g}_{\mathrm{dV},(2\left\{\frac{p_g}{2}\right\})}$};
\draw (0) -- (stdf);
\draw (stdf) -- (stdnf);
\draw (LPS) -- (stdf);
\draw (LPS) -- (LPI);
\draw (LPI*) -- (LPI);
\draw (LPI) -- (0);

  \end{tikzpicture}
  \caption{$p_g\ge7$ Case where $p_g-2\not\in4\mathbb{Z}$}
  \label{fig:hasse-p_gge7}
\end{figure}
\begin{figure}[htbp]
  \centering
  \begin{tikzpicture}[scale=1, transform shape,
    every node/.style={draw, circle, minimum size=1cm},
    every edge/.style={draw, thick}]
    
    \node (LPS)     at ( 1,2) {$\mathscr{G}^{p_g}_{\mathrm{LPS}}$};
    \node (LPI*)     at ( -2,2) {$\mathscr{G}^{p_g}_{\mathrm{LPI}^*}$};
    \node (LPI)     at ( 0,4) {$\mathscr{G}^{p_g}_{\mathrm{LPI}}$};
    \node (stdnf)   at (4,2)  {$\mathscr{G}^{p_g}_{\mathrm{stdnF}}$};
    \node (stdf)   at (2,4)  {$\mathscr{G}^{p_g}_{\mathrm{stdF}}$};
    \node (0)   at (1,6)  {$\mathscr{G}^{p_g}_{\mathrm{dV},(0)}$};
    \node (2k+2)   at (-4,6)  {$\mathscr{G}^{p_g}_{\mathrm{dV},(\frac{p_g}{2}+1)}$};
\draw (0) -- (stdf);
\draw (stdf) -- (stdnf);
\draw (LPS) -- (stdf);
\draw (LPS) -- (LPI);
\draw (LPI*) -- (LPI);
\draw (LPI) -- (0);

  \end{tikzpicture}
  \caption{$p_g\ge7$ Case where $p_g-2\in4\mathbb{Z}$}
  \label{fig:hasse-p_gge10mod4}
\end{figure}
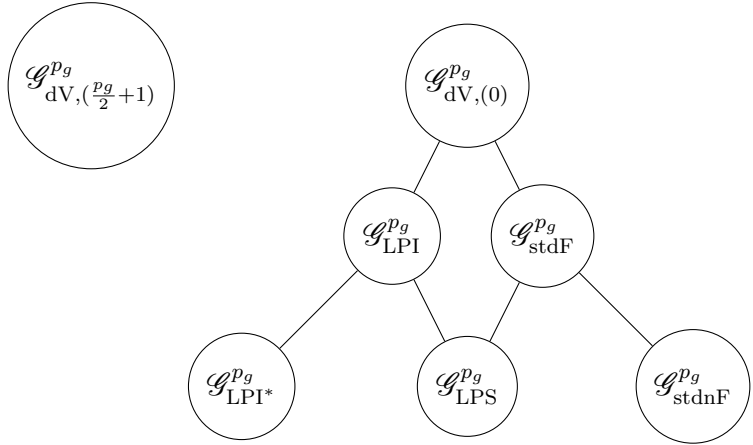

\end{thm}

\begin{proof}
First, the case where $p_g\ge7$ follows from Theorem \ref{intro-thm:Q2}, Propositions \ref{prop--partial--smoothing--from--non-Fano-to-Fano-with-anti-P}, \ref{prop--Lee-Park-from-I*-to-I}, \ref{prop--partial--smoothing--result--for--special}, \ref{prop--smoothability-of-special--cusp--k_1k_2not0}, \ref{prop--p_g=5-special-case}, \ref{prop--LPI-to-LPS}, Corollaries \ref{cor--properties--anti-P-standard-ii}, and \ref{cor--properties--anti-P-special}.

    We now treat the case where $p_g=3$.
    The most difficult part of this case has already been handled in Lemma \ref{lem--p_g=3--part--sm}.
    The remaining cases follow easily from Lemma \ref{lem--part--smoothing--again--smoothable} and \cite[Proposition 3.1]{HP}.
    
For $p_g=4$, the result follows from similar arguments as in the case $p_g\ge7$, together with Corollary \ref{cor--m(X)--down--special--type}, Lemmas \ref{lem--standard-for-low-p_g-smoothable}, \ref{lem--part--smoothing--again--smoothable}, \ref{lem--cone--part--smooth}, \cite[Proposition 3.1]{HP}, and Proposition \ref{prop--smoothability-of-special--low--p_g}.

In the case $p_g=5$, the result follows from arguments similar to those used for $p_g=4$ or $p_g \ge7$, as well as Lemma \ref{lem--p_g=5--ordinary}, Propositions \ref{prop--p_g=5-special-case} and \ref{prop--LPI-to-LPS}.

Finally, the case $p_g=6$ follows a similar argument to the case where $p_g=3,4,5$ or $\ge7$, along with Proposition \ref{prop--no--defo--spin}.
\end{proof}

\begin{proof}[{Proof of Theorem \ref{thm--connectedness}}]
    This follows from Theorem \ref{thm--stratification} and Corollary \ref{cor--p_g=10--connected}.
\end{proof}

\appendix
\section{List of singular fibers of elliptic surfaces with an involution}\label{app:sing_fiber_inv}

In this Appendix, we classify the involutions of germs of an elliptic singular fiber with a certain condition.
Let $f:Y \to T$ be a germ of an elliptic surface over an open disk $0 \in T$ with a central fiber $F_{0}$.
We consider the involution $\sigma:Y \to Y \in \mathrm{Aut}_T(Y)$ that satisfies the following conditions:
\begin{itemize}
     \item[($\star$)] The general fiber of the induced fibration $Y / \sigma \to T$ is a smooth rational curve.
\end{itemize}
Let $\rho:\widetilde{Y} \to Y$ be a resolution of the isolated fixed point of the involution $\sigma$ and let $\widetilde{\sigma} \in \mathrm{Aut}_T(\widetilde{Y})$ be the induced involution on $\widetilde{Y}$ by $\sigma$.
Since $\widetilde{\sigma}$ has no isolated fixed point, 
the quotient $\widetilde{W}:=\widetilde{Y}/\sigma$ is smooth.
Let $\widetilde{\varphi}:\widetilde{W} \to T$ be the natural ruling.
We denote by $\widetilde{B}$ the branch locus of the double covering $\widetilde{Y} \to \widetilde{W}$.
\begin{lem}$($cf. \cite[Lemma~5]{horikawa-genus2}$)$
\label{lem:normalized_branch}
There exists a relatively minimal model $\varphi:W \to T$ of $\widetilde{\varphi}$ that satisfies the following:
Let $\psi:\widetilde{W} \to W$ be the contraction map,
let $B :=\psi_{*}\widetilde{B}$ and let $\Gamma_{0}$ be a cetral fiber of $\varphi$.
For any point $x \in \Gamma_0$, it holds that $\mathrm{mult}_{x}B_{h} \leq 2$, where $B_{h}$ be the sub divisor consisting horizontal components of $B$.
\end{lem} 
We take a relatively minimal model $\varphi:W \to T$ of $\widetilde{\varphi}$ as Lemma~\ref{lem:normalized_branch}.
We have the following diagram:
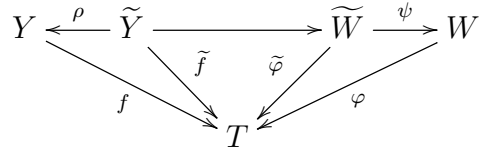
\[\xymatrix{
Y  \ar[rrd]_{f}  & \ar[l]_{\rho} \widetilde{Y} \ar[rr]  \ar[rd]^{\widetilde{f}} &     & \widetilde{W} \ar[ld]_{\widetilde{\varphi}} \ar[r]^{\psi} & W \ar[lld]^{\varphi}\\
 & & T &   &   \\
}\]

We decompose this into a sequence of blow-ups $\psi=\psi_{1}\circ \cdots \circ\psi_{N}\colon \widetilde{W}=W_{N}\to W_{N-1}\to \cdots \to W_{0}=W$.
Considering the even resolution process for $W$ and $B$,
we can define the double covering $X_i \to W_i$ and the branch locus $B_{i}$ for $i=1,\cdots, N$, inductively (Definition~\ref{defn--evenresol}).
Then it hold that $B_{N}=\widetilde{B}$.
Therefore, we see that an involution of an elliptic singular fiber satisfying condition $(\star)$ can be constructed from a ruled surface and a branch locus on it satisfying the conditions in Lemma~\ref{lem:normalized_branch}.

Let $W$, $B$ and $\Gamma_0$ be as Lemma~\ref{lem:normalized_branch} and
let $B_{\mathrm{hor}} = \sum_{i} B_{i}$ be the decomposition into irreducible components.
The branch locus $B$ satisfies one of the following around $\Gamma_0$
\begin{itemize}
    \item[$(1)$] $\Gamma_0 \not \subset B$.
    $B_{\mathrm{hor}}$ is smooth, there exist a point $q \in W$ with $(B_{\mathrm{hor}}\cdot \Gamma_0)_q \ge 3$
    \begin{itemize}
        \item[$(1\text{--}1)$] $B_{\mathrm{hor}}=B_1+B_2$. 
        $B_1$ intersect $\Gamma_0$ transversally and $(B_2\cdot \Gamma_0)_q = 3$
        \item[$(1\text{--}2)$] $B_{\mathrm{hor}}=B_1$ and $(B_1 \cdot \Gamma_0)_q  = 4$
    \end{itemize}
    \item[$(1)^*$] $\Gamma_0 \subset B$.
     $B_{\mathrm{hor}}$ is smooth, there exist a point $q \in W$ with $(B_{\mathrm{hor}}\cdot \Gamma_0)_q \ge 3$.
    As in the case of $(1)$, each $B_{\mathrm{hor}}$ is classified into two types $(1\text{--}1)^*$ and $(1\text{--}2)^*$.
    \item[$(2)$] $\Gamma_0 \not \subset B$.
    For any point $q \in B_{\mathrm{hor}} \cap \Gamma_0$, $(B_{\mathrm{hor}}\cdot\Gamma_0)_{q} \le 2$. 
    \begin{itemize}
        \item[$(2\text{-}1)$] 
         $B_{\mathrm{hor}}=B_1+B_2+B_3+B_4$.
         Each irreducible component $B_i$ intersects $\Gamma_0$ transversally.
         There exist two non-negative even integers $k_1$ and $k_2$ such that $B_1$ and $B_2$ consists a $A_{k_1-1}$-singularity, $B_3$ and $B_4$ consists a $A_{k_2-1}$-singularity.
         $B_1+B_2$ and $B_3+B_4$ are disjoint.
        \item[$(2\text{-}2)$]
         $B_{\mathrm{hor}}=B_1+B_2+B_3$.
         There exist a positive odd integer $k_1$ such that 
         $B_1$ is a $A_{k_1-1}$-singularity whose tangent cone meets $\Gamma_0$ transversally,
         $B_2$ and $B_3$ intersect $\Gamma_0$ transversally.
         There exist a positive a non-negative even integer $k_2$ such that $B_2$ and $B_3$ consists a $A_{k_2-1}$-singularity.
         $B_1$ and $B_2+B_3$ are disjoint.
        \item[$(2\text{-}3)$] 
        $B_{\mathrm{hor}}=B_1+B_2$.
        There exist two positive integers $k_1$ and $k_2$ that satisfies the following:
        $B_1$ is a $A_{k_1-1}$-singularity whose tangent cone meets $\Gamma_0$ transversally. 
        $B_2$ is a $A_{k_2-1}$-singularity whose tangent cone meets $\Gamma_0$ transversally.
         $B_1$ and $B_2$ are disjoint. 
    \end{itemize}
    \item[$(2)^{*}$] 
    $\Gamma_0 \subset B$.
    For any point $q \in B_{\mathrm{hor}} \cap \Gamma_0$, $(B_{\mathrm{hor}}\cdot\Gamma_0)_{q} \le 2$. 
    As in the case of $(2)$, each $B_{\mathrm{hor}}$ is classified into three types $(2\text{--}1)^*$--$(2\text{--}3)^*$.
    \item[$(3)$] $\Gamma_0 \not \subset B$.
    There exists a singular points $q \in B_{\mathrm{hor}}$ such that $(B_{\mathrm{hor}}\cdot\Gamma_0)_{q}=3$.
     \begin{itemize}
        \item[$(3\text{--}1)$] 
         $B_{\mathrm{hor}}=B_1+B_2+B_3$.
         $B_1$ and $B_2$ intersects $\Gamma_0$ transversally.
         $B_3$ is tangent to $\Gamma_0$ at the intersection of
         $B_2$ intersects $\Gamma_0$ with intersection multiplicity $2$.
         $B_1$ and $B_2+B_3$ are disjoint.
        \item[$(3\text{--}2)$]
         $B_{\mathrm{hor}}=B_1+B_2$.
         $B_1$ intersects $\Gamma_0$ transversally.
         $B_2$ is a $A_{1}$-singularity whose tangent cone coinsides with $\Gamma_0$.
         $B_1$ and $B_2$ are disjoint.
    \end{itemize}
    \item[$(3)^*$] $\Gamma_0 \subset B$.
    There exists a singular points $q \in B_{\mathrm{hor}}$ such that $(B_{\mathrm{hor}}\cdot\Gamma_0)_{q}=3$.
    As in the case of $(3)$, each $B_{\mathrm{hor}}$ is classified into two types $(3\text{--}1)^*$ and $(3\text{--}2)^*$.
    \item[$(4)$] $\Gamma_0 \not \subset B$.
    There exists a singular points $q \in B_{\mathrm{hor}}$ with $(B_{\mathrm{hor}}\cdot\Gamma_0)_{q}=4$.
\begin{itemize}    
\item[$(4\text{--}1)$] $B_{\mathrm{hor}} = B_1 + B_2$.  
$\mathrm{mult}_q(B_{\mathrm{hor}})=2$ and $(B_i\cdot \Gamma_0)_q=2$ for $i=1,2$.
If we blow up twice at the singular point of $B_{\mathrm{hor}}$, that is, the intersection point of $B_1$ and $B_2$, then the proper transform of $B_{\mathrm{hor}}$ has an $A_{k-1}$-singularity for some non-negative even integer $k$.
This $A_{k-1}$-singularity is not contact to the exceptional curve.

\item[$(4\text{--}2)$] $B_{\mathrm{hor}} = B_1$.
$\mathrm{mult}_q(B_{\mathrm{hor}})=2$ and $(B_1\cdot \Gamma_0)_q$=4.  
If we blow up twice at the singular point of $B_1$, then the proper transform of $B_1$ has an $A_{k-1}$-singularity for 
some positive odd integer $k$.
This $A_{k-1}$-singularity is not contact to the exceptional curve.

\end{itemize}
    
    \item[$(4)^*$] $\Gamma_0  \subset B$.
    There exists a singular points $q \in B_{\mathrm{hor}}$ with $(B_{\mathrm{hor}}\cdot\Gamma_0)_{q}=4$.
    As in the case of $(4)$, each $B_{\mathrm{hor}}$ is classified into two types $(4\text{--}1)^*$ and $(4\text{--}2)^*$.
\end{itemize}

We classify the fiber germs of elliptic surfaces arising from the above branches and their involutions.  
The classification is distinguished by the following data, and we also provide an explanation of the symbols used in the figures.

\begin{itemize}
    \item Irreducible components of $F_0$ that are exchanged by the involution.  
    In the figure,red lines represent irreducible components of $F_0$ fixed by $\sigma$.
    \item Intersection points between the horizontal $1$-dimensional $\sigma$-fixed locus and $F_0$.  
    Solid red points in the figure represent such intersection points.
    \item Isolated $\sigma$-fixed points that always appear on $F_0$.  
    Red circles in the figure represent these isolated $\sigma$-fixed points.
\end{itemize}

\begin{prop}\label{prop_inv_ell}
Let $f:Y \to T$ be a germ of an elliptic surface, and let $\sigma \in \mathrm{Aut}_T(Y)$ be an involution satisfying the condition $(\star)$.  
Then the central fiber $F_0$ and the involution $\sigma$ are one of the following.

\begin{itemize}
    \item[$(1\text{-}1)$]
    $F_0$ is of type $\mathrm{II}$. The involution in this case is as in the Table.
    \item[$(1\text{-}2)$] 
     $F_0$ is of type $\mathrm{III}$. The involution interchanges the two irreducible components. The configuration of the fixed locus is depicted in the figure in the table.
    \item[$(1\text{-}1)^*$]
    $F_0$ is of type $\mathrm{IV}^*$. The involution interchanges the two branches of length $2$.
    The configuration of the fixed locus is depicted in the figure in the table.
    \item[$(1\text{-}2)^*$]
    $F_0$ is of type $\mathrm{III}^*$. 
    The involution interchanges the two branches of length $3$.
    The configuration of the fixed locus is depicted in the figure in the table.
    \item[$(2\text{-}1)$]
    $F_0$ is of type $\mathrm{I}_{k_1+k_2}$. 
    The fixed locus is depicted in the figure in the table.
    The involution fixes the four red points in the figure and acts as a reflection that exchanges the remaining components symmetrically.
    \item[$(2\text{-}2)$]
    $F_0$ is of type $\mathrm{I}_{k_1+k_2}$. 
    The fixed locus is depicted in the figure in the table.
    The involution fixes the three red points in the figure and acts as a reflection that exchanges the remaining components symmetrically.
     \item[$(2\text{-}3)$]
    $F_0$ is of type $\mathrm{I}_{k_1+k_2}$. 
    The fixed locus is depicted in the figure in the table.
    The involution fixes the two red points in the figure and acts as a reflection that exchanges the remaining components symmetrically.
    \item[$(2\text{-}1)^*$]
    $F_0$ is of type $\mathrm{I}_{k_1+k_2}$. 
    The fixed locus is depicted in the figure in the table.
    \item[$(2\text{-}2)^*$]
    $F_0$ is of type $\mathrm{I}_{k_1+k_2}$. 
    The fixed locus is depicted in the figure in the table.
    The involution interchanges the two branches of length $1$.
    \item[$(2\text{-}3)^*$]
    $F_0$ is of type $\mathrm{I}_{k_1+k_2}$. 
    The fixed locus is depicted in the figure in the table.
    The involution exchanges the two branches located on either side of the figure.
    \item[$(3\text{-}1)$]
    $F_0$ is of type $\mathrm{III}$. 
    The fixed locus is depicted in the figure in the table.
    \item[$(3\text{-}2)$]
    $F_0$ is of type $\mathrm{IV}$. 
    The fixed locus is depicted in the figure in the table.
    The involution exchanges the two components.
    \item[$(3\text{-}1)^*$]
    $F_0$ is of type $\mathrm{III}^*$. 
    The fixed locus is depicted in the figure in the table.
      The involution interchanges the two branches of length $3$.
    \item[$(3\text{-}2)^*$]
    $F_0$ is of type $\mathrm{II}^*$. 
    The fixed locus is depicted in the figure in the table.
    \item[$(4\text{-}1)$]
    $F_0$ is of type $\mathrm{I}^*_{k}$. 
    The fixed locus is depicted in the figure in the table.
    The involution fixes the two red points in the figure and acts as a reflection that exchanges the remaining components symmetrically.
    \item[$(4\text{-}2)$]
    $F_0$ is of type $\mathrm{I}^*_{k}$. 
    The fixed locus is depicted in the figure in the table.
    The involution fixes the the red point in the figure and acts as a reflection that exchanges the remaining components symmetrically.
    \item[$(4\text{-}1)^*$]
    $F_0$ is of type $2\mathrm{I}_{k}$. 
    The fixed locus is depicted in the figure in the table.
    The involution fixes the two red points and two red circles in the figure and acts as a reflection that exchanges the remaining components symmetrically.
    \item[$(4\text{-}2)^*$]
    $F_0$ is of type $2\mathrm{I}_{k}$. 
    The fixed locus is depicted in the figure in the table.
    The involution fixes the red point and two red circles in the figure and acts as a reflection that exchanges the remaining components symmetrically.
\end{itemize}

\begin{table}[H]
  \begin{tabular}{|c|c|c|} 
    \hline  & branch locus $B$ & fiber $F$ and its involution  \\ \hline
$(1\text{--}1)$  &
\begin{minipage}{70mm}
\begin{center}
\begin{tikzpicture}[line cap=round,line join=round,>=triangle 45,x=1cm,y=1cm]
\clip(-3,-1) rectangle (3,1);
\draw [line width=1pt] (2.5,0)-- (-2.5,0);
\draw[line width=1pt,color=black,smooth,samples=100,domain=0.5:2.5,red] plot(\x,{(\x - 1.5)^2});
\draw [line width=1pt,red] (-1.5,-1)-- (-1.5,1);
\begin{scriptsize}
\draw (2.75,0) node {$\Gamma_0$};
\draw (-1.25,-0.25) node {$B_1$};
\draw (1.5,-0.25) node {$B_2$};
\end{scriptsize}
\end{tikzpicture}
\end{center}
\end{minipage} &
\begin{minipage}{50mm}
\begin{center}
\begin{tikzpicture}[line cap=round,line join=round,>=triangle 45,x=1cm,y=1cm]
\clip(-1,-1) rectangle (1,1);
\draw[line width=1pt,smooth,samples=100,domain=0.0000051148379696432314:4.09933597466865] plot(\x,{sqrt((\x)^(5))});
\draw[line width=1pt,smooth,samples=100,domain=0.0000051148379696432314:4.09933597466865] plot(\x,{0-sqrt((\x)^(5))});
\draw [fill=red] (0,0) circle (2.5pt);
\draw [fill=red] (0.69,0.4) circle (2.5pt);
\end{tikzpicture}
\end{center}
\end{minipage} \\ \hline
$(1\text{--}2)$ &
\begin{minipage}{70mm}
\begin{center}
\begin{tikzpicture}[line cap=round,line join=round,>=triangle 45,x=1cm,y=1cm]
\clip(-3,-1) rectangle (3,1);
\draw [line width=1pt] (2.5,0)-- (-2.5,0);
\draw[line width=1pt,color=black,smooth,samples=100,domain=-1:1,red] plot(\x,{(\x)^2});
\begin{scriptsize}
\draw (2.75,0) node {$\Gamma_0$};
\draw (0,-0.25) node {$B_1$};
\end{scriptsize}
\end{tikzpicture}
\end{center}
\end{minipage} &
\begin{minipage}{70mm}
\begin{center}
\begin{tikzpicture}[line cap=round,line join=round,>=triangle 45,x=0.75cm,y=0.75cm]
\clip(-3,-1.5) rectangle (3,1.5);
\draw [shift={(0,-2)},line width=1pt]  plot[domain=-1.5707963267948966:1.5707963267948966,variable=\t]({1*2*sin(\t r)+0*2*cos(\t r)},{0*2*sin(\t r)+1*2*cos(\t r)});
\draw [shift={(0,2)},line width=1pt]  plot[domain=1.5707963267948966:4.71238898038469,variable=\t]({1*2*sin(\t r)+0*2*cos(\t r)},{0*2*sin(\t r)+1*2*cos(\t r)});
\draw[<->] [shift={(-2,0)},line width=1pt]  plot[domain=3.141592653589793:6.283185307179586,variable=\t]({1*1*sin(\t r)+0*1*cos(\t r)},{0*1*sin(\t r)+1*1*cos(\t r)});
\begin{scriptsize}
\draw [fill=red] (0,0) circle (2.5pt);
\end{scriptsize}
\end{tikzpicture}
\end{center}
\end{minipage} \\ \hline
$(1 \text{--} 1)^*$ &
\begin{minipage}{70mm}
\begin{center}
\begin{tikzpicture}[line cap=round,line join=round,>=triangle 45,x=1cm,y=1cm]
\clip(-3,-1) rectangle (3,1);
\draw [line width=1pt,red] (2.5,0)-- (-2.5,0);
\draw[line width=1pt,color=black,smooth,samples=100,domain=0.5:2.5,red] plot(\x,{(\x - 1.5)^2});
\draw [line width=1pt,red] (-1.5,-1)-- (-1.5,1);
\begin{scriptsize}
\draw (2.75,0) node {$\Gamma_0$};
\draw (-1.25,-0.25) node {$B_1$};
\draw (1.5,-0.25) node {$B_2$};
\end{scriptsize}
\end{tikzpicture}
\end{center}
\end{minipage} &
\begin{minipage}{70mm}
\begin{center}
\begin{tikzpicture}[line cap=round,line join=round,>=triangle 45,x=0.75cm,y=0.75cm]
\clip(-4,-1) rectangle (5,2);
\draw [line width=1pt] (-3,1.5)-- (0,0);
\draw [line width=1pt,red] (-1,0)-- (1,1);
\draw [line width=1pt] (0,1)-- (3,-0.5);
\draw [line width=1pt] (1,0)-- (3,1);
\draw [line width=1pt] (2,-0.5)-- (4,0.5);
\draw [line width=1pt] (3.2439920935869733,0.380011184692855)-- (4,0);
\draw [line width=1pt] (2,1)-- (2.7323131647256136,0.6320321496544213);
\draw[<->] [shift={(3.9944114602591707,1.0889270716418986)},line width=1pt]  plot[domain=-0.6936192170105402:2.447973436579253,variable=\t]({1*0.6429802831071187*cos(\t r)+0*0.6429802831071187*sin(\t r)},{0*0.6429802831071187*cos(\t r)+1*0.6429802831071187*sin(\t r)});
\begin{scriptsize}
\draw [fill=red] (2,0) circle (2.5pt);
\draw [fill=red] (-2,1) circle (2.5pt);
\end{scriptsize}
\end{tikzpicture}
\end{center}
\end{minipage} \\ \hline
$(1 \text{--}2)^*$ &
\begin{minipage}{70mm}
\begin{center}
\begin{tikzpicture}[line cap=round,line join=round,>=triangle 45,x=1cm,y=1cm]
\clip(-3,-1) rectangle (3,1);
\draw [line width=1pt, red] (2.5,0)-- (-2.5,0);
\draw[line width=1pt,color=black,smooth,samples=100,domain=-1:1,red] plot(\x,{(\x)^2});
\begin{scriptsize}
\draw (2.75,0) node {$\Gamma_0$};
\draw (0,-0.25) node {$B_1$};
\end{scriptsize}
\end{tikzpicture}
\end{center}
\end{minipage} &
\begin{minipage}{70mm}
\begin{center}
\begin{tikzpicture}[line cap=round,line join=round,>=triangle 45,x=0.85cm,y=0.75cm]
\clip(-1.,-1.5) rectangle (6.5,1.5);
\draw [line width=1pt] (-1,1)-- (0.5,-0.5);
\draw [line width=1pt ] (0,-0.5)-- (1.5,1);
\draw [line width=1pt] (1,1)-- (2.5,-0.5);
\draw [line width=1pt] (1.9949545971809697,-0.3144062445991428)-- (3.4968674700510833,-0.3144062445991428);
\draw [line width=1pt] (3,-0.5)-- (4.5,1);
\draw [line width=1pt] (4,1)-- (5.5,-0.5);
\draw [line width=1pt] (5,-0.5)-- (6.5,1);
\draw [line width=1pt,red] (2.738643205403547,0.9783762368802792)-- (2.747547123946322,-0.4996742412203729);
\draw[<->] [shift={(2.75,-0.5)},line width=1pt]  plot[domain=3.141592653589793:6.283185307179586,variable=\t]({1*0.75*cos(\t r)+0*0.75*sin(\t r)},{0*0.75*cos(\t r)+1*0.75*sin(\t r)});
\begin{scriptsize}
\draw [fill=red] (2.5504089322370063,-0.3144062445991428) circle (2.5pt);
\end{scriptsize}
\end{tikzpicture}
\end{center}
\end{minipage} \\ \hline
\end{tabular}
\end{table}

\begin{table}[H]
  \begin{tabular}{|c|c|c|} 
    \hline  & branch locus & fiber $F$ and its involution \\ \hline
 $(2\text{--}1)$   &
    \begin{minipage}{70mm}
\begin{center}
\begin{tikzpicture}[line cap=round,line join=round,>=triangle 45,x=1cm,y=1cm]
\clip(-3,-1) rectangle (3,1);
\draw [line width=1pt] (-2.5,0)--(2.5,0);
\draw[line width=1pt,color=red,smooth,samples=100,domain=1.5:2.5] plot(\x,{sqrt((\x)-1.5)});
\draw[line width=1pt,color=red,smooth,samples=100,domain=1.5:2.5] plot(\x,{0-sqrt((\x)-1.5)});
\draw[line width=1pt,color=red,smooth,samples=100,domain=0.5:1.5] plot(\x,{sqrt(-((\x)-1.5))});
\draw[line width=1pt,color=red,smooth,samples=100,domain=0.5:1.5] plot(\x,{0-sqrt(-((\x)-1.5))});
\draw[line width=1pt,color=red,smooth,samples=100,domain=-1.5:-0.5] plot(\x,{sqrt((\x)+1.5)});
\draw[line width=1pt,color=red,smooth,samples=100,domain=-2.5:-1.5] plot(\x,{sqrt(-((\x)+1.5))});
\draw[line width=1pt,color=red,smooth,samples=100,domain=-1.5:-0.5] plot(\x,{0-sqrt((\x)+1.5)});
\draw[line width=1pt,color=red,smooth,samples=100,domain=-2.5:-1.5] plot(\x,{0-sqrt(-((\x)+1.5))});
\begin{scriptsize}
\draw (2.75,0) node {$\Gamma_0$};
\draw (-2,-0.25) node {$B_1$};
\draw (-1,-0.25) node {$B_2$};
\draw (1,-0.25) node {$B_3$};
\draw (2,-0.25) node {$B_4$};
\end{scriptsize}
\end{tikzpicture}
\end{center}
    \end{minipage} &
    \begin{minipage}{70mm}
\begin{center}
\begin{tikzpicture}[line cap=round,line join=round,>=triangle 45,x=1cm,y=0.375cm]
\clip(-3,-3.5) rectangle (3,3.5);
\draw [line width=1pt] (-1.396346686723849,2.9902261522629088)-- (1.3812368566545277,2.9902261522629088);
\draw [line width=1pt] (-0.7920188864908967,3.199416544651239)-- (-2.1982431908791127,1.8048139287290361); 
\draw [line width=1pt] (0.7885307448875938,3.199416544651239)-- (2.171511672343773,1.8396789941270912); 
\draw [line width=1pt] (-2.000674486956801,2.385898352029954)-- (-2.000674486956801,0.631023393661182); 
\draw [line width=1pt] (1.9971863453534984,2.385898352029954)-- (1.9971863453534984,0.6077800167291453);
\draw [line width=1pt] (-2.000674486956801,-0.6124972722027823)-- (-1.9774311100247646,-2.402237295969609);
\draw [line width=1pt] (-2.209864879345131,-1.7979094957366546)-- (-0.8152622634229333,-3.204133800124876);
\draw [line width=1pt] (-1.396346686723849,-3.006565096202564)-- (1.4044802335865643,-3.006565096202564);
\draw [line width=1pt] (0.8001524333536121,-3.1925121116588575)-- (2.2063767377418277,-1.7979094957366546);
\draw [line width=1pt] (2.0088080338195167,-2.4138589844356275)-- (2.0088080338195167,-0.6241189606688006);
\draw[<->] [line width=1pt] (-1.3847249982578307,0)-- (1.392858545120546,0);
\begin{scriptsize}
\draw [fill=red] (0.5,3) circle (2.5pt);
\draw [fill=red] (-0.5,3) circle (2.5pt);
\draw [fill=red] (0.5,-3) circle (2.5pt);
\draw [fill=red] (-0.5,-3) circle (2.5pt);
\draw [color=black] (-2,0) node {$\vdots$};
\draw [color=black] (2,0) node {$\vdots$};
\end{scriptsize}
\end{tikzpicture}
\end{center}
 \end{minipage} \\ \hline
 $(2\text{--}2)$ &
\begin{minipage}{70mm}
\begin{center}
\begin{tikzpicture}[line cap=round,line join=round,>=triangle 45,x=1cm,y=1cm]
\clip(-3,-1) rectangle (3,1);
\draw [line width=1pt] (2.5,0)-- (-2.5,0);
\draw[line width=1pt,color=red,smooth,samples=100,domain=1.5:2.5] plot(\x,{sqrt((\x)-1.5)});
\draw[line width=1pt,color=red,smooth,samples=100,domain=1.5:2.5] plot(\x,{0-sqrt((\x)-1.5)});
\draw[line width=1pt,color=red,smooth,samples=100,domain=0.5:1.5] plot(\x,{sqrt(-((\x)-1.5))});
\draw[line width=1pt,color=red,smooth,samples=100,domain=0.5:1.5] plot(\x,{0-sqrt(-((\x)-1.5))});
\draw[line width=1pt,color=red,smooth,samples=100,domain=-1.5:-0.5] plot(\x,{sqrt((\x)+1.5)});
\draw[line width=1pt,color=red,smooth,samples=100,domain=-2.5:-1.5] plot(\x,{sqrt(-((\x)+1.5))});
\begin{scriptsize}
\draw (2.75,0) node {$\Gamma_0$};
\draw (-1.5,-0.25) node {$B_1$};
\draw (1,-0.25) node {$B_2$};
\draw (2,-0.25) node {$B_3$};
\end{scriptsize}
\end{tikzpicture}
\end{center}
\end{minipage} &
\begin{minipage}{70mm}
\begin{center}
\begin{tikzpicture}[line cap=round,line join=round,>=triangle 45,x=0.75cm,y=0.375cm]
\clip(-3,-3.5) rectangle (3,3.5);
\draw [line width=1pt] (0.4034053251062811,3.192327577892131)-- (-2.1982431908791127,1.8048139287290361);
\draw [line width=1pt] (-0.4027259424068129,3.2029345682541455)-- (2.171511672343773,1.8396789941270912);
\draw [line width=1pt] (-2.000674486956801,2.385898352029954)-- (-2.000674486956801,0.631023393661182);
\draw [line width=1pt] (1.9971863453534984,2.385898352029954)-- (1.9971863453534984,0.6077800167291453);
\draw [line width=1pt] (-2.000674486956801,-0.6124972722027823)-- (-1.9774311100247646,-2.402237295969609);
\draw [line width=1pt] (-2.209864879345131,-1.7979094957366546)-- (-0.8152622634229333,-3.204133800124876);
\draw [line width=1pt] (-1.396346686723849,-3.006565096202564)-- (1.4044802335865643,-3.006565096202564);
\draw [line width=1pt] (0.8001524333536121,-3.1925121116588575)-- (2.2063767377418277,-1.7979094957366546);
\draw [line width=1pt] (2.0088080338195167,-2.4138589844356275)-- (2.0088080338195167,-0.6241189606688006);
\draw[<->] [line width=1pt] (-1.3847249982578307,0)-- (1.392858545120546,0);
\begin{scriptsize}
\draw [fill=red] (0,3) circle (2.5pt);
\draw [fill=red] (0.5,-3) circle (2.5pt);
\draw [fill=red] (-0.5,-3) circle (2.5pt);
\draw [color=black] (-2,0) node {$\vdots$};
\draw [color=black] (2,0) node {$\vdots$};
\end{scriptsize}
\end{tikzpicture}
\end{center}
\end{minipage} \\ \hline
$(2\text{--}3)$ & 
\begin{minipage}{70mm}
\begin{center}
\begin{tikzpicture}[line cap=round,line join=round,>=triangle 45,x=1cm,y=1cm]
\clip(-3,-1) rectangle (3,1);
\draw [line width=1pt] (2.5,0)-- (-2.5,0);
\draw[line width=1pt,color=red,smooth,samples=100,domain=1.5:2.5] plot(\x,{sqrt((\x)-1.5)});
\draw[line width=1pt,color=red,smooth,samples=100,domain=0.5:1.5] plot(\x,{sqrt(-((\x)-1.5))});
\draw[line width=1pt,color=red,smooth,samples=100,domain=-1.5:-0.5] plot(\x,{sqrt((\x)+1.5)});
\draw[line width=1pt,color=red,smooth,samples=100,domain=-2.5:-1.5] plot(\x,{sqrt(-((\x)+1.5))});
\begin{scriptsize}
\draw (2.75,0) node {$\Gamma_0$};
\draw (-1.5,-0.25) node {$B_1$};
\draw (1.5,-0.25) node {$B_2$};
\end{scriptsize}
\end{tikzpicture}
\end{center}
\end{minipage} &
\begin{minipage}{70mm}
\begin{center}
\begin{tikzpicture}[line cap=round,line join=round,>=triangle 45,x=0.75cm,y=0.375cm]
\clip(-3,-3.5) rectangle (3,3.5);
\draw [line width=1pt] (0.4034053251062811,3.192327577892131)-- (-2.1982431908791127,1.8048139287290361);
\draw [line width=1pt] (-0.4027259424068129,3.2029345682541455)-- (2.171511672343773,1.8396789941270912);
\draw [line width=1pt] (-2.000674486956801,2.385898352029954)-- (-2.000674486956801,0.631023393661182);
\draw [line width=1pt] (1.9971863453534984,2.385898352029954)-- (1.9971863453534984,0.6077800167291453);
\draw [line width=1pt] (-2.000674486956801,-0.6124972722027823)-- (-1.9774311100247646,-2.402237295969609);
\draw [line width=1pt] (-2.209864879345131,-1.7979094957366546)-- (0.38960179286642116,-3.1824138618571163);
\draw [line width=1pt] (-0.39733263152312553,-3.1967217604823808)-- (2.2063767377418277,-1.7979094957366546);
\draw [line width=1pt] (2.0088080338195167,-2.4138589844356275)-- (2.0088080338195167,-0.6241189606688006);
\draw[<->] [line width=1pt] (-1.3847249982578307,0)-- (1.392858545120546,0);
\begin{scriptsize}
\draw [fill=red] (0,3) circle (2.5pt);
\draw [fill=red] (0,-3) circle (2.5pt);
\draw [color=black] (-2,0) node {$\vdots$};
\draw [color=black] (2,0) node {$\vdots$};
\end{scriptsize}
\end{tikzpicture}
\end{center}
\end{minipage} \\ \hline
$(2\text{--}1)^*$  &
\begin{minipage}{70mm}
\begin{center}
\begin{tikzpicture}[line cap=round,line join=round,>=triangle 45,x=1cm,y=1cm]
\clip(-3,-1) rectangle (3,1);
\draw [line width=1pt,red] (-2.5,0)--(2.5,0);
\draw[line width=1pt,color=red,smooth,samples=100,domain=1.5:2.5] plot(\x,{sqrt((\x)-1.5)});
\draw[line width=1pt,color=red,smooth,samples=100,domain=1.5:2.5] plot(\x,{0-sqrt((\x)-1.5)});
\draw[line width=1pt,color=red,smooth,samples=100,domain=0.5:1.5] plot(\x,{sqrt(-((\x)-1.5))});
\draw[line width=1pt,color=red,smooth,samples=100,domain=0.5:1.5] plot(\x,{0-sqrt(-((\x)-1.5))});
\draw[line width=1pt,color=red,smooth,samples=100,domain=-1.5:-0.5] plot(\x,{sqrt((\x)+1.5)});
\draw[line width=1pt,color=red,smooth,samples=100,domain=-2.5:-1.5] plot(\x,{sqrt(-((\x)+1.5))});
\draw[line width=1pt,color=red,smooth,samples=100,domain=-1.5:-0.5] plot(\x,{0-sqrt((\x)+1.5)});
\draw[line width=1pt,color=red,smooth,samples=100,domain=-2.5:-1.5] plot(\x,{0-sqrt(-((\x)+1.5))});
\begin{scriptsize}
\draw (2.75,0) node {$\Gamma_0$};
\draw (-2,-0.25) node {$B_1$};
\draw (-1,-0.25) node {$B_2$};
\draw (1,-0.25) node {$B_3$};
\draw (2,-0.25) node {$B_4$};
\end{scriptsize}
\end{tikzpicture}
\end{center}
\end{minipage} &
\begin{minipage}{70mm}
\begin{center}
\begin{tikzpicture}[line cap=round,line join=round,>=triangle 45,x=0.75cm,y=1cm]
\clip(-2,-1.5) rectangle (4,1.5);
\draw [line width=1pt,red] (-1,1)-- (-1,-1);
\draw [line width=1pt] (-0.796405425080071,0.5940678833808372)-- (-2,1);
\draw [line width=1pt] (-0.8037318293379739,-0.5928096063994598)-- (-2,-1);
\draw [line width=1pt] (-1.2408739500595178,0)-- (0,0.5);
\draw [line width=1pt,red] (3,1)-- (3,-1);
\draw [line width=1pt] (2.7935326612923945,0.6013942876387403)-- (4,1);
\draw [line width=1pt] (2.800859065550297,-0.6001360106573629)-- (4,-1);
\draw [line width=1pt,red] (-0.5,0.5)-- (0.705507447790042,0.0006291384906886806);
\draw [line width=1pt,red] (1.3062725969380873,0)-- (2.5,0.5);
\draw [line width=1pt] (2,0.5)-- (3.5,0);
\draw[<->] [shift={(4.5,0)},line width=1pt]  plot[domain=-1.5707963267948966:1.5707963267948966,variable=\t]({1*1*cos(\t r)+0*1*sin(\t r)},{0*1*cos(\t r)+1*1*sin(\t r)});
\begin{scriptsize}
\draw [fill=black] (1,0) node {$\cdots$};
\draw [fill=red] (3.504089322370063,0.8361556731937827) circle (2.5pt);
\draw [fill=red] (3.5217486324743974,-0.8405229158363116) circle (2.5pt);
\draw [fill=red] (-1.613562300367173,-0.8711874334557244) circle (2.5pt);
\draw [fill=red] (-1.613562300367173,0.8711874334557244) circle (2.5pt);
\end{scriptsize}
\end{tikzpicture}
\end{center}
\end{minipage} \\ \hline
$(2\text{--}2)^*$ &
\begin{minipage}{70mm}
\begin{center}
\begin{tikzpicture}[line cap=round,line join=round,>=triangle 45,x=1cm,y=1cm]
\clip(-3,-1) rectangle (3,1);
\draw [line width=1pt,red] (2.5,0)-- (-2.5,0);
\draw[line width=1pt,color=red,smooth,samples=100,domain=1.5:2.5] plot(\x,{sqrt((\x)-1.5)});
\draw[line width=1pt,color=red,smooth,samples=100,domain=1.5:2.5] plot(\x,{0-sqrt((\x)-1.5)});
\draw[line width=1pt,color=red,smooth,samples=100,domain=0.5:1.5] plot(\x,{sqrt(-((\x)-1.5))});
\draw[line width=1pt,color=red,smooth,samples=100,domain=0.5:1.5] plot(\x,{0-sqrt(-((\x)-1.5))});
\draw[line width=1pt,color=red,smooth,samples=100,domain=-1.5:-0.5] plot(\x,{sqrt((\x)+1.5)});
\draw[line width=1pt,color=red,smooth,samples=100,domain=-2.5:-1.5] plot(\x,{sqrt(-((\x)+1.5))});
\begin{scriptsize}
\draw (2.75,0) node {$\Gamma_0$};
\draw (-1.5,-0.25) node {$B_1$};
\draw (1,-0.25) node {$B_2$};
\draw (2,-0.25) node {$B_3$};
\end{scriptsize}
\end{tikzpicture}
\end{center}
\end{minipage} &
\begin{minipage}{70mm}
\begin{center}
\begin{tikzpicture}[line cap=round,line join=round,>=triangle 45,x=0.75cm,y=0.75cm]
\clip(-4,-2) rectangle (4.5,2);
\draw [line width=1pt] (-1,1)-- (-1,-1);
\draw [line width=1pt] (-0.796405425080071,0.5940678833808372)-- (-2,1);
\draw [line width=1pt] (-0.8037318293379739,-0.5928096063994598)-- (-2,-1);
\draw [line width=1pt,red] (-1.2408739500595178,0)-- (0,0.5);
\draw [line width=1pt,red] (3,1)-- (3,-1);
\draw [line width=1pt] (2.7935326612923945,0.6013942876387403)-- (4,1);
\draw [line width=1pt] (2.800859065550297,-0.6001360106573629)-- (4,-1);
\draw [line width=1pt] (-0.5,0.5)-- (0.705507447790042,0.0006291384906886806);
\draw [line width=1pt,red] (1.3062725969380873,0)-- (2.5,0.5);
\draw [line width=1pt] (2,0.5)-- (3.5,0);
\draw[<->] [shift={(-2.5,0)},line width=1pt]  plot[domain=1.5707963267948966:4.71238898038469,variable=\t]({1*1*cos(\t r)+0*1*sin(\t r)},{0*1*cos(\t r)+1*1*sin(\t r)});
\begin{scriptsize}
\draw [fill=black] (1,0) node {$\cdots$};
\draw [fill=red] (3.504089322370063,0.8361556731937827) circle (2.5pt);
\draw [fill=red] (3.5217486324743974,-0.8405229158363116) circle (2.5pt);
\draw [fill=red] (-1,-0.35) circle (2.5pt);
\end{scriptsize}
\end{tikzpicture}
\end{center}
\end{minipage} \\ \hline
$(2 \text{--} 3)^*$ &
\begin{minipage}{70mm}
\begin{center}
\begin{tikzpicture}[line cap=round,line join=round,>=triangle 45,x=1cm,y=1cm]
\clip(-3,-1) rectangle (3,1);
\draw [line width=1pt,red] (2.5,0)-- (-2.5,0);
\draw[line width=1pt,color=red,smooth,samples=100,domain=1.5:2.5] plot(\x,{sqrt((\x)-1.5)});
\draw[line width=1pt,color=red,smooth,samples=100,domain=0.5:1.5] plot(\x,{sqrt(-((\x)-1.5))});
\draw[line width=1pt,color=red,smooth,samples=100,domain=-1.5:-0.5] plot(\x,{sqrt((\x)+1.5)});
\draw[line width=1pt,color=red,smooth,samples=100,domain=-2.5:-1.5] plot(\x,{sqrt(-((\x)+1.5))});
\begin{scriptsize}
\draw (2.75,0) node {$\Gamma_0$};
\draw (-1.5,-0.25) node {$B_1$};
\draw (1.5,-0.25) node {$B_2$};
\end{scriptsize}
\end{tikzpicture}
\end{center}
\end{minipage} &
\begin{minipage}{75mm}
\begin{center}
\begin{tikzpicture}[line cap=round,line join=round,>=triangle 45,x=0.75cm,y=0.75cm]
\clip(-4,-2) rectangle (6,2);
\draw [line width=1pt] (-1,1)-- (-1,-1);
\draw [line width=1pt] (-0.796405425080071,0.5940678833808372)-- (-2,1);
\draw [line width=1pt] (-0.8037318293379739,-0.5928096063994598)-- (-2,-1);
\draw [line width=1pt,red] (-1.2408739500595178,0)-- (0,0.5);
\draw [line width=1pt] (3,1)-- (3,-1);
\draw [line width=1pt] (2.7935326612923945,0.6013942876387403)-- (4,1);
\draw [line width=1pt] (2.800859065550297,-0.6001360106573629)-- (4,-1);
\draw [line width=1pt] (-0.5,0.5)-- (0.705507447790042,0.0006291384906886806);
\draw [line width=1pt] (1.3062725969380873,0)-- (2.5,0.5);
\draw [line width=1pt,red] (2,0.5)-- (3.5,0);
\draw[<->] [shift={(4.5,0)},line width=1pt]  plot[domain=-1.5707963267948966:1.5707963267948966,variable=\t]({1*1*cos(\t r)+0*1*sin(\t r)},{0*1*cos(\t r)+1*1*sin(\t r)});
\draw[<->] [shift={(-2.5,0)},line width=1pt]  plot[domain=1.5707963267948966:4.71238898038469,variable=\t]({1*1*cos(\t r)+0*1*sin(\t r)},{0*1*cos(\t r)+1*1*sin(\t r)});
\begin{scriptsize}
\draw [fill=black] (1,0) node {$\cdots$};
\draw [fill=red] (-1,-0.35) circle (2.5pt);
\draw [fill=red] (3,-0.35) circle (2.5pt);
\end{scriptsize}
\end{tikzpicture}
\end{center}
\end{minipage} \\ \hline
\end{tabular}
\end{table}

\newpage

\begin{table}[H]
  \begin{tabular}{|c|c|c|} 
    \hline  & branch locus & fiber $F$ and its involution \\ \hline
 $(3\text{--}1)$   &
    \begin{minipage}{70mm}
\begin{center}
\begin{tikzpicture}[line cap=round,line join=round,>=triangle 45,x=1cm,y=1cm]
\clip(-3,-1) rectangle (3,1);
\draw [line width=1pt, red] (-1.5,-1)-- (-1.5,1);
\draw [line width=1pt, red] (1.5,-1)-- (1.5,1);
\draw [line width=1pt] (2.5,0)-- (-2.5,0);
\draw[line width=1pt,color=black,smooth,samples=100,domain=0.5:2.5,red] plot(\x,{(\x - 1.5)^2});
\begin{scriptsize}
\draw (2.75,0) node {$\Gamma_0$};
\draw (-1.25,-0.25) node {$B_1$};
\draw (1.75,-0.25) node {$B_2$};
\draw (2,0.75) node {$B_3$};
\end{scriptsize}
\end{tikzpicture}
\end{center}
    \end{minipage} &
    \begin{minipage}{70mm}
    \begin{center}
\begin{tikzpicture}[line cap=round,line join=round,>=triangle 45,x=1cm,y=0.75cm]
\clip(-3,-1.5) rectangle (3,1.5);
\draw [shift={(0,-2)},line width=1pt]  plot[domain=-1.5707963267948966:1.5707963267948966,variable=\t]({1*2*sin(\t r)+0*2*cos(\t r)},{0*2*sin(\t r)+1*2*cos(\t r)});
\draw [shift={(0,2)},line width=1pt]  plot[domain=1.5707963267948966:4.71238898038469,variable=\t]({1*2*sin(\t r)+0*2*cos(\t r)},{0*2*sin(\t r)+1*2*cos(\t r)});
\begin{scriptsize}
\draw [fill=red] (0,0) circle (2.5pt);
\draw [fill=red] (1.6,0.8) circle (2.5pt);
\draw [fill=red] (1.6,-0.8) circle (2.5pt);
\end{scriptsize}
\end{tikzpicture}
   \end{center}
 \end{minipage} \\ \hline
 $(3\text{--}2)$ &
\begin{minipage}{70mm}
\begin{center}
\begin{tikzpicture}[line cap=round,line join=round,>=triangle 45,x=1cm,y=1cm]
\clip(-3,-1) rectangle (3,1);
\draw [line width=1pt, red] (-1.5,-1)-- (-1.5,1);
\draw [line width=1pt] (2.5,0)-- (-2.5,0);
\draw[line width=1pt,color=black,smooth,samples=100,domain=1:2,red] plot(\x,{(\x - 1)^2});
\draw[line width=1pt,color=black,smooth,samples=100,domain=1:2,red] plot(\x,{0-(\x - 1)^2});
\begin{scriptsize}
\draw (2.75,0) node {$\Gamma_0$};
\draw (-1.25,-0.25) node {$B_1$};
\draw (1.25,-0.5) node {$B_2$};
\end{scriptsize}
\end{tikzpicture}
\end{center}
\end{minipage} &
\begin{minipage}{70mm}
\begin{center}
\begin{tikzpicture}[line cap=round,line join=round,>=triangle 45,x=1cm,y=0.5cm]
\clip(-3,-2.9871467618634933) rectangle (3,3);
\draw [line width=1pt] (-2,2)-- (2,-2);
\draw [line width=1pt] (2,2)-- (-2,-2);
\draw [line width=1pt] (0,2.5)-- (0,-2.5);
\draw[<->] [shift={(-2,0)},line width=1pt]  plot[domain=3.141592653589793:6.283185307179586,variable=\t]({1*1*sin(\t r)+0*1*cos(\t r)},{0*1*sin(\t r)+1*1*cos(\t r)});
\begin{scriptsize}
\draw [fill=red] (0,0) circle (2.5pt);
\draw [fill=red] (0,2) circle (2.5pt);
\end{scriptsize}
\end{tikzpicture}
\end{center}
\end{minipage} \\ \hline
$(3\text{--}1)^*$ & 
\begin{minipage}{70mm}
\begin{center}
\begin{tikzpicture}[line cap=round,line join=round,>=triangle 45,x=1cm,y=1cm]
\clip(-3,-1) rectangle (3,1);
\draw [line width=1pt, red] (-1.5,-1)-- (-1.5,1);
\draw [line width=1pt, red] (1.5,-1)-- (1.5,1);
\draw [line width=1pt, red] (2.5,0)-- (-2.5,0);
\draw[line width=1pt,color=black,smooth,samples=100,domain=0.5:2.5,red] plot(\x,{(\x - 1.5)^2});
\begin{scriptsize}
\draw (2.75,0) node {$\Gamma_0$};
\draw (-1.25,-0.25) node {$B_1$};
\draw (1.75,-0.25) node {$B_2$};
\draw (2,0.75) node {$B_3$};
\end{scriptsize}
\end{tikzpicture}
\end{center}
\end{minipage} &
\begin{minipage}{70mm}
\begin{center}
\begin{tikzpicture}[line cap=round,line join=round,>=triangle 45,x=0.75cm,y=0.75cm]
\clip(-1.5,-2) rectangle (7,1.5);
\draw [line width=1pt] (-1,1)-- (0.5,-0.5);
\draw [line width=1pt ,red] (0,-0.5)-- (1.5,1);
\draw [line width=1pt] (1,1)-- (2.5,-0.5);
\draw [line width=1pt, red] (1.9949545971809697,-0.3144062445991428)-- (3.4968674700510833,-0.3144062445991428);
\draw [line width=1pt] (3,-0.5)-- (4.5,1);
\draw [line width=1pt, red] (4,1)-- (5.5,-0.5);
\draw [line width=1pt] (5,-0.5)-- (6.5,1);
\draw [line width=1pt] (2.738643205403547,0.9783762368802792)-- (2.747547123946322,-0.4996742412203729);
\draw[<->] [shift={(2.7,-0.5)},line width=1pt]  plot[domain=3.141592653589793:6.283185307179586,variable=\t]({1*1*cos(\t r)+0*1*sin(\t r)},{0*1*cos(\t r)+1*1*sin(\t r)});
\begin{scriptsize}
\draw [fill=red] (2.739984431439439,0.7557327149223082) circle (2.5pt);
\draw [fill=red] (-0.75,0.75) circle (2.5pt);
\draw [fill=red] (6.206763099570006,0.7067630995700059) circle (2.5pt);
\end{scriptsize}
\end{tikzpicture}
\end{center}
\end{minipage} \\ \hline
$(3\text{--}2)^*$  &
\begin{minipage}{70mm}
\begin{center}
\begin{tikzpicture}[line cap=round,line join=round,>=triangle 45,x=1cm,y=1cm]
\clip(-3,-1) rectangle (3,1);
\draw [line width=1pt, red] (-1.5,-1)-- (-1.5,1);
\draw [line width=1pt, red] (2.5,0)-- (-2.5,0);
\draw[line width=1pt,color=black,smooth,samples=100,domain=1:2,red] plot(\x,{(\x - 1)^2});
\draw[line width=1pt,color=black,smooth,samples=100,domain=1:2,red] plot(\x,{0-(\x - 1)^2});
\begin{scriptsize}
\draw (2.75,0) node {$\Gamma_0$};
\draw (-1.25,-0.25) node {$B_1$};
\draw (1.25,-0.5) node {$B_2$};
\end{scriptsize}
\end{tikzpicture}
\end{center}
\end{minipage} &
\begin{minipage}{70mm}
\begin{center}
\begin{tikzpicture}[line cap=round,line join=round,>=triangle 45,x=0.5cm,y=1cm]
\clip(-3,-0.5) rectangle (8.5,2);
\draw [line width=1pt] (-2.9748745370407925,1.48726981236248)-- (0,0);
\draw [line width=1pt, red] (-1,0)-- (1,1);
\draw [line width=1pt] (0,1)-- (2,0);
\draw [line width=1pt ,red] (1,0)-- (3,1);
\draw [line width=1pt] (2,1)-- (4,0);
\draw [line width=1pt, red] (3,0)-- (6.0013030000667555,1.4751234286857444);
\draw [line width=1pt] (4,1)-- (6,0);
\draw [line width=1pt] (4.871689318130352,1.4751234286857444)-- (6.815110706408036,0.5398518855771103);
\draw [line width=1pt, red] (5.98915661639002,0.5762910366073168)-- (7.993309923051381,1.6451728001600416);
\begin{scriptsize}
\draw [fill=red] (5.4450029028739655,0.27749854856301703) circle (2.5pt);
\draw [fill=red] (-2,1) circle (2.5pt);
\end{scriptsize}
\end{tikzpicture}
\end{center}
\end{minipage} \\ \hline
\end{tabular}
\end{table}

\begin{table}[H]
  \begin{tabular}{|c|c|c|} 
    \hline  & branch locus & fiber $F$ and its involution \\ \hline
 $(4\text{--}1)$   &
    \begin{minipage}{70mm}
\begin{center}
\begin{tikzpicture}[line cap=round,line join=round,>=triangle 45,x=1cm,y=1cm]
\clip(-3,-1) rectangle (3,1);
\draw [line width=1pt] (2.5,0)-- (-2.5,0);
\draw[line width=1pt,color=black,smooth,samples=100,domain=-1:1,red] plot(\x,{(\x)^2});
\draw[line width=1pt,color=black,smooth,samples=100,domain=-1:1,red] plot(\x,{0-(\x)^2});
\begin{scriptsize}
\draw (2.75,0) node {$\Gamma_0$};
\draw (1,0.5) node {$B_1$};
\draw (1,-0.5) node {$B_2$};
\end{scriptsize}
\end{tikzpicture}
\end{center}
    \end{minipage} &
    \begin{minipage}{70mm}
\begin{center}
\begin{tikzpicture}[line cap=round,line join=round,>=triangle 45,x=0.6cm,y=0.75cm]
\clip(-6,-2) rectangle (8,2);
\draw [line width=1pt] (-4,1)-- (-4,-1);
\draw [line width=1pt] (-3.810310855199785,0.5901754861508494)-- (-5.01896645566569,1.4850454980342638);
\draw [line width=1pt] (-3.7986891667337668,-0.6301018027810792)-- (-5.007344767199672,-1.5249718146644935);
\draw [line width=1pt] (-4.403016966966719,-0.1768559526063629)-- (-2.9967926625785033,0.32287665143242694);
\draw [line width=1pt] (-3.403551758889144,0.32287665143242694)-- (-1.9973274545009285,-0.16523426414034453);
\draw [line width=1pt] (2.0237767547414074,-0.22334270647043636)-- (3.4300010591296233,0.31125496296640853);
\draw [line width=1pt] (2.999998585886946,0.34612002836446365)-- (4.406222890275162,-0.20009932953839962);
\draw [line width=1pt] (4,1)-- (4,-1);
\draw [line width=1pt] (3.7786517131101727,0.6366622400149229)-- (5.022172378974132,1.4850454980342638);
\draw [line width=1pt] (3.7786517131101727,-0.5719933604509874)-- (4.9873073135760775,-1.4901067492664384);
\draw [line width=1pt] (-1.0327273118214086,0.1834163898402065)-- (0.9894464812657777,0.19503807830622488);
\draw [line width=1pt] (-0.40515613465641975,0.4042284706945555)-- (-1.3929996542679763,-0.22334270647043636);
\draw [line width=1pt] (0.40836205796486214,0.4042284706945555)-- (1.3845838891104003,-0.18847764107238127);
\draw[<->] [shift={(0,-0.5)},line width=1pt]  plot[domain=3.141592653589793:6.283185307179586,variable=\t]({1*1*cos(\t r)+0*1*sin(\t r)},{0*1*cos(\t r)+1*1*sin(\t r)});
\begin{scriptsize}
\draw [fill=red] (-0.4,0.19) circle (2.5pt);
\draw [fill=red] (0.4,0.19) circle (2.5pt);
\draw [color=black] (-1.75,-0.2) node {$\cdots$};
\draw [color=black] (1.75,-0.2) node {$\cdots$};
\end{scriptsize}
\end{tikzpicture}
\end{center}
 \end{minipage} \\ \hline
 $(4\text{--}2)$ &
\begin{minipage}{70mm}
\begin{center}
\begin{tikzpicture}[line cap=round,line join=round,>=triangle 45,x=1cm,y=1cm]
\clip(-3,-1) rectangle (3,1);
\draw [line width=1pt] (2.5,0)-- (-2.5,0);
\draw[line width=1pt,color=black,smooth,samples=100,domain=0:1,red] plot(\x,{(\x)^2});
\draw[line width=1pt,color=black,smooth,samples=100,domain=0:1,red] plot(\x,{0-(\x)^2});
\begin{scriptsize}
\draw (2.75,0) node {$\Gamma_0$};
\draw (0,-0.25) node {$B_1$};
\end{scriptsize}
\end{tikzpicture}
\end{center}
\end{minipage} &
\begin{minipage}{70mm}
\begin{center}
\begin{tikzpicture}[line cap=round,line join=round,>=triangle 45,x=0.6cm,y=0.75cm]
\clip(-6,-2) rectangle (7,2);
\draw [line width=1pt] (-4,1)-- (-4,-1);
\draw [line width=1pt] (-3.810310855199785,0.5901754861508494)-- (-5.01896645566569,1.4850454980342638);
\draw [line width=1pt] (-3.7986891667337668,-0.6301018027810792)-- (-5.007344767199672,-1.5249718146644935);
\draw [line width=1pt] (-4.403016966966719,-0.1768559526063629)-- (-2.9967926625785033,0.32287665143242694);
\draw [line width=1pt] (-3.403551758889144,0.32287665143242694)-- (-1.9973274545009285,-0.16523426414034453);
\draw [line width=1pt] (2.000533377809371,-0.20009932953839962)-- (3.4300010591296233,0.31125496296640853);
\draw [line width=1pt] (2.999998585886946,0.34612002836446365)-- (4.406222890275162,-0.20009932953839962);
\draw [line width=1pt] (4,1)-- (4,-1);
\draw [line width=1pt] (3.7786517131101727,0.6366622400149229)-- (5.022172378974132,1.4850454980342638);
\draw [line width=1pt] (3.7786517131101727,-0.5719933604509874)-- (4.9873073135760775,-1.4901067492664384);
\draw [line width=1pt] (-1.009483934889372,-0.1768559526063629)-- (0.39674036949884384,0.3926067822285371);
\draw [line width=1pt] (-0.40515613465641975,0.4158501591605739)-- (0.9894464812657777,-0.1768559526063629);
\draw[<->] [shift={(0,-0.5)},line width=1pt]  plot[domain=3.141592653589793:6.283185307179586,variable=\t]({1*1*cos(\t r)+0*1*sin(\t r)},{0*1*cos(\t r)+1*1*sin(\t r)});
\begin{scriptsize}
\draw [fill=red] (0,0.25) circle (2.5pt);
\draw [color=black] (-1.5,-0.2) node {$\cdots$};
\draw [color=black] (1.5,-0.2) node {$\cdots$};
\end{scriptsize}
\end{tikzpicture}
\end{center}
\end{minipage} \\ \hline
$(4\text{--}1)^*$ & 
\begin{minipage}{70mm}
\begin{center}
\begin{tikzpicture}[line cap=round,line join=round,>=triangle 45,x=1cm,y=1cm]
\clip(-3,-1) rectangle (3,1);
\draw [line width=1pt, red] (2.5,0)-- (-2.5,0);
\draw[line width=1pt,color=black,smooth,samples=100,domain=-1:1,red] plot(\x,{(\x)^2});
\draw[line width=1pt,color=black,smooth,samples=100,domain=-1:1,red] plot(\x,{0-(\x)^2});
\begin{scriptsize}
\draw (2.75,0) node {$\Gamma_0$};
\draw (1,0.5) node {$B_1$};
\draw (1,-0.5) node {$B_2$};
\end{scriptsize}
\end{tikzpicture}
\end{center}
\end{minipage} &
\begin{minipage}{70mm}
\begin{center}
\begin{tikzpicture}[line cap=round,line join=round,>=triangle 45,x=0.5cm,y=0.375cm]
\clip(-3,-3.5) rectangle (3,3.5);
\draw [line width=1pt] (-1.396346686723849,2.9902261522629088)-- (1.3812368566545277,2.9902261522629088);
\draw [line width=1pt] (-0.7920188864908967,3.199416544651239)-- (-2.1982431908791127,1.8048139287290361);
\draw [line width=1pt] (0.7885307448875938,3.199416544651239)-- (2.171511672343773,1.8396789941270912);
\draw [line width=1pt] (-2.000674486956801,2.385898352029954)-- (-2.000674486956801,0.631023393661182);
\draw [line width=1pt] (1.9971863453534984,2.385898352029954)-- (1.9971863453534984,0.6077800167291453);
\draw [line width=1pt] (-2.000674486956801,-0.6124972722027823)-- (-1.9774311100247646,-2.402237295969609);
\draw [line width=1pt] (-2.209864879345131,-1.7979094957366546)-- (-0.8152622634229333,-3.204133800124876);
\draw [line width=1pt] (-1.396346686723849,-3.006565096202564)-- (1.4044802335865643,-3.006565096202564);
\draw [line width=1pt] (0.8001524333536121,-3.1925121116588575)-- (2.2063767377418277,-1.7979094957366546);
\draw [line width=1pt] (2.0088080338195167,-2.4138589844356275)-- (2.0088080338195167,-0.6241189606688006);
\draw[<->] [line width=1pt] (-1.3847249982578307,0)-- (1.392858545120546,0);
\begin{scriptsize}
\draw [fill=red] (0.5,3) circle (2.5pt);
\draw [fill=red] (-0.5,3) circle (2.5pt);
\draw [color=red] (0.5,-3) circle (2.5pt);
\draw [color=red] (-0.5,-3) circle (2.5pt);
\draw [color=black] (-2,0) node {$\vdots$};
\draw [color=black] (2,0) node {$\vdots$};
\end{scriptsize}
\end{tikzpicture}
\end{center}
\end{minipage} \\ \hline
$(4\text{--}2)^{*}$  &
\begin{minipage}{70mm}
\begin{center}
\begin{tikzpicture}[line cap=round,line join=round,>=triangle 45,x=1cm,y=1cm]
\clip(-3,-1) rectangle (3,1);
\draw [line width=1pt, red] (2.5,0)-- (-2.5,0);
\draw[line width=1pt,color=black,smooth,samples=100,domain=0:1,red] plot(\x,{(\x)^2});
\draw[line width=1pt,color=black,smooth,samples=100,domain=0:1,red] plot(\x,{0-(\x)^2});
\begin{scriptsize}
\draw (2.75,0) node {$\Gamma_0$};
\draw (0,-0.25) node {$B_1$};
\end{scriptsize}
\end{tikzpicture}
\end{center}
\end{minipage} &
\begin{minipage}{70mm}
\begin{center}
\begin{tikzpicture}[line cap=round,line join=round,>=triangle 45,x=0.5cm,y=0.375cm]
\clip(-3,-3.5) rectangle (3,3.5);
\draw [line width=1pt] (0.4034053251062811,3.192327577892131)-- (-2.1982431908791127,1.8048139287290361);
\draw [line width=1pt] (-0.4027259424068129,3.2029345682541455)-- (2.171511672343773,1.8396789941270912);
\draw [line width=1pt] (-2.000674486956801,2.385898352029954)-- (-2.000674486956801,0.631023393661182);
\draw [line width=1pt] (1.9971863453534984,2.385898352029954)-- (1.9971863453534984,0.6077800167291453);
\draw [line width=1pt] (-2.000674486956801,-0.6124972722027823)-- (-1.9774311100247646,-2.402237295969609);
\draw [line width=1pt] (-2.209864879345131,-1.7979094957366546)-- (-0.8152622634229333,-3.204133800124876);
\draw [line width=1pt] (-1.396346686723849,-3.006565096202564)-- (1.4044802335865643,-3.006565096202564);
\draw [line width=1pt] (0.8001524333536121,-3.1925121116588575)-- (2.2063767377418277,-1.7979094957366546);
\draw [line width=1pt] (2.0088080338195167,-2.4138589844356275)-- (2.0088080338195167,-0.6241189606688006);
\draw[<->] [line width=1pt] (-1.3847249982578307,0)-- (1.392858545120546,0);
\begin{scriptsize}
\draw [fill=red] (0,3) circle (2.5pt);
\draw [color=red] (0.5,-3) circle (2.5pt);
\draw [color=red] (-0.5,-3) circle (2.5pt);
\draw [color=black] (-2,0) node {$\vdots$};
\draw [color=black] (2,0) node {$\vdots$};
\end{scriptsize}
\end{tikzpicture}
\end{center}
\end{minipage} \\ \hline
\end{tabular}
\end{table}

\end{prop}

\section{Smoothability of cusp singularities}\label{app:cusp}

This appendix is devoted to giving a necessary and sufficient condition for certain cusp singularities to be smoothable.

We begin by presenting two theorems about smoothable cusp singularities that play a crucial role in this appendix.

\begin{thm}[\cite{Ste}]\label{thm:Steenbrink}
    Let $(p\in X)$ be a cusp singularity whose exceptional curves form a cycle $[b_1,\dots,b_n]^\circ$.
    If $(p\in X)$ is smoothable, then it holds that
    \[
    n + 9 + \sum_{i=1}^n (2-b_i) \geq 0.
    \]
\end{thm}

\begin{thm}[{\cite{Loo}, \cite[Theorem 7.13]{GHK}, \cite{Engel}}]\label{thm:Looijenga_conjecture}
    Let $(p\in X)$ be a cusp singularity whose exceptional curves form a cycle 
    \[
    [a_1+3,2^{a_2},a_3+3,2^{a_4},\dots,a_{2n-1}+3,2^{a_{2n}}]^\circ \qquad (a_i\geq0).
    \]
    Then, $(p\in X)$ is smoothable if and only if the dual cycle 
    \[
    [2^{a_1},a_2+3,2^{a_3},a_4+3,\dots,2^{a_{2n-1}},a_{2n}+3]^\circ
    \]
    can be realized as an anticanonical divisor on a smooth rational surface. 
\end{thm}

We examine whether the cusp singularity given by one of the following cycle is smoothable: 
\begin{itemize}
    \item $[\beta+3,2^{n-1}]^\circ$ for $n\geq1$, $\beta\geq0$ (Corollary \ref{cor:cusp1})
    \item $[\chi-1,2^{n-\gamma-1},3,2^{\gamma-1}]^\circ$ for $0\leq\gamma\leq n$, $\chi\geq4$ (Theorem \ref{thm:cusp2})
    \item $[\chi-1,2^{n-\gamma-1},\chi-1,2^{\gamma-1}]^\circ$ for $0\leq\gamma\leq n$, $\chi\geq4$ (Corollary \ref{cor:cusp3})
    \item $[\chi-1,2^{k_1-1},\chi-1,2^{k_2-1}]^\circ$ for $\chi\geq4$, $k_1\geq0$, $k_2\geq0$, $k_1+k_2>0$ (Theorem \ref{thm:cusp4})
\end{itemize}
Note that the first cycle appears in the amulet $[\mathrm{I}_n \mathrm{E}]_\beta$ (Lemma \ref{modifyI}), and the second and the third cycles appear in the graph $[\mathrm{I}_n\mathrm{E2}]$ or the graph $[\mathrm{I}_n\mathrm{E1}]_{\beta_1,\beta_2,\gamma}$ with suitable $\beta_1$ and $\beta_2$ (Section \ref{subsec:0,2,0}).
The last cycle naturally appears as the exceptional curves of the singularity obtained by taking the double cover of a cone singularity of type $(\chi-3;k_1,k_2)$; see Proposition \ref{prop:classification_cone_sing} (5).

Before proceeding, we introduce the following lemma.
We will use this without explicit reference.

\begin{lem}
    Let $Y$ be a smooth rational surface, and let $D\in|-K_Y|$ be a cycle of rational curves.
    Let $Y'\to Y$ be a blow-up at $p\in D$.
    \begin{itemize}
        \item[$(1)$]
        If $p\in D$ is a smooth point, then the proper transform of $D$ is an anticanonical divisor of $Y'$.
        \item[$(2)$]
        If $p\in D$ is a node, then the reduced total transform of $D$ is an anticanonical divisor of $Y'$.
    \end{itemize}
\end{lem}

\begin{thm}\label{thm:cusp2}
    The cusp singularity $(p\in X)$ given by $[\chi-1,2^{n-\gamma-1},3,2^{\gamma-1}]^\circ$ ($0\leq\gamma\leq n$, $\chi\geq4$) is smoothable if and only if $n \geq \chi - 11$.
\end{thm}

\begin{proof}
    If the singularity $(p\in X)$ is smoothable, Theorem \ref{thm:Steenbrink} shows $n\geq\chi-11$.

    Henceforth, we assume $n\geq\chi-11$ and prove that $(p\in X)$ is smoothable, which, according to Theorem \ref{thm:Looijenga_conjecture}, is equivalent to proving that $[2^{\chi-4},n-\gamma+2,\gamma+2]^\circ$ arises as a divisor $D\in|-K_Y|$ where $Y$ is a smooth rational surface.

    If $n = \chi - 11$, we can construct $(Y,D)$ with 
    \[
    D = [2^{\chi-4},n-\gamma+2,\gamma+2]^\circ = [2^{\chi-4},\chi-\gamma-9,\gamma+2]^\circ
    \]
    by blow-up $\Sigma_{\gamma+1}$ with an anticanonical divisor $D_0=S+\Delta_0\in|-K_{\Sigma_{\gamma+1}}|$, which forms the cycle $[-\gamma-5,\gamma+1]^\circ$, as follows.
    We first take a blow-up at a node of $D_0$, and obtain the cycle $D_1=[-\gamma-4,1,\gamma+2]^\circ$.
    We take the $(n+6)$-times blow-ups at the intersection points of the strict ransforms of $S$ and the unique $(-1)$-curve in the anticanonical divisors.
    This gives the anticanonical divisor of the form $D'=[2^{\chi-5},1,\chi-\gamma-9,\gamma+2]^\circ$.
    Finally, we obtain $[2^{\chi-4},\chi-\gamma-9,\gamma+2]^\circ$ as the anticanonical divisor by taking a blow-up at any point on the $(-1)$-curve in $D'$ that is not a node of $D'$.
    
    Suppose that $n>\chi-11$.
    The above argument shows that the cycle $D=[2^{\chi-4},\chi-\gamma-9,\gamma+2]^\circ$ can be realized as an anticanonical divisor of a smooth rational surface.
    Let $E$ denote the irreducible component of $D$ whose self-intersection number is $-\chi+\gamma+9$.
    By taking the $(n-\chi+11)$-times blow-ups at general $(n-\chi+11)$ points on $E$, we obtain the anticanonical divisor of the form $[2^{\chi-4},n-\gamma+2,\gamma+2]^\circ$.
\end{proof}

\begin{cor}\label{cor:cusp1}
    The cusp singularity given by $[\beta+3,2^{n-1}]^\circ$ ($n\geq1$, $\beta\geq0$) is smoothable if and only if $\beta \leq n+8$.
\end{cor}

From now on, we focus on Theorem \ref{thm:cusp4}.
To prove it, we prepare some lemmas.

\begin{lem}\label{lem:useful1}
Let $Y$ be a rational surface.
Let $D \in |-K_{Y}|$ be a cycle of the form $[\alpha_1,2^{\beta_1},\alpha_2,2^{\beta_2}]^\circ$, where $\alpha_1 + \alpha_2 + 3 = \beta_1 + \beta_2$, $\alpha_i \ge 3$, and $\beta_i \ge 3$.
Then, any $(-1)$-curve on $Y$ does not meet the components corresponding to $\alpha_1$ and $\alpha_2$.
\end{lem}

\begin{proof}
We prove this by contradiction.  
We may assume that there exists a $(-1)$-curve $E$ on $Y$ that meets the component corresponding to $\alpha_1$.
To simplify the argument, we assume $\alpha_1 \ge 4$.
Let $Y \to Y'$ be the contraction of the $(-1)$-curve $E$.  
Then $Y'$ has a anticanonical divisor whose dual graph is of type $[\alpha_1 -1, 2^{\beta_1}, \alpha_2, 2^{\beta_2}]^\circ$.  
By Theorem~\ref{thm:Looijenga_conjecture},  
the cusp singularity of type $[\beta_1+3,2^{\alpha_1-4},\beta_2+3,2^{\alpha_2-3}]^\circ$ is smoothable.  
Therefore, by Theorem~\ref{thm:Steenbrink},  
we have $\alpha_1 + \alpha_2 + 2 \ge \beta_1 + \beta_2$.  
This contradicts the assumption $\alpha_1 + \alpha_2 + 3 = \beta_1 + \beta_2$.
\end{proof}

\begin{lem}\label{lem:useful2}
Let $Y$ be a rational surface, and let $D \in |-K_{Y}|$ be a cycle of the form $[\alpha_1,2^{\beta_1},\alpha_2,2^{\beta_2}]^\circ$, where $\alpha_i \ge 3$ and $\beta_i \ge 3$ ($i=1,2$).
Let $D'$ be the reduced subdivisor of $D$ consisting of the two curves corresponding to $\alpha_1$ and $\alpha_2$, together with the four $(-2)$-curves adjacent to them.
Then, there exists a $(-1)$-curve on $Y$ that meets $D'$.
\end{lem}

\begin{proof}
    First, we observe that there exists at least one $(-1)$-curve on $Y$ because $Y$ is rational and the number of components of $D$ exceeds $4$.
    Assume to the contrary that there is no $(-1)$-curve meeting $D'$.
    Then we can take a $(-1)$-curve $E$ meeting $D-D'$.
    By contracting $E$ and resulting $(-1)$-curves in the pushforwards of $D$, twice, we obtain a rational surface $Y_1$ with the anticanonical divisor $D_1$ that forms a cycle and contains a $0$-curve $\Gamma$.
    Then, by \cite[V, Proposition 4.3]{BHPV}, the linear system $|\Gamma|$ induces a ruling $Y_1\to\PP^1$.
    Again, there exists a $(-1)$-curve on $Y_1$ away from $\Gamma$ since $Y_1$ is rational and the number of components of $D_1$ exceeds $4$.
    By carrying out a similar procedure, we see that a relatively minimal model of $Y_1$ cannot be a Hirzebruch surface.
    This is a contradiction.
\end{proof}

\begin{thm}\label{thm:cusp4}
    The cusp singularity $(p\in X)$ given by $[\chi-1,2^{k_1-1},\chi-1,2^{k_2-1}]^\circ$ ($\chi\geq4$, $k_1\geq0$, $k_2\geq0$, $k_1+k_2>0$) is smoothable if and only if $k_1+k_2\geq2\chi-14$, $(k_1,k_2)=(0,2\chi-15)$, $(k_1,k_2)=(\chi-10,\chi-5)$,  $(k_1,k_2)=(\chi-9,\chi-6)$, $(k_1,k_2)=(\chi-6,\chi-9)$, $(k_1,k_2)=(\chi-5,\chi-10)$, or $(k_1,k_2)=(2\chi-15,0)$.
\end{thm}

\begin{proof}
    Since the cycle associated to $p$ is symmetric in $k_1$ and $k_2$, we may assume that $k_1\leq k_2$.
    Furthermore, we may assume that $k_1>0$, as otherwise one sees that $p\in X$ is smoothable if and only if $k_2\geq2\chi-15$ by Corollary \ref{cor:cusp1}.
    
    We demonstrate constructions of the dual cycle $D=[2^{\chi-4},k_1+2,2^{\chi-4},k_2+2]^\circ$ as an anticanonical divisor on a smooth rational surface separately for the cases
    \begin{itemize}
        \item $k_1+k_2\geq2\chi-14$,
        \item $(k_1,k_2)=(\chi-10,\chi-5)$,
        \item $(k_1,k_2)=(\chi-9,\chi-6)$.
    \end{itemize}
    
    Suppose $k_1+k_2\geq2\chi-14$.
    Let $D_0=S+\Delta_0$ be an anticanonical divisor on $\Sigma_{k_1}$ that forms a cycle $[-k_1-4,k_1]^\circ$.
    For each node $q$ of $D_0$, we take the blow-up at $q$ and the $(\chi-5)$-times blow-ups at the node over $q$ that is on the proper transforms of $S$.
    This gives the anticanonical divisor of the form $D'=[1,2^{\chi-5},k_1+2,2^{\chi-5},1,2\chi-k_1-12]^\circ$.
    Note that $2\chi-k_1-12\leq k_2+2$ from the assumption.
    Finally, by taking blow-ups at smooth points on $D'$ properly, we obtain the anticanonical divisor of the form $D$.
    
    Suppose $(k_1,k_2)=(\chi-10,\chi-5)$.
    Let $D_0=S+\Delta_0$ be an anticanonical divisor on $\Sigma_1$ that forms a cycle $[-5,1]^\circ$.
    Choose a node $q$ on $D_0$, and take the blow-up at $q$ and the $(\chi-4)$-times blow-ups at the nodes over $q$ that is on the proper transform of $S$. 
    Then, we obtain the anticanonical divisor $D'$ of the form $[\chi-9,2^{\chi-4},1]^\circ$.
    Next, we take the blow-up at the node $q'$, the intersection point of the proper transform of $S$ and the exceptional $(-1)$-curve $E$, and the $(\chi-5)$-times blow-ups at the nodes over $q'$ that is on the proper transform of $E$.
    This construction produces the anticanonical divisor of the form $[\chi-8,2^{\chi-4},\chi-3,1,2^{\chi-5}]^\circ$.
    Finally, by taking the blow-up at a smooth point on each $(-1)$-curve, we obtain the anticanonical divisor of the form $D$.

    Suppose $(k_1,k_2)=(\chi-9,\chi-6)$.
    Let $D_0=L+Q$ be an anticanonical divisor on $\PP^2$ where $L$ is a line and $Q$ is a conic intersecting with $L$ transeversely.
    Note that $D_0$ forms a cycle $[-4,-1]^\circ$.
    Choose a node $q$ on $D_0$, and take the blow-up at $q$ and the $(\chi-4)$-times blow-ups at the nodes over $q$ that is on the proper transform of $Q$.
    Then, we obtain the anticanonical divisor $D'$ of the form $[\chi-8,1,2^{\chi-5},0]^\circ$.
    Next, we take the blow-up at the node $q'$, the intersection point of the proper transforms of $L$ and $Q$, and the $(\chi-5)$-times blow-ups at the nodes over $q'$ that is on the proper transform of $L$.
    This construction produces the anticanonical divisor of the form $[\chi-7,1,2^{\chi-5},\chi-4,1,2^{\chi-5}]^\circ$.
    Finally, by taking the blow-up at a smooth point on each $(-1)$-curve, we obtain the anticanonical divisor of the form $D$. 
    
    Theorem \ref{thm:Steenbrink} shows that $k_1+k_2\geq2\chi-15$ whenever $(p\in X)$ is smoothable.
    We still need to show that a divisor of the form $D$ cannot sit on a smooth rational surface as an anticanonical divisor if $k_1+k_2=2\chi-15$, $0<k_1\leq k_2$ and $k_1\neq\chi-10,\chi-9$.
    
    To derive a contradiction, assume that there exists a rational surface $Y$ with an anticanonical divisor of the form $D$ for some $(k_1,k_2)$ such that $k_1+k_2=2\chi-15$, $0<k_1<k_2$ and $k_1\neq\chi-10,\chi-9$.
    By the assumption $k_1 + k_2 = 2\chi - 15$, the divisor $D$ is of the form $[2^{\chi-4}, k_2 + 2, 2^{\chi-4}, 2\chi - k_2 - 13]^\circ$.
    From Lemmas~\ref{lem:useful1} and~\ref{lem:useful2}, there exists a $(-1)$-curve on $Y$ that meets one of the $(-2)$-curves adjacent to the two components corresponding to $k_2 + 2$ and $2\chi - k_2 - 13$.

We assume that the $(-1)$-curve on $Y$ meets one of the $(-2)$-curves adjacent to the component corresponding to $2\chi - k_2 - 13$.
We show that this leads to a contradiction unless $k_2 = \chi - 6$.
We also remark that if, instead, we assume the $(-1)$-curve on $Y$ meets one of the $(-2)$-curves adjacent to the component corresponding to $k_2 + 2$, then a similar argument leads to a contradiction except when $k_2 = \chi - 5$.

Let $Y'$ be the surface obtained by first contracting this $(-1)$-curve, and then successively blowing down the $(-1)$-curves that newly appear, $(2\chi - k_2 - 13)$ times.
Then $Y'$ has an anticanonical divisor whose dual graph is of type $[2^{\chi-5}, k_2+2, 2^{k_2 - \chi + 8}, 1, 0, 2]^\circ$.  
We denote by $S_1$, $S_2$, and $\Gamma'$ the components in the cycle $[2^{\chi-5}, k_2+2, 2^{k_2-\chi+8}, 1, 0, 2]^\circ$ corresponding to the rightmost $2$, $1$, and $0$, respectively.
Then, by \cite[V, Proposition 4.3]{BHPV}, the linear system $|\Gamma'|$ induces a ruling $p': Y' \to \PP^1$.  
The curves $S_1$ and $S_2$ are sections of the ruling $p'$.  
The subchain $[2^{\chi-5}, k_2+2, 2^{k_2 - \chi + 8}]$ in the cycle $[2^{\chi-5}, k_2+2, 2^{k_2 - \chi + 8}, 1, 0, 2]^\circ$ is vertical with respect to $p'$.
Since $Y' \to \PP^1$ is not relatively minimal, there exists a $p'$-vertical $(-1)$-curve on $Y'$ that meets the chain $[2^{\chi-5}, k_2+2, 2^{k_2 - \chi + 8}]$.
By Lemma~\ref{lem:useful1}, this $(-1)$-curve does not meet the component corresponding to $k_2 + 2$ in the chain.
Since the chain $[2^{\chi-5}, k_2+2, 2^{k_2 - \chi + 8}]$ has length at least $4$, the $(-1)$-curve can only meet either one of the two terminal $(-2)$-curves, or one of the two $(-2)$-curves adjacent to the component corresponding to $k_2 + 2$.
By examining which of the above $(-2)$-curves the $(-1)$-curve intersects,  
we carry out a case-by-case analysis and obtain a contradiction in all cases except when $k_2 = \chi - 6$.
\end{proof}

\begin{cor}\label{cor:cusp3}
    The cusp singularity given by $[\chi-1,2^{n-\gamma-1},\chi-1,2^{\gamma-1}]^\circ$ $(0\leq\gamma\leq n,\chi\geq4)$ is smoothable if and only if $n\geq 2\chi - 14$, $(n,\gamma)=(2\chi-15,2\chi-15)$, $(2\chi-15,\chi-5)$, $(2\chi-15,\chi-6)$, $(2\chi-15,\chi-9)$, $(2\chi-15,\chi-10)$ or $(2\chi-15,0)$.
\end{cor}

\section{Drills on extended T-chains}\label{app:drill}

In this appendix, we summarize lists of certain P-admissible chains used in Sections \ref{sec:non-std_Horikawa} and \ref{sec--deformation-arekore-honke}.

We introduce the following operations for a chain.

\begin{defn}\label{def:chain_blowup}
\phantom{A}
\begin{itemize}
    \item[$(1)$] For a chain $C = [c_1, \cdots, c_r]$ with $c_i \ge 2$, the operation that associates the chain
    \[
    [c_1, \dots, c_i]R_1-1-L_1[c_{i+1},\cdots, c_r]
    \]
    is called the \emph{blow-up between $c_i$ and $c_{i+1}$}.

    \item[$(2)$] For a chain $T = [a_1, \cdots, a_r]-1-[b_1, \cdots, b_s]$ with $a_i, b_j \ge 2$, we define the following operations that associate new chains:
    \begin{itemize}
        \item[(L)] $[a_1, \cdots, a_r]R_1-1-L_2[b_1, \cdots, b_s]$
        \item[(R)] $[a_1, \cdots, a_r]R_2-1-L_1[b_1, \cdots, b_s]$
    \end{itemize}
    We refer to operation (L) as a \emph{left blow-up} and (R) as a \emph{right blow-up}.
    The above three blow-up operations are defined in the same way for chain of the form
    \[
    [b_{1,1},\ldots,b_{1,r_{1}}]-1-\cdots-1-[b_{l,1},\ldots,b_{l,r_{l}}].
    \]
\end{itemize}
\end{defn}

By the very definition, a P-admissible chain is a chain that becomes an ample T-train after repeated applications of the three operations described above.

We state some useful lemmas.
The following lemma provides a criterion for checking whether a given T-train is ample.

\begin{lem}\label{lem_ampleness}
Let $C_1-1-C_2$ be a T-train.
Then the following are equivalent.
\begin{itemize}
    \item $C_1-1-C_2$ satisfies the ample condition.
    \item $L_2C_1R_1-1-L_2C_2R_1$ satisfies the ample condition.
    \item $L_1C_1R_2-1-L_1C_2R_2$ satisfies the ample condition.
\end{itemize}
\end{lem}

\begin{proof}
By symmetry, it suffices to verify the equivalence between the first and the second conditions.

Let $\alpha_1$ (resp. $\alpha_2$) denote the log discrepancy of the rightmost component of $C_1$ (resp. the leftmost component of $C_2$).
By definition, the chain $C_1-1-C_2$ satisfies the ample condition if and only if $\alpha_1+\alpha_2<1$.
On the other hand, the chain $L_2C_1R_1-1-L_2C_2R_1$ satisfies the ample condition if and only if
\[
\frac{\alpha_1}{1+\alpha_1} + \frac{1}{2-\alpha_2} < 1.
\]
The equivalence of the two inequalities can be checked by direct computation.
\end{proof}

As a corollary, we see that the ample condition on T-trains imposes some condition on cores and blow-up processes.

\begin{cor}\label{cor:leftmost/rightmost_core}
    Let $\widetilde{C}$ be an ample T-train associated with a P-admissible chain $C$.
    Then, any core is not contracted under the blow-down process $\widetilde{C}\to C$.
    Moreover, the leftmost/rightmost core is greater than $2$ in $C$.
\end{cor}

\begin{proof}
    If $\widetilde{C}$ consists of exactly two vehicles, the claims follow from Lemma~\ref{lem_ampleness} by a simple computation.
    If $\widetilde{C}$ has more than two vehicles, applying the claims to all pairs of adjacent vehicles in $\widetilde{C}$ yields the desired result.
\end{proof}

\begin{rem}
    The above proof also works when $\widetilde{C}$ is a \emph{nef} T-train.
    In fact, Figueroa-Rana-Urz\'{u}a \cite[Corollary 2.4]{FRU} showed the first half of Corollary~\ref{cor:leftmost/rightmost_core} in the nef setting.
\end{rem}

The following corollary follows immediately from Corollary~\ref{cor:leftmost/rightmost_core}. 

\begin{cor}\label{cor_drill_2}
Let $[2^\beta, a_1 \cdots, a_r]$ be a P-admissible chain for some $\beta \ge 1$.
Then, during the process of transforming into an ample T-train, no blow-up occurs between the components $2^n$ and $2^{\beta - n}$ (for all $1 \le n \le \beta$).
\end{cor}

As we will see later, the following lemma reduces the classification of certain P-admissible chains.

\begin{lem}\label{lem_no_blowup_1}
Let $[2^\alpha, 3,2^\beta, a_1, \cdots, a_r]$ be a P-admissible chain for $\alpha \ge 0, \beta \ge -1$ with an associated ample T-train $\widetilde{C}$.
If a blow-up (in the sense of Definition \ref{def:chain_blowup}~(1)) occurs on the subchain $[2^\alpha, 3,2^\beta, a_1]$, then it occurs exactly once on $[2^\alpha, 3,2^\beta, a_1]$ and the leftmost vehicle of $\widetilde{C}$ becomes the chain $[2^\alpha,3,2^{k-1},3+\alpha]$ for some $k\geq0$ after the blow-ups.
\end{lem}

\begin{proof}
We assume that a blow-up occurs in the subchain $[2^\alpha, 3,2^{\beta}, a_1]$.
By Corollary~\ref{cor_drill_2}, no blow-up can occur within the subchain $[2^\alpha, 3]$, so a blow-up occurs in the subchain $[3,2^{\beta},a_1]$.

Suppose a blow-up occurs between $[2^{k},2^{\beta-k}]$.
Then we obtain
\[
[2^\alpha,3,2^{k-1},3]-1-[3,2^{\beta-k-1},a_1,\cdots, a_r].
\]
Performing $\alpha$-times left blow-ups yields the chain
\[
[2^\alpha,3,2^{k-1},3+\alpha]-1-[2^\alpha,3,2^{\beta-k-1},a_1,\cdots, a_r].
\]
No further blow-up occurs in the chain $[2^\alpha,3,2^{\beta-k-1},a_1]$ of the second chain, as otherwise the second vehicle has a form similar to the first one, contradicting the ample condition.
Hence, $[2^\alpha,3,2^{\beta-k-1},a_1,\cdots, a_r]$ is the desired extended T-chain satisfying the required condition. 
\end{proof}

Thanks to Lemma~\ref{lem_no_blowup_1}, when classifying P-admissible chains of the form 
\[
[2^\alpha, 3, 2^\beta, a_1, \cdots, a_r,2^\gamma,3,2^\delta],
\]
we may assume that no blow-up occurs on the subchains $[2^\alpha, 3,2^{\beta}, a_1]$ and $[a_r,2^{\gamma},3,2^\delta]$.
This observation simplifies computations for the classification.

In what follows, we apply the tools established above to classify certain P-admissible chains.

\begin{lem}\label{lem:extT_2m2}
    The chain $C=[2^{\alpha},\chi,2^{\beta}]$ ($\alpha, \beta \ge 0$) is P-admissible if and only if $C$ is T-chain and either $C=[2^{\chi-4},\chi]$ ($\alpha=0$, $\beta=\chi-4$) or $C=[\chi,2^{\chi-4}]$ ($\alpha=\chi-4$, $\beta=0$).
\end{lem}

\begin{proof}
By the Corollary~\ref{cor_drill_2},
$[2^\alpha,\chi,2^\beta]$ is a T-chain.
Therefore, we get $(\alpha,\beta)=(0,\chi-4)$ or $(\chi-4,0)$.
\end{proof}

\begin{lem} \label{lem:extT_I}
Let $\alpha_i \ge 0$, $\beta_i \ge 0$ and $n \ge 0$.
\begin{itemize}
    \item[$(1)$] The chain $C=[3,2^{\beta_{1}-1},3,2^{n-1},3,2^{\beta_2 -1},3]$ is not P-admissible.
    \item[$(2)$]  The chain $C=[2,4,2^{n-1},4,2]$ is a P-admissible chain and its associated T-chain is
    \[
    [2,5]-1-[3,2^{n-2},3]-1-[5,2].
    \]
     \item[$(3)$]  The chain $C=[2^{\alpha_1},3,2^{\beta_1 -1},3,2^{\alpha_2}]$ is P-admissible
     if and only if $\alpha_1=\alpha_2=0$.
     In this case, $C$ is a T-chain 
     $[3,2^{\beta_1 -1},3]$.
    \item[$(4)$]  The chain  
     $C = [2^{\alpha_1}, 3, 2^{\beta_1 - 1}, 3, 2^{\beta_2 - 1}, 3, 2^{\alpha_2}]$ is P-admissible if and only if $(\alpha_1, \alpha_2) = (1, 0)$ or $(0, 1)$,  
     and the associated ample T-train is $[2,3,2^{\beta_1 - 1}, 4] - 1 - [3, 2^{\beta_2 - 2}, 3]$ in the former case,  
     or $[3, 2^{\beta_1 - 2}, 3] - 1 - [2,3,2^{\beta_2 - 1}, 4]$ in the latter case.
     If $\beta_2 = 0$, we regard $[2,3,2^{\beta_1 - 1}, 4] - 1 - [3, 2^{\beta_2 - 2}, 3]$ as the T-chain $[2,3,2^{\beta_1 - 1}, 4]$.  Similarly, if $\beta_1 = 0$, we regard $[3, 2^{\beta_1 - 2}, 3] - 1 - [2,3,2^{\beta_2 - 1}, 4]$ as the T-chain $[2,3,2^{\beta_2 - 1}, 4]$.
\end{itemize}
\end{lem}
\begin{proof}
We will only show the case $(4)$.
By Lemma~\ref{lem_no_blowup_1},
to classify the P-admissible chain of the form 
\[\
C=[2^{\alpha_1},3,2^{\beta_{1}-1},3,2^{\beta_2 -1},3,2^{\alpha_2}],
\]
we may assume $C$ is a T-chain.
At first, we assume $\beta_1=\beta_2=0$.
In this case, $C=[2^{\alpha_2},5,2^{\alpha_2}]$.
So we get $(\alpha_1,\alpha_2)=(1,0)$ or $(0,1)$.
Therefore,
if $(\alpha_1,\beta_1,\beta_2,\alpha_2)=(1,0,0,0)$ or $(0,0,0,1)$,
then $C$ is a T-chain.
Next, we assume $\beta_1 \ge 1$ and $\beta_2=0$.
In this case, $C=[2^{\alpha_1},3,2^{\beta_1-1},4,2^{\alpha_2}]$.
If $\alpha_1=0$, then $\alpha_2=1$ and $[2^{\beta_1},4]$ must be a T-chain.
But it is impossible by the assumption $\beta_1 \ge 1$.
It follows that $\alpha_1 \ge1$.  
So we get $\alpha_2=0$.
Since $[3,2^{\beta_{1}-1},4-\alpha_1]$ must be a T-chain, we get $\alpha_1=1$.
Therefore, if 
$(\alpha_1,\beta_1,\beta_2,\alpha_2)=(1,\beta_1,0,0)$,
then $C$ is a T-chain.
Similarly, it follows that in the case of $\beta_2 \ge 1$ and $\beta_1=0$, the condition $(\alpha_1,\beta_1,\beta_2,\alpha_2)=(0,0,\beta_2,1)$ implies that $C$ becomes a T-chain.
We assume $\beta_1 \ge 1$ and $\beta_2 \ge 1$.
In this case, 
it must be $\alpha_1=0$ or $\alpha_2=0$.
If $\alpha_1=0$, then it must be $\alpha_2=1$.
So $[3, 2^{\beta_1-1},3,2^{\beta_2 -1},3, 2]$ must be a T-chain.
But it is impossible by the assumption $\beta_1 \ge 1$ and $\beta_2 \ge 1$.
Similarly, it can be seen that the case $\alpha_2 = 0$ is also impossible.
Summarizing the above discussion, $C$ is a T-chain if and only if 
$(\alpha_1,\beta_1,\beta_2,\alpha_2) = (1,\beta_1,0,0)$ or $(0,0,\beta_2,1)$.
\end{proof}

The following lemma can be proved by a similar argument as the previous one.

\begin{lem} \label{lem:extT_CS}
The chain $C=[2^{\alpha},4+\beta,2^{\gamma-1},3,2^{\beta}]$ with $\alpha, \beta, \gamma \ge 0$ is P-admissible if and only if
$(\alpha,\beta,\gamma)=(1,0,n)$ with $n \ge 0$ and its associated ample T-train is $[2,5]-1-[3,2^{n-2},3]$.
\end{lem}

\begin{lem}[Theorem~\ref{thm:one-section-T}] \label{lem:extT_step3}
For each of the chains, the following holds. 
\begin{itemize}
    \item[$(1)$] The chain $C=[2^\alpha,\chi,3+\beta,2^{\gamma-2},3,2^{\beta}]$ ($\alpha \ge0$, $\chi \geq 2$, $\gamma \geq 1$, $\beta \geq 0$) is P-admissible if and only if one of the following holds:
\begin{itemize}
    \item[$\bullet$] $(\alpha,\chi,\beta,\gamma)=(0,\chi,\chi-2,\chi-2+n)$ and $C$ is a P-admissible chain $[\chi,\chi+1,2^{\chi-4+n},3,2^{\chi-2}]$
    and its associated ample T-train is $[\chi,\chi+1,2^{\chi-4},3,2^{\chi-2}]-1-[\chi+1,2^{n-2},3,2^{\chi-2}]$. 
    \item[$\bullet$] $(\alpha,\chi,\beta,\gamma)=(2,5,0,n+1)$, C is a P-admissible chain $C=[2^2,5,3,2^{n-1},3]$ and its associated ample T-train is $[2^2,5,4]-1-[3,2^{n-2},3]$.
    \item[$\bullet$] $(\alpha,\chi,\beta,\gamma)=(1,\chi,0,\chi-2+n)$, $C$ is a P-admissible chain $[2,\chi,3,2^{\chi+n-4},3]$ and its associated ample T-train is $[2,\chi,3,2^{\chi-4},3]-1-[3,2^{n-2},3]$.
    \item[$\bullet$] $(\alpha,\chi,\beta,\gamma)=(\chi-2,\chi,0,n+1)$, $C$ is a P-admissible chain $C=[2^{\chi-2},\chi,3,2^{n-1},3]$ and its associated ample T-train is $[2^{\chi-2},\chi+2]-1-[2,5]-1-[3,2^{n-2},3]$.
\end{itemize}

\item[$(2)$] The chain $[2^\alpha,\chi,\beta,2^{\beta-4}]$ $(\beta=5,6,7,8)$ is P-admissible if and only if it is a T-chain $[3,5,2]$ ($\alpha=0$, $\chi=3$, $\beta=5$).

\item[$(3)$] The chain $C=[2^\alpha,\chi+1,2^{\beta-1},3,2]$ $(\alpha \ge 0, \chi \ge 1, \beta \ge 0)$ is P-admissible if and only if $C$ is a T-chain and  $C=[4,2^{\beta-1},3,2]$ ($\alpha=0$, $\chi=3$, $\forall \beta \ge 0$).

\item[$(4)$]
The chain $C=[2^\alpha,3,2^{\gamma-2},3+\alpha,\chi,3+\beta,2^{\delta-2},3,2^{\beta}]$ ($\alpha \ge0$, $\gamma \ge 1$ $\chi \geq 2$, $\delta \geq 1$, $\beta \geq 0$) is P-admissible if and only if one of the following holds:
\begin{itemize}
    \item[$\bullet$] $(\alpha,\gamma,\chi,\delta,\beta)=(1,m+2,3,n+2,0)$, $C$ is a P-admissible chain \\
    $[2,3,2^m,4,3,3,2^n,3]$ and its associated ample T-train is $[2,3,2^{m-2},4]-1-[2,3,4,3,3,3]-1-[3,2^{n-2},3]$.
    \item[$\bullet$] $(\alpha,\gamma,\chi,\delta,\beta)=(1,m+1,2,n+3,0)$,
      $C$ is a P-admissible chain \\
      $[2,3,2^{m-1},4,2,3,2^{n+1},3]$ 
      and its associated ample T-train is $[2,3,2^{m-2},4]-1-[2,5,2,3,2,3]-1-[3,2^{n-2},3]$. 
    \item[$\bullet$] $(\alpha,\gamma,\chi,\delta,\beta)=(1,m+1,5,n+4,0)$,
      $C$ is a P-admissible chain\\
      $[2,3,2^{m-1},4,5,3,2^{n+2},3]$ 
      and its associated ample T-train is $[2,3,2^{m-2},4]-1-[2,5,5,3,2^2,3]-1-[3,2^{n-2},3]$.
    \item[$\bullet$] $(\alpha,\gamma,\chi,\delta,\beta)=(2,m+\chi-1,\chi,n+1,0)$,
      $C$ is a P-admissible chain \\
      $[2^2,3,2^{m+\chi-3},5,\chi,3,2^{n-1},3]$
      and its associated ample T-train is $[2^2,3,2^{m-2},5]-1-[2^2,3,2^{\chi-3},5,\chi,4]-1-[3,2^{n-2},3]$.
      \item[$\bullet$] 
      $(\alpha,\gamma,\chi,\delta,\beta)=(0,m+4,4,n+1,1)$,
      $C$ is a P-admissible chain \\
      $[3,2^{m+2},3,4,4,2^{n-1},3,2]$ 
      and its associated ample T-train is $[3,2^{m-2},3]-1-[3,2^2,7,2]-1-[3,2^2,5,5,2]-1-[4,2^{n-2},3,2]$. 
      \item[$\bullet$]
      $(\alpha,\gamma,\chi,\delta,\beta)=(\chi,m+\chi,\chi,n+1,0)$,
      $C$ is a P-admissible chain\\
      $[2^\chi,3,2^{\chi-2+m},\chi+3,\chi,3,2^{n-1},3]$ 
      and its associated ample T-train is $[2^\chi,3,2^{m-2},\chi+3]-1-[2^\chi,3,2^{\chi-2},\chi+3,\chi+2]-1-[2,5]-1-[3,2^{n-2},3]$. 
      \item[$\bullet$]
      $(\alpha,\gamma,\chi,\delta,\beta)=(1,m+2,\chi,\chi+n,0)$,
      $C$ is a P-admissible chain \\
      $[2,3,2^m,4,\chi,3,2^{\chi-2+n},3]$ and its associated ample T-train is $[2,3,2^{m-2},4]-1-[2,3,5,3]-1-[2,\chi+2,3,2^{\chi-2},3]-1-[3,2^{n-2},3]$. 
      \item[$\bullet$]
      $(\alpha,\gamma,\chi,\delta,\beta)=(2,\chi+m,\chi,n+1,0)$,
      $C$ is a P-admissible chain \\
      $[2^2,3,2^{m+\chi-2},5,\chi,3,2^{n-1},3]$ and its associated ample T-train is $[2^2,3,2^{m-2},5]-1-[2^2,3,2^{\chi-2},\chi+3,4]-1-[2^2,3,2^{\chi-4},\chi+1,4]-1-[3,2^{n-2},3]$. 
      \item[$\bullet$]
      $(\alpha,\gamma,\chi,\delta,\beta)=(1,m+3,4,n+2,0)$,
      $C$ is a P-admissible chain \\
      $[2,3,2^{m+1},4,4,3,2^n,3]$ and its associated ample T-train is $[2,3,2^{m-2},4]-1-[2,3,2,6,3]-1-[2,3,5,3,3]-1-[3,2^{n-2},3]$. 
      \item[$\bullet$] 
      $(\alpha,\gamma,\chi,\delta,\beta)=(1,m+1,\chi,\chi+n+1,0)$,
      $C$ is a P-admissible chain \\
      $[2,3,2^{m-1},4,\chi,3,2^{\chi-1+n},3]$ and its associated ample T-train is $[2,3,2^{m-2},4]-1-[2,6,2,3]-1-[2,\chi+3,3,2^{\chi-1},3]-1-[3,2^{n-2},3]$. 
      \item[$\bullet$]
      $(\alpha,\gamma,\chi,\delta,\beta)=(0,m+2,\chi,\chi+n,1)$,
      $C$ is a P-admissible chain \\
      $[3,2^m,3,\chi,4,2^{\chi+n-2},3,2]$ and its associated ample T-train is $[3,2^{m-2},3]-1-[3,5,2]-1-[3,\chi+2,2^{\chi-3},3,2]-1-[3,\chi+3,2^{\chi-2},3,2]-1-[4,2^{n-2},3,2]$. 
\end{itemize}
\item[$(5)$] 
The chain $C=[2^{\beta-4},\beta,\chi,3+\gamma,2^{\delta-2},3,2^{\gamma}]$ $(\beta \in \{5,6,7,8\}, \gamma \ge 0, \delta \ge 1)$ is P-admissible if and only if one of the following holds:
\begin{itemize}
      \item[$\bullet$] $(\beta,\gamma,\chi,\delta)=(6,0,2,n+1)$,
      $C$ is a P-admissible chain $[2^2,6,2,3,2^{n-1},3]$
      and its associated ample T-train is $[2^2,6,2,4]-1-[3,2^{n-2},3]$.
      \item[$\bullet$] $(\beta,\gamma,\chi,\delta)=(5,0,2,n+3)$,
      $C$ is a P-admissible chain $[2,5,2,3,2^{n+1},3]$ and its associated ample T-train is $[2,5,2,3,2,3]-1-[3,2^{n-2},3]$. 
      \item[$\bullet$] $(\beta,\gamma,\chi,\delta)=(5,0,5,n+4)$,
      $C$ is a P-admissible chain $[2,5,5,3,2^{n+2},3]$ and its associated ample T-train is $[2,5,5,3,2^2,3]-1-[3,2^{n-2},3]$. 
      \item[$\bullet$] $(\beta,\gamma,\chi,\delta)=(5,0,4,n+4)$,
      $C$ is a P-admissible chain $[2,5,4,3,2^{n+2},3]$ and its associated ample T-train is $[2,5,5,2^2,3]-1-[2,7,2^2,3]-1-[3,2^{n-2},3]$. 
      \item[$\bullet$] $(\beta,\gamma,\chi,\delta)=(5,0,\chi,\chi+n+1)$,
      $C$ is a P-admissible chain $[2,5,\chi,3,2^{\chi+n-1},3]$ and its associated ample T-train is $[2,6,2,3]-1-[2,\chi+3,3,2^{\chi-1},3]-1-[3,2^{n-2},3]$. 
\end{itemize}

\item[$(6)$]
The chain $C=[2,3,2^{\alpha-1},\chi+1,3+\beta,2^{\gamma-2},3,2^{\beta}]$ $(\alpha\in\{0,1\},\chi \ge 1,\gamma \geq 1, \beta \ge 0)$ is P-admissible if and only if one of the following holds:

\begin{itemize}
      \item[$\bullet$] $(\alpha,\beta,\chi,\gamma)=(0,0,2,n+2)$,
      $C$ is a P-admissible chain $[2,4,3,2^n,3]$ and its associated ample T-train is $[2,4,3,3]-1-[3,2^{n-2},3]$.
      \item[$\bullet$] $(\alpha,\beta,\chi,\gamma)=(1,0,2,n+2)$,
      $C$ is a P-admissible chain $[2,3,3,3,2^n,3]$ and its associated ample T-train is $[2,5]-1-[2,4,3,3]-1-[3,2^{n-2},3]$.
      \item[$\bullet$] $(\alpha,\beta,\chi,\gamma)=(1,0,3,n+2)$,
      $C$ is a P-admissible chain $[2,3,4,3,2^n,3]$ and its associated ample T-train is $[2,3,5,3]-1-[2,5,3]-1-[3,2^{n-2},3]$. 
      \item[$\bullet$] $(\alpha,\beta,\chi,\gamma)=(1,0,4,n+2)$,
      $C$ is a P-admissible chain $[2,3,5,3,2^n,3]$ and its associated ample T-train is $[2,3,5,3,3]-1-[3,2^{n-2},3]$. 
\end{itemize}

\item[$(7)$]
The chain $C=[2^{\alpha-4},\alpha,\chi,\beta,2^{\beta-4}]$  is not P-admissible for any $\alpha, \beta \in\{5,6,7,8 \}$ and $\chi \ge 2$.

\item[$(8)$]
The chain $C=[2^{\alpha-4},\alpha,\chi+1,2^{\beta-1},3,2]$  is not P-admissible for any $\alpha \in\{5,6,7,8 \}, \beta \ge 0 $ and $\chi \ge 1$.

\item[$(9)$]
The chain $C=[2,3,2^{\alpha-1},\chi+2,2^{\beta-1},3,2]$ is not P-admissible for any $\alpha, \beta, \chi \ge 0$.

\end{itemize}
\end{lem}
\begin{proof}
We will only show the case $(1)$.
By Lemma~\ref{lem_no_blowup_1},
to classify the P-admissible chain of the form
\[
C=[2^\alpha,\chi,3+\beta,2^{\gamma-2},3,2^{\beta}],
\]
we may assume that blow-up only occur between $\chi$ and $3+\beta$.
First, we consider the case where $C = [2^\alpha, \chi, 3+\beta, 2^{\gamma-2}, 3, 2^{\beta}]$ is already a T-chain.  
When $\gamma=1$, we have $C = [2^\alpha, \chi, 4+\beta, 2^{\beta}]$.  
Therefore, by a simple computation, we find that $C$ is a T-chain if and only if  
$(\alpha, \chi, \beta, \gamma) = (0, 3, 1, 1)$, $(1, 3, 0, 1)$, or $(2, 5, 0, 1)$.  
When $\gamma \ge 2$, we have $C = [2^\alpha, \chi, 3+\beta, 2^{\gamma-2}, 3, 2^{\beta}]$.  
Again, by a simple computation, we find that $C$ is a T-chain if and only if  
$(\alpha, \chi, \beta, \gamma) = (1, \chi, 0, \chi - 2)$ or $(0, \chi, \chi - 2, \chi - 2)$.

Next, we consider the case where a blow-up occurs between $\chi$ and $3+\beta$.
Performing a blow-up yields the chain
\[
[2^{\alpha}, \chi+1] - 1 - [4+\beta, 2^{\gamma-2}, 3, 2^{\beta}].
\]
We must operate some blow-ups until the chain $[2^{\alpha}, \chi+1]$ becomes a T-chain.
By repeating left or right blow-ups, we can transform $[2^{\alpha}, \chi+1]$ into one of the following two types of T-chains:
\begin{itemize}
    \item[(i)] $L_2^{\alpha}[4]R_{1}^{\chi-3}R_{1}^{\alpha-(\chi-3)}$
    \item[(ii)] $L_2^{\alpha}L_1^{\chi-3}[4]R_{1}^{k}R_{2}^{k+\chi-3}R_1^{\alpha}$
\end{itemize}
Here, the former occurs if $\alpha \ge \chi - 3$, and in the latter, $k$ is a non-negative integer.

We consider case $(\mathrm{i})$. For the chain
\[
[2^{\alpha}, \chi+1] - 1 - [4+\beta, 2^{\gamma-2}, 3, 2^{\beta}],
\]
performing $(\alpha - (\chi - 3))$-times left blow-ups yields the chain
\[
[2^{\alpha}, \chi+1] R_{1}^{\alpha - (\chi - 3)} - 1 - L_{2}^{\alpha - (\chi - 3)} [4+\beta, 2^{\gamma-2}, 3, 2^{\beta}].
\]
Since $[2^{\alpha}, \chi+1] R_{1}^{\alpha - (\chi - 3)}$ is already a completed T-chain,  
we must complete
\[
L_{2}^{\alpha - (\chi - 3)} [4+\beta, 2^{\gamma-2}, 3, 2^{\beta}]
= [2^{\alpha - (\chi - 3)}, 4+\beta, 2^{\gamma-2}, 3, 2^{\beta}]
\]
as a P-admissible chain.
There exists no core in $[2^{\alpha - (\chi - 3)}]$  
and that no blow-up occurs within $[4+\beta, 2^{\gamma-2}, 3, 2^{\beta}]$,  
the chain $[2^{\alpha - (\chi - 3)}, 4+\beta, 2^{\gamma-2}, 3, 2^{\beta}]$ must be a T-chain.  
This is equivalent to $(\alpha, \chi, \beta, \gamma) = (\chi - 2, \chi, 0, 1)$.

We consider case ($\mathrm{ii}$). 
To realize this, we first perform $k$-times left blow-ups, then $k+\chi-3$-times right blow-ups, and finally $\alpha$-times left blow-ups.
This gives  
\[
L_2^{\alpha}L_1^{\chi-3}[4]R_{1}^{k}R_{2}^{k+\chi-3}R_1^{\alpha}-1-
L_2^{\alpha}L_1^{k+\chi-3}L_2^{k}[4+\beta,2^{\gamma-2},3,2^{\beta}].
\]  
Since $L_2^{\alpha}L_1^{\chi-3}[4]R_{1}^{k}R_{2}^{k+\chi-3}R_1^{\alpha}$ is already completed T-chain,
\[
L_2^{\alpha}L_1^{k+\chi-3}L_2^{k}[4+\beta,2^{\gamma-2},3,2^{\beta}]
=[2^\alpha,\chi-1+k,2^{k-1},4+\beta,2^{\gamma-2},3,2^{\beta}]
\]
must be P-admissible.

We show that $k = 0$ by contradiction.  
Assume $k \ge 1$.  
There are no core entries in $[2^\alpha, \chi - 1 + k, 2^{k-1}]$by Corollary~\ref{cor:leftmost/rightmost_core}, and no blow-up occurs within $[4 + \beta, 2^{\gamma-2}, 3, 2^{\beta}]$.  
Therefore,  
\[
[2^\alpha, \chi - 1 + k, 2^{k-1}, 4 + \beta, 2^{\gamma - 2}, 3, 2^{\beta}]
\]  
is a T-chain whose core contains $4 + \beta$.  
Thus, there exists a T-chain $T_2$ such that  
\[
[2^\alpha, \chi - 1 + k, 2^{k-1}, 4 + \beta, 2^{\gamma - 2}, 3, 2^{\beta}] =
L_2^{\alpha} L_1^{k + \chi - 3} L_2^k T_2 R_1^k R_2^{k + \chi - 3} R_1^\alpha.
\]  
The chain  
\[
L_2^{\alpha} L_1^{\chi - 3} [4] R_1^k R_2^{k + \chi - 3} R_1^{\alpha} -1- 
L_2^{\alpha} L_1^{k + \chi - 3} L_2^k T_2 R_1^k R_2^{k + \chi - 3} R_1^{\alpha}
\]  
satisfies the ample condition.  
This is equivalent to  
\[
[4] -1- L_2^k T_2 R_1^k
\]  
satisfying the ample condition.  
Hence, it follows that $k = 0$,  
which contradicts the assumption $k \ge 1$.  
Therefore, we conclude that $k = 0$, and that 
$
[2^\alpha,\chi+\beta+1,2^{\gamma - 2}, 3, 2^{\beta}]
$
is a T-chain. 

Since the T-train  
\[
[2^\alpha,\chi+1,2^{\chi-4},2+\alpha] - 1 - [2^\alpha,\chi+\beta+1,2^{\gamma - 2}, 3, 2^{\beta}]
\]    
satisfies the ample condition,
we must have $\alpha = 0$ and, in addition, either $\beta \ge 1$ or $\chi = 2$.
   
Since a blow-up occurs between $\chi$ and $3 + \beta$ in the chain $C$,  
we have $\chi \ge 3$ by Lemma~\ref{cor_drill_2}.  
From $\alpha = 0$, $\beta \ge 1$, and $\chi \ge 3$,  
we find that  
\[
[2^\alpha,\chi+\beta+1,2^{\gamma - 2}, 3, 2^{\beta}]
\]  
is not a T-chain.  
Hence, in case (ii), the chain $C$ cannot be P-admissible.

We summarize the above argument.  
Assuming that blow-ups occur only between $\chi$ and $3+\beta$,  
we find that the chain $C = [2^\alpha, \chi, 3+\beta, 2^{\gamma-2}, 3, 2^{\beta}]$  
is an extended T-chain with the ample condition  
if and only if $(\alpha, \chi, \beta, \gamma)$ is one of the following:
\[
(2,5,0,1),\quad (1,\chi,0,\chi-2),\quad (0,\chi,\chi-2,\chi-2),\quad (\chi-2,\chi,0,1).
\]
If $(\alpha, \chi, \beta, \gamma)=(2,5,0,1), (1,\chi,0,\chi-2), (0,\chi,\chi-2,\chi-2)$, then $C$ is a T-chain $[2^2,5,4]$, $[2^{\chi-2},3,2^{\chi-4},3]$, $[\chi,\chi+1,2^{\chi-4},3,2^{\chi-2}]$, respectively.
If $(\alpha, \chi, \beta, \gamma)= (\chi-2,\chi,0,1)$,
then $C$ is a P-admissible chain and its associated ample T-train is $[2^{\chi-2},\chi+2]-1-[2,5]$.
Therefore, by Lemma~\ref{lem_no_blowup_1},  
we conclude that the chain $C = [2^\alpha, \chi, 3+\beta, 2^{\gamma-2}, 3, 2^{\beta}]$  
is P-admissible  
if and only if $(\alpha, \chi, \beta, \gamma)$ is one of the following:
\[
(2,5,0,n+1),\quad (1,\chi,0,\chi-2+n),\quad (0,\chi,\chi-2,\chi-2+n),\quad (\chi-2,\chi,0,n+1).
\]
It also follows from Lemma~\ref{lem_no_blowup_1} that the associated ample T-train is determined.
\end{proof}

The following can be shown in the same manner as in Lemma~\ref{lem:extT_step3}.

\begin{lem}[Lemma~\ref{lem:amulet-b}]
\label{drill:lem:amulet-b}
Let $a \ge 0$, $b \in \{2,3,4,6\}$.
For each of the chains, the following holds.
\begin{itemize} 
    \item[$(1)$] The chain $C=[2^{b-2},3,2^{a-1},4+c,2^{d-1},3,2^{c}]$ $(c,d \ge 0)$ is P-admissible if and only if one of the following holds:
    \begin{itemize}
        \item[$\bullet$] $(b,a,c,d)=(2,m+1,1,n+1)$ $(m,n \geq 0)$,
        $C$ is a P-admissible chain $[3,2^m,5,2^n,3,2]$
        and its associated ample T-train is $[3,2^{m-2},3]-1-[3,5,3,2]-1-[4,2^{n-2},3,2]$.
        \item[$\bullet$] $(b,a,c,d)=(3,m,0,n+1)$ $(m,n \geq 0)$,
        $C$ is a P-admissible chain \\
        $[2,3,2^{m-1},4,2^n,3]$ and its associated ample T-train is $[2,3,2^{m-2},4]-1-[2,5,3]-1-[3,2^{n-2},3]$.
    \end{itemize}
    \item[$(2)$] The chain $C=[2^{b-2},3,2^{a-1},c+1,2^{c-4}]$ $(c \ge 4)$ is P-admissible if and only if one of the following holds:
     \begin{itemize}
     \item[$\bullet$] $(b,a,c)=(2,m+2,5)$ $(m \geq 0)$,
        $C$ is a P-admissible chain $[3,2^{m+1},6,2]$
        and its associated ample T-train is $[3,2^{m-2},3]-1- [3,2,6,2]$.
     \item[$\bullet$] $ (b,a,c)=(4,0,4)$ and 
        $C$ is a T-chain $[2^{2},6]$.
     \end{itemize}
    \item[$(3)$] The chain $C=[2^{b-2},3,2^{a-1},4,2]$ is P-admissible if and only if 
     $(b,a)=(2,m)$, $C$ is a P-admissible chain $[3,2^{m-1},4,2]$ and its associated ample T-train is $[3,2^{m-2},3]-1-[5,2]$.
     
\end{itemize}
\end{lem}

\begin{lem}\label{lem:drill_anti-P}
    Let $C=[d+1,2^{l-1},3]-1-[2]$ $(d\geq4,l\geq0)$.
    A T-train can be obtained from $C$ by performing an appropriate sequence of left blow-ups and right blow-ups if and only if $l=d-4$.
\end{lem}

\begin{proof}
    If $l=d-4$, there is a T-train induced by $C$:
    \[
    [d+1,2^{l-1},3]R_1^{d-3}R_2^{d-1}-1-L_1^{d-1}L_2^{d-3}[2] = [d+1,2^{d-5},d,2^{d-1}]-1-[d+1,2^{d-3}].
    \]
    One sees that there are countably infinitely many such T-trains if $l=d-4$:
    \[
    [d+1,2^{d-5},d+2,2^{d-3},d-2,2^{d-1}]-1-[d+1,2^{d-5},d,2^{d-1}],
    \]
    \[
    [d+1,2^{d-5},d+2,2^{d-5},d,2^{d-1},d-2,2^{d-1}]-1-[d+1,2^{d-5},d+2,2^{d-3},d-2,2^{d-1}],
    \]
    and so on.
    Conversely, it is straightforward to see that the only T-chains induced by $C$ are these.
\end{proof}

\end{document}